\documentclass[12pt,reqno]{amsbook}

\usepackage{amscd,amssymb}
\usepackage{epsfig}
\usepackage{float}
\usepackage{url}
\usepackage[all]{xy}
\usepackage{stmaryrd}
\usepackage{imakeidx}
\makeindex
\renewcommand{\theequation}{\mbox{\arabic{chapter}.\arabic{equation}}}

\renewcommand{\thefigure}{\mbox{\arabic{chapter}.\arabic{figure}}}

\newcommand{\N}{{\mathbb N}}
\newcommand{\Z}{{\mathbb Z}}
\newcommand{\Q}{{\mathbb Q}}
\newcommand{\C}{{\mathbb C}}
\newcommand{\D}{{\mathbb D}}
\newcommand{\R}{{\mathbb R}}
\newcommand{\F}{{\mathbb F}}
\renewcommand{\P}{{\mathbb P}}
\renewcommand{\H}{{\mathbb H}}

\newcommand{\BB}{{\mathcal B}}

\newcommand{\DD}{{\mathcal D}}
\newcommand{\EE}{{\mathcal E}}
\newcommand{\FF}{{\mathcal F}}
\newcommand{\GG}{{\mathcal G}}
\newcommand{\HH}{{\mathcal H}}
\newcommand{\KK}{{\mathcal K}}
\newcommand{\LL}{{\mathcal L}}
\newcommand{\MM}{{\mathcal M}}
\newcommand{\NN}{{\mathcal N}}

\newcommand{\OO}{{\mathcal O}}
\newcommand{\PP}{{\mathcal P}}

\newcommand{\RR}{{\mathcal R}}
\newcommand{\SSS}{{\mathcal S}}
\newcommand{\TT}{{\mathcal T}}

\newcommand{\VV}{{\mathcal V}}
\newcommand{\WW}{{\mathcal W}}
\newcommand{\XX}{{\mathcal X}}

\newcommand{\aaa}{{\bf a}}

\newcommand{\ddd}{{\rm d}}
\newcommand{\www}{\widetilde}
\newcommand{\whh}{\widehat}
\newcommand{\oooo}{\overline}
\newcommand{\uuuu}{\underline}

\newcommand{\mult}{{\rm mult}}
\newcommand{\paa}{\partial}

\newcommand{\pr}{{\rm pr}}

\newcommand{\Br}{{\rm Br}}

\DeclareMathOperator{\Aut}{Aut}
\DeclareMathOperator{\diag}{diag}

\DeclareMathOperator{\End}{End}
\DeclareMathOperator{\Hom}{Hom}
\DeclareMathOperator{\id}{id}
\DeclareMathOperator{\Imm}{Im}
\DeclareMathOperator{\Isom}{Isom}
\DeclareMathOperator{\Lie}{Lie}

\DeclareMathOperator{\mmod}{mod}
\DeclareMathOperator{\Ord}{Ord}
\DeclareMathOperator{\rk}{rk}
\DeclareMathOperator{\Rad}{Rad}
\DeclareMathOperator{\Ree}{Re}
\DeclareMathOperator{\sign}{sign}

\DeclareMathOperator{\Stab}{Stab}

\DeclareMathOperator{\tr}{tr}

\begin{document}

\theoremstyle{plain}
\newtheorem{lemma}{Lemma}[chapter]
\newtheorem{definition/lemma}[lemma]{Definition/Lemma}
\newtheorem{theorem}[lemma]{Theorem}
\newtheorem{proposition}[lemma]{Proposition}
\newtheorem{corollary}[lemma]{Corollary}
\newtheorem{conjecture}[lemma]{Conjecture}
\newtheorem{conjectures}[lemma]{Conjectures}

\theoremstyle{definition}
\newtheorem{definition}[lemma]{Definition}
\newtheorem{withouttitle}[lemma]{}
\newtheorem{remark}[lemma]{Remark}
\newtheorem{remarks}[lemma]{Remarks}
\newtheorem{example}[lemma]{Example}
\newtheorem{examples}[lemma]{Examples}
\newtheorem{notations}[lemma]{Notations}
\newtheorem{remarksandnotations}[lemma]{Remarks and Notations}

\title[Induced structures from upper triangular matrices]
{Unimodular bilinear lattices,\\ automorphism groups,\\ 
vanishing cycles,\\ monodromy groups,\\ distinguished bases,\\ 
braid group actions\\ and moduli spaces from \\upper triangular matrices}

\author[C. Hertling and K. Larabi]
{Claus Hertling and Khadija Larabi}

\address{Claus Hertling\\
Lehrstuhl f\"ur algebraische Geometrie, 
Universit\"at Mannheim,
B6 26, 68159 Mannheim, Germany}

\email{hertling@mail.uni-mannheim.de}

\address{Khadija Larabi\\
Lehrstuhl f\"ur algebraische Geometrie, 
Universit\"at Mannheim,
B6 26, 68159 Mannheim, Germany}

\email{klarabi@mail.uni-mannheim.de}

\date{December 21, 2024}

\subjclass[2020]{06B15, 20F55, 20F36, 14D05, 57M10, 32S30}

\keywords{Unimodular bilinear lattice, upper triangular matrix, Seifert form,
even and odd intersection form, even and odd monodromy group, vanishing cycle, 
braid group action, distinguished basis, manifold of
Stokes regions, isolated hypersurface singularity}

\thanks{This work was funded by the Deutsche 
Forschungsgemeinschaft (DFG, German Research Foundation) 
-- 494849004}


\begin{abstract}
This monograph starts with an upper triangular matrix
with integer entries and 1's on the diagonal. It develops from this
a spectrum of structures, which appear in different contexts,
in algebraic geometry, representation theory and the theory of
irregular meromorphic connections. 
It provides general tools to study these structures,
and it studies sytematically the cases of rank 2 and 3.
The rank 3 cases lead already to a rich variety of phenomena
and give an idea of the general landscape.
Their study takes up a large part of the monograph. 

Special cases are related to Coxeter groups, 
generalized Cartan lattices and exceptional sequences, 
or to isolated hypersurface singularities, their Milnor lattices
and their distinguished bases. 
But these make only a small part of all cases.
One case in rank 3 which is beyond them, is related to
quantum cohomology of $\P^2$ and to Markov triples.

The first structure associated to the matrix is 
a $\Z$-lattice with unimodular bilinear form (called Seifert form) 
and a triangular basis. It leads immediately to an even and
an odd intersection form, reflections and transvections,
an even and an odd monodromy group, even and odd vanishing cycles.
Braid group actions lead to braid group orbits of distinguished
bases and of upper triangular matrices.

Finally, complex manifolds, which consist of correctly glued Stokes regions, 
are associated to these braid group orbits.
A report on the case of isolated hypersurface singularities
concludes the monograph.
\end{abstract}

\maketitle

\tableofcontents

\setcounter{chapter}{0}

\chapter{Introduction}\label{s1}
\setcounter{equation}{0}
\setcounter{figure}{0}

This book develops many structures, starting from
a single upper triangular $n\times n$ matrix $S$
with integer entries and 1's on the diagonal.
The structures are introduced and are called playing 
characters in section \ref{s1.1}.

Such a matrix $S$ and induced structures appear in 
very different mathematical areas. 
Section \ref{s1.2} gives a panorama of these areas.
Though these areas are not subject of the book,
except for the area of isolated hypersurface singularities. 

Section \ref{s1.3} tells about the results in this book.
The book provides general tools and facts. 
It treats the cases $n=2$ and $n=3$ systematically. 
In chapter \ref{s10} it overviews the area of isolated
hypersurface singularities.

\section{Playing characters}
\label{s1.1}

$H_\Z$ will always be a {\it $\Z$-lattice}, so a free 
$\Z$-module of some  finite rank $n\in\N=\{1,2,3...\}$. 
Then $L$ will always be a nondegenerate bilinear form
$L:H_\Z\times H_\Z\to\Z$. It is called a {\it Seifert form}.
The pair $(H_\Z,L)$ is called a {\it bilinear lattice}. 
If for some $\Z$-basis $\uuuu{e}\in M_{1\times n}(H_\Z)$
of $H_\Z$ the determinant $\det L(\uuuu{e}^t,\uuuu{e})$
is $1$ the pair $(H_\Z,L)$ is called a 
{\it unimodular bilinear lattice}.
The notion {\it bilinear lattice} is from \cite{HK16}.
In chapter \ref{s2} we develop the following structures
for bilinear lattices, following \cite{HK16}.
Though in this introduction and in the chapters 
\ref{s3}--\ref{s10} we restrict to unimodular bilinear lattices. 
\begin{eqnarray*}
T^{uni}_n(\Z)&:=&\{ S\in M_{n\times n}(\Z)\,|\, 
S_{ij}=0\textup{ for }i>j,\ S_{ii}=1\}
\end{eqnarray*}
denotes the set of all upper triangular matrices with 
integer entries and 1's on the diagonal. 

Let $(H_\Z,L)$ be a unimodular bilinear lattice of rank $n$. 
A basis $\uuuu{e}$ of $H_\Z$ is called {\it triangular}
if $L(\uuuu{e}^t,\uuuu{e})^t\in T^{uni}_n(\Z)$. 
The transpose in the matrix is motivated by the case of
isolated hypersurface singularities. 
The set of triangular bases is called $\BB^{tri}$.
By far not every unimodular bilinear lattice has triangular
bases. But here we care only about those which have.

For fixed $n\in\N$ there is an obvious 1-1 correspondence
between the set of isomorphism classes of 
unimodular bilinear lattices $(H_\Z,L,\uuuu{e})$ with 
triangular bases and the set $T^{uni}_n(\Z)$, given
by the map 
$(H_\Z,L,\uuuu{e})\mapsto S:=L(\uuuu{e}^t,\uuuu{e})^t$.
To a given matrix $S\in T^{uni}_n(\Z)$ we always associate
the corresponding triple $(H_\Z,L,\uuuu{e})$. 

Let a unimodular bilinear lattice $(H_\Z,L,\uuuu{e})$ with a 
triangular basis $\uuuu{e}$ be given, with matrix 
$S=L(\uuuu{e}^t,\uuuu{e})^t\in T^{uni}_n(\Z)$. The following
objects are associated to this triple canonically. The names
are motivated by the case of isolated hypersurface
singularities.
\begin{list}{}{}
\item[(i)] A symmetric bilinear form 
$I^{(0)}:H_\Z\times H_\Z\to\Z$ and a skew-symmetric bilinear
form $I^{(1)}:H_\Z\times H_\Z\to\Z$ with 
\begin{eqnarray*}
I^{(0)}&=&L^t+L,\quad\textup{so } 
I^{(0)}(\uuuu{e}^t,\uuuu{e})=S+S^t,\\
I^{(1)}&=&L^t-L,\quad\textup{so } 
I^{(1)}(\uuuu{e}^t,\uuuu{e})=S-S^t,
\end{eqnarray*}
which are called {\it even} respectively {\it odd 
intersection form}.
\item[(ii)]
An automorphism $M:H_\Z\to H_\Z$ which is defined by
\begin{eqnarray*}
L(Ma,b)=L(b,a),\quad\textup{so }M(\uuuu{e}) = \uuuu{e}\cdot
S^{-1}S^t,
\end{eqnarray*}
and which is called {\it the monodromy}. It respects $L$
(and $I^{(0)}$ and $I^{(1)}$) because
$L(Ma,Mb)=L(Mb,a)=L(a,b)$. 
\item[(iii)] 
Six automorphism groups
\begin{eqnarray*}
O^{(k)}&:=& \Aut(H_\Z,I^{(k)})\qquad\textup{for }k\in\{0;1\},\\
G_\Z^M&:=& \Aut(H_\Z,M):=\{g:H_\Z\to H_\Z 
\textup{ automorphism }\,|\, gM=Mg\},\\
G_\Z^{(k)}&:=& \Aut(H_\Z,I^{(0)},M)=O^{(k)}\cap G_\Z^M
\qquad\textup{for }k\in\{0;1\},\\
G_\Z&:=&\Aut(H_\Z,L)=\Aut(H_\Z,L,I^{(0)},I^{(1)},M).
\end{eqnarray*}
\item[(iv)]
The set of {\it roots}
\begin{eqnarray*}
R^{(0)}:=\{a\in H_\Z\,|\, L(a,a)=1\},
\end{eqnarray*}
and the set
\begin{eqnarray*}
R^{(1)}:= H_\Z.
\end{eqnarray*}
\item[(v)] For $k\in\{0;1\}$ and $a\in R^{(k)}$ the 
{\it reflection} (if $k=0$) respectively {\it transvection} 
(if $k=1$) $s^{(k)}_a\in O^{(k)}$ with 
\begin{eqnarray*}
s^{(k)}_a(b):= b-I^{(k)}(a,b)a\quad\textup{for }b\in H_\Z.
\end{eqnarray*}
\item[(vi)]
For $k\in\{0;1\}$ the {\it even} (if $k=0$) respectively 
{\it odd} (if $k=1$) {\it monodromy group}
\begin{eqnarray*}
\Gamma^{(k)}:=\langle s^{(k)}_{e_1},...,s^{(k)}_{e_n}
\rangle \subset O^{(k)}.
\end{eqnarray*}
\item[(vii)]
For $k\in\{0;1\}$ the set of {\it even} (if $k=0)$ 
respectively {\it odd} (if $k=1)$ {\it vanishing cycles}
\begin{eqnarray*}
\Delta^{(k)}:=\Gamma^{(k)}\{\pm e_1,...,\pm e_n\}
\subset R^{(k)}.
\end{eqnarray*}
\end{list}
The definitions of all these objects require only
$S\in SL_n(\Z)$ and $e_1,...,e_n\in R^{(0)}$, 
not $S\in T^{uni}_n(\Z)$. But the formula (Theorem \ref{t2.6}) 
\begin{eqnarray*}
s^{(k)}_{e_1}...s^{(k)}_{e_n}=(-1)^{k+1}M
\quad\textup{for }k\in\{0;1\}
\end{eqnarray*}
depends crucially on $S\in T^{uni}_n(\Z)$. 

The even data $I^{(0)},O^{(0)},\Gamma^{(0)}$ and $\Delta^{(0)}$
are in many areas more important and are usually better understood
than the odd data $I^{(1)},O^{(1)},\Gamma^{(1)}$ and $\Delta^{(1)}$.
But in the area of isolated hypersurface singularities both turn up.

For $k\in\{0;1\}$ the group $\Gamma^{(k)}$ contains all 
reflections/transvections $s^{(k)}_a$ with $a\in\Delta^{(k)}$.
In the case of a bilinear lattice which is not unimodular this
holds for $k=0$, but not for $k=1$ (Remark \ref{t2.9} (iii)).
This is one reason why we restrict in the chapters 
\ref{s3}--\ref{s10} to unimodular bilinear lattices.

Section \ref{s3.2} gives an action of a semidirect product
$\Br_n\ltimes\{\pm 1\}^n$ of the braid group $\Br_n$ of
braids with $n$ strings and of a sign group $\{\pm 1\}^n$
on the set $(R^{(k)})^n$ for $k\in\{0;1\}$. It is compatible
with the Hurwitz action of $\Br_n$ on $(\Gamma^{(k)})^n$
with connecting map 
$$(R^{(k)})^n\to (\Gamma^{(k)})^n,\quad 
\uuuu{v}=(v_1,...,v_n)\mapsto (s^{(k)}_{v_1},...,s^{(k)}_{v_n}).$$
Both actions restrict to the same action on $\BB^{tri}$.
Especially, one obtains the orbit
$$\BB^{dist}:=\Br_n\ltimes\{\pm 1\}^n(\uuuu{e})\subset \BB^{tri}$$
of {\it distinguished bases} of $H_\Z$. 
The triple $(H_\Z,L,\BB^{dist})$ (up to isomorphism) is in
many cases a canonical object, whereas the choice of a 
distinguished basis $\uuuu{e}\in\BB^{dist}$ is a true choice.
The question whether $\BB^{dist}=\BB^{tri}$ or 
$\BB^{dist}\subsetneqq \BB^{tri}$ is usually a difficult question.
The subgroup
\begin{eqnarray*}
G_\Z^{\BB}:=\{g\in G_\Z\,|\, g(\BB^{dist})=\BB^{dist}\}
\end{eqnarray*}
is in many important cases equal to $G_\Z$. 
But if $\BB^{dist}\subsetneqq \BB^{tri}$, then 
$G_\Z^{\BB}\subsetneqq G_\Z$ is possible.

The action of $\Br_n\ltimes\{\pm 1\}^n$ on $\BB^{tri}$
is compatible with an action on $T^{uni}_n(\Z)$.
The orbit of $S$ is called
$$\SSS^{dist}:=\Br_n\ltimes\{\pm 1\}^n(S)\subset T^{uni}_n(\Z),$$
the matrices in it are called {\it distinguished matrices}.
As $\{\pm 1\}^n$ is the normal subgroup in the semidirect
product $\Br_n\ltimes\{\pm 1\}^n$, one can first divide out
the action of $\{\pm 1\}^n$. One obtains actions of $\Br_n$
on $\BB^{tri}/\{\pm 1\}^n$ and on $T^{uni}_n(\Z)/\{\pm 1\}^n$.
It will be interesting to determine the stabilizers
$(\Br_n)_{\uuuu{e}/\{\pm 1\}^n}$ of $\uuuu{e}/\{\pm 1\}^n$
and $(\Br_n)_{S/\{\pm 1\}^n}$ of $S/\{\pm 1\}^n$. 

Finally, chapter \ref{s8} introduces several complex manifolds
which are in fact {\it semisimple F-manifolds with 
Euler fields} (Definition \ref{t8.8}). First, the following
manifolds are universal and depend only on $n\in\N$,
\begin{eqnarray*}
C_n^{pure}&:=& \{\uuuu{u}=(u_1,...,u_n)\in\C^n\,|\, 
u_i\neq u_j\textup{ for all }i\neq j\},\\
C_n^{conf}&:=& C_n^{pure}/S_n \subset \C^n/S_n\cong \C^n,\\
C_n^{univ}&:=& \textup{the universal covering of }C_n^{conf}
\textup{ and }C_n^{pure}.
\end{eqnarray*}
$C_n^{conf}$ is the configuration space of $\Br_n$, with
$\pi_1(C_n^{conf})=\Br_n$. The subset of $C_n^{pure}$
\begin{eqnarray*}
F_n:=\{\uuuu{u}\in C_n^{pure}\,|\, \Imm(u_1)<...<\Imm(u_n)\}
\end{eqnarray*}
is a fundamental domain of the action of the symmetric
group $S_n$ on
$C_n^{pure}$ and maps under the projection 
$C_n^{pure}\to C_n^{conf}$ almost bijectively to $C_n^{conf}$.

Second, two manifolds $C_n^{\uuuu{e}/\{\pm 1\}^n}$ and
$C_n^{S/\{\pm 1\}^n}$ between $C_n^{conf}$ and $C_n^{univ}$
are constructed in section \ref{s8.1}. They are those coverings
of $C_n^{conf}$ whose fundamental groups embed into
$\Br_n=\pi_1(C_n^{conf})$ as the stabilizer 
$(\Br_n)_{\uuuu{e}/\{\pm 1\}^n}$ respectively the stabilizer  
$(\Br_n)_{S/\{\pm 1\}^n}$. In a more naive way, they are
obtained by glueing copies of $F_n$ in an appropriate way,
one copy for each element of $\BB^{dist}/\{\pm 1\}^n$
respectively $\SSS^{dist}/\{\pm 1\}^n$. 
These manifolds, especially $C_n^{\uuuu{e}/\{\pm 1\}^n}$,
are related to manifolds from algebraic geometry which
are associated to triples $(H_\Z,L,\uuuu{e})$,
and they have to compared with these manifolds.

\section
{The playing characters in action: areas, literature}
\label{s1.2}

\begin{figure}
\begin{xy}
\xymatrix{ 
\fbox{\parbox{3.5cm}{$\Br_n\ltimes\{\pm 1\}^n$ orbit
in $T^{uni}_n(\Z)$}} \ar@{-}[d]\ar@{-}[drr] & &\\
\fbox{\parbox{5cm}{Milnor lattice of an isola\-ted 
hypersurface singularity, distinguished basis}} \ar@{-}[d] & & 
\fbox{\parbox{4.5cm}{generalized Cartan lattice, 
crystallographic Coxeter group}} \ar@{-}[d] \\
\fbox{\parbox{5cm}{hom. mirror symmetry for invertible qh. 
singularities: Berglund-H\"ubsch duality}} \ar@{-}[drr] 
\ar@{-}[d]& & 
\fbox{\parbox{4.2cm}{non-crossing partition,
representation theory,
Grothendieck group of a hereditary algebra}} \ar@{-}[d] \ar@{-}[d]\\
\fbox{\parbox{4.5cm}{$\Z$-lattice of Lefschetz thimbles
for a tame function}} \ar@{-}[d]
\ar@{-}^{\textup{hom. mirror}}_{\textup{symmetry}}[rr] & & 
\fbox{\parbox{4.2cm}{semiorthogonal decom\-position in 
derived algebraic geometry}} \ar@{-}[d]\\
\fbox{\parbox{4.5cm}{moduli spaces from algebraic geometry}}
\ar@{-}[d] \ar@{-}[drr] 
\ar@{-}[rr]^{\textup{mirror}}_{\textup{symmetry}} 
& & 
\fbox{\parbox{4cm}{moduli spaces from quantum cohomology}}
\ar@{-}[d] \\ 
\fbox{\parbox{4.5cm}{isomonodromic family of connections 
with pole of Poincar\'e rank 2}} \ar@{-}[d] \ar@{-}[rr]
 & & \fbox{\parbox{4cm}{semisimple Frobenius manifold}} 
 \ar@{-}[dll] \ar@{-}[d]\\
\fbox{\parbox{4.5cm}{irregular meromorphic connection
with semisimple pole of order 2}} \ar@{-}[d] \ar@{-}[rr]
 & & \fbox{\parbox{4.5cm}{$N=2$ topological field theory,
 $tt^*$ geometry,\\
e.g. Toda equations}} \\ 
\fbox{\parbox{3cm}{$S\in T^{uni}_n(\C)\quad \& \\
\uuuu{u}\in C_n^{pure}\subset\C^n$}} 
\ar@{-}[rr] &  &
\fbox{\parbox{3.5cm}{$\Br_n\ltimes\{\pm 1\}^n$-action
on $T^{uni}_n(\C)$}} }
\end{xy}
\caption[Figure 1.1]{Mathematical areas where upper 
triangular (integer) matrices with 1's in the diagonal
appear}
\label{Fig:1.1}
\end{figure}
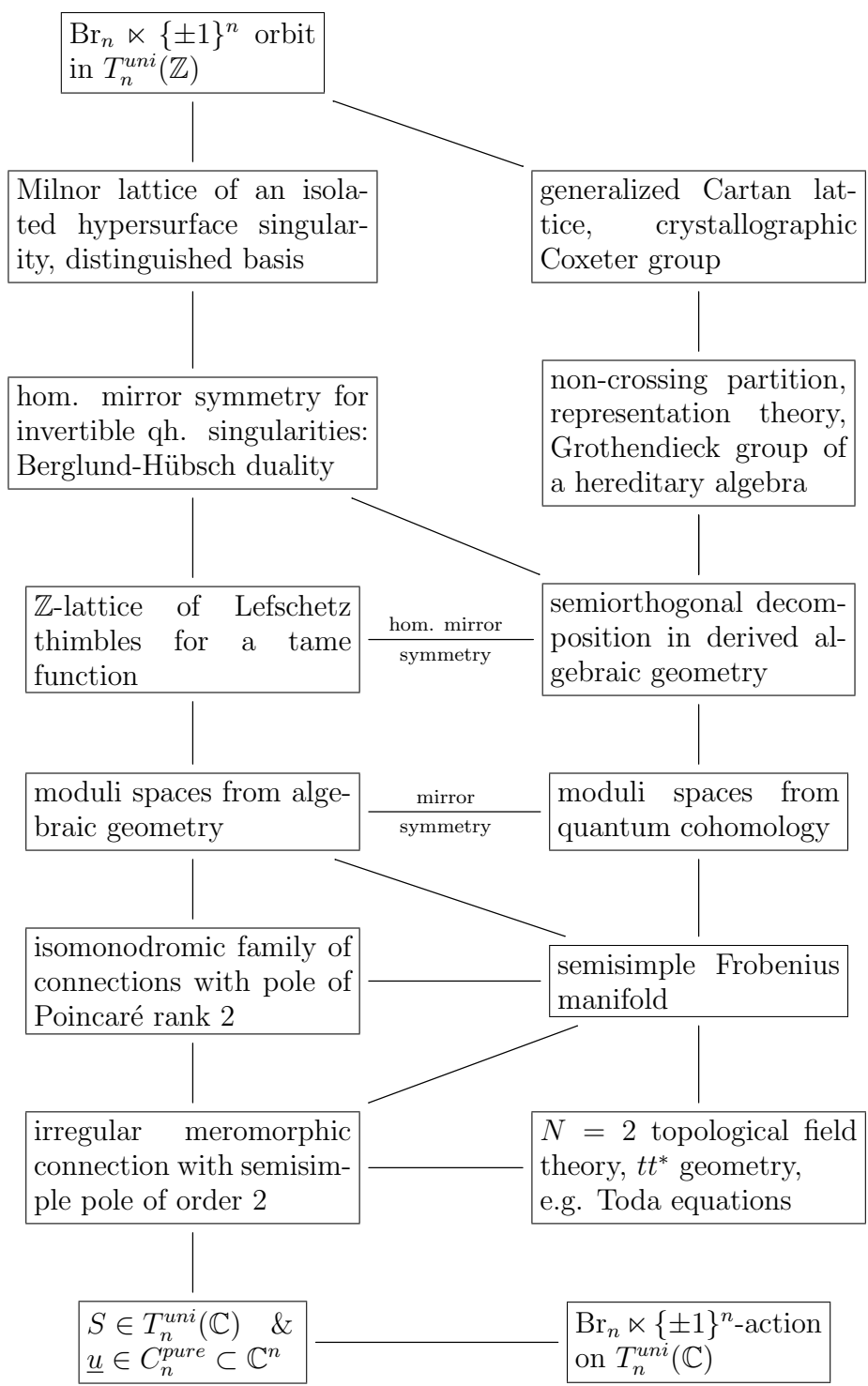

Upper triangular integer matrices $S$ with 1's on the diagonal
and the induced structures, which are described in section
\ref{s1.1}, arise in many different contexts.
Figure \ref{Fig:1.1} offers a landscape of mathematical areas
where they arise. The lines between the boxes indicate 
connections. Of course, the chart is very incomplete and 
subjective and does not at all give all connections between 
the areas.

In the following, for each of the boxes in Figure \ref{Fig:1.1} 
(and for a few more points),
some words are said and a few references are given. 
One can start with these references and use them to find 
more  references.

The mathematical areas in these boxes are not subject of 
this book, except for the area of isolated hypersurface 
singularities, which is treated in chapter \ref{s10}. 
(Therefore a reader could skip the rest of this section
\ref{s1.2}.) 

Historically the first appearance of
unimodular bilinear lattices $(H_\Z,L,\uuuu{e})$ with
triangular bases $\uuuu{e}$ is in the theory of isolated
hypersurfaces singularities. There they arise as
Milnor lattices with chosen distinguished bases \cite[Appendix]{Br70}.
All distinguished bases together form a $\Br_n\ltimes\{\pm 1\}^n$
orbit $\BB^{dist}$. A single distinguished basis $\uuuu{e}$ is a choice.
The triple $(H_\Z,L,\BB^{dist})$ is a canonical object associated 
to a singularity. These structures were studied a lot
\cite{AGV88}\cite{AGLV98}\cite{Eb01}\cite{Eb20}. 
Nevertheless it is not clear how to characterize the triples
$(H_\Z,L,\uuuu{e})$ or (equivalently) the matrices 
$S\in T^{uni}_n(\Z)$ which come from singularities with Milnor number $n$. 
See chapter \ref{s10} for more details and more references.

Closely related are the triples $(H_\Z,L,\uuuu{e})$ which come
by essentially the same construction from tame functions on
affine manifolds \cite{Ph83}\cite{Ph85}. 
Though they were less in the focus of people from singularity theory.

More recently they arose on one side of a version of mirror 
symmetry which connects such functions on the B-side with 
quantum cohomology of toric orbifolds on the A-side.
Iritani considered the integral structures on both sides \cite{Ir09}.
On the A-side this leads to the K-group of orbifold vector bundles
on toric orbifolds with the Mukai pairing.

This is part of a general appearance of triples
$(H_\Z,L,\uuuu{e})$ in derived algebraic geometry.
Gorodentsev \cite{Go94-1}\cite{Go94-2} studied them also
abstractly. Chapter 4 in the book \cite{CDG24} gives an
excellent survey on Gorodentsev's results.
It starts with a general point of view, the Grothendieck
group of a small triangulated category. Later it considers
the bounded derived category $\DD^b(X)$ of a complex 
projective manifold and its Grothendieck group.
Chapter 3 in the same book gives a survey on helices
in triangulated categories, which are behind the appearance
of triples $(H_\Z,L,\uuuu{e})$ in derived algebraic geometry.

Gorodentsev also classified pairs $(H_\C,L)$ over $\C$.
They split into indecomposable objects, which can be 
classified. A finer classification of pairs $(H_\R,L)$
over $\R$ was done by Nemethi \cite{Ne98} and in a more
explicit way in \cite{BH19}. Also pairs $(H_\R,L)$ split
into indecomposable objects, which can be classified.
Pairs $(H_\R,L)$ up to isomorphism are in 1-1 
correspondence with tuples of {\it spectral pairs} 
where the first entry is modulo $2\Z$.

A homological mirror symmetry for polynomials with 
singularities beyond the one by Iritani
was proposed by Takahashi. It compares structures above
pairs of invertible polynomials. One result for
chain type singularities is in \cite{AT20}.

The appearance of triples $(H_\Z,L,\uuuu{e})$ in 
derived algebraic geometry is related to the appearance
of such triples in representation theory. 
This motivated Hubery and Krause \cite{HK16} to consider
bilinear lattices $(H_\Z,L)$ with triangular bases
$\uuuu{e}$. Chapter \ref{s2} in this book takes up their
definition. Though the chapters \ref{s3} to \ref{s10} in this book
restrict to {\it unimodular} bilinear lattices with triangular
bases. One reason is that we consider, contrary to \cite{HK16},
also the odd monodromy group $\Gamma^{(1)}$ and the set 
$\Delta^{(1)}$ of odd vanishing cycles. They do not seem to
behave well in the case of a not unimodular bilinear lattice
with triangular basis (Remark \ref{t2.9} (iii)). 

\cite{HK16} is interested especially in the case of
a {\it generalized Cartan lattice}. There a matrix
$S\in M_{n\times n}(\Z)$, with which one starts, has entries 
$S_{ij}=0$ for $i>j$, $S_{ij}\leq 0$ for $i<j$, 
and diagonal entries $S_{ii}\in\N$ with 
$\frac{S_{ij}}{S_{ii}},\frac{S_{ij}}{S_{jj}}\in\Z_{\leq 0}$ 
for $i< j$, and any such matrix works. This is easier
than in the case of isolated hypersurface singularities.
The intersection of both cases consists only of 
the ADE-singularities.

Every generalized Cartan lattice arises as the 
Grothendieck group of a hereditary artin algebra
\cite[Proposition 4.6]{HK16}. A crucial result here
(cited in Theorem \ref{t3.6})
is the surjectivity of the braid group action on the set
of those $n$-tuples of reflections in the Coxeter group of a 
Coxeter system, whose product is the monodromy 
\cite{IS10}\cite{BDSW14}.
Using this and additional work, \cite{HK16}
can reprove results of Ingalls-Thomas,
Igusa-Schiffler, Crawley-Boevey and Ringel
on exceptional sequences and non-crossing partitions
related to finite dimensional modules of 
finite dimensional hereditary algebras. 
\cite{BBGKK19} gives a survey on non-crossing partitions.

\cite{BWY19} is between the points of view
of \cite{HK16} and of \cite{AT20}.
The last result in it concerns the hyperbolic singularities.

The homological mirror symmetry is related to a
more down to earth version of mirror symmetry,
which compares moduli spaces on the A-side and B-side.
On the B-side one has moduli spaces of (generalized)
complex structures, on the A-side they are related to
symplectic data, namely Gromov-Witten invariants.
In good cases, both sides lead to Frobenius manifolds
which are then by mirror symmetry isomorphic \cite{RS15}.
These Frobenius manifolds are generically semisimple
and also in other aspects rather special, as
they come from algebraic geometry.

This leads to the last appearance of upper triangular
matrices $S\in T^{uni}_n(\Z)$ which we want to mention,
namely as the carrier of the Stokes structure of a 
holomorphic vector bundle on the unit disk $\D$
with a (flat) meromorphic connection with a semisimple
order 2 pole at $0\in\D$ and pairwise different
eigenvalues of the pole part. After some choice this
structure is given by a tuple 
$\uuuu{u}\in C_n^{pure}\subset\C^n$ of eigenvalues
of the pole part, two upper triangular matrices in
$T^{uni}_n(\C)$, which are called {\it Stokes matrices},
and a diagonal matrix carrying the exponents of the
formal connection
\cite{Ma83}\cite{Sa02}\cite{HS11}. 

In cases from algebraic geometry, the diagonal matrix carrying
the exponents is the unit matrix times 
a constant in $\frac{1}{2}\Z$, 
and the two Stokes matrices coincide. 
The Stokes matrix is usually in $T^{uni}_n(\Z)$ 
as there the connection comes from oscillating integrals
over a distinguished basis of Lefschetz thimbles, 
which therefore give a splitting of the Stokes structure.

Varying $\uuuu{u}$ in a suitable covering of $C_n^{conf}$
(the best is the covering $C_n^{\uuuu{e}/\{\pm 1\}^n}$ 
in the notation of chapter \ref{s8}) 
gives rise to a holomorphic bundle on 
$\D\times C_n^{\uuuu{e}/\{\pm 1\}^n}$ with a meromorphic
connection with a pole of Poincar\'e rank 1 along
$\{0\}\times C_n^{\uuuu{e}/\{\pm 1\}^n}$ \cite{Sa02}.

Some additional choices/enrichments allow to construct
from this isomonodromic family of connections 
the structure of a Frobenius manifold on the covering
$C_n^{\uuuu{e}/\{\pm 1\}^n}$ \cite{Sa02}\cite{He03}\cite{Du99}
(section \ref{s8} describes the underlying F-manifold). 
The construction of Frobenius manifolds from quantum
cohomology is described in \cite{Ma99}. 

A different way (without additional choices) to enrich
the isomonodromic family leads to a TEZP structure
in the sense of \cite{He03} on the manifold
$C_n^{\uuuu{e}/\{\pm 1\}^n}$. This formalizes the
$tt^*$ geometry of \cite{CV93}. It induces on the manifold
$C_n^{\uuuu{e}/\{\pm 1\}^n}$ interesting 
differential geometric structure.
\cite{CV93} starts with families of 
$N=2$ topological field theories, which are objects
from mathematical physics. But the heart of them
are real analytic isomonodromic families of holomorphic 
bundles on $\P^1$ with meromorphic connections 
with poles of order 2 at 0 and $\infty$ 
which both come from the same matrix $S\in T^{uni}_n(\Z)$. 

As concrete cases, \cite{CV93} classifies and interprets 
those  $\Br_3\ltimes\{\pm 1\}^3$ orbits of matrices 
$S\in T^{uni}_3(\Z)$ whose monodromy matrices 
$S^{-1}S^t$ have eigenvalues in $S^1$. We refine this
in Theorem \ref{t4.6}. \cite{GL12} studies certain
higher rank cases, which are related to the Toda equations.
Also \cite{GL12} considers especially the cases related to
matrices $S\in T^{uni}_n(\Z)$.

$tt^*$ geometry makes also sense for matrices 
$S\in T^{uni}_n(\R)$. Frobenius manifolds can also
be constructed from matrices $S\in T^{uni}_n(\C)$
(and additional choices). The group 
$\Br_n\ltimes\{\pm 1\}^n$ acts also on $T^{uni}_n(\C)$. 
This manifold or coverings or other enrichments of it
parametrize Frobenius manifolds \cite{Ma99} or
meromorphic connections (their monodromy data)
and come equipped with interesting differential geometric
structure. But in this book we care only about the
points in the discrete subset 
$T^{uni}_n(\Z)\subset T^{uni}_n(\C)$.

\section{Results}
\label{s1.3}

Section \ref{s1.1} associated to each matrix $S\in T^{uni}_n(\Z)$
an impressive list of algebraic-combinatorial data and also a
few complex manifolds. For a given matrix $S$ there are many 
natural questions which all aim at controlling parts of these data,
for example:
\begin{list}{}{}
\item[(i)]
What can one say about the $\Z$-lattice $(H_\Z,I^{(0)})$ with the
even intersection form, e.g. its signature?
\item[(ii)]
What can one say about the $\Z$-lattice $(H_\Z,I^{(1)})$ with the
odd intersection form?
\item[(iii)]
What are the eigenvalues and the Jordan block structure of 
the monodromy $M$?
\item[(iv)] 
How big are the groups $G_\Z$, $G_\Z^{(0)}$, $G_\Z^{(1)}$
and $G_\Z^M$?
\item[(v)]
How good can one understand the even monodromy group $\Gamma^{(0)}$?
Is it determined by the pair $(H_\Z,I^{(0)})$ alone? 
\item[(vi)]
How good can one understand the odd monodromy group $\Gamma^{(1)}$?
Is it determined by the pair $(H_\Z,I^{(1)})$ alone? 
\item[(vii)]
Is $\Delta^{(0)}=R^{(0)}$ or $\Delta^{(0)}\subsetneqq R^{(0)}$?
How explicitly can one control these two sets?
How explicitly can one control $\Delta^{(1)}$?
\item[(viii)]
Is there an easy description of the set $\BB^{dist}$ of 
distinguished bases? Is $\BB^{dist}=\BB^{tri}$ or
$\BB^{dist}\subsetneqq \BB^{tri}$? 
\item[(ix)]
Is $G_\Z^{\BB}=G_\Z$ or $G_\Z^{\BB}\subsetneqq G_\Z$?
\item[(x)]
How do the manifolds $C_n^{\uuuu{e}/\{\pm 1\}^n}$ and
$C_n^{S/\{\pm 1\}^n}$ look like? Do they have natural
partial compactifications?
\end{list}

In the cases from isolated hypersurface singularities, 
$S$ is quite special.
There some of the questions have beautiful answers
((i), (ii), (iii) and (v), partially (vi) and (vii)), others are
wide open. The even case is better understood than the odd case.
Chapter \ref{s10} overviews this.

In the case of a generalized Cartan lattice, some $S$ is easy to
give. The even case is well understood ((i), (iii) and (v), 
partially (vii) and (viii)), the odd case is wide open and 
has probably much less beautiful structure.

In this book we concentrate on general tools
and on the cases of rank 2 and rank 3. 
The cases of rank 2 are already interesting, but still very special.
The cases of rank 3 are still in some sense small, 
but they show already a big variety of different types and 
phenomena. We consider them as sufficiently general to give an 
idea of the landscape for arbitrary rank $n\in\N$. 
The large number of pages of this book is due to the 
systematic study of all cases of rank 3.

Here the singularity cases form just two cases ($A_3$, $A_2A_1$), 
and also the cases from generalized Cartan lattices form a subset 
which one can roughly estimate as one third of all cases, 
not containing some of the most interesting cases ($\HH_{1,2}$, $\P^2$).

\begin{examples}\label{t1.1}
In the following examples, some matrices in $T^{uni}_2(\Z)$ 
and $T^{uni}_3(\Z)$ are distinguished. They cover the most 
important cases in $T^{uni}_2(\Z)$ and $T^{uni}_3(\Z)$. This will
be made precise in Theorem \ref{t1.2}, which gives results
on the braid group action on $T^{uni}_3(\Z)$. 
\begin{eqnarray*}
\begin{array}{ccccc}
S(A_1^2) & S(A_2) & S(\P^1) & S(x)\textup{ for }x\in\Z & 
S(A_1^3) \\
\begin{pmatrix}1&0\\0&1\end{pmatrix} & 
\begin{pmatrix}1&-1\\0&1\end{pmatrix} &
\begin{pmatrix}1&-2\\0&1\end{pmatrix} &
\begin{pmatrix}1&x\\0&1\end{pmatrix} &
\begin{pmatrix}1&0&0\\0&1&0\\0&0&1\end{pmatrix} 
\end{array}
\end{eqnarray*}
\begin{eqnarray*}
\begin{array}{cccc}
S(\P^2) & S(A_2A_1) & S(A_3) & S(\P^1A_1) \\
\begin{pmatrix}1&-3&3\\0&1&-3\\0&0&1\end{pmatrix} &
\begin{pmatrix}1&-1&0\\0&1&0\\0&0&1\end{pmatrix} &
\begin{pmatrix}1&-1&0\\0&1&-1\\0&0&1\end{pmatrix} & 
\begin{pmatrix}1&-2&0\\0&1&0\\0&0&1\end{pmatrix} 
\end{array}
\end{eqnarray*}
\begin{eqnarray*}
\begin{array}{cccc}
 & & S(-l,2,-l) & S(x_1,x_2,x_3) \\
S(\widehat{A}_2) & S(\HH_{1,2}) & \textup{for }l\geq 3 & 
\textup{for }x_1,x_2,x_3\in\Z\\
\begin{pmatrix}1&-1&-1\\0&1&-1\\0&0&1\end{pmatrix} &
\begin{pmatrix}1&-2&2\\0&1&-2\\0&0&1\end{pmatrix} &
\begin{pmatrix}1&-l&2\\0&1&-l\\0&0&1\end{pmatrix} & 
\begin{pmatrix}1&x_1&x_2\\0&1&x_3\\0&0&1\end{pmatrix} 
\end{array}
\end{eqnarray*}
The notations $A_1^2$
\index{$A_1^2,\ A_2,\, A_1^3,\ A_2A_1,\ A_3,\ \whh{A}_2$}, $A_2$, 
$A_1^3$, $A_2A_1$, $A_3$
and $\whh{A}_2$ are due the facts that $(H_\Z,I^{(0)})$ 
is in these cases the corresponding root lattice respectively
in the case $\whh{A}_2$ the affine root lattice of type 
$\whh{A}_2$. The notations $\P^1$\index{$\P^1,\ \P^2$} 
and $\P^2$ come from
the quantum cohomology of $\P^1$ and $\P^2$: The matrix
$S$ is in these cases one Stokes matrix of the associated
meromorphic connection with a semisimple pole of order 2
\cite{Gu99}\cite[Example 4.4 and (4.97)]{Du99}.
The notation $\HH_{1,2}$\index{$\HH_{1,2}$} is related to a Hurwitz space:
Here $(H_\Z,L)$ comes from a $\Z$-lattice of Lefschetz
thimbles for a branched covering of degree 2 from an elliptic
curve to $\P^1$ with 4 simple branch points,
one of them above $\infty$ \cite[Lecture 5]{Du96}.
A different point of view is that here $(H_\Z,I^{(0)})$
is an extended affine root lattice of type
$A_1^{(1,1)*}$. We prefer the notation $\HH_{1,2}$. 
\end{examples}

A large part of this book is devoted to answering the 
questions above for the cases of rank 2 and 3. 
Though the chapters \ref{s2}, \ref{s3}, \ref{s8}
and the sections \ref{s5.1}, \ref{s6.1} and \ref{s7.1}
offer also a lot of background material and tools.
Chapter \ref{s10} overviews the state of the art in 
the case of isolated hypersurface singularities.
In the following, we present some key results from the
chapters \ref{s4} to \ref{s9}. 

The action of $\Br_3\ltimes\{\pm 1\}^3$ on $T^{uni}_3(\Z)$
boils down to an action of $PSL_2(\Z)\ltimes G^{sign}$
on $T^{uni}_3(\Z)$ where $G^{sign}\cong \{\pm 1\}^2$ comes
from the action of the sign group $\{\pm 1\}^3$. 
As the action of $PSL_2(\Z)$ is partially nonlinear, 
it is good to write it as a semidirect product $PSL_2(\Z)\cong
G^{phi}\rtimes\langle\gamma\rangle$ where $\gamma$
acts cyclically and linearly of order 3 and $G^{phi}$
is a free Coxeter group with 3 generators which act
nonlinearly.

The sections \ref{s4.2}--\ref{s4.4} analyze the action
on $T^{uni}_3(\Z)$ carefully. The first result Theorem
\ref{t4.6} builds on coarser classifications of
Kr\"uger \cite[\S 12]{Kr90} and Cecotti-Vafa 
\cite[Ch. 6.2]{CV93}. The following theorem 
gives a part of Theorem \ref{t4.6}.

\begin{theorem}\label{t1.2} (Part of Theorem \ref{t4.6})

(a) The characteristic polynomial of $S^{-1}S^t$ and of the
monodromy $M$ of $(H_\Z,L,\uuuu{e})$ for 
$S=S(\uuuu{x})\in T^{uni}_3(\Z)$ with $\uuuu{x}\in\Z^3$ is 
\begin{eqnarray*}
p_{ch,M}&=&(t-1)(t^2-(2-r(\uuuu{x}))t+1),\\
\textup{where}&& r:\Z^3\to\Z,\quad
\uuuu{x}\mapsto x_1^2+x_2^2+x_3^2-x_1x_2x_3.
\end{eqnarray*} 
The characteristic polynomial and $r(\uuuu{x})$ are invariants
of the $\Br_3\ltimes\{\pm 1\}^3$ orbit of $S(\uuuu{x})$. 
All eigenvalues of $p_{ch,M}$ are unit roots if and only
if $r(\uuuu{x})\in\{0,1,2,3,4\}$. 

(b)  For $\rho\in\Z-\{4\}$ the fiber $r^{-1}(\rho)\subset\Z^3$
consists only of finitely many $\Br_3\ltimes\{\pm 1\}^3$ orbits.
The following table gives the symbols in Example \ref{t1.1}
for the fibers over $r\in\{0,1,2,3,4\}$, so there are only
seven orbits plus one series of orbits over $r\in\{0,1,2,3,4\}$. 
\begin{eqnarray*}
\begin{array}{c|c|c|c|c|c}
r(\uuuu{x}) & 0 & 1 & 2 & 3 & 4 \\  \hline 
 & A_1^3,\P^2 & A_2A_1 & A_3 & - & \P^1A_1,\whh{A}_2,\HH_{1,2}, 
S(-l,2,-l)\textup{ with  }l\geq 3
\end{array}
\end{eqnarray*}
\end{theorem}

With the help of certain (beautiful) graphs, 
in Theorem \ref{t4.13} the stabilizers
$(\Br_3)_{S/\{\pm 1\}^3}$ are calculated for certain representatives
of all $\Br_3\ltimes\{\pm 1\}^3$-orbits in $T^{uni}_3(\Z)$. 
We work with 14 graphs and 24 sets of representatives. 

Lemma \ref{t5.8} gives informations on the characteristic
polynomial and the signature of $I^{(0)}$ in all rank 3 cases.

\begin{lemma}\label{t1.3} (Lemma \ref{t5.8} (b)) 

Consider $\uuuu{x}\in\Z^3$ with $r=r(\uuuu{x})<0$
or $>4$ or with $S(\uuuu{x})$ one of the cases in the table
in Theorem \ref{t1.2}. Then 
$p_{ch,M}=(t-\lambda_1)(t-\lambda_2)\Phi_1$ and $\sign I^{(0)}$
are as follows ($\Phi_m=$ the cyclotomic polynomial of 
$m$-th primitive unit roots). 
\begin{eqnarray*}
\begin{array}{llll}
r(\uuuu{x}) & p_{ch,M} & \sign I^{(0)}\hspace*{1cm} 
& S(\uuuu{x}) \\
r<0 & \lambda_1,\lambda_2>0 & (+--) & S(\uuuu{x}) \\
r=0 & \Phi_1^3 & (+++) & S(A_1^3)\\
r=0 & \Phi_1^3 & (+--) & S(\P^2)\\
r=1 & \Phi_6\Phi_1 & (+++) & S(A_2A_1) \\
r=2 & \Phi_4\Phi_1 & (+++) & S(A_3) \\
r=4 & \Phi_2^2\Phi_1 & (++\ 0)
& S(\P^1A_1) \\
r=4 & \Phi_2^2\Phi_1 & (++\ 0)
& S(\whh{A}_2) \\
r=4 & \Phi_2^2\Phi_1 & (+\ 0\ 0)
& S(\HH_{1,2}) \\
r=4 & \Phi_2^2\Phi_1 & (+\ 0\ -)
& S(-l,2,-l) \textup{ with }l\geq 3 \\
r>4 & \lambda_1,\lambda_2<0 & (++-) & S(\uuuu{x})
\end{array}
\end{eqnarray*}
\end{lemma}

Chapter \ref{s5} analyzes the groups $G_\Z,G_\Z^{(0)},G_\Z^{(1)}$
and $G_\Z^M$ in all rank 3 cases. This leads into an intricate
case discussion. The case $\HH_{1,2}$ is different from all
other cases as it is the only case where $G_\Z$ is not 
abelian and where the subgroup $\{\pm M^m\,|\, m\in\Z\}$
does not have finite index in $G_\Z$.
The automorphism $Q\in G_\Q:=\Aut(H_\Q,L)$ is defined for 
$r(\uuuu{x})\neq 0$. It is $\id$ on $\ker(M-\id)$ and 
$-\id$ on $\ker(M^2-(2-r)M+\id)$ (Definition \ref{t5.9}).
It is only in a few cases in $G_\Z$ (Theorem \ref{t5.11}).

\begin{theorem}\label{t1.4} (Part of the Theorems 
\ref{t5.11}, \ref{5.13}, \ref{t5.14}, \ref{t5.16}, \ref{t5.18},
\ref{t3.28})

(a) In the $\Br_3\ltimes\{\pm 1\}^3$ orbit of $S(\HH_{1,2})$ 
$$G_\Z\cong SL_2(\Z)\times\{\pm 1\},\quad M=Q,$$ 
and the subgroup
$\{\pm M^m\,|\, m\in\Z\}=\{\pm \id,\pm Q\}$ 
has infinite index in $G_\Z$. 

(b) In all other rank 3 cases the subgroup
$\{\pm M^m\,|\, m\in\Z\}$ has finite index in $G_\Z$
and $G_\Z$ is abelian. 
Then one of the five possibilities holds, 
\begin{eqnarray}\label{1.1}
G_\Z&=&O_3(\Z),\\
G_\Z&=&\{\id,Q\}\times\{\pm (M^{root})^m\,|\, m\in\Z\},
\label{1.2}\\
G_\Z&=&\{\pm (M^{root})^m\,|\, m\in\Z\},\label{1.3}\\
G_\Z&=&\{\id,Q\}\times\{\pm M^m\,|\, m\in\Z\},\label{1.4}\\
G_\Z&=&\{\pm M^m\,|\, m\in\Z\},\label{1.5}
\end{eqnarray}
where $M^{root}$ is a root of $\pm M$ or of $MQ$.
The following table gives the index 
$[G_\Z:\{\pm M^m\,|\, m\in\Z\}]\in\N$ and informations on $M^{root}$.

\begin{eqnarray*}
\begin{array}{llll}
 & \textup{matrix} & \textup{index} & M^{root} \\ 
\hline
\eqref{1.1} & S(A_1^3) & 24 &  \\ \hline 
\eqref{1.2} & S(x,0,0)\textup{ with }x<0 & 4 
& (M^{root})^2=MQ \\
\eqref{1.2} & S(-l,2,-l)\textup{ with }l\textup{ even} & l^2-4 
& (M^{root})^{l^2/2-2}=MQ \\
\eqref{1.2} & S(4,4,4)\textup{ and }S(5,5,5) & 6 
& (M^{root})^3=-M \\
\eqref{1.2} & S(4,4,8) & 4 & (M^{root})^2=M \\ \hline 
\eqref{1.3} & S(\P^2) & 3 & (M^{root})^3=M \\
\eqref{1.3} & S(\whh{A}_2) \textup{ and }S(x,x,x)& 3 
& (M^{root})^3=-M \\
 & \textup{with }x\in\Z-\{-1,0,...,5\} & & \\
\eqref{1.3} & S(-l,2,-l) \textup{ with }l\textup{ odd} 
& l^2-4 & (M^{root})^{l^2-4}=-M \\
\eqref{1.3} & S(2y,2y,2y^2) \textup{ with }x\in\Z_{\geq 3} 
& 2 & (M^{root})^2=M \\ \hline 
\eqref{1.4} & S(3,3,4)\textup{ and }S(x,x,0) & 2 & \\
& \textup{with }x\in\Z_{\geq 2} & & \\ \hline 
\eqref{1.5} & S(A_3)\textup{ and }S(\uuuu{x}) & 1 & \\
& \textup{in other }\Br_3\ltimes\{\pm 1\}^3\textup{ orbits} & &  
\end{array}
\end{eqnarray*}

(c) $G_\Z=G_\Z^{\BB}$ holds for all rank 3 cases except four cases,
the $\Br_3\ltimes\{\pm 1\}^3$ orbits of $S(\uuuu{x})$ with
$$\uuuu{x}\in\{(3,3,4),(4,4,4),(5,5,5),(4,4,8)\}.$$
In these four cases $Q\in G_\Z-G_\Z^{\BB}$. 
\end{theorem}

Though in higher rank it is easier to construct matrices
$S$ with $G_\Z^\BB\subsetneqq G_\Z$ (Remarks \ref{t3.29}). 

Chapter \ref{s6} studies the even and odd monodromy groups
and the sets of even and of odd vanishing cycles in the
rank 2 and rank 3 cases. The following theorem catches
some of the results on the even monodromy group 
$\Gamma^{(0)}$ and the set $\Delta^{(0)}$ of even vanishing cycles
for the rank 3 cases. 
The group $O^{(0),*}$ (Definition \ref{t6.4}) is a certain subgroup
of $O^{(0)}$ which is determined only by $(H_\Z,I^{(0)})$
(so independently of $\BB$).
Part (b) discusses only the (in general difficult) problem
whether $\Delta^{(0)}=R^{(0)}$ or $\Delta^{(0)}\subsetneqq R^{(0)}$.
Theorem \ref{t6.14} contains many more informations on 
$\Delta^{(0)}$. Theorem \ref{t6.11} contains many more
informations on $\Gamma^{(0)}$ than part (a) below. Remarkably,
$\Gamma^{(0)}\cong G^{fCox,3}$ (the free Coxeter group with
three generators) holds not only for the Coxeter cases
$\uuuu{x}\in\Z^3_{\leq -2}$ (which all satisfy $r(\uuuu{x})>4$), 
but also in all cases $\uuuu{x}\in\Z^3$ with $r(\uuuu{x})<0$ and
in the case $\P^2$.

\begin{theorem}\label{t1.5}
(a) (Part of Lemma \ref{t2.11} and Theorem \ref{t6.11}) 

(i) (Part of Lemma \ref{t2.11}) The case $A_1^n$, $n\in\N$: 
\begin{eqnarray*}
\Gamma^{(0)}\cong \{\pm 1\}^n,\quad
\Gamma^{(1)}=\{\id\},\quad 
\Delta^{(0)}=R^{(0)}=\Delta^{(1)}=\{\pm e_1,...,\pm e_n\}.
\end{eqnarray*}

(ii) The cases with $r(\uuuu{x})>0$ and the cases
$A_3,\whh{A}_2,A_2A_1,\P^1A_1$: They contain all reducible
rank 3 cases except $A_1^3$. Then $\Gamma^{(0)}$
is a Coxeter group. If $\uuuu{x}\in\Z_{\leq -2}^3$ then
$\Gamma^{(0)}\cong G^{fCox,3}$. 

(iii) The cases $A_3,\whh{A}_2,\HH_{1,2}$: Then
$\Gamma^{(0)}= O^{(0),*}$.

(iv) The cases $S(-l,2,-l)$ with $l\geq 3$: Then
$\Gamma^{(0)}\stackrel{1:l}{\subset} O^{(0),*}.$

(v) The cases $\P^2$ and $\uuuu{x}\in\Z^3$ with 
$r(\uuuu{x})<0$: Then $\Gamma^{(0)}\cong G^{fCox,3}$. 

(b) (Part of Theorem \ref{t6.14})

(i) $\Delta^{(0)}=R^{(0)}$ holds in the following cases:
$A_3,\whh{A}_2,\P^2$, all $S(\uuuu{x})$ with 
$\uuuu{x}\in\{0,-1,-2\}$, all reducible cases.

(ii) $\Delta^{(0)}\subsetneqq R^{(0)}$ holds in the following
cases: $\HH_{1,2}$, all $S(-l,2,-l)$ with $l\geq 3$,
$S(3,3,4),S(4,4,4),S(5,5,5),S(4,4,8)$.

(iii) In the cases of the other $\Br_3\ltimes\{\pm 1\}^3$ orbits
in $T^{uni}_3(\Z)$, we do not know whether
$\Delta^{(0)}=R^{(0)}$ or $\Delta^{(0)}\subsetneqq R^{(0)}$
holds.
\end{theorem}

Let $E_n$ denote the $n\times n$ unit matrix.
Given $S\in T^{uni}_n(\Z)$ with associated triple
$(H_\Z,L,\uuuu{e})$, consider the
matrix $\www{S}:=2E_n-S\in T^{uni}_n(\Z)$ with the
associated triple $(H_\Z,\www{L},\uuuu{e})$.
Then $\www{L}$, $\www{I}^{(0)}$ and $\www{M}$ are far from
$L$, $I^{(0)}$ and $M$, but 
$\www{I}^{(1)}=-I^{(1)}$, 
$\www{\Gamma}^{(1)}=\Gamma^{(1)}$ and 
$\www{\Delta}^{(1)}=\Delta^{(1)}$ (see the Remarks \ref{t4.17}).
For example the cases $A_3$ and $\whh{A}_2$ are related
in this way, and also the Coxeter case $(-2,-2,-2)$ 
and the case $\HH_{1,2}$ are related in this way
(in both cases after an action of $\Br_3\ltimes\{\pm 1\}^3$).

This motivates in the rank 3 cases to consider the action
of the bigger group 
$(G^{phi}\ltimes\www{G}^{sign})\rtimes\langle\gamma\rangle$
on $\Z^3$ 
where $\www{G}^{sign}$ is generated by $G^{sign}$ and
the total sign change $\delta^\R:\uuuu{x}\mapsto -\uuuu{x}$.
Lemma \ref{t4.18} gives representatives for all orbits
of this action on $\Z^3$ (respectively $T^{uni}_3(\Z)$).
Still it is difficult to see for a given triple
$\uuuu{x}\in\Z^3$ in which orbit it is. 

We had for some time the hope that the beautiful facts
on the even monodromy group $\Gamma^{(0)}$ for the Coxeter cases
$\uuuu{x}\in\Z^3_{\leq 0}$ would have analoga for the
odd monodromy group $\Gamma^{(1)}$, but this does not
hold in general. In the case $(-2,-2,-2)$ 
$\Gamma^{(0)}\cong G^{fCox,3}$, but in the case $\HH_{1,2}$
not, and in both cases together 
$\Gamma^{(1)}\not\cong G^{free,3}$. 
On the other hand, $\Gamma^{(1)}\cong G^{free,3}$ for 
$\uuuu{x}\in B_1$, where $B_1\subset\Z^3$ is as follows.
\begin{eqnarray*}
B_1&:=& (G^{phi}\ltimes \www{G}^{sign})\rtimes 
\langle\gamma\rangle (\{\uuuu{x}\in\Z^3-\{(0,0,0)\}
\,|\, r(\uuuu{x})\leq 0 \}),\\
B_2&:=& \{\uuuu{x}\in\Z^3-\{(0,0,0)\}\,|\, S(\uuuu{x})
\textup{ is reducible}\},\\
B_3&:=& \{(0,0,0)\}.
\end{eqnarray*}
Though the set $B_1$ is difficult to understand
(see the Examples \ref{t4.20}). It contains 
$(3,3,3)$, so the orbit of $\P^2$. 
$B_2\cup B_3$ consists of the triples $\uuuu{x}$ 
with reducible $S(\uuuu{x})$, so with two or three
zero entries. 

Consider $\uuuu{x}\in\Z^3-B_3$.
The radical $\Rad I^{(1)}$ has rank 1, so the quotient
lattice $\oooo{H_\Z}^{(1)}:= H_\Z/\Rad I^{(1)}$ has rank 2.
Denote by $\Gamma^{(1)}_s$ the image of $\Gamma^{(1)}$
under the natural homomorphism
$\Gamma^{(1)}\to \Aut(\oooo{H_\Z}^{(1)})$ and by
$\Gamma^{(1)}_u$ the kernel of it. There is an exact sequence
$$\{1\}\to\Gamma^{(1)}_u\to \Gamma^{(1)}\to 
\Gamma^{(1)}_s\to\{1\}.$$
Denote by $\oooo{\Delta^{(1)}}\subset \oooo{H_\Z}^{(1)}$
the image of $\Delta^{(1)}$ in $\oooo{H_\Z}^{(1)}$. 
Often $\oooo{\Delta^{(1)}}$ is easier to describe
than $\Delta^{(1)}$. 

The long Theorems \ref{t6.18} and \ref{t6.21} offer 
detailed results about $\Gamma^{(1)}$ and $\Delta^{(1)}$
for the representatives in Lemma \ref{t4.18} of the
$(G^{phi}\ltimes\www{G}^{sign})\rtimes\langle\gamma\rangle$
orbits in $\Z^3$. The next theorem gives only a rough
impression.

\begin{theorem}\label{t1.6}
Consider $S=S(\uuuu{x})\in T^{uni}_3(\Z)$ and the
associated triple $(H_\Z,L,\uuuu{e})$. 

(a) (Part of Theorem \ref{t6.18}) 
Consider $\uuuu{x}\neq (0,0,0)$. 
\begin{eqnarray*}
\Gamma^{(1)}\cong G^{free,3}&\iff& \uuuu{x}\in B_1,\\
\Gamma^{(1)}_u=\{\id\}&\iff& \uuuu{x}\in B_1\cup B_2,\\
\Gamma^{(1)}_u\cong\Z^2&\iff& \uuuu{x}\in\Z^3-
(B_1\cup B_2\cup B_3),
\end{eqnarray*}
\begin{eqnarray*}
\Gamma^{(1)}_s&\cong& \textup{one of the groups }
SL_2(\Z),G^{free,2},G^{free,2}\times\{\pm 1\}\\
&&\textup{for }\uuuu{x}\in \Z^3-(B_1\cup B_3).
\end{eqnarray*}

(b) (Part of Theorem \ref{t6.21}) 

(i) In the cases of $A_3$ and $\whh{A}_2$ 
$\oooo{\Delta^{(1)}}=\oooo{H_\Z}^{(1),prim}$, so
$\oooo{\Delta^{(1)}}$ is the set of primitive vectors
in $\oooo{H_\Z}^{(1)}$, and $\Delta^{(1)}$ is the
full preimage in $H_\Z$ of $\oooo{\Delta^{(1)}}$. 

(ii) Though in many other cases 
$\oooo{\Delta^{(1)}}\not\subset \oooo{H_\Z}^{(1),prim}$,
and $\Delta^{(1)}$ is not the full preimage in $H_\Z$
of $\oooo{\Delta^{(1)}}$, but each fiber has infinitely
many elements.

(iii) But for $\uuuu{x}\in B_1$ the map 
$\Delta^{(1)}\to\oooo{\Delta^{(1)}}$ is a bijection. 
Especially for $\P^2$ $\oooo{\Delta^{(1)}}$ is easy
to describe (Theorem \ref{t6.21} (h)), but $\Delta^{(1)}$ not. 
\end{theorem}

Chapter \ref{s7} studies the set 
$\BB^{dist}=\Br_n\ltimes\{\pm 1\}^n(\uuuu{e})$ of 
distinguished bases for a given triple $(H_\Z,L,\uuuu{e})$. 
In general, it is difficult to characterize this orbit
in easy terms. We know that the inclusions in 
\eqref{3.3} and \eqref{3.4} hold. We are interested
when they are equalities.
\begin{eqnarray*}
\BB^{dist}\subset\{\uuuu{v}\in(\Delta^{(0)})^n\,|\, 
s_{v_1}^{(0)}...s_{v_n}^{(0)}=-M\},\hspace*{2cm}(3.3)\\
\BB^{dist}\subset\{\uuuu{v}\in(\Delta^{(1)})^n\,|\, 
s_{v_1}^{(1)}...s_{v_n}^{(1)}=M\}.\hspace*{2.4cm}(3.4)
\end{eqnarray*}
In general, this is a difficult question. In the rank 3 
cases Theorem \ref{t7.3} and Theorem \ref{t7.7}
give our results for \eqref{3.3} and \eqref{3.4}.

\begin{theorem}\label{t1.7}
Consider $S(\uuuu{x})\in T^{uni}_3(\Z)$ and the
associated triple $(H_\Z,L,\uuuu{e})$. 

(a) (Part of Theorem \ref{t7.3}) 

\eqref{3.3} is an equality for all cases except for
$\uuuu{x}$ in the 
$\Br_3\ltimes\{\pm 1\}^3$ orbit of $\HH_{1,2}$. 
There the right hand side of \eqref{3.3} consists of
countably many $\Br_3\ltimes\{\pm 1\}^3$ orbits.

(b) (Part of Theorem \ref{t7.7})

(i) The inclusion in \eqref{3.4} is an equality
$\iff\uuuu{x}\in B_1\cup B_2$. 

(ii) The cases $A_3$, $\whh{A}_2$, $\HH_{1,2}$
and $S(-l,2,l)$ with $l\geq 3$ are not in $B_1$.
But there the inclusion in \eqref{3.4} becomes an
equality if one adds on the right hand side of \eqref{3.4}
the condition $\sum_{i=1}^n \Z v_i=H_\Z$. 
\end{theorem}

The last section \ref{s7.4} of chapter \ref{s7}
builds on Theorem \ref{t4.16} which determines for a 
representative
$S\in T^{uni}_3(\Z)$ of each $\Br_3\ltimes\{\pm 1\}$ orbit
in $T^{uni}_3(\Z)$ the stabilizer $(Br_3)_{S/\{\pm 1\}^3}$. 
Theorem \ref{t7.11} determines in each of these cases
the stabilizer $(\Br_3)_{\uuuu{e}/\{\pm 1\}^3}$. 
The graphs $\GG_1,...,\GG_{14}$ in section \ref{s4.4}
used the groups $G^{phi}\rtimes\langle\gamma\rangle$.
At the end of section \ref{s7.4} different graphs,
which use the group $\Br_3$,
are introduced for the orbits of matrices as well as the
orbits of triangular bases. For the cases of finite orbits
and for the case $\whh{A}_2$ the graphs are given explicitly.
In the case of $A_3$ the orbit $\SSS^{dist}/\{\pm 1\}^3$ has four
elements and the orbit $\BB^{dist}/\{\pm 1\}^3$ has 16 elements.

The stabilizers $(\Br_3)_{\uuuu{e}/\{\pm 1\}^3}$ 
and $(\Br_3)_{S/\{\pm 1\}^3}$ are used in section \ref{s9}
for the construction of the manifolds 
$C_3^{\uuuu{e}/\{\pm 1\}^3}$ and 
$C_3^{S/\{\pm 1\}^3}$ in the rank 3 cases. 
Recall that they are obtained by glueing Stokes regions,
one for each element of an orbit $\BB^{dist}/\{\pm 1\}^3$
or $\SSS^{dist}/\{\pm 1\}^3$. 
Theorem \ref{t9.3} gives complete results, it covers all cases.
Remarkably, in all cases where $\Gamma^{(0)}\cong G^{fCox,3}$
or $\Gamma^{(1)}\cong G^{free,3}$, the manifold 
$C^{\uuuu{e}/\{\pm 1\}^3}$ is just $C_3^{univ}\cong\C^2\times\H$.
This holds especially for the case $\P^2$. 
In this aspect, the case $\P^2$ is easier than the case $A_3$,
where $C_3^{\uuuu{e}/\{\pm 1\}^3}$ consists of 16 
Stokes regions and is isomorphic to $\C^3$ minus two
hypersurfaces (a caustic and a Maxwell stratum).

Corollary \ref{t9.8} discusses natural partial compactifications.
Section \ref{s9.1} prepares this by giving explicitly
deck transformations of $C_3^{univ}$ which generate the
group $\Br_3$ of deck transformations of the covering
$\pr_3^{u,c}: C_3^{univ}\to C_3^{conf}$.
The formulas for the deck transformations use a
Schwarzian triangle function which is equivalent to the
$\lambda$-function and a certain lift of the logarithm.

Chapter \ref{s8} gives general background for these manifolds:
the configuration space $C_n^{conf}$ of $\Br_n$
and the related manifolds $C_n^{pure}$ and $C_n^{univ}$,
the notion of a (semisimple) F-manifold with Euler field,
distinguished systems of paths, two a priori different,
but closely related braid group actions on 
(homotopy classes) of distinguished systems of paths,
and a construction of $\Z$-lattice bundles on
$C_n^{univ}$, $C_n^{\uuuu{e}/\{\pm 1\}^n}$ and
$C_n^{S/\{\pm 1\}^n}$. They are the door to more 
transcendent structures on these manifolds, like
Dubrovin-Frobenius manifolds. 

We refer to the beginning of section \ref{s10} for an
introduction to the material of section \ref{s10},
which presents known results in the case of isolated 
hypersurface singularities. 

Appendix \ref{sa} recalls properties of the
hyperbolic plane and subgroups of its group of isometries.
The upper half plane, M\"obius transformations from
matrices in $SL_2(\R)$ and hyperbolic polygons are mentioned. Three special cases of the Poincar\'e-Maskit theorem 
are made explicit. Later also a model of the hyperbolic
plane from a cone of positive vectors in $\R^3$ with
a metric with signature $(+--)$ is explained.
This connects groups of $3\times 3$ matrices with 
groups of isometries. Appendix \ref{sa} is mainly
used in chapter \ref{s6} for the determination of some
groups $\Gamma^{(0)}$ and $\Gamma^{(1)}$, but also in
chapter \ref{s4} for the group $G^{phi}$.

Appendix \ref{sb} is used for the explicit formulas for
some generating deck transformations of the covering
$\pr_3^{u,c}:C_3^{univ}\to C_3^{conf}$ in section \ref{s9.1}.
It recalls properties of the first congruence subgroups
$\Gamma(n)\subset SL_2(\Z)$, those with $n\in\{2,3,4,5\}$,
with special emphasis on the most important case $n=2$.
In this case, the Schwarzian triangle function
$T:\H\to \C-\{0;1\}\cong\H/\Gamma(2)$, which we use, is
obtained by composing the $\lambda$-function with some
automorphism of $\C-\{0;1\}$. Later we introduce
a lift $\kappa:\H\to\C$ of the restriction 
$\ln:\C-\{0;1\}\dashrightarrow\C$ of the logarithm
$\ln:\C-\{0\}\dashrightarrow\C$, lifting it with $T$.

Appendix \ref{sc} contains Lemma \ref{tc.1} which determines
the unit groups in two families of rings of quadratic
algebraic integers. The lemma is used at many places,
especially in chapter \ref{s5} for the determination of
the groups $G_\Z,G_\Z^{(0)},G_\Z^{(1)},G_\Z^M$.
In order to present an elegant proof of Lemma \ref{t6.1},
we use the theory of continued fractions.
The results are certainly known, but we did not find a
completely satisfying reference.

Also Appendix \ref{sd} is used in chapter \ref{s5}.
It studies powers of units in rings of quadratic algebraic
integers.

{\bf Acknowledgements.}
We thank Martin Guest and Atsushi Takahashi for discussions
on topics related to this book.

\chapter{Bilinear lattices and induced structures}\label{s2}
\setcounter{equation}{0}
\setcounter{figure}{0}

This chapter fixes the basic notions,
a bilinear lattice and its associated data, namely 
a Seifert form, an even and an odd intersection form,
a monodromy, the roots, the triangular bases, 
an even and an odd monodromy group, 
the even and the odd vanishing cycles.
The notion of a bilinear lattice and the even part
of the associated data are considered in \cite{HK16}.
The more special case of a unimodular bilinear lattice
and even and odd data are considered since long time
in singularity theory \cite{AGV88}\cite{Eb01}.
In this paper we are mainly interested in unimodular
bilinear lattices. Only this chapter \ref{s2}
treats the general case, partially following \cite{HK16}.

\begin{notations}\label{t2.1}
In these notations, 
$R$ will be either the ring $\Z$ or
one of the fields $\Q$, $\R$ or $\C$. 
Later we will work mainly with $\Z$. 
If $R=\Z$ write $\www{R}:=\Q$, else write
$\www{R}:=R$. 

In the whole paper, \index{$H_R$}
$H_R\supsetneqq\{0\}$ is a finitely generated free $R$-module,
so a $\Z$-lattice\index{$\Z$-lattice} if $R=\Z$, 
and a finite dimensional $R$-vector space if $R$ is $\Q$, $\R$ or $\C$. 
Its rank will usually be called $n\in\N=\{1,2,3,...\}$
\index{$\N=\{1,2,3,...\}$}
(it is its dimension if $R$ is $\Q$, $\R$ or $\C$).
If $R_1$ and $R_2$ are both in the
list $\Z,\Q,\R,\C$ and $R_1$ is left of $R_2$ and
$H_{R_1}$ is given, then 
$H_{R_2}:=H_{R_1}\otimes_{R_1}R_2$.

In the whole paper, 
$L:H_R\times H_R\to R$ will be a nondegenerate
$R$-bilinear form. 
If $U\subset H_R$ is an $R$-submodule, then 
$U^{\perp}:=\{b\in H_R\,|\, L(U,b)=0\}$ and
${}^{\perp}U:=\{a\in H_R\,|\, L(a,U)=0\}$.
In the case $R=\Z$, $U^{\perp}$ and ${}^{\perp}U$
are obviously primitive $\Z$-submodules of $H_\Z$. 

In Lemma \ref{t2.2} we will start with $H_R$ and a symmetric
$R$-bilinear form $I^{[0]}:H_R\times H_R\to R$
or a skew-symmetric $R$-bilinear form 
$I^{[1]}:H_R\times H_R\to R$.
With the square brackets in the index we distinguish them
from the bilinear forms $I^{(0)}$ and $I^{(1)}$, which
are induced in Definition \ref{t2.3} 
by a given bilinear form $L$.
Though later they will be identified.

Suppose that $M:H_R\to H_R$ is an automorphism.
Then $M_s,M_u,N:H_{\www{R}}\to H_{\www{R}}$ denote the
semisimple part\index{semisimple part}, 
the unipotent part\index{unipotent part} 
and the \index{nilpotent part}nilpotent part 
of $M$ with $M=M_sM_u=M_uM_s$ and $N=\log M_u,e^N=M_u$.
Denote 
$H_\lambda:=\ker(M_s-\lambda\cdot \id):H_\C\to H_\C$.

For $m\in\N$ denote by $\Phi_m\in\Z[t]$ 
\index{$\Phi_m$} the cyclotomic polynomial 
\index{cyclotomic polynomial} 
whose zeros are the primitive $m$-th unit roots.
\end{notations}

The following lemma is elementary and classical. We skip the 
proof. 

\begin{lemma}\label{t2.2}
Let $R\in\{\Q,\R,\C\}$ and let $H_R$ be an $R$-vector space
of dimension $n\in\N$. 

(a) Let $I^{[0]}:H_R\times H_R\to R$ be a symmetric bilinear form.
Consider $a\in H_R$ with $I^{[0]}(a,a)\neq 0$. The map
\begin{eqnarray*}
s_a^{[0]}:H_R\to H_R,\quad 
s_a^{[0]}(b):=b-\frac{2I^{[0]}(a,b)}{I^{[0]}(a,a)}a,
\end{eqnarray*}
is a {\sf \index{reflection}reflection}, so it is in $\Aut(H_r,I^{[0]})$,
it fixes the codimension 1 subspace $\{b\in H_R\,|\, 
I^{[0]}(a,b)=0\}$ and it maps $a$ to $-a$. 
Especially $(s_a^{[0]})^2=\id$. 

(b) Let $I^{[1]}:H_R\times H_R\to R$ be a skew-symmetric 
bilinear form. Consider $a\in H_R$. The map
\begin{eqnarray*}
s_a^{[1]}:H_R\to H_R,\quad 
s_a^{[1]}(b):=b-I^{[1]}(a,b)a,
\end{eqnarray*}
is in $\Aut(H_R,I^{[1]})$ with 
$$(s_a^{[1]})^{-1}(b)=b+I^{[1]}(a,b)a.$$
It is $\id$ if $a\in \Rad(I^{[1]})$.
If $a\notin \Rad(I^{[1]})$ then 
it fixes the codimension 1 subspace $\{b\in H_R\,|\, 
I^{[1]}(a,b)=0\}$, and $s_a^{[1]}-\id$ is nilpotent 
with a single $2\times 2$ Jordan block.
Then it is called a {\sf \index{transvection}transvection}. 

(c) Fix $k\in\{0;1\}$ and consider $I^{[k]}$ as in (a) or (b).
An element $g\in\Aut(H_R,I^{[k]})$ and an element $a\in H_R$
with $I^{[0]}(a,a)\neq 0$ if $k=0$ satisfy
\begin{eqnarray*}
g\, s_a^{[k]}\, g^{-1}=s_{g(a)}^{[k]}.
\end{eqnarray*}

\end{lemma}

\begin{definition}\label{t2.3}
(a) \cite[ch. 2]{HK16} A 
{\it \index{bilinear lattice}bilinear lattice} is a pair
\index{$H_\Z$}\index{$(H_\Z,L)$}
$(H_\Z,L)$ with $H_\Z$ a $\Z$-lattice of some rank
$n\in\N$ together with a nondegenerate bilinear form
$L:H_\Z\times H_\Z\to \Z$. \index{$L$: Seifert form}
If $\det L(\uuuu{e}^t,\uuuu{e})=1$ for some
$\Z$-basis $\uuuu{e}=(e_1,...,e_n)$ of $H_\Z$ then
$L$ and the pair $(H_\Z,L)$ are called 
{\it \index{unimodular bilinear lattice}unimodular}.
The bilinear form is called {\it \index{Seifert form}Seifert form} in this paper.

(b) A bilinear lattice induces several structures:
\begin{list}{}{}
\item[(i)] 
\cite[ch. 2]{HK16} A symmetric bilinear form
\begin{eqnarray*}
I^{(0)}=L^t+L:H_\Z\times H_\Z\to\Z,\quad 
\textup{so }I^{(0)}(a,b)=L(b,a)+L(a,b),
\end{eqnarray*}
\index{$I^{(0)},\ I^{(1)}$} which is called 
{\it \index{even intersection form}even intersection form}.

\item[(ii)] 
A skew-symmetric bilinear form
\begin{eqnarray*}
I^{(1)}=L^t-L:H_\Z\times H_\Z\to\Z,\quad 
\textup{so }I^{(1)}(a,b)=L(b,a)-L(a,b),
\end{eqnarray*}
which is called {\it \index{odd intersection form}odd intersection form}.

\item[(iii)]
\cite[ch. 2]{HK16} An automorphism $M:H_\Q\to H_\Q$ 
which is defined by \index{$M$: monodromy} 
\begin{eqnarray*}
L(Ma,b)=L(b,a)\quad\textup{for }a,b\in H_\Q,
\end{eqnarray*}
which is called {\it \index{monodromy}monodromy}.

\item[(iv)]
Six \index{automorphism group}automorphism groups
\index{$O^{(0)},\ O^{(1)}$}
\index{$G_\Z^M,\ G_\Z^{(0)},\ G_\Z^{(1)},\ G_\Z$}
\begin{eqnarray*}
O^{(k)}&:=& \Aut(H_\Z,I^{(k)})\qquad\textup{for }k\in\{0;1\},\\
G_\Z^M&:=& \Aut(H_\Z,M):=\{g:H_\Z\to H_\Z 
\textup{ automorphism }\,|\, gM=Mg\},\\
G_\Z^{(k)}&:=& \Aut(H_\Z,I^{(0)},M)=O^{(k)}\cap G_\Z^M
\qquad\textup{for }k\in\{0;1\},\\
G_\Z&:=&\Aut(H_\Z,L).
\end{eqnarray*}

\item[(v)]
\cite[ch. 2]{HK16} The set $R^{(0)}\subset H_\Z$ of 
\index{$R^{(0)}$}{\it \index{root}roots},
\begin{eqnarray*}
R^{(0)}:=\{a\in H_\Z\,|\, L(a,a)>0; \frac{L(a,b)}{L(a,a)},\frac{L(b,a)}{L(a,a)}\in \Z
\textup{ for all }b\in H_\Z\}.
\end{eqnarray*}

\item[(vi)]
\cite[ch. 2]{HK16} The set $\BB^{tri}$ \index{$\BB^{tri}$} 
of {\it \index{triangular basis}triangular bases},
\begin{eqnarray*}
\BB^{tri}:=\{\uuuu{e}=(e_1,...,e_n)\in (R^{(0)})^n\,|\,
\bigoplus_{i=1}^n\Z e_i=H_\Z, L(e_i,e_j)=0\textup{ for }i<j\}.
\end{eqnarray*}
\end{list}

(c) Let $n\in\N$ and $R\in\{\Z,\Q,\R,\C\}$. The sets 
$T^{tri}_n$ and $T^{uni}_n(R)$ of \index{$T^{uni}_n(\Z)$} 
\index{upper triangular matrix}upper triangular matrices
are defined by 
\begin{eqnarray*}
T^{uni}_n(R):= \{S=(s_{ij})\in M_{n\times n}(R)&|& 
s_{ii}=1,s_{ij}=0\textup{ for }i>j\},\\
T^{tri}_n:= \{S=(s_{ij})\in M_{n\times n}(\Z)&|& 
s_{ii}\in \N, s_{ij}=0\textup{ for }i>j, \\
&& \frac{s_{ij}}{s_{ii}},\frac{s_{ji}}{s_{ii}}\in\Z
\textup{ for }i\neq j\}.
\end{eqnarray*}
Obviously $T^{uni}_n(\Z)\subset T^{tri}_n$. 
\end{definition}

\begin{remarks}\label{t2.4}
(i) There are bilinear lattices with $\BB^{tri}=\emptyset$.
We are interested only in bilinear lattices with 
$\BB^{tri}\neq\emptyset$. 

(ii) A triangular basis $\uuuu{e}\in\BB^{tri}$ is called
in \cite{HK16} a {\it complete exceptional sequence}.

(iii) In the case $\BB^{tri}\neq\emptyset$,
\cite{HK16} considers the bilinear form $L^t$
(with $L^t(a,b)=L(b,a)$). Our choice $L$ is motivated
by singularity theory. Also the names for 
$L$, $I^{(0)}$, $I^{(1)}$ and $M$, namely 
{\it Seifert form, even intersection form, odd intersection form} 
and {\it monodromy} are motivated by singularity theory.
The roots in $R^{(0)}$ are in \cite{HK16} also called
{\it pseudo-real roots}.

(iv) In this paper we are mainly interested in the cases of
unimodular bilinear lattices with $\BB^{tri}\neq\emptyset$.
Singularity theory leads to such cases.

(v) \cite{HK16} is mainly interested in the cases of 
{\it \index{generalized Cartan lattice}generalized Cartan lattices}. 
A generalized Cartan lattice
is a triple $(H_\Z,L,\uuuu{e})$ with $(H_\Z,L)$ a
bilinear lattice and $\uuuu{e}\in \BB^{tri}$ with
$L(e_i,e_j)\leq 0$ for $i>j$.
\end{remarks}

\begin{remarks}\label{t2.5}
(i) The classification of pairs $(H_\R,L)$ and pairs
$(H_\C,L)$ with $L$ a nondegenerate bilinear form on
$H_\R$ respectively $H_\C$ is well understood. 
Such a pair decomposes into an orthogonal sum
of irreducible pairs. This and the 
classification of the irreducible pairs
over $\R$ is carried out in \cite{Ne98} and, more explicitly,
in \cite{BH19}. 

In both references it is also proved that
a pair $(H_\R,L)$ of rank $n\in\N$ 
up to isomorphism is uniquely determined by an unordered
tuple of $n$ spectral pairs modulo $2\Z$, i.e. by
$n$ pairs $([\alpha_1],l_1),...,([\alpha_n],l_n)
\in\R/2\Z\times\Z$. Here $\alpha_1,...,\alpha_n\in\R$.
The eigenvalues of the monodromy $M$ are the numbers
$e^{-2\pi i\alpha_1},...,e^{-2\pi i \alpha_n}$.
The numbers $l_1,...,l_n$ determine the Jordan block
structure, see \cite{BH19} for details.

The classification over $\C$ follows easily. Though 
it was carried out before in \cite{Go94-1}\cite{Go94-2},
and it is formulated also in \cite[Theorem 4.22]{CDG24}.

(ii) A unimodular bilinear lattice $(H_\Z,L)$ is called in 
\cite{CDG24} a {\it Mukai pair}. In \cite[4.1--4.4]{CDG24}
basic results of Gorodentsev for $R=\Z$ or $R=\C$ are 
rewritten. The monodromy is there called {\it canonical
operator}. A triangular basis
is there called {\it exceptional}.

(iii) The classification over $\Z$, so of unimodular
bilinear lattices $(H_\Z,L)$, is wide open for larger $n$.
The case $n=3$ is treated in great detail in this book.
\end{remarks}

\begin{lemma}\label{t2.6}
(a) Let $(H_\Z,L)$ be a bilinear lattice of rank $n\in\N$.
\begin{list}{}{}
\item[(i)]
Let $\uuuu{e}=(e_1,...,e_n)$ be a $\Z$-basis of $H_\Z$.
Define $S:=L^t(\uuuu{e}^t,\uuuu{e})=L(\uuuu{e}^t,\uuuu{e})^t
\in  M_{n\times n}(\Z)\cap GL_n(\Q)$. Then
\begin{eqnarray*}
I^{(0)}(\uuuu{e}^t,\uuuu{e})=S+S^t,\ 
I^{(1)}(\uuuu{e}^t,\uuuu{e})=S-S^t,\ 
M(\uuuu{e})=\uuuu{e}S^{-1}S^t.
\end{eqnarray*}
\item[(ii)]
\begin{eqnarray*}
I^{(0)}(a,b)=L((M+\id)a,b),\quad 
\Rad I^{(0)}=\ker((M+\id):H_\Z\to H_\Z),\\\
I^{(1)}(a,b)=L((M-\id)a,b),\quad 
\Rad I^{(1)}=\ker((M-\id):H_\Z\to H_\Z).
\end{eqnarray*}
\item[(iii)]
\begin{eqnarray*}
G_\Z=\Aut(H_\Z,L,I^{(0)},I^{(1)},M)\subset
\left\{\begin{array}{c}G^{(0)}_\Z\\ G^{(1)}_\Z\end{array}
\right\}\subset G^M_\Z.
\end{eqnarray*}
\item[(iv)]
$M\in G_\Z$ if $(H_\Z,L)$ is unimodular or if 
$\BB^{tri}\neq\emptyset$. 
\item[(v)]
If $a\in R^{(0)}$ then 
\index{$s_a^{(0)},\ s_a^{(1)}$} 
\begin{eqnarray*}
s_a^{(0)}:=s_a^{[0]}\textup{ and }
s_a^{(1)}:=s_{a/\sqrt{L(a,a)}}^{[1]},\quad
\textup{so }s_a^{(1)}(b)=b-\frac{I^{(1)}(a,b)}{L(a,a)}a,
\end{eqnarray*}
are in $O^{(0)}$ respectively $O^{(1)}$. 
\item[(vi)]
\cite[Lemma 2.1]{HK16} If $a,b\in R^{(0)}$ then $s_a^{(0)}(b)\in R^{(0)}$
(but not necessarily $s_a^{(1)}(b)\in R^{(0)}$).
\item[(vii)]
If $a,b\in R^{(0)}$ with $L(a,b)=0$ then
\begin{eqnarray*}
L(s_a^{(1)}b,s_a^{(1)}b)=L(b,b),\quad 
s_a^{(1)}(b)\in R^{(0)},\quad
s_{s_a^{(1)}(b)}^{(1)}=s_a^{(1)} s_b^{(1)}(s_a^{(1)})^{-1}.
\end{eqnarray*}
\end{list}

(b) The map 
\begin{eqnarray*}
\{(H_\Z,L,\uuuu{e})\,|\, \begin{array}{ll}(H_\Z,L)
\textup{ is a bilinear}\\
\textup{lattice of rank }n,\ \uuuu{e}\in\BB^{tri}\end{array}
\}/\textup{isomorphism} \to T^{tri}_n
\end{eqnarray*}
is a bijection and restricts to a bijection 
\begin{eqnarray*}
\{(H_\Z,L,\uuuu{e})\,|\, \begin{array}{ll}(H_\Z,L)
\textup{ is a unimodular}\\ 
\textup{bilinear lattice of rank }n,\ \uuuu{e}\in\BB^{tri}\end{array}
\}/\textup{isom.} \to T^{uni}_n(\Z).
\end{eqnarray*}

(c) \cite[Lemma 3.10]{HK16} 
Let $(H_\Z,L)$ be a unimodular bilinear lattice with
$\BB^{tri}\neq\emptyset$. Then
\begin{eqnarray*}
R^{(0)}=\{a\in H_\Z\,|\, L(a,a)=1\}.
\end{eqnarray*}

(d) Let $(H_\Z,L)$ be a unimodular bilinear lattice with
$\BB^{tri}\neq\emptyset$. Define for $a\in H_\Z$
$s_a^{(1)}:=s_a^{[1]}\in O^{(1)}$. This definition is 
compatible with the definition of $s_a^{(1)}$ for $a\in R^{(0)}$
in part (a) (v). Furthermore now for $a,b\in H_\Z$
$$s^{(1)}_{s^{(1)}_a(b)}=s^{(1)}_a s^{(1)}_b (s^{(1)}_a)^{-1}.$$
\end{lemma}

{\bf Proof:}
(a) (i) The defining equation for $M$ can be written as
$L((M\uuuu{e})^t,\uuuu{e})=L(\uuuu{e},\uuuu{e})^t$, which
implies $M\uuuu{e}=\uuuu{e}S^{-1}S^t$. The rest is trivial.

(ii) Trivial.

(iii) $g\in G_\Z$ commutes with $M$ because
\begin{eqnarray*}
L(gMa,gb)=L(Ma,b)=L(b,a)=L(gb,ga)=L(Mga,gb).
\end{eqnarray*}
Of course it respects $I^{(0)}$ and $I^{(1)}$.
Therefore $G_\Z=\Aut(H_\Z,L,I^{(0)},I^{(1)},M)$. The rest is trivial.

(iv) The calculation $L(Ma,Mb)=L(Mb,a)=L(a,b)$ shows that $M$
respects $L$. It remains to show $M\in \Aut(H_\Z)$. 

This is clear if $(H_\Z,L)$ is unimodular. Suppose
$\BB^{tri}\neq\emptyset$ and $(H_\Z,L)$ not unimodular. 
Consider $\uuuu{e}\in \BB^{tri}$,
$S:=L(\uuuu{e}^t,\uuuu{e})^t\in T^{tri}_n$ and
$D:=\textup{diag}(s_{11},...,s_{nn})\in M_{n\times n}(\Z)$.
Then $D^{-1}S,SD^{-1}\in T^{uni}_n(\Z)$ and
\begin{eqnarray*}
S^{-1}S^t=S^{-1}DD^{-1}S^t=(D^{-1}S)^{-1}(SD^{-1})^t\in GL_n(\Z),
\end{eqnarray*}
so $M\in\Aut(H_\Z)$. 

(v) If $a\in R^{(0)}$ and $b\in H_\Z$ then 
$\frac{L(a,b)}{L(a,a)},\frac{L(b,a)}{L(a,a)}\in \Z$, so
\begin{eqnarray*}
\frac{2 I^{(0)}(a,b)}{I^{(0)}(a,a)},\frac{I^{(1)}(a,b)}{L(a,a)}
\in\Z\textup{ and }s^{(0)}_a(b),s^{(1)}_a(b),(s^{(1)}_a)^{-1}(b)
\in H_\Z.
\end{eqnarray*}

(vi) $L(s_a^{(0)}(b),s_a^{(0)}(b))=L(b,b)$ because 
$s_a^{(0)}\in G_\Z^{(0)}$ and $I^{(0)}=L^t+L$
(in general $L(s_a^{(1)}(b),s_a^{(1)}(b))\neq L(b,b)$).
For $c\in H_\Z$
\begin{eqnarray*}
\frac{L(s_a^{(0)}(b),c)}{L(s_a^{(0)}(b),s_a^{(0)}(b))}
&=& \frac{L(b-\frac{L(a,b)+L(b,a)}{L(a,a)}a,c)}{L(b,b)}\\
&=&\frac{L(b,c)}{L(b,b)}
-\frac{L(a,b)+L(b,a)}{L(b,b)} \frac{L(a,c)}{L(a,a)}\in\Z,
\end{eqnarray*}
and analogously 
$\frac{L(c,s_a^{(0)}(b))}{L(s_a^{(0)}(b),s_a^{(0)}(b))}\in\Z$.

(vii) $I^{(1)}(a,b)=L(b,a)$ because of $L(a,b)=0$. 
\begin{eqnarray*}
L(s_a^{(1)}(b),s_a^{(1)}(b))
&=& L(b-\frac{L(b,a)}{L(a,a)}a,b-\frac{L(b,a)}{L(a,a)}a)\\
&=& L(b,b)-\frac{L(b,a)}{L(a,a)}L(b,a)-0+
(-\frac{L(b,a)}{L(a,a)})^2L(a,a)\\
&=& L(b,b).
\end{eqnarray*}
For $c\in H_\Z$
\begin{eqnarray*}
\frac{L(s_a^{(1)}(b),c)}{L(s_a^{(1)}(b),s_a^{(1)}(b))}
&=& \frac{L(b-\frac{L(b,a)}{L(a,a)}a,c)}{L(b,b)}\\
&=& \frac{L(b,c)}{L(b,b)}-\frac{L(b,a)}{L(b,b)}
\frac{L(a,c)}{L(a,a)}\in\Z,
\end{eqnarray*}
so $s_a^{(1)}(b)\in R^{(0)}$. Finally
\begin{eqnarray*}
s_a^{(1)}s_b^{(1)}(s_a^{(1)})^{-1}
&=& s_a^{(1)}s_{b/\sqrt{L(b,b)}}^{[1]}(s_a^{(1)})^{-1}\\
&\stackrel{\textup{Lemma \ref{t2.2} (c)}}{=}&
s_{s_a^{(1)}(b/\sqrt{L(b,b)})}^{[1]}
=s_{s_a^{(1)}(b)/\sqrt{L(b,b)}}^{[1]}\\
&=& s_{s_a^{(1)}(b)/\sqrt{L(s_a^{(1)}(b),s_a^{(1)}(b))}}^{[1]}
=s_{s_a^{(1)}(b)}^{(1)}.
\end{eqnarray*}

(b) Starting with $S\in T^{tri}_n$, one can define
$H_\Z:=M_{n\times 1}(\Z)$ with standard $\Z$-basis 
$\uuuu{e}=(e_1,...,e_n)$, and one can define
$L:H_\Z\times H_\Z\to\Z$ by $L(\uuuu{e}^t,\uuuu{e})=S^t$.
Then $\uuuu{e}\in \BB^{tri}$.

If $(H_\Z,L)$ is unimodular then $\pm 1=\det 
L(\uuuu{e}^t,\uuuu{e})=L(e_1,e_1)... L(e_n,e_n)$
and $L(e_i,e_i)\in \N$, so $L(e_i,e_i)=1$ and
$L(\uuuu{e}^t,\uuuu{e})^t\in T^{uni}_n(\Z)$. 
The rest is trivial.

(c) The inclusion $R^{(0)}\supset\{a\in H_\Z\,|\, L(a,a)=1\}$
is obvious. Consider $\uuuu{e}\in\BB^{tri}$.
By part (b), the matrix
$S:=L(\uuuu{e}^t,\uuuu{e})^t$ is in $T^{uni}_n(\Z)$.
Consider a root 
$a=\sum_{i=1}^n \alpha_ie_i\in R^{(0)}$. Then
\begin{eqnarray*}
(L(a,e_1),...,L(a,e_n))&=&(\alpha_1,...,\alpha_n)S,\\
\textup{ so }
\gcd(L(a,e_1),...,L(a,e_n))&=&\gcd(\alpha_1,...,\alpha_n).
\end{eqnarray*}
But $L(a,a)$ divides $\gcd(L(a,e_1),...,L(a,e_n))$ because 
$a$ is a root. Therefore $L(a,a)$ divides each $\alpha_i$.
Thus $L(a,a)^2$ divides 
$\sum_{i=1}^n\sum_{j=1}^n\alpha_is_{ij}\alpha_j=L(a,a)$,
so $L(a,a)=1$. 

(d) By part (c) $L(a,a)=1$ for $a\in R^{(0)}$. 
\hfill$\Box$ 

\bigskip
Up to now, the triangular shape of the matrix 
$L(\uuuu{e}^t,\uuuu{e})^t\in T^{tri}_n$ has not been used.
It leads to the result in Theorem \ref{t2.7}. In algebraic geometry
and the theory of meromorphic differential equations, this result 
is well known, it is a piece of Picard-Lefschetz theory.
In the frame of singularity theory, it is treated in 
\cite{AGV88} and \cite{Eb01}. 
An elementary direct proof for a unimodular lattice is given
in \cite{BH20}. The case $k=0$ is proved in \cite{HK16}.

\begin{theorem}\label{t2.7}
Let $(H_\Z,L,\uuuu{e})$ be a bilinear lattice with a triangular
basis $\uuuu{e}$. Let $k\in\{0,1\}$. Then
\begin{eqnarray*}
s_{e_1}^{(k)}...s_{e_n}^{(k)}=(-1)^{k+1}M.
\end{eqnarray*}
\end{theorem}

{\bf Proof:}
The case $k=0$ is a special case
of Proposition 2.4 in \cite{HK16}. The case $k=1$ can be 
proved by an easy modification of Lemma 2.3 (5) and 
Proposition 2.4 in \cite{HK16}. Both cases are proved for
a unimodular bilinear lattice in \cite[Theorem 4.1]{BH20}.
\hfill$\Box$ 

\bigskip
In Picard-Lefschetz theory and singularity theory
also the following notions are standard.

\begin{definition}\label{t2.8}
Let $(H_\Z,L,\uuuu{e})$ be a bilinear lattice with a 
triangular basis $\uuuu{e}$. It induces several structures:
\index{$\Gamma^{(0)},\ \Gamma^{(1)}$}
\index{$\Delta^{(0)},\ \Delta^{(1)}$} 
\begin{list}{}{}
\item[(a)]
The {\it \index{even monodromy group}even monodromy group} 
$\Gamma^{(0)}:=\langle s_{e_1}^{(0)},...,s_{e_n}^{(0)}\rangle
\subset O^{(0)}.$
\item[(b)]
The {\it \index{odd monodromy group}odd monodromy group} 
$\Gamma^{(1)}:=\langle s_{e_1}^{(1)},...,s_{e_n}^{(1)}\rangle
\subset O^{(1)}.$
\item[(c)]
The set $\Delta^{(0)}:=\Gamma^{(0)}\{\pm e_1,...,\pm e_n\}
\subset H_\Z$ of {\it \index{even vanishing cycle}even vanishing cycles}.
\item[(d)]
The set $\Delta^{(1)}:=\Gamma^{(1)}\{\pm e_1,...,\pm e_n\}
\subset H_\Z$ of {\it \index{odd vanishing cycle}odd vanishing cycles}.
\end{list}
\end{definition}

\begin{remarks}\label{t2.9}
(i) The even vanishing cycles are roots, i.e. 
$\Delta^{(0)}\subset R^{(0)}$, because of Lemma \ref{t2.6} (a) 
(vi). In general $\Delta^{(1)}\not\subset R^{(0)}$. 
The name {\it \index{vanishing cycle}vanishing cycles} for the elements of $\Delta^{(0)}$
and $\Delta^{(1)}$ and the name {\it \index{monodromy group}monodromy group} 
stem from singularity theory.
In \cite{HK16} the elements of $\Delta^{(0)}$ are called
{\it real roots}. $\Gamma^{(1)}$ and $\Delta^{(1)}$ are not
considered in \cite{HK16}.

(ii) A matrix $S\in T^{tri}_n$ or $T^{uni}_n(\Z)$ determines
by Lemma \ref{t2.6} (b) a bilinear lattice $(H_\Z,L,\uuuu{e})$
with a triangular basis (up to isomorphism). This leads to the
program to determine for a given matrix $S$ the data $I^{(0)},I^{(1)},G_\Z,G_\Z^{(0)},
G_\Z^{(1)},G_\Z^M,\Gamma^{(0)},\Gamma^{(1)},\Delta^{(0)}$
and $\Delta^{(1)}$. One should start with relevant invariants 
like $\sign I^{(0)}$, $\Rad I^{(0)}$, 
$\Rad I^{(1)}$, the characteristic polynomial and the Jordan
normal form of $M$. 

(iii) The odd monodromy group $\Gamma^{(1)}$ arises naturally in 
many cases where $(H_\Z,L)$ is a unimodular bilinear lattice,
for example in cases from isolated hypersurface singularities.
But it is not clear whether it is natural in cases where
$(H_\Z,L)$ is a bilinear lattice which is not unimodular.
Theorem \ref{t2.7} is positive evidence. But the following
is negative evidence. The monodromy group 
\begin{eqnarray*}
\Gamma^{(1)}=\langle s^{[1]}_{e_i/\sqrt{L(e_i,e_i)}}\,|\,
i\in\{1,...,n\}\rangle
\end{eqnarray*}
contains because of Lemma \ref{t2.2} (c) all transvections
$s^{[1]}_{g(e_i)/\sqrt{L(e_i,e_i)}}$ for $g\in \Gamma^{(1)}$.
Only in the unimodular cases these coincide with
the transvections $s^{[1]}_a$ for $a\in\Delta^{(1)}$. 
We will only consider the unimodular cases.

(iv) We will work on this program rather systematically
in the chapters \ref{s5} and \ref{s6} 
for $S\in T^{uni}_2(\Z)$ and $S\in T^{uni}_3(\Z)$.
We will discuss in chapter \ref{s9} 
known results for matrices $S\in T^{uni}_n(\Z)$
from singularity theory, with emphasis on the simple
singularities (ADE-type, i.e. the ADE root lattices)
and the simple elliptic singularities (i.e. 
$\www{E_6},\www{E_7},\www{E_8}$).
\end{remarks}

Definition \ref{t2.10} and Lemma \ref{t2.11} discuss the case
when a unimodular bilinear lattice $(H_\Z,L,\uuuu{e})$ with
triangular basis is {\it reducible}. Then also the 
monodromy groups, the set of roots 
and the sets of vanishing cycles split.
But beware that here reducibility involves not only
$(H_\Z,L)$, but also $\uuuu{e}$.

\begin{definition}\label{t2.10}
(a) Let $(H_\Z,L,\uuuu{e})$ be a unimodular bilinear lattice
of rank $n\in\N$ with a triangular basis $\uuuu{e}$.
Let $\{1,...,n\}=I_1\ \dot\cup\  I_2$ be a decomposition into 
disjoint subsets such that
\begin{eqnarray*}
L(e_i,e_j)=L(e_j,e_i)=0\quad\textup{for }i\in I_1,\ j\in I_2.
\end{eqnarray*}
Then the triple $(H_\Z,L,\uuuu{e})$ is called 
{\it \index{reducible triple}reducible}.
If such a decomposition does not exist the triple is called
{\it \index{irreducible triple}irreducible}.

(b) A matrix $S\in T^{uni}_n(\Z)$ is called reducible if
the triple $(H_\Z,L,\uuuu{e})$ (which is unique up to isomorphism) 
is reducible, where $(H_\Z,L)$ is a unimodular bilinear lattice 
and $\uuuu{e}$ is a triangular basis with $S=L(\uuuu{e}^t,\uuuu{e})^t$.
\end{definition}

\begin{lemma}\label{t2.11}
Keep the situation of Definition \ref{t2.10}.
For $l\in\{1,2\}$ let $\sigma_l:\{1,2,...,|I_l|\}\to I_l$
be the unique bijection with $\sigma_l(i)<\sigma_l(j)$ for
$i<j$. Define
\begin{eqnarray*}
\uuuu{e}_l&:=&(e_{1,l},e_{2,l},...,e_{|I_l|,l}):=
(e_{\sigma_l(1)},e_{\sigma_l(2)},...,e_{\sigma_l(|I_l|}),\\
H_{\Z,l}&:=&\bigoplus_{i=1}^{|I_l|}\Z\cdot e_{i,l},\quad
L_l:=L|_{H_{\Z,l}}.
\end{eqnarray*}
Then $(H_{\Z,l},L_l,\uuuu{e}_l)$ is a unimodular bilinear
lattice with triangular basis. The decomposition
$H_\Z=H_{\Z,1}\oplus H_{\Z,2}$ is left and right 
\index{orthogonal}$L$-orthogonal.
Denote by $\Gamma^{(0)}_l$, $\Gamma^{(1)}_l$, $\Delta^{(0)}_l$,
$\Delta^{(1)}_l$ and $R^{(0)}_l$ the monodromy groups and sets of
vanishing cycles and roots of $(H_{\Z,l},L_l,\uuuu{e}_l)$.
Denote by $\www{M}_l$ the automorphism of
$H_\Z$ which extends the monodromy $M_l$ on $H_{\Z,l}$
by the identity on $H_{\Z,m}$, where $\{l,m\}=\{1,2\}$. Then 
\begin{eqnarray*}
\Gamma^{(k)}&=& \Gamma^{(k)}_1\times \Gamma^{(k)}_2,\\
R^{(0)}&=& R^{(0)}_1\ \dot\cup\ R^{(0)}_2,\\
\Delta^{(k)}&=& \Delta^{(k)}_1\ \dot\cup \ \Delta^{(k)}_2,\\
M&=& \www{M}_1\www{M}_2=\www{M}_2\www{M}_1.
\end{eqnarray*}
\end{lemma}

The proof is trivial. Because of this lemma, we will study
the monodromy groups and the sets of vanishing cycles only for
irreducible triples. In the Examples \ref{t1.1} this
excludes the cases $S(A^2_1)$, $S(A^3_1)$, $S(A_2A_1)$, 
$S(\P^1A_1)$ and all cases $S(x_1,x_2,x_3)$ where two of the
three numbers $x_1,x_2,x_3$ are zero. 

The following lemma treats the cases $S(A_1^n):=E_n$ for
$n\in\N$. It is a trivial consequence of the special case
$S(A_1)=(1)\in M_{1\times 1}(\Z)$ and Lemma \ref{t2.11}, 
but worth to be stated.

\begin{lemma}\label{t2.12}
The case \index{$A_1^n$}$A_1^n$ for any $n\in\N$: 
\begin{eqnarray*}
H_\Z=\bigoplus_{i=1}^n\Z\cdot e_i,\quad S=S(A_1^n):=E_n,\quad
I^{(0)}=2L,\quad I^{(1)}=0,
\end{eqnarray*}
the reflections 
$s_{e_i}^{(0)}\textup{ with }
s_{e_i}^{(0)}(e_j)=\left\{\begin{array}{ll}
e_j&\textup{if }j\neq i,\\ -e_i&\textup{if }j=i,\end{array}\right\}$
commute, the transvections $s_{e_i}^{(1)}$ are 
$s_{e_i}^{(1)}=\id$,
\begin{eqnarray*}
\Gamma^{(0)}&=&\{ \prod_{i=1}^n
(s_{e_i}^{(0)})^{l_i}\,|\, (l_1,...,l_n)\in\{0;1\}^n\}
\cong \{\pm 1\}^n,\\
\Gamma^{(1)}&=&\{\id\},\\
\Delta^{(0)}&=&R^{(0)}=\{\pm e_1,...,\pm e_n\}=\Delta^{(1)}.
\end{eqnarray*}
\end{lemma}


\chapter{Braid group actions}\label{s3}
\setcounter{equation}{0}
\setcounter{figure}{0}

In the sections \ref{s3.2}--\ref{s3.4} a unimodular
bilinear lattice $(H_\Z,L)$ of some rank $n\geq 2$ is 
considered. The braid group 
$\Br_n$ is introduced in
section \ref{s3.1}. It acts on several sets of $n$-tuples
and of matrices associated to $(H_\Z,L)$. 

Section \ref{s3.1} starts with the Hurwitz action
on $G^n$ where $G$ is a group. Results
of Artin, Birman-Hilden and Igusa-Schiffler for
$G$ a free group, a free Coxeter group or any Coxeter
group are cited and applied.
This is relevant as many of the monodromy groups
$\Gamma^{(1)}$ and $\Gamma^{(0)}$ of rank 2 or rank 3
unimodular bilinear lattices with triangular bases 
are such groups.

It turns out that the Hurwitz action of $\Br_n$ on $(O^{(k)})^n$ 
lifts to an action of a semidirect product 
$\Br_n\ltimes\{\pm 1\}^n$ on sets of certain 
$n$-tuples of cycles in $H_\Z$. This is discussed
in section \ref{s3.2}. 

Most important is the set $\BB^{tri}$ of triangular
bases of $(H_\Z,L)$ (if this set is not empty)
and the subset $\BB^{dist}=\Br_n\ltimes\{\pm 1\}^n(\uuuu{e})$
of a chosen triangular basis $\uuuu{e}$.
Section \ref{s3.3} poses questions on the characterization
of such a set $\BB^{dist}$ of {\it distinguished bases} 
which will guide our work in chapter \ref{s7}. 
It also offers several examples with quite different properties.

Section \ref{s3.4} connects the group $\Br_n\ltimes\{\pm 1\}^n$
via its action on the orbit of a triangular basis $\uuuu{e}$
with the group $G_\Z$. There is a group antihomomorphism
$Z:(\Br_n\ltimes\{\pm 1\}^n)_S\to G_\Z$, where 
$(\Br_n\ltimes\{\pm 1\}^n)_S$ denotes the stabilizer of a
matrix $S\in T^{uni}_n(\Z)$. In this way certain braids induce
automorphisms in $G_\Z$, and in many cases these automorphisms
generate $G_\Z$, i.e. $Z$ is surjective.
Theorem \ref{t3.26} (b) states the well known fact
$Z((\delta^{1-k}\sigma^{root})^n)=(-1)^{k+1}M$ for $k\in\{0;1\}$. 
Theorem \ref{t3.26} (c) gives a condition when
$Z(\delta^{1-k}\sigma^{root})$ is in $G_\Z$ and thus an
$n$-th root of $(-1)^{k+1}M$. 
Theorem \ref{t3.28} states that for almost all cases with 
rank $\leq 3$ the map $Z$ is surjective. The exceptions
are only four cases.

\section[The braid group and the Hurwitz action]
{The braid group and the Hurwitz action, some classical results}
\label{s3.1}

Choose $n\in\Z_{\geq 2}$. 
The \index{braid group}braid group \index{$\Br_n$}$\Br_n$ 
of braids with $n$ strings was introduced by Artin \cite{Ar25}. 
It is the fundamental group of a configuration space. 
We will come to this geometric point of view in chapter \ref{s8}.
Here we take a purely algebraic point of view.
Artin \cite[Satz 1]{Ar25} showed that $\Br_n$ is generated 
by $n-1$ elementary braids $\sigma_1,...,\sigma_{n-1}$ 
\index{$\sigma_1,...,\sigma_{n-1}$} 
and that all relations come from the relations
\begin{eqnarray*}
\sigma_i\sigma_j&=&\sigma_j\sigma_i\quad
\textup{for }i,j\in\{1,...,n-1\}\textup{ with }|i-j|\geq 2,\\
\sigma_i\sigma_{i+1}\sigma_i&=&\sigma_{i+1}\sigma_i\sigma_{i+1}
\quad\textup{for }i\in\{1,...,n-2\}.
\end{eqnarray*}
He also showed \cite[Theorem 19]{Ar47} that the 
\index{center}center of $\Br_n$ is
\begin{eqnarray*}
\textup{Center}(\Br_n)=\langle \sigma^{mon}\rangle,
\end{eqnarray*}
where \index{$\sigma^{root},\ \sigma^{mon}$} 
\begin{eqnarray*}
\sigma^{root}:=\sigma_{n-1}\sigma_{n-2}...\sigma_2\sigma_1,\quad
\sigma^{mon}:=(\sigma^{root})^n.
\end{eqnarray*}

An important action of $\Br_n$ is the 
{\it \index{Hurwitz action}Hurwitz action}
on the $n$-th power $G^n$ for any group $G$. 
The braid group $\Br_n$ acts via
\begin{eqnarray*}
\sigma_i(g_1,...,g_n)&:=& (g_1,...,g_{i-1},g_ig_{i+1}g_i^{-1},g_i, 
g_{i+2},...,g_n),\\
\sigma_i^{-1}(g_1,...,g_n)&:=& (g_1,...,g_{i-1},g_{i+1},
g_{i+1}^{-1}g_ig_{i+1},g_{i+2},...,g_n).
\end{eqnarray*}
The fibers of the map
\begin{eqnarray*}
\pi_n:G^n\to G,\quad \uuuu{g}=(g_1,...,g_n)\mapsto g_1...g_n,
\end{eqnarray*}
are invariant under this action,
\begin{eqnarray*}
\pi_n(\uuuu{g})=\pi_n(\sigma_i\uuuu{g})
=\pi_n(\sigma_i^{-1}\uuuu{g}).
\end{eqnarray*}
We will study this action for $n=3$ in the cases of the 
monodromy groups for the rank 3 unimodular bilinear lattices.
The following results in Theorem \ref{t3.2} of Artin \cite{Ar25} 
and Birman-Hilden \cite{BH73} will be relevant.

\begin{definition}\label{t3.1}
(a) Let $G^{free,n}$ be the \index{free group}
\index{$G^{free,n},\ G^{fCox,n}$}
free group 
with $n$ generators $x_1,...,x_n$. Let
$$\Delta(G^{free,n}):=\bigcup_{i=1}^n\{wx_iw^{-1}\,|\, 
w\in G^{free,n}\}$$
be the set of elements conjugate to $x_1,...,x_n$. 
Obviously $\Br_n((x_1,...,x_n))\subset\Delta(G^{free,n})^n$.

(b) Let $G^{fCox,n}$ be the 
\index{free Coxeter group}free Coxeter group with $n$
generators $x_1,...,x_n$, so all relations are generated
by the relations $x_1^2=...=x_n^2=e$. Let
$$\Delta(G^{fCox,n}):=\bigcup_{i=1}^n\{wx_iw^{-1}\,|\, 
w\in G^{fCox,n}\}$$
be the set of elements conjugate to $x_1,...,x_n$. 
Obviously $\Br_n((x_1,...,x_n))\subset\Delta(G^{fCox,n})^n$.
\end{definition}

\begin{theorem}\label{t3.2}
(a) \cite[Satz 7 and Satz 9]{Ar25}
$\Br_n$ acts simply transitively on the set of
tuples $$\{(w_1,...,w_n)\in \Delta(G^{free,n})^n\,|\, 
w_1...w_n=x_1...x_n\}.$$

(b) \cite[Theorem 7]{BH73} 
$\Br_n$ acts simply transitively on the set of
tuples $$\{(w_1,...,w_n)\in \Delta(G^{fCox,n})^n\,|\, 
w_1...w_n=x_1...x_n\}.$$
\end{theorem}

\begin{remarks}\label{t3.3}
(i) Both results were reproved by Kr\"uger in 
\cite[Satz 7.6]{Kr90}.

(ii) Theorem 1.31 in \cite{KT08} gives a weaker version
of Artin's result Theorem \ref{t3.2} (a). Theorem 1.31 in 
\cite{KT08} is equivalent to the statement that $\Br_n$ acts
simply transitively on the set of tuples
\begin{eqnarray*}
&&\{(w_1,...,w_n)\in \Delta(G^{free,n})^n\,|\, 
w_1...w_n=x_1...x_n,\\
&&\hspace*{1cm} w_1,...,w_n\textup{ generate }G^{free,n},\ 
\textup{ a permutation}\\
&&\hspace*{1cm} \sigma\in S_n\textup{ exists with }
w_i\textup{ conjugate to }x_{\sigma(i)}\}.
\end{eqnarray*}

(iii) The formulation of Theorem 1.31 
in \cite{KT08} is different. 
There a group automorphism $\varphi$ of $G^{free,n}$ is called a
{\it \index{braid automorphism}braid automorphism} 
if $\varphi(x_1...x_n)=x_1...x_n$ 
and if a permutation $\sigma\in S_n$ with 
$\varphi(x_i)$ conjugate to $x_{\sigma(i)}$ exists. 
The group of all braid automorphisms is called
$\www{\Br}_n$. Theorem 1.31 in \cite{KT08} states that the map
\begin{eqnarray*}
\www{Z}:\{\sigma_1^{\pm 1},...,\sigma_{n-1}^{\pm 1}\}
\to \www{Br}_n\quad\textup{with}\\
(\www{Z}(\sigma_i)(x_1),...,\www{Z}(\sigma_i)(x_n))
= \sigma_i^{-1}(x_1,...,x_n),\\
(\www{Z}(\sigma_i^{-1})(x_1),...,\www{Z}(\sigma_i^{-1})(x_n))
= \sigma_i(x_1,...,x_n),
\end{eqnarray*}
extends to a group isomorphism $\Br_n\to \www{\Br}_n$. 

(iv) In order to understand the equivalence of the statements
in (ii) and (iii), it is crucial to see that the extension
$\www{Z}:\Br_n\to\www{\Br}_n$ which is defined by 
\begin{eqnarray*}
(\www{Z}(\beta)(x_1),...,\www{Z}(\beta)(x_n))
=\beta^{-1}(x_1,...,x_n)\quad\textup{for }\beta\in\Br_n
\end{eqnarray*}
is a group homomorphism. This follows from the equations
\begin{eqnarray*}
\www{Z}(\beta\sigma_i)=\www{Z}(\beta)\www{Z}(\sigma_i)
\quad\textup{and}\quad 
\www{Z}(\beta\sigma_i^{-1})=\www{Z}(\beta)\www{Z}(\sigma_i^{-1})
\end{eqnarray*}
for $\beta\in\Br_n$ and $i\in\{1,....n-1\}$. 
The first equation holds because of
\begin{eqnarray*}
&&(\www{Z}(\beta\sigma_i)(x_1),...,
\www{Z}(\beta\sigma_i)(x_n))\\
&=& (\beta\sigma_i)^{-1}(x_1,...,x_n)\\
&=&\sigma_i^{-1}(\beta^{-1}(x_1,...,x_n))\\
&=& \sigma_i^{-1}(\www{Z}(\beta)(x_1),...,\www{Z}(\beta)(x_n))\\
&=& (\www{Z}(\beta)(x_1),...,\www{Z}(\beta)(x_{i-1}),
\www{Z}(\beta)(x_{i+1}),\\
&& (\www{Z}(\beta)(x_{i+1}))^{-1}\www{Z}(\beta)(x_i)
\www{Z}(\beta)(x_{i+1}),
\www{Z}(\beta)(x_{i+2}),...,\www{Z}(\beta)(x_n))\\
&=& (\www{Z}(\beta)(x_1),...,\www{Z}(\beta)(x_{i-1}),
\www{Z}(\beta)(x_{i+1}),\\
&& \www{Z}(\beta)(x_{i+1}^{-1}x_ix_{i+1}),
\www{Z}(\beta)(x_{i+2}),...,\www{Z}(\beta)(x_n))\\
&=& (\www{Z}(\beta)\www{Z}(\sigma_i)(x_1),...,
\www{Z}(\beta)\www{Z}(\sigma_i)(x_n)).
\end{eqnarray*}
The second equation is proved by a similar calculation.
\end{remarks}

\begin{example}\label{t3.4}
Theorem \ref{t3.2} will be applied in the following situation,
which in fact arises quite often. 

Consider a unimodular bilinear lattice $(H_\Z,L,\uuuu{e})$
of rank $n\geq 2$ with triangular basis $\uuuu{e}$ such that 
for some $k\in\{0,1\}$ the following holds:
\begin{eqnarray*}
\Gamma^{(k)}=\left\{\begin{array}{lll}
G^{fCox,n}&\textup{ with generators }s_{e_1}^{(0)},...,
s_{e_n}^{(0)}&\textup{ if }k=0,\\
G^{free,n}&\textup{ with generators }s_{e_1}^{(1)},...,
s_{e_n}^{(1)}&\textup{ if }k=1.
\end{array}\right.
\end{eqnarray*}
Then in the notation of Definition \ref{t3.1}
\begin{eqnarray*}
\Delta(G^{fCox,n})=
\{s_v^{(0)}\,|\,v\in\Delta^{(0)}\}&& \textup{if }k=0,\\
\Delta(G^{free,n})=
\{s_v^{(1)}\,|\,v\in\Delta^{(1)}\}&& \textup{if }k=1.
\end{eqnarray*}
By Theorem \ref{t3.2}, two statements hold:
\begin{list}{}{}
\item[(1)] 
The set 
\begin{eqnarray*}
\{(s_{v_1}^{(k)},...,s_{v_n}^{(k)})\,|\, v_1,...,v_n\in\Delta^{(k)}, 
s_{v_1}^{(k)}...s_{v_n}^{(k)}=(-1)^{k+1}M\}
\end{eqnarray*}
is a single orbit under the Hurwitz action of $\Br_n$. 
\item[(2)]
The stabilizer of any such tuple $(s_{v_1}^{(k)},...,s_{v_n}^{(k)})$
under the Hurwitz action of $\Br_n$ is $\{\id\}$. 
\end{list}
\end{example}

Theorem \ref{t3.2} (b) concerns a free Coxeter group with $n$ generators. 
The transitivity of the action generalizes to arbitrary Coxeter groups
and can be applied to generalize the statement (1) in Example \ref{t3.4},
as is explained in the following.

\begin{definition}\label{t3.5}
(Classical, e.g \cite[5.1]{Hu90}) A 
{\it \index{Coxeter system}Coxeter system}
$(W,S^{gen})$ consists of a group $W$ and a finite set
$S^{gen}=\{s_1,...,s_n\}\subset W$ for some $n\in\N$ of
generators of the group 
such that there are generating relations as follows.
There is a subset $I\subset\{(i,j)\in\{1,...,n\}^2\,|\, i<j\}$
and a map $a:I\to\Z_{\geq 2}$ such that the generating relations
are
\begin{eqnarray*}
s_1^2=...=s_n^2=1,\quad 1=(s_is_j)^{a(i,j)}\textup{ for }
(i,j)\in I.
\end{eqnarray*}
The group $W$ is then called a {\it \index{Coxeter group}Coxeter group}.
\end{definition}

The following theorem was proved by Deligne \cite{De74}
for the ADE Weyl groups and in general by Igusa and Schiffler
\cite{IS10}. A short proof was given by Baumeister,
Dyer, Stump and Wegener \cite{BDSW14}.

\begin{theorem}\label{t3.6}
\cite{De74}\cite{IS10}\cite{BDSW14}
Let $(W,S^{gen})$ with $S^{gen}=\{s_1,...,s_n\}$ be a Coxeter
system with $n\geq 2$. Define 
$\Delta(W,S^{gen}):=\bigcup_{i=1}^n \{ws_iw^{-1}\,|\, w\in W\}$. The set
\begin{eqnarray*}
\{(w_1,...,w_n)\in \Delta(W,S^{gen})^n\,|\, w_1...w_n=s_1...s_n\}
\end{eqnarray*}
is a single orbit under the Hurwitz action of $\Br_n$.
\end{theorem}

Part (a) of the following theorem is classical if 
$S_{ij}\in\{0,-1,-2\}$ for $i<j$ and due to Vinberg
\cite{Vi71} in the general case.

\begin{theorem}\label{t3.7}
Let $(H_\Z,L)$ be a unimodular bilinear lattice of rank
$n\geq 2$ and let $\uuuu{e}$ be a triangular basis
such that the matrix $S=L(\uuuu{e}^t,\uuuu{e})^t\in T^{uni}_n(\Z)$
satisfies $S_{ij}\leq 0$ for $i<j$.

(a) (Classical for $S_{ij}\in\{0,-1,-2\}$, 
\cite[Proposition 6, Theorem 1, Theorem 2,
Proposition 17]{Vi71} for $S_{ij}\leq 0$) 
The pair $(\Gamma^{(0)},\{s_{e_1}^{(0)},...,s_{e_n}^{(0)}\})$
is a Coxeter system with
\begin{eqnarray*}
I&=&\{(i,j)\in\{1,...,n\}^2\,|\, i<j,S_{ij}\in\{0,-1\}\}
\quad\textup{and}\\
a(i,j)&=&\left\{\begin{array}{ll}
2&\textup{ if }S_{ij}=0,\\ 3&\textup{ if }S_{ij}=-1.
\end{array}\right.
\end{eqnarray*}

(b) The set 
\begin{eqnarray*}
\{(g_1,...,g_n)\in\bigl(\{s_v^{(0)}\,|\, 
v\in\Delta^{(0)}\}\bigr)^n\,|\, g_1...g_n=-M\}
\end{eqnarray*}
is a single orbit under the Hurwitz action of $\Br_n$. 
\end{theorem}

{\bf Proof of part (b):} Observe
\begin{eqnarray*}
\Delta(\Gamma^{(0)},\{s_{e_1}^{(0)},...,s_{e_n}^{(0)}\})
=\{s_v^{(0)}\,|\, v\in\Delta^{(0)}\}.
\end{eqnarray*}
Apply Theorem \ref{t3.6}.
\hfill$\Box$ 

\begin{remarks}\label{t3.8}
(i) The transitivity result in part (b) holds also for a 
bilinear lattice $(H_\Z,L)$ which is not necessarily unimodular,
if it comes equipped with a triangular basis $\uuuu{e}$
with $L(e_i,e_j)\leq 0$ for $i>j$. This is the case of a generalized
Cartan lattice (Remark \ref{t2.4} (v)). This is crucial in 
\cite{HK16}. 

(ii) Especially in the case of a unimodular bilinear lattice
$(H_\Z,L,\uuuu{e})$ with triangular basis $\uuuu{e}$ and matrix
$S=L(\uuuu{e}^t,\uuuu{e})^t\in T^{uni}_n(\Z)$ with $S_{ij}\leq -2$
for all $i<j$ we have $\Gamma^{(0)}=G^{fCox,n}$ with generators 
$s_{e_1}^{(0)},...,s_{e_n}^{(0)}$. 

(iii) Theorem \ref{t6.11} (g) gives in the case $n=3$
$\Gamma^{(0)}=G^{fCox,3}$ with generators 
$s_{e_1}^{(0)},s_{e_2}^{(0)},s_{e_3}^{(0)}$ also in the following
cases: if $S_{ij}\geq 3$
for $i<j$ and if additionally 
$$2S_{12}\leq S_{13}S_{23},\quad 
2S_{13}\leq S_{12}S_{23},\quad
2S_{23}\leq S_{12}S_{13}.$$

(iv) Theorem \ref{t6.18} (g) gives in the situation of part (iii)
also $\Gamma^{(1)}=G^{free,3}$ with generators 
$s_{e_1}^{(1)},s_{e_2}^{(1)},s_{e_3}^{(1)}$.

(v) Though in the situation of part (ii) there are cases with
$\Gamma^{(1)}= G^{free,n}$, and there are cases with
$\Gamma^{(1)}\neq G^{free,n}$. 
The odd cases are more complicated than the even cases. 
For the cases with $n=3$ see the Remarks \ref{t4.17},
Lemma \ref{t4.18} and Theorem \ref{t6.18}.
\end{remarks}

\section{Braid group action on tuples of cycles}
\label{s3.2}

Consider a unimodular bilinear lattice $(H_\Z,L)$ of some rank
$n\geq 2$ and the groups $O^{(k)}=\Aut(H_\Z,I^{(k)})$ 
for $k\in\{0;1\}$. Recall that here the set of roots $R^{(0)}$ 
is $$R^{(0)}=\{\delta\in H_\Z\,|\, L(\delta,\delta)=1\}.$$ 
In order to treat the even case $k=0$ and the odd case $k=1$
uniformly, we define  \index{$R^{(1)}$} 
$$R^{(1)}:=H_\Z.$$
The Hurwitz action of $\Br_n$ on $(O^{(k)})^n$ restricts 
because of 
\begin{eqnarray}\label{3.1}
s_a^{(k)}s_b^{(k)}(s_a^{(k)})^{-1}
=s_{s_a^{(k)}(b)}^{(k)}\quad\textup{for}\quad a,b\in R^{(k)}
\end{eqnarray}
(Lemma \ref{t2.2} (c)) 
to an action on the subset $(\{s_v^{(k)}\,|\, v\in R^{(k)}\})^n$.

It turns out that this action has a natural lift to an action
of a certain semidirect product $\Br_n\ltimes\{\pm 1\}^n$
on the set $(R^{(k)})^n$. Here the sets $(R^{(k)})^n$ and
$(\{s_v^{(k)}\,|\, v\in R^{(k)}\})^n$ are related by the map
\begin{eqnarray*}
\pi_n^{(k)}: (R^{(k)})^n&\to& (\{s_v^{(k)}\,|\, v\in R^{(k)}\})^n
\subset (O^{(k)})^n, \\
\quad \uuuu{v}=(v_1,...,v_n)&\mapsto& (s_{v_1}^{(k)},...,s_{v_n}^{(k)}).
\end{eqnarray*}
\index{$\pi_n^{(k)},\ \pi_n$}
Recall also the map 
$$\pi_n: (O^{(k)})^n\to O^{(k)},\quad (g_1,...,g_n)\mapsto g_1...g_n$$
which was defined for an arbitrary group before Definition \ref{t3.1}.

Furthermore it turns out that both actions, for $k=0$ and for
$k=1$, restrict to the same action on the set $\BB^{tri}$
of triangular bases if this set is not empty.
This is the action in which we are interested most. 
In the case of a unimodular bilinear lattice from singularity theory,
it is well known \cite[5.7]{Eb01}
\cite[\S 1.9]{AGV88} and has been studied by 
A'Campo, Brieskorn, Ebeling, Gabrielov, Gusein-Zade, 
Kr\"uger and others.
In fact it works also for bilinear lattices with
$\BB^{tri}\neq\emptyset$ which are not necessarily unimodular,
see Remark \ref{t3.14}. 

Finally, the actions induce actions of $\Br_n\ltimes \{\pm 1\}^n$
on several spaces of matrices. 
The purpose of this section is to fix all these well known actions.

Lemma \ref{t3.9} presents the semidirect product 
$\Br_n\ltimes\{\pm 1\}^n$.
Lemma \ref{t3.10} gives its action on $(R^{(k)})^n$.
Lemma \ref{t3.11} gives its action on $\BB^{tri}$ if this
set is not empty.

\begin{lemma}\label{t3.9}
Fix $n\in\Z_{\geq 2}$. 

(a) The multiplicative group $\{\pm 1\}^n$ is called
{\it sign group}. It is generated by the elements
$\delta_j=((-1)^{\delta_{ij}})_{i=1,...,n}\in\{\pm 1\}^n$
(here $\delta_{ij}$ is the Kronecker symbol)
for $j\in\{1,...,n\}$.

(b) The following relations define a 
\index{semidirect product}semidirect product 
\index{$\Br_n\ltimes\{\pm 1\}^n$}$\Br_n\ltimes\{\pm 1\}^n$ 
of $\Br_n$ and $\{\pm 1\}^n$
with $\{\pm 1\}^n$ as normal subgroup,
\begin{eqnarray*}
\sigma_j\delta_i\sigma_j^{-1}&=& \delta_i\quad
\textup{for }i\in\{1,...,n\}-\{j,j+1\},\\
\sigma_j\delta_j\sigma_j^{-1}&=&\delta_{j+1},\quad
\sigma_j\delta_{j+1}\sigma_j^{-1}=\delta_j.
\end{eqnarray*}
In the following $\Br_n\ltimes\{\pm 1\}^n$ always means
this semidirect product.
\end{lemma}

{\bf Proof:} 
Part (a) is a notation. Part (b) requires a proof.
We have the exact sequence
\begin{eqnarray*}
\{1\}\to \Br_n^{pure}\to\Br_n\to S_n\to\{1\}
\end{eqnarray*}
where $\Br_n^{pure}\subset \Br_n$ is the normal subgroup
of pure braids (and $\sigma_i\in \Br_n$ maps to the 
transposition $(i\ i+1)\in S_n$). See Remark \ref{t8.4} (vii)
for this group and this exact sequence. 
The natural action of $S_n$ on $\{\pm 1\}^n$,
\begin{eqnarray*}
S_n\owns \alpha: \uuuu{\varepsilon}
=(\varepsilon_1,...,\varepsilon_n)\mapsto
(\varepsilon_{\alpha^{-1}(1)},...,\varepsilon_{\alpha^{-1}(n)})
=:\alpha . \uuuu{\varepsilon}
\end{eqnarray*}
lifts to an action of $\Br_n$ on $\{\pm 1\}^n$,
$\sigma:\uuuu{\varepsilon}\mapsto \sigma .\uuuu{\varepsilon}.$
This action can be used to define a semidirect product of
$\Br_n$ and $\{\pm 1\}^n$ by 
$\sigma\uuuu{\varepsilon}\sigma^{-1}:=\sigma .\uuuu{\varepsilon}$.
It is the semidirect product in part (b).\hfill$\Box$

\begin{lemma}\label{t3.10}
Let $(H_\Z,L)$ be a unimodular bilinear lattice of rank
$n\geq 2$. Fix $k\in\{0,1\}$.

(a) The following formulas define an action of the 
semidirect product $\Br_n\ltimes \{\pm 1\}^n $ from
Definition \ref{t3.9} (b) on the set $(R^{(k)})^n$,
\begin{eqnarray*}
\sigma_j(\uuuu{v})&=& (v_1,...,v_{j-1},
s_{v_j}^{(k)}(v_{j+1}),v_j,v_{j+2},...,v_n),\\
\sigma_j^{-1}(\uuuu{v})&=& (v_1,...,v_{j-1},
v_{j+1},(s_{v_{j+1}}^{(k)})^{-1}(v_j),v_{j+2},...,v_n),\\
\delta_j(\uuuu{v})&=& (v_1,...,v_{j-1},-v_j,v_{j+1},...,v_n),
\end{eqnarray*}
for $\uuuu{v}=(v_1,...,v_n)\in (R^{(k)})^n$.

(b) The map $\pi_n^{(k)}:(R^{(k)})^n\to (O^{(k)})^n$
is compatible with the action of $\Br_n\ltimes \{\pm 1\}^n$
on $(R^{(k)})^n$ from part (a) and the Hurwitz action
of $\Br_n$ on $(O^{(k)})^n$, so the diagram
\begin{eqnarray*}
\begin{CD}
(R^{(k)})^n @>{\sigma_j}>>  (R^{(k)})^n\\
@V{\pi^{(k)}_n}VV  @VV{\pi^{(k)}_n}V \\
(O^{(k)})^n @>{\sigma_j}>> (O^{(k)})^n 
\end{CD}
\end{eqnarray*}
commutes. Here the sign group $\{\pm 1\}^n$
acts trivially on $(O^{(k)})^n$.
Especially, each orbit in $(R^{(k)})^n$ is contained in one fiber
of the projection $\pi_n\circ\pi_n^{(k)}:(R^{(k)})^n\to O^{(k)}$.
\end{lemma}

{\bf Proof:} (a) 
We denote the actions in part (a) by 
$\sigma_j^{(k)},(\sigma_j^{-1})^{(k)}$ and $\delta_j^{(k)}$
(of course $\delta_j^{(0)}=\delta_j^{(1)}$). The identities
\begin{eqnarray*}
\sigma_j^{(k)}(\sigma_j^{-1})^{(k)}
=(\sigma_j^{-1})^{(k)}\sigma_j^{(k)}=\id&&
\textup{for }j\in\{1,...,n-1\}\\
\sigma_i^{(k)}\sigma_j^{(k)}=\sigma_j^{(k)}\sigma_i^{(k)}&&
\textup{for }|i-j|\geq 2,\\
\sigma_j^{(k)}\delta_i^{(k)}(\sigma_j^{-1})^{(k)}=\delta_i^{(i)}&&
\textup{for }i\in\{1,...,n\}-\{j,j+1\},\\
\sigma_j^{(k)}\delta_j^{(k)}(\sigma_j^{-1})^{(k)}
=\delta_{j+1}^{(k)}&&\textup{and}\\
\sigma_j^{(k)}\delta_{j+1}^{(k)}(\sigma_j^{-1})^{(k)}
=\delta_j^{(k)}&&\textup{for }j\in\{1,...,n-1\}
\end{eqnarray*}
are obvious or easy to see. The identities
\begin{eqnarray*}
\sigma_i^{(k)}\sigma_{i+1}^{(k)}\sigma_i^{(k)}
=\sigma_{i+1}^{(k)}\sigma_i^{(k)}\sigma_{i+1}^{(k)}
&&\textup{for }i\in\{1,...,n-2\}
\end{eqnarray*}
are proved by the following calculation with $\uuuu{v}\in (R^{(k)})^n$, 
\begin{eqnarray*}
&& \sigma_i^{(k)}\sigma_{i+1}^{(k)}\sigma_i^{(k)}(\uuuu{v})\\
&=& \sigma_i^{(k)}\sigma_{i+1}^{(k)}(...,v_{i-1},s^{(k)}_{v_i}(v_{i+1}),v_i,
v_{i+2},v_{i+3}...)\\
&=& \sigma_i^{(k)}(...,v_{i-1},s^{(k)}_{v_i}(v_{i+1}),
s^{(k)}_{v_i}(v_{i+2}),v_i,v_{i+3},...)\\
&=& (...,v_{i-1},s^{(k)}_{s^{(k)}_{v_i}(v_{i+1})}
(s^{(k)}_{v_i}(v_{i+2})),s^{(k)}_{v_i}(v_{i+1}),v_i,v_{i+3},...)\\
&\stackrel{\eqref{3.1}}{=}&(...,v_{i-1},s^{(k)}_{v_i}s^{(k)}_{v_{i+1}}
(s^{(k)}_{v_i})^{-1}(s^{(k)}_{v_i}(v_{i+2})),
s^{(k)}_{v_i}(v_{i+1}),v_i,v_{i+3},...)\\
&=& (...,v_{i-1},s^{(k)}_{v_i}(s^{(k)}_{v_{i+1}}(v_{i+2})),
s^{(k)}_{v_i}(v_{i+1}),v_i,v_{i+3},...)\\
&=&\sigma_{i+1}^{(k)}(...,v_{i-1},s^{(k)}_{v_i}
(s^{(k)}_{v_{i+1}}(v_{i+2})),v_i,v_{i+1},v_{i+3},...)\\
&=&\sigma_{i+1}^{(k)}\sigma_i^{(k)}(...,v_{i-1},v_i,
s^{(k)}_{v_{i+1}}(v_{i+2}),v_{i+1},v_{i+3},...)\\
&=&\sigma_{i+1}^{(k)}\sigma_i^{(k)}\sigma_{i+1}^{(k)}(\uuuu{v}).
\end{eqnarray*}
The maps $\sigma_j^{(k)}$ and $\delta_i^{(k)}$ satisfy all relations
between the generators $\sigma_j$ and $\delta_i$ of the group
$\Br_n\ltimes\{\pm 1\}^n$. Therefore the formulas in part (a)
define an action of this group on the set $(R^{(k)})^n$.

(b) The actions are compatible because of \eqref{3.1}.
The sign group acts trivially on $(O^{(k)})^n$ because
$s_v^{(k)}=s_{-v}^{(k)}$ for $v\in R^{(k)}$. Each orbit of the
Hurwitz action on $(O^{(k)})^n$ is contained in one fiber of the
map $\pi_n$, as was remarked in section \ref{s3.1}.\hfill$\Box$

\begin{lemma}\label{t3.11}
Let $(H_\Z,L)$ be a unimodular bilinear lattice with nonempty set
$\BB^{tri}$ of triangular bases. The actions in Lemma \ref{t3.10}
of $\Br_n\ltimes\{\pm 1\}^n$ on $(R^{(0)})^n$ and on $(R^{(1)})^n$ 
both restrict to the same action on $\BB^{tri}$. 
This action can also be written as follows,
\begin{eqnarray*}
\sigma_j(\uuuu{v})&=& (v_1,...,v_{j-1},
v_{j+1}-L(v_{j+1},v_j)v_j,v_j,v_{j+2},...,v_n),\\
\sigma_j^{-1}(\uuuu{v})&=& (v_1,...,v_{j-1},
v_{j+1},v_j-L(v_{j+1},v_j)v_{j+1},v_{j+2},...,v_n),\\
\delta_j(\uuuu{v})&=& (v_1,...,v_{j-1},-v_j,v_{j+1},...,v_n),
\end{eqnarray*}
for $\uuuu{v}=(v_1,...,v_n)\in \BB^{tri}$.
\end{lemma}

{\bf Proof:}
$\uuuu{v}\in\BB^{tri}$ implies $L(v_j,v_{j+1})=0$ and
$2L(v_j,v_j)=2=I^{(0)}(v_j,v_j)$. Recall $I^{(0)}=L+L^t$ and
$I^{(1)}=L^t-L$. Therefore 
\begin{eqnarray*}
s_{v_j}^{(k)}(v_{j+1})
&=& v_{j+1}-I^{(k)}(v_j,v_{j+1})v_j
=v_{j+1}-L(v_{j+1},v_j)v_j,\\
(s_{v_{j+1}}^{(k)})^{-1}(v_j)
&=& v_j-(-1)^kI^{(k)}(v_{j+1},v_j)v_{j+1}
=v_j-L(v_{j+1},v_j)v_{j+1}.
\end{eqnarray*}
So $\sigma_j(\uuuu{v})$ and $\sigma_j^{-1}(\uuuu{v})$ are given
by the formulas in Lemma \ref{t3.11}. It remains to see that the
images are again in $\BB^{tri}$. They are in $(R^{(0)})^n$
because of the even case $k=0$. 
They form $\Z$-bases of $H_\Z$ because $\uuuu{v}$ is a $\Z$-basis
of $H_\Z$. They are triangular bases because
\begin{eqnarray*}
L(\sigma_j(\uuuu{v})_j,\sigma_j(\uuuu{v})_{j+1})
&=& L(v_{j+1}-L(v_{j+1},v_j)v_j,v_j)=0,\\
L(\sigma_j^{-1}(\uuuu{v})_j,\sigma_j^{-1}(\uuuu{v})_{j+1})
&=& L(v_{j+1},v_j-L(v_{j+1},v_j)v_{j+1})=0.
\end{eqnarray*}
Of course $\delta_j(\uuuu{v})\in\BB^{tri}$.\hfill$\Box$

\begin{definition}\label{t3.12}
Fix $n\in\N$ and $R\in\{\Z,\Q,\R,\C\}$.

(a) Recall the definition of the set $T^{uni}_n(R)$ 
of upper triangular $n\times n$ matrices with entries 
in $R$ and 1's on the diagonal in Definition \ref{t2.3} (c). 
Additionally we define the
sets of symmetric and skewsymmetric matrices
\begin{eqnarray*} 
T^{(0)}_n(R)&:=& \{A\in M_{n\times n}(R)\,|\, A^t=A,
A_{ii}=2\},\\
T^{(1)}_n(R)&:=& \{A\in M_{n\times n}(R)\,|\, A^t=-A\}.
\end{eqnarray*}

(b) For $a\in\Z$ define the $n\times n$ matrix
\begin{eqnarray*}
C_{n,j}(a)&=&
\begin{pmatrix}
1 & & & & & \\ & \ddots & & & & \\ 
& & a & 1 & & \\
& & 1 & 0 & & \\ & & & & \ddots & \\ & & & & & 1 \end{pmatrix}
\end{eqnarray*}
which differs from the unit matrix only in the positions
$(j,j),(j,j+1),(j+1,j),(j+1,j+1)$. 
Its inverse is
\begin{eqnarray*}
C_{n,j}^{-1}(a)&=&
\begin{pmatrix}
1 & & & & & \\ & \ddots & & & & \\
& & 0 & 1 & & \\ 
& & 1 & -a & & \\
& & & & \ddots & \\ & & & & & 1 \end{pmatrix}
\end{eqnarray*}
with $-a$ at the position $(j+1,j+1)$. 
\end{definition}

\begin{lemma}\label{t3.13}
Fix $n\in\Z_{\geq 2}$ and $R\in\{\Z,\Q,\R,\C\}$. 

(a) The following formulas define an action of the semidirect
product $\Br_n\ltimes\{\pm 1\}^n$ on each of the sets of matrices
$T^{(0)}_n(R)$, $T^{(1)}_n(R)$ and $T^{uni}_n(R)$,
\begin{eqnarray*}
\sigma_j(A)&=& C_{n,j}(-A_{j,j+1})\cdot A\cdot C_{n,j}(-A_{j,j+1})
\quad\textup{for }j\in\{1,...,n-1\},\\
\sigma_j^{-1}(A)&=& C_{n,j}^{-1}(A_{j,j+1})\cdot A\cdot C_{n,j}^{-1}(A_{j,j+1})
\quad\textup{for }j\in\{1,...,n-1\},\\
\delta_j(A)&=& \diag(((-1)^{\delta_{ij}})_{i=1,...,n})\cdot A\cdot 
\diag(((-1)^{\delta_{ij}})_{i=1,...,n})\\
&&\hspace*{4cm} \textup{for }j\in\{1,...,n\}.
\end{eqnarray*}

(b) Let $(H_\Z,L)$ be a unimodular bilinear lattice of rank $n$. 

(i) Fix $k\in\{0;1\}$. The map
\begin{eqnarray*}
(R^{(k)})^n\to T^{(k)}_n(\Z),\quad \uuuu{v}\mapsto 
I^{(k)}(\uuuu{v}^t,\uuuu{v}),
\end{eqnarray*}
is compatible with 
the actions of $\Br_n\ltimes\{\pm 1\}^n$ on $(R^{(k)})^n$ 
and on $T^{(k)}_n(\Z)$.

(ii) Suppose that $\BB^{tri}$ is not empty. The map
\begin{eqnarray*}
\BB^{tri}\to T^{uni}_n(\Z),\quad \uuuu{v}\mapsto 
L(\uuuu{v}^t,\uuuu{v})^t,
\end{eqnarray*}
is compatible with the actions of $\Br_n\ltimes\{\pm 1\}^n$ 
on $\BB^{tri}$ and on $T^{uni}_n(\Z)$.
\end{lemma}

{\bf Proof:}
(a) For $A\in T^{(k)}_n(\Z)$ there is a unique matrix
$S\in T^{uni}_n(\Z)$ with $A=S+(-1)^kS^t$.
Consider a unimodular bilinear lattice $(H_\Z,L,\uuuu{e})$
with a triangular basis $\uuuu{e}$ with matrix
$S=L(\uuuu{e}^t,\uuuu{e})^t$. Then 
\begin{eqnarray*}
\sigma_j(\uuuu{e})&=& (...,e_{j-1},e_{j+1}-S_{j,j+1}e_j,e_j,e_{j+2},...)\\
&=& \uuuu{e}\cdot C_{n,j}(-S_{j,j+1}),\\
\sigma_j^{-1}(\uuuu{e})&=& 
(...,e_{j-1},e_{j+1},e_j-S_{j,j+1}e_{j+1},e_{j+2},...)\\
&=& \uuuu{e}\cdot C_{n,j}^{-1}(S_{j,j+1}),\\
\delta_j(\uuuu{e})&=& (...,e_{i-1},-e_i,e_{i+1},...)\\
&=& \uuuu{e}\cdot \diag(((-1)^{\delta_{ij}})_{i=1,...,n}).
\end{eqnarray*}
Observe that $C_{n,j}(a)$ for $a\in\Z$ is symmetric. Therefore for $k\in\{0;1\}$
\begin{eqnarray*}
I^{(k)}(\sigma_j(\uuuu{e})^t,\sigma_j(\uuuu{e}))
&=& C_{n,j}(-S_{j,j+1})\cdot I^{(k)}(\uuuu{e}^t,\uuuu{e})
\cdot C_{n,j}(-S_{j,j+1}),\\
I^{(k)}(\sigma_j^{-1}(\uuuu{e})^t,\sigma_j^{-1}(\uuuu{e}))
&=& C_{n,j}^{-1}(S_{j,j+1})\cdot I^{(k)}(\uuuu{e}^t,\uuuu{e})
\cdot C_{n,j}^{-1}(S_{j,j+1}),
\end{eqnarray*}
similarly for $L$ instead of $I^{(k)}$, 
and also similarly for the action of $\delta_j$.
This shows part (a) for $R=\Z$. 
Changing the set of scalars does not change the matrix
identities which say that the group $\Br_n\ltimes\{\pm 1\}^n$ acts. 

(b) This follows from the proof of part (a).\hfill$\Box$

\begin{remarks}\label{t3.14}
Let $(H_\Z,L)$ be a bilinear lattice, not necessarily unimodular.

(i) The action of $\Br_n\ltimes\{\pm 1\}^n$ on $(R^{(0)})^n$ in 
Lemma \ref{t3.10} (a) works also in this case.
It restricts as in Lemma \ref{t3.11} to an action on $\BB^{tri}$,
if this set is not empty, though here for $\uuuu{v}\in\BB^{tri}$
\begin{eqnarray*}
\sigma_j(\uuuu{v})&=& (v_1,...,v_{j-1},
v_{j+1}-\frac{L(v_{j+1},v_j)}{L(v_j,v_j)}v_j,v_j,v_{j+2},...,v_n),\\
\sigma_j^{-1}(\uuuu{v})&=& (v_1,...,v_{j-1},
v_{j+1},v_j-\frac{L(v_{j+1},v_j)}{L(v_{j+1},v_{j+1})}
v_{j+1},v_{j+2},...,v_n).
\end{eqnarray*}

(ii) The action of $\Br_n\ltimes\{\pm 1\}^n$ 
in Lemma \ref{t3.10} on $(R^{(1)})^n$ does not generalize.
In Lemma \ref{t2.6} (a) (v) we defined in the case of a general
bilinear lattice $s_a^{(1)}$ only for $a\in R^{(0)}$.
We defined $s_a^{(1)}$ for any $a\in R^{(1)}$ in Lemma \ref{t2.6} (d)
only in the case of a unimodular bilinear lattice.

(iii) On the other hand, part (a) (vii) of Lemma \ref{t2.6} says 
that the action in Lemma \ref{t3.10} for $k=1$ works for 
$\uuuu{v}\in \BB^{tri}$. 
Though at the end this is just the action in (i) above. 

(iv) The action in (i) on $\BB^{tri}$ is compatible with an
action on the set $T^{tri}_n$ of matrices in Lemma \ref{t2.3} (c), 
which generalizes the action in Lemma \ref{t3.13} (a). Here 
$C_{n,j}(-S_{j,j+1})$ and $C_{n,j}^{-1}(S_{j,j+1})$ 
in Lemma \ref{t3.13} (a) have to be replaced by
$C_{n,j}(-\frac{S_{j,j+1}}{S_{j,j}})$ and 
$C_{n,j}^{-1}(\frac{S_{j,j+1}}{S_{j+1,j+1}})$.
\end{remarks}

$\{\pm 1\}^n$ is the normal subgroup in the semidirect
product $\Br_n\ltimes\{\pm 1\}^n$. Therefore, if 
$\Br_n\ltimes\{\pm 1\}^n$ acts on some set $\Sigma$,
the group $\Br_n$ acts on the quotient $\Sigma/\{\pm 1\}^n$.
Often it is good to consider this quotient and the action
of $\Br_n$ on it.

\begin{lemma}\label{t3.15}
Let $(H_\Z,L)$ be a unimodular bilinear lattice of some rank
$n\geq 2$. Fix $k\in\{0;1\}$.

(a) The map
\begin{eqnarray*}
\pi_n^{(k)}:(R^{(k)})^n\to (\{s_v^{(k)}\,|\, v\in R^{(k)}\})^n
\subset (O^{(k)})^n,\quad 
\uuuu{v}\mapsto (s^{(k)}_{v_1},...,s^{(k)}_{v_n}),
\end{eqnarray*}
factors into maps
\begin{eqnarray*}
\begin{CD}
(R^{(k)})^n @>>> (R^{(k)})^n/\{\pm 1\}^n
@>{\pi_n^{(k)}/\{\pm 1\}^n}>>
(\{s_v^{(k)}\,|\, v\in R^{(k)}\})^n.
\end{CD}
\end{eqnarray*}
$\Br_n$ acts on the quotient $(R^{(k)})^n/\{\pm 1\}^n$, and 
the second map $\pi_n^{(k)}/\{\pm 1\}^n$ is $\Br_n$ equivariant.
The image of $\uuuu{v}$ in $(R^{(k)})^n/\{\pm 1\}^n$ is
denoted by $\uuuu{v}/\{\pm 1\}^n$. 

(b) The second map $$\pi_n^{(k)}/\{\pm 1\}^n:
(R^{(k)})^n/\{\pm 1\}^n\to (\{s_v^{(k)}\,|\, v\in R^{(k)}\})^n$$ 
in part (a) is a bijection if 
$k=0$ or if $k=1$ and $\Rad I^{(1)}=\{0\}$.

(c) Consider the case $k=1$ and $\Rad I^{(1)}\supsetneqq\{0\}$.
Consider a triangular basis $\uuuu{e}\in\BB^{tri}$ and the 
induced set $\Delta^{(1)}$ of odd vanishing cycles. 
The second map restricts to a $\Br_n$ equivariant bijection
\begin{eqnarray*}
(\Delta^{(1)})^n/\{\pm 1\}^n\to (\{s_v^{(1)}\,|\, v\in\Delta^{(1)}\})^n,\quad
\uuuu{v}\mapsto (s_{v_1}^{(1)},...,s_{v_n}^{(1)}), 
\end{eqnarray*}
if $(H_\Z,L,\uuuu{e})$ is either irreducible or reducible with
at most one summand of type $A_1$. 
\end{lemma}

{\bf Proof:}
Part (a) is trivial. (b) Suppose $k=0$ or ($k=1$ and $\Rad I^{(1)}=\{0\}$).
If $k=1$ and $v=0$ then $s_v^{(1)}=\id$. If $k=0$ and $v\in R^{(0)}$ or
if $k=1$ and $v\in R^{(1)}-\{0\}$ then  
$v\notin\Rad I^{(k)}$ and $s_v^{(k)}\neq\id$ 
for any $v\in R^{(k)}$. Then one can recover $\pm v$ from $s_v^{(k)}$,
essentially because of
\begin{eqnarray*}
\{0\}\subsetneqq (s_v^{(k)}-\id)(H_\Z)\subset\Z v.
\end{eqnarray*}

(c) If $(H_\Z,L,\uuuu{e})$ is irreducible then 
$\Delta^{(1)}\cap \Rad I^{(1)}=\emptyset$, 
and the argument of part (b) holds. 
If $(H_\Z,L,\uuuu{e})$ is reducible with only one summand
of type $A_1$ then for a unique $j\in\{1,...,n\}$
$e_j\in \Rad I^{(1)}$, and then 
$\Delta^{(1)}\cap \Rad I^{(1)}=\{\pm e_j\}$. 
Then $v\in\Delta^{(1)}$ satisfies $s_v^{(1)}=\id$
if and only if $v=\pm e_j$. So also then
one can recover $\pm v$ from $s_v^{(1)}$ for 
any $v\in\Delta^{(1)}$.\hfill$\Box$

\begin{remarks}\label{t3.16}
(i) Consider the action of $\Br_n$ on $T^{uni}_n(\Z)$. 
The elementary braid $\sigma_j$ maps
$S=(S_{ij})\in T^{uni}_n(\Z)$ to 
\begin{eqnarray*}
\sigma_j(S)&=&
C_{n,j}(-S_{j,j+1})\cdot S\cdot C_{n,j}(-S_{j,j+1})\\
\textup{ with }&&
\begin{pmatrix} \sigma_j(S)_{jj} & \sigma_j(S)_{j,j+1} \\
\sigma_j(S)_{j+1,j} & \sigma_j(S)_{j+1,j+1}\end{pmatrix}
=\begin{pmatrix} 1 & -S_{j,j+1} \\ 0 & 1 \end{pmatrix}.
\end{eqnarray*}

(ii) Especially, in the case $n=2$ $\delta_1$, $\delta_2$
and $\sigma_1$ all map
$S=\begin{pmatrix}1 & x\\0&1\end{pmatrix}$ to 
$\begin{pmatrix} 1& -x\\0&1\end{pmatrix}$,
so the $\Br_2\ltimes\{\pm 1\}^2$ orbit equals the 
$\Br_2$ orbit and the $\langle \delta_1\rangle$ orbit
and consists only of $S=\begin{pmatrix}1 & x\\0&1\end{pmatrix}$ and 
$\begin{pmatrix}1 &-x\\0&1\end{pmatrix}$. 

(iii) Under rather special circumstances, also in higher
rank $n$ the sign group action is eaten up by the braid group
action. Ebeling proved the following lemma.
\end{remarks}

\begin{lemma}\label{t3.17}
Let $(H_\Z,L)$ be a unimodular bilinear lattice and 
$\uuuu{e}\in\BB^{tri}$ a triangular basis with
$S=L(\uuuu{e}^t,\uuuu{e})^t\in T^{uni}_n(\Z)$. 

(a) \cite[proof of Prop. 2.2]{Eb83} 
Suppose $S_{j,j+1}=\varepsilon\in\{\pm 1\}$ for some
$j\in\{1,...,n-1\}$. Then
\begin{eqnarray*}
\sigma_j^{3\varepsilon}(\uuuu{e})=\delta_j(\uuuu{e}),\quad
\sigma_j^{-3\varepsilon}(\uuuu{e})=\delta_{j+1}(\uuuu{e}).
\end{eqnarray*}

(b) \cite[Prop. 2.2]{Eb83}
Suppose $S_{ij}\in\{0,1,-1\}$ for all $i,j$ and that
$(H_\Z,L)$ is irreducible. 
Then the orbit of $\uuuu{e}$ under $\Br_n$ coincides with
the orbit of $\uuuu{e}$ under $\Br_n\ltimes\{\pm 1\}^n$. 
\end{lemma}

\section{Distinguished bases}
\label{s3.3}

In this section we continue the discussion of the braid group
action on $\BB^{tri}$. Now we fix one triangular basis $\uuuu{e}$.
Definition \ref{t3.18} gives notations.
The Remarks \ref{t3.19} pose questions on the orbit of 
$\uuuu{e}$ under $\Br_n\ltimes\{\pm 1\}^n$.
The questions will guide the work which will be done in chapter \ref{s7}.

\begin{definition}\label{t3.18}
Let $(H_\Z,L)$ be a unimodular lattice of rank $n\geq 2$ 
with nonempty set $\BB^{tri}$ of triangular bases.

Given a triangular basis $\uuuu{e}\in\BB^{tri}$ 
and $k\in\{0,1\}$, we are interested in the following orbits:
\index{$\BB^{dist}$}\index{$\RR^{(k),dist}$}\index{$\SSS^{dist}$}
\begin{eqnarray*}
\textup{the set }\BB^{dist}
&:=&\Br_n\ltimes\{\pm 1\}^n(\uuuu{e})
\textup{ of {\it \index{distinguished basis}distinguished bases}},\\
\textup{the set }\RR^{(k),dist}
&:=&\Br_n(\pi_n^{(k)}(\uuuu{e}))
\textup{ of {\it 
\index{distinguished tuple of reflections}distinguished tuples of}}\\
&&\textup{\it reflections or transvections},
\\
\textup{the set }\SSS^{dist}
&:=&\Br_n\ltimes\{\pm 1\}^n(S)
\textup{ of {\it \index{distinguished matrix}distinguished matrices},}\\
&&\textup{where }S=L(\uuuu{e}^t,\uuuu{e})^t\in T^{uni}_n(\Z),
\end{eqnarray*}
the quotient sets $\BB^{dist}/\{\pm 1\}^n$ and 
$\SSS^{dist}/\{\pm 1\}^n$, which are $\Br_n$ orbits.
We are also interested in the \index{stabilizer}stabilizers 
in $\Br_n$ of the points $\uuuu{e}/\{\pm 1\}^n\in \BB^{dist}/\{\pm 1\}^n$
and $S/\{\pm 1\}^n\in \SSS^{dist}/\{\pm 1\}^n$, namely the groups
\index{$(\Br_n)_{\uuuu{e}/\{\pm 1\}^n},\ (\Br_n)_{S/\{\pm 1\}^n}$}
\begin{eqnarray*}
(\Br_n)_{\uuuu{e}/\{\pm 1\}^n}
\subset (\Br_n)_{S/\{\pm 1\}^n} \subset \Br_n.
\end{eqnarray*}
\end{definition}

\begin{remarks}\label{t3.19}
In the situation of Definition \ref{t3.18}
the following constraints on the set $\BB^{dist}$ 
of distinguished bases are clear from what
has been said,
\begin{eqnarray}\nonumber
\BB^{dist}&\subset& \BB^{tri}\cap (\Delta^{(0)})^n 
\cap (\Delta^{(1)})^n \cap 
\{\uuuu{v}\in (H_\Z)^n\,|\, \sum_{i=1}^n\Z v_i=H_\Z\}\\
&&\cap \ (\pi_n\circ\pi_n^{(0)})^{-1}(-M)
\cap (\pi_n\circ\pi_n^{(1)})^{-1}(M),\label{3.2}
\end{eqnarray}
where $\pi_n\circ\pi_n^{(k)}:(R^{(k)})^n\to O^{(k)},\quad
\uuuu{v}\mapsto s_{v_1}^{(k)}...s_{v_n}^{(k)}$. 

An interesting problem is which - if any - 
of these constraints
are sufficient to characterize the orbit $\BB^{dist}$.
We are most interested in the questions whether the inclusions
\begin{eqnarray}\label{3.3}
\BB^{dist}&\subset& \{\uuuu{v}\in (\Delta^{(0)})^n\,|\, 
\pi_n\circ\pi_n^{(0)}(\uuuu{v})=-M\},\\
\BB^{dist}&\subset& \{\uuuu{v}\in (\Delta^{(1)})^n\,|\, 
\pi_n\circ\pi_n^{(1)}(\uuuu{v})=M\},\label{3.4}
\end{eqnarray}
are equalities. 
We will study this problem systematically in chapter \ref{s7} 
for $n=2$ and $n=3$.
In this section \ref{s3.3} we give some examples.
\end{remarks}

\begin{remarks}\label{t3.20}
Consider a unimodular bilinear lattice $(H_\Z,L)$ 
of rank $n\geq 2$ and $k\in\{0;1\}$. 

(i) Two basic invariants of the $\Br_n\ltimes \{\pm 1\}^n$
orbit of a tuple $\uuuu{v}\in (R^{(k)})^n$ are the product
$(\pi_n\circ\pi_n^{(k)})(\uuuu{v}))
=s_{v_1}^{(k)}...s_{v_n}^{(k)}\in O^{(k)}$ 
and the sublattice $\sum_{i=1}^n\Z v_i\subset H_\Z$, 
namely 
\begin{eqnarray*}
(\pi_n\circ\pi_n^{(k)})(\sigma_j({\uuuu{v}}))
=(\pi_n\circ\pi_n^{(k)})(\uuuu{v})
\quad\textup{and}\quad
\sum_{i=1}^n\Z \sigma_j(\uuuu{v})_i=\sum_{i=1}^n\Z v_i.
\end{eqnarray*}

(ii) A triangular basis $\uuuu{e}\in \BB^{tri}$ induces
the even and odd monodromy groups $\Gamma^{(0)}$ and $\Gamma^{(1)}$
and the sets $\Delta^{(0)}$ and $\Delta^{(1)}$ 
of even and odd vanishing cycles. 
Each distinguished basis $\uuuu{v}\in\BB^{dist}$
induces the same even and odd monodromy groups 
$\Gamma^{(0)}$ and $\Gamma^{(1)}$ and 
the same sets $\Delta^{(0)}$ and $\Delta^{(1)}$ 
of even and odd vanishing cycles. 
This is obvious from the action of
$\Br_n\ltimes \{\pm 1\}^n$ on $\BB^{dist}$,
the Hurwitz action of $\Br_n$ on $\RR^{(k),dist}$
and the definition of $\Gamma^{(k)}$ and $\Delta^{(k)}$.
So they are invariants of the set $\BB^{dist}$ of distinguished bases. 

We did not mention the monodromy $M$, because it is by Theorem
\ref{t2.7} an invariant of $(H_\Z,L)$ if $\BB^{tri}\neq\emptyset$, 
so it does not depend on the choice of a $\Br_n\ltimes\{\pm 1\}^n$
orbit in $\BB^{tri}$. 

(iii) A matrix $S\in T^{uni}_n(\Z)$ determines a triple 
$(H_\Z,L,\uuuu{e})$ with $(H_\Z,L)$ a unimodular bilinear lattice
and $\uuuu{e}$ a triangular basis with $S=L(\uuuu{e}^t,\uuuu{e})^t$
up to isomorphism. For a second matrix $\www{S}$ in the
$\Br_n\ltimes\{\pm 1\}^n$ orbit of $S$ then a triangular 
basis $\www{\uuuu{e}}$ of $(H_\Z,L)$ with 
$\www{S}=L(\www{\uuuu{e}}^t,\www{\uuuu{e}})^t$ exists
(but it is not unique in general). Therefore the triple
$(H_\Z,L,\BB^{dist})$ and all induced data depend only
on the $\Br_n\ltimes\{\pm 1\}^n$ orbit $\SSS^{dist}$ of $S$.
These induced data comprise $R^{(0)}$, $\Gamma^{(0)}$,
$\Gamma^{(1)}$, $\Delta^{(0)}$ and $\Delta^{(1)}$. 

(iv) Choose a triangular basis $\uuuu{e}$. 
Then $s_\delta^{(k)}\in \Gamma^{(k)}$ for $\delta\in\Delta^{(k)}$. 
This follows from the definition of a vanishing cycle and 
from formula \ref{3.2}. It also implies that the set 
$(\Delta^{(k)})^n$ is invariant under the action of 
$\Br_n\ltimes\{\pm 1\}^n$ on $(R^{(k)})^n$. 
\end{remarks}

\begin{remarks}\label{t3.21}
(i) Consider a unimodular bilinear lattice which is not 
a lattice of type $A_1^n$ 
and a fixed triangular basis $\uuuu{e}$. 
Then $\Delta^{(0)}\subset R^{(0)}$, but
$\Delta^{(1)}\not\subset R^{(0)}$
by Corollary \ref{t6.22} (a).
Nevertheless $\BB^{dist}\subset (\Delta^{(0)}\cap\Delta^{(1)})^n$, 
so many odd vanishing cycles do not turn up 
in bases in the braid group orbit of $\uuuu{e}$,
i.e. in distinguished bases. 

(ii) In all cases except $A_1^n$ 
$\Delta^{(1)}\not\subset\Delta^{(0)}$ because
$\Delta^{(1)}\not\subset R^{(0)}$.
In some cases $\Delta^{(0)}\subset \Delta^{(1)}$,
in many cases $\Delta^{(0)}\not\subset \Delta^{(1)}$. 
See Corollary \ref{t6.22}.
\end{remarks}

Given a unimodular bilinear lattice $(H_\Z,L)$, 
any element $g\in O^{(k)}$ acts on $(R^{(k)})^n$ by
$g(\uuuu{v}):=(g(v_1),...,g(v_n))$. Part (a) of the
next Lemma \ref{t3.22} says especially that this action
commutes with the action of $\Br_n\ltimes\{\pm 1\}^n$ 
on $(R^{(k)})^n$. Part (b) gives implications, 
which will be used to construct the
interesting Examples \ref{t3.23} (i) and (ii).

\begin{lemma}\label{t3.22}
Let $(H_\Z,L)$ be a unimodular bilinear lattice. 
Fix $k\in \{0,1\}$.

(a) If $g\in O^{(k)}$ and $(\alpha,\varepsilon)\in 
\Br_n\ltimes \{\pm 1\}^n$ then for $\uuuu{v}\in (R^{(k)})^n$
\begin{eqnarray*}
g((\alpha,\varepsilon)(\uuuu{v})) &=&
(\alpha,\varepsilon)(g(\uuuu{v})),\\
(\pi_n\circ\pi_n^{(k)})(g(\uuuu{v})) &=& 
g\circ (\pi_n\circ \pi_n^{(k)})(\uuuu{v})\circ g^{-1}.
\end{eqnarray*}
Here $g(\uuuu{v})$ means $(g(v_1),...,g(v_n))$, and similarly
for $g((\alpha,\varepsilon)(\uuuu{v}))$.

(b) If $g\in G_\Z^{(k)}-G_\Z$ and $\uuuu{e}\in\BB^{tri}$
then 
\begin{eqnarray*}
g(\uuuu{e})&\notin& \BB^{tri},\\ \textup{so especially}\quad 
g(\uuuu{e})&\notin& \Br_n\ltimes\{\pm 1\}^n(\uuuu{e}),\\
\textup{but}\quad
(\pi_n\circ\pi_n^{(k)})(g(\uuuu{e}))&=&(\pi_n\circ\pi_n^{(k)})
(\uuuu{e})=(-1)^{k+1}M.
\end{eqnarray*}
\end{lemma}

{\bf Proof}
(a) $g((\id,\varepsilon)(\uuuu{v}))
=(\id,\varepsilon)(g(\uuuu{v}))$
is trivial. Consider $(\alpha,\varepsilon)=(\sigma_j,
(1,...,1))=\sigma_j$.
\begin{eqnarray*}
g(\sigma_j(\uuuu{v}))
&=&(g(v_1),...,g(v_{j-1}),gs_{v_j}^{(k)}(v_{j+1}),
g(v_j),g(v_{j+2}),...,g(v_n))\\
&=& (g(v_1),...,g(v_{j-1}),s_{g(v_j)}^{(k)}(gv_{j+1}),
g(v_j),g(v_{j+2}),...,g(v_n))\\
&=&\sigma_j(g(\uuuu{v})),
\end{eqnarray*}
because of $gs_{v_j}^{(k)}g^{-1}=s_{g(v_j)}^{(k)}$ (Lemma 
\ref{t2.2} (c)). 
\begin{eqnarray*}
(\pi_n\circ\pi_n^{(k)})(g(\uuuu{v})) &=& 
s_{g(v_1)}^{(k)}...s_{g(v_n)}^{(k)}
=\bigl(gs_{v_1}^{(k)}g^{-1}\bigr)...
\bigl(gs_{v_n}^{(k)}g^{-1}\bigr)\\
&=& gs_{v_1}^{(k)}...s_{v_n}^{(k)}g^{-1}
=g\circ (\pi_n\circ \pi_n^{(k)})(\uuuu{v})\circ g^{-1}.
\end{eqnarray*}

(b) Suppose $g\in G_\Z^{(k)}-G_\Z$ and $\uuuu{e}\in\BB^{tri}$.
Then  $gMg^{-1}=M$ and 
\begin{eqnarray*}
(\pi_n\circ \pi_n^{(k)})(g(\uuuu{e}))
&\stackrel{(a)}{=}&
g\circ(\pi_n\circ\pi_n^{(k)})(\uuuu{e})\circ g^{-1}\\
&=&g\bigl( (-1)^{k+1}M\bigr)g^{-1}
=(-1)^{k+1}M=(\pi_n\circ\pi_n^{(k)})(\uuuu{e}).
\end{eqnarray*}
Furthermore $S=L(\uuuu{e}^t,\uuuu{e})^t\in T^{uni}_n(\Z)$,
$I^{(k)}=L^t+(-1)^{k}L$ and 
\begin{eqnarray*}
I^{(k)}(g(\uuuu{e})^t,g(\uuuu{e}))
&=& I^{(k)}(\uuuu{e}^t,\uuuu{e})=S+(-1)^{k}S^t
\quad\textup{because }g\in G_\Z^{(k)},\\
L(g(\uuuu{e})^t,g(\uuuu{e}))^t&\neq& L(\uuuu{e}^t,\uuuu{e})^t=S
\quad\textup{because }g\notin G_\Z,\\
\textup{so }L(g(\uuuu{e})^t,g(\uuuu{e}))&\notin& T^{uni}_n(\Z),
\end{eqnarray*}
so $g(\uuuu{e})\notin \BB^{tri}$, so 
$g(\uuuu{e})\notin \Br_n\ltimes \{\pm 1\}^n(\uuuu{e})$.
\hfill$\Box$

\begin{examples}\label{t3.23}
Let $(H_\Z,L)$ be a unimodular bilinear lattice of rank $n\geq 2$,
$\uuuu{e}\in \BB^{tri}$ a triangular basis,
$S=L(\uuuu{e}^t,\uuuu{e})^t\in T^{uni}_n(\Z)$
and $k\in\{0,1\}$. 

(i) $k=1$, $n=3$, $S=S(-3,3,-3)=S(\P^2)$, the odd case $\P^2$.
To carry out this example we need two results which will be proved
later, Theorem \ref{t5.14} (b) and Theorem \ref{t6.21} (h). 

By Theorem \ref{t5.14} (b) (i) there is a root $M^{root}\in G_\Z$
of the monodromy with $(M^{root})^3=M$ and 
\begin{eqnarray*}
G_\Z^{(1)}&=&\{\pm (M^{root})^l(\id +a(M^{root}-\id)^2)\,|\,
l,a\in\Z\}\\
&\supsetneqq& G_\Z=\{\pm (M^{root})^l\,|\, l\in \Z\}.
\end{eqnarray*}
Also, here $\Rad I^{(1)}=\Z f_3$ with $f_3=e_1+e_2+e_3$. 
The shape of $M^{root}$ in Theorem \ref{t5.14} (b) shows  
\begin{eqnarray*}
(M^{root}-\id)^2(\uuuu{e})=\uuuu{e}
(\begin{pmatrix}3&-3&1\\1&0&0\\0&1&0\end{pmatrix}-E_3)^2
=f_3(1,-2,1).
\end{eqnarray*}
For example $g:=\id +(M^{root}-\id)^2\in G_\Z^{(1)}-G_\Z$ 
satisfies 
\begin{eqnarray*}
g(\uuuu{e})&=&\uuuu{e}+f_3(1,-2,1),\\
I^{(1)}(g(\uuuu{e})^t,g(\uuuu{e}))
&=& I^{(1)}(\uuuu{e}^t,\uuuu{e})=S-S^t
=\begin{pmatrix}0&-3&3\\3&0&-3\\-3&3&0\end{pmatrix},\\
L(g(\uuuu{e})^t,g(\uuuu{e}))^t
&=&\begin{pmatrix}3&-7&5\\-4&9&-7\\2&-4&3\end{pmatrix},\\
s_{g(e_1)}^{(1)}s_{g(e_2)}^{(1)}s_{g(e_3)}^{(1)}
&=&M,\\
g(e_1),g(e_2),g(e_3)&\notin& R^{(0)},\quad\textup{because }
3\neq 1\textup{ and }9\neq 1,\\
g(e_1),g(e_2),g(e_3)&\notin& \Delta^{(1)}.
\end{eqnarray*}
The last claim $g(e_j)\notin \Delta^{(1)}$
holds because of $g(e_j)\neq e_j$, but 
$g(e_j)\in e_j+\Z f_3$, and because the projection
$H_\Z\to H_\Z/\Z f_3$ restricts to an injective map
$\Delta^{(1)}\to H_\Z/\Z f_3$ by Theorem \ref{t6.21} (h).

(ii) $k=0$, $n=3$, $S=S(-2,2,-2)=S(\HH_{1,2})$, the even case
$\HH_{1,2}$. 
To carry out this example we need two results which will be proved
later, Theorem \ref{t5.14} (a) and Theorem \ref{t6.14} (e). 

Recall Theorem \ref{t5.14} (b) (i),
\begin{eqnarray*}
(H_\Z,L)&=& (H_{\Z,1},L_1)\oplus (H_{\Z,2},L_2)\\
\textup{with }H_{\Z,1}&=& \Z f_1\oplus \Z f_2=\ker\Phi_2(M),\quad
H_{\Z,2}=\Z f_3=\ker \Phi_1(M),\\
(f_1,f_2,f_3)&=& \uuuu{e}
\begin{pmatrix}1&0&1\\1&1&1\\0&1&1\end{pmatrix},\\
G_\Z^{(0)}&=& G_{\Z,1}^{(0)}\times G_{\Z,2}^{(0)}
\supsetneqq G_\Z=G_{\Z,1}\times G_{\Z,2}\\
\textup{with }
G_{\Z,1}^{(0)}&=&\Aut(H_{\Z,1})\supsetneqq 
G_{\Z,1}=\{g\in\Aut(H_{\Z,1})\,|\, \det g=1\},\\
G_{\Z,2}^{(0)}&=&G_{\Z,2}=\{\pm\id|_{\Z f_3}\}.
\end{eqnarray*}
For example $g:=((f_1,f_2,f_3)\mapsto (f_1,-f_2,f_3))
\in G_\Z^{(0)}-G_\Z$ satisfies
\begin{eqnarray*}
g(\uuuu{e})&=& g(f_3-f_2,-f_3+f_1+f_2,f_3-f_1)\\
&=& (f_3+f_2,-f_3+f_1-f_2,f_3-f_1)\\
&=&\uuuu{e}+f_2(2,-2,0),\\
s_{g(e_1)}^{(0)}s_{g(e_2)}^{(0)}s_{g(e_3)}^{(0)}&=&-M,
\end{eqnarray*}
\begin{eqnarray*}
I^{(0)}(g(\uuuu{e})^t,g(\uuuu{e}))
&=&I^{(0)}(\uuuu{e}^t,\uuuu{e})=S+S^t,\\
L(g(\uuuu{e})^t,g(\uuuu{e}))^t&=& 
\begin{pmatrix}1&0&0\\-2&1&0\\2&-2&1\end{pmatrix}
=S^t\neq S 
=L(\uuuu{e}^t,\uuuu{e})^t,\\
g(\uuuu{e})&\notin& \BB^{tri},
\textup{ so }g(\uuuu{e})\notin \BB^{dist},\\
g(e_1),g(e_2),g(e_3)&\in\Delta^{(0)}.
\end{eqnarray*}
The last claim $g(e_j)\in\Delta^{(0)}$ holds because of
Theorem \ref{t6.14} (e). Therefore
\begin{eqnarray*}
g(\uuuu{e})\in 
\{\uuuu{v}\in (\Delta^{(0)})^3\,|\, 
(\pi_3\circ\pi_3^{(0)})(\uuuu{v})=-M,
\sum_{i=1}^3\Z v_i=H_\Z\}.
\end{eqnarray*}
Especially, here the inclusion in \eqref{3.3} is not an
equality. 
In a certain sense, this example is the worst case within
all cases $S(\uuuu{x})$ with $k=0$, $n=3$ and eigenvalues
of $M$ unit roots. See Theorem \ref{t7.3} (b).

(iii) $n\in\N$, $S=E_n$, the even and odd case $A_1^n$. 
Compare Lemma \ref{t2.12}. 
\begin{eqnarray*}
\Delta^{(0)}&=& R^{(0)}=\Delta^{(1)}=\{\pm e_1,...,\pm e_n\},\\
M&=&\id,\\
\BB^{dist}&=& \{(\varepsilon_1e_{\sigma(1)},...,
\varepsilon_ne_{\sigma(n)})\,|\, \varepsilon_1,...,\varepsilon_n
\in \{\pm 1\},\sigma\in S_n\},\\
&=&\{\uuuu{v}\in(\Delta^{(0)})^n\,|\, 
(\pi_n\circ\pi_n^{(0)})(\uuuu{e})=-M=-\id\}
\end{eqnarray*}
The last equality follows from $-M=-\id$ and
\begin{eqnarray*}
s_{e_i}^{(0)}|_{\Z e_i}=-\id,\quad 
s_{e_i}^{(0)}|_{\sum_{j\neq i}\Z e_j}=\id.
\end{eqnarray*}
Here the inclusion in \eqref{3.3} is an equality. 
On the contrary, the inclusion \eqref{3.4} is not 
an equality, the set 
\begin{eqnarray*}
\{\uuuu{v}\in(\Delta^{(1)})^n\,|\, 
(\pi_n\circ\pi_n^{(1)})(\uuuu{v})=M=\id\}
\end{eqnarray*}
is much bigger than $\BB^{dist}$ 
if $n\geq 2$, it consists of many 
$\Br_n\ltimes \{\pm 1\}^n$ orbits. Each of these orbits
contains a unique one of the following tuples,
\begin{eqnarray*}
(e_{l(1)},e_{l(2)},...,e_{l(n)})\quad\textup{with}\quad 
1\leq l(1)\leq l(2)\leq ...\leq l(n)\leq n.
\end{eqnarray*}
This follows from $s_{e_i}^{(1)}=\id$. 

The braids act just by permutations on the $\Br_n$ orbit 
$\BB^{dist}/\{\pm 1\}^n$. Therefore the stabilizer of
$\uuuu{e}/\{\pm 1\}^n$ is the group
$\Br_n^{pure}$ of pure braids
(see Remark \ref{t8.2} (vii) for this group). 
The stabilizer of $S/\{\pm 1\}^n$ is the whole group $\Br_n$. 

(iv) Reconsider Example \ref{t3.4}, so a case where
\begin{eqnarray*}
\Gamma^{(k)}=\left\{\begin{array}{ll}
G^{fCox,n}&\textup{ with generators }s_{e_1}^{(0)},...,
s_{e_n}^{(0)}\textup{ if }k=0,\\
G^{free,n}&\textup{ with generators }s_{e_1}^{(1)},...,
s_{e_n}^{(1)}\textup{ if }k=1.
\end{array}\right.
\end{eqnarray*}
In the notation of Definition \ref{t3.1}
$\Delta(G^{fCox,n})=\{s_\delta^{(0)}\,|\,\delta\in\Delta^{(0)}\}$
if $k=0$ and
$\Delta(G^{free,n})=\{s_\delta^{(1)}\,|\,\delta\in\Delta^{(1)}\}$
if $k=1$. By Theorem \ref{t3.2}
\begin{eqnarray*}
\RR^{dist,(k)}=\{(s_{v_1}^{(k)},...,s_{v_n}^{(k)})\,|\, 
\uuuu{v}\in(\Delta^{(k)})^n,s_{v_1}^{(k)}...s_{v_n}^{(k)}
=(-1)^{k+1}M\}.
\end{eqnarray*}
Furthermore, the shape of $\Gamma^{(k)}$ shows that 
$(H_\Z,L,\uuuu{e})$ is irreducible. Lemma \ref{t3.15} (b) or (c) 
applies. Therefore
\begin{eqnarray*}
\BB^{dist}=\{\uuuu{v}\in (\Delta^{(k)})^n\,|\, 
s_{v_1}^{(k)}...s_{v_n}^{(k)}=(-1)^{k+1}M\}.
\end{eqnarray*}
So in this case only the constraints $\uuuu{v}\in(\Delta^{(k)})^n$
and $s_{v_1}^{(k)}...s_{v_n}^{(k)}=(-1)^{k+1}M$ in
the Remarks \ref{t3.19} are needed in order to characterize the orbit
$\BB^{dist}$. The inclusions in \eqref{3.3} and
\eqref{3.4} are here equalities.

By Theorem \ref{t3.2} the stabilizer of 
$\pi_n^{(k)}(\uuuu{e})$ and of $\uuuu{e}/\{\pm 1\}^n$
is $\{\id\}\subset\Br_n$. 
The size of the stabilizer $(\Br_n)_{S/\{\pm 1\}^n}$ depends
on the case. Theorem \ref{t7.11} gives cases with $n=3$ where
it is $\langle \sigma_2\sigma_1\rangle$ or
$\langle \sigma_2\sigma_1^2\rangle$ or
$\langle \sigma^{mon}\rangle$. 

(v) Suppose $S_{ij}\leq 0$ for $i<j$, so $(H_\Z,L,\uuuu{e})$ 
is a generalized Cartan lattice as in Theorem \ref{t3.7}.
By Theorem \ref{t3.7} (b)
\begin{eqnarray*}
\RR^{dist,(0)}&=& \{(g_1,...,g_n)\in\bigl(\{s_\delta^{(0)}\,|\, 
\delta\in\Delta^{(0)}\}\bigr)^n\,|\, g_1...g_n=-M\}.
\end{eqnarray*}
With Lemma \ref{t3.15} (b) this implies
\begin{eqnarray*}
\BB^{dist}&=& \{\uuuu{v}\in (\Delta^{(0)})^n\,|\, 
s_{v_1}^{(k)}...s_{v_n}^{(k)}=-M\}.
\end{eqnarray*}
Also here only the two constraints $\uuuu{v}\in(\Delta^{(0)})^n$
and $s_{v_1}^{(k)}...s_{v_n}^{(k)}=-M$ in
the Remarks \ref{t3.19} are needed in order to characterize the orbit
$\BB^{dist}$. The inclusion in \eqref{3.3} is here an equality.
\end{examples}

\section{From $\Br_n\ltimes\{\pm 1\}^n$ to $G_\Z$}
\label{s3.4}

Definition \ref{t3.24} gives a map 
$Z:\Br_n\ltimes\{\pm 1\}^n\to \Aut(H_\Z)$, which restricts to a 
group antihomomorphism 
$Z:(\Br_n\ltimes\{\pm 1\}^n)_S\to G_\Z$. 
The definition and the restriction to a group antihomomorphism 
are classical. 
Lemma \ref{t3.25} provides basic facts around this map. 
Also Theorem \ref{t3.26} (b) is classical. It states that
$Z((\delta_n^{1-k}\sigma^{root})^n)=(-1)^{k+1}M$.

Theorem \ref{t3.26} (c) gives a condition when
$Z(\delta_n^{1-k}\sigma^{root})\in G_\Z$. 
Then this is an $n$-th root of $(-1)^{k+1}M$. 
Theorem \ref{t3.26} (c) embraces Theorem 4.5 (a)+(b) in 
\cite{BH20}. It gives more because in \cite{BH20} no braids are
considered. The braids allow a new and more elegant proof
than the one in \cite{BH20}. 
Theorem \ref{t3.26} (c) will be used in the discussion
of the groups $G_\Z$ in many cases in chapter \ref{s5}.

\begin{definition}\label{t3.24}
Let $(H_\Z,L,\uuuu{e})$ be a unimodular bilinear lattice of rank
$n\geq 2$ with a triangular basis $\uuuu{e}$.
For $(\alpha,\varepsilon)\in\Br_n\ltimes \{\pm 1\}^n$ define
an automorphism $Z((\alpha,\varepsilon))\in Aut(H_\Z)$ by the 
following action on the $\Z$-basis $\uuuu{e}$ of $H_\Z$,
\index{$Z$}
\begin{eqnarray*}
Z:\Br_n\ltimes\{\pm 1\}^n&\to& \Aut(H_\Z),\\
(\alpha,\varepsilon)&\mapsto& Z((\alpha,\varepsilon))
=(\uuuu{e}\to (\alpha,\varepsilon)(\uuuu{e})).
\end{eqnarray*}
\end{definition}

\begin{lemma}\label{t3.25}
Let $(H_\Z,L,\uuuu{e})$ be a unimodular bilinear lattice of rank
$n\geq 2$ with a triangular basis $\uuuu{e}$.

(a) For $(\alpha,\varepsilon)\in\Br_n\ltimes\{\pm 1\}^n$
\begin{eqnarray*}
Z((\alpha,\varepsilon))\in G_\Z
\iff Z((\alpha,\varepsilon))\in G_\Z^{(0)}
\iff Z((\alpha,\varepsilon))\in G_\Z^{(1)}.
\end{eqnarray*}

(b) The stabilizer of $S$ in $\Br_n\ltimes\{\pm 1\}^n$ is
\begin{eqnarray*}
(\Br_n\ltimes\{\pm 1\}^n)_S=\{(\alpha,\varepsilon)
\in\Br_n\ltimes\{\pm 1\}^n\,|\, Z((\alpha,\varepsilon))\in G_\Z\}.
\end{eqnarray*}

(c) The restriction of the map $Z$ to the stabilizer
$(\Br_n\ltimes\{\pm 1\}^n)_S$ is also denoted $Z$, 
\begin{eqnarray*}
Z:(\Br_n\ltimes \{\pm 1\}^n)_S\to G_\Z.
\end{eqnarray*}
It is a group antihomomorphism with kernel the stabilizer
$(\Br_n\ltimes\{\pm 1\}^n)_{\uuuu{e}}$ of $\uuuu{e}$.

(d) The triple $(H_\Z,L,\BB^{dist})$ with 
$\BB^{dist}=\Br_n\ltimes\{\pm 1\}^n(\uuuu{e})$ 
the set of distinguished bases 
(Definition \ref{t3.18}) gives rise to the subgroup 
$G_\Z^{\BB}$ of $G_\Z$, \index{$G_\Z^{\BB}$}
\begin{eqnarray*}
G_\Z^{\BB}:=\Aut(H_\Z,L,\BB^{dist})
:=\{g\in G_\Z\,|\, g(\BB^{dist})=\BB^{dist}\}\subset G_\Z.
\end{eqnarray*}
It does not depend on $\uuuu{e}$, but only on the triple 
$(H_\Z,L,\BB^{dist})$. Then $G_\Z^\BB$ is the image of
$(\Br_n\ltimes\{\pm 1\}^n)_S$ under $Z$ in $G_\Z$, 
\begin{eqnarray*}
G_\Z^{\BB}=Z((\Br_n\ltimes\{\pm 1\}^n)_S)\subset G_\Z.
\end{eqnarray*}

(e) The subgroup $Z((\{\pm 1\}^n)_S)$ of $G_\Z^{\BB}$ 
is a normal subgroup of $G_\Z^{\BB}$, and the group 
antihomomorphism $Z$ in part (c) induces a group 
antihomomorphism
\begin{eqnarray*}
\oooo{Z}:(\Br_n)_{S/\{\pm 1\}^n}\to 
G_\Z^{\BB}/Z((\{\pm 1\}^n)_S)
\end{eqnarray*}
with kernel $(\Br_n)_{\uuuu{e}/\{\pm 1\}^n}$,
which is isomorphic to 
$(\Br_n\ltimes \{\pm 1\}^n)_{\uuuu{e}}$. 

(f) Suppose that $(H_\Z,L,\uuuu{e})$ is irreducible
(Definition \ref{t2.10} (a)). Then 
\begin{eqnarray*}
(\{\pm 1\}^n)_S
&=&\{(1,...,1),(-1,...,-1)\},\\
Z((-1,...,-1))&=&-\id\in G_\Z,\\
Z((\{\pm 1\}^n)_S)&=& \{\pm \id\}.
\end{eqnarray*}
$\{\pm \id\}$ is a normal subgroup of $G_\Z$.
The group antihomorphism $\oooo{Z}$ in part (d) becomes
\begin{eqnarray*}
\oooo{Z}:(\Br_n)_{S/\{\pm 1\}^n}\to G_\Z/\{\pm \id\}
\end{eqnarray*}
with kernel $(\Br_n)_{\uuuu{e}/\{\pm 1\}^n}$ and image
$G_\Z^{\BB}/\{\pm \id\}$. 
\end{lemma}

{\bf Proof:}
(a) Fix $k\in\{0,1\}$ and $(\alpha,\varepsilon)\in 
\Br_n\ltimes\{\pm 1\}^n$. Then
\begin{eqnarray*}
Z((\alpha,\varepsilon))\in G_\Z^{(k)}
\iff I^{(k)}((\alpha,\varepsilon)(\uuuu{e})^t,
(\alpha,\varepsilon)(\uuuu{e}))=S+(-1)^{k}S^t.
\end{eqnarray*}
If this equality holds then $I^{(k)}=L^t+(-1)^{k}L$ and 
$L((\alpha,\varepsilon)(\uuuu{e})^t,
(\alpha,\varepsilon)(\uuuu{e}))^t\in T^{uni}_n(\Z)$ imply
$L((\alpha,\varepsilon)(\uuuu{e})^t,
(\alpha,\varepsilon)(\uuuu{e}))^t=S$, so
$Z((\alpha,\varepsilon))\in G_\Z$.

(b) Trivial with the compatibility of the actions of
$\Br_n\ltimes\{\pm 1\}^n$ on $\BB^{dist}$ and on
$T^{uni}_n(\Z)$ in Theorem \ref{t3.4} (d).

(c) The following calculation shows that the map
$Z:(\Br_n\ltimes \{\pm 1\}^n)_S\to G_\Z$ 
is a group antihomomorphism,
\begin{eqnarray*}
Z((\alpha,\varepsilon)(\beta,\www{\varepsilon}))(\uuuu{e})
&=& (\alpha,\varepsilon)(\beta,\www{\varepsilon})(\uuuu{e})\\
&=& (\alpha,\varepsilon)(Z((\beta,\www{\varepsilon}))(e_1),...,
Z((\beta,\www{\varepsilon}))(e_n))\\
&\stackrel{\textup{\ref{t3.22} (a)}}{=}&
Z((\beta,\www{\varepsilon}))(\alpha,\varepsilon)(\uuuu{e})
=Z((\beta,\www{\varepsilon}))Z((\alpha,\varepsilon))(\uuuu{e}).
\end{eqnarray*}
It is trivial that the kernel of this map $Z$ is 
$(\Br_n\ltimes\{\pm 1\}^n)_{\uuuu{e}}$.

(d) $G_\Z\subset Z((\Br_n\ltimes\{\pm 1\}^n)_S)$: 
Consider $g\in G_\Z^{\BB}$. Then $g(\uuuu{e})\in\BB^{dist}$
comes with same matrix $S$ as $\uuuu{e}$ because $g$ respects $L$.
There is a pair 
$(\alpha,\varepsilon)\in \Br_n\ltimes\{\pm 1\}^n$ with 
$Z((\alpha,\varepsilon))(\uuuu{e})
=(\alpha,\varepsilon)(\uuuu{e})=g(\uuuu{e})$.
Therefore $g=Z((\alpha,\varepsilon))$. 

$G_\Z\supset Z((\Br_n\ltimes\{\pm 1\}^n)_S)$: 
Consider $(\alpha,\varepsilon)\in (\Br_n\ltimes\{\pm 1\}^n)_S$,
$(\beta,\www{\varepsilon})\in\Br_n\ltimes\{\pm 1\}^n$
and $\uuuu{v}:=(\beta,\www{\varepsilon})(\uuuu{e})$. We have to show
$Z((\alpha,\varepsilon))(\uuuu{v})\in \BB^{dist}$. 
This is rather obvious with the commutativity of the actions of
$O^{(k)}$ and $\Br_n\ltimes\{\pm 1\}^n$ in Lemma \ref{t3.22} (a), 
\begin{eqnarray*}
Z((\alpha,\varepsilon))(\uuuu{v})
&=&Z((\alpha,\varepsilon))(\beta,\www{\varepsilon})(\uuuu{e})\\
&\stackrel{\ref{t3.22} (a)}{=}&
(\beta,\www{\varepsilon})Z((\alpha,\varepsilon))(\uuuu{e})\\
&=& (\beta,\www{\varepsilon})(\alpha,\varepsilon)(\uuuu{e}).
\end{eqnarray*}

(e) Elementary group theory gives the group isomorphisms
\begin{eqnarray*}
(\Br_n)_{S/\{\pm 1\}^n}&\cong& 
\frac{(\Br_n\ltimes\{\pm 1\}^n)_S}
{(\{\pm 1\}^n)_S},\\
(\Br_n)_{\uuuu{e}/\{\pm 1\}^n}&\cong& 
\frac{(\Br_n\ltimes\{\pm 1\}^n)_{\uuuu{e}}}
{(\{\pm 1\}^n)_{\uuuu{e}}}
\cong (\Br_n\ltimes\{\pm 1\}^n)_{\uuuu{e}}
\end{eqnarray*}
Here use $(\{\pm 1\}^n)_{\uuuu{e}}=\{(1,...,1)\}$. 
$\{\pm 1\}^n$ is a normal subgroup of 
$\Br_n\ltimes\{\pm 1\}^n$. 
Therefore $(\{\pm 1\}^n)_S$ is a normal subgroup of
$(\Br_n\ltimes\{\pm 1\}^n)_S$.
Therefore $Z((\{\pm 1\}^n)_S)$ is a normal subgroup of
$G_\Z^{\BB}$. Therefore $\oooo{Z}$ is well defined. 
Its kernel is still $(\Br_n)_{\uuuu{e}/\{\pm 1\}^n}
\cong (\Br_n\ltimes \{\pm 1\}^n)_{\uuuu{e}}$ 
because $(\Br_n\ltimes \{\pm 1\}^n)_{\uuuu{e}}\cap 
(\{\pm 1\}^n)_{S}=(\{\pm 1\}^n)_{\uuuu{e}}=\{(1,...,1)\}$.


(f) If $(H_\Z,L,\uuuu{e})$ is irreducible,
then the following graph is connected: its vertices 
are $e_1,...,e_n$, and it has 
an edge between $e_i$ and $e_j$ for $i<j$ if 
$S_{ij}(=L(e_j,e_i))\neq 0$. 
Therefore then 
$(\{\pm 1\}^n)_S=\{(1,...,1),(-1,...,-1)\}$. 
Everything else follows from this and from part (d). 
\hfill$\Box$ 

\bigskip
The antihomomorphism 
$Z:(\Br_n\ltimes \{\pm 1\}^n)_S\to G_\Z$ is not always surjective,
but in many cases. See Theorem \ref{t3.28}, the Remarks \ref{t3.29}
and the Remarks \ref{t9.1}. Theorem \ref{t3.26} (b) writes $(-1)^{k+1}M$
as an image of a braid by $Z$. Theorem \ref{t3.26} (c) gives conditions when
it has an $n$-th root which is also an image of a braid by $Z$.

\begin{theorem}\label{t3.26}
Let $(H_\Z,L,\uuuu{e})$ be a unimodular bilinear lattice of
rank $n\geq 2$ with a triangular basis $\uuuu{e}$
and matrix $S=L(\uuuu{e}^t,\uuuu{e})^t\in T^{uni}_n(\Z)$. 
Fix $k\in\{0,1\}$. Recall from chapter \ref{s3.1}
\begin{eqnarray*}
\sigma^{root}&:=&\sigma_{n-1}\sigma_{n-2}...\sigma_2\sigma_1
\in \Br_n,\\
\sigma^{mon}&:=&(\sigma^{root})^n,\\ 
\textup{center}(\Br_n)&=&\langle \sigma^{mon}\rangle.
\end{eqnarray*}

(a) 
\begin{eqnarray*}
\delta_n^{1-k}\sigma^{root}(\uuuu{e})
&=& Z(\delta^{1-k}_n\sigma^{root})(\uuuu{e}) \\
&=&(s_{e_1}^{(k)}(e_2),s_{e_1}^{(k)}(e_3),...,
s_{e_1}^{(k)}(e_n),s_{e_1}^{(k)}(e_1))\\
&=& \uuuu{e}\cdot R
\end{eqnarray*}
\begin{eqnarray*}
\textup{with }R&:=&\left(\begin{array}{ccc|c}
-q_{n-1}& \dots & -q_1  & -q_0 \\ \hline
 & & & \\
 & E_{n-1} & & \\
 & & & \end{array}\right)\\
\textup{so }R_{ij}&=&\left\{\begin{array}{ll}
-q_{n-j}& \textup{if }i=1, \\ \delta_{i-1,j}& \textup{if }i\geq 2,
\end{array}\right. 
\end{eqnarray*}
where $q_0=(-1)^k$, $q_{n+1-j}=S_{1j}$ for $j\in\{2,...,n\}$.

(b) 
\begin{eqnarray*}
Z((\delta^{1-k}_n\sigma^{root})^n)&=&(-1)^{k+1}M,\\
\textup{so especially }
(\delta^{1-k}_n\sigma^{root})^n&\in& (\Br_n\ltimes\{\pm 1\}^n)_S.
\end{eqnarray*}

(c) Write $q_0=(-1)^k$, $q_{n+1-j}=S_{1j}$ for $j\in\{2,...,n\}$ 
as in part (a) and additionally $q_n:=1$. Suppose
$q_{n-j}=q_0q_j$ for $j\in\{1,...,n-1\}$, and suppose that $S$
has the following shape, 
\begin{eqnarray*}
S&:=&\begin{pmatrix}1&q_{n-1}& \dots & q_1\\
 & \ddots & \ddots & \vdots \\
 &  & \ddots & q_{n-1}\\ 
 &  & & 1 \end{pmatrix}, \\
\textup{so }S_{ij}&=&\left\{\begin{array}{ll}
0& \textup{if }i>j, \\ q_{n-(j-i)}& \textup{if }i\leq j,
\end{array}\right.
\end{eqnarray*}
Then
\begin{eqnarray*}
M^{root}&:=&Z(\delta^{1-k}_n\sigma^{root})\in G_\Z\\
\textup{ with }(M^{root})^n&=&(-1)^{k+1}M.
\end{eqnarray*}
\index{root of the monodromy}\index{$M^{root}$}
$M^{root}$ is regular and cyclic and has the characteristic
polynomial $q(t):=\sum_{i=0}^n q_it^i\in\Z[t]$. 
\end{theorem}

{\bf Proof:}
(a) The second line follows from the definition of the action 
of $\delta^{1-k}_n\sigma^{root}$ on $\uuuu{e}$. For the third 
line observe $s_{e_1}^{(k)}(e_j)=e_j-S_{1j}e_1$ for 
$j\geq 2$ and $s_{e_1}^{(k)}(e_1)=-q_0e_1$. 

(b) Use part (a) and 
\begin{eqnarray*}
s_{s_{e_1}^{(k)}(e_2)}^{(k)}(s_{e_1}^{(k)}(e_j))
=s_{e_1}^{(k)}s_{e_2}^{(k)}(s_{e_1}^{(k)})^{-1}s_{e_1}^{(k)}(e_j)
=s_{e_1}^{(k)}s_{e_2}^{(k)}(e_j)
\end{eqnarray*}
to find 
\begin{eqnarray*}
(\delta^{1-k}_n\sigma^{root})^2(\uuuu{e})
=(s_{e_1}^{(k)}s_{e_2}^{(k)}(e_3),...,
s_{e_1}^{(k)}s_{e_2}^{(k)}(e_n),
s_{e_1}^{(k)}s_{e_2}^{(k)}(e_1),
s_{e_1}^{(k)}s_{e_2}^{(k)}(e_2)).
\end{eqnarray*}
One continues inductively and finds
\begin{eqnarray*}
(\delta^{1-k}_n\sigma^{root})^n(\uuuu{e})
=(s_{e_1}^{(k)}...s_{e_n}^{(k)}(e_1),...,
s_{e_1}^{(k)}...s_{e_n}^{(k)}(e_n))
=(-1)^{k+1}M(\uuuu{e}),
\end{eqnarray*}
so $Z((\delta^{1-k}_n\sigma^{root})^n)=(-1)^{k+1}M$. 

(c) If $S$ is as in part (c) then
\begin{eqnarray*}
I^{(k)}(\delta^{1-k}_n\sigma^{root}(\uuuu{e})^t,
\delta^{1-k}_n\sigma^{root}(\uuuu{e}))
&\stackrel{(1)}{=}& 
I^{(k)}((e_2\ e_3\ ...\ e_n\ e_1)^t,(e_2\ e_3\ ...\ e_n\ e_1))\\
&\stackrel{(2)}{=}&S+(-1)^kS^t 
=I^{(k)}(\uuuu{e}^t,\uuuu{e}).
\end{eqnarray*}
Here $\stackrel{(1)}{=}$ uses $s_{e_1}^{(k)}\in O^{(k)}$,
and $\stackrel{(2)}{=}$ uses that 
$I^{(k)}(\uuuu{e}^t,\uuuu{e})=S+(-1)^kS^t$ 
and that $S$ is as in part (c).

Therefore $M^{root}:=Z(\delta^{1-k}_n\sigma^{root})
\in G_\Z^{(k)}$, so by Lemma \ref{t3.25} (a)
\begin{eqnarray*}
M^{root}\in G_\Z\quad\textup{and}\quad\delta^{1-k}_n\sigma^{root}
\in (\Br_n\ltimes\{\pm 1\}^n)_S.
\end{eqnarray*}
Also
\begin{eqnarray*}
(M^{root})^n=(Z(\delta^{1-k}_n\sigma^{root}))^n
=Z((\delta^{1-k}_n\sigma^{root})^n)
=(-1)^{k+1}M.
\end{eqnarray*}

Let $\uuuu{e}^*$ be the $\Z$-basis
of $H_\Z$ which is left $L$-dual to the $\Z$-basis $\uuuu{e}$,
so with $L((\uuuu{e}^*)^t,\uuuu{e})=E_n$. Remark 4.8 in 
\cite{BH20} says
\begin{eqnarray*}
M^{root}\uuuu{e}^* = \uuuu{e}^*R^{-t}=\uuuu{e}^*\cdot
\left(\begin{array}{ccc|c}
 & & & -q_0 \\ \hline & & & -q_1 \\ & E_{n-1} & & \vdots \\
 & & & -q_{n-1}\end{array}\right).
\end{eqnarray*}
The matrix $R^{-t}$ is the companion matrix of the polynomial
$q(t)$. Therefore $M^{root}$ is regular, cyclic with generating
vector $c=e_1^*$ and has the characteristic polynomial $q(t)$.
\hfill$\Box$

\begin{remarks}\label{t3.27}
The main part of part (c) of Theorem \ref{t3.26} 
has also the following matrix version: 
For $q(t)=\sum_{i=0}^nq_it^i\in\Z[t]$ with 
$q_n=1$, $q_0=(-1)^k$ for some $k\in\{0;1\}$ and
$q_{n-j}=q_0q_j$ the matrix $R$ in part (a) and
the matrix $S$ in part (c) of Theorem \ref{t3.26} 
satisfy $$R^n=(-1)^{k+1}S^{-1}S^t.$$
A proof using matrices of this version of part (c)
of Theorem \ref{t3.26} was given in 
\cite[Theorem 4.5 (a)+(b)]{BH20}.
The proof here with the braid group action is more elegant.
\end{remarks}

The antihomomorphism $Z:(\Br_n\ltimes\{\pm 1\}^n)_S\to G_\Z$
is not surjective in general. A simple example with $n=4$
is given in the Remarks \ref{t3.29}. 
But it is surjective in the case $n=1$, in all cases with $n=2$
and in almost all cases with $n=3$. 
Theorem \ref{t3.28} gives precise statements.
Its proof requires first a good control of the braid group
action on $T^{uni}_3(\Z)$, which is the subject of chapter \ref{s4}
and second complete knowledge of the group $G_\Z$ for all cases
with $n\leq 3$, which is the subject of chapter \ref{s5}.
Theorem \ref{t3.28} is proved within the theorems in chapter
\ref{s5} which treat the different cases with $n\in\{1,2,3\}$,
namely Lemma \ref{t5.4} (the cases $A_1^n$), 
Theorem \ref{t5.5} (the rank 2 cases),
Theorem \ref{t5.13} (the reducible rank 3 cases),
Theorem \ref{t5.14} (the irreducible rank 3 cases with all eigenvalues
in $S^1$),
Theorem \ref{t5.16} (some special other rank 3 cases),
Theorem \ref{t5.18} (the rest of the rank 3 cases).

\begin{theorem}\label{t3.28}
Let $(H_\Z,L,\uuuu{e})$ be a unimodular bilinear lattice of 
rank $n\leq 3$ with triangular basis $\uuuu{e}$ and
matrix $S=L(\uuuu{e}^t,\uuuu{e})^t\in T^{uni}_n(\Z)$.
The group antihomomorphism $Z:(\Br_n\ltimes\{\pm 1\}^n)_S
\to G_\Z$ is not surjective in the four cases with $n=3$ 
where $S$ is in the $\Br_3\ltimes\{\pm 1\}^3$ orbit of $S(\uuuu{x})$
with $$\uuuu{x}\in\{(3,3,4),(4,4,4),(5,5,5),(4,4,8)\},$$
so then $G_\Z\supsetneqq G_\Z^{\BB}$.
It is surjective in all other cases with $n\leq 3$,
so then $G_\Z=G_\Z^{\BB}$. 
\end{theorem}

\begin{remarks}\label{t3.29}
(i) Contrary to Theorem \ref{t3.28} for the cases 
$n=1,2,3$, it is in the cases $n\geq 4$ easy to find matrices
$S\in T^{uni}_n(\Z)$ such that the group antihomomorphism
$Z:(\Br_n\ltimes \{\pm 1\}^n)_S\to G_\Z$ is not surjective.
Though the construction which we propose in part (ii) and
carry out in one example in part (iii) leads to matrices
which are rather particular. 
For a given matrix $S$ it is in general not easy to see whether
$Z$ is surjective or not. 

(ii) Consider a reducible unimodular bilinear lattice
$(H_\Z,L,\uuuu{e})$ of rank $n$ with triangular basis
$\uuuu{e}$ and matrix $S=L(\uuuu{e}^t,\uuuu{e})^t\in 
T^{uni}_n(\Z)$. There are an $L$-orthogonal decomposition
$H_\Z=\bigoplus_{j=1}^l H_{\Z,j}$ with $l\geq 2$ and 
$\textup{rank}\, H_{\Z,j}\geq 1$ and a surjective map 
$\alpha:\{1,...,n\}\to\{1,...,l\}$ with 
$e_i\in H_{\Z,\alpha(i)}$. Then for $k\in\{0,1\}$
\begin{eqnarray*}
\Gamma^{(k)}\{e_i\}\subset H_{\Z,\alpha(i)},\\
\textup{so}\quad \Delta^{(k)}\subset \bigcup_{j=1}^l H_{\Z,j}.
\end{eqnarray*}
Especially any $g\in G_\Z^{\BB}\subset G_\Z$ maps each $e_i$ to
an element of $\bigcup_{j=1}^l H_{\Z,j}$. This does not 
necessarily hold for any $g\in G_\Z$. Part (iii) gives an 
example.

(iii) Consider the unimodular bilinear lattice 
$(H_\Z,L,\uuuu{e})$ of rank 4 with triangular basis
$\uuuu{e}$ and matrix 
$$S=L(\uuuu{e}^t,\uuuu{e})^t=
\begin{pmatrix}1&2&0&0\\0&1&0&0\\0&0&1&2\\0&0&0&1\end{pmatrix}
\in T^{uni,4}(\Z).$$
Then $H_\Z= H_{\Z,1}\oplus H_{\Z,2}$ with 
$H_{\Z,1}=\Z e_1\oplus \Z e_2$ and
$H_{\Z,2}=\Z e_3\oplus \Z e_4$. The $\Z$-linear map
$g:H_\Z\to H_\Z$  with 
\begin{eqnarray*}
(g(e_1),g(e_2),g(e_3),g(e_4))
=(e_1+(e_3-e_4),e_2+(e_3-e_4),\\
e_3+(e_1-e_2),e_4+(e_1-e_2))
\end{eqnarray*}
is not in $G_\Z^{\BB}$ because 
$g(e_1)$, $g(e_2)$, $g(e_3)$, $g(e_4)\notin H_{\Z,1}\cup
H_{\Z,2}$. But $g\in G_\Z$ because
\begin{eqnarray*}
L(e_1-e_2,e_1-e_2)=0=L(e_3-e_4,e_3-e_4),\\
\textup{so}\quad L(g(e_i),g(e_j))=L(e_i,e_j)\quad\textup{for }
\{i,j\}\subset\{1,2\}\textup{ or }\{i,j\}\subset\{3,4\},
\end{eqnarray*}
and also
\begin{eqnarray*}
L(g(e_i),g(e_j))=L(e_i,e_j)\quad\textup{for }
(i,j)\in(\{1,2\}\times\{3,4\})\cup
(\{3,4\}\times\{1,2\}).
\end{eqnarray*}
So here $g\in G_\Z-G_\Z^{\BB}$, so $G_\Z^{\BB}\subsetneqq G_\Z$. 
\end{remarks}

\chapter{Braid group action on upper triangular
$3\times 3$ matrices }\label{s4}
\setcounter{equation}{0}
\setcounter{figure}{0}

The subject of this chapter is the case $n=3$ of the
action in Lemma \ref{t3.13} of $\Br_n\ltimes\{\pm 1\}^n$
on the matrices in $T^{uni}_n(\Z)$.

In section \ref{s4.1} the action on $T^{uni}_3(\R)$ is 
made concrete. The (quotient) group of actions is given
in new generators. It is 
$$ (G^{phi}\ltimes 
G^{sign})\rtimes \langle\gamma\rangle\cong 
(G^{phi}\rtimes \langle\gamma\rangle)\ltimes G^{sign},$$
where $G^{phi}$
is a free Coxeter group with three generators,
$G^{sign}$ is the group of actions in $T^{uni}_3(\R)$ which
the sign group $\{\pm 1\}^3$ induces, 
and $\gamma$ acts cyclically of order 3.
In fact, $G^{phi}\rtimes \langle\gamma\rangle
\cong PSL_2(\Z),$
so we have a nonlinear action of $PSL_2(\Z)$, but this
way to look at it is less useful than the presentation 
as $G^{phi}\rtimes\langle\gamma\rangle$.

The action on $T^{uni}_3(\Z)$ had been studied already 
by Kr\"uger \cite[\S 12]{Kr90} and by Cecotti-Vafa
\cite[Ch. 6.2]{CV93}. Section \ref{4.2} recovers and
refines their results.  Like them, it puts emphasis
on the cases where the monodromy of a corresponding
unimodular bilinear lattice has eigenvalues in $S^1$.
Section \ref{s4.2} follows largely Kr\"uger 
\cite[\S 12]{Kr90}.

Section \ref{s4.3} uses {\it pseudo-graphs} to 
systematically study all cases, not only those
where the monodromy has eigenvalues in $S^1$.
This goes far beyond Kr\"uger and Cecotti-Vafa. 

The results of section \ref{s4.3} and the pseudo-graphs
are used in section \ref{s4.4} to determine in all cases
the stabilizer $(\Br_3\ltimes\{\pm 1\}^3)_S$
respectively the stabilizer $(\Br_3)_{S/\{\pm 1\}^3}$. 
Section \ref{s7.4} will build on this and determine 
the stabilizer $(\Br_3)_{\uuuu{e}/\{\pm 1\}^3}$
of a distingushed basis $\uuuu{e}\in\BB^{dist}$
for any unimodular bilinear lattice of rank 3
with a fixed triangular basis. 

Section \ref{s4.5} starts with the observation
that a matrix $S\in T^{uni}_n(\Z)$ and the
matrix $\www{S}\in T^{uni}_n(\Z)$ with 
$\www{S}_{ij}=-S_{ij}$ for $i<j$ lead to
unimodular bilinear lattices with the same 
odd monodromy groups and the same odd vanishing cycles.
This motivates to study the action on $T^{uni}_n(\Z)$
which extends the action of $\Br_n\ltimes\{\pm 1\}^n$
by this global sign change.
Section \ref{s4.5} carries this out in the case $n=3$
and gives standard representatives for each orbit.
Though examples show that the action is rather wild.
Similar looking triples in $\Z^3$ are in the orbits
of very different standard representatives.

\section[Real upper triangular $3\times 3$ matrices]
{Braid group action on real upper triangular
$3\times 3$ matrices}
\label{s4.1}

The action of $\Br_3\ltimes \{\pm 1\}^3$ on $T^{uni}_3(\Z)$ 
will be studied in the next sections. It extends to an action on 
$T^{uni}_3(\R)\cong\R^3$ which will be studied here.
By Theorem \ref{t3.4} (d), $\sigma_1$ acts on $T^{uni}_3(\Z)$ by
\begin{eqnarray*}
\sigma_1:\begin{pmatrix}1&x_1&x_2\\0&1&x_3\\0&0&1\end{pmatrix}
&\mapsto& \begin{pmatrix}-x_1&1&0\\1&0&0\\0&0&1\end{pmatrix}
\begin{pmatrix}1&x_1&x_2\\0&1&x_3\\0&0&1\end{pmatrix}
\begin{pmatrix}-x_1&1&0\\1&0&0\\0&0&1\end{pmatrix}\\
&=&\begin{pmatrix}1&-x_1&x_3-x_1x_2\\0&1&x_2\\0&0&1\end{pmatrix}.
\end{eqnarray*}
It extends to an action on $T^{uni}_3(\R)$.
With the isomorphism
\begin{eqnarray*}
T^{uni}_3(R)\stackrel{\cong}{\longrightarrow}R^3,\quad
\begin{pmatrix}1&x_1&x_2\\0&1&x_3\\0&0&1\end{pmatrix}
\mapsto (x_1,x_2,x_3)
\quad\textup{for }R\in\{\Z,\Q,\R,\C\}
\end{eqnarray*}
this gives the action
\begin{eqnarray*}
\sigma_1^{\R}:\R^3\to\R^3,\quad
(x_1,x_2,x_3)\mapsto (-x_1,x_3-x_1x_2,x_2).
\end{eqnarray*}
Analogously \index{$\sigma_j^\R:\R^3\to\R^3$}
\index{$\delta_j^\R:\R^3\to\R^3$}
\begin{eqnarray*}
(\sigma_1^{\R})^{-1}:\R^3\to\R^3,&& 
(x_1,x_2,x_3)\mapsto (-x_1,x_3,x_2-x_1x_3), \\
\sigma_2^{\R}:\R^3\to\R^3,&& 
(x_1,x_2,x_3)\mapsto (x_2-x_1x_3,x_1,-x_3), \\
(\sigma_2^{\R})^{-1}:\R^3\to\R^3,&&  
(x_1,x_2,x_3)\mapsto (x_2,x_1-x_2x_3,-x_3), \\
\delta_1^{\R}:\R^3\to\R^3,&& 
(x_1,x_2,x_3)\mapsto (-x_1,-x_2,x_3), \\
\delta_2^{\R}:\R^3\to\R^3,&& 
(x_1,x_2,x_3)\mapsto (-x_1,x_2,-x_3), \\
\delta_3^{\R}:\R^3\to\R^3,&& 
(x_1,x_2,x_3)\mapsto (x_1,-x_2,-x_3).
\end{eqnarray*}
One sees \index{$G^{sign}$}
\begin{eqnarray*}
\delta_3^{\R}=\delta_1^{\R}\delta_2^{\R}\quad
\textup{and}\quad 
G^{sign}:=\langle\delta_1^{\R},\delta_2^{\R}\rangle
\cong\{\pm 1\}^2.
\end{eqnarray*}
The group $\langle \sigma_1^{\R},\sigma_2^{\R}\rangle
\ltimes G^{sign}\subset\Aut_{pol}(\R^3)$ of polynomial
automorphisms of $\R^3$ will be more transparent in other 
generators.

\begin{definition}\label{t4.1}
Define the polynomial automorphisms of $\R^3$
\index{$\varphi_j:\R^3\to\R^3$}\index{$\gamma:\R^3\to\R^3$} 
\begin{eqnarray*}
\varphi_1:\R^3\to\R^3,&& 
(x_1,x_2,x_3)\mapsto (x_2x_3-x_1,x_3,x_2), \\
\varphi_2:\R^3\to\R^3,&&  
(x_1,x_2,x_3)\mapsto (x_3,x_1x_3-x_2,x_1), \\
\varphi_3:\R^3\to\R^3,&&  
(x_1,x_2,x_3)\mapsto (x_2,x_1,x_1x_2-x_3), \\
\gamma:\R^3\to\R^3,&& 
(x_1,x_2,x_3)\mapsto (x_3,x_1,x_2),
\end{eqnarray*}
and the group $G^{phi}:=\langle \varphi_1,\varphi_2,\varphi_3
\rangle\subset\Aut_{pol}(\R^3)$. \index{$G^{phi}$} 
\end{definition}

\begin{theorem}\label{t4.2}
(a) The group $G^{phi}$ is a free Coxeter group
with the three generators $\varphi_1,\varphi_2,\varphi_3$,
so $G^{phi}\cong G^{fCox,3}$. 

(b) $\langle\gamma\rangle\cong \Z/3\Z\cong A_3\subset S_3$.

(c) $\langle \sigma_1^{\R},\sigma_2^{\R}\rangle
\ltimes G^{sign} = 
(G^{phi}\ltimes G^{sign})\rtimes \langle\gamma\rangle.$
\end{theorem}

{\bf Proof:}
(a) $\varphi_1^2=\varphi_2^2=\varphi_3^2=\id$ is obvious,
and also that $(2,2,2)\in\R^3$ is a fixed point of $G^{phi}$.
We will show that the group 
$\langle d_{(2,2,2)}\varphi_1,d_{(2,2,2)}\varphi_2,
d_{(2,2,2)}\varphi_3\rangle$ of induced actions on the 
tangent space $T_{(2,2,2)}\R^3$ is a free Coxeter group 
with three generators. This will imply
$G^{phi}\cong G^{fCox,3}$.

Affine linear coordinates $(\www{x}_1,\www{x}_2,\www{x}_3)$ on
$\R^3$ which vanish at $(2,2,2)$ with
\begin{eqnarray*}
(x_1,x_2,x_3)=(2+\www{x}_1,2+\www{x}_2,2+\www{x}_3)
=(2,2,2)+(\www{x}_1,\www{x}_2,\www{x}_3)
\end{eqnarray*}
are also linear coordinates on $T_{(2,2,2)}\R^3$. We have 
\begin{eqnarray*}
\varphi_1((2,2,2)+(\www{x}_1,\www{x}_2,\www{x}_3))
=(2,2,2)+(\www{x}_2\www{x}_3+2\www{x}_2+2\www{x}_3-\www{x}_1,
\www{x}_3,\www{x}_2),\\
d_{(2,2,2)}\varphi_1(\www{x}_1,\www{x}_2,\www{x}_3)
=(\www{x}_1,\www{x}_2,\www{x}_3)
\begin{pmatrix}-1 &0&0&\\2&0&1\\2&1&0\end{pmatrix},
\end{eqnarray*}
and analogously
\begin{eqnarray*}
d_{(2,2,2)}\varphi_2(\www{x}_1,\www{x}_2,\www{x}_3)
=(\www{x}_1,\www{x}_2,\www{x}_3)
\begin{pmatrix}0&2&1\\0&-1&0\\1&2&0\end{pmatrix},\\
d_{(2,2,2)}\varphi_3(\www{x}_1,\www{x}_2,\www{x}_3)
=(\www{x}_1,\www{x}_2,\www{x}_3)
\begin{pmatrix}0&1&2\\1&0&2\\0&0&-1\end{pmatrix}.
\end{eqnarray*}
The group $G^{phi}$ respects the fibers of the map
\begin{eqnarray*}
r_\R:\R^3\to \R,\quad (x_1,x_2,x_3)\mapsto 
x_1^2+x_2^2+x_3^2-x_1x_2x_3.
\end{eqnarray*}
The group $\langle d_{(2,2,2)}\varphi_1,d_{(2,2,2)}\varphi_2,
d_{(2,2,2)}\varphi_3\rangle$ respects the tangent cone
at $(2,2,2)$ of the fiber $r_\R^{-1}(4)$. 
This tangent cone is the zero set of the quadratic form 
\begin{eqnarray*}
q_{\R^3,(2,2,2)}:\R^3&\to& \R,\\
(\www{x}_1,\www{x}_2,\www{x}_3)&\mapsto&
-\www{x}_1^2-\www{x}_2^2-\www{x}_3^2+2\www{x}_1\www{x}_2
+2\www{x}_1\www{x}_3+2\www{x}_2\www{x}_3.
\end{eqnarray*}
This quadratic form is indefinite with signature $(+,-,-)$.
As in Theorem \ref{ta.4}, its cone of positive vectors is 
called $\KK$. Consider the six vectors
\begin{eqnarray*}
v_1=(1,1,0),\ v_2=(1,0,1),\ v_3=(0,1,1),\\
w_1=v_1+v_2,\ w_2=v_1+v_3,\ w_3=v_2+v_3.
\end{eqnarray*}
Then
\begin{eqnarray*}
v_1,v_2,v_3\in\paa\KK,\quad w_1,w_2,w_3\in\KK,\\
d_{(2,2,2)}\varphi_1:v_1\leftrightarrow v_2,\ w_1\mapsto w_1,\\
d_{(2,2,2)}\varphi_2:v_1\leftrightarrow v_3,\ w_2\mapsto w_2,\\
d_{(2,2,2)}\varphi_3:v_2\leftrightarrow v_3,\ w_3\mapsto w_3.
\end{eqnarray*}
Compare Theorem \ref{ta.4}. In the model $\KK/\R^*$ of the
hyperbolic plane, $d_{(2,2,2)}\varphi_i$ for $i\in\{1,2,3\}$
gives a rotation with angle $\pi$ and elliptic fixed point
$\R^* w_i$, which maps the hyperbolic line with euclidean 
boundary points $\R^*v_1\& \R^*v_2$ respectively
$\R^*v_1\& \R^*v_3$ respectively $\R^*v_2\& \R^*v_3$ to 
itself 

\begin{figure}[H]
\includegraphics[width=0.5\textwidth]{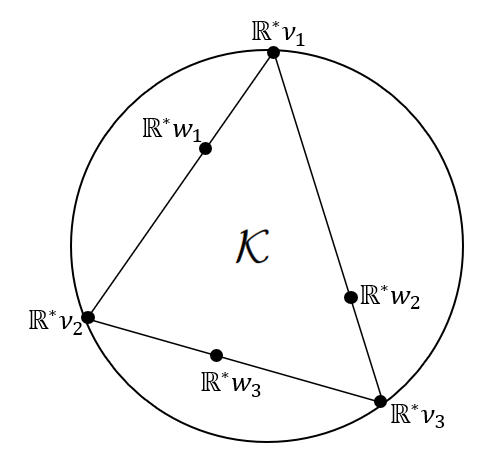}
\caption[Figure 4.1]{$G^{fCox,3}$ generated by 3 elliptic
M\"obius transformations, an application of 
Theorem \ref{ta.2} (b)}
\label{Fig:4.1}
\end{figure}

Theorem \ref{ta.2} (b) applies and shows
$\langle d_{(2,2,2)}\varphi_i\,|\, i\in\{1,2,3\}\rangle
\cong G^{fCox,3}$. Figure \ref{Fig:4.1} illustrates this.

(b) Trivial.

(c) The equality of groups
\begin{eqnarray*}
\langle \sigma_1^{\R},\sigma_2^{\R}\rangle\ltimes G^{sign}
=\langle\varphi_1,\varphi_2,\varphi_3,\gamma,
\delta_1^{\R},\delta_2^{\R}\rangle
\end{eqnarray*}
follows from 
\begin{eqnarray}\label{4.1}
\gamma &=& \delta_3^{\R}\sigma_2^{\R}\sigma_1^{\R},\\
\varphi_1&=& \delta_1^{\R}\gamma^{-1}(\sigma_2^{\R})^{-1},
\label{4.2}\\
\varphi_2&=& \delta_1^{\R}\gamma \sigma_2^{\R}
=\delta_3^{\R}\gamma^{-1}(\sigma_1^{\R})^{-1},
\label{4.3}\\
\varphi_3&=& \delta_3^{\R}\gamma\sigma_1^{\R}.\label{4.4}
\end{eqnarray}
and
\begin{eqnarray}\label{4.5}
\sigma_1^{\R}=\gamma^{-1}\delta_3^{\R}\varphi_3,\quad
\sigma_2^{\R}=\gamma^{-1}\delta_1^{\R}\varphi_2.
\end{eqnarray}
$G^{phi}$ fixes $(2,2,2)$, and therefore
$G^{phi}\cap G^{sign}=\{\id\}$. As $G^{sign}$ is a normal
subgroup of $\langle \sigma_1^{\R},\sigma_2^{\R}\rangle
\ltimes G^{sign}$, it is also a normal subgroup of
$\langle \varphi_1,\varphi_2,\varphi_3,\delta_1^{\R},
\delta_2^{\R}\rangle$, so 
$\langle\varphi_1,\varphi_2,\varphi_3,\delta_1^{\R},
\delta_2^{\R}\rangle =G^{phi}\ltimes G^{sign}$. 
More precisely
\begin{eqnarray*}
\varphi_i\delta_j^{\R}\varphi_i^{-1}=\delta_k^{\R}
\quad\textup{for }(i,j,k)\in\{(1,3,3),(2,2,2),(3,1,1), \\
(1,1,2),(1,2,1),(2,1,3),(2,3,1),(3,2,3),(3,3,2)\}.\end{eqnarray*}
We claim $\gamma\notin G^{phi}\ltimes G^{sign}$.
If $\gamma$ were in $G^{phi}\ltimes G^{sign}$ then
$\gamma\in G^{phi}$ as $\gamma$ fixes $(2,2,2)$. But all
elements of finite order in $G^{phi}\cong G^{fCox,3}$
have order two, though $\gamma$ has order three.
Hence $\gamma\notin G^{phi}$ and $\gamma\notin 
G^{phi}\ltimes G^{sign}$.

The claim and 
\begin{eqnarray}\label{4.6}
\gamma\varphi_1\gamma^{-1}=\varphi_2,\ 
\gamma\varphi_2\gamma^{-1}=\varphi_3,\
\gamma\varphi_3\gamma^{-1}=\varphi_1,\\ 
\gamma\delta_1^{\R}\gamma^{-1}=\delta_3^{\R},\ 
\gamma\delta_2^{\R}\gamma^{-1}=\delta_1^{\R},\ 
\gamma\delta_3^{\R}\gamma^{-1}=\delta_2^{\R},\label{4.7}
\end{eqnarray}
show
\begin{eqnarray*}
\langle \varphi_1,\varphi_2,\varphi_3,\delta_1^{\R},
\delta_2^{\R},\gamma\rangle
=(G^{phi}\ltimes G^{sign})\rtimes \langle\gamma\rangle.
\hspace*{2cm}\Box
\end{eqnarray*}

\section[Integer upper triangular $3\times 3$ matrices]
{Braid group action on integer upper triangular
$3\times 3$ matrices}
\label{s4.2}

In this section we will give a partial classification
of the orbits of the action of $\Br_3\ltimes\{\pm 1\}^3$
on $T^{uni}_3(\Z)$. This refines results which were obtained
independently by Kr\"uger \cite[\S 12]{Kr90} and Cecotti-Vafa
\cite[Ch. 6.2]{CV93} (building on Mordell \cite[p 106 ff]{Mo69}).

The refinement consists in the following.
By Theorem \ref{t4.2} the action of $\Br_3\ltimes\{\pm 1\}^3$
on $T^{uni}_3(\Z)$ coincides with the action of 
$(G^{phi}\ltimes G^{sign})\rtimes\langle\gamma\rangle$.
For reasons unknown to us, Kr\"uger and Cecotti-Vafa
considered the action of the slightly larger group
$(G^{phi}\ltimes G^{sign})\rtimes\langle \gamma,\gamma_2\rangle$
with 
\begin{eqnarray*}
\gamma_2:\R^3\to\R^3,\quad (x_1,x_2,x_3)\mapsto (x_2,x_1,x_3),
\quad\textup{so }\langle \gamma,\gamma_2\rangle\cong S_3.
\end{eqnarray*}
Thus they obtained a slightly coarser classification.
Nevertheless Theorem \ref{t4.6} is essentially due to them
(and Mordell \cite[p 106 ff]{Mo69}). 
The following definition and lemma prepare it.
They are due to Kr\"uger \cite[\S 12]{Kr90}.

\begin{definition}\label{t4.3} \cite[Def. 12.2]{Kr90}
For $\uuuu{x}=(x_1,x_2,x_3)\in\R^3$ we set as usual
$\|\uuuu{x}\|:=\sqrt{x_1^2+x_2^2+x_3^2}$. A
tuple $\uuuu{x}\in\R^3$ is called a
{\it local minimum} if \index{local minimum}
\begin{eqnarray*}
\| \uuuu{x}\| \leq\min(\|\sigma_1^{\R}(\uuuu{x})\|,
\|(\sigma_1^{\R})^{-1}(\uuuu{x})\|, 
\|\sigma_2^{\R}(\uuuu{x})\|,
\|(\sigma_2^{\R})^{-1}(\uuuu{x})\|).
\end{eqnarray*}
This is obviously equivalent to 
\begin{eqnarray*}
\| \uuuu{x}\| \leq\min(\|\varphi_1(\uuuu{x})\|,
\|\varphi_2(\uuuu{x})\|, 
\|\varphi_3(\uuuu{x})\|).
\end{eqnarray*}
\end{definition}

\begin{lemma}\label{t4.4} \cite[Lemma 12.3]{Kr90}
$\uuuu{x}\in\R^3$ is a local minimum if and only if it satisfies
(i) or (ii),
\begin{list}{}{}
\item[(i)] 
$x_1x_2x_3\leq 0$,
\item[(ii)]
$x_1x_2x_3>0,\ 2|x_1|\leq |x_2x_3|,\ 2|x_2|\leq |x_1x_3|,\ 
2|x_3|\leq |x_1x_2|$.
\end{list}
In the case (ii) also $|x_1|\geq 2$, $|x_2|\geq 2$ and 
$|x_3|\geq 2$ hold.
\end{lemma}

{\bf Proof:}
$\uuuu{x}\in\R^3$ is a local minimum if for all $i,j,k$ with
$\{i,j,k\}=\{1,2,3\}$ 
\begin{eqnarray*}
x_i^2+x_j^2+x_k^2\leq x_i^2+x_j^2+(x_k-x_ix_j)^2
\end{eqnarray*}
holds, which is equivalent to 
\begin{eqnarray*}
2x_ix_jx_k\leq x_i^2x_j^2.
\end{eqnarray*}

{\bf 1st case,} $x_1x_2x_3\leq 0$: Then $\uuuu{x}$ is a local
minimum.

{\bf 2nd case,} $x_1x_2x_3>0$: Then the condition
$2x_ix_jx_k\leq x_i^2x_j^2$ is equivalent to 
\begin{eqnarray*}
2|x_k|\leq |x_ix_j|.
\end{eqnarray*}
These three conditions together imply
\begin{eqnarray*}
4|x_k|\leq 2|x_i||x_j|\leq |x_i||x_i||x_k|,\quad\textup{so }
4\leq |x_i|^2,\quad\textup{so }2\leq |x_i|.
\hspace*{1cm}\Box
\end{eqnarray*}

\bigskip
The square $\|.\|^2:\Z^3\to\Z_{\geq 0}$ of the norm 
has on $\Z^3$ values in $\Z_{\geq 0}$. Therefore
each $\Br_3\ltimes\{\pm 1\}^3$ orbit in $\Z^3$ has local minima.
Kr\"uger showed that the only $\Br_3\ltimes\{\pm 1\}^3$ orbits
in $\R^3$ without local minimal are of the following shape.
We will not use this result, but we find it interesting enough
to cite it.

\begin{theorem}\label{t4.5} \cite[Theorem 12.6]{Kr90}
Let $\uuuu{x}\in\R^3$ whose $\Br_3\ltimes\{\pm 1\}^3$ orbit
does not contain a local minimum. Then
\begin{eqnarray*}
x_1x_2x_3>0,\quad  
2<\min(|x_1|,|x_2|,|x_3|),\\
4=r_\R(\uuuu{x})(=x_1^2+x_2^2+x_3^2-x_1x_2x_3).
\end{eqnarray*}
Furthermore, there is a sequence $(\psi_n)_{n\in\N}$ with
$\psi_n\in\{\varphi_1,\varphi_2,\varphi_3\}$ with 
$\psi_n\neq\psi_{n+1}$ such that the sequence
$(\uuuu{x}^{(n)})_{n\in\N\cup\{0\}}$ with 
$\uuuu{x}^{(0)}=\uuuu{x}$ and 
$\uuuu{x}^{(n+1)}=\psi_n(\uuuu{x}^{(n)})$ satisfies
\begin{eqnarray*}
\| \uuuu{x}^{(n+1)}\| < \| \uuuu{x}^{(n)}\| \quad
\textup{for all }n\in\N\cup\{0\},\\
\lim_{n\to\infty} (|x^{(n)}_1|,|x^{(n)}_2|,|x^{(n)}_3|)=(2,2,2).
\end{eqnarray*}
\end{theorem}

Now we come to the classification of $\Br_3\ltimes\{\pm 1\}^3$
orbits in $\Z^3$. The following result is except for its part (f) 
a refinement of
\cite[Theorem 12.7]{Kr90} and of \cite[Ch 6.2]{CV93}. 
The proof below follows (except for the part (f)) 
the proof in \cite{Kr90}. 
Recall that each $\Br_3\ltimes\{\pm 1\}^3$ orbit in $\R^3$
is contained in one fiber of the map
$r_\R:\R^3\to\R$, $\uuuu{x}\mapsto x_1^2+x_2^2+x_3^2-x_1x_2x_3$.

\begin{theorem}\label{t4.6}
(a) Each fiber of \index{$r,\ r_\R$}
$$r:\Z^3\to\Z,\ 
r(\uuuu{x})=x_1^2+x_2^2+x_3^2-x_1x_2x_3$$ 
except the fiber $r^{-1}(4)$ 
contains only finitely many local minima.

(b) Each fiber of $r:\Z^3\to\Z$ except the fiber $r^{-1}(4)$
consists of only finitely many $\Br_3\ltimes\{\pm 1\}^3$ orbits.

(c) For $\rho\in\Z_{<0}$, each local minimum 
$\uuuu{x}\in r^{-1}(\rho)$ satisfies $x_1x_2x_3>0$
and $|x_1|\geq 3$, $|x_2|\geq 3$, $|x_3|\geq 3$.

(d) For $\rho\in\N-\{4\}$, each local minimum 
$\uuuu{x}\in r^{-1}(\rho)$ satisfies $x_1x_2x_3\leq 0$.

(e) The following table gives all local minima in
$r^{-1}(\{0,1,2,3,4\})$. The local minima in one 
$\Br_3\ltimes\{\pm 1\}^3$ orbit are in one line.
The last entry in each line is one matrix in the corresponding
orbit in $T^{uni}_3(\Z)$.
\begin{eqnarray*}
\begin{array}{rll}
r=3 & - & - \\
r=0 & (0,0,0) & S(A_1^3) \\
r=0 & (3,3,3),(-3,-3,3),(-3,3,-3),(3,-3,-3) & S(\P^2) \\
r=1 & (\pm 1,0,0),(0,\pm 1,0),(0,0,\pm 1) & S(A_2A_1) \\
r=2 & (\pm 1,\pm 1,0),(\pm 1,0,\pm 1),(0,\pm 1,\pm 1) & S(A_3)\\
r=4 & (\pm 2,0,0),(0,\pm 2,0),(0,0,\pm 2) & S(\P^1A_1) \\
r=4 & (-1,-1,-1),(1,1,-1),(1,-1,1),(-1,1,1) & S(\whh{A}_2)\\
r=4 & (2,2,2),(-2,-2,2),(-2,2,-2),(2,-2,-2) & S(\HH_{1,2})\\
\left\{\begin{array}{c}r=4 \\ l\in\Z_{\geq 3}\end{array}\right\}
& \left\{\begin{array}{l} 
(\varepsilon_1 2,\varepsilon_2 l,\varepsilon_3 l),
(\varepsilon_1 l,\varepsilon_2 2,\varepsilon_3 l),
(\varepsilon_1 l,\varepsilon_2 l,\varepsilon_3 2)\\
\textup{for }\varepsilon_1,\varepsilon_2,\varepsilon_3\in\{\pm 1\}
\textup{ with }\varepsilon_1\varepsilon_2\varepsilon_3=1
\end{array}\right\} & S(-l,2,-l)
\end{array}
\end{eqnarray*}
So there are seven single $\Br_3\ltimes \{\pm 1\}^3$ orbits
and one series with parameter $l\in\Z_{\geq 3}$ of
$\Br_3\ltimes\{\pm 1\}^3$ orbits with
$r\in\{0,1,2,3,4\}$.
These are the most interesting orbits as the monodromy matrix
$S(\uuuu{x})^{-1}S(\uuuu{x})^t$ for $\uuuu{x}\in\R^3$
has eigenvalues in $S^1$ if and only if 
$r(\uuuu{x})\in[0,4]$.

(f) For a given local minimum $\uuuu{x}\in\Z^3$ the set of
all local minima in the $\Br_3\ltimes\{\pm 1\}^3$ orbit
of $\uuuu{x}$ is either the set 
$G^{sign}\rtimes \langle\gamma\rangle(\uuuu{x})$ or the set
$G^{sign}\rtimes \langle\gamma,\gamma_2\rangle(\uuuu{x})$
(see Theorem \ref{t4.13} (b) for details). 
\end{theorem}

{\bf Proof:}
(a) Fix $\rho\in\Z-\{4\}$. Let $\uuuu{x}\in \Z^3$ be a local
minimum with $r(\uuuu{x})=\rho$. 

{\bf 1st case,} $x_1x_2x_3\leq 0$:
Then $\rho=r(\uuuu{x})=x_1^2+x_2^2+x_3^2-x_1x_2x_3\geq \|\uuuu{x}\|^2$,
so $\rho\geq 0$. 
The closed ball of radius $\sqrt{\rho}$ around $0$ in
$\R^3$ intersects $\Z^3$ only in finitely many points.

{\bf 2nd case,} $x_1x_2x_3>0$: 
We can suppose $x_i>0$ for $i\in\{1,2,3\}$ because of the
action of $G^{sign}$ on $\R^3$ and $\Z^3$. 
Lemma \ref{t4.4} says $2x_1\leq x_2x_3$, $2x_2\leq x_1x_3$, 
$2x_3\leq x_1x_2$, $x_i\geq 2$ for $i\in\{1,2,3\}$.

We can suppose $x_1=\min(x_1,x_2,x_3)$ (the other cases are
analogous). If $x_1=2$ then
$4\neq \rho=r(\uuuu{x})=4+x_2^2+x_3^2-2x_2x_3=4+(x_2-x_3)^2$,
so $x_2\neq x_3$, which is a contradiction to
$2x_2\leq x_1x_3=2x_3$, $2x_3\leq x_1x_2=2x_2$.
Therefore $x_1\geq 3$. 

We can suppose $x_1\leq x_2\leq x_3$ (the other cases are
analogous). 
\begin{eqnarray*}
\rho&=& r(\uuuu{x}) = x_1^2+x_2^2+(x_3-\frac{1}{2}x_1x_2)^2
-\frac{1}{4}x_1^2x_2^2\\
&\leq & x_1^2+x_2^2+(x_2-\frac{1}{2}x_1x_2)^2 -
\frac{1}{4}x_1^2x_2^2
\quad(\textup{because }x_2\leq x_3\leq \frac{1}{2}x_1x_2)\\
&=& x_1^2+2x_2^2-x_1x_2^2
=(x_1-2)(x_1+2-x_2^2)+4 \\
&\leq& (x_1-2)(x_1+2-x_1^2) + 4 
=-(x_1-2)(x_1-2)(x_1+1) + 4 \\
&\leq& \left\{\begin{array}{ll}
-(3-2)(3-2)(3+1)+4 =0, \\
-(x_1-2)^3+4.\end{array}\right. 
\end{eqnarray*}
This implies $\rho\leq 0$ and $x_1\leq 2+\sqrt[3]{4-\rho}$,
so $x_1$ is one of the finitely many values in 
$\Z\cap [3,2+\sqrt[3]{4-\rho}]$. 

The inequality $\rho\leq (x_1-2)(x_1+2-x_2^2)+4$ implies
\begin{eqnarray*}
x_2^2\leq \frac{4-\rho}{x_1-2}+x_1+2,
\end{eqnarray*}
so $x_2$ is one of the finitely many values in
$\Z\cap[x_1,\sqrt{\frac{4-\rho}{x_1-2}+x_1+2}]$.

Because of $x_3\leq \frac{1}{2}x_1x_2$ also $x_3$ can take
only finitely many values. 

(b) Each $\Br_3\ltimes\{\pm 1\}^3$ orbit in $\Z^3$ is 
mapped by $\|.\|^2$ to a subset of $\Z_{\geq 0}$.
A preimage in this orbit of the minimum of this subset is a local
minimum. Therefore (a) implies (b).

(c) Suppose $\rho<0$ and $\uuuu{x}\in r^{-1}(\rho)$ is a 
local minimum. $0>\rho=\|\uuuu{x}\|^2-x_1x_2x_3$ implies
$x_1x_2x_3>0$. Lemma \ref{t4.4} gives $|x_1|\geq 2$. 

If $x_1=2\varepsilon$ with $\varepsilon\in\{\pm 1\}$ then
$\rho=r(\uuuu{x})=4+(x_2-\varepsilon x_3)^2\geq 4$,
a contradiction. So $|x_1|\geq 3$.
Analogously $|x_2|\geq 3$ and $|x_3|\geq 3$.

(d) Suppose $\rho\in\N-\{4\}$ and $\uuuu{x}\in r^{-1}(\rho)$
is a local minimum. In the second case $x_1x_2x_3>0$ in the
proof of part (a) we concluded $\rho\leq 0$. 
Therefore we are in the first case in the proof of part (a),
so $x_1x_2x_3\leq 0$. 

(e) Suppose $\rho\in\{0,1,2,3,4\}$, and $\uuuu{x}\in r^{-1}(\rho)$
is a local minimum.

In the cases $\rho\in\{1,2,3\}$ by part (d) $x_1x_2x_3\leq 0$
and $\rho=r(\uuuu{x})=x_1^2+x_2^2+x_3^2+|x_1x_2x_3|$, 
so in these cases all $x_i\neq 0$ is impossible, so some
$x_i=0$, so $\rho=r(\uuuu{x})=x_j^2+x_k^2$ where
$\{i,j,k\}=\{1,2,3\}$. 

{\bf The case $\rho=3$:}
$3=x_j^2+x_k^2$ is impossible, the case $\rho=3$ is impossible,
$r^{-1}(3)=\emptyset$.

{\bf The case $\rho=1$:}
$1=x_j^2+x_k^2$ is solved only by 
$(x_j,x_k)\in\{(\pm 1,0),(0,\pm 1)\}$. 
The six local minima $(\pm 1,0,0),(0,\pm 1,0),(0,0,\pm 1)$
are in one orbit of $\Br_3\ltimes\{\pm 1\}^3$ because
$\gamma(1,0,0)=(0,1,0)$, $\gamma(0,1,0)=(0,0,1)$. 

{\bf The case $\rho=2$:}
$x_j^2+x_k^2=2$ is solved only by 
$(x_j,x_k)\in\{(\pm 1,\pm 1)\}$.
The twelve local minimal $(\pm 1,\pm 1,0),(\pm 1,0,\pm 1),
(0,\pm 1,\pm 1)$ are in one orbit of $\Br_3\ltimes\{\pm 1\}^3$
because $\gamma(1,1,0)=(0,1,1)$, $\gamma(0,1,1)=(1,0,1)$.

{\bf The case $\rho=0$:}
We use the proof of part (a).

{\bf 1st case,} $x_1x_2x_3\leq 0$: 
$0=\rho=\|x\|^2-x_1x_2x_3\geq\|x\|^2$, so $\uuuu{x}=(0,0,0)$.
Its $\Br_3\ltimes\{\pm 1\}^3$ orbit consists only of
$(0,0,0)$. 

{\bf 2nd case,} $x_1x_2x_3>0$: 
We can suppose $x_i>0$ for each $i\in\{1,2,3\}$.
Suppose $x_i\leq x_j\leq x_k$ for $\{i,j,k\}=\{1,2,3\}$.
The proof of part (a) gives
\begin{eqnarray*}
3\leq x_i\leq 2+\sqrt[3]{4-\rho}=2+\sqrt[3]{4},
\textup{ so }x_i=3
\end{eqnarray*}
and
\begin{eqnarray*}
3=x_i\leq x_j\leq \sqrt{\frac{4-\rho}{x_i-2}+x_i+2}=3,
\textup{ so }x_j=3.
\end{eqnarray*}
$0=r(\uuuu{x})=9+9+x_k^2-9x_k=(x_k-3)(x_k-6)$ and
$x_k\leq \frac{1}{2}x_ix_j=\frac{9}{2}$ show $x_k=3$.
The four local minima $(3,3,3)$, $(-3,-3,3)$, $(-3,3,-3)$
and $(3,-3,-3)$ are in one $\Br_3\ltimes \{\pm 1\}^3$ orbit
because of the action of $G^{sign}$. 

{\bf The case $\rho=4$:}

{\bf 1st case,} some $x_i=0$: Then with 
$\{i,j,k\}=\{1,2,3\}$ $4=r(\uuuu{x})=x_j^2+x_k^2$.
This is solved only by $(x_j,x_k)\in\{(\pm 2,0),(0,\pm 2)\}$.
The six local minima $(\pm 2,0,0)$, $(0,\pm 2,0)$, $(0,0,\pm 2)\}$
are in one $\Br_3\ltimes\{\pm 1\}^3$ orbit because
$\gamma((2,0,0))=(0,2,0)$, $\gamma((0,2,0))=(0,0,2)$.

{\bf 2nd case,} all $x_i\neq 0$ and $x_1x_2x_3<0$: 
$4=r(\uuuu{x})=x_1^2+x_2^2+x_3^2+|x_1x_2x_3|$, so 
$(x_1,x_2,x_3)\in \{(-1,-1,-1),(1,1,-1),(1,-1,1),(-1,1,1)\}$.
These four local minima are in one $\Br_3\ltimes\{\pm 1\}^3$
orbit because of the action of $G^{sign}$. 

{\bf 3rd case,} all $x_i\neq 0$ and $x_1x_2x_3>0$: 
We can suppose $x_i>0$ for each $i\in\{1,2,3\}$ and
$x_i\leq x_j\leq x_k$ for some $i,j,k$ with 
$\{i,j,k\}=\{1,2,3\}$. 
As in the proof of part (a) we obtain the estimate
\begin{eqnarray*}
4=\rho=r(\uuuu{x})\leq -(x_i-2)^3+4,\quad\textup{so }x_i=2,
\end{eqnarray*}
and
\begin{eqnarray*}
4=\rho=r(\uuuu{x})=4+(x_j-x_k)^2,\quad\textup{so }
l:=x_j=x_k\geq 2.
\end{eqnarray*}
For $l=2$ the four local minima
\begin{eqnarray*}
(2,2,2),(-2,-2,2),(-2,2,-2),(2,-2,-2)
\end{eqnarray*}
and for $l\geq 3$ the 24 local minima 
\begin{eqnarray*}
(\varepsilon_1 2,\varepsilon_2 l,\varepsilon_3 l),
(\varepsilon_1 l,\varepsilon_2 2,\varepsilon_3 l),
(\varepsilon_1 l,\varepsilon_2 l,\varepsilon_3 2)\\
\textup{ with }\varepsilon_1,\varepsilon_2,\varepsilon_3
\in\{\pm 1\},
\varepsilon_1\varepsilon_2\varepsilon_3=1,
\end{eqnarray*}
are in one $\Br_3\ltimes\{\pm 1\}^3$ orbit because of the
action of $G^{sign}$ and $\gamma$.

It remains to see that local minima in different lines in
the list in part (e) are in different $\Br_3\ltimes\{\pm 1\}^3$
orbits. One reason is part (f). Another way to argue is given
in the Remarks \ref{t4.7}.

(f) See Lemma \ref{t4.10} (e). 
\hfill$\Box$

\begin{remarks}\label{t4.7}
Part (f) of Theorem \ref{t4.6} is strong and allows easily to
see when $\Br_3\ltimes\{\pm 1\}^3$ orbits are separate. 
Nevertheless it is also interesting
to find invariants of the orbits which separate them.

Now we discuss several invariants which help to prove the claim 
that local minima in different lines in the list in part (e) 
are in different  $\Br_3\ltimes\{\pm 1\}^3$ orbits. 

The number $r(\uuuu{x})\in\Z$ is such an invariant.
Furthermore the set $\{(0,0,0)\}$ is a single orbit and thus
different from the orbit of $S(\P^2)$. 
Therefore the claim is true for the lines with $r\in\{0,1,2,3\}$. 

It remains to consider the $3+\infty$ lines with $r=4$.
Certainly the reducible case $S(\P^1A_1)$ is separate from the other
cases, which are all irreducible. The signature of $I^{(0)}$ is an
invariant. It is given in Lemma \ref{t5.7}. It allows to see that
the orbits of $S(\whh{A}_2)$, $S(\HH_{1,2})$ are different from one another
and from the orbits of $S(-l,2,-l)$ for $l\geq 3$. 
In order to see that the orbits of $S(-l,2,-l)$ for $l\geq 3$
are pairwise different, we can offer Lemma \ref{t7.10},
which in fact allows to separate all the lines with $r=4$. 
It considers the induced monodromy on the quotient lattice
$H_\Z/\Rad I^{(1)}$. 
\end{remarks}

\begin{remarks}\label{t4.8}
Kr\"uger \cite[\S 12]{Kr90} and Cecotti-Vafa
\cite[Ch. 6.2]{CV93} considered the action of the 
group $(G^{phi}\ltimes G^{sign})\rtimes
\langle\gamma,\gamma_2\rangle$ 
with $\gamma_2:\R^3\to\R^3$, $(x_1,x_2,x_3)\mapsto (x_2,x_1,x_3)$,
\index{$\gamma_2:\R^3\to\R^3$} 
so $\langle \gamma,\gamma_2\rangle\cong S_3$ 
which is slightly larger than
$(G^{phi}\ltimes G^{sign})\rtimes\langle\gamma\rangle$.
Because of
\begin{eqnarray*}
\gamma_2\varphi_1\gamma_2^{-1}=\varphi_2, \
\gamma_2\varphi_2\gamma_2^{-1}=\varphi_1, \
\gamma_2\varphi_3\gamma_2^{-1}=\varphi_3, \
\gamma_2\gamma\gamma_2^{-1}=\gamma^{-1},
\end{eqnarray*}
we
have
\begin{eqnarray*}
\Br_3\ltimes\{\pm 1\}^3\bigl(\gamma_2(\uuuu{x})\bigr)
=\gamma_2\bigl(\Br_3\ltimes\{\pm 1\}^3(\uuuu{x})\bigr).
\end{eqnarray*}
Especially, the $\Br_3\ltimes \{\pm 1\}^3$ orbit of 
$\uuuu{x}$ coincides with the 
$(G^{phi}\ltimes G^{sign})\rtimes \langle\gamma,\gamma_2\rangle$
orbit of $\uuuu{x}$ in the following cases:
\begin{list}{}{}
\item[(i)]
if $x_i=x_j$ for some $i\neq j$,
\item[(ii)]
if $x_i=0$ for some $i$ (observe $\delta_3^{\R}\gamma
\varphi_1(x_1,x_2,0)=(x_2,x_1,0)$), 
\item[(iii)]
if $\uuuu{x}=(x_1,x_2,\frac{1}{2}x_1x_2)$ with
$|x_i|\geq 3$ and $|x_1|\neq |x_2|$
(observe $\varphi_3(\uuuu{x})=(x_2,x_1,\frac{1}{2}x_1x_2)$). 
\end{list}
In Lemma \ref{t4.12} 24 sets $C_1,...,C_{24}$ of local
minima are considered. The only local minima 
$\uuuu{x}\in\bigcup_{i=1}^{24} C_i$ which satisfy
none of the conditions (i), (ii) and (iii) are those in
$C_{16}\cup C_{22}\cup C_{24}$. Theorem \ref{t4.13} (b) shows
that in these cases the $\Br_3\ltimes\{\pm 1\}^3$ 
orbits of $\uuuu{x}$ and of $\gamma_2(\uuuu{x})$ are 
indeed disjoint.

Especially all orbits in the fibers $r^{-1}(\rho)$
with $\rho\in\{0,1,2,3,4\}$ contain local minima
which satify (i), (ii) or (iii), 
so there the classifications in Theorem \ref{t4.6}
and the classification by Kr\"uger and Cecotti-Vafa
coincide. 

We do not know whether for $\uuuu{x}\in C_{16}\cup
C_{22}\cup C_{24}$ and $\gamma_2(\uuuu{x})$ the 
corresponding unimodular bilinear lattices with sets
of distinguished bases are isomorphic or not. 

For $\uuuu{x}$ and $\gamma_2(\uuuu{x})$
in one $\Br_3\ltimes\{\pm 1\}^3$ orbit they are isomorphic,
see Remark \ref{t3.20} (iii). 
\end{remarks}

\section{A classification of the $\Br_3\ltimes \{\pm 1\}^3$ orbits in $\Z^3$}
\label{s4.3}

This section refines the results of section \ref{s4.2}
on the braid group action on integer upper triangular
$3\times 3$ matrices. Using pseudo-graphs, it gives a 
classification of all orbits of $\Br_3$ on $\Z^3/\{\pm 1\}^3$.
Definition \ref{t4.9} makes precise what is meant here by
a pseudo-graph, and it defines a pseudo-graph $\GG(\uuuu{x})$ 
for any local minimum $\uuuu{x}\in\Z^3$. 

As $\Br_3\ltimes\{\pm 1\}^3$ and $(G^{phi}\ltimes\{\pm 1\}^3)
\rtimes\langle \gamma\rangle$ are semidirect products with
normal subgroups $\{\pm 1\}^3$, the groups $\Br_3$ and
$G^{phi}\rtimes\langle\gamma\rangle$ act on
$\Z^3/\{\pm 1\}^3$.

\begin{definition}\label{t4.9}
(a) For any set $\VV$, $\PP(\VV)$ denotes its power set,
so the set of all its subsets, and $\PP_k(\VV)$ for some
$k\in\N$ denotes the set of all subsets with $k$ elements.
We will use only $\PP_1(\VV)$ and $\PP_2(\VV)$.

(b) A {\it pseudo-graph} \index{pseudo-graph}
is here a tuple
$\GG=(\VV,\VV_0,\VV_1,\VV_2,v_0,\EE_1,\EE_2,\EE_3,\EE_\gamma)$
with the following ingredients:

$\VV$ is a non-empty 
finite or countably infinite set of vertices.

$\VV_0,\VV_1,\VV_2\subset\VV$ are pairwise disjoint subsets,
$\VV_0$ is not empty (the sets $\VV_1$ and $\VV_2$ may be empty,
the union $\VV_0\cup\VV_1\cup\VV_2$ 
can be equal to $\VV$ or a proper subset of $\VV$).

$v_0\in\VV_0$ is a distinguished vertex in $\VV_0$.

$\EE_1,\EE_2,\EE_3\subset \PP_1(\VV)\cup\PP_2(\VV)$ are 
sets of undirected edges.
A subset of $\VV$ with two elements means an edge between
the two vertices. A subset of $\VV$ with one element means
a loop from the vertex to itself.

$\EE_\gamma=\{(v_0,v_1),(v_2,v_0)\}$ for some $v_1,v_2\in\VV_0$
is a set of two or one directed edges, one only if
$v_1=v_2=v_0$, and then it is a directed loop.

(c) An isomorphism between two pseudo-graphs $\GG$ and 
$\www{\GG}$ is a bijection $\phi:\VV\to\www{\VV}$ with
$\phi(v_0)=\www{v_0}$ which induces bijections
$\phi:\VV_i\to \www{\VV_i}$, $\phi:\EE_j\to\www{\EE_j}$ and 
$\phi:\EE_\gamma\to\www{\EE_\gamma}$.

(d) $\GG|_{\VV_0\cup\VV_1}$ denotes the restriction of a 
pseudo-graph $\GG$ to the vertex set $\VV_0\cup\VV_1$,
so one deletes all vertices in $\VV-(\VV_0\cup\VV_1)$
and all edges with at least one end in 
$\VV-(\VV_0\cup\VV_1)$. Analogously, $\GG|_{\VV_0}$ denotes
the restriction of a pseudo-graph $\GG$ to the vertex set $\VV_0$.

(e) Define
\begin{eqnarray*}
\LL_0&:=& \{\uuuu{x}/\{\pm 1\}^3\,|\, \uuuu{x}\in\Z^3
\textup{ is a local minimum}\},\\
\LL_1&:=& \{\uuuu{y}/\{\pm 1\}^3\,|\,\uuuu{y}\in\Z^3\textup{ with } 
|y_i|=1\textup{ for some }i\}-\LL_0,\\
\LL_2&:=& \Z^3/\{\pm 1\}^3-(\LL_0\cup\LL_1).
\end{eqnarray*}

(f) A pseudo-graph $\GG(\uuuu{x})$ is associated to a local
minimum $\uuuu{x}\in\Z^3$ in the following way:
\begin{eqnarray*}
\VV&:=&\Br_3(\uuuu{x}/\{\pm 1\}^3)\subset\Z^3/\{\pm 1\}^3,\\
\VV_0&:=& \VV\cap \LL_0\textup{ is the set of sign classes
in }\VV\textup{ of local minima},\\
\VV_1&:=& \VV\cap \LL_1,\\
\VV_2&:=& \{w\in\VV\cap\LL_2\,|\, \textup{an }
i\in\{1,2,3\}\textup{ with }
\varphi_i(w)\in\VV_0\cup\VV_1\textup{ exists}\},\\
v_0&:=& \uuuu{x}/\{\pm 1\}^3,\\
\EE_i&:=& \{\{w,\varphi_i(w)\}\,|\, w\in \VV\} 
\quad\textup{ for }i\in\{1,2,3\},\\
\EE_\gamma&:=& \{(v_0,\gamma(v_0)),(\gamma^{-1}(v_0),v_0)\}.
\end{eqnarray*}

(g) An {\it infinite tree} \index{infinite tree} 
$(\WW,\FF)$ consists of a 
countably infinite set $\WW$ of vertices and a set 
$\FF\subset\PP_2(\WW)$ of undirected
edges such that the graph is connected and has no cycles.
A $(2,\infty\times 3)$-tree 
\index{$(2,\infty\times 3)$-tree}
is an infinite tree with a 
distinguished vertex with two neighbours such that any other
vertex has three neighbours.
\end{definition}

The next lemma gives already structural results about the
pseudo-graphs $\GG(\uuuu{x})$ for the local minima 
$\uuuu{x}\in\Z^3$. Theorem \ref{t4.13} and the Remarks
\ref{t4.14} will give a complete classification of all
isomorphism classes of pseudo-graphs $\GG(\uuuu{x})$ for
the local minima $\uuuu{x}\in\Z^3$.

\begin{lemma}\label{t4.10}
Let $\uuuu{x}\in\Z^3$ be a local minimum with pseudo-graph
$\GG(\uuuu{x})=(\VV,\VV_0,\VV_1,\VV_2,v_0,\EE_1,\EE_2,\EE_3,
\EE_\gamma)$. 

(a) ($\uuuu{x}$ and $\GG(\uuuu{x})$ are not used in part (a))
For $w\in\LL_2$, there are $i,j,k$ with $\{i,j,k\}=\{1,2,3\}$
and $\|\varphi_i(w)\|<\|w\|$, $\|\varphi_j(w)\|>\|w\|$,
$\|\varphi_k(w)\|>\|w\|$, $\varphi_j(w)\in\LL_2$,
$\varphi_k(w)\in\LL_2$ and $\varphi_j(w)\neq\varphi_k(w)$.

(b) Let $w\in\VV_2$, so especially $w\in\LL_2$.
Choose $i,j,k$ as in part (a). The edge which connects
$w\in\VV_2$ with $\VV_0\cup\VV_1$ is in $\EE_i$. After
deleting this edge, the component of the remaining 
pseudo-graph which contains $w$ is a 
$(2,\infty\times 3)$-tree with distinguished vertex $w$
and all vertices in $\LL_2$. 

(c) The pseudo-graph $\GG(\uuuu{x})$ is connected. 

(d) The pseudo-graph $\GG(\uuuu{x})|_{\VV_0\cup \VV_1}$
is connected.

(e) The pseudo-graph $\GG(\uuuu{x})|_{\VV_0\cup\VV_1}$ is finite,
and 
\begin{eqnarray*}
\VV_0=\langle\gamma\rangle (v_0)\quad\textup{or}\quad 
\VV_0=\langle\gamma,\gamma_2\rangle(v_0).
\end{eqnarray*}
\end{lemma}

{\bf Proof:}
(a) Consider $w=\uuuu{y}/\{\pm 1\}^3\in\LL_2$.
Then $|y_i|\geq 2$ for $i\in\{1,2,3\}$ because 
$w\notin\LL_0\cup\LL_1$. $w\notin\LL_0$ implies
$y_1y_2y_3>0$. We can suppose $y_1,y_2,y_3\in\Z_{\geq 2}$. 
Observe for $i,j,k$ with $\{i,j,k\}=\{1,2,3\}$
\begin{eqnarray*}
&&\|\varphi_j(\uuuu{y})\|^2-\|\varphi_i(\uuuu{y})\|^2\\
&=& (y_i^2+(y_iy_k-y_j)^2+y_k^2)
-((y_jy_k-y_i)^2+y_j^2+y_k^2)\\
&=& (y_i^2-y_j^2)y_k^2.
\end{eqnarray*}
Consider the case $2\leq y_1\leq y_2\leq y_3$.
The other cases are analogous. Then
\begin{eqnarray*}
\|\varphi_3(\uuuu{y})\|\leq \|\varphi_2(\uuuu{y})\|
\leq \|\varphi_1(\uuuu{y})\|.
\end{eqnarray*}
Because $w\notin\LL_0$ $\|\varphi_3(\uuuu{y})\|<\|\uuuu{y}\|$. 
Also $\varphi_2(\uuuu{y})=(y_3,y_1y_3-y_2,y_1)$ with
\begin{eqnarray*}
\varphi_2(\uuuu{y})_2=y_1y_3-y_2
\geq (y_1-1)y_3\left\{\begin{array}{ll}
\geq 2y_3>y_2&\textup{ if }y_1>2,\\
= y_3>y_2&\textup{ if }y_1=2.\end{array}\right.
\end{eqnarray*}
Here $y_3>y_2$ if $y_1=2$ because $\uuuu{y}=(2,y_2,y_3)$
is not a local minimum.

Therefore $\|\varphi_2(\uuuu{y})\|>\|\uuuu{y}\|$, and also
$\|\varphi_1(\uuuu{y})\| \geq \|\varphi_2(\uuuu{y})\|
>\|\uuuu{y}\|$. 
Especially $\varphi_2(\uuuu{y})$ and $\varphi_1(\uuuu{y})$
are not local minima, so 
$\varphi_2(w)\notin\LL_0$ and $\varphi_1(w)\notin\LL_0$.

Obviously $\varphi_2(\uuuu{y})_i\geq 2$ and  
$\varphi_1(\uuuu{y})_i\geq 2$ for $i\in\{1,2,3\}$, so 
$\varphi_2(w)\in\LL_2$ and $\varphi_1(w)\in\LL_2$.
The inequality $\varphi_2(w)\neq \varphi_1(w)$ follows from 
\begin{eqnarray*}
\varphi_2(\uuuu{y})_2\geq 2y_3>y_3=\varphi_1(\uuuu{y})_2
&&\textup{if }y_1>2,\\
\varphi_2(\uuuu{y})_2=y_1y_3-y_2>(y_1-1)y_3=y_3=\varphi_1(\uuuu{y})_2
&&\textup{if }y_1=2.
\end{eqnarray*}

(b) Because $\varphi_j(w),\varphi_k(w)\in\LL_2$, 
the edge which connects $w$ to $\VV_0\cup\VV_1$ cannot be in
$\EE_j$ or $\EE_k$, so it is in $\EE_i$. 

Using part (a) again and again, one sees that the component
of $\GG(\uuuu{x})-\textup{(this edge)}$ which contains 
$w$ is a $(2,\infty\times 3)$ tree.

(c) Any vertex $w\in\VV=(G^{phi}\rtimes \langle\gamma\rangle)
(v_0)$ is obtained from $v_0$ by applying an element 
$\psi\gamma^\xi$ with $\psi\in G^{phi}$ and $\xi\in\{0,\pm 1\}$.
As $G^{phi}$ is a free Coxeter group with generators
$\varphi_1$, $\varphi_2$, $\varphi_3$, applying
$\psi\gamma^\xi$ to $v_0$ yields a path in $\GG(\uuuu{x})$
from $v_0$ to $w$. 

(d) This follows from (b) and (c).

(e) Consider $w=\uuuu{y}/\{\pm 1\}^3\in \VV_1$. 
Because $w\notin \VV_0$, $y_1y_2y_3>0$. 
We can suppose $y_1,y_2,y_3\in\N$, and one of them is equal to $1$. 
Suppose $1=y_1\leq y_2\leq y_3$. Then
\begin{eqnarray}
\varphi_3(\uuuu{y})&=&(y_2,1,y_2-y_3),\textup{ so }
y_2\cdot 1\cdot (y_2-y_3)\leq 0,\nonumber\\
&&\textup{so }\varphi_3(w)\in\VV_0,\label{4.8}\\
\varphi_2(\uuuu{y})&=&(y_3,y_3-y_2,1),\nonumber\\
&&\textup{so }\varphi_2(w)\in\VV_0\cup\VV_1\ 
(\textup{in }\VV_0\textup{ only if }y_3=y_2),\label{4.9}\\
\varphi_1\varphi_2(\uuuu{y})&=& (-y_2,1,y_3-y_2),\nonumber\\
&&\textup{so }\varphi_1\varphi_2(w)=\varphi_3(w)\in\VV_0.
\label{4.10}
\end{eqnarray}
Especially, each vertex in $\VV_1$ is connected by an edge to a 
vertex in $\VV_0$. 

Therefore the main point is to show
$\VV_0=\langle \gamma\rangle(v_0)$ or
$\VV_0=\langle \gamma,\gamma_2\rangle(v_0)$.
Then $\VV_0$ and $\VV_1$ are finite. 

First case, the restricted pseudo-graph $\GG(\uuuu{x})|_{\VV_0}$
is connected: Then $\|w\|=\|v_0\|$ for each $w\in\VV_0$.
This easily implies 
$\VV_0=\langle\gamma\rangle(v_0)$ or 
$\VV_0=\langle\gamma,\gamma_2\rangle$.

Second case, the restricted pseudo-graph $\GG(\uuuu{x})|_{\VV_0}$
is not connected: We will lead this to a contradiction.
Within all paths in $\GG(\uuuu{x})|_{\VV_0\cup\VV_1}$ which
connect vertices in different components of $\GG(\uuuu{x})|_{\VV_0}$,
consider a shortest path. It does not contain an edge in $\EE_\gamma$,
because else one could go over to a path of the same length with
an edge in $\EE_\gamma$ at one end and between the same vertices, 
but dropping that edge would lead to a shorter path. 
Because each vertex in $\VV_1$ is connected by an edge to a 
vertex in $\VV_0$, a shortest path contains either one or two
vertices in $\VV_1$. The observations \eqref{4.8}--\eqref{4.10}
lead in both cases to the vertices at the end of the path 
being in the same component of $\GG(\uuuu{x})|_{\VV_0}$, 
so to a contradiction.
\hfill$\Box$

\begin{examples}\label{t4.11}
The following 14 figures show the pseudo-graphs 
$\GG_1,...,\GG_{14}$ 
\index{$\GG_1,...,\GG_{14}$}
for 14 values
$v_0=\uuuu{x}/\{\pm 1\}^3$ with $\uuuu{x}\in\Z^3$ a local minimum.
The ingredients of the figures have the following meaning.

$\bullet$\qquad a vertex in $\VV_0$,

{\tiny{$\otimes$}} \qquad a vertex in $\VV_1$,

\noindent \hspace*{0.14cm}
\includegraphics[width=0.04\textwidth]{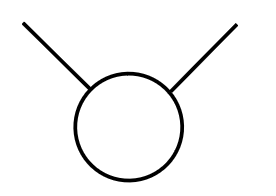}
\qquad a vertex in $\VV_2$ together with the 
$(2,\infty\times 3)$ tree (compare Lemma \ref{t4.10} (b)),

$\stackrel{i}{\mbox{-----}}$ \qquad an edge in $\EE_i$,

$\stackrel{\gamma}{\rightarrow\hspace*{-0.3cm}\mbox{---}}$ 
\qquad an edge in $\EE_\gamma$. 

The pseudo-graphs are enriched in the following way. 
Each vertex is labeled with a value $\uuuu{y}$ of its sign
class $\uuuu{y}/\{\pm 1\}^3$. We have chosen
$\uuuu{y}\in\N^3$ if $y_1y_2y_3>0$ 
(this holds for all $\uuuu{y}/\{\pm 1\}^3\in\VV_1\cup\VV_2$
and some $\uuuu{y}/\{\pm 1\}^3\in\VV_0)$
and $\uuuu{y}\in\Z_{\leq 0}^3$
if $y_1y_2y_3\leq 0$ (this holds for some 
$\uuuu{y}/\{\pm 1\}^3\in \VV_0$). 
The vertex $v_0$ can be recognized by the edges in $\EE_\gamma$
leading to and from it. 
The sets $C_i$ are defined in Lemma \ref{t4.12}.
The relations $\GG_j\,:\, C_i$ are explained in 
Theorem \ref{t4.13}.
\end{examples}

\begin{figure}[H]
\includegraphics[width=0.3\textwidth]{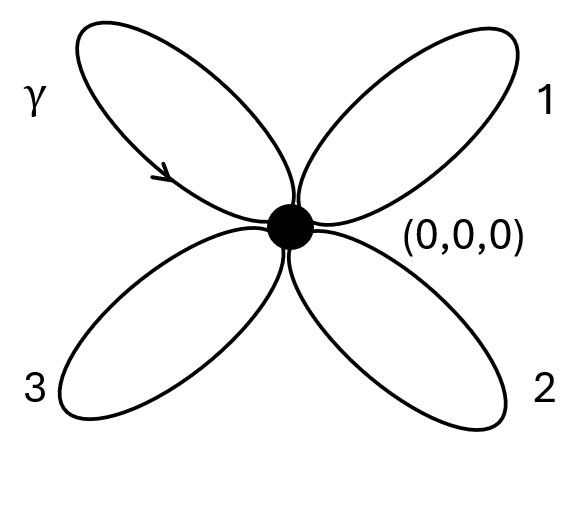}
\caption[Figure 4.2]{$G_1:\ C_1\ (A_1^3),\ C_2\ (\HH_{1,2}).
\quad \textup{Here }\uuuu{x}=(0,0,0)\in C_1.$}
\label{Fig:4.2}
\end{figure}

\begin{figure}[H]
\includegraphics[width=0.5\textwidth]{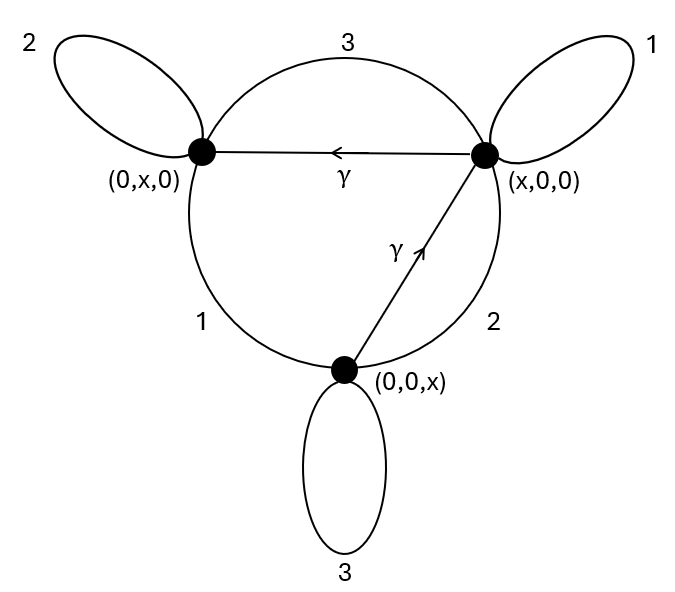}
\caption[Figure 4.3]{$G_2:\ C_3, C_4, C_5$, so the reducible
cases without $A_1^3$. 
\quad Here $\uuuu{x}=(x,0,0)\in C_3\cup C_4\cup C_5,\ x<0.$}
\label{Fig:4.3}
\end{figure}

\begin{figure}[H]
\includegraphics[width=0.5\textwidth]{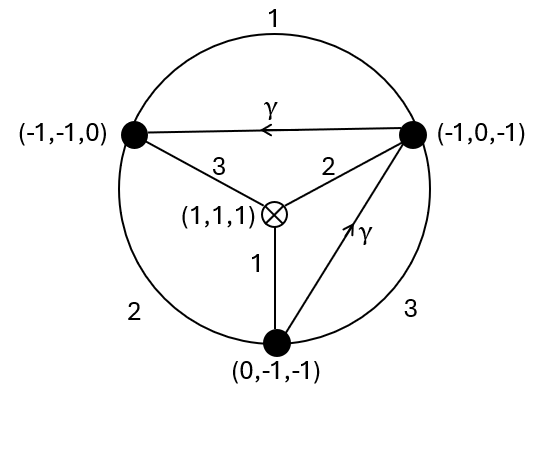}
\caption[Figure 4.4]{$G_3:\ C_6\ (A_3)$. 
\quad Here $\uuuu{x}=(-1,0,-1)$.}
\label{Fig:4.4}
\end{figure}

\begin{figure}[H]
\includegraphics[width=0.5\textwidth]{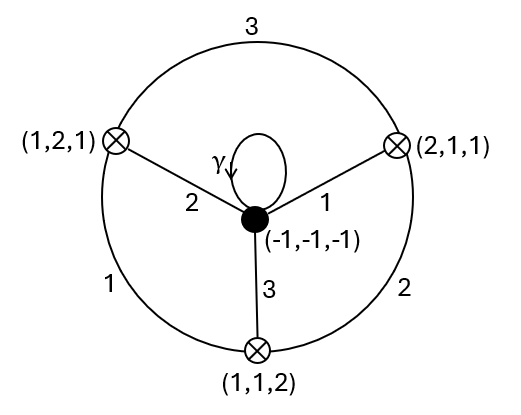}
\caption[Figure 4.5]{$G_4:\ C_7\ (\widehat{A}_2)$.
\quad Here $\uuuu{x}=(-1,-1,-1)$.}
\label{Fig:4.5}
\end{figure}

\begin{figure}[H]
\includegraphics[width=0.7\textwidth]{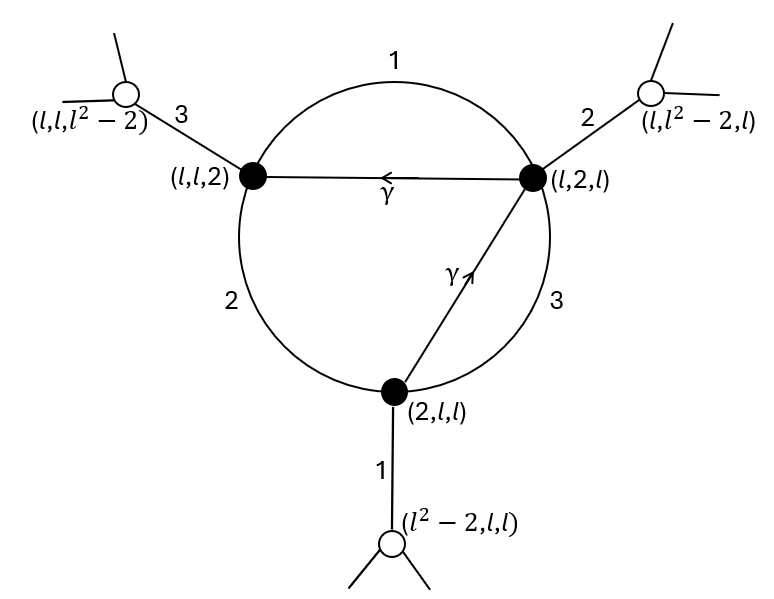}
\caption[Figure 4.6]{$G_5:\ C_8,C_9$.
\quad Here $\uuuu{x}=(l,2,l)\sim (-l,2,-l)$ with $l\geq 3$.}
\label{Fig:4.6}
\end{figure}

\begin{figure}[H]
\includegraphics[width=0.5\textwidth]{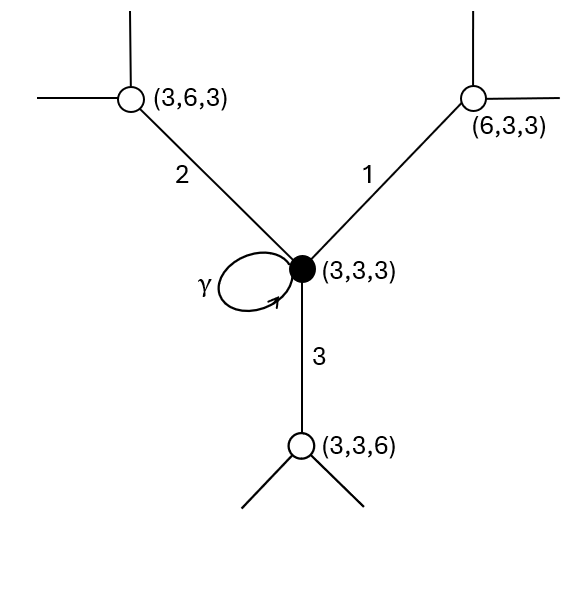}
\caption[Figure 4.7]{$G_6:\ C_{10},C_{11},C_{12}.$
\quad Here $\uuuu{x}=(3,3,3)\in C_{10}$.}
\label{Fig:4.7}
\end{figure}

\begin{figure}[H]
\includegraphics[width=0.7\textwidth]{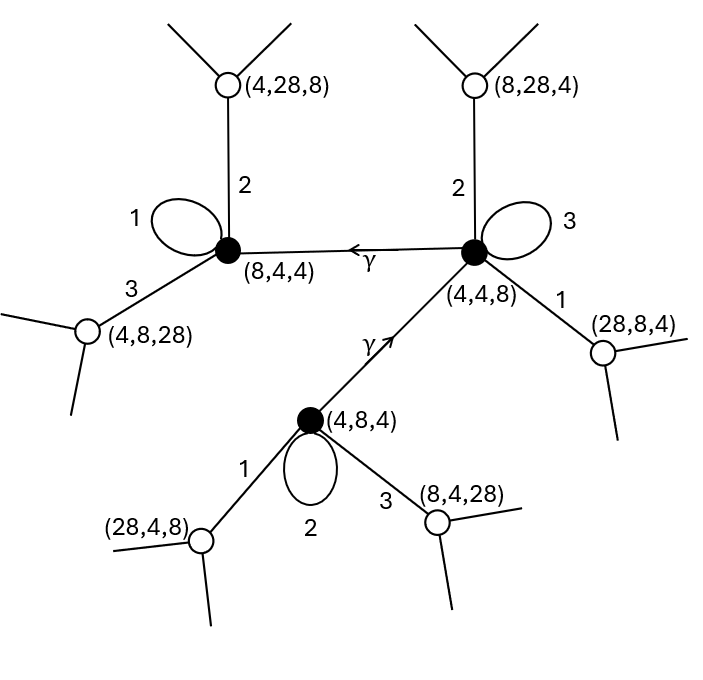}
\caption[Figure 4.8]{$G_7:\ C_{13}.$
\quad Here $\uuuu{x}=(4,4,8)$.}
\label{Fig:4.8}
\end{figure}

\begin{figure}[H]
\includegraphics[width=0.85\textwidth]{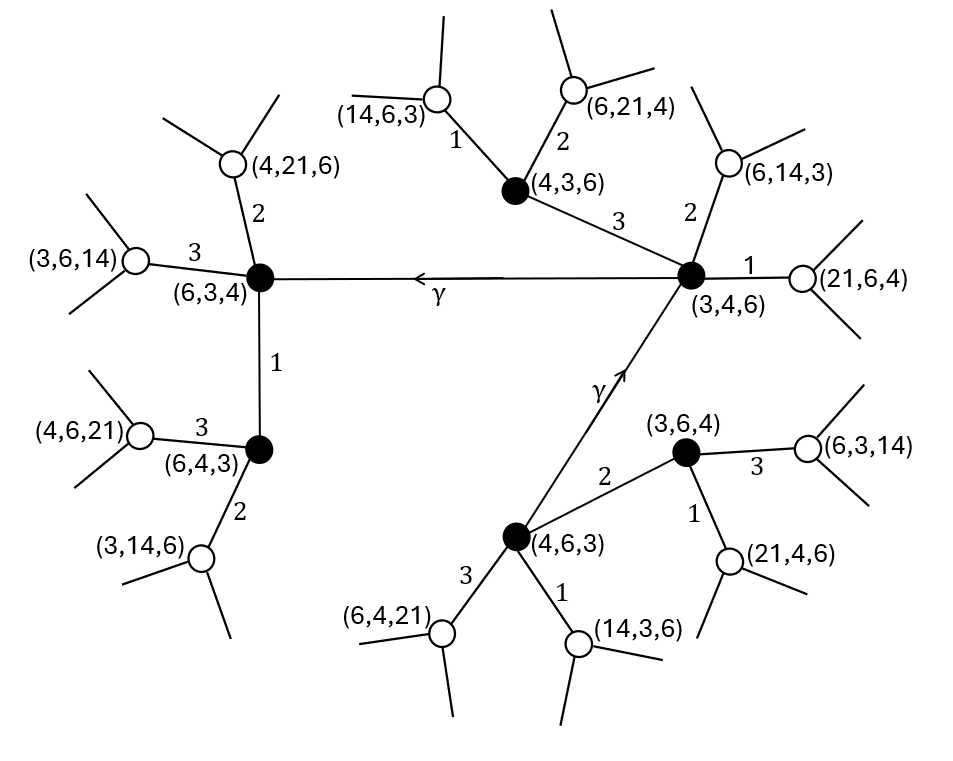}
\caption[Figure 4.9]{$G_8:\ C_{14}.$
\quad Here $\uuuu{x}=(3,4,6)$.}
\label{Fig:4.9}
\end{figure}

\begin{figure}[H]
\includegraphics[width=0.85\textwidth]{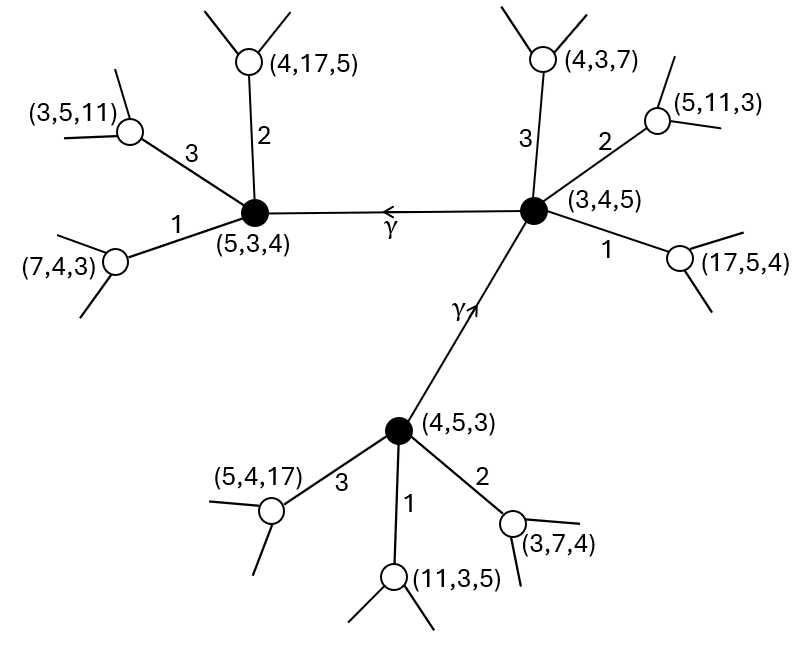}
\caption[Figure 4.10]{$G_9:\ C_{15},C_{16},C_{23}, C_{24}.$
\quad Here $\uuuu{x}=(3,4,5)\in C_{16}$.}
\label{Fig:4.10}
\end{figure}

\begin{figure}[H]
\includegraphics[width=0.75\textwidth]{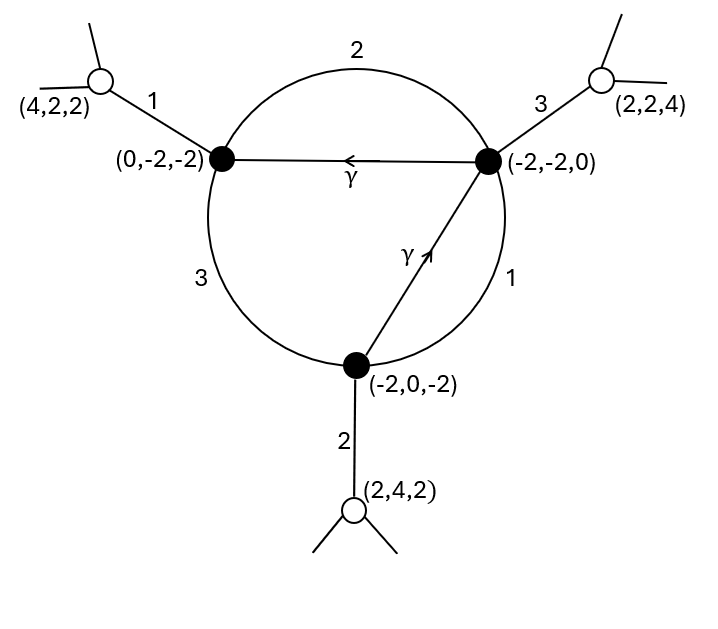}
\caption[Figure 4.11]{$G_{10}:\ C_{17}.$
\quad Here $\uuuu{x}=(-2,-2,0)$.}
\label{Fig:4.11}
\end{figure}

\begin{figure}[H]
\includegraphics[width=0.65\textwidth]{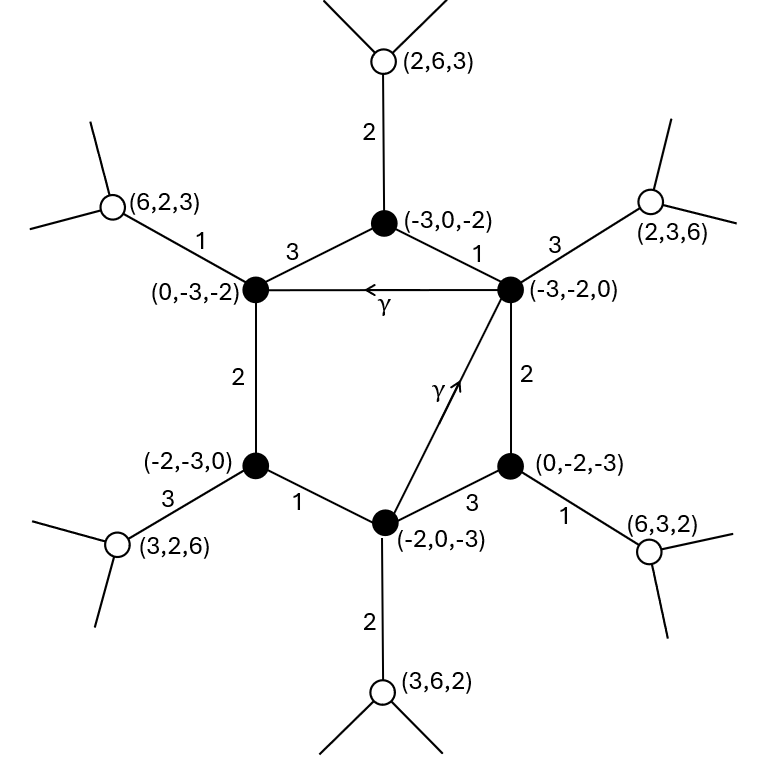}
\caption[Figure 4.12]{$G_{11}:\ C_{18}.$
\quad Here $\uuuu{x}=(-3,-2,0)$.}
\label{Fig:4.12}
\end{figure}

\begin{figure}[H]
\includegraphics[width=0.8\textwidth]{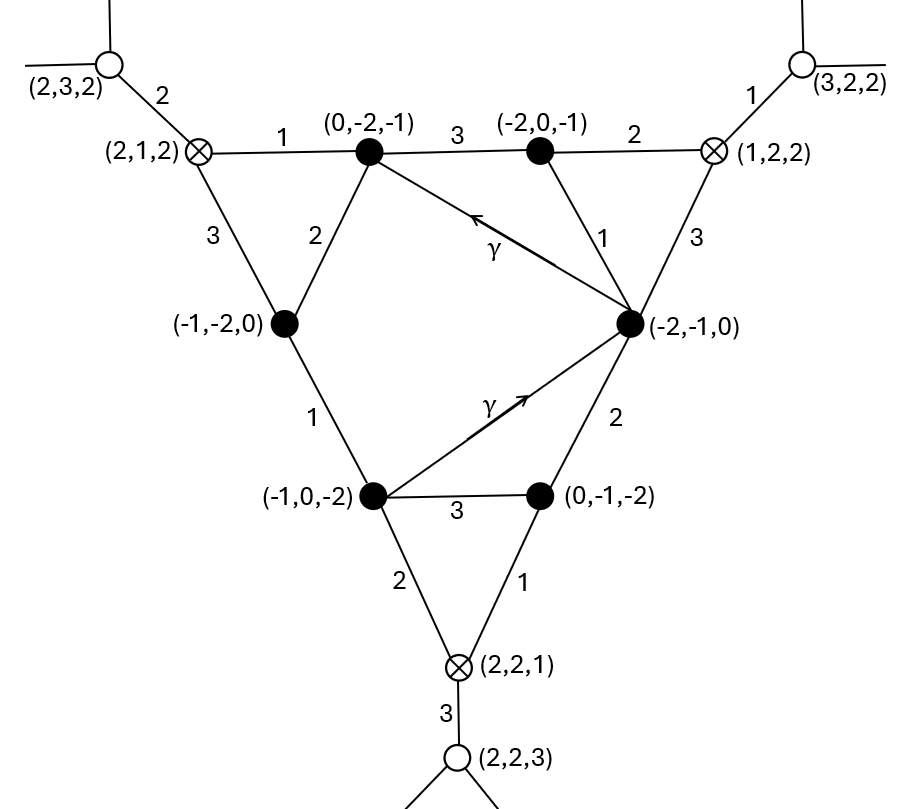}
\caption[Figure 4.13]{$G_{12}:\ C_{19}.$
\quad Here $\uuuu{x}=(-2,-1,0)$.}
\label{Fig:4.13}
\end{figure}

\begin{figure}[H]
\includegraphics[width=0.85\textwidth]{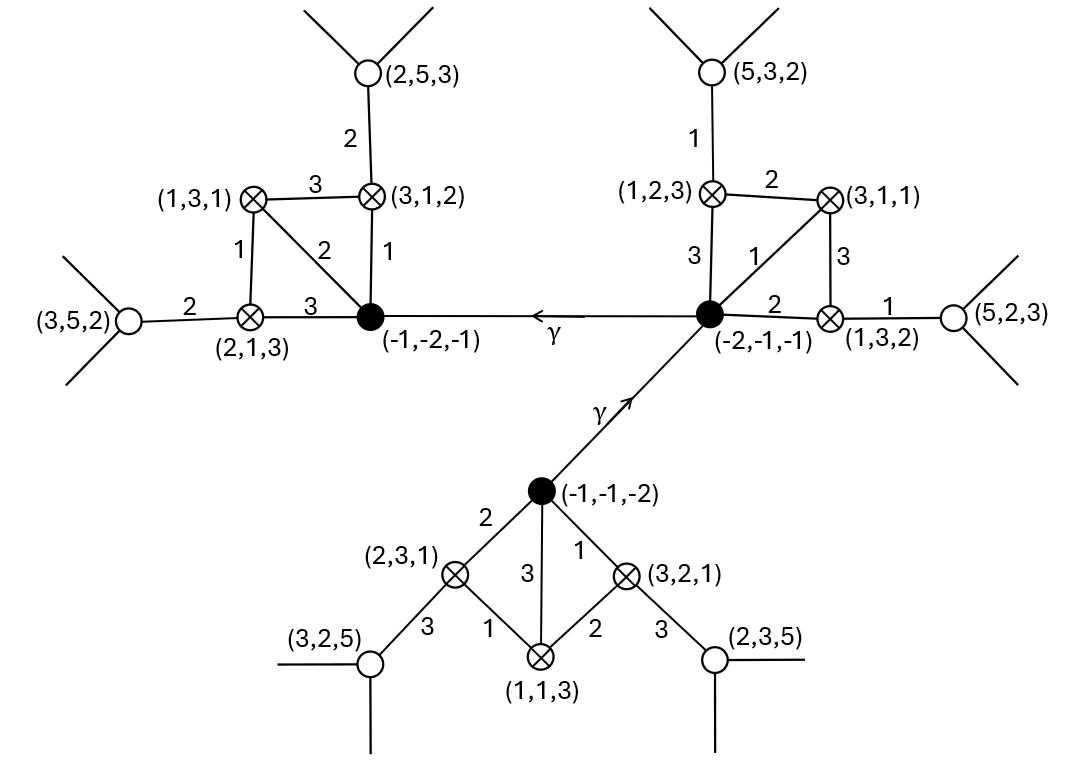}
\caption[Figure 4.14]{$G_{13}:\ C_{20}.$
\quad Here $\uuuu{x}=(-2,-1,-1)$.}
\label{Fig:4.14}
\end{figure}

\begin{figure}[H]
\includegraphics[width=0.85\textwidth]{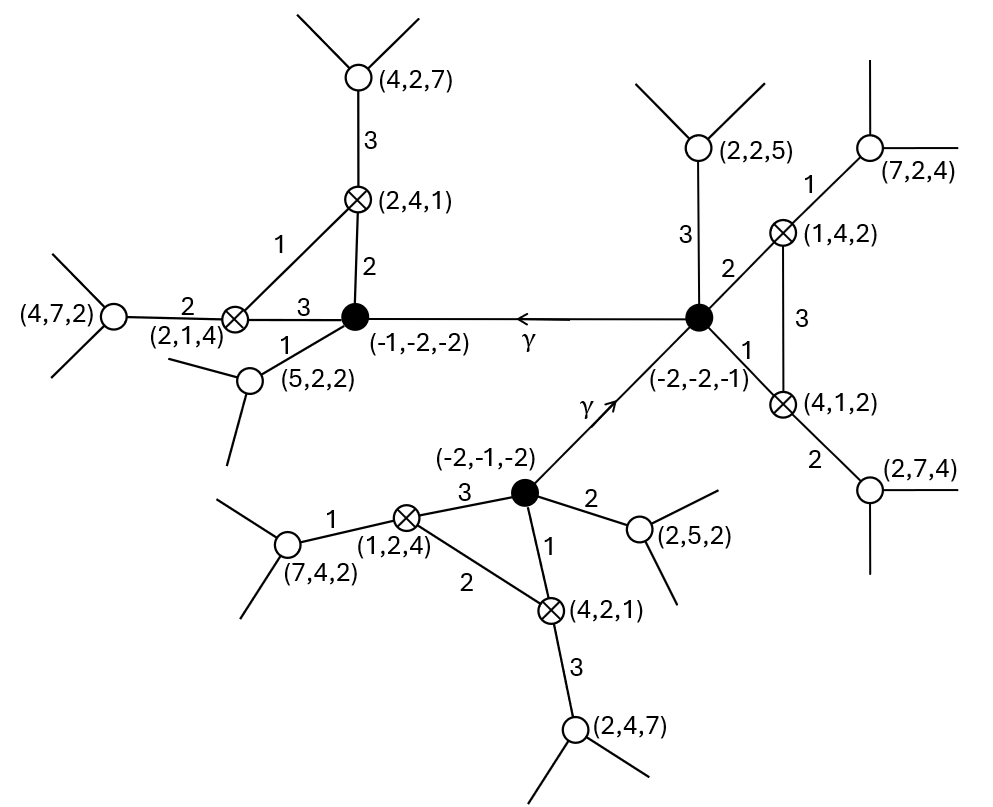}
\caption[Figure 4.15]{$G_{14}:\ C_{21},C_{22}.$
\quad Here $\uuuu{x}=(-2,-2,-1)\in C_{21}$.}
\label{Fig:4.15}
\end{figure}

\begin{lemma}\label{t4.12}
Consider the following 24 sets $C_i$, $i\in\{1,2,...,24\}$, of 
triples in $\Z^3$.
\index{$C_1,...,C_{24}$} 
\begin{eqnarray*}
C_1&=& \{(0,0,0)\}\quad (A_1^3),\\
C_2&=& \{(2,2,2)\}\quad (\HH_{1,2}),\\
C_3&=& \{(-1,0,0)\}\quad (A_2A_1),\\
C_4&=& \{(-2,0,0)\}\quad (\P^1A_1),\\
C_5&=& \{(x,0,0)\,|\, x\in\Z_{\leq -3}\}\\
&& (\textup{the reducible cases without }
A_1^3,A_2A_1,\P^1A_1),\\
C_6&=& \{(-1,0,-1)\}\quad (A_3),\\
C_7&=& \{(-1,-1,-1)\}\quad (\widehat{A}_2),\\
C_8&=& \{(-l,2,-l)\,|\, l\in\Z_{\geq 3}\textup{ odd}\},\\
C_9&=& \{(-l,2,-l)\,|\, l\in\Z_{\geq 4}\textup{ even}\},\hspace*{6cm}\\
C_{10}&=& \{(3,3,3)\}\quad (\P^2),\\
C_{11}&=& \{(x,x,x)\,|\, x\in\Z_{\geq 4}\},\\
C_{12}&=& \{(x,x,x)\,|\, x\in\Z_{\leq -2}\},\\
C_{13}&=& \{(2y,2y,2y^2)\,|\, y\in\Z_{\geq 2}\},\\
C_{14}&=& \{(x_1,x_2,\frac{1}{2}x_1x_2\,|\, 3\leq x_1<x_2,
x_1x_2\textup{ even}\},\\
C_{15}&=& \{(x_1,x_1,x_3)\,|\, 3\leq x_1<x_3<\frac{1}{2}x_1^2\}
\\
&& \cup\ 
\{(x_1,x_2,x_2)\,|\, 3\leq x_1<x_2\},\\
C_{16}&=& \{(x_1,x_2,x_3)\,|\, 3\leq x_1<x_2<x_3
<\frac{1}{2}x_1x_2\}\\
&&\cup\ 
\{(x_1,x_2,x_3)\,|\, 3\leq x_2<x_1<x_3<\frac{1}{2}x_1x_2\},\\
C_{17}&=& \{(x,x,0)\,|\, x\in\Z_{\leq -2}\},\\
C_{18}&=& \{(x_1,x_2,0)\,|\, x_1<x_2\leq -2\},\\
C_{19}&=& \{(x,-1,0)\,|\, x\in\Z_{\leq -2}\},\\
C_{20}&=& \{(x,-1,-1)\,|\, x\in\Z_{\leq -2}\},\\
C_{21}&=& \{(x,x,-1)\,|\, x\in\Z_{\leq -2}\},\\
C_{22}&=& \{(x_1,x_2,-1)\,|\, x_1<x_2\leq -2\}\\
&&\cup\ \{(x_1,x_2,-1)\,|\, x_2<x_1\leq -2\},\\
C_{23}&=& \{(x_1,x_1,x_3)\,|\, x_1<x_3\leq -2\}\\
&&\cup\ \{(x_1,x_2,x_2)\,|\, x_1<x_2\leq -2\},\\ 
C_{24}&=& \{(x_1,x_2,x_3)\,|\, x_1<x_2<x_3\leq -2\}\\
&&\cup\ \{(x_1,x_2,x_3)\,|\, x_2<x_1<x_3\leq -2\}.
\end{eqnarray*}

(a) Each triple in $\bigcup_{i=1}^{24}C_i$ is a local minimum.
All $\uuuu{x}$ in one set $C_i$ have the 
value in the following table 
or satisfy the inequality in the following table,
\begin{eqnarray*}
\begin{array}{l|l}
\rho & i\textup{ with }r(\uuuu{x})=\rho\textup{ for }
\uuuu{x}\in C_i\\ \hline 
0 & 1,\ 10\\
1 & 3 \\
2 & 6 \\
4 & 2,\ 4,\ 7,\ 8,\ 9\\
<0 & 11,\ 13,\ 14,\ 15,\ 16 \\
>4 & 5,\ 12,\ 17,\ 18,\ 19,\ 20,\ 21,\ 22,\ 23,\ 24
\end{array}
\end{eqnarray*}

(b) The following table makes statements about the 
$\langle\gamma\rangle$ orbits and the 
$\langle\gamma,\gamma_2\rangle$ orbits of
$v_0=\uuuu{x}/\{\pm 1\}^3$ with 
$\uuuu{x}\in \bigcup_{i=1}^{24}C_i$,
\begin{eqnarray*}
\begin{array}{l|l}
i\in\{1,2,7,10,11,12\} & \langle\gamma\rangle(v_0)
\textup{ is }\gamma_2\textup{-invariant}\\ 
 & \qquad\textup{ and has size }1\\ \hline 
i\in\{3,4,5,6,8,9,13,15,17,20,21,23\}& \langle\gamma\rangle(v_0)
\textup{ is }\gamma_2\textup{-invariant}\\
 & \qquad\textup{ and has size }3\\ \hline
i\in\{14,16,18,19,22,24\}& \langle\gamma\rangle(v_0)
\textup{ is not }\gamma_2\textup{-invariant}\\
 & \qquad\textup{ and has size }3,\\
  & \langle\gamma,\gamma_2\rangle(v_0)
  \textup{ has size }6
\end{array}
\end{eqnarray*}

(c) The set of all local minima in $\Z^3$ is the following 
disjoint union,
\begin{eqnarray*}
&&\Bigl(\dot\bigcup_{i\in\{1,...,24\}-\{14,18,19\}}
\dot\bigcup_{\uuuu{x}\in C_i}
(G^{sign}\rtimes\langle\gamma\rangle)\{\uuuu{x}\}\Bigr)\\
&\dot\cup&
\Bigl(\dot\bigcup_{i\in\{14,18,19\}}
\dot\bigcup_{\uuuu{x}\in C_i}
(G^{sign}\rtimes\langle\gamma,\gamma_2\rangle)\{\uuuu{x}\}
\Bigr).
\end{eqnarray*}
\end{lemma}

{\bf Proof:} Part (b) is trivial. The parts (a) and (c) follow
with the characterization of local minima in Lemma 
\ref{t4.4} and Theorem \ref{t4.6} (c)--(e). \hfill$\Box$

\begin{theorem}\label{t4.13}
(a) $\Z^3$ is the disjoint union 
$$\dot\bigcup_{i\in\{1,...,24\}}\dot\bigcup_{\uuuu{x}\in C_i}
(\Br_3\ltimes\{\pm 1\}^3)(\uuuu{x}).$$

(b) For $v_0=\uuuu{x}/\{\pm 1\}^3$ with 
$\uuuu{x}\in \bigcup_{i=1}^{24}C_i$, the set 
$\VV_0=\Br_3(v_0)\cap\LL_0$ of sign classes of local minima in
the $\Br_3\ltimes\{\pm 1\}^3$ orbit of $\uuuu{x}$ is as follows,
\begin{eqnarray*}
\VV_0= \langle\gamma\rangle(v_0)&\textup{ if }i\in 
\{1,...,24\}-\{14,18,19\},\\
\VV_0= \langle\gamma,\gamma_2\rangle(v_0)&\textup{ if }i\in 
\{14,18,19\}.
\end{eqnarray*}

(c) The set 
$\{\GG(\uuuu{x})\,|\, \uuuu{x}\in \bigcup_{i=1}^{24}C_i\}$
of pseudo-graphs $\GG(\uuuu{x})$ for 
$\uuuu{x}\in\bigcup_{i=1}^{24}C_i$ consists of the 14 
isomorphism classes $\GG_1,...,\GG_{14}$ in the Examples
\ref{t4.11}. All $\uuuu{x}$ in one set $C_i$ have the same
pseudo-graph. The first and second column in the following
table give for each of the 14 pseudo-graphs $\GG_j$ the set
or the sets $C_i$ with $\GG(\uuuu{x})=\GG_j$ for 
$\uuuu{x}\in C_i$. The third and fourth column in the following
table are subject of Theorem \ref{t4.16}.
\begin{eqnarray*}
\begin{array}{l|l|l|l}
 & \textup{sets} & 
(G^{phi}\rtimes \langle \gamma\rangle )_{\uuuu{x}/\{\pm 1\}^3} 
& (\Br_3)_{\uuuu{x}/\{\pm 1\}^3} \\ \hline 
\GG_1 & C_1\ (A_1^3),\ C_2\ (\HH_{1,2})  
& G^{phi}\rtimes\langle\gamma\rangle & \Br_3 \\
\GG_2 & C_3,C_4,C_5\textup{ (red. cases)} 
& \langle \varphi_1,\gamma^{-1}\varphi_3\rangle 
& \langle \sigma_1,\sigma_2^2\rangle \\
\GG_3 & C_6\ (A_3) 
& \langle \gamma\varphi_3,\varphi_2\varphi_1\varphi_3\rangle 
& \langle\sigma_1\sigma_2,\sigma_1^3\rangle \\
\GG_4 & C_7 \ (\widehat{A}_2) 
& \langle\gamma,\varphi_2\varphi_1\varphi_3\rangle
& \langle\sigma_2\sigma_1,\sigma_1^3\rangle\\
\GG_5 & C_8,\ C_9\ ((-l,2,-l))
& \langle\gamma^{-1}\varphi_1\rangle
&\langle\sigma^{mon},\sigma_1^{-1}\sigma_2^{-1}\sigma_1\rangle\\
\GG_6 & C_{10} (\P^2),\ C_{11},\ C_{12} 
& \langle\gamma\rangle
& \langle\sigma_2\sigma_1\rangle\\
\GG_7 & C_{13}\ (\textup{e.g. }(4,4,8)) 
& \langle \varphi_3\rangle 
& \langle \sigma_2\sigma_1^2\rangle \\
\GG_8 & C_{14} (\textup{e.g. }(3,4,6)) 
& \langle\id\rangle 
& \langle \sigma^{mon}\rangle \\
\GG_9 & C_{15},\ C_{16},\ C_{23},\ C_{24} 
& \langle\id\rangle 
&\langle \sigma^{mon}\rangle \\
\GG_{10} & C_{17}\ (\textup{e.g. }(-2,-2,0)) 
& \langle \gamma^{-1}\varphi_2\rangle 
& \langle \sigma^{mon},\sigma_2\rangle \\
\GG_{11} & C_{18}\ (\textup{e.g. }(-3,-2,0)) 
& \langle \gamma^{-1}\varphi_3\varphi_1\rangle
& \langle \sigma^{mon},\sigma_2^2\rangle \\
\GG_{12} & C_{19}\ (\textup{e.g. }(-2,-1,0)) 
& \langle \gamma^{-1}\varphi_3\varphi_1,
\varphi_3\varphi_2\varphi_1\rangle 
& \langle \sigma^{mon},\sigma_2^2,\sigma_2\sigma_1^3\sigma_2^{-1}
\rangle \\
\GG_{13} & C_{20}\ (\textup{e.g. }(-2,-1,-1)) 
& \langle \varphi_2\varphi_3\varphi_1,\varphi_3\varphi_2\varphi_1
\rangle 
& \langle \sigma^{mon},\sigma_2^3,\sigma_2\sigma_1^3\sigma_2^{-1}
\rangle \\
\GG_{14} & C_{21},\ C_{22}\ (\textup{e.g. }(-2,-2,-1))
& \langle \varphi_2\varphi_3\varphi_1\rangle 
& \langle \sigma^{mon},\sigma_2^3\rangle
\end{array}
\end{eqnarray*}
\end{theorem}

{\bf Proof:}
We start with part (c). It can be seen rather easily for all
$\uuuu{x}$ in one family $C_i$ simultaneously. We do not give
more details.

(b) The pseudo-graphs $G_8,G_{11}$ and $G_{12}$ are the only
of the 14 pseudo-graphs with $|\VV_0|=6$.
By inspection of them or by Lemma \ref{t4.10} (e), for them
$\VV_0=\langle\gamma,\gamma_2\rangle(v_0)$. The table in part (c)
gives the correspondence 
$G_8\leftrightarrow C_{14}$, $G_{11}\leftrightarrow C_{18}$,
$G_{12}\leftrightarrow C_{19}$. The other 11 pseudo-graphs
satisfy $|\VV_0|=1$ or $|\VV_0|=3$, so in any case
$\VV_0=\langle\gamma\rangle(v_0)$. 

(a) Part (c) of Lemma \ref{t4.12} alone shows already that $\Z^3$
is the union given. Part (b) of Theorem \ref{t4.13} adds only 
the fact that this is a disjoint union.
\hfill$\Box$

\begin{remarks}\label{t4.14}
(i) We have 14 pseudo-graphs $\GG_1,...,\GG_{14}$, but 
24 sets $C_1,...,C_{24}$ because the separation into the sets
shall be fine enough for the table in Theorem \ref{t4.13} (c)
and the tables in Lemma \ref{t4.12} (a) and (b).

(ii) In the pseudo-graphs $\GG_j$ with $|\VV_0|=3$ or $|\VV_0|=6$
one can choose another distinguished vertex $\www{v_0}\in\VV_0$
and change the set $\EE_\gamma$ to a set $\www{\EE_\gamma}$
accordingly. 
This gives a pseudo-graph $\www{\GG_j}$ which is not equal to
$\GG_j$, but closely related. 

The graphs $\GG_5$ and $\GG_{10}$ are related by such a change. 

The following table shows for each $\GG_j$ except $\GG_{10}$
the number of isomorphism classes of pseudo-graphs obtained in 
this way (including the original pseudo-graphs). 
In the cases $\GG_8$ and $\GG_{11}$ there is $|\VV_0|=6$,
but because of some symmetry of the pseudo-graph without
$\EE_\gamma$, there are only 3 related pseudo-graphs.
\begin{eqnarray*}
\begin{array}{c|cccccccccccccc|c}
\GG_i,\ i & 1 & 2 & 3 & 4 & 5 & 6 & 7 & 8 & 9 & 10 & 11 & 12 & 13 
& 14 & \sum \\ \hline 
\textup{related} & 1 & 3 & 3 & 1 & 3 & 1 & 3 & 3 & 3 & 
G_5 & 3 & 6 & 3 & 3 & 36\\
\textup{pseudo-graphs} &&&&&&&&&&&&&&&
\end{array}
\end{eqnarray*}
The total number 36 is the number of isomorphism classes
of pseudo-graphs $\GG(\uuuu{x})$ for $\uuuu{x}\in\Z^3$
a local minimum.
\end{remarks}

\section{The stabilizers of upper triangular 
$3\times 3$ matrices}
\label{s4.4}

The groups $G^{phi}\ltimes \langle\gamma\rangle$ and $\Br_3$
act on $\Z^3/\{\pm 1\}^3$.
The pseudo-graphs in the Examples \ref{t4.11} and Theorem
\ref{t4.13} offer a convenient way to determine the stabilizers
$(G^{phi}\ltimes\langle\gamma\rangle)_{v_0}$ and 
$(\Br_3)_{v_0}$ for $v_0=\uuuu{x}/\{\pm 1\}^3$ with
$\uuuu{x}\in\bigcup_{i=1}^{24}C_i$ a local minimum.

The stabilizers of $v_0$ depend only on the pseudo-graph 
$G_j$ with $G_j=\GG(\uuuu{x})$. The results are presented in 
Theorem \ref{t4.16}. The Remarks \ref{t4.15} prepare this.

\begin{remarks}\label{t4.15}
(i) First we recall some well known facts about the groups
$SL_2(\Z)$ and $PSL_2(\Z)$ and their relation to $\Br_3$.
The group $SL_2(\Z)$ is generated by the matrices
$$A_1:=\begin{pmatrix}1&-1\\0&1\end{pmatrix}
\quad\textup{and}\quad
A_2:=\begin{pmatrix}1&0\\1&1\end{pmatrix}.$$
\index{generators of $SL_2(\Z)$} 
Generating relations are
\begin{eqnarray*}
A_1A_2A_1=A_2A_1A_2\quad\textup{and}\quad
(A_2A_1)^6=E_2.
\end{eqnarray*}
The group $\Br_3$ is generated by the elementary braids
$\sigma_1$ and $\sigma_2$. The only generating relation is
$\sigma_1\sigma_2\sigma_1=\sigma_2\sigma_1\sigma_2.$
Therefore there is a surjective group homomorphism
\begin{eqnarray*}
\Br_3\to SL_2(\Z),\quad \sigma_1\mapsto A_1,
\quad\sigma_2\mapsto A_2,
\end{eqnarray*}
with kernel $\langle (\sigma_2\sigma_1)^6\rangle
=\langle (\sigma^{mon})^2\rangle$. 
It induces a surjective group homomorphism
\begin{eqnarray*}
\Br_3\to PSL_2(\Z),\quad \sigma_1\mapsto[A_1],
\quad \sigma_2\mapsto [A_2],
\end{eqnarray*}
with kernel $\langle\sigma^{mon}\rangle$ because
$(A_2A_1)^3=-E_2$. 

(ii) The action of $\Br_3\ltimes \{\pm 1\}^3$ on 
$T^{uni}_3(\Z)$ and on $\Z^3$ is fixed in the beginning
of section \ref{s4.1}. One sees that 
$\sigma^{mon}=(\sigma_2\sigma_1)^3$ acts trivially on
$T^{uni}_3(\Z)$ and $\Z^3$. This can be checked directly.
Or it can be seen as a consequence of the following two facts.
\begin{list}{}{}
\item[(1)]
The action of $\Br_3\ltimes\{\pm 1\}^3$ on 
$(\Br_3\ltimes\{\pm 1\}^3)(\uuuu{x})$ for some $\uuuu{x}\in\Z^3$
is induced by the action of $\Br_3\ltimes\{\pm 1\}^3$ on
the set $\BB^{dist}$ of distinguished bases of a triple
$(H_\Z,L,\uuuu{e})$ with $L(\uuuu{e}^t,\uuuu{e})^t=S(\uuuu{x})$
by $S((\alpha,\varepsilon)(\uuuu{x}))=
L((\alpha,\varepsilon)(\uuuu{e})^t,
(\alpha,\varepsilon)(\uuuu{e}))^t$ for 
$(\alpha,\varepsilon)\in \Br_3\ltimes\{\pm 1\}^3$.
\item[(2)]
For $(\alpha,\varepsilon)=(\sigma^{mon},(1,1,1))$
$(\alpha,\varepsilon)(\uuuu{e})=Z((\alpha,\varepsilon))
(\uuuu{e})=M(\uuuu{e})$ by Theorem \ref{t3.10},
and $L(M(\uuuu{e})^t,M(\uuuu{e}))=L(\uuuu{e}^t,\uuuu{e})$. 
\end{list}
In any case the action of $\Br_3\ltimes\{\pm 1\}^3$ on $\Z^3$
boils down to a nonlinear action of
$PSL_2(\Z)\ltimes G^{sign}$ where 
$G^{sign}=\langle\delta_1^\R,\delta_2^\R\rangle$ catches the
action of $\{\pm 1\}^3$ on $\Z^3$, see section \ref{s4.1}.

(iii)
The shape of this nonlinear action led us in Definition 
\ref{t4.1} and Theorem \ref{t4.2} to the group 
$(G^{phi}\ltimes G^{sign})\rtimes \langle\gamma\rangle
=(G^{phi}\rtimes \langle\gamma\rangle)\ltimes G^{sign}$.
In fact, $G^{phi}\rtimes\langle\gamma\rangle
\cong PSL_2(\Z)$. This can be seen as follows.

The formulas \eqref{4.1}--\eqref{4.4}
in the proof of Theorem \ref{t4.2} (c)
give lifts to $\Br_3\ltimes\{\pm 1\}^3$ of the 
generators $\varphi_1,\varphi_2,\varphi_3$ and $\gamma$ of
$G^{phi}\ltimes\langle\gamma\rangle$. 
Dropping the generators of the sign action in these lifts, 
we obtain the following lifts to $\Br_3$,
\begin{eqnarray}
\left.\begin{array}{l}
l(\gamma)=\sigma_2\sigma_1,
\quad l(\gamma^{-1})=\sigma_1^{-1}\sigma_2^{-1},\\
l(\varphi_1)=l(\gamma)^{-1}\sigma_2^{-1}
=\sigma_1^{-1}\sigma_2^{-2},\\
l(\varphi_2)=l(\gamma)\sigma_2=\sigma_2\sigma_1\sigma_2
=\sigma_1\sigma_2\sigma_1,\\
l(\varphi_3)=l(\gamma)\sigma_1=\sigma_2\sigma_1^2.
\end{array}\right\}\label{4.11}
\end{eqnarray}
The equality of groups in Theorem \ref{t4.2} (c) boils down 
after dropping the sign action to an equality of groups
\begin{eqnarray}\label{4.12}
\langle \sigma_1^\R,\sigma_2^\R\rangle\cong 
G^{phi}\rtimes\langle\gamma\rangle.
\end{eqnarray}
As $(\sigma_2^\R\sigma_1^\R)^3=\id$, we obtain a surjective
group homomorphism 
$PSL_2(\Z)\to G^{phi}\rtimes\langle\gamma\rangle$ with 
$[A_i]\mapsto \sigma_i^\R$. 
The subgroup 
$\langle [A_1]^{-1}[A_2]^{-2},[A_2][A_1][A_2],
[A_2][A_1]^2\rangle$ of $PSL_2(\Z)$ is mapped to $G^{phi}$.
One easily calculates that this subgroup is 
the free Coxeter group
with three generators which was considered in Remark 
\ref{t6.12} (iv) and which has index three in $PSL_2(\Z)$.
As $G^{phi}$ is also a free Coxeter group with three generators
and has index 3 in $G^{phi}\rtimes\langle\gamma\rangle$,
the map $PSL_2(\Z)\to G^{phi}\rtimes\langle\gamma\rangle$
is a group isomorphism.

(iv) For use in the proof of Theorem \ref{t4.16} we recall
the formulas \eqref{4.6}
\begin{eqnarray}
\left.\begin{array}{ccccccccc}
\gamma\varphi_1&=&\varphi_2\gamma,& 
\gamma\varphi_2&=&\varphi_3\gamma,& 
\gamma\varphi_3&=&\varphi_1\gamma,\\
\varphi_1\gamma^{-1}&=&\gamma^{-1}\varphi_2,&
\varphi_2\gamma^{-1}&=&\gamma^{-1}\varphi_3,&
\varphi_3\gamma^{-1}&=&\gamma^{-1}\varphi_1,
\end{array}\right\}\label{4.13}
\end{eqnarray}
from the proof of Theorem \ref{t4.2}.

(v) The relation 
$\sigma_1\sigma_2\sigma_1=\sigma_2\sigma_1\sigma_2$
is equivalent to each of the two relations
\begin{eqnarray*}
\sigma_1\sigma_2\sigma_1^{-1}=\sigma_2^{-1}\sigma_1\sigma_2
\quad\textup{and}\quad
\sigma_1^{-1}\sigma_2\sigma_1=\sigma_2\sigma_1\sigma_2^{-1}
\end{eqnarray*}
and induces for any $m\in\Z$ the relations
\begin{eqnarray}\label{4.14}
\sigma_1\sigma_2^m\sigma_1^{-1}=\sigma_2^{-1}\sigma_1^m\sigma_2
\quad\textup{and}\quad
\sigma_1^{-1}\sigma_2^m\sigma_1=\sigma_2\sigma_1^m\sigma_2^{-1}.
\end{eqnarray}
Also this will be useful in the proof of Theorem
\ref{t4.16}.
\end{remarks}

\begin{theorem}\label{t4.16}
Consider $v_0:=\uuuu{x}/\{\pm 1\}^3$ with
$\uuuu{x}\in\bigcup_{i=1}^{24}C_i\subset\Z^3$
a local minimum, and consider the pseudo-graph $\GG_j$ with
$\GG_j=\GG(\uuuu{x})$. The entries in the third and fourth
column in the table in Theorem \ref{t4.13}, which are in the
line of $\GG_j$, give the stabilizers 
$(G^{phi}\rtimes\langle\gamma\rangle)_{v_0}$ and
$(\Br_3)_{v_0}$ of $v_0$. 
\end{theorem}

{\bf Proof:}
First we treat the stabilizer 
$(G^{phi}\rtimes\langle\gamma\rangle)_{v_0}$. 

The {\it total set} $|\GG_j|$ of the pseudo-graph $\GG_j$
means the union of vertices and edges in an embedding of
the pseudo-graph in the real plane $\R^2$ 
as in the figures in the Examples \ref{t4.11}. 
The fundamental group $\pi_1(|\GG_j|,v_0)$ is a free group
with 0, 1, 3, 4, 5 or 6 generators. The number of 
generators is the number of compact components in
$\R^2-|\GG_j|$. A generator is the class of a closed path 
which starts and ends at $v_0$ and turns once around one of 
these compact components. Any such closed path induces a word
in $\varphi_1,\varphi_2,\varphi_3,\gamma$ and $\gamma^{-1}$.
This word gives an element of 
$(G^{phi}\rtimes\langle\gamma\rangle)_{v_0}$. We obtain a 
group homomorphism
\begin{eqnarray*}
\pi_1(|\GG_j|,v_0)\to 
(G^{phi}\rtimes\langle\gamma\rangle)_{v_0}.
\end{eqnarray*}
It is surjective because any element of 
$(G^{phi}\rtimes\langle\gamma\rangle)_{v_0}$ can be written
as $\psi\gamma^\xi$ with $\psi\in G^{phi}$ and 
$\xi\in\{0,\pm 1\}$. The element $\psi\gamma^\xi$ leads
to and comes from a closed path in $|\GG_j|$ which starts
and ends at $v_0$. 

In fact, this shows that we could restrict to closed paths
which run never or only once at the beginning through an edge
in $\EE_\gamma$. But we will not use this fact.

The following list gives for each of the 14 pseudo-graphs
$\GG_1,...,\GG_{14}$ in the first line one word in
$\varphi_1,\varphi_2,\varphi_3,\gamma$ and $\gamma^{-1}$ for
each compact component of $\R^2-|\GG_j|$. In the following
lines the relations \eqref{4.13} are used to show that all
these words are generated in 
$G^{phi}\rtimes \langle\gamma\rangle$ by the generators in the
third column in the table in Theorem \ref{t4.13}.
The generators are underlined.

\begin{list}{}{}
\item[$\GG_1:$]
$\uuuu{\varphi_1}$, $\uuuu{\varphi_2}$, 
$\uuuu{\varphi_3}$, $\uuuu{\gamma}$. 
\item[$\GG_2:$]
$\uuuu{\varphi_1}$, $\varphi_3\gamma$, $\gamma\varphi_2$, 
$\gamma^{-1}\varphi_2\gamma$, $\gamma\varphi_3\gamma^{-1}$, 
$\gamma\varphi_1\gamma$. 
\begin{eqnarray*}
\gamma\varphi_1\gamma&=&\gamma^2\varphi_3
=\uuuu{\gamma^{-1}\varphi_3},\\
\gamma\varphi_2&=& \varphi_3\gamma 
=(\gamma^{-1}\varphi_3)^{-1},\\
\gamma^{-1}\varphi_2\gamma&=&\varphi_1,\\
\gamma\varphi_3\gamma^{-1}&=&\varphi_1.
\end{eqnarray*}
\item[$\GG_3:$]
$\varphi_1\gamma$, $\varphi_2\varphi_3\gamma$, 
$\gamma\varphi_2\gamma$, $\gamma\varphi_1\varphi_2$, 
$\uuuu{\gamma\varphi_3}$.
\begin{eqnarray*}
\varphi_1\gamma&=& \gamma\varphi_3,\\
\varphi_2\varphi_3\gamma&=& \varphi_2\gamma\varphi_2
=\gamma\varphi_1\varphi_2 = (\gamma\varphi_3)
(\uuuu{\varphi_2\varphi_1\varphi_3})^{-1},\\
\gamma\varphi_2\gamma&=& \gamma^2\varphi_1
=\gamma^{-1}\varphi_1 =\varphi_3\gamma^{-1}
=(\gamma\varphi_3)^{-1}.
\end{eqnarray*}
\item[$\GG_4:$]
$\uuuu{\gamma}$, $\uuuu{\varphi_2\varphi_1\varphi_3}$, 
$\varphi_3\varphi_2\varphi_1$, 
$\varphi_1\varphi_3\varphi_2$.
\begin{eqnarray*}
\gamma\varphi_2\varphi_1\varphi_3\gamma^{-1}
&=&\gamma\varphi_2\varphi_1\gamma^{-1}\varphi_1
=\gamma\varphi_2\gamma^{-1}\varphi_2\varphi_1
=\varphi_3\varphi_2\varphi_1,\\
\gamma^{-1}\varphi_2\varphi_1\varphi_3\gamma
&=&\gamma^{-1}\varphi_2\varphi_1\gamma\varphi_2
=\gamma^{-1}\varphi_2\gamma\varphi_3\varphi_2
=\varphi_1\varphi_3\varphi_2.
\end{eqnarray*}
\item[$\GG_5:$]
$\uuuu{\gamma^{-1}\varphi_1}$, $\gamma\varphi_2\gamma$, 
$\gamma\varphi_3$. 
\begin{eqnarray*}
\gamma\varphi_2\gamma&=& \gamma^2\varphi_1
=\gamma^{-1}\varphi_1,\\
\gamma\varphi_3&=& \varphi_1\gamma 
=(\gamma^{-1}\varphi_1)^{-1}.
\end{eqnarray*}
\item[$\GG_6:$]
$\uuuu{\gamma}$. 
\item[$\GG_7:$]
$\uuuu{\varphi_3}$, $\gamma^{-1}\varphi_1\gamma$, 
$\gamma\varphi_2\gamma^{-1}$.
\begin{eqnarray*}
\gamma^{-1}\varphi_1\gamma&=& \varphi_3,\\
\gamma\varphi_2\gamma^{-1}&=& \varphi_3.
\end{eqnarray*}
\item[$\GG_8:$]
$\uuuu{\id}$. 
\item[$\GG_9:$]
$\uuuu{\id}$. 
\item[$\GG_{10}:$]
$\uuuu{\gamma^{-1}\varphi_2}$, 
$\gamma\varphi_3\gamma$, $\gamma\varphi_1$.
\begin{eqnarray*}
\gamma\varphi_3\gamma&=&\gamma^2\varphi_2
=\gamma^{-1}\varphi_2,\\
\gamma\varphi_1&=&\varphi_2\gamma
=(\gamma^{-1}\varphi_2)^{-1}.
\end{eqnarray*}
\item[$\GG_{11}:$]
$\uuuu{\gamma^{-1}\varphi_3\varphi_1}$, 
$\gamma\varphi_1\varphi_2\gamma$,
$\varphi_2\varphi_3\gamma^{-1}$.
\begin{eqnarray*}
\gamma\varphi_1\varphi_2\gamma
&=& \gamma\varphi_1\gamma\varphi_1
=\gamma^2\varphi_3\varphi_1=\gamma^{-1}\varphi_3\varphi_1,\\
\varphi_2\varphi_3\gamma^{-1}
&=& \varphi_2\gamma^{-1}\varphi_1=\gamma^{-1}\varphi_3\varphi_1.
\end{eqnarray*}
\item[$\GG_{12}:$]
$\uuuu{\varphi_3\varphi_2\varphi_1}$, 
$\uuuu{\gamma^{-1}\varphi_3\varphi_1}$,
$\gamma^{-1}\varphi_1\varphi_3\varphi_2\gamma$,
$\gamma\varphi_1\varphi_2\gamma$, 
$\varphi_2\varphi_3\gamma^{-1}$, 
$\gamma\varphi_2\varphi_1\varphi_3\gamma^{-1}$.
\begin{eqnarray*}
\gamma^{-1}\varphi_1\varphi_3\varphi_2\gamma
&=& \gamma^{-1}\varphi_1\varphi_3\gamma\varphi_1
=\gamma^{-1}\varphi_1\gamma\varphi_2\varphi_1
=\varphi_3\varphi_2\varphi_1,\\
\gamma\varphi_1\varphi_2\gamma
&=& \gamma\varphi_1\gamma\varphi_1
=\gamma^2\varphi_3\varphi_1=\gamma^{-1}\varphi_3\varphi_1,\\
\varphi_2\varphi_3\gamma^{-1}
&=& \varphi_2\gamma^{-1}\varphi_1
=\gamma^{-1}\varphi_3\varphi_1,\\
\gamma\varphi_2\varphi_1\varphi_3\gamma^{-1}
&=& \gamma\varphi_2\varphi_1\gamma^{-1}\varphi_1
=\gamma\varphi_2\gamma^{-1}\varphi_2\varphi_1
=\varphi_3\varphi_2\varphi_1.
\end{eqnarray*}
\item[$\GG_{13}:$]
$\uuuu{\varphi_2\varphi_3\varphi_1}$, 
$\uuuu{\varphi_3\varphi_2\varphi_1}$, 
$\gamma^{-1}\varphi_3\varphi_1\varphi_2\gamma$,
$\gamma^{-1}\varphi_1\varphi_3\varphi_2\gamma$,
$\gamma\varphi_1\varphi_2\varphi_3\gamma^{-1}$,
$\gamma\varphi_2\varphi_1\varphi_3\gamma^{-1}$.
\begin{eqnarray*}
\gamma^{-1}\varphi_3\varphi_1\varphi_2\gamma
&=& \gamma^{-1}\varphi_3\varphi_1\gamma\varphi_1
=\gamma^{-1}\varphi_3\gamma\varphi_3\varphi_1
=\varphi_2\varphi_3\varphi_1,\\
\gamma^{-1}\varphi_1\varphi_3\varphi_2\gamma
&=& \gamma^{-1}\varphi_1\varphi_3\gamma\varphi_1
=\gamma^{-1}\varphi_1\gamma\varphi_2\varphi_1
=\varphi_3\varphi_2\varphi_1,\\
\gamma\varphi_1\varphi_2\varphi_3\gamma^{-1}
&=& \gamma\varphi_1\varphi_2\gamma^{-1}\varphi_1
=\gamma\varphi_1\gamma^{-1}\varphi_3\varphi_1
=\varphi_2\varphi_3\varphi_1,\\
\gamma\varphi_2\varphi_1\varphi_3\gamma^{-1}
&=& \gamma\varphi_2\varphi_1\gamma^{-1}\varphi_1
=\gamma\varphi_2\gamma^{-1}\varphi_2\varphi_1
=\varphi_3\varphi_2\varphi_1.
\end{eqnarray*}
\item[$\GG_{14}:$]
$\uuuu{\varphi_2\varphi_3\varphi_1}$, 
$\gamma^{-1}\varphi_3\varphi_1\varphi_2\gamma$,
$\gamma\varphi_1\varphi_2\varphi_3\gamma^{-1}$. 
\begin{eqnarray*}
\gamma^{-1}\varphi_3\varphi_1\varphi_2\gamma
&=& \varphi_2\varphi_3\varphi_1\quad 
(\textup{see }\GG_{13}),\\
\gamma\varphi_1\varphi_2\varphi_3\gamma^{-1}
&=& \varphi_2\varphi_3\varphi_1 \quad
(\textup{see }\GG_{13}).
\end{eqnarray*}
\end{list}
Therefore the stabilizer 
$(G^{phi}\rtimes\langle\gamma\rangle)_{v_0}$
is as claimed in the third column of the table in Theorem
\ref{t4.13}. 

Now we treat the stabilizer $(\Br_3)_{v_0}$. 
It is the preimage in $\Br_3$ of
$(G^{phi}\rtimes\langle\gamma\rangle)_{v_0}$
under the surjective group homomorphism
$\Br_3\to G^{phi}\rtimes \langle\gamma\rangle$
with kernel $\langle \sigma^{mon}\rangle$. 
So if $(G^{phi}\rtimes \langle\gamma\rangle)_{v_0}
=\langle g_1,...,g_m\rangle$ and $h_1,...,h_m$ are
any lifts to $\Br_3$ of $g_1,...,g_m$ then 
$(\Br_3)_{v_0}=\langle \sigma^{mon},h_1,...,h_m\rangle$. 

For any word in $\varphi_1,\varphi_2,\varphi_3,\gamma$
and $\gamma^{-1}$ we use the lifts in \eqref{4.11} 
to construct a lift of this word.

The following list gives for each of the 14 
pseudo-graphs $\GG_1,...,\GG_{14}$ for each generator of
$(G^{phi}\times \langle\gamma\rangle)_{v_0}$ in the 
third column of the table in Theorem \ref{t4.13}
this lift and rewrites it using the relations \eqref{4.14}.
The generators in the fourth column of the table in
Theorem \ref{t4.13} are underlined. 

\begin{eqnarray*}
\GG_1:& G^{phi}&\rightsquigarrow \Br_3,\\
\GG_2:& \varphi_1 &\rightsquigarrow 
\sigma_1^{-1}\sigma_2^{-2}
=\sigma_2(\sigma_2^{-1}\sigma_1^{-1}\sigma_2^{-1})\sigma_2^{-1}
=\sigma_2(\sigma_1^{-1}\sigma_2^{-1}\sigma_1^{-1})\sigma_2^{-1}
\\
&&=(\sigma_2^2\sigma_1)(\sigma_1^{-1}\sigma_2^{-1})^3
=\uuuu{\sigma_2^2}\sigma_1(\sigma^{mon})^{-1} ,\\
      & \gamma^{-1}\varphi_3 &\rightsquigarrow 
      (\sigma_1^{-1}\sigma_2^{-1})(\sigma_2\sigma_1^2)
      =\uuuu{\sigma_1} ,\\
\GG_3:& \gamma\varphi_3 &\rightsquigarrow 
(\sigma_2\sigma_1)(\sigma_2\sigma_1^2)
=(\sigma_2\sigma_1\sigma_2)\sigma_1^2
=(\uuuu{\sigma_1\sigma_2})\sigma_1^3 ,\\
      & \varphi_2\varphi_1\varphi_3 &\rightsquigarrow 
(\sigma_1\sigma_2\sigma_1)(\sigma_1^{-1}\sigma_2^{-2})
(\sigma_2\sigma_1^2)=\uuuu{\sigma_1^3} ,\\
\GG_4:& \gamma &\rightsquigarrow 
\uuuu{\sigma_2\sigma_1} ,\\
      & \varphi_2\varphi_1\varphi_3 &\rightsquigarrow 
\uuuu{\sigma_1^3} ,\\
\GG_5:& \gamma^{-1}\varphi_1 &\rightsquigarrow 
(\sigma_1^{-1}\sigma_2^{-1})(\sigma_1^{-1}\sigma_2^{-2})
=(\sigma_1^{-1}\sigma_2^{-1})^3\sigma_2\sigma_1\sigma_2^{-1}\\
&&\stackrel{\eqref{4.14}}{=}(\sigma^{mon})^{-1}(\uuuu{\sigma_1^{-1}\sigma_2^{-1}
\sigma_1})^{-1},\\
\GG_6:& \gamma &\rightsquigarrow 
\uuuu{\sigma_2\sigma_1},\\
\GG_7:& \varphi_3 &\rightsquigarrow 
\uuuu{\sigma_2\sigma_1^2},
\end{eqnarray*}
\begin{eqnarray*}
\GG_8:& \id &\rightsquigarrow \uuuu{\id} ,\\
\GG_9:& \id &\rightsquigarrow \uuuu{\id} ,\\
\GG_{10}:& \gamma^{-1}\varphi_2 &\rightsquigarrow 
(\sigma_1^{-1}\sigma_2^{-1})(\sigma_2\sigma_1\sigma_2)
=\uuuu{\sigma_2} ,\\
\GG_{11}:& \gamma^{-1}\varphi_3\varphi_1 &\rightsquigarrow 
(\sigma_1^{-1}\sigma_2^{-1})(\sigma_2\sigma_1^2)
(\sigma_1^{-1}\sigma_2^{-2})
=\sigma_2^{-2}=(\uuuu{\sigma_2^2})^{-1} ,\\
\GG_{12}:& \gamma^{-1}\varphi_3\varphi_1 &\rightsquigarrow 
(\uuuu{\sigma_2^2})^{-1} ,\\
         & \varphi_3\varphi_2\varphi_1 &\rightsquigarrow 
(\sigma_2\sigma_1^2)(\sigma_1\sigma_2\sigma_1)
(\sigma_1^{-1}\sigma_2^{-2})
=\uuuu{\sigma_2\sigma_1^3\sigma_2^{-1}} ,\\
\GG_{13}:& \varphi_2\varphi_3\varphi_1 &\rightsquigarrow 
(\sigma_1\sigma_2\sigma_1)(\sigma_2\sigma_1^2)
(\sigma_1^{-1}\sigma_2^{-2})
=(\sigma_1\sigma_2)^3\sigma_2^{-3}\\
&&=\sigma^{mon}(\uuuu{\sigma_2^3})^{-1} ,\\
         & \varphi_3\varphi_2\varphi_1 &\rightsquigarrow 
\uuuu{\sigma_2\sigma_1^3\sigma_2^{-1}} ,\\
\GG_{14}:& \varphi_2\varphi_3\varphi_1 &\rightsquigarrow 
\sigma^{mon}(\uuuu{\sigma_2^3})^{-1}. 
\end{eqnarray*}
Observe
$$(\sigma_2\sigma_1)^3=\sigma^{mon},\quad 
(\sigma_2\sigma_1^2)^2=\sigma^{mon}.$$
Therefore the stabilizer $(\Br_3)_{v_0}$ is as claimed in the
fourth column of the table in Theorem \ref{t4.13}.
\hfill$\Box$

\section{A global sign change, relevant for the odd case}
\label{s4.5}

\begin{remarks}\label{t4.17}
(i) For $S\in T^{uni}_n(\Z)$ consider a unimodular bilinear
lattice $(H_\Z,L)$ with a triangular basis
$\uuuu{e}=(e_1,...,e_n)$ with $L(\uuuu{e}^t,\uuuu{e})^t=S$.
Consider also the matrix $\www{S}\in T^{uni}_n(\Z)$ with
$\www{S}_{ij}=-S_{ij}$ for $i<j$. 

On the same lattice $H_\Z$ and the same basis $\uuuu{e}$
we define a second unimodular bilinear form $\www{L}$
by $\www{L}(\uuuu{e}^t,\uuuu{e})^t=\www{S}$
and denote all objects associated to $(H_\Z,\www{L},\uuuu{e})$
with a tilde, 
$\www{I}^{(k)}$, $\www{\BB}^{tri}$, $\www{\Gamma}^{(k)}$, 
$\www{\Delta}^{(k)}$.

Most of them differ a lot from the objects associated to 
$(H_\Z,L,\uuuu{e})$. But the odd intersection forms differ
only by the sign. Therefore the monodromies $M$ and $\www{M}$
are different (in general), but the odd monodromy groups and
the sets of odd vanishing cycles coincide:
\begin{eqnarray*}
\www{I}^{(1)}&=& -I^{(1)},\\
\textup{so }\www{s}_{e_i}^{(1)}&=&(s_{e_i}^{(1)})^{-1},\\
\www{M}&=& \www{s}_{e_1}^{(1)}\circ ...\circ 
\www{s}_{e_n}^{(1)}
\stackrel{\textup{in general}}{\neq} M
=s_{e_1}^{(1)}\circ ...\circ s_{e_n}^{(1)},\\
\textup{but }\www{\Gamma}^{(1)}
&=&\langle \www{s}_{e_1}^{(1)},...,\www{s}_{e_n}^{(1)}\rangle
=\langle s_{e_1}^{(1)},...,s_{e_n}^{(1)}\rangle
=\Gamma^{(1)},\\
\www{\Delta}^{(1)}&=& \www{\Gamma}^{(1)}\{\pm e_1,...,\pm e_n\}
={\Gamma}^{(1)}\{\pm e_1,...,\pm e_n\}=\Delta^{(1)}.
\end{eqnarray*}
Because of $\www{\Gamma}^{(1)}=\Gamma^{(1)}$ and 
$\www{\Delta}^{(1)}=\Delta^{(1)}$ the global sign change
from $S$ to $\www{S}$ is interesting.

(ii) In this section we will study the action on
$T^{uni}_3(\Z)$ of the extension of the action of
$\Br_3\ltimes\{\pm 1\}^3$ by this global sign change.
\index{global sign change}
Define \index{$\delta^\R:\R^3\to\R^3$} 
\begin{eqnarray*}
\delta^{\R}:\R^3\to\R^3,\ (x_1,x_2,x_3)\mapsto (-x_1,-x_2,-x_3)
\end{eqnarray*}
and
\begin{eqnarray*}
\www{G}^{sign}:=\langle \delta_1^{\R},\delta_2^{\R},
\delta^{\R}\rangle\cong\{\pm 1\}^3.
\end{eqnarray*}
It is easy to see that the double semidirect product
$(G^{phi}\ltimes G^{sign})\rtimes \langle\gamma\rangle$
extends to the double semidirect product
$(G^{phi}\ltimes \www{G}^{sign})\rtimes\langle\gamma\rangle$.
\end{remarks}

\begin{lemma}\label{t4.18}
Each $(G^{phi}\ltimes\www{G}^{sign})\rtimes 
\langle\gamma\rangle$ orbit in $\Z^3$ contains at least
one local minimum of one of the following types:
\begin{list}{}{}
\item[(a)]
$\uuuu{x}\in\Z^3_{\geq 3}$ with $2x_i\leq x_jx_k$ for 
$\{i,j,k\}=\{1,2,3\}$.
\item[(b)]
$(-l,2,-l)$ for some $l\in\Z_{\geq 2}$.
\item[(c)]
$(x_1,x_2,0)$ for some $x_1,x_2\in\Z_{\geq 0}$ with
$x_1\geq x_2$.
\end{list}
\end{lemma}

{\bf Proof:} In (c) we can restrict to $x_1\geq x_2$ because
$\delta^\R_3\gamma\varphi_1(x_1,x_2,0)=(x_2,x_1,0)$.

Each $(G^{phi}\ltimes \www{G}^{sign})\rtimes
\langle\gamma\rangle$ orbit consists of one or several
$(G^{phi}\ltimes G^{sign})\rtimes\langle\gamma\rangle$
orbits and thus contains local minima.

Suppose that $\uuuu{x}\in\Z^3$ is such a local minimum
and is not obtained with $G^{sign}$ from a local minimum
in (a), (b) or (c).
Then either $\uuuu{x}$ is a local minimum associated to
$S(\whh{A}_2)$, so $\uuuu{x}\in \{(-1,-1,-1),(1,1,-1),
(1,-1,1),(-1,1,1)\}$ or $x_1x_2x_3<0$ and
$r(\uuuu{x})>4$.

In the first case $\delta^{\R}(-1,-1,-1)=(1,1,1)
=\varphi_3(1,1,0)$, so the orbit contains a local minimum in (c).

In the second case $r(\delta^{\R}(\uuuu{x}))=r(-\uuuu{x})
=r(\uuuu{x})+2x_1x_2x_3<r(\uuuu{x})$. 
We consider a local minimum $\uuuu{x}^{(1)}$ in the 
$\Br_3\ltimes\{\pm 1\}^3$ orbit of $-\uuuu{x}$.
If it is not obtained with $G^{sign}$ from a local minimum
in (a), (b) or (c), then $\uuuu{x}^{(1)}$ is a local minimum
associated to $S(\whh{A}_2)$ or 
$x_1^{(1)}x_2^{(1)}x_3^{(1)}<0$ and $r(\uuuu{x}^{(1)})>4$.

We repeat this procedure until we arrive at a local minimum
obtained with $G^{sign}$ from one in (a), (b) or (c).
This stops after finitely many steps because
$r(\uuuu{x}^{(1)})=r(-\uuuu{x})<r(\uuuu{x})$.
\hfill$\Box$

\begin{remark}\label{t4.19}
Corollary \ref{t6.23} will say that the 
$(G^{phi}\ltimes \www{G}^{sign})\rtimes \langle\gamma\rangle$
orbits of the local minima in the parts (b) and (c) of
Lemma \ref{t4.18} are pairwise different and also different
from the orbits of the local minima in part (a). 
\end{remark}

\begin{examples}\label{t4.20}
Given an element $\uuuu{x}\in\Z^3$, 
it is not obvious which local minimum of a type in (a), (b) 
or (c) is contained in the $(G^{phi}\ltimes \www{G}^{sign})
\rtimes\langle\gamma\rangle$ orbit of $\uuuu{x}$. 
We give four families of examples.
An arrow $\mapsto$ between two elements of $\Z^3$ means that
these two elements are in the same orbit. 

(i) Start with $(x_1,x_2,-1)$ with $x_1\geq x_2>0$.
\begin{eqnarray*}
(x_1,x_2,-1)&\stackrel{\www{G}^{sign}}{\mapsto}&
(x_1,x_2,1)\stackrel{\varphi_1}{\mapsto}
(x_2-x_1,1,x_2)\stackrel{\www{G}^{sign}}{\mapsto}
(x_1-x_2,1,x_2)\\
&\mapsto & ... \mapsto (\gcd(x_1,x_2),1,0).
\end{eqnarray*}

(ii) Start with $(x_1,x_2,-2)$ with $x_1\geq x_2\geq 2$.
\begin{eqnarray*}
(x_1,x_2,-2)&\stackrel{\www{G}^{sign}}{\mapsto}&
(x_1,x_2,2)\stackrel{\varphi_1}{\mapsto}
(2x_2-x_1,2,x_2)\mapsto ...\\
& \mapsto&  \left\{\begin{array}{l}
(\gcd(x_1,x_2),2,0) \quad\textup{ if }x_1\textup{ and }x_2\\
\hspace*{2cm}\textup{contain different powers of }2,\\
(-\gcd(x_1,x_2),2,-\gcd(x_1,x_2))\quad \textup{ if }
x_1\textup{ and }x_2\\
\hspace*{2cm}\textup{contain the same power of }2
\end{array}\right.
\end{eqnarray*}
In order to understand the case discussion, observe
that $2x_2-x_1$ and $x_2$ contain the same power of $2$
if and only if $x_1$ and $x_2$ contain the same power of $2$.
Furthermore $0$ is divisible by an arbitrarily large power of $2$. 

(iii) The examples in part (ii) lead in the special case
$\gcd(x_1,x_2)=1$ to the following, 
\begin{eqnarray*}
(x_1,x_2,-2)\mapsto ... \mapsto 
&\mapsto& (1,2,0) \textup{ or }(-1,2,-1) \\
&\mapsto& (2,1,0)\textup{ or }(1,1,0),
\end{eqnarray*}
again depending on whether $x_1$ and $x_2$ contain different
powers of $2$ or the same power of $2$. 

(iv) Start with $(-3,-3,-l)$ for some $l\in\Z_{\geq 2}$.
\begin{eqnarray*}
(-3,-3,-l)&\stackrel{\delta^{\R}}{\mapsto}& 
(3,3,l)\mapsto (3,3,\pm (l-9))\\
&\mapsto& (3,3,\www{l})\textup{ for some }
\www{l}\in\{0,1,2,3,4\}\\
&\mapsto& (3,3,0),\ (3,1,0),\ (-3,2,-3),\ 
(3,3,3)\textup{ or }(4,3,3).
\end{eqnarray*}
$(3,3,0)$ and $(3,1,0)$ are in part (c),
$(-3,2,-3)$ is in part (b),
and $(3,3,3)$ and $(4,3,3)$ are in part (a) of Lemma
\ref{t4.18}. 
\end{examples}

\chapter{Automorphism groups}\label{s5}
\setcounter{equation}{0}
\setcounter{figure}{0}

A unimodular bilinear lattice $(H_\Z,L)$ comes equipped with
four automorphism groups $G_\Z$, $G_\Z^{(0)}$, $G_\Z^{(1)}$
and $G_\Z^M$ of $H_\Z$, which all respect the monodromy $M$
and possibly some bilinear form. 
They are the subject of this chapter. 

Section \ref{s5.1} gives two basic observations which serve
as general tools to control these groups under reasonable
conditions. Another very useful tool is Theorem \ref{t3.26} (c),
which gives in favourable situations an $n$-th root of
$(-1)^{k+1}M$ in $G_\Z$. Section \ref{s5.1} treats also the
cases $A_1^n$. 

Section \ref{s5.2} takes care of the rank 2 cases.
It makes use of some statements on quadratic units in 
Lemma \ref{tc.1} (a) in Appendix \ref{sc}. 

All further sections \ref{s5.3}--\ref{s5.7} are devoted
to the rank 3 cases. Section \ref{s5.3} discusses the 
setting and the basic data. It also introduces the special
automorphism $Q\in\Aut(H_\Q,L)$ which is $\id$ on 
$\ker(M-\id)$ and $-\id$ on $\ker(M^2-(2-r)M+\id)$ 
and determines the (rather few) cases where $Q$ is in $G_\Z$. 

The treatment of the reducible rank 3 cases in section
\ref{s5.4} builds on the rank 2 cases and is easy.

The irreducible rank 3 cases with all eigenvalues in $S^1$
form the four single cases $A_3$, $\whh{A}_2$, $\P^2$,
$\HH_{1,2}$ and the series $S(-l,2,-l)$ with $l\in\Z_{\geq 3}$. 
Here in section \ref{s5.5} third roots of the monodromy 
and in the series $S(-l,2,-l)$ even higher roots 
of the monodromy turn up.

The sections \ref{s5.6} and \ref{s5.7} treat all irreducible
rank 3 cases with eigenvalues not all in $S^1$
(1 is always one eigenvalue). Section \ref{s5.6}
takes care of those families of cases where 
$G_\Z\supsetneqq \{\pm M^m\,|\, m\in\Z\}$,
section \ref{s5.7} of all others.

Section \ref{s5.6} is rather long. Here again roots 
of the monodromy turn up, and statements on quadratic units
in Lemma \ref{tc.1} are used. 

Section \ref{s5.7} is very long. The main result
$G_\Z= \{\pm M^m\,|\, m\in\Z\}$ in these cases 
requires an extensive case discussion. 

This chapter determines $G_\Z$ in all cases with rank $n\leq 3$.
An application is a proof of Theorem \ref{t3.28}, 
which says that in almost all cases with rank $n\leq 3$ the map 
$Z:(\Br_n\ltimes\{\pm 1\}^n)_S\to G_\Z$ is surjective, the
exception being four cases in section \ref{s5.6}. 

\section{Basic observations}
\label{s5.1}

Given a bilinear lattice $(H_\Z,L)$, the most important
of the four automorphisms groups 
$G_\Z^M,G_\Z^{(0)},G_\Z^{(1)}$ and $G_\Z$ in Definition
\ref{t2.3} (b) (iv) and Lemma \ref{t2.6} (a) (iii) is the
smallest group $G_\Z$. But the key to it is often the largest
group $G_\Z^M$. We collect some elementary observations 
on these groups.

\begin{lemma}\label{t5.1}
Let $H_\Z$ be a $\Z$-lattice of some rank $n\in\N$ and let
$M:H_\Z\to H_\Z$ be an automorphism of it.

(a) The characteristic polynomial $p_{ch,M}(t)\in\Z[t]$ 
of the automorphism $M$ is unitary. 
Each eigenvalue $\lambda\in\C$ of $M$ is 
an algebraic integer and a unit in the ring 
$\OO_{\Q[\lambda]}\subset\Q[\lambda]$ of algebraic integers
\index{$\OO_{\Q[\lambda]},\ \OO_{\Q[\lambda]}^*$} 
in $\Q[\lambda]$, so in $\OO_{\Q[\lambda]}^*$.
Also $\lambda^{-1}\in \OO_{\Q[\lambda]}^*$.

(b) Suppose that $M$ is {\sf regular}, 
\index{regular endomorphism}
that means, 
$M:H_\C\to H_\C$ has for each eigenvalue only one Jordan block.
\begin{list}{}{}
\item[(i)]
Then 
\begin{eqnarray*}
\Q[M]=\bigoplus_{i=0}^{n-1}\Q M^i\stackrel{!}{=}
\End(H_\Q,M)
:=\{g:H_\Q\to H_\Q\,|\, gM=Mg\}.
\end{eqnarray*}
\item[(ii)]
Consider a polynomial $p(t)=\sum_{i=0}^{n-1}p_it^i\in\Q[t]$
and the endomorphism $g=p(M)\in \End(H_\Q,M)$ 
of $H_\Q$ and $H_\C$. 
Then $g$ has the eigenvalue $p(\lambda)$ on the generalized
eigenspace $H_\lambda$ of $M$ with eigenvalue $\lambda$.
If $g\in \End(H_\Z,M)$ then 
$p(\lambda)\in\OO_{\Q[\lambda]}$ for each eigenvalue $\lambda$
of $M$. If $g\in G_\Z^M$ then 
$p(\lambda)\in \OO_{\Q[\lambda]}^*$ for each eigenvalue 
$\lambda$ of $M$.
\end{list}

(c) Suppose that $M$ is the monodromy of a bilinear lattice
$(H_\Z,L)$. Suppose that $M$ is regular. 
\begin{list}{}{}
\item[(i)]
Then 
\begin{eqnarray*}
G_\Z^{(0)}\cup G_\Z^{(1)}\subset \{p(M)\,|\, 
p(t)=\sum_{i=0}^{n-1}p_it^i\in\Q[t],\ p(M)\in\End(H_\Z),\\ \\
p(\lambda)p(\lambda^{-1})=1\textup{ for each eigenvalue }
\lambda\textup{ of }M\}.
\end{eqnarray*}
\item[(ii)] 
If $M$ is semisimple then $G_\Z=G_\Z^{(0)}=G_\Z^{(1)}$,
and this set is equal to the set on the right hand side of (i).
\end{list}
\end{lemma}

{\bf Proof:}
(a)  Trivial.

(b) (i) As $M$ is regular, one can choose a vector $c\in H_\Q$
with $H_\Q=\bigoplus_{i=0}^{n-1}\Q M^ic$. The inclusion
$\Q[M]\subset\End(H_\Q,M)$ is clear.
Suppose $g\in\End(H_\Q,M)$. Write 
$gc=p(M)c$ for some polynomial $p(t)=\sum_{i=0}^{n-1}p_it^i
\in\Q[t]$. As
$$gM^kc=M^kgc=M^kp(M)c=p(M)M^kc$$
for each $k\in\{0,1,...,n-1\}$, $g=p(M)$. Thus
$\Q[M]=\End(H_\Q,M)$.

(ii) Similar to part (a).

(c) (i) Suppose $g=p(M)\in G_\Z^{(k)}$ for some $k\in\{0;1\}$
with $p(t)=\sum_{i=0}^{n-1}p_it^i\in\Q[t]$. 
Recall $\Rad I^{(k)}=\ker(M-(-1)^{k+1}\id:H_\Z\to H_\Z)$.

For $\lambda\neq (-1)^{k+1}$ the generalized eigenspaces
$H_\lambda$ and $H_{\lambda^{-1}}$ are dual to one another
with respect to $I^{(k)}$, and $g$ has eigenvalue $p(\lambda)$
on $H_\lambda$ and eigenvalue $p(\lambda^{-1})$ on 
$H_{\lambda^{-1}}$. That $g$ respects $I^{(k)}$ implies
$p(\lambda)p(\lambda^{-1})=1$. 

For $\lambda=(-1)^{k+1}$, $g$ restricts to an automorphism
of the sublattice $H_\lambda\cap H_\Z$ of $H_\Z$ with 
determinant 
$\pm 1=\det (g|_{H_\lambda\cap H_\Z})
=p(\lambda)^{\dim H_\lambda}$, so $p(\lambda)=\pm 1$. 

(ii) Suppose additionally that $M$ is semisimple and
that $g=p(M)$ with $p(t)=\sum_{i=0}^{n-1}p_it^i\in\Q[t]$ 
satifies $g\in\End(H_\Z)$ and 
$p(\lambda)p(\lambda^{-1})=1$ for each eigenvalue $\lambda$
of $M$. As $M$ is regular and semisimple, each eigenvalue has
multiplicity 1. The (1-dimensional) 
eigenspaces $H_\lambda$ and $H_{\lambda^{-1}}$ are dual
to one another with respect to $L$. The conditions
$p(\lambda)p(\lambda^{-1})=1$ imply that $g$ respects $L$
and also that 
$\det g=\prod_{\lambda\textup{ eigenvalue}}p(\lambda)=\pm 1$. 
Together with
$g\in \End(H_\Z)$ this shows $g\in G_\Z$.
\hfill$\Box$ 

\bigskip
The situation in the following Lemma \ref{t5.2} arises
surprisingly often. One reason is Theorem \ref{t3.26} (c).
See the Remarks \ref{t5.3} below.

\begin{lemma}\label{t5.2}
Let $H_\Z$ be a $\Z$-lattice of some rank $n\in\N$,
let $M:H_\Z\to H_\Z$ and $M^{root}:H_\Z\to H_\Z$ be
automorphisms of $H_\Z$, and let $l\in\N$ and 
$\varepsilon\in\{\pm 1\}$ be such that the following holds:
$M$ is regular, 
$$(M^{root})^l=\varepsilon M,$$
and $M^{root}$ is {\sf cyclic}, 
\index{cyclic endomorphism}
that means, a vector
$c\in H_\Z$ with $\bigoplus_{i=0}^{n-1}\Z(M^{root})^ic=H_\Z$
exists.

(a) Then $M^{root}$ is regular, and 
\begin{eqnarray*}
\End(H_\Z,M) & = & \End(H_\Z,M^{root})=\Z[M^{root}],\\
\Aut(H_\Z,M)&=&\Aut(H_\Z,M^{root})=\{p(M^{root})\,|\,
p(t)=\sum_{i=0}^{n-1}p_it^i\in\Z[t],\\
&& p(\kappa)\in(\Z[\kappa])^*\textup{ for each 
eigenvalue }\kappa\textup{ of }M^{root}\}.
\end{eqnarray*}

(b) Suppose that $M$ is the monodromy of a bilinear lattice 
$(H_\Z,L)$ and that the set of eigenvalues of $M^{root}$ is
invariant under inversion, that means, with $\kappa$
an eigenvalue of $M^{root}$ also $\kappa^{-1}$ is an eigenvalue
of $M^{root}$. 
\begin{list}{}{}
\item[(i)]
Then $M^{root}\in G_\Z$ and
\begin{eqnarray*}
G_\Z^{(0)}\cup G_\Z^{(1)}&\subset& \{p(M^{root})\,|\, 
p(t)=\sum_{i=0}^{n-1}p_it^i\in\Z[t],\\
&&p(\kappa)p(\kappa^{-1})=1\textup{ for each eigenvalue }
\kappa\textup{ of }M^{root}\}.
\end{eqnarray*}
\item[(ii)] 
If $M$ is semisimple then $G_\Z=G_\Z^{(0)}=G_\Z^{(1)}$,
and this set is equal to the set on the right hand side if (i).
\end{list}
\end{lemma}

{\bf Proof:}
(a) $M^{root}$ is regular as $M$ is regular and 
$(M^{root})^l=\varepsilon M$. This equation also implies
$\Q[M]\subset \Q[M^{root}]$. As these $\Q$-vector spaces have 
both dimension $n$, 
$$\End(H_\Q,M)=\Q[M]=\Q[M^{root}]=\End(H_\Q,M^{root}).$$
Then $\End(H_\Z,M)=\End(H_\Z)\cap\Q[M^{root}]$.
Consider $g\in \End(H_\Z,M)$. Choose a cyclic generator $c\in H_\Z$
with $\bigoplus_{i=0}^{n-1}\Z(M^{root})^ic=H_\Z$.
Write $g(c)=p(M^{root})c$ for some polynomial 
$p(t)=\sum_{i=0}^{n-1}p_it^i\in\Z[t]$. As in the proof of
Lemma \ref{t5.1}, one finds $g=p(M^{root})$, so 
$\End(H_\Z,M)=\Z[M^{root}]$. The element $g=p(M^{root})$ above
is in $\Aut(H_\Z,M)$ if and only if $\det g=\pm 1$, 
and this holds if and only if the algebraic integer 
$p(\kappa)$ is a unit in $\Z[\kappa]$ for each 
eigenvalue $\kappa$ of $M^{root}$. 

(b) (i) The main point is to show $M^{root}\in G_\Z$.
As $M$ and $M^{root}$ are both regular, the map
$\kappa\mapsto \varepsilon \kappa^l$ is a bijection from
the set of eigenvalues of $M^{root}$ to the set of 
eigenvalues of $M$. For $\lambda=\varepsilon\kappa^l$,
the generalized eigenspaces $H_\lambda$ and $H_{\lambda^{-1}}$ 
of $M$ are the generalized eigenspaces of $M^{root}$ with
eigenvalues $\kappa$ and $\kappa^{-1}$. These two spaces
are dual to one another with respect to $I^{(0)}$
(if $\lambda\neq -1$), $I^{(1)}$ (if $\lambda\neq 1$) and $L$.

Consider the decomposition of $M^{root}$ into the commuting
semisimple part $M^{root}_s$ and unipotent part 
$M^{root}_u$ with nilpotent part $N^{root}$ with
$\exp(N^{root})=M^{root}_u$, and also the decomposition
$M=M_sM_u$ with nilpotent part $N$ with $\exp N=M_u$.

$M^{root}_s$ and $M_s$ respect $L$  
because they have eigenvalue $\kappa$ and $\lambda$ on $H_\lambda$
and eigenvalue $\kappa^{-1}$ and $\lambda^{-1}$ 
on $H_{\lambda^{-1}}$.

As $M$ and $M_s$ respect $L$, also $M_u$ respects $L$.
Therefore $N$ is an infinitesimal isometry.
Because $N=lN^{root}$, also $N^{root}$ is an infinitesimal
isometry. Therefore $M^{root}_u$ respects $L$.
Thus also $M^{root}=M^{root}_sM^{root}_u$ respects $L$,
so $M^{root}\in G_\Z$.

Part (ii) and the rest of part (i) are proved as 
part (b) in Lemma \ref{t5.1}. \hfill$\Box$

\begin{remarks}\label{t5.3}
(i) The pair $(H_\Z,M^{root})$ in Lemma \ref{t5.2} is an
{\it Orlik block} if $M^{root}$ is of finite order.
Orlik blocks will be defined and discussed in 
the beginning of section \ref{s10.3}.
They are important building blocks in the unimodular
bilinear lattices $(H_\Z,L)$ for many isolated hypersurface
singularities. 

(ii) If the matrix $S=L(\uuuu{e}^t,\uuuu{e})^t$ of a unimodular
bilinear lattice $(H_\Z,L)$ with a triangular basis 
$\uuuu{e}$ has the special shape in Theorem \ref{t3.26} (c),
then the monodromy $(-1)^{k+1} M$ has by Theorem \ref{t3.26} (c) 
a specific $n$-th root
$M^{root}\in G_\Z$. This situation is special, but it
arises surprisingly often, in singularity theory and 
in the cases in part (ii). 

(iii) Theorem \ref{t3.26} (c) applies to all matrices 
$S(x)=\begin{pmatrix}1&x\\0&1\end{pmatrix}$ with $x\in\Z$
and to all matrices 
$S=\begin{pmatrix}1&x&\varepsilon x\\ 0&1&x\\ 0&0&1\end{pmatrix}$
with $x\in\Z$ and $\varepsilon\in\{\pm 1\}$.

It applies especially to the matrices $S(A_1^3)$, $S(\widehat{A}_2)$,
$S(\HH_{1,2})$ and $S(\P^2)$ in the Examples \ref{t1.1}
and to the matrix $S(-1,1,-1)$ in the $\Br_3\ltimes\{\pm 1\}^3$ 
orbit of $S(A_3)$. Though it is not useful
in the cases $S(A_1^3)$ and $S(\HH_{1,2})$
because their monodromies are not regular. 
In the case $S(A_3)$ we will not use it as there the monodromy
itself is cyclic. 
\end{remarks}

The completely reducible cases $A_1^n$ for $n\in\N$ can be treated
easily, building on Lemma \ref{t2.12}.

\begin{lemma}\label{t5.4}
Fix $n\in\N$ and consider the case $A_1^n$ with $S=S(A_1^n)=E_n$.
Then \index{$A_1^n$}
\begin{eqnarray*}
G_\Z=G_\Z^{(0)}&\cong& O_n(\Z)=\{A\in GL_n(\{0;\pm 1\})\,|\, 
\exists\ \sigma\in S_n\\
&&\exists\ \varepsilon_1,...,\varepsilon_n\in\{\pm 1\}
\textup{ such that }A_{ij}=\varepsilon_i\delta_{i\sigma(j)}\},\\
G_\Z^{(1)}=G_\Z^M&=&\Aut(H_\Z)\cong GL_n(\Z).
\end{eqnarray*}
The map $Z:(\Br_n\ltimes\{\pm 1\}^n)_{S}=\Br_n\ltimes\{\pm 1\}^n
\to G_\Z$ is surjective.
\end{lemma}

{\bf Proof:} The groups $G_\Z^{(1)}$ and $G_\Z^M$ are as claimed because
$M=\id$ and $I^{(1)}=0$. The groups $G_\Z$ and $G_\Z^{(0)}$ map
the set $R^{(0)}=\{\pm e_1,...,\pm e_n\}$ to itself and are therefore
also as claimed. 

The stabilizer of $S=E_n$ is the whole group $\Br_n\ltimes\{\pm 1\}^n$. 
The subgroup $\{\pm 1\}^n$ gives all sign changes of the basis 
$\uuuu{e}=(e_1,...,e_n)$. The subgroup $\Br_n$ gives under $Z$ all permutations
of the elements of the tuple $(e_1,...,e_n)$. Therefore $Z$ is surjective.
\hfill$\Box$

\section{The rank 2 cases}\label{s5.2}

For $x\in\Z$ consider the matrix 
$S=S(x)=\begin{pmatrix}1&x\\0&1\end{pmatrix}\in T^{uni}_2(\Z)$,
and consider a unimodular bilinear lattice $(H_\Z,L)$ with
a triangular basis $\uuuu{e}=(e_1,e_2)$ with
$L(\uuuu{e}^t,\uuuu{e})^t=S$. Then
\begin{eqnarray*}
M\uuuu{e}&=&\uuuu{e}S^{-1}S^t=\uuuu{e}
\begin{pmatrix}1-x^2 & -x\\ x & 1\end{pmatrix},\\
p_{ch,M}(t)&=& t^2-(2-x^2)t+1\\ 
\textup{ with zeros }
\lambda_{1/2}&=&\frac{2-x^2}{2}\pm \frac{1}{2}x\sqrt{x^2-4}.
\end{eqnarray*}
Theorem \ref{t3.26} (c) applies with $(n,k,q_0,q_1)=(2,0,1,x)$, namely
\begin{eqnarray*}
\delta_2\sigma^{root}&=&\delta_2\sigma_1\in 
(\Br_2\ltimes\{\pm 1\}^2)_S,\\
M^{root}&:=&Z(\delta_2\sigma_1)\in G_\Z\\
\textup{with }
M^{root}\uuuu{e}&=&\uuuu{e}\begin{pmatrix}-x&-1\\1&0\end{pmatrix}\\
\textup{and }(M^{root})^2&=&-M.
\end{eqnarray*}
$M^{root}$ is regular and cyclic,
\begin{eqnarray*}
p_{ch,M^{root}}(t) &=& t^2+xt+1\\
\textup{ with zeros }
\kappa_{1/2}&=&-\frac{x}{2}\pm \frac{1}{2}\sqrt{x^2-4}\\
\textup{with }\kappa_i^2&=&-\lambda_i=-x\kappa_i-1,
\ \kappa_1+\kappa_2=-x,
\ \kappa_1\kappa_2=1.
\end{eqnarray*}
$M$ is regular if $x\neq 0$. $M$ and $M^{root}$ are 
semisimple if $x\neq \pm 2$. 
If $x=\pm 1$, $M$ has eigenvalues $e^{\pm 2\pi i/6}$
and $M^{root}$ has eigenvalues $e^{\pm 2\pi i /3}$
respectively $e^{\pm 2\pi i/6}$.
If $|x|\geq 3$, $M$ and $M^{root}$ have real eigenvalues and
infinite order. 
If $x=\pm 2$,  they have a $2\times 2$
Jordan block with eigenvalue $-1$ respectively 
$-\frac{x}{2}$.

\begin{theorem}\label{t5.5}
(a) If $x\neq 0$ then
\begin{eqnarray*}
G_\Z&=& G_\Z^{(0)}=G_\Z^{(1)}=\{\pm (M^{root})^l\,|\, l\in\Z\},\\
G_\Z^M&=& \left\{\begin{array}{ll} G_\Z & \textup{if }x\neq \pm 3,\\
\{\pm (M^{root}+\frac{x}{|x|}\id)^l\,|\, l\in\Z\} & \textup{if }x=\pm 3.
\end{array}\right.
\end{eqnarray*}
If $x=0$ then
\begin{eqnarray*}
&& G_\Z=G_\Z^{(0)}\cong O_2(\Z)=
\{\begin{pmatrix}\varepsilon_1 & 0\\ 0 &\varepsilon_2\end{pmatrix},
\begin{pmatrix} 0& \varepsilon_1 \\ \varepsilon_2&0\end{pmatrix}
\,|\, \varepsilon_1,\varepsilon_2\in\{\pm 1\}\},\\
&& G_\Z^{(1)}=G_\Z^M =\Aut(H_\Z)\cong GL_2(\Z).
\end{eqnarray*}
In all cases $G_\Z=G_\Z^{\BB}$, so $Z:(\Br_2\ltimes\{\pm 1\}^2)_S\to G_\Z$
is surjective. 

(b) Properties of $I^{(0)}$ and $I^{(1)}$: 
\begin{eqnarray*}
x=0:&& I^{(1)}=0,\ \Rad I^{(1)}=H_\Z,\ L(\uuuu{e}^t,\uuuu{e})^t=E_2,\ 
I^{(0)}(\uuuu{e}^t,\uuuu{e})=2E_2.\\
x\neq 0:&& \Rad I^{(1)}=\{0\}.\\
|x|\leq 1:&& I^{(0)}\textup{ is positive definite.}\\
|x|=2:&& I^{(0)}\textup{ is positive semi-definite},\ 
\Rad I^{(0)}=\Z(e_1-\frac{x}{|x|}e_2).\\
|x|>2:&& I^{(0)}\textup{ is indefinite},\ \Rad I^{(0)}=\{0\}.
\end{eqnarray*}
\end{theorem}

{\bf Proof:}
Part (b) is obvious. 

The case $x=0$ is the case $A_1^2$.
It is covered by Lemma \ref{t5.4}. 

Consider the cases $x\neq 0$ in part (a). We can restrict
to $x<0$ because of 
$L((e_1,-e_2)^t,(e_1,-e_2))^t=\begin{pmatrix} 1& -x\\0&1\end{pmatrix}$.
So suppose $x<0$. 

We know $\{\pm (M^{root})^l\,|\, l\in\Z\}\subset G_\Z\subset G_\Z^M$.
By Lemma \ref{t5.2} (a)
\begin{eqnarray*}
G_\Z^M=\{p(M^{root})\,|\, p(t)=p_1t+p_0\in\Z[t], \ 
p(\kappa_1)p(\kappa_2)\in\{\pm 1\}\}.
\end{eqnarray*}
The map
\begin{eqnarray*}
Q_2:\Z^2\to\Z,\quad (p_1,p_0)\mapsto 
p(\kappa_1)p(\kappa_2)=p_1^2-p_1p_0x+p_0^2,
\end{eqnarray*}
is a quadratic form. Lemma \ref{t5.2} (a) shows
\begin{eqnarray*}
G_\Z^M=\{p_1M^{root}+p_0\id\,|\, 
(p_1,p_0,\varepsilon_1)\in\Z^2\times\{\pm 1\},
Q_2(p_1,p_0)=\varepsilon\}.
\end{eqnarray*}
Lemma \ref{t5.2} (b) (ii) shows for 
$x\neq -2$ 
\begin{eqnarray*}
G_\Z=G_\Z^{(0)}=G_\Z^{(1)}=\{p_1M^{root}+p_0\id\,|\, 
(p_1,p_0)\in\Z^2, Q_2(p_1,p_0)=1\}.
\end{eqnarray*}
For $x=-2$ it shows only 
\begin{eqnarray*}
G_\Z,G_\Z^{(0)},G_\Z^{(1)}
\subset \{p_1M^{root}+p_0\id\,|\, (p_1,p_0)\in\Z^2, Q_2(p_1,p_0)=1\}.
\end{eqnarray*}
We discuss the cases $x=-1$, $x=-2$ and $x\leq -3$ separately.

{\bf The case $x=-1$:}
$Q_2$ is positive definite, so $Q_2(p_1,p_0)=-1$ is impossible, 
and 
\begin{eqnarray*}
\{(p_1,p_0)\,|\, Q_2(p_1,p_0)=1\}
=\{\pm (0,1),\pm(1,0),\pm(1,-1)\},\\
G_\Z=G_\Z^{(0)}=G_\Z^{(1)}=G_\Z^M=\{\pm\id,\pm M^{root},\pm (M^{root})^2\}.
\end{eqnarray*}
Because of $M^2=-M^{root}$ this equals
$\{\pm \id,\pm M,\pm M^2\}$. 

{\bf The case $x=-2$:} This follows from Lemma \ref{t5.6} below.

Remark: Here $Q_2$ is positive semidefinite with 
$Q_2(p_1,p_0)=(p_1+p_0)^2$. The solution 
$(p_1,p_0)=(p_1,-p_1+\varepsilon_2)$ with $\varepsilon_2\in\{\pm 1\}$ 
of $Q_2(p_1,p_0)=1$ corresponds to
\begin{eqnarray*}
p_1M^{root}+(-p_1+\varepsilon_2)\id 
&=&\varepsilon_2(\id+\varepsilon_2p_1 (M^{root}-\id))\\
&=&\varepsilon_2(\id+(M^{root}-\id))^{\varepsilon_2p_1}\\
&=&\varepsilon_2 (M^{root})^{\varepsilon_2p_1}.
\end{eqnarray*}

{\bf The cases $x\leq -3$:}
The arguments above show
\begin{eqnarray*}
&&G_\Z^M\cong (\Z[\kappa_1])^*\\
&\supset &G_\Z\cong\{p_1\kappa_1+p_0\in\Z[\kappa_1]\,|\, 
(p_1\kappa_1+p_0)(p_1\kappa_2+p_0)=1\}.
\end{eqnarray*}
So we need to understand the unit group of $\Z[\kappa_1]$
and the subgroup of elements with norm 1.
Both are treated in Lemma \ref{tc.1} (a) in Appendix C.

It remains to show $G_\Z=G_\Z^{\BB}$ in all cases $x\in\Z_{\leq -1}$. 
This follows from $G_\Z=\{\pm (M^{root})^l\,|\, l\in\Z\}$ and
$$Z(\delta_1\delta_2)=-\id,\quad Z(\delta_2\sigma_1)=M^{root}.
\hspace*{2cm}\Box$$

\begin{lemma}\label{t5.6}
Let $H_\Z$ be a $\Z$-lattice of rank 2, and let $\www{M}:H_\Z\to H_\Z$
be an automorphism with a $2\times 2$ Jordan block and eigenvalue
$\lambda\in\{\pm 1\}$. 

(a) Then a cyclic automorphism $\www{M}^{root}:H_\Z\to H_\Z$
with eigenvalue 1 
and a number $l\in\N$ with $(\www{M}^{root})^l=\lambda M$ exist.
They are unique. 
$$\Aut(H_\Z,\www{M})=\{\pm (\www{M}^{root})^l\,|\, l\in\Z\}.$$

(b) If $\www{I}:H_\Z\times H_\Z\to\Z$ is an $\www{M}$-invariant bilinear form
then it is also $\www{M}^{root}$-invariant
and 
$$\Aut(H_\Z,\www{M},\www{I})=\Aut(H_\Z,\www{M})
=\{\pm (\www{M}^{root})^l\,|\, l\in\Z\}.$$
\end{lemma}

{\bf Proof:} (a) There is a $\Z$-basis $\uuuu{f}=(f_1,f_2)$ of $H_\Z$
and an $l\in\N$ with 
$$\www{M}\uuuu{f}=\uuuu{f}\lambda\begin{pmatrix}1&l\\0&1\end{pmatrix}.$$
Here $f_1$ is a generator of the rank 1 $\Z$-lattice 
$\ker(\www{M}-\lambda\id)\subset H_\Z$. It is unique up to the sign. 
It is a primitive element of $H_\Z$. 
An element $f_2$ with $H_\Z=\Z f_1\oplus \Z f_2$ exists.
It is unique up to sign and up to adding a multiple of $f_1$.
The sign is fixed by $l$ in the matrix above being positive.
$l$ is unique. 

Define $\www{M}^{root}:H_\Z\to H_\Z$ by
$$\www{M}^{root}\uuuu{f}=\uuuu{f}\begin{pmatrix}1&1\\0&1\end{pmatrix}.$$
Obviously $(\www{M}^{root})^l=\lambda \www{M}$.

Any $g\in \Aut(H_\Z,\www{M})$ must fix $\Z f_1=\ker(\www{M}-\lambda\id)$.
Therefore it must be up to the sign a power of $\www{M}^{root}$. 

(b) That $\www{M}^{root}$ respects $\www{I}$ follows 
by the same arguments as $M^{root}\in G_\Z$ in 
the proof of Lemma \ref{t5.2} (b) (i)
(but now the situation is simpler, as $\www{M}$ and $\www{M}^{root}$
have a single $2\times 2$ Jordan block).
The rest follows with part (a). \hfill$\Box$

\section{Generalities on the rank 3 cases}\label{s5.3}

For $\uuuu{x}=(x_1,x_2,x_3)\in\Z^3$ consider the matrix 
$S=S(\uuuu{x})=
\begin{pmatrix}1&x_1&x_2\\0&1&x_3\\0&0&1\end{pmatrix}
\in T^{uni}_3(\Z)$,
and consider a unimodular bilinear lattice $(H_\Z,L)$ with
a triangular basis $\uuuu{e}=(e_1,e_2,e_3)$ with
$L(\uuuu{e}^t,\uuuu{e})^t=S$. Then
\begin{eqnarray*}
S^{-1}&=&\begin{pmatrix}1&-x_1&x_1x_3-x_2\\0&1&-x_3\\0&0&1\end{pmatrix},\\
S^{-1}S^t&=&\begin{pmatrix}1-x_1^2-x_2^2+x_1x_2x_3&-x_1-x_2x_3+x_1x_3^2
&x_1x_3-x_2\\x_1-x_2x_3&1-x_3^2&-x_3\\x_2&x_3&1\end{pmatrix},\\
M\uuuu{e}&=&\uuuu{e}S^{-1}S^t,
\end{eqnarray*}
\begin{eqnarray*}
I^{(0)}(\uuuu{e}^t,\uuuu{e})&=&S+S^t=
\begin{pmatrix}2&x_1&x_2\\x_1&2&x_3\\x_2&x_3&2\end{pmatrix},\\
I^{(1)}(\uuuu{e}^t,\uuuu{e})&=&S-S^t=
\begin{pmatrix}0&x_1&x_2\\-x_1&0&x_3\\-x_2&-x_3&0\end{pmatrix},\\
p_{ch,M}(t)&=& (t-1)(t^2-(2-r(\uuuu{x}))t+1),
\end{eqnarray*}
where 
\begin{eqnarray*}
r:\Z^3\to\Z,\quad 
\uuuu{x}=(x_1,x_2,x_3)\mapsto x_1^2+x_2^2+x_3^2-x_1x_2x_3.
\end{eqnarray*}
\index{$r,\ r_\R$} 
For $(x_1,x_2,x_3)\neq (0,0,0)$ define
\begin{eqnarray*}
f_3&:=&\uuuu{e}\, 
\frac{1}{\gcd(x_1,x_2,x_3)}\begin{pmatrix}-x_3\\x_2\\-x_1 \end{pmatrix}.
\end{eqnarray*}
\index{$f_3$}
This is a primitive vector in $H_\Z$. 
\begin{eqnarray*}
\Rad I^{(1)} \stackrel{\textup{2.5 (a)(ii)}}{=}
\ker (M-\id)=\left\{\begin{array}{ll}
\Z f_3&\textup{if }(x_1,x_2,x_3)\neq (0,0,0),\\
H_\Z&\textup{if }(x_1,x_2,x_3)=(0,0,0).\end{array}\right.
\end{eqnarray*}
Also
\begin{eqnarray*}
p_{ch,S+S^t}(t)&=& t^3-6t^2+(12-x_1^2-x_2^2-x_3^2)t-2(4-r),\\
L(f_3,f_3)&=&\frac{r(\uuuu{x})}{\gcd(x_1,x_2,x_3)^2},\quad 
I^{(0)}(f_3,f_3)=2L(f_3,f_3).
\end{eqnarray*}
The eigenvalues of $M$ are called
\begin{eqnarray*}
\lambda_{1/2}=\frac{2-r}{2}\pm \frac{1}{2}\sqrt{r(r-4)},\quad 
\lambda_3=1
\end{eqnarray*}
with $\lambda_1+\lambda_2=2-r$, $\lambda_1\lambda_2=1$.

The following Lemma gives implicitly precise information on
$p_{ch,M}$ and $\sign I^{(0)}$ for all $\uuuu{x}\in\Z^3$.
{\it Implicitly}, because one has to determine with the tools
from section \ref{s4.2} in the cases
$r(\uuuu{x})\in\{0,1,2,4\}$ in which $\Br_3\ltimes \{\pm 1\}^3$
orbit in Theorem \ref{t4.6} (e) the matrix $S(\uuuu{x})$  is.

\begin{lemma}\label{t5.7}
(a) $r^{-1}(3l)=\emptyset$ for $l\in\Z-3\Z$.

(b) Consider $\uuuu{x}\in\Z^3$ with $r=r(\uuuu{x})<0$
or $>4$ or with $S(\uuuu{x})$ one of the cases in
Theorem \ref{t4.6} (e). Then $p_{ch,M}$ and $\sign I^{(0)}$
are as follows.

\begin{eqnarray*}
\begin{array}{llll}
 & p_{ch,M} & \sign I^{(0)}\hspace*{1cm} & S(\uuuu{x}) \\
r<0 & \lambda_1,\lambda_2>0 & (+--) & S(\uuuu{x}) \\
r=0 & \Phi_1^3 & (+++) & S(A_1^3)\\
r=0 & \Phi_1^3 & (+--) & S(\P^2)\\
r=1 & \Phi_6\Phi_1 & (+++) & S(A_2A_1) \\
r=2 & \Phi_4\Phi_1 & (+++) & S(A_3) \\
r=4 & \Phi_2^2\Phi_1 & (++\ 0)
& S(\P^1A_1) \\
r=4 & \Phi_2^2\Phi_1 & (++\ 0)
& S(\whh{A}_2) \\
r=4 & \Phi_2^2\Phi_1 & (+\ 0\ 0)
& S(\HH_{1,2}) \\
r=4 & \Phi_2^2\Phi_1 & (+\ 0\ -)
& S(-l,2,-l) \textup{ with }l\geq 3 \\
r>4 & \lambda_1,\lambda_2<0 & (++-) & S(\uuuu{x})
\end{array}
\end{eqnarray*}
\end{lemma}

{\bf Proof:}
(a) If $(3| x_1,3| x_2,3| x_3)$ then $9| r$.

If $(3| x_1,3| x_2,3\nmid x_3)$ then 
$3| (r-1),3\nmid r$.

If $(3| x_1,3\nmid x_2,3\nmid x_3)$ then 
$3| (r-2),3\nmid r$.

If $(3\nmid x_1,3\nmid x_2,3\nmid x_3)$ then 
$3| (x_1^2+x_2^2+x_3^2),3\nmid r$.

(b) The statements on $p_{ch,M}$ are obvious. 
$$\Rad I^{(0)}\stackrel{\textup{2.5 (a)(ii)}}{=}
\ker (M+\id)\supsetneqq\{0\}
\iff \Phi_2|p_{ch,M}\iff r=4.$$
In the cases with $r=4$, one calculates $p_{ch,S+S^t}(t)$
and reads off $\sign I^{(0)}$ from the zeros of 
$p_{ch,S+S^t}(t)$. The case $S(A_1^3)=S(0,0,0)$ is trivial.

Consider the cases with $r\neq 4$ and $\uuuu{x}\neq (0,0,0)$. 
Then $I^{(0)}$ is nondegenerate.
The product of the signs in the signature of $I^{(0)}$ is the
sign of $\det(S+S^t)=2(4-r)$. Because of the $2$'s on the
diagonal of $S+S^t$, $I^{(0)}$ cannot be negative definite.
Also recall $I^{(0)}(f_3,f_3)=2r(\gcd(x_1,x_2,x_3))^{-2}$.
This shows $\sign I^{(0)}=(+--)$ for $r<0$,
$\sign I^{(0)}=(++-)$ for $r>4$ and 
$\sign I^{(0)}=(+++)$ or $(+--)$ for $r\in\{0,1,2\}$.

The classification of $\Br_3\ltimes\{\pm 1\}$ orbits
in $T^{uni}_3(\Z)$ in Theorem \ref{t4.6} (e) says that for each
of the cases $r\in\{0,1,2\}$ there is only one orbit
(with $\uuuu{x}\neq (0,0,0)$ in the case $r=0$),
namely $S(\P^2)$, $S(A_2A_1)$ and $S(A_3)$. 
One checks the claims on $\sign I^{(0)}=\sign (S+S^t)$
immediately.
\hfill$\Box$

\begin{remarks}\label{t5.8}
(i) It is very remarkable that the fibers $r^{-1}(1)$
and $r^{-1}(2)\subset\Z^3$ of $r:\Z^3\to\Z$ consist each
of only one orbit. If one looks at the fibers of the
real map
\begin{eqnarray*}
r_\R:\R^3\to\R,\quad \uuuu{x}\mapsto x_1^2+x_2^2+x_3^2-x_1x_2x_3,
\end{eqnarray*}
this does not hold. Each real fiber $r_\R^{-1}(\rho)$
with $\rho\in(0,4)$ has five components, one compact
(homeomorphic to a 2-sphere), four non-compact
(homeomorphic to $\R^2$). The four non-compact components
are related by the action of $G^{sign}$.
It is remarkable that the fibers $r_\R^{-1}(1)$ and 
$r_\R^{-1}(2)\subset\R^3$ intersect $\Z^3$ only in the central piece.

(ii) By $p_{ch,M}=(t-1)(t^2-(2-r(\uuuu{x}))t+1)$, the monodromy
matrix $S^{-1}S^t$ for $\uuuu{x}\in\R^3$ and $S=S(\uuuu{x})$ 
has all eigenvalues in $S^1$ if and only if $r_\R(\uuuu{x})\in[0,4]$. 

(iii) The semialgebraic subvariety $r_\R^{-1}([0,4])\subset\R^3$ 
was studied in \cite[5.2]{BH20}. It has a central piece 
which is $G^{sign}$ invariant and which looks like a 
tetrahedron with smoothened edges  
and four other pieces which are permuted by $G^{sign}$.
Each other piece is homeomorphic to $[0,1]\times\R^2$
and is glued in one point (its only singular point) to 
one of the vertices of the central piece. 
The four vertices are $(2,2,2),(2,-2,-2),(-2,2,-2),(-2,-2,2)$,
so the elements of the $\Br_3\ltimes\{\pm 1\}^3$ orbit of
$(2,2,2)$. 

For $\uuuu{x}$ inside the central piece
$S+S^t$ is positive definite, on its boundary except
the vertices $\sign (S+S^t)=(++0)$, at a vertex
$\sign(S+S^t)=(+00)$. 
On those boundary components of the other four pieces
which contain one of the vertices $\sign(S+S^t)=(+0-)$
(except at the vertex). 
On the interior of $r_\R^{-1}(-\infty,4]$ except the
central piece $\sign(S+S^t)=(+--)$. On the exterior
of $r_\R^{-1}(-\infty,4]$ $\sign(S+S^t)=(++-)$. 

(iv) Due to Lemma \ref{t5.7} (b) and Theorem \ref{t4.6}, the 
seven cases
$S(A_1^3)$, $S(\P^2)$, $S(A_2A_1)$, $S(A_3)$, $S(\P^1A_1)$,
$S(\widehat{A}_2)$, $S(\HH_{1,2})$ and the series 
$S(-l,2,-l)$ for $l\geq 3$ give the only rank 3 unimodular
bilinear lattices where all eigenvalues of the monodromy
are unit roots. In the sections \ref{s5.4} and \ref{s5.5} 
we will focus on the reducible cases and these cases. 
In the sections \ref{s5.6} and \ref{s5.7}
we will treat the other cases. 
\end{remarks}

The following definition presents a special automorphism
$Q$ in $\Aut(H_\Q,L)$. Theorem \ref{t5.11} will say in which
cases $Q$ is in $G_\Z$ and in which cases not. 
The determination of the group $G_\Z$ in all irreducible rank 3 cases
in the sections \ref{s5.5}--\ref{s5.7} will build on this result. 
It is preceded by Lemma \ref{t5.10} which provides notations
and estimates which will be used in the proof of Theorem
\ref{t5.11} and also later.

\begin{definition}\label{t5.9}
Consider $\uuuu{x}\in\Z^3$ with $r(\uuuu{x})\neq 0$.
Then $H_\Q=H_{\Q,1}\oplus H_{\Q,2}$ with
$H_{\Q,1}:=\ker(M^2-(2-r(\uuuu{x}))M+\id:H_\Q\to H_\Q)$ and
$H_{\Q,2}:=\ker(M-\id:H_\Q\to H_\Q)$. 
This decomposition is left and right $L$-orthogonal.
Then $Q:H_\Q\to H_\Q$ denotes the automorphism with
$Q|_{H_{\Q,1}}=-\id$ and $Q|_{H_{\Q,2}}=\id$. 
It is in $\Aut(H_\Q,L)$. \index{$Q$}
\end{definition}

\begin{lemma}\label{t5.10}
(a) For $\uuuu{x}\in \Z^3-\{(0,0,0)\}$ write $r:=r(\uuuu{x})$
and define \index{$g=g(\uuuu{x})=\gcd(x_1,x_2,x_3)$}
\index{$\www{\uuuu{x}}=(\www{x}_1,\www{x}_2,\www{x}_3)=g^{-1}\uuuu{x}$} 
\begin{eqnarray*} 
g:=g(\uuuu{x})&:=&\gcd(x_1,x_2,x_3)\in\N,\nonumber\\
\uuuu{\www{x}}:=(\www{x}_1,\www{x}_2,\www{x}_3)
&:=&g^{-1}\uuuu{x}\in\Z^3.
\end{eqnarray*}
Then $f_3=-\www{x}_3e_1+\www{x}_2e_2-\www{x}_1e_3$,
$\gcd(\www{x}_1,\www{x}_2,\www{x}_3)=1$ and 
\begin{eqnarray}
g^2\,|\, r,\quad 
\frac{r}{g^2}= \www{x}_1^2+\www{x}_2^2+\www{x}_3^2-g\www{x}_1\www{x}_2\www{x}_3.
\label{5.1}  
\end{eqnarray}

(b) Consider a local minimum (Definition \ref{t4.3}) 
$\uuuu{x}\in\Z^3_{\geq 3}$
with $x_i\leq x_j\leq x_k$ for some $i,j,k$ with
$\{i,j,k\}=\{1,2,3\}$. 
Then 
\begin{eqnarray}
\www{x}_i&\leq& \frac{2+(4-r)^{1/3}}{g},\label{5.2}\\
x_j^2&\leq & \frac{4-r}{x_i-2}+x_i+2,\label{5.3}\\
x_k&\leq& \frac{1}{2}x_ix_j 
\quad\textup{and}\quad
\www{x}_k\leq \frac{g}{2}\www{x}_i\www{x}_j.\label{5.4}
\end{eqnarray}
\end{lemma}

{\bf Proof:}
(a) Trivial.

(b) Lemma \ref{t4.4} shows $x_k\leq\frac{1}{2}x_ix_j$,
which is \eqref{5.4}. This is equivalent to 
$\www{x}_k\leq\frac{g}{2}\www{x}_i\www{x}_j$. 
We also have $\www{x}_i\leq \www{x}_j\leq\www{x}_k$,
and we know $r\leq 0$ from Theorem \ref{4.6}.

The proof of \eqref{5.2} is similar to the second case
in the proof of Theorem \ref{t4.6} (a). 
\begin{eqnarray*}
\frac{r}{g^2}&=& \www{x}_i^2+\www{x}_j^2+
(\www{x}_k-\frac{g}{2}\www{x}_i\www{x}_j)^2
-\frac{g^2}{4}\www{x}_i^2\www{x}_j^2\\
&\leq& \www{x}_i^2+\www{x}_j^2+
(\www{x}_j-\frac{g}{2}\www{x}_i\www{x}_j)^2
-\frac{g^2}{4}\www{x}_i^2\www{x}_j^2
\quad(\textup{because }\www{x}_j\leq\www{x}_k\leq 
\frac{g}{2}\www{x_i}\www{x}_j)\\
&=& \www{x}_i^2+2\www{x}_j^2-g\www{x}_i\www{x}_j^2\\
&=& (\www{x}_i-g\www{x}_j^2+\frac{2}{g})(\www{x}_i-\frac{2}{g})
+\frac{4}{g^2}.
\end{eqnarray*}
If $\www{x}_i<\frac{2}{g}$, then \eqref{5.2} holds
anyway. If $\www{x}_i\geq \frac{2}{g}$ then we can further estimate
the last formula using $-g\www{x}_j^2\leq -g\www{x}_i^2$.
Then we obtain 
\begin{eqnarray*}
\frac{r}{g^2}&\leq& (\www{x}_i-g\www{x}_i^2+\frac{2}{g})
(\www{x}_i-\frac{2}{g})+\frac{4}{g^2}\\
&=& -g(\www{x}_i-\frac{2}{g})^2(\www{x}_i+\frac{1}{g})
+\frac{4}{g^2}\\
&\leq& -g(\www{x}_i-\frac{2}{g})^3+\frac{4}{g^2},\\
\textup{so }
(\www{x}_i-\frac{2}{g})^3&\leq& \frac{4-r}{g^3}.
\end{eqnarray*}
This shows \eqref{5.2}. The inequality \eqref{5.3}
was proved within the second case in the proof of Theorem
\ref{t4.6} (a).\hfill$\Box$

\begin{theorem}\label{t5.11}
Consider $\uuuu{x}\in\Z^3$ with $r\neq 0$. The automorphism
$Q\in \Aut(H_\Q,L)$ which was defined in Definition
\ref{t5.9} can be written in two interesting ways,
\begin{eqnarray}\label{5.5}
Q&=& (\uuuu{e}\mapsto -\uuuu{e}+\frac{2g^2}{r}f_3
(-\www{x}_3,\www{x}_2-g\www{x}_1\www{x}_3,-\www{x}_1)),\\
Q&=& \id + 2(M-\id)+\frac{2}{r}(M-\id)^2.\label{5.6}
\end{eqnarray}
We have
\begin{eqnarray}
Q\in G_\Z\iff \frac{2g^2}{r}\in\Z 
\iff \frac{r}{g^2}\in\{\pm 1,\pm 2\}.\label{5.7}
\end{eqnarray}
This holds if and only if $\uuuu{x}$ is in the 
$\Br_3\ltimes \{\pm 1\}^3$ orbit of a triple in the following
set:
\begin{eqnarray*}
&&\{(x,0,0)\,|\, x\in\N\} \quad 
(\textup{these are the reducible cases except }A_1^3),\\
&\cup& \{(x,x,0)\,|\, x\in\N\}\quad
\textup{(these cases include }A_3),\\
&\cup& \{(-l,2,-l)\,|\, l\geq 2\textup{ even}\}\quad
\textup{(these cases include }\HH_{1,2}),\\
&\cup&\{(3,3,4),(4,4,4),(5,5,5),(4,4,8)\}.
\end{eqnarray*}
So, within the cases with $r\in\{0,1,2,3,4\}$, 
$Q$ is not defined for $A_1^3$ and $\P^2$,
and $Q\notin G_\Z$ for $\widehat{A}_2$ and
$S(-l,2,-l)$ with $l\geq 3$ odd. 
\end{theorem}

{\bf Proof:} First we prove \eqref{5.5}. 
The 2-dimensional subspace
$H_{\Q,1}=\ker(M^2-(2-r)M+\id)\subset H_\Q$, on which
$Q$ is $-\id$, can also be characterized as the right
$L$-orthogonal subspace $H_{\Q,1}=(\Q f_3)^{\perp}$ to 
$\Q f_3$. For $b=\uuuu{e}\cdot \uuuu{y}^t\in H_\Q$ with 
$\uuuu{y}\in\Q^3$ 
\begin{eqnarray*}
L(f_3,b)&=&(-\www{x}_3, \www{x}_2, -\www{x}_1)
\begin{pmatrix}1&0&0\\x_1&1&0\\x_2&x_3&1\end{pmatrix}
\begin{pmatrix}y_1\\y_2\\y_3\end{pmatrix}\\
&=&(-\www{x}_3, \www{x}_2-g\www{x}_1\www{x}_3,-\www{x}_1)
\begin{pmatrix}y_1\\y_2\\y_3\end{pmatrix},
\end{eqnarray*}
so
\begin{eqnarray}
H_{\Q,1}=\{\uuuu{e}\cdot \uuuu{y}^t\,|\, \uuuu{y}\in\Q^3,
0=(-\www{x}_3, \www{x}_2-g\www{x}_1\www{x}_3,-\www{x}_1)
\begin{pmatrix}y_1\\y_2\\y_3\end{pmatrix}\}.\label{5.8}
\end{eqnarray}
Denote the endomorphisms on the right hand sides of \eqref{5.5}
and \eqref{5.6} by $Q^{\eqref{5.5}}$ respectively
$Q^{\eqref{5.6}}$. The formulas \eqref{5.5} and \eqref{5.8}
show $Q^{\eqref{5.5}}|_{H_{\Q,1}}=-\id$. Also
\begin{eqnarray*}
Q^{\eqref{5.5}}(f_3)=-f_3+f_3\frac{2g^2}{r}
(\www{x}_3^2+(\www{x}_2-g\www{x}_1\www{x}_3)\www{x}_2
+\www{x}_1^2)
=f_3.
\end{eqnarray*}
Therefore $Q^{\eqref{5.5}}=Q$, so \eqref{5.5} holds. 

Now we prove \eqref{5.6}.
Because $(M-\id)(f_3)=0$, we have $Q^{\eqref{5.6}}(f_3)=f_3$. 
Consider $b\in H_{\Q,1}$. Then
\begin{eqnarray*}
0&=&(M^2-(2-r)M+\id)(b),\quad\textup{so}\\
(M-\id)^2(b)&=&-rM(b),\quad\textup{so}\\
Q^{\eqref{5.6}}(b)
&=& (\id+2(M-\id)+\frac{2}{r}(-rM))(b)
=(-id)(b)=-b.
\end{eqnarray*}
Therefore $Q^{\eqref{5.6}}=Q$, so \eqref{5.6} holds. 

\eqref{5.7} can be proved either with \eqref{5.6} and
Lemma \ref{t5.17} below or with \eqref{5.5}, which is easier
and which we do now. Observe $\gcd(-\www{x}_3,
\www{x}_2-g\www{x}_1\www{x}_3,-\www{x}_1)=
\gcd(\www{x}_1,\www{x}_2,\www{x}_3)=1$.
Also, $f_3$ is a primitive vector in $H_\Z$. This shows that
$Q(e_1),Q(e_2),Q(e_3)$ are all in $H_\Z$ if and only if
$\frac{2g^2}{r}\in\Z$. This shows \eqref{5.7}. 

It is easy to see that all triples in the set in Theorem 
\ref{t5.11} satisfy $\frac{r}{g^2}\in\{\pm 1,\pm 2\}$:
\begin{eqnarray*}
\frac{r}{g^2}((x,0,0))=\frac{x^2}{x^2}=1,\quad
\frac{r}{g^2}((x,x,0))=\frac{2x^2}{x^2}=2,\\
\textup{for even }l\quad \frac{r}{g^2}((-l,2,-l))
=\frac{4}{4}=1,
\end{eqnarray*}
\begin{eqnarray*}
\frac{r}{g^2}((3,3,4))= \frac{-2}{1} = -2,\quad  
\frac{r}{g^2}((4,4,4))= \frac{-16}{16} = -1,\\
\frac{r}{g^2}((5,5,5))= \frac{-50}{25} = -2,\quad
\frac{r}{g^2}((4,4,8))= \frac{-32}{16} = -2.
\end{eqnarray*}

The difficult part is to see that there are no other
$\uuuu{x}\in\Z^3$ with $r\neq 0$
and $\frac{r}{g^2}\in\{\pm 1,\pm 2\}$.
It is sufficient to consider local minima
(Definition \ref{t4.3}). 
The calculations 
\begin{eqnarray*}
\frac{r}{g^2}((-l,2,-l))=\frac{4}{1}=4
\quad\textup{for odd }l\geq 3\\
\textup{and}\quad \frac{r}{g^2}((-1,-1,-1))=\frac{4}{1}=4
\end{eqnarray*}
deal with the other cases with $r\in\{1,2,3,4\}$, see
Theorem \ref{t4.6} (e). 

Consider $\uuuu{x}\in\Z^3_{\leq 0}$ with $x_i\leq x_j\leq x_k$
for some $i,j,k$ with $\{i,j,k\}=\{1,2,3\}$.
Then $\frac{r}{g^2}=\www{x}_1^2+\www{x}_2^2+\www{x}_3^2
+g|\www{x}_1\www{x}_2\www{x}_3|$ can be $1$ or $2$ only
if $\www{x}_k=0$ and $\www{x}_i=\www{x}_j=-1$.
Then $\uuuu{x}=(-g,-g,0)$. This is in the 
$\Br_3\ltimes\{\pm 1\}^3$ orbit of $(g,g,0)$. 

Consider a local minimum $\uuuu{x}\in\Z^3_{\geq 3}$ with
$x_i\leq x_j\leq x_k$ for some $i,j,k$ with 
$\{i,j,k\}=\{1,2,3\}$. Suppose $\frac{r}{g^2}\in\{-1,-2\}$. 
We have to show
$\uuuu{x}\in\{(3,3,4),(4,4,4),(5,5,5),(4,4,8)\}$. 
Of course $\uuuu{x}\neq (3,3,3)$ because $r\neq 0$.

If $g=1$ then $r\in\{-1,-2\}$. One sees easily that 
$r=-1$ is impossible and that $r=-2$ is only satisfied
for $\uuuu{x}=(3,3,4)$. 

From now on suppose $g\geq 2$. 
Write $\rho:=|\frac{r}{g^2}|\in\{1,2\}$.
\eqref{5.2} takes the shape
\begin{eqnarray*}
\www{x}_i\leq 
\frac{2}{g}+(\frac{\rho}{g}+\frac{4}{g^3})^{1/3}.
\end{eqnarray*}
The only pairs $(\www{x}_i,g)\in\N\times\Z_{\geq 2}$ 
which satisfy this and $x_i=\www{x}_ig\geq 3$ 
are in the following two tables,
\begin{eqnarray*}
\rho=1: \begin{array}{l|l|l|l|l}
g & 2 & 3 & 4 & 5 \\ \hline 
\www{x}_i & 2 & 1 & 1 & 1 \end{array},\quad
\rho=2: \begin{array}{l|l|l|l|l|l}
g & 2 & 3 & 4 & 5 & 6 \\ \hline 
\www{x}_i & 2 & 1 & 1 & 1 & 1\end{array}
\end{eqnarray*}
The following table discusses these nine cases. 
Three of them lead to $(4,4,4)$, $(5,5,5)$ and $(4,4,8)$,
six of them are impossible. The symbol $\circledast$
denotes {\it impossible}. The inequalities 
$x_i\leq x_j\leq x_k$ and \eqref{5.3} and \eqref{5.4}
are used, and also $x_i=g\www{x}_i$ and $x_j,x_k\in g\N$.
\begin{eqnarray*}
\begin{array}{l|l|l|l|l|l|l|l|l}
\rho & g & -r & \www{x}_i & x_i & 
x_j^2\leq \frac{4-r}{x_i-2}+x_i+2 & x_j & 
x_k\leq \frac{1}{2}x_ix_j & x_k \\ \hline 
1 & 2 & 4 & 2 & 4 & x_j^2 \leq 10 & \circledast & & \\
1 & 3 & 9 & 1 & 3 & x_j^2 \leq 18 & 3 & x_k\leq \frac{9}{2} & 
3 \ \circledast \\
1 & 4 & 16 & 1 & 4 & x_j^2 \leq 16 & 4 & x_k\leq 8 & 4 \\
1 & 5 & 25 & 1 & 5 & x_j^2 \leq 16+\frac{2}{3} & \circledast & & \\
2 & 2 & 8 & 2 & 4 & x_j^2 \leq 12 & \circledast & & \\
2 & 3 & 18 & 1 & 3 & x_j^2 \leq 27 & 3 & x_k\leq \frac{9}{2} &
3 \ \circledast \\
2 & 4 & 32 & 1 & 4 & x_j^2 \leq 24 & 4 & x_k\leq 8 & 8 \\
2 & 5 & 50 & 1 & 5 & x_j^2 \leq 25 & 5 & x_k\leq 12+\frac{1}{2} & 
5 \\
2 & 6 & 72 & 1 & 6 & x_j^2 \leq 27 & \circledast & & 
\end{array}
\end{eqnarray*}
This finishes the proof of Theorem \ref{t5.11}.\hfill$\Box$

\section{The reducible rank 3 cases}
\label{s5.4}

Definition \ref{t2.10} proposed the notion of reducible triple
$(H_\Z,L,\uuuu{e})$ where $(H_\Z,L)$ is a unimodular bilinear
lattice and $\uuuu{e}$ is a triangular basis. The following
Remarks propose the weaker notion when a unimodular bilinear
lattice $(H_\Z,L)$ (without a triangular basis) is reducible.
Then the groups $G_\Z,G_\Z^{(0)},G_\Z^{(1)}$ and $G_\Z^M$
split accordingly if also the eigenvalues of $M$ split
in a suitable sense.

\begin{remarks}\label{t5.12}
(i) Suppose that a unimodular bilinear lattice $(H_\Z,L)$
splits into a direct sum $H_{\Z,1}\oplus H_{\Z,2}$
which is left and right $L$-orthogonal. Then the restrictions
of $L,I^{(0)},I^{(1)}$ and $M$ to $H_{\Z,i}$ are called
$L_i,I^{(0)}_i,I^{(1)}_i$ and $M_i$ for $i\in\{1;2\}$. 
We say that $(H_\Z,L)$ is {\it reducible} and that it splits into
the direct sum $(H_{\Z,1},L_1)\oplus (H_{\Z,2},L_2)$. 

(ii) In the situation of (i), suppose that the eigenvalues
of $M_1$ are pairwise different from the eigenvalues of $M_2$. 
Then any element of $G_\Z^M$ respects the splitting.
For $\{i,j\}=\{1,2\}$ write $G_{\Z,i}^M:=G_\Z^M(H_{\Z,i},L_i)$.
Then
\begin{eqnarray*}
G_\Z^M = G_{\Z,1}^M\times G_{\Z,2}^M,
\end{eqnarray*}
and with analogous notations
\begin{eqnarray*}
G_\Z = G_{\Z,1}\times G_{\Z,2},
\ G_\Z^{(0)} = G_{\Z,1}^{(0)}\times 
G_{\Z,2}^{(0)},
\ G_\Z^{(1)} = G_{\Z,1}^{(1)}\times 
G_{\Z,2}^{(1)}.
\end{eqnarray*}

(iii) There is only one unimodular bilinear lattice of rank 1.
We call it $A_1$-lattice and denote the matrix
$S=S(A_1)=(1)\in M_{1\times 1}(\Z)$. 
Here $G_\Z=G_\Z^{(0)}=G_\Z^{(1)}=G_\Z^M=\{\pm\id\}$.

(iv) Suppose that the characteristic polynomial 
$p_{ch,M}(t)\in\Z[t]$ of the monodromy $M$ of a unimodular
bilinear lattice $(H_\Z,L)$ 
splits into a product $p_{ch,M}=p_1p_2$
of non-constant polynomials $p_1$ and $p_2$
with $\gcd(p_1,p_2)=1$. Then 
$\ker p_1(M)\oplus \ker p_2(M)$ is a sublattice of finite
index in $H_\Z$, and the summands are left and right
$L$-orthogonal to one another. If the index is 1, 
we are in the situation of (ii). Theorem \ref{t5.14}
will show that this applies to the cases
$S(\HH_{1,2})$ and $S(-l,2,-l)$ with $l\geq 4$ even, 
but not to the cases $S(A_3)$, $S(\widehat{A}_2)$ and 
$S(-l,2,-l)$  with $l\geq 3$ odd.
\end{remarks}

These remarks apply especially to the reducible $3\times 3$ cases
except $A_1^3$ (which is part of Lemma \ref{t5.4}). 
This includes the two reducible cases $A_2A_1$ and
$\P^1A_1$ with eigenvalues in $S^1$.

\begin{theorem}\label{t5.13}
Consider $\uuuu{x}=(x,0,0)\in\Z^3$ with $x\neq 0$
and the unimodular bilinear lattice $(H_\Z,L,\uuuu{e})$
with triangular basis $\uuuu{e}$ with 
$L(\uuuu{e}^t,\uuuu{e})^t=S(\uuuu{x})\in T^{uni}_3(\Z)$. 

Then $(H_\Z,L,\uuuu{e})$ is reducible with the summands 
$(H_{\Z,1},L_1,(e_1,e_2))$ and $(H_{\Z,2},L_2,e_3)$ with
$H_{\Z,1}=\Z e_1\oplus \Z e_2$ and $H_{\Z,3}=\Z e_3$. 
The first summand is an irreducible rank two unimodular
bilinear lattice with triangular basis. Its groups 
$G_{\Z,1},G_{\Z,1}^{(0)},G_{\Z,1}^{(1)}$ and $G_{\Z,1}^M$
are treated in Theorem \ref{t5.5}.
The second summand is of type $A_1$. 
See Remark \ref{t5.12} (iii) for its groups. 

The decompositions in Remark \ref{t5.12} (ii)
hold for the groups $G_\Z$, $G_\Z^{(0)}$, $G_\Z^{(1)}$ and $G_\Z^M$.
Here $G_\Z^M=G_\Z$ if $x\neq \pm 3$ and 
$G_\Z=G_\Z^{(0)}=G_\Z^{(1)}$ always.

The map $Z:(\Br_3\ltimes\{\pm 1\}^3)_S\to G_\Z$ is surjective.
\end{theorem}

{\bf Proof:} The first point to see is that Remark \ref{t5.12}
(ii) applies. It does because the characteristic polynomials
of the monodromies $M_1$ and $M_2$ of the two summands 
are $t^2-(2-x^2)t+1$ and $t-1$, and here $x\neq 0$,
so that the eigenvalues of $M_1$ are not equal to the eigenvalue $1$
of $M_2$. 

The second point to see is the surjectivity of the map $Z$. 
This follows from the surjectivity of the map $Z$ in the irreducible 
rank 2 cases in Theorem \ref{t5.5} and in the case $A_1$ in Lemma
\ref{t5.4}.
\hfill$\Box$

\section{The irreducible rank 3 cases with all eigenvalues in $S^1$}
\label{s5.5}

Theorem \ref{t5.14} is the only point in this section.
It treats the irreducible rank 3 cases with all eigenvalues in $S^1$.

\begin{theorem}\label{t5.14}
Consider for each of the matrices
$S(\P^2)$, $S(A_3)$, $S(\widehat{A}_2)$, 
$S(\HH_{1,2})$ and $S(-l,2,-l)$ for $l\geq 3$
in the Examples \ref{t1.1} a unimodular bilinear lattice
$(H_\Z,L)$ with a triangular basis $\uuuu{e}$ with
$L(\uuuu{e}^t,\uuuu{e})^t=S$. 

(a) The cases $S(\HH_{1,2})$ 
and $S(-l,2,-l)$ for $l\geq 4$ even:
Then $(H_\Z,L)$ is reducible 
(in the sense of Remark \ref{t5.12} (i)), 
$H_\Z=H_{\Z,1}\oplus H_{\Z,2}$ with
\begin{eqnarray*}
H_{\Z,1}&:=&\ker(M+\id)^2\textup{ of rank 2}\\
\textup{and }H_{\Z,2}&:=&\ker(M-\id)\textup{ of rank 1.}
\end{eqnarray*}
$(H_{\Z,2},L_2)$ is an $A_1$-lattice. In all cases
the decompositions in Remark \ref{t5.12} (ii) hold
for the groups $G_\Z^M$, $G_\Z^{(0)}$, $G_\Z^{(1)}$ and $G_\Z$. 
The groups $G_{\Z,1}^M$, $G_{\Z,1}^{(0)}$, $G_{\Z,1}^{(1)}$
and $G_{\Z,1}$ are as follows.
\begin{list}{}{}
\item[(i)]
$S(\HH_{1,2})$: $H_{\Z,1}$ has a $\Z$-basis 
$\uuuu{f}=(f_1,f_2)$ with 
\begin{eqnarray*}
L(\uuuu{f}^t,\uuuu{f})^t=\begin{pmatrix}0&-1\\1&0\end{pmatrix},
\quad 
I^{(0)}(\uuuu{f}^t,\uuuu{f})=\begin{pmatrix}0&0\\0&0
\end{pmatrix},\\
I^{(1)}(\uuuu{f}^t,\uuuu{f})=\begin{pmatrix}0&-2\\2&0
\end{pmatrix},\quad
M\uuuu{f}=\uuuu{f}\begin{pmatrix}-1&0\\0&-1\end{pmatrix},
\end{eqnarray*}
\begin{eqnarray*}
G_{\Z,1}^M&=&G_{\Z,1}^{(0)}=\Aut(H_{\Z,1})\cong GL_2(\Z),\\
G_{\Z,1}&=&G_{\Z,1}^{(1)}=\{g\in\Aut(H_{\Z,1})\,|\, \det g=1\}
\cong SL_2(\Z).
\end{eqnarray*}
\item[(ii)]
$S(-l,2,-l)$ for $l\geq 4$ even: $H_{\Z,1}$ has a $\Z$-basis 
$\uuuu{f}=(f_1,f_2)$ with 
\begin{eqnarray*}
L(\uuuu{f}^t,\uuuu{f})^t
=\begin{pmatrix}0&-1\\1&1-\frac{l^2}{4}\end{pmatrix},
\quad 
I^{(0)}(\uuuu{f}^t,\uuuu{f})
=\begin{pmatrix}0&0\\0&2-\frac{l^2}{2}\end{pmatrix},\\
I^{(1)}(\uuuu{f}^t,\uuuu{f})
=\begin{pmatrix}0&-2\\2&0
\end{pmatrix},\quad
M\uuuu{f}=\uuuu{f}
\begin{pmatrix}-1&2-\frac{l^2}{2}\\0&-1\end{pmatrix},
\end{eqnarray*}
Define $M_1^{root}\in \Aut(H_{\Z,1})$ by 
$M_1^{root}\uuuu{f}=\uuuu{f}\begin{pmatrix}1&1\\0&1\end{pmatrix}$.
Then 
\begin{eqnarray*}
(M_1^{root})^{l^2/2-2}=-M_1\qquad \textup{ and }\\
G_{\Z,1}=G_{\Z,1}^{(0)}=G_{\Z,1}^{(1)}=G_{\Z,1}^M
=\{\pm (M_1^{root})^m\,|\, m\in\Z\}.
\end{eqnarray*}
\end{list}

(b) The cases $S(\P^2)$, $S(A_3)$, $S(\widehat{A}_2)$, 
$S(-l,2,-l)$ with $l\geq 3$ odd: 
Then $(H_\Z,L)$ is irreducible.
\begin{eqnarray*}
G_\Z=G_\Z^{(0)}=\{\pm(M^{root})^m\,|\, m\in\Z\}.
\end{eqnarray*}
Here $M^{root}$ is defined by \index{$M^{root}$}
\begin{eqnarray}
M^{root}&:=& Z(\sigma^{root})
\quad\textup{ in the case }S(\P^2),\label{5.9}\\
M^{root}&:=& M= Z(\sigma^{mon})
\quad\textup{ in the case }S(A_3),\label{5.10}\\
M^{root}&:=& Z(\delta_3\sigma^{root})
\quad\textup{ in the case }S(\whh{A}_2),\label{5.11}\\
M^{root}&:=& (-M)\circ 
Z(\delta_3\sigma_1^{-1}\sigma_2^{-1}\sigma_1)^{(5-l^2)/2}
\nonumber \\
&&\textup{ in the case }S(-l,2,-l)\textup{ for odd }l\geq 3.
\label{5.12}
\end{eqnarray}
It satisfies 
\begin{eqnarray*}
M^{root}\uuuu{e}=\uuuu{e} M^{root,mat}\quad\textup{and}\quad
(M^{root})^m=\varepsilon M
\end{eqnarray*}
where 
$M^{root,mat}$, $m$ and $\varepsilon$ are as follows:
\begin{eqnarray*}
\begin{array}{ccc}
 & S(\P^2) & S(A_3) \\
M^{root,mat} & 
\begin{pmatrix}3&-3&1\\1&0&0\\0&1&0 \end{pmatrix} & 
\begin{pmatrix}0&0&1\\-1&0&1\\0&-1&1 \end{pmatrix} \\
(m,\varepsilon) & (3,1) & (1,1)\\ \hline \hline 
 & S(\widehat{A}_2) & 
 S(-l,2,-l)\textup{ with }l\geq 3\textup{ odd}\\
M^{root,mat} & 
\begin{pmatrix}1&1&-1\\1&0&0\\0&1&0 \end{pmatrix} & 
\frac{1}{2}
\begin{pmatrix}1-l^2&l^3-l&-1-l^2\\
-2l&2l^2-2&-2l\\ 1-l^2&l^3-3l&3-l^2 \end{pmatrix} \\
(m,\varepsilon) &  (3,-1) & (l^2-4,-1)
\end{array}
\end{eqnarray*}
$M$ and $M^{root}$ are regular. $M^{root}$ is cyclic.
In the cases $S(A_3)$, $S(\widehat{A}_2)$ and 
$S(-l,2,-l)$ for $l\geq 3$ odd
\begin{eqnarray*}
G_\Z=G_\Z^{(0)}=G_\Z^{(1)}=G_\Z^M
=\{\pm (M^{root})^m\,|\, m\in\Z\}.
\end{eqnarray*}
Some additional information:
\begin{list}{}{}
\item[(i)]
$S(\P^2)$: $\sign I^{(0)}=(+--)$, 
$p_{ch,M}=p_{ch,M^{root}}=\Phi_1^3$,
$M$ and $M^{root}$ have a $3\times 3$ Jordan block,
\begin{eqnarray*}
G_\Z^{(1)}=G_\Z^M
=\{\pm (M^{root})^m(\id+a(M^{root}-\id)^2)\,|\, m,
a\in\Z\}\supsetneqq G_\Z.
\end{eqnarray*}
\item[(ii)]
$S(A_3)$: $\sign I^{(0)}=(+++)$, 
$p_{ch,M}=\Phi_4\Phi_1$, $M=M^{root}$, 
$|G_\Z|=8$.
\item[(iii)]
$S(\widehat{A}_2)$: $\sign I^{(0)}=(++\, 0)$, 
$p_{ch,M}=\Phi_2^2\Phi_1$, $
p_{ch,M^{root}}=\Phi_1^2\Phi_2$, 
$M$ and $M^{root}$ have a 
$2\times 2$ Jordan block with eigenvalue $-1$
respectively $1$. 
\item[(iv)]
$S(-l,2,-l)$ with $l\geq 3$ odd: $\sign I^{(0)}=(+\, 0\, -)$,
$p_{ch,M}=\Phi_2^2\Phi_1$, $p_{ch,M^{root}}=\Phi_1^2\Phi_2$,  
$M$ and $M^{root}$ have a $2\times 2$ Jordan block with 
eigenvalue $-1$ respectively $1$. 
\end{list}

(c) In all cases in this theorem the map
$Z:(\Br_3\ltimes\{\pm 1\}^3)_S\to G_\Z$ is surjective,
so $G_\Z=G_\Z^{\BB}$. 
\end{theorem}

{\bf Proof:}
(a) Recall 
\begin{eqnarray*}
H_{\Z,2}=\ker(M-\id)=\Rad I^{(1)}=\Z f_3,\ 
f_3=-\www{x}_3e_1+\www{x}_2e_2-\www{x}_1e_3.
\end{eqnarray*}
We will choose a $\Z$-basis $\uuuu{f}=(f_1,f_2)$ of 
$H_{\Z,1}:=\ker(M+\id)^2$.
Denote $\www{\uuuu{f}}:=(f_1,f_2,f_3)$. 
\index{$f_1,\ f_2,\ f_3$}
In all cases it will be easy to see that
$\www{\uuuu{f}}$ is a $\Z$-basis of $H_\Z$.
Therefore in all cases $H_\Z=H_{\Z,1}\oplus H_{\Z,2}$.

(i) $S(\HH_{1,2})$: Recall $p_{ch,M}=\Phi_2^2\Phi_1$, 
\begin{eqnarray*}
S=\begin{pmatrix}1&-2&2\\0&1&-2\\0&0&1\end{pmatrix},\ 
S^{-1}S^t=\begin{pmatrix}1&-2&2\\2&-3&2\\2&-2&1\end{pmatrix},\ 
f_3=\uuuu{e}\begin{pmatrix}1\\1\\1\end{pmatrix}.
\end{eqnarray*}
Define
\begin{eqnarray*}
f_1:=\uuuu{e}\begin{pmatrix}1\\1\\0\end{pmatrix},\ 
f_2:=\uuuu{e}\begin{pmatrix}0\\1\\1\end{pmatrix}.
\end{eqnarray*}
Then 
\begin{eqnarray*}
L(\www{\uuuu{f}}^t,\www{\uuuu{f}})^t
=\begin{pmatrix}0&-1&0\\1&0&0\\0&0&1\end{pmatrix},\
M\www{\uuuu{f}}=\www{\uuuu{f}}
\begin{pmatrix}-1&0&0\\0&-1&0\\0&0&1\end{pmatrix},\ 
H_{\Z,1}=\Z f_1\oplus \Z f_2.
\end{eqnarray*}
The claims on the groups $G_{\Z,1}^M,G_{\Z,1}^{(1)},G_{\Z,1}^{(0)}$ 
and $G_{\Z,1}$ follow from the shape of the matrices of
$M,I^{(1)},I^{(0)}$ and $L$ with respect to the basis
$\uuuu{f}$ of $H_{\Z,1}$. 

(ii) $S(-l,2,-l)$ with $l\geq 4$ even: 
 Recall $p_{ch,M}=\Phi_2^2\Phi_1$, 
\begin{eqnarray*}
S=\begin{pmatrix}1&-l&2\\0&1&-l\\0&0&1\end{pmatrix},\ 
S^{-1}S^t=\begin{pmatrix}l^2-3&-l^3+3l&l^2-2\\
l&-l^2+1&l\\2&-l&1\end{pmatrix},\ 
f_3=\uuuu{e}
\begin{pmatrix}\frac{l}{2}\\1\\ \frac{l}{2}\end{pmatrix}.
\end{eqnarray*}
Define
\begin{eqnarray*}
f_1:=\uuuu{e}\begin{pmatrix}1\\0\\-1\end{pmatrix},\ 
f_2:=\uuuu{e}
\begin{pmatrix}\frac{l^2-2}{2}\\ \frac{l}{2}\\0\end{pmatrix}.
\end{eqnarray*}
Then 
\begin{eqnarray*}
L(\www{\uuuu{f}}^t,\www{\uuuu{f}})^t
=\begin{pmatrix}0&-1&0\\1&1-\frac{l^2}{4}&0\\0&0&1\end{pmatrix},
\ M\www{\uuuu{f}}=\www{\uuuu{f}}
\begin{pmatrix}-1&2-\frac{l^2}{2}&0\\0&-1&0\\0&0&1\end{pmatrix},\\ 
H_{\Z,1}=\Z f_1\oplus \Z f_2,\ 
H_\Z=H_{\Z,1}\oplus H_{\Z,2},
\end{eqnarray*}
$M_1=M|_{H_{\Z,1}}$ and $M_1^{root}$ have each a $2\times 2$ Jordan block,
$M_1^{root}$ is cyclic and $(M_1^{root})^{l^2/2-2}=-M_1$.
Lemma \ref{t5.6} shows
\begin{eqnarray*}
G_{\Z,1}^M=G_{\Z,1}^{(0)}=G_{\Z,1}^{(1)}=G_{\Z,1}
=\{\pm (M_1^{root})^m\,|\, m\in\Z\}.
\end{eqnarray*}

(b) Recall 
\begin{eqnarray*}
\begin{array}{ccc}
 & S(\P^2) & S(A_3) \\
S^{-1}S^t & 
\begin{pmatrix}10&-15&6\\6&-8&3\\3&-3&1 \end{pmatrix} & 
\begin{pmatrix}0&0&1\\-1&0&1\\0&-1&1 \end{pmatrix} \\
p_{ch,M} & \Phi_1^3 & \Phi_4\Phi_1 \\ \hline \hline 
 & S(\widehat{A}_2) & 
 S(-l,2,-l)\textup{ with }l\geq 3\textup{ odd}\\
S^{-1}S^t & 
\begin{pmatrix}-2&-1&2\\-2&0&1\\-1&-1&1 \end{pmatrix} & 
\begin{pmatrix}l^2-3&-l^3+3l&l^2-2\\
l&-l^2+1&l\\2&-l&1 \end{pmatrix} \\
p_{ch,M} &  \Phi_2^2\Phi_1 & \Phi_2^2\Phi_1
\end{array}
\end{eqnarray*}

The case $S(-l,2,-l)$ with $l\geq 3$ odd will be treated
separately below. 
Theorem \ref{t3.26} (c) applies in the case $S(\P^2)$ with $k=1$
and in the case $S(\whh{A}_2)$ with $k=0$. It shows 
in these cases $(M^{root})^3=\varepsilon M$. 
By Theorem \ref{t3.26} (c) the matrices $M^{root,mat}$
are as claimed in the cases $S(\P^2)$ and $S(\whh{A}_2)$. 
In the case $A_3$ by definition $M^{root}=M$. 

In all three cases $S(\P^2)$, $S(A_3)$ and $S(\whh{A}_2)$ 
$M^{root}$ is cyclic with cyclic generator $e_1$.
In the cases $S(\P^2)$ and  $S(A_3)$ 
$p_{ch,M^{root}}=p_{ch,M}$. In the case
$S(\whh{A}_2)$ $p_{ch,M^{root}}=\phi_1^2\phi_2$ 
and $p_{ch,M}=\phi_2^2\phi_1$.
Lemma \ref{t5.2} (a) shows in all three cases
\begin{eqnarray*}
G_\Z^M&=& \{p(M^{root})\,|\, 
p(t)=\sum_{i=0}^2p_it^i\in \Z[t],\\
&&p(\kappa)\in(\Z[\kappa])^*
\textup{ for each eigenvalue }\kappa\textup{ of }M^{root}\}.
\end{eqnarray*}

(i) $S(\P^2)$: $M^{root}-\id$ is nilpotent with 
$(M^{root}-\id)^2\neq 0$, $(M^{root}-\id)^3=0$.
An element of $\Z[M^{root}]$ can be written in the form
\begin{eqnarray*}
q_0\id+q_1(M^{root}-\id)+q_2(M^{root}-\id)^2\textup{ with }
q_0,q_1,q_2\in\Z.
\end{eqnarray*}
It is in $\Aut(H_\Z)$ if and only if $q_0\in\{\pm 1\}$. 
Then it can be written as
\begin{eqnarray*}
q_0(M^{root})^{q_0q_1}(\id + \www{q}_2(M^{root}-\id)^2)
\textup{ for some }\www{q}_2\in\Z.
\end{eqnarray*}
Therefore $G_\Z^M$ is as claimed. 

Because $(M^{root}-\id)^3=0$, 
$(M^{root}-\id)^2(H_\Z)\subset\ker (M^{root}-\id)
=\Rad I^{(1)}$.
Therefore $\id+\www{q}_2(M^{root}-\id)^2$ and thus any
element of $G_\Z^M$ respects $I^{(1)}$,
so $G_\Z^M=G_\Z^{(1)}$.

On the other hand, one easily checks that 
$\id+\www{q}_2(M^{root}-\id)^2$ respects $I^{(0)}$
only if $\www{q}_2=0$. Therefore
\begin{eqnarray*}
G_\Z=G_\Z^{(0)}=\{\pm (M^{root})^m\,|\, m\in\Z\}
\subsetneqq G_\Z^{(1)}=G_\Z^M.
\end{eqnarray*}

(ii) $S(A_3)$: For $p(t)=\sum_{i=0}^2p_it^i\in\Z[t]$ write 
$\mu_j:=p(\lambda_j)$ for $j\in\{1,2,3\}$
for the eigenvalues of the element
$p(M)\in\Z[M]=\End(H_\Z,M)$ where 
$\lambda_1=i,\lambda_2=-i,\lambda_3=1$ are the eigenvalues
of $M$. Because of $(\Z[i])^*=\{\pm 1,\pm i\}$, 
one can multiply a given element of $G_\Z^M$ with a 
suitable power of $M$ and obtain an element
with $\mu_1=\mu_2=1$. Therefore
\begin{eqnarray*}
G_\Z^M=\{M^{m_1}(\id+m_2\Phi_4(M))&|& m_1\in\{0,1,2,3\},
m_2\in\Z, \\
&& 1+m_2\Phi_4(1)\in\{\pm 1\}\}.
\end{eqnarray*}
This forces $m_2\in\{0,-1\}$. The case $m_2=-1$ gives
$\id+(-1)(M^2+\id)=-M^2$. Therefore
\begin{eqnarray*}
\{\pm M^m\,|\, m\in\{0,1,2,3\}\}
=\{\pm M^m\,|\, m\in\Z\}
=G_\Z^M=G_\Z^{(0)}=G_\Z^{(1)}=G_\Z.
\end{eqnarray*}

(iii) $S(\widehat{A}_2)$: Write
$\uuuu{f}=(f_1,f_2):=\uuuu{e}
\begin{pmatrix}1&2\\1&1\\1&0\end{pmatrix}$. Then
\index{$f_1,\ f_2,\ f_3$} 
\begin{eqnarray*}
M^{root}\uuuu{f}&=&\uuuu{f}
\begin{pmatrix}1&1\\0&1\end{pmatrix},\\
H_{\Z,1}&=&\ker\Phi_2^2(M)=\ker\Phi_1^2(M^{root})
=\Z f_1\oplus \Z f_2.
\end{eqnarray*}
Lemma \ref{t5.6} implies 
\begin{eqnarray*}
\{g|_{H_{\Z,1}}\,|\, g\in G_\Z^M\}
=\{\pm (M^{root}|_{H_{\Z,1}})^m\,|\, m\in\Z\},\\
G_\Z^M=\{\pm (M^{root})^m\,|\, m\in\Z\}\times 
\{g\in G_\Z^M\,|\, g|_{H_{\Z,1}}=\id\}.
\end{eqnarray*}
But $p(t)=1+q\Phi_1^2(t)$ with $q\in\Z$ satisfies 
$p(-1)=1+q\cdot 2^2\in\{\pm 1\}$ only if $q=0$. Therefore
\begin{eqnarray*}
\{\pm (M^{root})^m\,|\, m\in\Z\}
=G_\Z^M=G_\Z^{(0)}=G_\Z^{(1)}=G_\Z.
\end{eqnarray*} 

(iv) $S(-l,2,-l)$ with $l\geq 3$ odd: Recall $f_3$ and define
$(f_1,\www{f}_2)$ and $\uuuu{d}=(d_1,d_2,d_3)$ with
\begin{eqnarray*}
(f_1,\www{f}_2,f_3)=\uuuu{e}\begin{pmatrix}
1&l^2-2&l\\0&l&2\\-1&0&l\end{pmatrix},\quad
\uuuu{d}=\uuuu{e}\begin{pmatrix}
\frac{l^2-l-2}{2}&\frac{l^2-1}{2}&\frac{l^2-l}{2}\\
\frac{l-1}{2}&\frac{l+1}{2}&\frac{l-1}{2}\\
0&\frac{l-1}{2}&-1\end{pmatrix}.
\end{eqnarray*}
The matrix which expresses $(f_1,\www{f}_2,f_3)$ with
$\uuuu{e}$ has determinant 4, the matrix which expresses
$\uuuu{d}$ with $\uuuu{e}$ has determinant 1.
Therefore $\Z f_1\oplus\Z \www{f}_2\oplus \Z f_3$ is a 
sublattice of index 4 in $H_\Z$, and $\uuuu{d}$ is
a $\Z$-basis of $H_\Z$. 
One calculates
\begin{eqnarray*}
M(f_1,\www{f}_2,f_3)&=&(f_1,\www{f}_2,f_3)\begin{pmatrix}
-1&-l^2+4&0\\0&-1&0\\0&0&1\end{pmatrix}.
\end{eqnarray*}
Especially 
\begin{eqnarray*}
\Z f_1\oplus \Z \www{f}_2=\ker(M^2+\id)^2\supset\Z f_1
=\ker(M+\id)=\Rad I^{(0)}.
\end{eqnarray*}
Observe  
\begin{eqnarray}
\delta_3\sigma_1^{-1}\sigma_2^{-1}\sigma_1\in 
(\Br_3\ltimes\{\pm 1\}^3)_S.\label{5.13}
\end{eqnarray}
Define
\begin{eqnarray}\label{5.14}
\www{M}:=Z(\delta_3\sigma_1^{-1}\sigma_2^{-1}\sigma_1)\in G_\Z.
\end{eqnarray}
Then $M^{root}=(-M)\circ \www{M}^{(5-l^2)/2}\in G_\Z$. 
One calculates
\begin{eqnarray}
\www{M}\uuuu{e}&=&\uuuu{e}+f_1(-1,l,-1),\label{5.15}\\
\www{M}(f_1,\www{f}_2,f_3)&=&(f_1,\www{f}_2,f_3)\begin{pmatrix}
1&2&0\\0&1&0\\0&0&1\end{pmatrix},\label{5.16}\\
M^{root}(f_1,\www{f}_2,f_3)&=&(f_1,\www{f}_2,f_3)\begin{pmatrix}
1&1&0\\0&1&0\\0&0&-1\end{pmatrix},\nonumber\\
(M^{root})^{l^2-4}&=&-M,\nonumber\\
(M^{root})^2&=& \www{M}.\nonumber 
\end{eqnarray}

Finally, one calculates
\begin{eqnarray*}
M^{root}\uuuu{d}=\uuuu{d}\begin{pmatrix}
0&0&-1\\1&0&1\\0&1&1\end{pmatrix}.
\end{eqnarray*}
Therefore $M^{root}$ is cyclic with cyclic generator $d_1$
and regular. Lemma \ref{t5.2} (a) shows
\begin{eqnarray*}
G_\Z^M=\{p(M^{root})\,|\, p(t)=\sum_{i=0}^2p_it^i\in\Z[t],
p(1),p(-1)\in\{\pm 1\}\}.
\end{eqnarray*}
As in the case $S(\widehat{A}_2)$ one finds with Lemma \ref{t5.6}
\begin{eqnarray*}
\{g|_{H_{\Z,1}}\,|\, g\in G_\Z^M\}
=\{\pm (M^{root}|_{H_{\Z,1}})^m\,|\, m\in\Z\},\\
G_\Z^M=\{\pm (M^{root})^m\,|\, m\in\Z\}\times 
\{g\in G_\Z^M\,|\, g|_{H_{\Z,1}}=\id\}.
\end{eqnarray*}
But $p(t)=1+q\Phi_1^2(t)$ with $q\in\Z$ satisfies 
$p(-1)=1+q\cdot 2^2\in\{\pm 1\}$ only if $q=0$. Therefore
\begin{eqnarray*}
\{\pm (M^{root})^m\,|\, m\in\Z\}
=G_\Z^M=G_\Z^{(0)}=G_\Z^{(1)}=G_\Z.
\end{eqnarray*} 

(c) Of course $Z(\delta_1\delta_2\delta_3)=-\id$. In the cases
in part (b) $G_\Z=\{\pm (M^{root})^m\,|\, m\in\Z\}$.
The definitions \eqref{5.9}--\eqref{5.12} show that $Z$ is
surjective.

The case $\HH_{1,2}$: $G_\Z\cong SL_2(\Z)\times\{\pm 1\}$. 
The group $G_\Z$ is generated by $-\id$, $h_1$ and $h_2$ with
\begin{eqnarray*}
h_1:=(\www{\uuuu{f}}\mapsto \www{\uuuu{f}}
\begin{pmatrix}1&-1&0\\0&1&0\\0&0&1\end{pmatrix})
=(\uuuu{e}\mapsto \uuuu{e}
\begin{pmatrix}2&-1&0\\1&0&0\\0&0&1\end{pmatrix})
=Z(\delta_2\sigma_1),\\
h_2:=(\www{\uuuu{f}}\mapsto \www{\uuuu{f}}
\begin{pmatrix}1&0&0\\1&1&0\\0&0&1\end{pmatrix})
=(\uuuu{e}\mapsto \uuuu{e}
\begin{pmatrix}1&0&0\\0&2&-1\\0&1&0\end{pmatrix})
=Z(\delta_3\sigma_2).
\end{eqnarray*}

The cases $S(-l,2,-l)$ with $l\geq 4$ even: 
\eqref{5.13}--\eqref{5.16} hold also for
even $l$. With respect to the $\Z$-basis $\www{\uuuu{f}}=(f_1,f_2,f_3)
=(f_1,\frac{1}{2}\www{f}_2,f_3)$ of $H_\Z$ 
\begin{eqnarray*}
\www{M}\www{\uuuu{f}}=\www{\uuuu{f}}
\begin{pmatrix}1&1&0\\0&1&0\\0&0&1\end{pmatrix}.
\end{eqnarray*}
One sees
\begin{eqnarray*}
G_\Z&=&\langle -\id,\www{M},Q\rangle
\end{eqnarray*}
with 
\begin{eqnarray*}
Q&=&(\www{\uuuu{f}}\mapsto \www{\uuuu{f}}
\begin{pmatrix}-1&0&0\\0&-1&0\\0&0&1\end{pmatrix})
= M\circ \www{M}^{2-l^2/2}\\
&=& Z(\sigma^{mon})\circ 
Z(\delta_3\sigma_1^{-1}\sigma_2^{-1}\sigma_1)^{2-l^2/2}.
\hspace*{2cm}\Box
\end{eqnarray*}

\hfill$\Box$

\section{Special rank 3 cases with eigenvalues not all in $S^1$}
\label{s5.6}

This section starts with a general lemma for all rank 3 cases
with eigenvalues not all in $S^1$. It gives coarse information
on the four groups $G_\Z^M,G_\Z^{(0)},G_\Z^{(1)}$ and $G_\Z$. 
Afterwards Theorem \ref{t5.16} determines these groups
precisely for three series of cases and one exceptional case.
Theorem \ref{t5.18} in section \ref{s5.7} will treat the
other irreducible cases with eigenvalues not all in $S^1$.

\begin{lemma}\label{t5.15}
Fix $\uuuu{x}\in\Z^3$ with $r(\uuuu{x})<0$ or $r(\uuuu{x})>4$. Then
\begin{eqnarray*}
G_\Z^M\stackrel{(2:1)\textup{ or }(1:1)}
{\supset} G_\Z=G_\Z^{(0)}=G_\Z^{(1)}
\stackrel{(\textup{finite}:1)}{\supset}
\{\pm M^l\,|\, l\in\Z\}.
\end{eqnarray*}
\end{lemma}

{\bf Proof:}
For $p(t)=\sum_{i=0}^2p_it^i\in\Q[t]$ write
$\mu_j:=p(\lambda_j)$ for $j\in\{1,2,3\}$.
Then $\mu_1,\mu_2$ and $\mu_3$ are the eigenvalues of
$p(M)\in\End(H_\Q)$. The monodromy is because of
$r(\uuuu{x})\in \Z-\{0,1,2,3,4\}$ and Lemma \ref{t5.7}
semisimple and regular. Lemma \ref{t5.1} (c) (i) and (ii)
applies and gives   
\begin{eqnarray}\label{5.17}
G_\Z^M&=&\{p(M)\,|\, p(t)=\sum_{i=0}^2p_it^i\in\Q[t],\\
&&p(M)\in\End(H_\Z),
\mu_1\mu_2\in\{\pm 1\}, \mu_3\in\{\pm 1\}\}\nonumber\\
\supset G_\Z &=& G_\Z^{(0)}=G_\Z^{(1)}=\{p(M)\in G_\Z^M\,|\, 
\mu_1\mu_2=1\}\label{5.18}
\end{eqnarray}
Especially $h\in G_\Z^M\Rightarrow h^2\in G_\Z$. 
Also, $\End(H_\Q,M)=\Q[M]$, and the map
\begin{eqnarray*}
\End(H_\Q,M)=\Q[M]\to \Q[\lambda_1]\times\Q,\quad 
p(M)\mapsto (\mu_1,\mu_3)
\end{eqnarray*}
is an isomorphism of $\Q$-algebras. This is a special
case of the chinese remainder theorem. 
$Q$ is mapped to $(-1,1)$. 

Observe $(-\id)\in G_\Z\subset G_\Z^M$. Therefore the
subgroup $\{h\in G_\Z\,|\, \mu_3=1\}$ has index 2 in
$G_\Z$, and the subgroup $\{h\in G_\Z^M\,|\, \mu_3=1\}$
has index 2 in $G_\Z^M$. The map
\begin{eqnarray}
\{h\in G_\Z^M\,|\, \mu_3=1\}\to \OO_{\Q[\lambda_1]}^*,\quad
h\mapsto \mu_1,\label{5.19}
\end{eqnarray}
is injective. The element $-1\in \OO_{\Q[\lambda_1]}^*$ 
is in the image of the map in \eqref{5.19} if and only if 
$Q\in G_\Z^M$, and then it is the image of $Q$. 
By Dirichlet's unit theorem \cite[Ch. 2 4.3 Theorem 5]{BSh73}
the group $\OO_{\Q[\lambda_1]}^*$ is isomorphic to 
the group $\{\pm 1\}\times \Z$.
Therefore $\{h\in G_\Z^M\,|\, \mu_3=1\}$ is isomorphic
to $\{\pm 1\}\times \Z$ if $Q\in G_\Z^M$ and to $\Z$ if $Q\notin G_\Z^M$. 
If $Q\in G_\Z^M$ then because of $Q\in \Aut(H_\Q,L)$ 
also $Q\in G_\Z$. 
This and the implication
$h\in G_\Z^M\Rightarrow h^2\in G_\Z$ show
\begin{eqnarray*}
[G_\Z^M:G_\Z]
=[\{h\in G_\Z^M\,|\, \mu_3=1\}:\{h\in G_\Z\,|\, \mu_3=1\}]
\in\{1,2\}.
\end{eqnarray*}
The group $\{\pm M^l\,|\, l\in\Z\}\subset G_\Z
\subset G_\Z^M$ is isomorphic to $\{\pm 1\}\times\Z$,
so it has a free part of rank 1, just as $G_\Z$ and $G_\Z^M$. 
Therefore
$[G_\Z:\{\pm M^l\,|\, l\in\Z\}]<\infty$.
\hfill$\Box$

\bigskip
Later we will see precisely how much bigger 
$G_\Z^M$ and $G_\Z$ are than $\{\pm M^l\,|\, l\in\Z\}$.
In the majority of the cases they are not bigger, but 
$G_\Z^M=G_\Z=\{\pm M^l\,|\, l\in\Z\}$.

The following theorem determines the groups 
$G_\Z$ and $G_\Z^M\supset G_\Z$ for three series of triples
and the exceptional case $(3,3,4)$. 
In Theorem \ref{t5.18} in section \ref{s5.7} 
we will see that the $\Br_3\ltimes\{\pm 1\}^3$
orbits of these three series and of the triple $(3,3,4)$ are
the only triples $\uuuu{x}$ with $r(\uuuu{x})\in\Z_{<0}\cup\Z_{>4}$,
$(H_\Z,L,\uuuu{e})$ irreducible 
and {\it not} $G_\Z^M=G_\Z=\{\pm M^l\,|\, l\in\Z\}$.

\begin{theorem}\label{t5.16}
For each $\uuuu{x}\in\Z^3$ below fix also the associated triple
$(H_\Z,L,\uuuu{e})$.

(a) Consider $\uuuu{x}=(x,x,x)$ with $x\in\Z-\{-1,0,1,2,3\}$
and $S=S(\uuuu{x})$. 
Then $\delta_3\sigma_2\sigma_1\in (\Br_3\ltimes\{\pm 1\}^3)_S$
and by Theorem \ref{t3.26} (c) (with $k=0$) 
\begin{eqnarray*}
M^{root,3}:=Z(\delta_3\sigma_2\sigma_1)\in G_\Z
\quad\textup{with}\quad (M^{root,3})^3=-M.
\end{eqnarray*}
$M^{root,3}$ is cyclic with $M^{root,3}(f_3)=-f_3$.
In the case $x=4$ define also 
\begin{eqnarray*}
M^{root,6}:=-(M^{root,3})^2-2M^{root,3}.
\end{eqnarray*}
Then
\begin{eqnarray*}
G_\Z&=&G_\Z^M=\{\pm (M^{root,3})^l\,|\, l\in\Z\}\quad\textup{if }
x\notin\{4,5\},\\
G_\Z&=&G_\Z^M=\{\id,Q\}\times \{\pm (M^{root,3})^l\,|\, l\in\Z\}\quad\textup{if }x=5,\\
G_\Z&=& G_\Z^{(0)}=G_\Z^{(1)}=
\{\id,Q\}\times \{\pm (M^{root,3})^l\,|\, l\in\Z\}\\
&\stackrel{1:2}{\subset}&
G_\Z^M=\{\id,Q\}\times \{\pm (M^{root,6})^l\,|\, l\in\Z\}
\quad\textup{ if }x=4.
\end{eqnarray*}
If $x=4$ then $(M^{root,6})^2=-M^{root,3}$. 

(b) Consider $\uuuu{x}=(2y,2y,2y^2)$ with $y\in\Z_{\geq 2}$ 
and $S=S(\uuuu{x})$. 
Then $\sigma_2\sigma_1^2\in (\Br_3\ltimes\{\pm 1\}^3)_S$.
By Lemma \ref{t3.25} (b) 
\begin{eqnarray*}
M^{root,2}:=Z(\sigma_2\sigma_1^2)\in G_\Z.
\end{eqnarray*}
It satisfies $(M^{root,2})^2=M$ and $M^{root,2}(f_3)=-f_3$.
In the case $y=2$ define also 
\begin{eqnarray*}
M^{root,4}:=-\frac{1}{4}M-2M^{root,2}-\frac{3}{4}\id.
\end{eqnarray*}
Then
\begin{eqnarray*}
G_\Z&=&G_\Z^M=\{\pm (M^{root,2})^l\,|\, l\in\Z\}
\quad\textup{if }
y\geq 3,\\
G_\Z&=&G_\Z^{(0)}=G_\Z^{(1)}=
\{\id,Q\}\times \{\pm (M^{root,2})^l\,|\, l\in\Z\}\\
&\stackrel{1:2}{\subset}&
G_\Z^M=\{\id,Q\}\times \{\pm (M^{root,4})^l\,|\, l\in\Z\}
\quad\textup{ if }y=2.
\end{eqnarray*}
If $y=2$ then $(M^{root,4})^2=-M^{root,2}$. 

(c) Consider $\uuuu{x}=(x,x,0)$ with $x\in\Z_{\geq 2}$
and $S=S(\uuuu{x})$. 
In the case $x=2$ define also 
\begin{eqnarray*}
M^{root,2}:=\frac{1}{2}M+\frac{1}{2}\id.
\end{eqnarray*}
Then
\begin{eqnarray*}
G_\Z&=&G_\Z^M=\{\id,Q\}\times\{\pm M^l\,|\, l\in\Z\}\quad\textup{if }
x\geq 3,\\
G_\Z&=&G_\Z^{(0)}=G_\Z^{(1)}=\{\id,Q\}\times \{\pm M^l\,|\, l\in\Z\}\\
&\stackrel{1:2}{\subset}&
G_\Z^M=\{\id,Q\}\times \{\pm (M^{root,2})^l\,|\, l\in\Z\}
\quad\textup{ if }x=2.
\end{eqnarray*}
If $x=2$ then $(M^{root,2})^2=M$. 

(d) Consider $\uuuu{x}=(3,3,4)$. Then
\begin{eqnarray*}
G_\Z=G_\Z^M=\{\id,Q\}\times\{\pm M^l\,|\, l\in\Z\}.
\end{eqnarray*}

(e) In all cases in (a)--(d) except for the four cases
$$\uuuu{x}\in\{(4,4,4),(5,5,5),(4,4,8),(3,3,4)\}$$
the map $Z:(\Br_3\ltimes\{\pm 1\}^3)_S\to G_\Z$ is surjective
so $G_\Z=G_\Z^{\BB}$. 
In the four exceptional cases $Q\in G_\Z-G_\Z^{\BB}$. 
\end{theorem}

{\bf Proof:} 
In all cases in this theorem $r(\uuuu{x})<0$ or $r(\uuuu{x})>4$,
so Lemma \ref{t5.15} applies, so
$G_\Z=G_\Z^{(0)}=G_\Z^{(1)}$. 

(a) $S$ is as in Theorem \ref{t3.26} (c) with $k=0$. 
Therefore $\delta_3\sigma_2\sigma_1$
is in the stabilizer of $S$ and $M^{root,3}$ is in $G_\Z$,
it is cyclic, and it satisfies $(M^{root,3})^3=-M$. 
Explicitly (by Theorem \ref{t3.26} (a))
$$ M^{root,3}(\uuuu{e})=\uuuu{e}\cdot 
\begin{pmatrix}-x&-x&-1\\1&0&0\\0&1&0\end{pmatrix}.$$
One sees $M^{root,3}(f_3)=-f_3$ where $f_3=-e_1+e_2-e_3$, 
so its third eigenvalue is $\kappa_3=-1$. 
The other two eigenvalues $\kappa_1$ and $\kappa_2$ are determined
by the trace $-x=\kappa_1+\kappa_2-1$ and the determinant 
$-1=\kappa_1\kappa_2(-1)$ of $M^{root,3}$.
The eigenvalues are 
$$\kappa_{1/2}=\frac{1-x}{2}\pm \frac{1}{2}\sqrt{x^2-2x-3},\quad
\kappa_3=-1.$$

Because $M$ and $M^{root,3}$ are regular, Lemma \ref{t5.1} 
applies. It gives an isomorphism of $\Q$-algebras
\begin{eqnarray*}
\End(H_\Q,M)=\End(H_\Q,M^{root,3})&&\\
=\{p(M^{root,3})\,|\, p(t)=\sum_{i=0}^2p_it^i\in\Q[t]\}
&\to& \Q[\kappa_1]\times \Q\\
p(M^{root,3})&\mapsto& (p(\kappa_1),p(-1)).
\end{eqnarray*}
The image of $-\id$ is $(-1,-1)$, 
the image of $Q$ is $(-1,1)$, 
the image of $M^{root,3}$ is $(\kappa_1,-1)$. 
The image of $G_\Z^M$ is a priori 
a subgroup of
$\OO_{\Q[\kappa_1]}^*\times\{\pm 1\}$. 
We have to find out which one. By Theorem \ref{t5.11}
$Q\in G_\Z^M$ (and then also $Q\in G_\Z$) only for $x\in\{4,5\}$.

Lemma \ref{t5.2} applies because $M^{root,3}$ is cyclic.
It shows
\begin{eqnarray*}
G_\Z^M&=& \{p(M^{root,3})\,|\, p(t)=p_2t^2+p_1t+p_0\in\Z[t]
\textup{ with }\\
&& (p(\kappa_1),p(-1))\in\Z[\kappa_1]^*\times\{\pm 1\}\},\\
G_\Z&=& \{p(M^{root,3})\,|\, p(t)=p_2t^2+p_1t+p_0\in\Z[t]
\textup{ with }\\
&& (p(\kappa_1),p(-1))\in\Z[\kappa_1]^*\times\{\pm 1\}
\textup{ and }p(\kappa_1)p(\kappa_2)=1\}.
\end{eqnarray*}
Now Lemma \ref{tc.1} (a) is useful. It says 
\begin{eqnarray*}
\Z[\kappa_1]^*=\left\{\begin{array}{ll}
\{\pm \kappa_1^l\,|\, l\in\Z\}&\textup{ for }x\notin\{4,-2\},\\
\{\pm (\kappa_1-1)^l\,|\, l\in\Z\}&\textup{ for }x=-2,\\
\{\pm (\kappa_1+1)^l\,|\, l\in\Z\}&\textup{ for }x=4.
\end{array}\right.
\end{eqnarray*}
with 
\begin{eqnarray*}
(\kappa_1-1)(\kappa_2-1)=-1\quad\textup{and}\quad 
(\kappa_1-1)^2=\kappa_1\quad\textup{for}\quad x=-2,\\
(\kappa_1+1)(\kappa_2+1)=-1\quad\textup{and}\quad 
(\kappa_1+1)^2=-\kappa_1\quad\textup{for}\quad x=4.
\end{eqnarray*}

For $x\notin\{4,5,-2\}$ the facts $-\id,M^{root,3}\in G_\Z$
and $Q\notin G_\Z$ show that the image of $G_\Z$ and $G_\Z^M$
in $\Z[\kappa_1]^*\times\{\pm 1\}=\{\pm \kappa_1^l\,|\, l\in\Z\}
\times\{\pm 1\}$ has index 2. Therefore
then $G_\Z=G_\Z^M$ is as claimed. 

For $x=5$ the facts $-\id,M^{root,3},Q\in G_\Z$ show that the 
image of $G_\Z$ and $G_\Z^M$ is $\Z[\kappa_1]^*\times\{\pm 1\}
= \{\pm \kappa_1^l\,|\, l\in\Z\} \times\{\pm 1\}$. Therefore
then $G_\Z=G_\Z^M$ is as claimed.

Consider the case $x=-2$. If an automorphism $p(M^{root,3})$
is in $G_\Z^M$ which corresponds to a pair $(\kappa_1-1,\pm 1)$,
then $p(\kappa_1)=\kappa_1-1$ means
$p(t)=t-1+l_2(t^2-3t+1)$ for some $l_2\in\Z$.
But then $p(-1)=-2+l_2\cdot 5\notin\{\pm 1\}$.
So $G_\Z^M$ does not contain such an automorphism.
$G_\Z^M$ and $G_\Z$ are as claimed, because $Q\notin G_\Z^M$.

Consider the case $x=4$. The polynomial
$p(t)=-t^2-2t$ satisfies $p(-1)=1$, 
$p(\kappa_1)=-\kappa_1^2-2\kappa_1=-(-3\kappa_1-1)-2\kappa_1
=\kappa_1+1$. Therefore $M^{root,6}=p(M^{root,3})\in G_\Z^M$.
Because $\kappa_1+1$ has norm $-1$, $M^{root,6}$ is not in
$G_\Z$. The groups $G_\Z^M$ and $G_\Z$ are as claimed. 

(b) Compare the sections \ref{s4.1} and \ref{s3.2} 
for the actions of $\Br_3\ltimes\{\pm 1\}^3$ on 
$\Z^3$ and on $\BB^{dist}$.
\begin{eqnarray*}
&&\sigma_2\sigma_1^2(2y,2y,2y^2)\\
&=& \sigma_2\sigma_1(-2y,2y^2-2y\cdot 2y,2y)
=\sigma_2\sigma_1(-2y,-2y^2,2y)\\
&=& \sigma_2(2y,2y-(-2y)(-2y^2),-2y^2)
=\sigma_2(2y,2y-4y^3,-2y^2)\\
&=& ((2y-4y^3)-2y\cdot (-2y^2),2y,2y^2)
=(2y,2y,2y^2),
\end{eqnarray*}
so $\sigma_2\sigma_1^2$ is in the stabilizer of $(2y,2y,2y^2)$,
so $M^{root,2}:=Z(\sigma_2\sigma_1^2)\in G_\Z$. 
\begin{eqnarray*}
&& M^{root,2}(\uuuu{e})=\sigma_2\sigma_1^2(\uuuu{e})\\
&=& \sigma_2\sigma_1(s_{e_1}^{(0)}(e_2),e_1,e_3)\\
&=&\sigma_2\sigma_1(e_2-2ye_1,e_1,e_3)\\
&=& \sigma_2(s_{e_2-2ye_1}^{(0)}(e_1),e_2-2ye_1,e_3)\\
&=&\sigma_2(e_1+2y(e_2-2ye_1),e_2-2ye_1,e_3)\\
&=& ((1-4y^2)e_1+2ye_2,s_{e_2-2ye_1}^{(0)}(e_3),e_2-2ye_1)\\
&=&((1-4y^2)e_1+2ye_2,e_3+2y^2(e_2-2ye_1),e_2-2ye_1)\\
&=& \uuuu{e}\cdot 
\begin{pmatrix} 1-4y^2&-4y^3&-2y\\2y&2y^2&1\\0&1&0\end{pmatrix}
=:\uuuu{e}\cdot M^{root,2,mat}.
\end{eqnarray*}
The map $Z:(\Br_3\ltimes\{\pm 1\}^3)_S\to G_\Z$ in Lemma 
\ref{t3.25} is a group antihomomorphism. By Theorem \ref{t3.26}
$Z(\sigma^{mon})=M$. Therefore 
\begin{eqnarray*}
(M^{root,2})^2&=& Z(\sigma_2\sigma_1^2)Z(\sigma_2\sigma_1^2)
=Z(\sigma_2\sigma_1(\sigma_1\sigma_2\sigma_1)\sigma_1)\\
&=& Z(\sigma_2\sigma_1(\sigma_2\sigma_1\sigma_2)\sigma_1)
=Z((\sigma_2\sigma_1)^3)=Z(\sigma^{mon})=M.
\end{eqnarray*}
One sees $M^{root,2}(f_3)=-f_3$, where $f_3=-ye_1+e_2-e_3$,
so its third eigenvalue is $\kappa_3=-1$. The other two
eigenvalues $\kappa_1$ and $\kappa_2$ are determined by
the trace $1-2y^2=\kappa_1+\kappa_2-1$ and the product
$\kappa_1\kappa_2=1$, which holds because of $M^{root,2}\in G_\Z$
(or because $\det M^{root,2}=-1$). The eigenvalues are
$$\kappa_{1,2}=(1-y^2)\pm y\sqrt{y^2-2},\quad \kappa_3=-1.$$

Because $M$ and $M^{root,2}$ are regular, Lemma \ref{t5.1} 
applies. It gives an isomorphism of
$\Q$-algebras
\begin{eqnarray*}
\End(H_\Q,M)=\End(H_\Q,M^{root,2})&&\\
=\{p(M^{root,2})\,|\, p(t)=\sum_{i=0}^2p_it^i\in\Q[t]\}
&\to& \Q[\kappa_1]\times \Q\\
p(M^{root,2})&\mapsto& (p(\kappa_1),p(-1)).
\end{eqnarray*}
The image of $-\id$ is $(-1,-1)$, 
the image of $Q$ is $(-1,1)$, 
the image of $M^{root,2}$ is $(\kappa_1,-1)$. 
The image of $G_\Z^M$ is a priori 
a subgroup of
$\OO_{\Q[\kappa_1]}^*\times\{\pm 1\}$. 
We have to find out which one. By Theorem \ref{t5.11}
$Q\in G_\Z^M$ (and then also $Q\in G_\Z$) only for $y=2$.

Consider the decomposition $H_\Q=H_{\Q,1}\oplus H_{\Q,2}$
as in Definition \ref{t5.9} and the primitive sublattices
$H_{\Z,1}=H_{\Q,1}\cap H_\Z$ and 
$H_{\Z,2}=H_{\Q,2}\cap H_\Z=\Z f_3$ in $H_\Z$.  
The sublattice $H_{\Z,1}$ is the right $L$-orthogonal subspace 
$(\Z f_3)^\perp\subset H_\Z$ of $\Z f_3$, see \eqref{5.8}:
\begin{eqnarray*}
H_{\Z,1}&=&\{\uuuu{e}\cdot \uuuu{z}^t\,|\, \uuuu{z}\in \Z^3,
0=(-y,1-2y^2,-1)\begin{pmatrix}z_1\\z_2\\z_3\end{pmatrix}\}\\
&=& \langle f_1,f_2\rangle_\Z
\quad\textup{with}\quad 
f_1=e_1-ye_3,\quad f_2=e_2+(1-2y^2)e_3.
\end{eqnarray*}
Write $\uuuu{f}=(f_1,f_2,f_3)=\uuuu{e}\cdot M(\uuuu{e},\uuuu{f})$, 
\begin{eqnarray*}
M^{root,2}(\uuuu{e})=\uuuu{e}\cdot M^{root,2,mat},\quad 
M^{root,2}(\uuuu{f})=\uuuu{f}M^{root,2,mat,\uuuu{f}}.
\end{eqnarray*}
Then 
\begin{eqnarray*}
M(\uuuu{e},\uuuu{f})
&=&\begin{pmatrix}1&0&-y\\0&1&1\\-y&1-2y^2&-1\end{pmatrix},\\
M(\uuuu{e},\uuuu{f})^{-1}
&=&\frac{1}{y^2-2}
\begin{pmatrix}2y^2-2&2y^3-y&y\\-y&-y^2-1&-1\\y&2y^2-1&1
\end{pmatrix},
\end{eqnarray*}
\begin{eqnarray*}
M^{root,2,mat,\uuuu{f}}
=M(\uuuu{e},\uuuu{f})^{-1}M^{root,2,mat}M(\uuuu{e},\uuuu{f})
= \begin{pmatrix} -2y^2+1&-2y&0\\y&1&0\\0&0&-1\end{pmatrix}.
\end{eqnarray*}

Any element $h=p(M)\in G_\Z^M$ with $p(t)\in\Q[t]$ 
restricts to an automorphism of
$H_{\Z,1}$ which commutes with $M^{root,2}|_{H_{\Z,1}}$, so its restriction
to $H_{\Z,1}$ has the shape 
$h|_{H_{\Z,1}}=a\id+b M^{root,2}|_{H_{\Z,1}}$ with $a,b\in\Q$ with 
\begin{eqnarray*}
a\begin{pmatrix}1&0\\0&1\end{pmatrix}
+b\begin{pmatrix}-2y^2+1&-2y\\y&1\end{pmatrix}\in GL_2(\Z).
\end{eqnarray*}
This implies $by\in\Z$ and $a+b\in\Z$. The eigenvalue 
$p(\kappa_1)$ is 
\begin{eqnarray*}
&&a+b\kappa_1=a+b(1-y^2+y\sqrt{y^2-2})\\
&=&(a+b)+by(-y+\sqrt{y^2-2})\in\Z[\sqrt{y^2-2}]^*.
\end{eqnarray*}
Now Lemma \ref{tc.1} (b) is useful. It says 
\begin{eqnarray*}
\Z[\sqrt{y^2-2}]^*=\left\{\begin{array}{ll}
\{\pm \kappa_1^l\,|\, l\in\Z\}&\textup{ if }y\geq 3,\\
\{\pm (1+\sqrt{2})^l\,|\, l\in\Z\}&\textup{ if }y=2.
\end{array}\right.
\end{eqnarray*}
Furthermore $\kappa_1$ has norm 1, $1+\sqrt{2}$ has norm $-1$,
and $(1+\sqrt{2})^2=-\kappa_2$ if $y=2$. 

Consider the cases $y\geq 3$. The map 
$$G_\Z^M\to \OO_{\Q[\kappa_1]}^*\times\{\pm 1\},
\quad p(M^{root,2})\mapsto (p(\kappa_1),p(-1)),$$
has because of $Q\notin G_\Z^M$ as image the index 2 subgroup 
of $\{\pm\kappa_1^l\,|\, l\in\Z\}\times\{\pm 1\}$ which is generated
by $(\kappa_1,-1)$ and $(-1,-1)$, because $Q\notin G_\Z^M$.
Therefore then $G_\Z^M=G_\Z=\{\pm (M^{root,2})^l\,|\, l\in\Z\}$. 

Consider the case $y=2$. Then $Q\in G_\Z\subset G_\Z^M$.
Therefore then 
$G_\Z=\{\id,Q\}\times\{\pm (M^{root,2})^l\,|\, l\in\Z\}$. 
The question remains whether $1+\sqrt{2}$ arises as
eigenvalue $p(\kappa_1)$ for an element $p(M^{root,2})\in G_\Z^M$.
It does. $M^{root,4}$ has the first eigenvalue
$$-\frac{1}{4}(-3+2\sqrt{2})^2-2(-3+2\sqrt{2})-\frac{3}{4}
=1-\sqrt{2}$$
with $(1-\sqrt{2})^2=3-2\sqrt{2}=-\kappa_1$
and the third eigenvalue $-\frac{1}{4}-2(-1)-\frac{3}{4}=1$. 
Therefore $(M^{root,4})^2=-M^{root,2}$.
$M^{root,4}$ is in $G_\Z^M$ because
\begin{eqnarray*}
M^{root,4}(\uuuu{e})
&=&\uuuu{e}\left(-\frac{1}{4}
\begin{pmatrix} 97 & 220 & 28 \\ -28 & -63 & -8 \\ 4 & 8 & 1\end{pmatrix}
-2\begin{pmatrix} -15 & -32 & -4 \\ 4 & 8 & 1 \\ 0&1&0\end{pmatrix}
-\frac{3}{4}E_3\right)\\
&=& \uuuu{e}\begin{pmatrix}5&9&1\\-1&-1&0\\-1&-4&-1\end{pmatrix}.
\end{eqnarray*}
Therefore then 
$G_\Z^M=\{\id,Q\}\times\{\pm (M^{root,4})^l\,|\, l\in\Z\}$.

(c) Observe $r=2x^2$. 
$M$ has the eigenvalues $\lambda_1,\lambda_2,\lambda_3$ with
$$\lambda_{1/2}=(1-x^2)\pm x\sqrt{x^2-2},\quad  \lambda_3=1.$$
Because $M$ is regular, Lemma \ref{t5.1} applies. 
It gives an isomorphism of $\Q$-algebras
\begin{eqnarray*}
\End(H_\Q,M)&&\\
=\{p(M)\,|\, p(t)=p_2t^2+p_1t+p_0\in\Q[t]\}
&\to& \Q[\lambda_1]\times \Q\\
p(M)&\mapsto& (p(\lambda_1),p(1)).
\end{eqnarray*}
The image of $-\id$ is $(-1,-1)$, 
the image of $Q$ is $(-1,1)$, 
the image of $M$ is $(\lambda_1,1)$. 
The image of $G_\Z^M$ is a priori a subgroup of
$\OO_{\Q[\lambda_1]}^*\times\{\pm 1\}$. 
We have to find out which one. By Theorem \ref{t5.11}
$Q\in G_\Z\subset G_\Z^M$. 

Consider the decomposition $H_\Q=H_{\Q,1}\oplus H_{\Q,2}$
as in Definition \ref{t5.9} and the primitive sublattices
$H_{\Z,1}=H_{\Q,1}\cap H_\Z$ and 
$H_{\Z,2}=H_{\Q,2}\cap H_\Z=\Z f_3$ in $H_\Z$.  
The sublattice $H_{\Z,1}$ is the right $L$-orthogonal subspace 
$(\Z f_3)^\perp\subset H_\Z$ of $\Z f_3$, see \eqref{5.8}:
\begin{eqnarray*}
H_{\Z,1}&=&\{\uuuu{e}\cdot \uuuu{z}^t\,|\, \uuuu{z}\in \Z^3,
0=(0,1,-1)\begin{pmatrix}z_1\\z_2\\z_3\end{pmatrix}\}\\
&=& \langle f_1,f_2\rangle_\Z
\quad\textup{with}\quad 
f_1=e_1,\quad f_2=e_2+e_3.
\end{eqnarray*}
Write $\uuuu{f}=(f_1,f_2,f_3)=\uuuu{e}\cdot M(\uuuu{e},\uuuu{f})$, 
\begin{eqnarray*}
M(\uuuu{e})=\uuuu{e}\cdot M^{mat},\quad 
M(\uuuu{f})=\uuuu{f}M^{mat,\uuuu{f}}.
\end{eqnarray*}
Then 
\begin{eqnarray*}
M(\uuuu{e},\uuuu{f})
&=&\begin{pmatrix}1&0&0\\0&1&1\\0&1&-1\end{pmatrix},\quad 
M(\uuuu{e},\uuuu{f})^{-1}
=\frac{1}{2}
\begin{pmatrix}2&0&0\\0&1&1\\0&1&-1\end{pmatrix},\\
M^{mat}&=&\begin{pmatrix}1-2x^2&-x&-x\\x&1&0\\x&0&1\end{pmatrix},\\
M^{mat,\uuuu{f}}
&=&M(\uuuu{e},\uuuu{f})^{-1}M^{mat}M(\uuuu{e},\uuuu{f})
= \begin{pmatrix} 1-2x^2&-2x&0\\x&1&0\\0&0&1\end{pmatrix}.
\end{eqnarray*}

The upper left $2\times 2$-matrix in $M^{mat,\uuuu{f}}$ 
coincides after identification
of $x$ and $y$ with the upper left $2\times 2$-matrix in
$M^{root,2,mat,\uuuu{f}}$ in the proof of part (b).
Therefore we can argue exactly as in the proof of part (b).
Lemma \ref{tc.1} (b) applies in the same way. 

We obtain for $x\geq 3$ 
$G_\Z^M=G_\Z=\{\id,Q\}\times \{\pm M^l\,|\, l\in\Z\}$
and for $x=2$
$G_\Z=\{\id,Q\}\times \{\pm M^l\,|\, l\in\Z\}$.

Consider the case $x=2$. Then $M^{root,2}$ has the first eigenvalue
$$\frac{1}{2}\lambda_1+\frac{1}{2}=-1+\sqrt{2}$$
with $(-1+\sqrt{2})^2=3-2\sqrt{2}=-\lambda_1$
and the third eigenvalue $\frac{1}{2}+\frac{1}{2}=1$. 
Therefore $(M^{root,2})^2=QM$. $M^{root,2}$ is in $G_\Z^M$
because
\begin{eqnarray*}
M^{root,2}(\uuuu{e})&=&\uuuu{e}\left(\frac{1}{2}
\begin{pmatrix}-7&-2&-2\\2&1&0\\2&0&1\end{pmatrix}
+\frac{1}{2}E_3\right)
= \uuuu{e}\begin{pmatrix}-3&-1&-1\\1&1&0\\1&0&1\end{pmatrix}.
\end{eqnarray*}
Therefore $G_\Z^M=\{\id,Q\}\times\{\pm (M^{root,2})^l\,|\, l\in\Z\}$
for $x=2$. 

(d) Here $r=-2$. $M$ has the eigenvalues $\lambda_1,\lambda_2,\lambda_3$ with
$$\lambda_{1/2}=2\pm \sqrt{3},\quad  \lambda_3=1.$$
It is well known and can be seen easily either elementarily or with
Theorem \ref{tc.6} that 
$$\OO_{\Q[\lambda_1]}^*=\{\pm \lambda_1^l\,|\, l\in\Z\}.$$
$Q\in G_\Z\subset G_\Z^M$ by Theorem \ref{t5.11}.
Recall the proof of Lemma \ref{t5.15}.
The restriction of the map in \eqref{5.19} to the map
$$\{\id,Q\}\times \{M^l\,|\, l\in\Z\}\to \OO_{\Q[\lambda_1]}^*$$ 
is an isomorphism. Therefore the map in \eqref{5.19} 
is an isomorphism and 
$$G_\Z^M=G_\Z=\{\id,Q\}\times\{\pm M^l\,|\, l\in\Z\}.$$

(e) Observe in part (c)
\begin{eqnarray*}
Z(\sigma_2)(\uuuu{e})&=&\sigma_2(\uuuu{e})
=(e_1,e_3,e_2),\\
\textup{so }Z(\sigma_2)(f_1,f_2,f_3)&=&(f_1,f_2,-f_3),\\
\textup{ so }Z(\sigma_2)&=&-Q.
\end{eqnarray*}
Now in all cases 
\begin{eqnarray*}
-\id= Z(\delta_1\delta_2\delta_3),\quad M=Z(\sigma^{mon})
\end{eqnarray*}
and 
\begin{eqnarray*}
&&\textup{ in part (a): }\left\{ \begin{array}{lll}
M^{root,3}&=&Z(\delta_3\sigma_2\sigma_1),\\ 
G_\Z&=&\{\pm (M^{root,3})^l\,|\, l\in\Z\}\textup{ for }
x\notin\{4,5\},\end{array}\right. \\
&&\textup{ in part (b): }\left\{ \begin{array}{lll}
M^{root,2}&=&Z(\sigma_2\sigma_1^2),\\ 
G_\Z&=&\{\pm (M^{root,2})^l\,|\, l\in\Z\}\textup{ for }
y\neq 2,\end{array}\right. \\
&&\textup{ in part (c): }\left\{ \begin{array}{lll}
Q&=&Z(\delta_1\delta_2\delta_3\sigma_2),\\ 
G_\Z&=&\{\id,Q\}\times \{\pm M^l\,|\, l\in\Z\}.
\end{array}\right.
\end{eqnarray*}
This shows $G_\Z=G_\Z^{\BB}$ in all but the four cases 
$\uuuu{x}\in \{(4,4,4),(5,5,5),(4,4,8),(3,3,4)\}$.
In these four cases $Q\in G_\Z$. 
It remains to see $Q\notin G_\Z^{\BB}$. We offer two proofs.

First proof: It uses that in these four cases
the stabilizer of $\uuuu{e}$ in 
$\Br_3\ltimes\{\pm 1\}^3$ is $\{\id\}$, which will be proved
as part of Theorem \ref{t7.11}. It also follows from
$\Gamma^{(1)}=G^{free,3}$ in Theorem \ref{t6.18} (g)
or $\Gamma^{(0)}=G^{fCox,3}$ in Theorem \ref{t6.11} (g)
and from Example \ref{t3.4} (respectively Theorem \ref{t3.2}
(a) or (b)).  This implies that here
$Z:(\Br_3\ltimes\{\pm 1\}^3)_S\to G_\Z$ is injective.
Observe furthermore $Q^2=\id$. 
If $Q=Z(\beta)$ for some braid $\beta$, then $\beta^2=\id$
as $Z$ is a group antihomomorphism. But there is no braid
of order two.

Second proof: By formula \eqref{5.5} in Theorem \ref{t5.11} 
in the four cases
\begin{eqnarray*}
Q(\uuuu{e})&=&-\uuuu{e} + 2f_3(1,3,1)
\quad \textup{in the case }\uuuu{x}=(4,4,4),\\
Q(\uuuu{e})&=&-\uuuu{e} + f_3(1,4,1)
\quad \textup{in the case }\uuuu{x}=(5,5,5),\\
Q(\uuuu{e})&=&-\uuuu{e} + f_3(2,7,1)
\quad \textup{in the case }\uuuu{x}=(4,4,8),\\
Q(\uuuu{e})&=&-\uuuu{e} + f_3(4,9,3)
\quad \textup{in the case }\uuuu{x}=(3,3,4).
\end{eqnarray*}
By Theorem \ref{t6.21} (g) the restriction to $\Delta^{(1)}$
of the projection $\pr^{H,(1)}:H_\Z\to \oooo{H_\Z}^{(1)}$
is injective. Therefore in all four cases 
$Q(e_i)\notin\Delta^{(1)}$ for $i\in\{1,2,3\}$.
But any automorphism in $G_\Z^{\BB}$ maps each $e_i$
to an odd vanishing cycle. Thus $Q\notin G_\Z^{\BB}$.
\hfill$\Box$

\section{General rank 3 cases with eigenvalues not all in $S^1$}
\label{s5.7}

Theorem \ref{t5.18} below will show 
$G_\Z=G_\Z^M=\{\pm M^l\,|\, l\in\Z\}$
in all irreducible rank 3 with eigenvalues not all in $S^1$ 
which have not been treated in Theorem \ref{t5.16}.
This result is simple to write down, but the proof is long.
It is a case discussion with many subcases.
It builds on part (b) of the technical Lemma \ref{t5.17}
which gives necessary and sufficient conditions when an
endomorphism in $\End(H_\Q)$ of a certain shape is in
$\End(H_\Z)$. 

Any element of $\{h\in G_\Z^M\,|\, \mu_3=1\}$ can be 
written as $h=q(M)$ with $q(t)=1+q_0(t-1)+q_1(t-1)^2$ with
unique coefficients $q_0,q_1\in\Q$, but not all values 
$q_0,q_1\in\Q$ give such an element. 
Part (b) of the following lemma says which integrality conditions on $q_0$
and $q_1$ are necessary for $q(M)\in \End(H_\Z)$.
Part (a) is good to know in this context.

\begin{lemma}\label{t5.17}
Fix $\uuuu{x}\in \Z^3-\{(0,0,0)\}$
and the associated triple $(H_\Z,L,\uuuu{e})$.
Recall $g=\gcd(x_1,x_2,x_3)$ and $\www{x}_i=g^{-1}x_i$. Define 
\index{$g_1,\ g_2$} 
\begin{eqnarray} 
g_1&:=&\gcd(2x_1-x_2x_3,2x_2-x_1x_3,2x_3-x_1x_2)
\in\N\cup\{0\},\nonumber\\
g_2&:=& \frac{g_1}{g}=\gcd(2\www{x}_1-g\www{x}_2\www{x}_3,
2\www{x}_2-g\www{x}_1\www{x}_3,2\www{x}_3-g\www{x}_1\www{x}_2)\in\N\cup\{0\}.\nonumber \\
&&\label{5.20}
\end{eqnarray}

(a) We separate three cases.\\
(i) Case (three or two of $x_1,x_2,x_3$ are odd): 
Then $g$ and $g_2$ are odd and
\begin{eqnarray*}
\gcd(g_2,\www{x}_i)=\gcd(g_2,g)=1, 
\quad g_2^2\,|\, (r-4).
\end{eqnarray*}
(ii) Case (exactly one of $x_1,x_2,x_3$ is odd): 
Then $g$ is odd, $g_2\equiv 2(4)$ and
\begin{eqnarray*}
\gcd(\frac{g_2}{2},\www{x}_i)=\gcd(\frac{g_2}{2},g)=1, 
\quad (\frac{g_2}{2})^2\,|\, (r-4).
\end{eqnarray*}
(iii) Case (none of $x_1,x_2,x_3$ is odd): Then $g$ and $g_2$ are even.
More precisely, $g_2\equiv 0(4)$ only if $\frac{g}{2}$ and 
$\www{x}_1,\www{x}_2,\www{x}_3$ are odd. Else $g_2\equiv 2(4)$. Always 
\begin{eqnarray*}
\gcd(\frac{g_2}{2},\www{x}_i)=\gcd(\frac{g_2}{2},\frac{g}{2})=1, 
\quad g_2^2\,|\, (r-4).
\end{eqnarray*}

(b) Consider $q_0,q_1\in\Q$, $q(t):=1+q_0(t-1)+q_1(t-1)^2\in\Q[t]$
and $h=q(M)\in \Q[M]$. Define $q_2:=q_0-2q_1\in\Q$. 
\index{$q_0,\ q_1,\ q_2$} 

Then $h\in \End(H_\Z,M)$ if and only if the following 
integrality conditions \eqref{5.21}--\eqref{5.24} are satisfied.
\begin{eqnarray}
q_2\cdot g^2\in\Z,&&\label{5.21}\\
q_1\cdot gg_1\in\Z,\label{5.22}\\
q_0x_i-q_1x_jx_k\in\Z&&\textup{for}\quad \{i,j,k\}=\{1,2,3\},\label{5.23}\\
q_1(x_i^2-x_j^2)\in\Z&&\textup{for}\quad \{i,j,k\}=\{1,2,3\}.\label{5.24}
\end{eqnarray}
If these conditions hold, then also the following holds,
\begin{eqnarray}
q_0\cdot g_1\in\Z.\label{5.25}
\end{eqnarray}

(c) In part (b) the eigenvalue of $q(M)$ on $H_{\C,\lambda_1}$ is 
\begin{eqnarray*}
\mu_1&:=&q(\lambda_1)= (1-rq_1)+(q_0-rq_1)(\lambda_1-1)\\
&=& (1-q_0)+(q_0-rq_1)\lambda_1.
\end{eqnarray*}
\end{lemma}

{\bf Proof:}
(a) $g$ is odd in the cases (i) and (ii) and even in case (iii).
Therefore $g_2$ is odd in case (i) and even in the cases (ii) and (iii),
and furthermore $\frac{g_2}{2}$ is odd in case (ii). 
Also $\frac{g_2}{2}$ odd in case (iii) almost always, namely except when
$\frac{g}{2}$ and $\www{x}_1,\www{x}_2,\www{x}_3$ are odd, as can
be seen from the definition \eqref{5.20} of $g_2$. 
Here observe that at least one of $\www{x}_1,\www{x}_2,\www{x}_3$ is odd
because $\gcd(\www{x}_1,\www{x}_2,\www{x}_3)=1$.

Now we consider first case (iii). A common divisor of $\frac{g_2}{2}$
and $\www{x}_1$ would be odd. Because of the second term 
$2\www{x}_2-g\www{x}_1\www{x}_3$ and the third term
$2\www{x}_3-g\www{x}_1\www{x}_2$ in \eqref{5.20} it would also divide
$\www{x}_2$ and $\www{x}_3$. This is impossible because of
$\gcd(\www{x}_1,\www{x}_2,\www{x}_3)=1$. Therefore 
$\gcd(\frac{g_2}{2},\www{x}_1)=1$. Analogously for $\www{x}_2$ and 
$\www{x}_3$. 

A common divisor of $\frac{g_2}{2}$ and $\frac{g}{2}$ would be odd. 
Because of all three terms in \eqref{5.20} it would divide
$\www{x}_1,\www{x}_2,\www{x}_3$. Therefore
$\gcd(\frac{g_2}{2},\frac{g}{2})=1$. 

For $i,j,k$ with $\{i,j,k\}=\{1,2,3\}$ observe
\begin{eqnarray}
2(2\www{x}_i-g\www{x}_j\www{x}_k) 
+ g\www{x}_k(2\www{x}_j-g\www{x}_i\www{x}_k) 
=\www{x}_i(4-x_k^2),\label{5.26} \\
4(r-4)= g^2(2\www{x}_i-g\www{x}_j\www{x}_k)^2-
(4-x_j^2)(4-x_k^2).\label{5.27}
\end{eqnarray}
\eqref{5.26} and $\gcd(\frac{g_2}{2},\www{x}_i)=1$ imply in case 
(iii) that $\frac{g_2}{2}$ divides $4^{-1}(4-x_k^2)$. 
This and \eqref{5.27} imply that $(\frac{g_2}{2})^2$ divides
$\frac{r-4}{4}$, so $g_2^2$ divides $r-4$. 

The claims for the cases (i) and (ii) follow similarly.

(b) Recall the shape of $M^{mat}\in M_{3\times 3}(\Z)$ with 
$M\uuuu{e}=\uuuu{e}M^{mat}$ from the beginning of Section 
\ref{s5.3}. It gives 
\begin{eqnarray*}
M^{mat}-E_3=
\begin{pmatrix}-x_1^2-x_2^2+x_1x_2x_3&-x_1-x_2x_3+x_1x_3^2
&x_1x_3-x_2\\x_1-x_2x_3&-x_3^2&-x_3\\x_2&x_3&0\end{pmatrix},\\
\begin{pmatrix}1&x_1&x_2\\0&1&x_3\\0&0&1\end{pmatrix}
(M^{mat}-E_3)
=\begin{pmatrix}0&-x_1&-x_2\\x_1&0&-x_3\\x_2&x_3&0\end{pmatrix},\\
(M^{mat}-E_3)
\begin{pmatrix}1&0&0\\0&1&0\\-x_2&-x_3&1\end{pmatrix}
=\begin{pmatrix}-x_1^2&-x_1&x_1x_3-x_2\\x_1&0&-x_3\\x_2&x_3&0\end{pmatrix}.
\end{eqnarray*}
Now 
$$q(M)\in \End(H_\Z)\iff q(M^{mat})-E_3\in M_{3\times 3}(\Z),$$ 
and this is equivalent to the following matrix being in $M_{3\times 3}(\Z)$, 
\begin{eqnarray*}
&&\begin{pmatrix}1&x_1&x_2\\0&1&x_3\\0&0&1\end{pmatrix}
(q(M^{mat})-E_3)
\begin{pmatrix}1&0&0\\0&1&0\\-x_2&-x_3&1\end{pmatrix}\\
&=& 
\begin{pmatrix}0&-x_1&-x_2\\x_1&0&-x_3\\x_2&x_3&0\end{pmatrix}\\
&&\cdot\left[
q_0\begin{pmatrix}1&0&0\\0&1&0\\-x_2&-x_3&1\end{pmatrix}
+q_1\begin{pmatrix}-x_1^2&-x_1&x_1x_3-x_2\\x_1&0&-x_3\\x_2&x_3&0\end{pmatrix}
\right]\\
&=& q_0
\begin{pmatrix} x_2^2 & -x_1+x_2x_3 & -x_2 \\
x_1+x_2x_3 & x_3^2 & -x_3 \\ x_2 & x_3 & 0\end{pmatrix} \\
&+& q_1
\begin{pmatrix} -x_1^2-x_2^2 & -x_2x_3 & x_1x_3 \\
-x_1^3-x_2x_3 & -x_1^2-x_3^2 & x_1^2x_3 -x_1x_2 \\
-x_1^2x_2 +x_1x_3 & -x_1x_2 & x_1x_2x_3-x_2^2-x_3^2\end{pmatrix}.
\end{eqnarray*}
This gives nine scalar conditions, which we denote by their
place $[a,b]$ with $a,b\in\{1,2,3\}$ in the matrix, so for example
$[2,1]$ is the condition $q_0(x_1+x_2x_3)+q_1(-x_1^3-x_2x_3)\in\Z$.
These nine conditions are sufficient and necessary for $q(M)\in \End(H_\Z)$.

The following trick allows an easy derivation of implied conditions.
Recall the cyclic action $\gamma:\Z^3\to\Z^3,\uuuu{x}\mapsto(x_3,x_1,x_2)$.
It lifts to an action of $\Br_3\ltimes\{\pm 1\}^3$ on triangular bases
of $(H_\Z,L)$. Therefore together with 
$M^{mat}=S(\uuuu{x})^{-1}S(\uuuu{x})^t$ also the matrix
$\www{M}^{mat}:= S(\gamma(\uuuu{x}))^{-1}S(\gamma(\uuuu{x}))^t$ 
is a monodromy matrix. 
Integrality of $q(M^{mat})$ is equivalent to integrality of
$q(\www{M}^{mat})$. Therefore if the nine conditions hold, also the
conditions hold which are obtained from the nine conditions by replacing
$(x_1,x_2,x_3)$ by $(x_3,x_1,x_2)$ or by $(x_2,x_3,x_1)$. 
In the following $[a,b]$ denotes all three so obtained conditions,
so for example $[2,1]$ denotes the conditions
\begin{eqnarray*}
q_0(x_i+x_jx_k)+q_1(-x_i^3-x_jx_k)\in\Z\\\
\textup{for } (i,j,k)\in\{(1,2,3),(3,1,2),(2,3,1)\}.
\end{eqnarray*}

We have to show the following equivalence: 
\begin{eqnarray*}
\textup{the conditions }[a,b]\textup{ for }a,b\in\{1,2,3\}
\iff \textup{ the conditions }\eqref{5.21}-\eqref{5.24}.
\end{eqnarray*}
$\Longrightarrow$: $[1,3]$ and $[3,2]$ are equivalent to one another
and to \eqref{5.23}.

$[3,3]$ says $q_1(x_i^2-r)\in\Z$. One derives $q_1(x_i^2-x_j^2)\in\Z$,
which is \eqref{5.24}. 

$[1,1]$ and \eqref{5.24} give $q_2x_i^2\in\Z$, so 
$q_2\gcd(x_1^2,x_2^2,x_3^2)=q_2g^2\in\Z$ which is \eqref{5.21}.

The derivation of $q_1gg_1\in\Z$ is laborious and goes as follows:\\
$[3,1]\& [1,3]$ imply $q_1x_i(2x_k-x_ix_j)\in\Z$.\\
$[3,2]\& [2,3]$ imply $q_1x_i(2x_j-x_ix_k)\in\Z$.\\
$[3,3]\& \eqref{5.24}$ imply $q_1x_i(2x_i-x_jx_k)\in\Z$. \\
One sees $q_1x_ig_1\in\Z$ and then $q_1gg_1\in\Z$.

$\Longleftarrow$: $\eqref{5.23}$ gives $[1,3]$ and $[3,2]$.

\eqref{5.21} and \eqref{5.24} give $[1,1]$ and $[2,2]$.

\eqref{5.23} reduces $[1,2]$, $[2,1]$, $[2,3]$ and $[3,1]$ to
$q_2x_jx_k$, $q_0x_jx_k-q_1x_i^3$, $q_1x_i(x_ix_k-2x_j)$
and $q_1x_i(-x_ix_j+2x_k)$. The first follows from \eqref{5.21},
the third and fourth follow from \eqref{5.22}.
The second reduces with \eqref{5.24} to 
$x_j(q_0x_k-q_1x_ix_j)$, which follows from \eqref{5.23}.

$[3,3]$ reduces with \eqref{5.24} to $q_1x_j(x_ix_k-2x_j)$
which follows from \eqref{5.22}. 

The equivalence of the conditions $[a,b]$ with the conditions
\eqref{5.21}--\eqref{5.24} is shown. 

It remains to show how
\eqref{5.21}-\eqref{5.24} imply \eqref{5.25}. One combines two times
\eqref{5.23}, $2q_0x_i-2q_1x_jx_k\in\Z$, 
with \eqref{5.21}, $(q_0-2q_1)x_jx_k\in\Z$, 
and obtains $q_0(2x_i-x_jx_k)\in\Z$.

(c) Recall 
\begin{eqnarray*}
(t-\lambda_1)(t-\lambda_2)=t^2-(2-r)t+1,\\
\textup{so}\quad \lambda_1+\lambda_2=2-r,\quad \lambda_1\lambda_2=1,\\
\lambda_1^2=(2-r)\lambda_1-1,\quad(\lambda_1-1)^2=(-r)\lambda_1,
\end{eqnarray*}
so
\begin{eqnarray*}
\mu_1&=&1+q_0(\lambda_1-1)+q_1(\lambda_1-1)^2\\
&=&(1-rq_1)+(q_0-rq_1)(\lambda_1-1).\hspace*{2cm}\Box
\end{eqnarray*}

\begin{theorem}\label{t5.18}
Consider a triple $\uuuu{x}\in\Z^3$ with 
$r(\uuuu{x})\in\Z_{<0}\cup\Z_{>4}$ which is neither reducible
nor in the $\Br_3\ltimes\{\pm 1\}^3$ orbit of a triple
in Theorem \ref{t5.16}. More explicitly, the triple $\uuuu{x}$
is any triple in $\Z^3$ which is not in the 
$\Br_3\ltimes\{\pm 1\}^3$ orbits of the triples in the following
set,
\begin{eqnarray*}
\{(x,0,0)\,|\, x\in\Z\}\cup
\{(x,x,x)\,|\, x\in\Z\}\cup
\{(2y,2y,2y^2\,|\, y\in\Z_{\geq 2}\}\\
\cup \{(x,x,0)\,|\, x\in\Z\}\cup
\{(-l,2,-l)\,|\, l\in\Z_{\geq 3}\}\cup
\{(3,3,4)\}.
\end{eqnarray*}
Consider the associated triple $(H_\Z,L,\uuuu{e})$
with $L(\uuuu{e}^t,\uuuu{e})^t=S(\uuuu{x})$. Then  
\begin{eqnarray*}
G_\Z=G_\Z^M=\{\pm M^l\,|\, l\in\Z\}.
\end{eqnarray*}
The map $Z:(\Br_3\ltimes\{\pm 1\}^3)_S\to G_\Z$ is surjective,
so $G_\Z=G_\Z^{\BB}$. 
\end{theorem}

{\bf Proof:}
The surjectivity of $Z$ follows from $G_\Z=\{\pm M^l\,|\, l\in\Z\}$, 
$-\id=Z(\delta_1\delta_2\delta_3)$ and $M=Z(\sigma^{mon})$.
The main point is to prove  
$G_\Z^M=\{\pm M^l\,|\, l\in\Z\}$.

Theorem \ref{t5.11} says for which $\uuuu{x}$ the automorphism
$Q$ of $H_\Q$ in Definition \ref{t5.9} is in $G_\Z^M$.
They are all excluded here. So here $Q\notin G_\Z^M$. 

We use the notation $g=g(\uuuu{x})$ in Lemma \ref{t5.10} 
and the notations from the beginning of section \ref{s5.3}. 
Especially $r:=r(\uuuu{x})\in\Z_{<0}\cup\Z_{>4}$, and
$\lambda_3=1$ and $\lambda_{1/2}=\frac{2-r}{2}\pm\sqrt{r(r-4)}$
are the eigenvalues of the monodromy. 

The proof of Lemma \ref{t5.15} gives a certain control on $G_\Z^M$
and $G_\Z$. Recall the notations there. For
$p(t)=\sum_{i=0}^2p_it^i\in\Q[t]$ write $\mu_j:=p(\lambda_j)$
for the eigenvalues of $p(M)\in \End(H_\Q)$. 
Recall \eqref{5.17}, \eqref{5.18}, the isomorphism of $\Q$-algebras 
\begin{eqnarray*}
\End(H_\Q,M)
=\{p(M)\,|\, p(t)=\sum_{i=0}^2p_it^i\in\Q[t]\}&\to& 
\Q[\lambda_1]\times \Q,\\
p(M)&\mapsto& (p(\lambda_1),p(1)),
\end{eqnarray*}
and its restriction in \eqref{5.19}, the injective group
homomorphism
\begin{eqnarray}\label{5.28}
\{h\in G_\Z^M\,|\, \mu_3=1\}\to \OO_{\Q[\lambda_1]}^*, \quad 
h=p(M)\mapsto \mu_1=p(\lambda_1).
\end{eqnarray}
The image of $Q\in\End(H_\Q,M)$ in $\Q[\lambda_1]\times \Q$
is $(-1,1)$. Because $Q\notin G_\Z^M$, the image in
\eqref{5.28} does not contain $-1$, so it is a cyclic group.
It contains $\lambda_1$ which is the image of $M$.
Therefore the group $\{h\in G_\Z^M\,|\, \mu_3=1\}$
is cyclic. It has two generators which are inverse to one 
another. We denote by $h_{gen}$  the generator
such that a positive power of it is $M$, namely
$(h_{gen})^{l_{gen}}=M$ for a unique number $l_{gen}\in\N$.
We have to prove $l_{gen}=1$. 

We will argue indirectly. We will assume the existence of
a root $h=p(M)\in G_\Z^M$ with $h^l=M$ for some $l\geq 2$,
first eigenvalue $\mu_1=p(\lambda_1)$ and third eigenvalue
$\mu_3=p(1)=1$. Then $\mu_1^l=\lambda_1$
and $\mu_1=(1-q_0)+(q_0-rq_1)\lambda_1$ for certain
$q_0,q_1\in\Q$ which must satisfy the properties in 
\eqref{5.21}--\eqref{5.24}. We will come to a contradiction.

We can restrict to $\uuuu{x}$ in the following set, as
the $\Br_3\ltimes\{\pm 1\}^3$ orbits of the elements in this
set are all $\uuuu{x}$ which we consider in this theorem:

Consider the two sets $Y_I$ and $Y_{II}\subset \Z^3$,
\begin{eqnarray*}
Y_I&:=&\{\uuuu{x}\in \Z^3_{\leq 0}\,|\, x_1\leq x_2\leq x_3\}\\
&-&\bigl[\{(x,0,0)\,|\, x\in\Z_{<0}\}\cup\ \{(x,x,x)\,|\, x\in\Z_{\leq 0}\}\\
&& \cup\ \{(x,x,0)\,|\, x\in\Z_{<0}\}\bigr],\\
Y_{II}&:=&\{\uuuu{x}\in \Z^3_{\geq 3}\,|\, 
x_1\leq x_2\leq x_3, 2x_3\leq x_1x_2\}\\
&-&\bigl[\{(x,x,x)\,|\, x\in\Z_{\geq 3}\}\\
&& \cup\ \{(2y,2y,2y^2)\,|\, y\in\Z_{\geq 2}\}
\cup\{(3,3,4)\}\bigr].
\end{eqnarray*}
All triples $\uuuu{x}$ in this theorem are in the
$\Br_3\ltimes\{\pm 1\}^3$ orbits of the triples in
$Y_I\cup Y_{II}\cup\{\uuuu{x}\,|\, (x_2,x_1,x_3)\in 
Y_I\cup Y_{II}\}$.
We will restrict to $\uuuu{x}\in Y_I\cup Y_{II}$. 
For $\uuuu{x}$ with $(x_2,x_1,x_3)\in Y_I\cup Y_{II}$, one 
can copy the following proof and exchange $x_1$ and $x_2$.

Because the triples $(x,x,x)$ are excluded, 
$x_1<x_3$ and $\www{x}_1<\www{x}_3$. 

We will assume the existence of a unit 
$\mu_1 \in \OO_{\Q[\lambda_1]}^*$ with $\mu_1^l=\lambda_1$
for some $l\geq 2$ and some norm $\NN(\mu_1)\in\{\pm 1\}$.
The integrality conditions in Lemma \ref{t5.17} (b)
for $q_0,q_1,q_2\in\Q$ with $\mu_1=(1-q_0)+(q_0-rq_1)\lambda_1$ 
and $q_2=q_0-2q_1$ will lead to a contradiction.
The proof is a case discussion. The cases split as follows.

Case I: $\uuuu{x}\in Y_I$. 

\hspace*{0.5cm} Subcase I.1: $l\geq 3$ odd.

\hspace*{0.5cm} Subcase I.2: $l=2$. 

Case II: $\uuuu{x}\in Y_{II}$. 

\hspace*{0.5cm} Subcase II.1: $l\geq 3$ odd.

\hspace*{0,5cm} Subcase II.2: $l=2$.

\hspace*{1cm} Subcase II.2.1: $\NN(\mu_1)=1$:

\hspace*{1.5cm} Subcase II.2.1.1: 
$\mu_1=\kappa_a$ for some $a\in\Z_{\geq 3}$.

\hspace*{1.5cm} Subcase II.2.1.2:
$\mu_1=-\kappa_a$ for some $a\in\Z_{\geq 3}$. 

\hspace*{1cm} Subcase II.2.2: $\NN(\mu_1)=-1$.

The treatment of the cases I, II.1 and II.2.2
will be fairly short. The treatment of the cases II.2.1 
will be laborious.

Lemma \ref{td.2} prepares all cases with 
$\NN(\mu_1)=1$. Consider such a case. 
Suppose $\mu_1=\kappa_a$ for some 
$a\in\Z_{\leq -3}\cup\Z_{\geq 3}$.
Compare Lemma \ref{td.2} (c) and Lemma \ref{t5.17} (c):
\begin{eqnarray*}
q_0=q_{0,l}(a),\quad q_1=q_{1,l}(a),\quad q_2=q_{2,l}(a),
\quad r=r_l(a). 
\end{eqnarray*}
The integrality condition \eqref{5.21} $q_2g^2\in\Z$ together
with \eqref{d.6} and \eqref{d.5} tells that 
$r/(2-a)=r_l(a)/(2-a)$ divides $g^2$.

For $l\geq 3$ odd $r/(2-a)$ is itself a square
by \eqref{d.4}, and $g$ can be written as 
$g=\gamma_1\gamma_3$ with $\gamma_1,\gamma_3\in\N$
and $\gamma_1^2=r/(2-a)$. 
For $l=2$ $r/(2-a)=a+2$, and $g$ can be written as
$g=\gamma_1\gamma_2\gamma_3$ with 
$\gamma_1,\gamma_2,\gamma_3\in\N$, $\gamma_2$ squarefree, and 
$a+2=\gamma_1^2\gamma_2$, so $g^2=(a+2)\gamma_2\gamma_3^2$.

On the other hand $g^2$ divides $r=r_l(a)$ by Lemma 1.3.
\eqref{5.1} takes the shape
\begin{eqnarray}\label{5.29}
\www{x}_1^2+\www{x}_2^2+\www{x}_3^2-g\www{x}_1\www{x}_2\www{x}_3
=\frac{r}{g^2} 
=\left\{\begin{array}{ll}
\frac{2-a}{\gamma_3^2}&\textup{ if }l\geq 3\textup{ is odd,}\\
\frac{2-a}{\gamma_2\gamma_3^2}&\textup{ if }l=2.
\end{array}\right.
\end{eqnarray}
This equation will be the key to contradictions in the cases
discussed below. The absolute value of the left hand side will
be large, the absolute value of the right hand side will be small.
Now we start the case discussion.

{\bf Case I.1}, $\uuuu{x}\in Y_I$, $l\geq 3$ odd:
$\NN(\mu_1)^l=\NN(\lambda_1)=1$ and $l$ odd imply $\NN(\mu_1)=1$. 
Here $\lambda_1\in(-1,0)$, so $\mu_1\in(-1,0)$,
so $\mu_1=\kappa_a$ for some $a\in\Z_{\leq -3}$.
By the discussion above  $g=\gamma_1\gamma_3$,  
$\gamma_1=|b_{(l+1)/2}+b_{(l-1)/2}|$, and \eqref{5.29} holds.

{\bf Case I.1.1}, all $x_i<0$: We excluded the triples $(x,x,x)$.
Therefore $\www{x}_1\leq -2$. For $l\geq 3$ 
$\gamma_1=|b_{(l+1)/2}+b_{(l-1)/2}|\geq |b_2+b_1|=|a|-1$ 
by Lemma \ref{td.2} (b). Now 
by \eqref{5.29}
\begin{eqnarray*}
|a|+2\geq \frac{2-a}{\gamma_3^2}
\geq 2g+4+1+1=2\gamma_1\gamma_3+6\geq 2(|a|-1)+6,
\end{eqnarray*}
a contradiction.

{\bf Case I.1.2}, $x_3=0$: 
The integrality conditions \eqref{5.23} and \eqref{5.24}
say here
\begin{eqnarray*}
q_0x_1,q_0x_2,q_1x_1x_2\in\Z,\quad 
q_1x_1^2,q_1x_2^2\in\Z,\quad\textup{so}\quad q_1g^2\in\Z
\end{eqnarray*}
(which is a bit stronger than \eqref{5.22}). 
$q_1g^2\in\Z$ means
\begin{eqnarray*}
\frac{(b_l(a)-b_{l-1}(a)-1)g^2}{rb_l(a)}
=\frac{b_l(a)-b_{l-1}(a)-1}{b_l(a)(2-a)/\gamma_3^2}\in\Z.
\end{eqnarray*}
But $\frac{b_l(a)-b_{l-1}(a)-1}{b_l(a)}\in (1,2)$
because of Lemma \ref{td.2} (b), and 
$\frac{2-a}{\gamma_3^2}\in\N$, a contradiction. 

{\bf Case I.2}, $\uuuu{x}\in Y_{I}$, $l=2$: Here $\lambda_1<0$. 
Therefore $\mu_1^2=\lambda_1$ is impossible.

{\bf Case II.1}, $\uuuu{x}\in Y_{II}$, $l\geq 3$ odd: 
$\NN(\mu_1)^l=\NN(\lambda_1)=1$ and $l$ odd imply $\NN(\mu_1)=1$. 
Here $\lambda_1>1$, so $\mu_1>1$,
so $\mu_1=\kappa_a$ for some $a\in\Z_{\geq 3}$.
By the discussion above  $g=\gamma_1\gamma_3$,  
$\gamma_1=b_{(l+1)/2}+b_{(l-1)/2}$, and \eqref{5.29} holds.
The proof of Lemma \ref{t5.10} (b) gives the first inequality
below,
\begin{eqnarray}\nonumber
\frac{a-2}{\gamma_3^2}&=&\frac{|r|}{g^2}
\geq g\www{x}_1\www{x}_2^2-\www{x}_1^2-2\www{x}_2^2
\geq \www{x}_2^2(g\www{x}_1-3),\\
\textup{so}\quad
a-2&\geq& \gamma_3^2\www{x}_2^2(\gamma_1\gamma_3\www{x}_1-3),
\nonumber\\
\textup{so}\quad
(\gamma_3^2\www{x}_2^2-1)3&\geq& 
(\gamma_3^2\www{x}_2^2-1)\gamma_1\gamma_3\www{x}_1
+\gamma_1\gamma_3\www{x}_1-(a+1).\label{5.30}
\end{eqnarray}
Observe with Lemma \ref{td.2} (b) 
$$\gamma_1=b_{(l+1)/2}+b_{(l-1)/2}\geq b_2+b_1=a+1\geq 4.$$
Therefore the inequality \eqref{5.30} can only hold if
$\gamma_3=\www{x}_2=\www{x}_1=1$ and $\gamma_1=a+1$, so $l=3$. 
Then also $g=\gamma_1\gamma_3=a+1$.
\begin{eqnarray*}
a-2&=&\frac{a-2}{\gamma_3^2}=\frac{|r|}{g^2}= 
g\www{x}_3-\www{x}_3^2-1-1=(a+1-\www{x}_3)\www{x}_3-2,\\
0&=& (\www{x}_3-1)(\www{x}_3-a).
\end{eqnarray*}
We excluded the triples $(x,x,x)$, so $\www{x}_3>1$.
But also $\www{x}_3\leq\frac{g}{2}\www{x}_1\www{x}_2
=\frac{a+1}{2}<a$. A contradiction. 

{\bf Case II.2}, $\uuuu{x}\in Y_{II}$, $l=2$: 
Then $\NN(\mu_1)=\varepsilon_1$ for some
$\varepsilon_1\in\{\pm 1\}$. 
Also $\lambda_1>1$ and $\varepsilon_2\mu_1>1$ for some
$\varepsilon_2\in\{\pm 1\}$. 
Then $\mu_1$ is a zero of a polynomial 
$t^2-\varepsilon_2at+\varepsilon_1$, namely 
\begin{eqnarray*}
\mu_1=\varepsilon_2(\frac{a}{2}
+\frac{1}{2}\sqrt{a^2-4\varepsilon_1})\quad
\textup{with}\left\{\begin{array}{ll}
a\in\Z_{\geq 3}&\textup{ if }\varepsilon_1=1,\\
a\in\N&\textup{ if }\varepsilon_1=-1,\end{array}\right.\\
\mu_1+\mu_1^{conj}=\varepsilon_2a,
\quad \mu_1\mu_1^{conj}=\varepsilon_1,\quad
\mu_1^2=\varepsilon_2 a\mu_1-\varepsilon_1.
\end{eqnarray*} 
Comparison with
\begin{eqnarray*}
\lambda_1&=&\frac{2-r}{2}+\frac{1}{2}\sqrt{r(r-4)}\\
&=& \mu_1^2=\frac{a^2-2\varepsilon_1}{2}+
\frac{a}{2}\sqrt{a^2-4\varepsilon_1}\\\
&=&\varepsilon_2a\mu_1-\varepsilon_1
\end{eqnarray*}
shows
\begin{eqnarray*}
r&=& -a^2+2(\varepsilon_1+1),\\
\mu_1&=& \frac{\varepsilon_1\varepsilon_2}{a}
+\frac{\varepsilon_2}{a}\lambda_1\\
&=& (1-q_0)+(q_0-rq_1)\lambda_1,
\end{eqnarray*}
with
\begin{eqnarray*}
q_0&=&\frac{a-\varepsilon_1\varepsilon_2}{a},\\
q_1&=&\frac{\varepsilon_2(\varepsilon_1+1)-a}
{a(a^2-2(\varepsilon_1+1))}.
\end{eqnarray*}

{\bf Case II.2.1}, $\NN(\mu_1)=1$: Then $\varepsilon_1=1$
and
\begin{eqnarray*}
r&=& -a^2+4=(2-a)(2+a),\\
q_0&=& \frac{a-\varepsilon_2}{a},\\
q_1&=&\frac{2\varepsilon_2-a}{a(a^2-4)}
=\frac{-1}{a(a+2\varepsilon_2)},\\\
q_2&=& q_0-2q_1=\frac{a+\varepsilon_2}{a+2\varepsilon_2}.
\end{eqnarray*}
Write $a+2\varepsilon_2=\gamma_1^2\gamma_2$ with 
$\gamma_1,\gamma_2\in\N$ and $\gamma_2$ squarefree.
The integrality condition \eqref{5.21} $q_2g^2\in\Z$ tells
\begin{eqnarray}\label{5.31}
g=\gamma_1\gamma_2\gamma_3\quad\textup{with}\quad 
a+2\varepsilon_2=\gamma_1^2\gamma_2, \quad
g^2=(a+2\varepsilon_2)\gamma_2\gamma_3^2,
\end{eqnarray} 
for some $\gamma_3\in\N$. The conditions $r<0$ and $g^2|r$ tell
\begin{eqnarray}\label{5.32}
\frac{a-2\varepsilon_2}{\gamma_2\gamma_3^2}
&=& \frac{-r}{g^2}=g\www{x}_1\www{x}_2\www{x}_3
-\www{x}_1^2-\www{x}_2^2-\www{x}_3^2\in\N.
\end{eqnarray}
$\gamma_2$ divides $a+2\varepsilon_2$ and $a-2\varepsilon_2$,
so it divides $4$. As it is squarefree, $\gamma_2\in\{1,2\}$. 
Also $\gcd(g,a)=\gcd(\gamma_1\gamma_2\gamma_3,a)\in\{1,2\}$ as 
$a+2\varepsilon_2=\gamma_1^2\gamma_2$ and 
$\gamma_2\gamma_3^3$ divides $a-2\varepsilon_2$.

The integrality conditions \eqref{5.23},
$q_0x_i-q_1x_jx_k\in\Z$ tell
\begin{eqnarray*}
\frac{a-\varepsilon_2}{a}x_i+\frac{1}{a(a+2\varepsilon_2)}x_jx_k
=x_i+\frac{1}{a}(-g\varepsilon_2\www{x}_i
+\gamma_2\gamma_3^2\www{x}_j\www{x}_k)\in\Z,\nonumber\\
\textup{so after multiplying with }\gamma_1\qquad 
\frac{\gamma_3}{a}(-(a+2\varepsilon_2)\varepsilon_2\www{x}_i
+g\www{x_j}\www{x}_k)\in\Z,\nonumber\\
\textup{so}\quad \frac{\gamma_3}{a}
(-2\www{x}_i+g\www{x}_j\www{x}_k)\in\Z,\nonumber\\
\textup{so}\quad \frac{\gamma_3}{a}g_2\in\Z.\nonumber
\end{eqnarray*}
This is a bit stronger than the integrality condition
\eqref{5.22} $q_1g^2g_2\in\Z$ which says 
$\frac{\gamma_2\gamma_3^2}{a}g_2\in\Z$. 
We can improve it even more, to 
\begin{eqnarray}\label{5.33}
\frac{g_2}{a}\in\Z,\quad 
\frac{1}{a}(-2\www{x}_i+g\www{x}_j\www{x}_k)\in\Z,
\end{eqnarray}
by the following case discussion: 
If $a\equiv 1(2)$, then $\gamma_3\equiv 1(2)$, so 
$\gcd(a,\gamma_3)=1$, so \eqref{5.33} holds.
It $a\equiv 0(4)$, then $a-2\varepsilon_2\equiv 2(4)$,
so $\gamma_3\equiv 1(2)$, so $\gcd(a,\gamma_3)=1$,
so \eqref{5.33} holds. 
If $a\equiv 2(4)$, then a priori 
$\frac{2}{a}g_2\in\Z$. But then $\frac{a}{2}\equiv 1(2)$
and $g\equiv 0(2)$, so $g_2\equiv 0(2)$, so \eqref{5.33} holds.

Finally, the integrality
conditions \eqref{5.24} $q_1(x_i^2-x_j^2)$ say
\begin{eqnarray}
\frac{\gamma_2\gamma_3^2}{a}(\www{x}_i^2-\www{x}_j^2)\in\Z.
\label{5.34}
\end{eqnarray}
The following estimate which arises from \eqref{5.2} 
will also be useful:
\begin{eqnarray}\nonumber
\www{x}_1^2&\leq& \frac{(2+(4-r)^{1/3})^2}{g^2}
= \frac{(2+a^{2/3})^2}{g^2}
= \frac{4+4a^{2/3}+a^{4/3}}{g^2}\\
&\leq& \frac{4+2a^{1/3}+2a+a^{4/3}}{g^2}
=\frac{(a+2)
(2+a^{1/3})}{(a+2\varepsilon_2)
\gamma_2\gamma_3^2}.\label{5.35}
\end{eqnarray}

{\bf Case II.2.1.1}, $\varepsilon_2=1$, $\mu_1=\kappa_a$
for some $a\in\Z_{\geq 3}$: 
The estimate \eqref{5.35} says 
\begin{eqnarray}\label{5.36}
\www{x}_1^2\leq\lfloor \frac{2+a^{1/3}}{\gamma_2\gamma_3^2}
\rfloor
\leq \lfloor 2+a^{1/3}\rfloor \leq a\quad
(\textup{recall }a\in\Z_{\geq 3}).
\end{eqnarray}
Recall that $\uuuu{x}$ is a local minimum, so 
$2\www{x}_3\leq g\www{x}_1\www{x}_2$ and that
$\www{x}_1\leq \www{x}_2\leq\www{x}_3$ and 
$1\leq \www{x}_1<\www{x}_3$ (as $(x,x,x)$ is excluded).

{\bf Case II.2.1.1.1}, $2\www{x}_3<g\www{x}_1\www{x}_2$:
Then \eqref{5.33} gives the existence of $\alpha\in\N$ with
$g\www{x}_1\www{x}_2=\alpha a+2\www{x}_3$. With this we go
into \eqref{5.32},
\begin{eqnarray*}
\frac{a-2}{\gamma_2\gamma_3^2} &=&
\alpha a\www{x}_3-\www{x}_1^2+(\www{x}_3^2-\www{x}_2^2)\\
&\stackrel{\eqref{5.36}}{\geq}& \alpha a\www{x}_3-a+0
\stackrel{\www{x}_3\geq 2}{\geq} a,
\end{eqnarray*}
a contradiction.

{\bf Case II.2.1.1.2}, $2\www{x}_3=g\www{x}_1\www{x}_2$,
$\www{x}_1\geq 2$: \eqref{5.32} takes the shape
\begin{eqnarray*}
\frac{a-2}{\gamma_2\gamma_3^2} &=&
\www{x}_3^2-\www{x}_1^2-\www{x}_2^2
=g^2(\frac{\www{x}_1}{2})^2\www{x}_2^2-\www{x}_1^2-\www{x}_2^2\\
&\geq& (g^2-2)\www{x}_2^2\geq ((a+2)-2)\www{x}_2^2=a\www{x}_2^2
\geq a,
\end{eqnarray*}
a contradiction. 

{\bf Case II.2.1.1.3,} $2\www{x}_3=g\www{x}_1\www{x}_2$,
$\www{x}_1=1$: Write $\gamma_4:=\gamma_2\gamma_3^2$. Then 
\eqref{5.32} takes the shape 
\begin{eqnarray*}
\frac{a-2}{\gamma_4}&=&
\www{x}_3^2-1-\www{x}_2^2=
(\frac{1}{4}(a+2)\gamma_4-1)\www{x}_2^2-1,\\
&\stackrel{\www{x}_2\geq 1}{\geq} & 
\frac{1}{4}(a+2)\gamma_4-2,\\ 
(\gamma_4^2-4)a&\leq& -2(\gamma_4^2-4)+8(\gamma_4-2).
\end{eqnarray*}
If $\gamma_4>2$ then $a\leq -2+\frac{8}{\gamma_4+2}\leq
-2+\frac{8}{5},$ a contradiction.
Therefore $\gamma_4\in\{1,2\}$. If $\gamma_4=1$ then 
\eqref{5.32} becomes 
\begin{eqnarray*}
a-2&=& \frac{a-2}{4}\www{x}_2^2-1,
\quad\textup{so}\quad 4=(a-2)(\www{x}_2^2-4), 
\end{eqnarray*}
which has no solution $(a,\www{x}_2)\in\Z_{\geq 3}\times\N$,
so a contradiction.

If $\gamma_4=2$ then \eqref{5.32} is solved with
$\www{x}_2=1$ and $a\in\Z_{\geq 3}$ arbitrary. Then
$\www{x}_1=1$, $\gamma_2=2$, $\gamma_3=1$, $g=2\gamma_1$,  
$\www{x}_3=\frac{g}{2}=\gamma_1$, 
$\uuuu{x}=(2\gamma_1,2\gamma_1,2\gamma_1^2)$.
These cases are excluded.

{\bf Case II.2.1.2}, $\varepsilon_2=-1$, $\mu_1=-\kappa_a$
for some $a\in\Z_{\geq 3}$: The estimates \eqref{5.35} say here
\begin{eqnarray}\label{5.37}
\www{x}_1^2&\leq& \lfloor \frac{(2+a^{2/3})^2}
{(a-2)\gamma_2\gamma_3^2} \rfloor \leq 
\lfloor \frac{(a+2)(2+a^{1/3})}{(a-2)\gamma_2\gamma_3^2}\rfloor.
\end{eqnarray}
This implies
\begin{eqnarray}\label{5.38}
\www{x}_1^2<a\quad\textup{if}\quad a\geq 8.
\end{eqnarray} 
We treat small $a$ first. Recall $\gamma_2\in\{1,2\}$
and $a=\gamma_1^2\gamma_2+2$. So if $a\leq 9$ then $a\in\{3,4,6\}$.
The following table lists constraints for $a\in\{3,4,6\}$.
For $\www{x}_1$ \eqref{5.37} gives an upper bound and 
$\frac{3}{g}\leq \frac{x_1}{g}=\www{x}_1$ gives a lower bound.
Recall the conditions $a-2=\gamma_1^2\gamma_2$ and 
$\frac{a+2}{\gamma_2\gamma_3^2}\in\N$ .
\begin{eqnarray*}
\begin{array}{c|c|c|c}
a & 3 & 4 & 6 \\
(\gamma_1,\gamma_2,\gamma_3) & (1,1,1) & (1,2,1) & (2,1,1)
\textup{ or }(2,1,2)\\
\gamma_2\gamma_3^2 & 1 & 2 & 1\textup{ or }4\\
g & 1 & 2 & 2 \textup{ or }4\\
\frac{3}{g} & 3 & \frac{3}{2} & \frac{3}{2} \textup{ or }\frac{3}{4} \\
\frac{(2+a^{2/3})^2}{a-2} 
 & 16,64.. & 10,21.. & 7,03.. \\
\www{x}_1 & 3\textup{ or }4 & 2 & 2 \textup{ or }1 
\end{array}
\end{eqnarray*}
This gives five cases $(a,\www{x}_1)\in\{(3,3),(3,4),(4,2),(6,2),(6,1)\}$
with $a\leq 9$. We treat these cases first and then all cases with $a\geq 10$. 
Because of \eqref{5.33} a number $\alpha\in\Z_{\geq 0}$ with
$g\www{x}_1\www{x}_2=\alpha a+2\www{x}_3$ exists. Then \eqref{5.32} becomes
\begin{eqnarray}\label{5.39}
\frac{a+2}{\gamma_2\gamma_3^2}=\alpha a\www{x}_3-\www{x}_1^2 
+(\www{x}_3^2-\www{x}_2^2).
\end{eqnarray}
Also recall \eqref{5.34} 
$\frac{\gamma_2\gamma_3^2}{a}(\www{x}_i^2-\www{x}_j^2)\in\Z$.

{\bf Case II.2.1.2.1}, $(a,\gamma_1,\gamma_2,\gamma_3,g,\www{x}_1)
=(3,1,1,1,1,3)$: \eqref{5.39} says
$$5=3 \alpha\www{x}_3-9+(\www{x}_3^2-\www{x}_2^2).$$ 
\eqref{5.34} says that $3$ divides $\www{x}_3^2-\www{x}_2^2$. 
A contradiction.

{\bf Case II.2.1.2.2}, $(a,\gamma_1,\gamma_2,\gamma_3,g,\www{x}_1)
=(3,1,1,1,1,4)$: \eqref{5.39} says
$$5=3 \alpha\www{x}_3-16+(\www{x}_3^2-\www{x}_2^2).$$ 
$g\www{x}_1\www{x}_2=4\www{x}_2=\alpha a+2\www{x}_3$ implies that
$\alpha$ is even. This and $\www{x}_3>\www{x}_1=4$ and \eqref{5.39}
show $\alpha=0$, so $2\www{x}_2=\www{x}_3$, so 
$5=0-16+3\www{x}_2^2$, a contradiction.

{\bf Case II.2.1.2.3}, $(a,\gamma_1,\gamma_2,\gamma_3,g,\www{x}_1)
=(4,1,2,1,2,2)$: \eqref{5.39} says
$$3=4 \alpha\www{x}_3-4+(\www{x}_3^2-\www{x}_2^2).$$ 
\eqref{5.34} says that $2$ divides $\www{x}_3^2-\www{x}_2^2$.
A contradiction.

{\bf Case II.2.1.2.4}, $(a,\gamma_1,\gamma_2,\gamma_3,g,
\www{x}_1)=(6,2,1,1,2,2)$: \eqref{5.39} says
$$8=6 \alpha\www{x}_3-4+(\www{x}_3^2-\www{x}_2^2).$$ 
$\www{x}_3>\www{x}_1=2$ and \eqref{5.39} imply $\alpha=0$, so
$12=\www{x}_3^2-\www{x}_2^2$. Only $\www{x}_2=2$ and $\www{x}_3=4$
satisfy this. But then $\gcd(\www{x}_1,\www{x}_2,\www{x}_3)=2\neq 1$,
a contradiction.

{\bf Case II.2.1.2.5}, $(a,\gamma_1,\gamma_2,\gamma_3,g,\www{x}_1)
=(6,2,1,2,4,1)$: \eqref{5.39} says
$$2=6 \alpha\www{x}_3-1+(\www{x}_3^2-\www{x}_2^2).$$ 
It implies $\alpha=0$ and $\www{x}_2=1$, $\www{x}_3=2$,
so $\uuuu{x}=(4,4,8)$. This case was excluded in 
Theorem \ref{t5.18}.

{\bf Case II.2.1.2.6}, $a\geq 10$: 

{\bf Case II.2.1.2.6.1}, $\alpha>0$: \eqref{5.38} gives $\www{x}_1^2<a$.
This and \eqref{5.39} and $\www{x}_3>\www{x}_1\geq 1$ show 
$\alpha=1$, $\www{x}_3=2$, $\www{x}_1=1$, $\gamma_2=\gamma_3=1$,
so $a+2=2a-1+(4-\www{x}_2^2)$. A contradiction to $a\geq 10$. 

{\bf Case II.2.1.2.6.2}, $\alpha=0$, $\www{x}_1\geq 2$: 
\eqref{5.39} says
\begin{eqnarray*}
\frac{a+2}{\gamma_2\gamma_3^2}
&=&\www{x}_3^2-\www{x}_1^2-\www{x}_2^2
=((a-2)\gamma_2\gamma_3^2\frac{1}{4}\www{x}_1^2-1)\www{x}_2^2
-\www{x}_1^2,\\
\textup{so}\quad a+2&\geq&
((a-2)-1)\www{x}_2^2-\www{x}_1^2\geq ((a-2)-2)4,\\
\textup{so}\quad 3a&\leq& 18,\quad a\leq 6,
\end{eqnarray*}
a contradiction.

{\bf Case II.2.1.2.6.3}, $\alpha=0$, $\www{x}_1=1$: 
Write $\gamma_4:=\gamma_2\gamma_3^2$. Then
\eqref{5.39} says
\begin{eqnarray}\nonumber
\frac{a+2}{\gamma_4}&=& \www{x}_3^2-1-\www{x}_2^2
=((a-2)\gamma_4\frac{1}{4}-1)\www{x}_2^2-1,\\
\textup{so}\quad 
\www{x}_2^2&=& \frac{4}{\gamma_4^2}\cdot 
\frac{a+2+\gamma_4}{a-2-\frac{4}{\gamma_4}}.\label{5.40}
\end{eqnarray}
The right hand side must be $\geq 1$. This means
\begin{eqnarray*}
\gamma_4^2(a-2-\frac{4}{\gamma_4})&\leq& 4(a+2+\gamma_4),\\
(\gamma_4^2-4)a&\leq& 2(\gamma_4^2-4)+8(\gamma_4+2).
\end{eqnarray*}
If $\gamma_4>2$ then $a\leq 2+\frac{8}{\gamma_4-2}$,
which is in contradiction to $a\geq 10$, as $\gamma_4=3$
would mean $\gamma_2=3$, which is impossible.
Therefore $\gamma_4\in\{1,2\}$. 

If $\gamma_4=1$ then \eqref{5.40} says 
$\www{x}_2^2=4\frac{a+3}{a-6}=4+\frac{36}{a-6}$. 
But the right hand side is not a square for any $a\geq 10$,
a contradiction.

If $\gamma_4=2$ then \eqref{5.40} says
$\www{x}_2^2=\frac{a+4}{a-4}$, which is also not a square
for any $a\geq 10$, a contradiction.

{\bf Case II.2.2}, $\NN(\mu_1)=\varepsilon_1=-1$:
Recall the formulas for $r,q_0$ and $q_1$ at the beginning
of case II.2. Now
\begin{eqnarray*}
r&=& -a^2,\\
q_0&=& \frac{a+\varepsilon_2}{a},\\
q_1&=& \frac{-1}{a^2},\\
q_2&=& \frac{a^2+\varepsilon_2 a+2}{a^2}.
\end{eqnarray*}
The integrality condition \eqref{5.21} $q_2g^2\in\Z$ says
$\frac{a}{2}\,|\, g$ if $a\equiv 2(4)$ and 
$a\,|\, g$ if $a\equiv 1(2)$ or $a\equiv 0(4)$. 
Also $g^2\,|\, r=-a^2$. Therefore $g=a$ or $g=\frac{a}{2}$,
and $g=\frac{a}{2}$ only if $a\equiv 2(4)$. 
The case $g=a$ means $\frac{r}{g^2}=-1$ which implies
by Lemma \ref{t5.11} $Q\in G_\Z$. But all such cases are 
excluded in Theorem \ref{t5.18}. 

Therefore $g=\frac{a}{2}$ and $a\equiv 2(4)$. 
The integrality condition \eqref{5.22} $q_1g^2g_2\in\Z$ says 
$\frac{g_2}{4}\in\Z$. But $g_2\equiv 0(4)$ and $g\equiv 1(2)$
are together impossible in view of the definition
$g_2=\gcd(2\www{x}_1-g\www{x}_2\www{x}_3,
2\www{x}_2-g\www{x}_1\www{x}_3,2\www{x}_3-g\www{x}_1\www{x}_2)$
and $\gcd(\www{x}_1,\www{x}_2,\www{x}_3)=1$. A contradiction.

Therefore in all cases the assumption that a nontrivial
root $\mu_1$ of $\lambda_1$ exists which satisfies the
integrality conditions \eqref{5.21}--\eqref{5.24}
leads to a contradiction. Theorem \ref{t5.18} is proved
\hfill$\Box$

\begin{remark}\label{t5.19}
The results in this chapter 
give complete results on $G_\Z$ and $G_\Z^M\supset G_\Z$
for all unimodular bilinear lattices of rank 3.

The reducible cases: Lemma \ref{t5.4}, Theorem \ref{t5.13}. 

The irreducible cases with $r\in\{0,1,2,4\}$: Theorem \ref{t5.14} 

The irreducible cases with $r\in\Z_{<0}\cup\Z_{>4}$
and $G_\Z^M\supsetneqq \{\pm M^l\,|\,l\in\Z\}$: Theorem \ref{t5.16}.

The irreducible cases with $r\in \Z_{<0}\cup\Z_{>4}$
and $G_\Z^M=\{\pm M^l\,|\,l\in\Z\}$: Theorem \ref{t5.18}.
\end{remark}

\chapter{Monodromy groups and vanishing cycles}\label{s6}
\setcounter{equation}{0}
\setcounter{figure}{0}

This chapter studies the monodromy groups 
$\Gamma^{(0)}$ and $\Gamma^{(1)}$ of the unimodular bilinear
lattices $(H_\Z,L,\uuuu{e})$ with triangular basis $\uuuu{e}$
which have rank 2 or 3. In rank 3 the even as well as the odd
cases split into many different case studies. They make the
chapter long. 

Section \ref{s6.1} considers for $k\in \{0;1\}$ 
the quotient lattice $\oooo{H_\Z}^{(1)}:= H_\Z/\Rad I^{(k)}$
and the induced bilinear form $\oooo{I}^{(k)}$ on it.
Because $\Gamma^{(k)}$ acts trivially on $\Rad I^{(k)}$,
it acts on this quotient lattice and respects $\oooo{I}^{(k)}$. 
The homomorphism $\Gamma^{(k)}\to \Aut(\oooo{H_\Z}^{(1)},
\oooo{I}^{(k)})$ has an image $\Gamma^{(k)}_s$, the 
{\it simple part} of $\Gamma^{(k)}$ and a kernel
$\Gamma^{(k)}_u$, the {\it unipotent part} of $\Gamma^{(k)}$.
There is the exact sequence
$$\{\id\}\to \Gamma^{(k)}_u\to\Gamma^{(k)}\to 
\Gamma^{(k)}_s\to\{\id\}.$$ 
We will study $\Gamma^{(k)}$ together with
$\Gamma^{(k)}_s$ and $\Gamma^{(k)}_u$. 
Also the natural homomorphism $j^{(k)}:H_\Z\to H_\Z^\sharp
:= \Hom_\Z(H_\Z,\Z)$ and in the even case the spinor norm
will be relevant.
Section \ref{s6.1} fixes more or less well known general facts.

Section \ref{s6.2} treats the rank 2 cases. The even cases
$A_1^2$ and $A_2$ are classical and easy. 
In the other even cases $\Gamma^{(0)}\cong G^{fCox,2}$. 
There we can characterize $\Delta^{(0)}$ arithmetically
and geometrically. In the odd case $A_2$ 
$\Gamma^ {(1)}\cong SL_2(\Z)$. In the other irreducible
odd cases $\Gamma^{(1)}\cong G^{free,2}$. The matrix
group $\Gamma^{(1),mat}\subset SL_2(\Z)$ is a Fuchsian
group of the second kind, but has infinite index in 
$SL_2(\Z)$ in most cases. We do not have a characterization
of $\Delta^{(1)}$ which is as nice as in the even cases.

The long Theorem \ref{t6.11} in section \ref{s6.3}
states our results on the even monodromy group 
$\Gamma^{(0)}$ in the rank 3 cases.
The results are detailed except for the local minima
$\uuuu{x}\in\Z^3_{\geq 3}$ with $r(\uuuu{x})\leq 0$
where we only state $\Gamma^{(0)}\cong \Gamma^{(0)}_s\cong G^{fCox,3}$
and $\Gamma^{(0)}_u=\{\id\}$. It is followed by
Theorem \ref{t6.14} which gives the set $\Delta^{(0)}$
of even vanishing cycles in many, but not all cases. 
Especially in the cases of the local minima 
$\uuuu{x}\in\Z^3_{\geq 3}$ with $r(\uuuu{x})\leq 0$ 
we know little and only state
$\Delta^{(0)}=R^{(0)}$ in the case $(3,3,3)$, but
$\Delta^{(0)}\subsetneqq R^{(0)}$ in the four cases
$(3,3,4)$, $(4,4,4)$, $(5,5,5)$ and $(4,4,8)$. 
The result $\Delta^{(0)}=R^{(0)}$ in the case $(3,3,3)$
seems to be new . Its proof is rather laborious. 

Section \ref{s6.4} treats the odd monodromy group $\Gamma^{(1)}$
and the set of odd vanishing cycles $\Delta^{(1)}$
in the rank 3 cases. The long Theorem \ref{t6.18}
fixes the results on $\Gamma^{(1)}$.
The even longer Theorem \ref{t6.21} fixes the results on
$\Delta^{(1)}$. Also their proofs are long.
They are preceded by two technical lemmas, the second one
helps to control $\Gamma^{(1)}_u$. Similar to the even rank 3 cases,
in the case of a local minimum $\uuuu{x}\in\Z^3_{\geq 3}$
with $r(\uuuu{x})\leq 0$ 
$\Gamma^{(1)}\cong \Gamma^{(1)}_s\cong G^{free,3}$
and $\Gamma^{(1)}_u=\{\id\}$. 
In the same case, interestingly, the map $\Delta^{(1)}
\to \oooo{H_\Z}^{(1)}$ is injective. 
This leads in this case to the problem how to recover 
an odd vanishing cycle from its image in $\oooo{H_\Z}^{(1)}$. 
One solution is offered in the most important case
$\uuuu{x}=(3,3,3)$ in Lemma \ref{t6.26}. 
One general application of the Theorems \ref{t6.18} and 
\ref{t6.21} is given in Corollary \ref{t6.23}.
It allows to separate many of the orbits
of the bigger group $(G^{phi}\ltimes \www{G}^{sign})
\rtimes\langle\gamma\rangle$ which acts on
$T^{uni}_3(\Z)$ and $\Z^3$ in Lemma \ref{t4.18}.

\section{Basic observations}\label{s6.1}

Let $(H_\Z,L,\uuuu{e})$ be a unimodular bilinear lattice
of rank $n\in\N$ with a triangular basis $\uuuu{e}$.
Definition \ref{t2.8} gave two monodromy groups
$\Gamma^{(0)}$ and $\Gamma^{(1)}$ and two sets
$\Delta^{(0)}$ and $\Delta^{(1)}$ of vanishing cycles.
Later in this chapter they shall be studied rather 
systematically in essentially all cases with $n=2$
or $n=3$. For that we need some notations and basic facts,
which are collected here. Everything in this section is
well known. Most of it is stated in the even case in
[Eb84] and in the odd case in [Ja83].

\begin{definition}\label{t6.1}
Let $(H_\Z,L)$ be a unimodular bilinear lattice
of rank $n\in\N$.
In the following $k\in\{0;1\}$. Denote
\index{dual lattice}\index{$H_\Z^\sharp=\Hom(H_\Z,\Z)$}
\index{$\oooo{H_\Z}^{(k)},\ \oooo{H_\Z}^{(k),\sharp}$} 
\index{$\OO^{(k),Rad},\ \OO^{(k),Rad}_u$}
\begin{eqnarray*}
O^{(k)}&:=& \Aut(H_\Z,I^{(k)})\quad\textup{the group of 
automorphisms of }H_\Z\\
&& \textup{ which respect }I^{(k)}.\\
H_\Z^\sharp &:=&\Hom(H_\Z,\Z)\quad\textup{the dual lattice}.\\
j^{(k)}&:&H_\Z\to H_\Z^\sharp,\quad 
a\mapsto(b\mapsto I^{(k)}(a,b)).\\
t^{(k)}&:&O^{(k)}\to \Aut(H_\Z^\sharp),\quad
g\mapsto (l\mapsto l\circ g^{-1}).
\end{eqnarray*}
\begin{eqnarray*}
\oooo{H_\Z}^{(k)}&:=& H_\Z/\Rad I^{(k)},\quad
\oooo{H_\R}^{(k)}:=H_\R/\Rad_\R I^{(k)}.\\
\pr^{H,(k)}&=&\oooo{(.)}^{(k)}:H_\Z\to \oooo{H_\Z}^{(k)},\quad
a\mapsto \oooo{a}^{(k)},\quad\textup{the projection}.\\
\oooo{I}^{(k)}&:&\oooo{H_\Z}^{(k)}\times \oooo{H_\Z}^{(k)}\to\Z
\quad\textup{the bilinear form on }\oooo{H_\Z}^{(k)}\\
&& \textup{ which is induced by }I^{(k)},\\
\oooo{H_\Z}^{(k),\sharp}&:=& \Hom(\oooo{H_\Z}^{(k)},\Z)
\quad\textup{the dual lattice}.\\
O^{(k),Rad}&:=&\{g\in O^{(k)}\,|\, g|_{\Rad I^{(k)}}=\id\}.\\
\pr^{A,(k)}&=&\oooo{(.)}:O^{(k),Rad}\to 
\Aut(\oooo{H_\Z}^{(k)},\oooo{I}^{(k)}),\quad g\mapsto\oooo{g},\\ 
&&\textup{the natural map to the set of induced automorphisms}.
\end{eqnarray*}
For any subgroup $G^{(k)}\subset O^{(k),Rad}$ define
\begin{eqnarray*}
G^{(k)}_s&:=& \pr^{A,(k)}(G^{(k)})\subset 
\Aut(\oooo{H_\Z}^{(k)},\oooo{I}^{(k)}),\\
G^{(k)}_u&:=& \ker(\pr^{A,(k)}:G^{(k)}\to 
\Aut(\oooo{H_\Z}^{(k)},\oooo{I}^{(k)})).
\end{eqnarray*}
$G^{(k)}_s$ is called the {\it simple part} of $G^{(k)}$, 
and $G^{(k)}_u$ is called the {\it unipotent part} of $G^{(k)}$. 
\index{simple part}\index{unipotent part}
\end{definition}

\begin{lemma}\label{t6.2}
Let $(H_\Z,L,\uuuu{e})$ be a unimodular bilinear lattice
of rank $n\in\N$ with a triangular basis $\uuuu{e}$.

(a) The map $t^{(k)}:O^{(k)}\to\Aut(H_\Z^\sharp)$ is a group
homomorphism. For $g\in O^{(k)}$, $t^{(k)}(g)$ maps 
$j^{(k)}(H_\Z)\subset H_\Z^\sharp$ to itself, so it induces
an automorphism $\tau^{(k)}(g)\in\Aut(H_\Z^\sharp/j^{(k)}(H_\Z))$.
The map \index{$\tau^{(k)}$} 
\begin{eqnarray*}
\tau^{(k)}:O^{(k)}\to\Aut (H_\Z^\sharp/j^{(k)}(H_\Z))
\end{eqnarray*}
is a group homomorphism.

(b) For $a\in R^{(0)}$ if $k=0$ and for $a\in H_\Z$ if $k=1$
the reflection or transvection $s^{(k)}_a\in O^{(k)}$ 
is in $\ker\tau^{(k)}$.
Therefore $\Gamma^{(k)}\subset \ker \tau^{(k)}$.

(c) $\ker\tau^{(k)}\subset O^{(k),Rad}$.

(d) The horizontal lines of the following diagram are
exact sequences. 
\begin{eqnarray*}
\begin{array}{ccccccccc}
\{\id\} & \to & \Gamma^{(k)}_u & \to & \Gamma^{(k)} 
& \to & \Gamma^{(k)}_s & \to & \{\id\} \\
\| & & \cap & & \cap & & \cap & & \| \\
\{\id\} & \to & (\ker \tau^{(k)})_u & \to & \ker \tau^{(k)} 
& \to & (\ker \tau^{(k)})_s & \to & \{\id\} \\
\| & & \cap & & \cap & & \cap & & \| \\
\{\id\} & \to & O^{(k),Rad}_u & \to & O^{(k),Rad} 
& \to & O^{(k),Rad}_s & \to & \{\id\}
\end{array}
\end{eqnarray*}
The second and third exact sequence split non-canonically. 
\begin{eqnarray*}
O^{(k),Rad}_s=\Aut(\oooo{H_\Z}^{(k)},\oooo{I}^{(k)}).
\end{eqnarray*}

(e) The map 
\index{$T:\oooo{H_\Z}^{(k),\sharp}\otimes \Rad I^{(k)} \to O^{(k),Rad}_u$}
\begin{eqnarray*}
T:\oooo{H_\Z}^{(k),\sharp}\otimes \Rad I^{(k)} &\to& 
O^{(k),Rad}_u,\\
\sum_{i\in I}l_i\otimes r_i&\mapsto& \bigl(a\mapsto 
a+\sum_{i\in I}l_i(\oooo{a}^{(k)})r_i\bigr),\\
\textup{ shorter: }h&\mapsto& 
\bigl(a\mapsto a+h(\oooo{a}^{(k)})\bigr),
\end{eqnarray*}
is an isomorphism between abelian groups with
\begin{eqnarray*}
T(h_1+h_2)=T(h_1)\circ T(h_2), \quad T(h)^{-1}=T(-h).
\end{eqnarray*}
It restricts to an isomorphism
\begin{eqnarray*}
T:\oooo{j}^{(k)}(\oooo{H_\Z}^{(k)})\otimes \Rad I^{(k)}\to
(\ker \tau^{(k)})_u,
\end{eqnarray*}where 
$\oooo{j}^{(k)}:\oooo{H_\Z}^{(k)}\to \oooo{H_\Z}^{(k),\sharp}$ 
is the map
\begin{eqnarray*}
a\mapsto \bigl(b\mapsto \oooo{I}^{(k)}(a,b)\bigr)\quad
\textup{for an arbitrary }b\in \oooo{H_\Z}^{(k)}.
\end{eqnarray*}

(f) For $g\in O^{(k),Rad}$, $a\in\oooo{H_\Z}^{(k)}$, 
$r\in\Rad I^{(k)}$
\begin{eqnarray*}
g\circ T(\oooo{j}^{(k)}(a)\otimes r)\circ g^{-1}
=T(\oooo{j}^{(k)}(\oooo{g}(a))\otimes r).
\end{eqnarray*}

(g) Analogously to $t^{(k)}$ and $\tau^{(k)}$ there are the 
group homomorphisms
\begin{eqnarray*}
\oooo{t}^{(k)}:O^{(k),Rad}_s=\Aut(\oooo{H_\Z}^{(k)},
\oooo{I}^{(k)})&\to& \Aut(\oooo{H_\Z}^{(k),\sharp}),\ 
g\mapsto (l\mapsto l\circ g^{-1}),\\\
\textup{and}\quad 
\oooo{\tau}^{(k)}:O^{(k),Rad}_s&\to& \Aut(\oooo{H_\Z}^{(k),\sharp}
/\oooo{j}^{(k)}(\oooo{H_\Z}^{(k)}).
\end{eqnarray*}
$\tau^{(k)}$ and $\oooo{\tau}^{(k)}$ satisfy
\begin{eqnarray*}
(\ker \tau^{(k)})_s =\ker \oooo{\tau}^{(k)}.
\end{eqnarray*}
\end{lemma}

{\bf Proof:}
(a) The map $t^{(k)}$ is a group homomorphism because
\begin{eqnarray*}
t^{(k)}(g_1g_2)(l)=l\circ(g_1g_2)^{-1}=l\circ g_2^{-1}\circ
g_1^{-1}=t^{(k)}(g_1)t^{(k)}(g_2)(l).
\end{eqnarray*}
$t^{(k)}(g)$ maps $j^{(k)}(H_\Z)$ to itself because
\begin{eqnarray*}
t^{(k)}(g)(j^{(k)}(a))&=&j^{(k)}(a)\circ g^{-1}=
I^{(k)}(a,g^{-1}(.))\\
&=&I^{(k)}(g(a),(.))=j^{(k)}(g(a)).
\end{eqnarray*}
$\tau^{(k)}$ is a group homomorphism because $t^{(k)}$ is one.

(b) Choose $l\in H_\Z^\sharp$ and $b\in H_\Z$. Then
\begin{eqnarray*}
(t^{(k)}(s^{(k)}_a)(l)-l)(b)&=&l\circ (s^{(k)}_a)^{-1}(b)-l(b)\\
&=& l(b-(-1)^k I^{(k)}(a,b)a)-l(b)\\
&=& (-1)^{k+1}I^{(k)}(a,b)l(a)\\
&=&j^{(k)}((-1)^{k+1}l(a)a)(b),\\
\textup{so }t^{(k)}(s^{(k)}_a)(l)-l&=&
j^{(k)}((-1)^{k+1}l(a)a)\in j^{(k)}(H_\Z),\\
\textup{so }\tau^{(k)}(s^{(k)}_a)&=&\id,
\end{eqnarray*}
so $s^{(k)}_a\in\ker \tau^{(k)}$.

(c) Let $g\in \ker\tau^{(k)}$ and let $r\in \Rad I^{(k)}$.
Also $g^{-1}\in\ker \tau^{(k)}$. Choose $l\in H_\Z^\sharp$.
Now $\tau^{(k)}(g^{-1})=\id$ implies
\begin{eqnarray*}
t^{(k)}(g^{-1})(l)-l = j^{(k)}(a)\quad\textup{for some }
a\in H_\Z,\\
0=I^{(k)}(a,r)=j^{(k)}(a)(r)=
\bigl(t^{(k)}(g^{-1})(l)-l)(r)
=l((g-\id)(r)).
\end{eqnarray*}
Because $l$ is arbitrary, $(g-\id)(r)=0$, so 
$g(r)=r$, so $g\in O^{(k),Rad}$. 

(d)  The exact sequences are obvious. Choose an arbitrary splitting
of $H_\Z$ as $\Z$-module into $\Rad I^{(k)}$ and a suitably chosen
$\Z$-module $\www{H_\Z}^{(k)}$,
\begin{eqnarray*}
H_\Z=\Rad I^{(k)}\oplus \www{H_\Z}^{(k)}.
\end{eqnarray*}
The projection $\pr^{H,(k)}:H_\Z\to \oooo{H_\Z}^{(k)}$ restricts
to an isomorphism
\begin{eqnarray*}
\pr^{H,(k)}:(\www{H_\Z}^{(k)},I^{(k)}|_{\www{H_\Z}^{(k)}})
\to (\oooo{H_\Z}^{(k)},\oooo{I}^{(k)}).
\end{eqnarray*}
Via this isomorphism, any element of 
$\Aut(\oooo{H_\Z}^{(k)},\oooo{I}^{(k)})$ lifts to an element
of $O^{(k),Rad}$. 
This shows $O^{(k),Rad}_s=\Aut(\oooo{H_\Z}^{(k)},
\oooo{I}^{(k)})$,
and it gives a non-canonical splitting of the third exact
sequence. The end of the proof of part (g) will show
that this splitting restricts to a non-canonical splitting of the
second exact sequence. 

(e) The fact $\oooo{r}^{(k)}=0$ for $r\in \Rad I^{(k)}$
easily implies that $T$ is a group homomorphism with
$T(h_1+h_2)=T(h_1)T(h_2)$ and with image in $O^{(k),Rad}_u$.

Consider $g\in O^{(k),Rad}_u$. Then $g|_{\Rad I^{(k)}}=\id$
and $(g-\id)(a)\in\Rad I^{(k)}$ for any $a\in H_\Z$.
If $b\in a+\Rad I^{(k)}$ then $(g-\id)(a-b)=0$, so
$(g-\id)(b)=(g-\id)(a)$. Thus there is an element 
$h\in \oooo{H_\Z}^{(k),\sharp}\otimes \Rad I^{(k)}$ with 
$h(\oooo{a}^{(k)})=(g-\id)(a)$ for any $a\in H_\Z$, so 
$T(h)(a)=a+h(\oooo{a}^{(k)})=g(a)$, so $T(h)=g$.

Choose a $\Z$-basis $r_1,...,r_m$ of $\Rad I^{(k)}$ and 
linear forms $l_1,...,l_m\in H_\Z^\sharp$ with
$l_i(r_j)=\delta_{ij}$. Then any 
$r\in \Rad I^{(k)}$ satisfies $r=\sum_{i=1}^ml_i(r)r_i$.

Consider $h\in \oooo{H_\Z}^{(k),\sharp}\otimes \Rad I^{(k)}$
with $T(h)\in (\ker\tau^{(k)})_u$. Then
\begin{eqnarray*}
t^{(k)}(T(h))(l_i)-l_i &=&j^{(k)}(a_i)\quad
\textup{for some }a_i\in H_\Z, \textup{ and also}\\
t^{(k)}(T(h))(l_i)-l_i
&=& l_i\circ T(h)^{-1}-l_i = l_i\circ T(-h)-l_i\\
&=& l_i\circ (\id -h(\oooo{(.)}^{(k)}))-l_i
=-l_i(h(\oooo{(.)}^{(k)})).
\end{eqnarray*}
For $b\in H_\Z$ $h(\oooo{b}^{(k)})\in\Rad I^{(k)}$, so
\begin{eqnarray*}
h(\oooo{b}^{(k)})&=& 
\sum_{i=1}^m l_i(h(\oooo{b}^{(k)}))r_i 
=-\sum_{i=1}^m j^{(k)}(a_i)(b)r_i,\\
\textup{so }h&\in& \oooo{j}^{(k)}(\oooo{H_\Z}^{(k)})
\otimes \Rad I^{(k)}.
\end{eqnarray*}

Going backwards through these arguments, one sees that any
$h\in \oooo{j}^{(k)}(\oooo{H_\Z}^{(k)})\otimes \Rad I^{(k)}$ 
satisfies $T(h)\in (\ker\tau^{(k)})_u$.

(f) For $b\in H_\Z$
\begin{eqnarray*}
(g\circ T(\oooo{j}^{(k)}(a)\otimes r)\circ g^{-1})(b)
&=& g(g^{-1}(b) + \oooo{I}^{(k)}(a,\oooo{g^{-1}(b)}^{(k)})r)\\
&=& b+\oooo{I}^{(k)}(\oooo{g}(a),\oooo{b}^{(k)})r\\
&=& T(\oooo{j}^{(k)}(\oooo{g}(a))\otimes r)(b).
\end{eqnarray*}

(g) The projection $\pr^{H,(k)}=\oooo{()}^{(k)}:H_\Z\to
\oooo{H_\Z}^{(k)}$ induces the embedding
\begin{eqnarray*}
i^{(k)}:\oooo{H_\Z}^{(k),\sharp}&\hookrightarrow& H_\Z^\sharp,\ 
l\mapsto l\circ \pr^{H,(k)},\\
\textup{with}\qquad\Imm(i^{(k)}) 
&=&\{l\in H_\Z^\sharp\,|\, l|_{\Rad I^{(k)}}=0\}\\
\textup{and}\qquad j^{(k)}(H_\Z)&=&i^{(k)}(\oooo{j}^{(k)}(\oooo{H_\Z}^{(k)}))\subset \Imm (i^{(k)}).
\end{eqnarray*}
The three lattices
\begin{eqnarray*}
H_\Z^{\sharp}\supset \Imm(i^{(k)})\supset j^{(k)}(H_\Z)
\end{eqnarray*}
have ranks $n$, $n-\rk \Rad I^{(k)}$, $n-\rk \Rad I^{(k)}$
and are for each $g\in O^{(k),\Rad}$ invariant under the map
$t^{(k)}(g)=(l\mapsto l\circ g^{-1})$.

This map acts trivially on the quotient
$H_\Z^\sharp /\Imm(i^{(k)})$. 
It acts trivially on the quotient $\Imm(i^{(k)})/j^{(k)}(H_\Z)$
if and only if $\oooo{g}^{(k)}\in\ker \oooo{\tau}^{(k)}$.
It acts trivially on the quotient $H_\Z^\sharp/j^{(k)}(H_\Z)$
if and only if $g\in \ker\tau^{(k)}$.
Therefore $(\ker\tau^{(k)})_s\subset\oooo{\tau}^{(k)}$.
It remains to find for each $\www{g}\in \ker\oooo{\tau}^{(k)}$
an element $g\in \ker\tau^{(k)}$ with $\oooo{g}^{(k)}=\www{g}$.

Choose a $\Z$-basis $r_1,...,r_n$ of $H_\Z$ such that
$r_1,...,r_m$ (with $m=\rk\Rad I^{(k)}$) is a $\Z$-basis of
$\Rad I^{(k)}$. Then 
$H_\Z=\Rad I^{(k)}\oplus \www{H_\Z}^{(k)}$ with 
$\www{H_\Z}^{(k)}=\bigoplus_{j=m+1}^n\Z \cdot r_j$ 
is a splitting of 
$H_\Z$ with $\www{H_\Z}^{(k)}\cong\oooo{H_\Z}^{(k)}$. 

Consider the dual $\Z$-basis $l_1,...,l_n$ of $H_\Z^\sharp$ with
$l_i(r_j)=\delta_{ij}$. Then
$\Imm(i^{(k)}) =\bigoplus_{j=m+1}^n \Z\cdot l_j \supset
j^{(k)}(H_\Z)$. 

An element $\www{g}\in O^{(k),Rad}_s$ has a unique lift to an 
element $g\in O^{(k),Rad}$ with 
$g(\www{H_\Z}^{(k)})=\www{H_\Z}^{(k)}$. 
This splitting of the third exact sequence in part (d)
was used already in the proof of part (d).
We claim that $g\in \ker\tau^{(k)}$ if 
$\www{g}\in \ker\oooo{\tau}^{(k)}$. We have 
\begin{eqnarray*}
l_j-l_j\circ g^{-1}&\in& j^{(k)}(H_\Z)\quad\textup{ for }
j\in\{m+1,...,n\}\\
&& \textup{ because of }
\www{g}\in \ker\oooo{\tau}^{(k)},\\
l_i-l_i\circ g^{-1}&=&0\quad\textup{ for }i\in\{1,...,m\}.
\end{eqnarray*}
This shows the claim. Therefore 
$(\ker\tau^{(k)})_s =\ker\oooo{\tau}^{(k)}$.
The claim also shows that the non-canonical splitting of the
third exact sequence in part (d) restricts to a non-canonical
splitting of the second exact sequence in part (d). 
\hfill$\Box$

\begin{remarks}\label{t6.3}
(i) The exact sequence 
$\{\id\}\to \Gamma^{(k)}_u\to\Gamma^{(k)}\to \Gamma^{(k)}_s\to\{\id\}$
splits sometimes, sometimes not. When it splits and when
$\Gamma^{(k)}_u$, $\Gamma^{(k)}_s$ and the splitting are known,
then also $\Gamma^{(k)}$ is known.

(ii) Suppose that one has a presentation of 
$\Gamma^{(k)}_s$, namely an isomorphism
\begin{eqnarray*}
\Gamma^{(k)}_s&\stackrel{\cong}{\longrightarrow}& 
\langle g_1,...,g_n\,|\, w_1(g_1,...,g_m),
...,w_m(g_1,...,g_n)\rangle,\\
\oooo{s^{(k)}_{e_i}}&\mapsto& g_i,
\end{eqnarray*}
where $w_1(g_1....,g_n),...,w_m(g_1,...,g_n)$ are certain words
in $g_1^{\pm 1},...,g_n^{\pm 1}$. Then the group $\Gamma^{(k)}_u$
is the normal subgroup of $\Gamma^{(k)}$ generated by the 
elements $w_1(s^{(k)}_{e_1},...,s^{(k)}_{e_n})$,..., 
$w_m(s^{(k)}_{e_1},...,s^{(k)}_{e_n})$. 
In many of the cases with $n=2$ or $n=3$
we have such a presentation. 

(iii) The symmetric bilinear form $\oooo{I}^{(0)}$ on 
$\oooo{H_\Z}^{(0)}$ is nondegenerate.
It is well known that for any 
$g\in\Aut(\oooo{H_\R}^{(0)},\oooo{I}^{(0)})$ some $m\in\N$
and elements $a_1,...,a_m\in\oooo{H_\R}^{(0)}$ with 
$\oooo{I}^{(0)}(a_i,a_i)\in\R^*$ and 
$g=\oooo{s}^{(0)}_{a_1}...\oooo{s}^{(0)}_{a_m}$ 
exist and that the sign
\begin{eqnarray*}
\oooo{\sigma}(g):=\prod_{i=1}^m\sign(\oooo{I}^{(0)}(a_i,a_i))
\in\{\pm 1\}
\end{eqnarray*}
is independent of $m$ and $a_1,...,a_m$. This sign
$\oooo{\sigma}(g)\in\{\pm 1\}$ 
is the {\it spinor norm} of $g$. The map $\oooo{\sigma}:
\Aut(\oooo{H_\R}^{(0)},\oooo{I}^{(0)})\to\{\pm 1\}$
is obviously a group homomorphism. \index{spinor norm}
\end{remarks}

\begin{definition}\label{t6.4}
Keep the situation of Definition \ref{t6.1}.
Define the {\it spinor norm homomorphism}
\begin{eqnarray*}
\sigma: O^{(0),Rad}\to \{\pm 1\},\quad \sigma(g):=
\oooo{\sigma}(\oooo{g}).
\end{eqnarray*}
Define the subgroup $O^{(k),*}$ of $O^{(k),Rad}$
\begin{eqnarray*}
O^{(k),*}&:=& \left\{\begin{array}{ll}
\ker\tau^{(1)}& \textup{if }k=1,\\
\ker\tau^{(0)}\cap\ker\sigma & \textup{if }k=0.
\end{array}\right. 
\end{eqnarray*}
\end{definition}

\begin{remarks}\label{t6.5}
(i) For $a\in R^{(0)}$ of course $\sigma(s^{(0)}_a)=1$.
Therefore $\Gamma^{(0)}\subset\ker\sigma$. Thus
\begin{eqnarray*}
\Gamma^{(k)}\subset O^{(k),*}\quad\textup{for }k\in\{0;1\}.
\end{eqnarray*}
(ii) For $g\in (\ker\tau^{(0)})_u$ $\sigma(g)=1$ because
$\oooo{g}=\id$. Therefore
\begin{eqnarray*}
O^{(k),*}_u=(\ker\tau^{(k)})_u\quad\textup{for }k\in\{0;1\}.
\end{eqnarray*}
\end{remarks}

\begin{remarks}\label{t6.6}
Finally we make some comments on the sets of vanishing cycles
in a unimodular bilinear lattice $(H_\Z,L,\uuuu{e})$ with
a triangular basis. 

(i) Let $H_\Z^{prim}$ denote the set of primitive vectors in
$H_\Z$, i.e. vectors $a\in H_\Z-\{0\}$ with $\Z a=\Q a\cap H_\Z$,
and analogously $\oooo{H_\Z}^{(k),prim}$.
Then
\begin{eqnarray*}
\Delta^{(0)}\subset R^{(0)}\subset H_\Z^{prim},\quad
\Delta^{(1)}\subset H_\Z^{prim},\\
\oooo{\Delta}^{(0)}\subset \oooo{R}^{(0)}\subset 
\oooo{H_\Z}^{(0),prim},\quad\textup{where }
\oooo{\Delta}^{(k)}:=\pr^{H,(k)}(\Delta^{(k)}).
\end{eqnarray*}
Here $R^{(0)}\subset H_\Z^{prim}$ and 
$\oooo{R}^{(0)}\subset\oooo{H_\Z}^{(0),prim}$ because of
$2=I^{(0)}(a,a)=\oooo{I}^{(0)}(\oooo{a}^{(0)},\oooo{a}^{(0)})$
for $a\in R^{(0)}$. Furthermore 
\begin{eqnarray*}
\oooo{e_i}^{(1)}\in \oooo{H_\Z}^{(1),prim}&\iff&
\Gamma^{(1)}_s\{\oooo{e_i}^{(1)}\}\subset \oooo{H_\Z}^{(1),prim}.
\end{eqnarray*}
Whether or not $\oooo{e_i}^{(1)}\in \oooo{H_\Z}^{(1),prim}$,
that depends on the situation.
$\oooo{\Delta}^{(1)}\subset\oooo{H_\Z}^{(1),prim}$ may hold
or not. 

(ii) In general, an element $a\in H_\Z$ satisfies
\begin{eqnarray*}
\oooo{a}^{(k)}\in \oooo{H_\Z}^{(k),prim}&\iff&
a+\Rad I^{(k)}\subset H_\Z^{prim}.
\end{eqnarray*}

(iii) The set $\oooo{\Delta}^{(k)}\subset \oooo{H_\Z}^{(k)}$
is often simpler to describe than the set $\Delta^{(k)}$.
The control of $\oooo{\Delta}^{(k)}$ is a step towards the 
control of $\Delta^{(k)}$. 

(iv) Given $a\in\Delta^{(k)}$, it is interesting to understand the
three sets
\begin{eqnarray*}
(a+\Rad I^{(k)})\cap\Delta^{(k)}
\stackrel{(1)}{\supset} (a+\Rad I^{(k)})\cap \Gamma^{(k)}\{a\}
\stackrel{(2)}{\supset} \Gamma^{(k)}_u\{a\}.
\end{eqnarray*}
The next lemma makes comments on the inclusions
$\stackrel{(1)}{\supset}$ and $\stackrel{(2)}{\supset}$. 
\end{remarks}

\begin{lemma}\label{t6.7}
Keep the situation in Remark \ref{t6.6} (iv).

(a) In $\stackrel{(1)}{\supset}$ equality holds if 
and only if the image in $\oooo{H_\Z}^ {(k)}$ of any
$\Gamma^{(k)}$ orbit different from $\Gamma^{(k)}\{a\}$
is different from $\Gamma^{(k)}_s\{\oooo{a}^{(k)}\}$. 

(b) The following is an inclusion of groups,
\begin{eqnarray*}
\Stab_{\Gamma^{(k)}}(\oooo{a}^{(k)})\stackrel{(3)}{\supset} 
\Gamma^{(k)}_u\cdot \Stab_{\Gamma^{(k)}}(a).
\end{eqnarray*}
In $\stackrel{(2)}{\supset}$ equality holds 
if and only if in $\stackrel{(3)}{\supset}$ equality holds.
\end{lemma}

{\bf Proof:} Trivial. \hfill$\Box$

\section{The rank 2 cases}\label{s6.2}

For $x\in\Z-\{0\}$ consider the matrix
$S=S(x)=\begin{pmatrix}1&x\\0&1\end{pmatrix}\in T^{uni}_2(\Z)$,
and consider a unimodular bilinear lattice $(H_\Z,L)$ with
a triangular basis $\uuuu{e}=(e_1,e_2)$ with
$L(\uuuu{e}^t,\uuuu{e})^t=S$.
Recall the formulas and the results in section \ref{s5.2},
especially $M^{root}:H_\Z\to H_\Z$ and its eigenvalues 
$\kappa_{1/2}=\frac{-x}{2}\pm\frac{1}{2}\sqrt{x^2-4}$.

We can restrict to $x<0$ because of
$L((e_1,-e_2)^t,(e_1,-e_2))^t=\begin{pmatrix}1&-x\\0&1\end{pmatrix}$.
So suppose $x<0$. 

First we consider the even cases. Then $\Gamma^{(0)}$ is a
Coxeter group, and $\Gamma^{(0)}$ and $\Delta^{(0)}$ are
well known. Still we want to document and derive the facts
in our way.

\begin{theorem}\label{t6.8}
(a) We have
\begin{eqnarray*}
s^{(0)}_{e_i}\uuuu{e}&=&\uuuu{e}\cdot s_{e_i}^{(0),mat}
\quad\textup{ with }\\
s^{(0),mat}_{e_1}&=&\begin{pmatrix}-1&-x\\0&1\end{pmatrix},\quad
s^{(0),mat}_{e_2}=\begin{pmatrix}1&0\\-x&-1\end{pmatrix},\\
\Gamma^{(0)}&\cong& \Gamma^{(0),mat}:=
\langle s^{(0),mat}_{e_1},s^{(0),mat}_{e_2}\rangle 
\subset GL_2(\Z),\\
R^{(0)}&=&\{y_1e_1+y_2e_2\in H_\Z\,|\, 1=y_1^2+xy_1y_2+y_2^2\}.
\end{eqnarray*}

(b) The case $x=-1$: $(H_\Z,I^{(0)})$ is the $A_2$ root lattice.
$\Gamma^{(0)}\cong D_6$ is a dihedral group with six elements,
the identity, three reflections and two rotations,
\index{dihedral group}
\begin{eqnarray*}
\Gamma^{(0)}&=& 
\langle\id,\ s^{(0)}_{e_1},\ s^{(0)}_{e_2},
\ s^{(0)}_{e_1}s^{(0)}_{e_2}s^{(0)}_{e_1},
\ s^{(0)}_{e_1}s^{(0)}_{e_2},\ s^{(0)}_{e_2}s^{(0)}_{e_1}\rangle
\cong\Gamma^{(0),mat}\\  &=&\langle
E_2,\begin{pmatrix}-1&1\\0&1\end{pmatrix},
\begin{pmatrix}1&0\\1&-1\end{pmatrix},
\begin{pmatrix}0&-1\\-1&0\end{pmatrix},
\begin{pmatrix}0&-1\\1&-1\end{pmatrix},
\begin{pmatrix}-1&1\\-1&0\end{pmatrix}\rangle.
\end{eqnarray*}

The set $\Delta^{(0)}$ of vanishing cycles coincides with
the set $R^{(0)}$ of roots and is 
$$\Delta^{(0)}=R^{(0)}=\{\pm e_1,\pm e_2,\pm(e_1+e_2)\}.$$
The following picture shows the action of $\Gamma^{(0)}$
on $\Delta^{(0)}$. One sees the action of $D_6$ on the vertices
of a regular 6-gon.
\begin{figure}[H]
\includegraphics[width=1.0\textwidth]{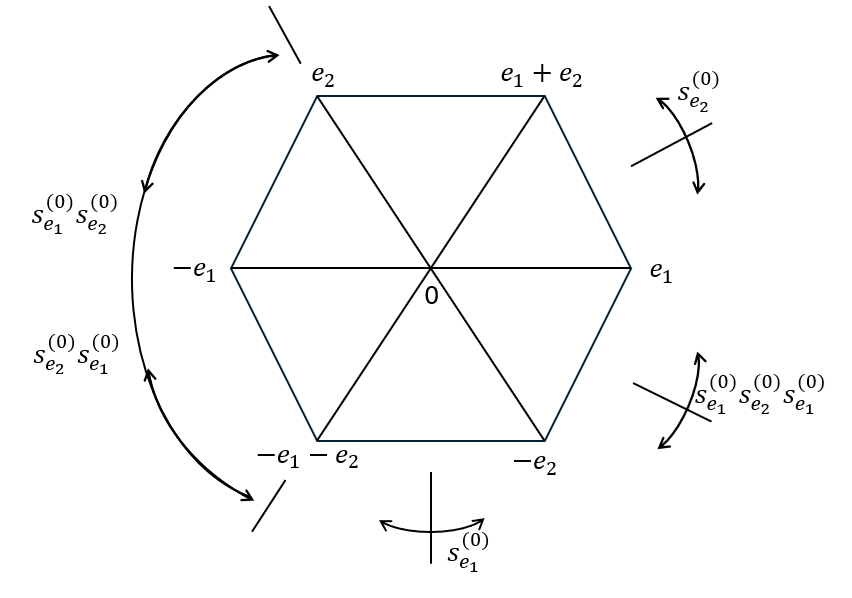}
\caption[Figure 6.1]{A regular 6-gon, actions of 
$D_6$ and $D_{12}$}
\label{Fig:6.1}
\end{figure}

One sees also the action of $D_{12}$ on the regular 6-gon.
This fits to the following.
\begin{eqnarray*}
D_6\cong \Gamma^{(0)}=\ker\tau^{(0)}=O^{(0),*}
\stackrel{1:2}{\subset}O^{(0)}\cong D_{12}.
\end{eqnarray*}

(c) The case $x=-2$: 
$\Gamma^{(0)}\cong G^{fCox,2}$ is a free Coxeter group with 
the two generators $s^{(0)}_{e_1}$ and $s^{(0)}_{e_2}$. Here
\begin{eqnarray*}
\Rad I^{(0)}=\Z f_1\quad\textup{with}\quad f_1=e_1+e_2.
\end{eqnarray*}
The set $\Delta^{(0)}$ of vanishing cycles coincides with the set
$R^{(0)}$ of roots and is

\begin{eqnarray*}
\Delta^{(0)}=R^{(0)}&=&\{y_1e_1+y_2e_2\in H_\Z\,|\,
1=(y_1-y_2)^2\}\\
&=& (e_1+\Z f_1)\ \dot\cup\ (e_2+\Z f_1).
\end{eqnarray*}
It splits into the two disjoint orbits
\begin{eqnarray*}
\Gamma^{(0)}\{e_1\}=\Gamma^{(0)}\{-e_1\}
&=& (e_1+\Z 2f_1)\ \dot\cup\ (-e_1+\Z 2f_1),\\
\Gamma^{(0)}\{e_2\}=\Gamma^{(0)}\{-e_2\}
&=& (e_2+\Z 2f_1)\ \dot\cup\ (-e_2+\Z 2f_1).
\end{eqnarray*}
$s^{(0)}_{e_1}$ acts on $\Delta^{(0)}$ by permuting 
vanishing cycles horizontally, so by adding $\pm 2e_1$.
$s^{(0)}_{e_2}$ acts on $\Delta^{(0)}$ by permuting 
vanishing cycles vertically, so by adding $\pm 2e_2$,
see the following formulas and Figure \ref{Fig:6.2}.
For $\varepsilon\in\{\pm 1\}$ and $m\in\Z$
\begin{eqnarray*}
s_{e_1}^{(0)}(\varepsilon e_1+mf_1)=-\varepsilon e_1+mf_1,\quad
s_{e_1}^{(0)}(\varepsilon e_2+mf_1)=-\varepsilon e_2+(m+2\varepsilon)f_1,\\
s_{e_2}^{(0)}(\varepsilon e_1+mf_1)=-\varepsilon e_1+(m+2\varepsilon)f_1,\quad
s_{e_2}^{(0)}(\varepsilon e_2+mf_1)=-\varepsilon e_2+mf_1.
\end{eqnarray*}

\begin{figure}[H]
\includegraphics[width=0.7\textwidth]{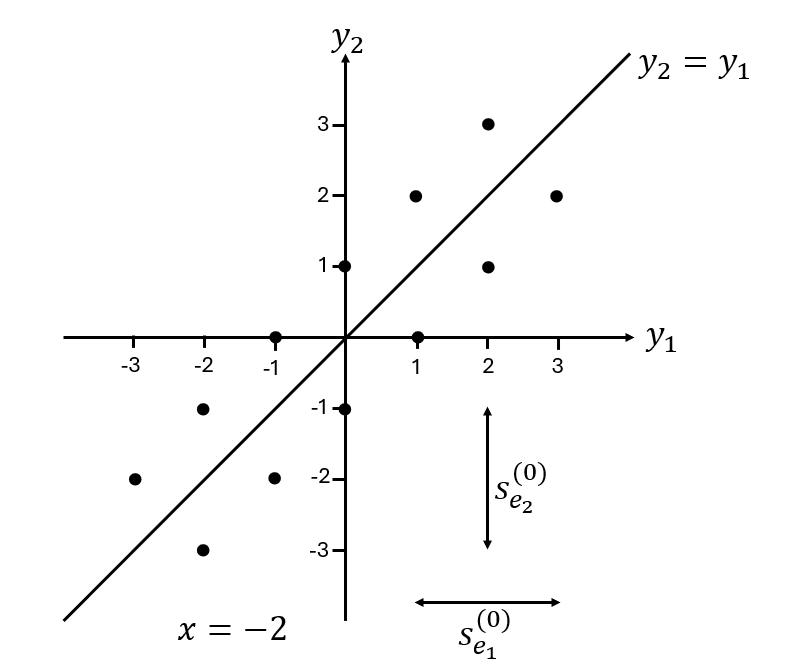}
\caption[Figure 6.2]{Even vanishing cycles in the case
$S=\begin{pmatrix}1&-2\\0&1\end{pmatrix}$}
\label{Fig:6.2}
\end{figure}

\noindent
The matrix group $\Gamma^{(0),mat}$ is given by the following formulas
for $m\in\Z$,
\begin{eqnarray*}
(s_{e_1}^{(0)}s_{e_2}^{(0)})^m(\uuuu{e}) &=& 
\uuuu{e}\begin{pmatrix}2m+1&-2m\\2m&-2m+1\end{pmatrix},\\
(s_{e_1}^{(0)}s_{e_2}^{(0)})^ms_{e_1}^{(0)}(\uuuu{e}) &=& 
\uuuu{e}\begin{pmatrix}-2m-1&2m+2\\-2m&2m+1\end{pmatrix}.
\end{eqnarray*}

(d) The cases $x\leq -3$: 
$\Gamma^{(0)}\cong G^{fCox,2}$ is a free Coxeter group with 
the two generators $s^{(0)}_{e_1}$ and $s^{(0)}_{e_2}$.
The set $\Delta^{(0)}$ of vanishing cycles coincides with the set
$R^{(0)}$ of roots. More information on $\Delta^{(0)}$:

(i) Recall Lemma \ref{tc.1} (a). The map
\begin{eqnarray*}
u:\Delta^{(0)}&\to& \{\textup{units in }\Z[\kappa_1]\textup{ with 
norm }1\} = \{\pm \kappa_1^l\,|\, l\in\Z\}\\
y_1e_1+y_2e_2&\mapsto & y_1-\kappa_1 y_2,
\end{eqnarray*}
is well defined and a bijection with
\begin{eqnarray*}
s^{(0)}_{e_1}(u^{-1}(\varepsilon\kappa_1^l))
=u^{-1}(-\varepsilon\kappa_1^{-l}),\quad 
s^{(0)}_{e_2}(u^{-1}(\varepsilon\kappa_1^l))
=u^{-1}(-\varepsilon\kappa_1^{2-l}),
\end{eqnarray*}
for $\varepsilon\in\{\pm 1\}$, $l\in\Z$. 
Especially, $\Delta^{(0)}$ splits into the two disjoint orbits
$\Gamma^{(0)}\{e_1\}$ and $\Gamma^{(0)}\{e_2\}$. 

(ii) The matrix group $\Gamma^{(0),mat}$ 
is given by the following formulas for $m\in\Z$,
\begin{eqnarray*}
(s_{e_1}^{(0)}s_{e_2}^{(0)})^m(\uuuu{e}) &=& 
\uuuu{e}\begin{pmatrix}y_1&-y_2\\y_2&y_1+xy_2\end{pmatrix}
=(u^{-1}(\kappa_1^{-2m}),u^{-1}(-\kappa_1^{-2m+1})),\\
&&\textup{where }u^{-1}(\kappa_1^{-2m})=y_1-\kappa_1y_2,\\
(s_{e_1}^{(0)}s_{e_2}^{(0)})^ms_{e_1}^{(0)}(\uuuu{e}) &=& 
\uuuu{e}\begin{pmatrix}y_1&xy_1+y_2\\y_2&-y_1\end{pmatrix}
=(u^{-1}(-\kappa_1^{-2m}),u^{-1}(\kappa_1^{-2m-1}))\\
&&\textup{where }u^{-1}(-\kappa_1^{-2m})=y_1-\kappa_1y_2.
\end{eqnarray*}

(iii) $\Delta^{(0)}\subset H_\R\cong\R^2$ is part of the
hyperbola
$\{y_1e_1+y_2e_2\in H_\R\,|\, (y_1-\kappa_1y_2)(y_1-\kappa_2y_2)
=1\}$ with asymptotic lines
$y_2=\kappa_2y_1$ and $y_2=\kappa_1y_1$. Both branches of this
hyperbola are strictly monotonously increasing. The lower right
branch is concave and contains the points
$u^{-1}(\kappa_1^l)$ for $l\in\Z$, the upper left branch is
convex and contains the points $u^{-1}(-\kappa^l)$ for $l\in\Z$.
The horizontal respectively vertical line through a 
vanishing cycle $a\in \Delta^{(0)}$ intersects the other branch
of the hyperbola in $s^{(0)}_{e_1}(a)$ respectively
$s^{(0)}_{e_2}(a)$. See Figure \ref{Fig:6.3}.
\begin{figure}
\includegraphics[width=1.0\textwidth]{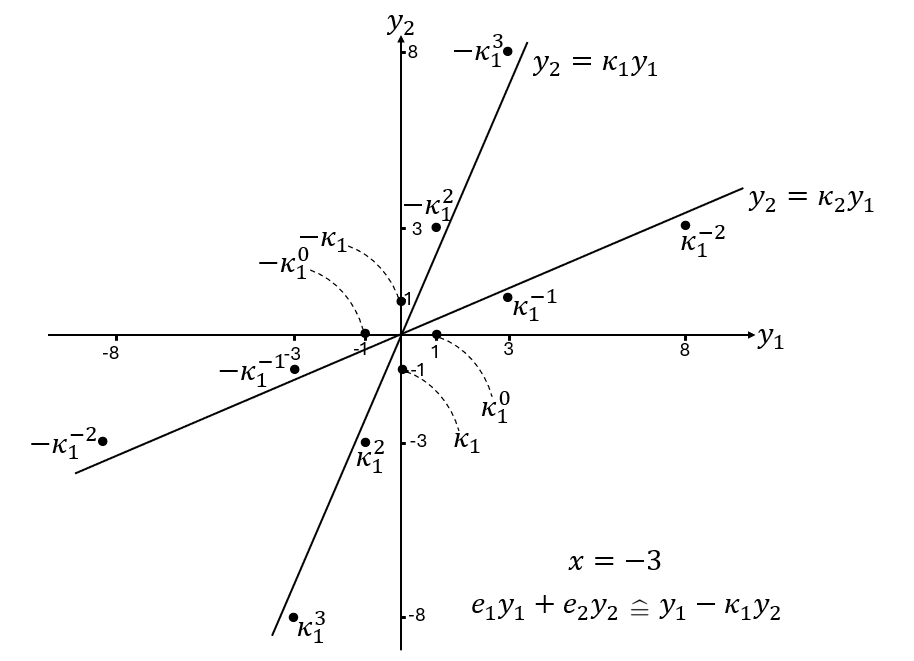}
\caption[Figure 6.3]{Even vanishing cycles in the case
$S=\begin{pmatrix}1&-3\\0&1\end{pmatrix}$}
\label{Fig:6.3}
\end{figure}

(iv) Denote by $\www{s}^{(0)}_{e_i}$, $\www{\Gamma}^{(0)}$
and $\www{\Delta}^{(0)}$ the objects for $x=-2$ and 
as usual by $s^{(0)}_{e_i}$, $\Gamma^{(0)}$ and $\Delta^{(0)}$
the objects for an $x\leq -3$. 

The map $\www{s}^{(0)}_{e_i}\mapsto s^{(0)}_{e_i}$ extends to 
a group isomorphism $\www{\Gamma}^{(0)}\to\Gamma^{(0)}$.
The map 
\begin{eqnarray*}
&&\www{\Delta}^{(0)}\to\Delta^{(0)},\\
&&\uuuu{e}\begin{pmatrix}1-l\\-l\end{pmatrix}\mapsto
u^{-1}(\kappa_1^l),\quad
\uuuu{e}\begin{pmatrix}l-1\\l\end{pmatrix}\mapsto
u^{-1}(-\kappa_1^l)\quad \textup{ for }l\in\Z,
\end{eqnarray*}
is a bijection. The bijections $\www{\Gamma}^{(0)}\to\Gamma^{(0)}$
and $\www{\Delta}^{(0)}\to\Delta^{(0)}$ are compatible
with the action of $\www{\Gamma}^{(0)}$ on $\www{\Delta}^{(0)}$
and of $\Gamma^{(0)}$ on $\Delta^{(0)}$.

(e) More on the cases $x\leq -2$: The automorphism
$g_{1,2}:H_\Z\to H_\Z$ with $g_{1,2}:e_1\leftrightarrow e_2$
is in $O^{(0)}$. The set $\{\pm \id,\pm g_{1,2}\}$ is a subgroup
of $O^{(0)}$ with
\begin{eqnarray*}
\Gamma^{(0)}=\ker\tau^{(0)}=O^{(0),*}
\stackrel{1:4}{\subset} O^{(0)}=\Gamma^{(0)}\rtimes
\{\pm \id,\pm g_{1,2}\}.
\end{eqnarray*}
\end{theorem}

{\bf Proof:}
(a) Everything except possibly the shape of $R^{(0)}$ is obvious.
\begin{eqnarray*}
R^{(0)}&=&\{y_1e_1+y_2e_2\in H_\Z\,|\,
2=I^{(0)}(y_1e_1+y_2e_2,y_1e_1+y_2e_2)\}\\
&=&\{y_1e_1+y_2e_2\in H_\Z\,|\, 
2=\begin{pmatrix}y_1&y_2\end{pmatrix}
\begin{pmatrix}2&x\\x&2\end{pmatrix}
\begin{pmatrix}y_1\\y_2\end{pmatrix}\}\\
&=&\{y_1e_1+y_2e_2\in H_\Z\,|\, 1=y_1^2+xy_1y_2+y_2^2\}.
\end{eqnarray*}

(b) This is classical and elementary. 
$R^{(0)}=\{\pm e_1,\pm e_2,\pm (e_1+e_2)\}$.
The actions of $s^{(0)}_{e_1}$ and $s^{(0)}_{e_2}$
on this set extend to the action of the dihedral group $D_6$ 
on the vertices of a regular 6-gon. Therefore
$\Delta^{(0)}=R^{(0)}$ and $\Gamma^{(0)}\cong D_6$.

$O^{(0)}\cong D_{12}$ is obvious as the vanishing cycles 
form the vertices of a regular 6-gon in $(H_\Z,I^{(0)})$.
It remains to show for some element $g\in O^{(0)}-\Gamma^{(0)}$
$g\notin\ker\tau^{(0)}$. 

Consider the reflection $g\in O^{(0)}$ with 
$g(\uuuu{e})=(e_1,-e_1-e_2)$ and the linear form 
$l:H_\Z\to\Z$ with $l(\uuuu{e})=(1,0)$. Then
\begin{eqnarray*}
t^{(0)}(g)(l)=l\circ g^{-1},\quad 
(l\circ g^{-1})(\uuuu{e})=(l(e_1),l(-e_1-e_2))=(1,-1),\\
(l-l\circ g^{-1})(\uuuu{e})=(0,1),\\
l-l\circ g^{-1}\notin j^{(0)}(H_\Z)=
\langle (\uuuu{e}\mapsto (2,-1)),(\uuuu{e}\mapsto (-1,2))\rangle,
\end{eqnarray*}
so $g\notin\ker\tau^{(0)}$. 

(c) and (d) The group $\Gamma^{(0)}$ for $x\leq -2$:
Recall the Remarks and Notations \ref{ta.1}.
The matrices
$s_{e_1}^{(0),mat}=\begin{pmatrix}-1&-x\\0&1\end{pmatrix}$ and
$s_{e_2}^{(0),mat}=\begin{pmatrix}1&0\\-x&-1\end{pmatrix}$
have the eigenvectors $\begin{pmatrix}1\\0\end{pmatrix}$ 
respectively $\begin{pmatrix}0\\1\end{pmatrix}$ with eigenvalue 
$-1$ and the eigenvectors $\begin{pmatrix}-x/2\\1\end{pmatrix}$
respectively $\begin{pmatrix}-2/x\\1\end{pmatrix}$ with 
eigenvalue $1$

\begin{figure}[H]
\includegraphics[width=0.9\textwidth]{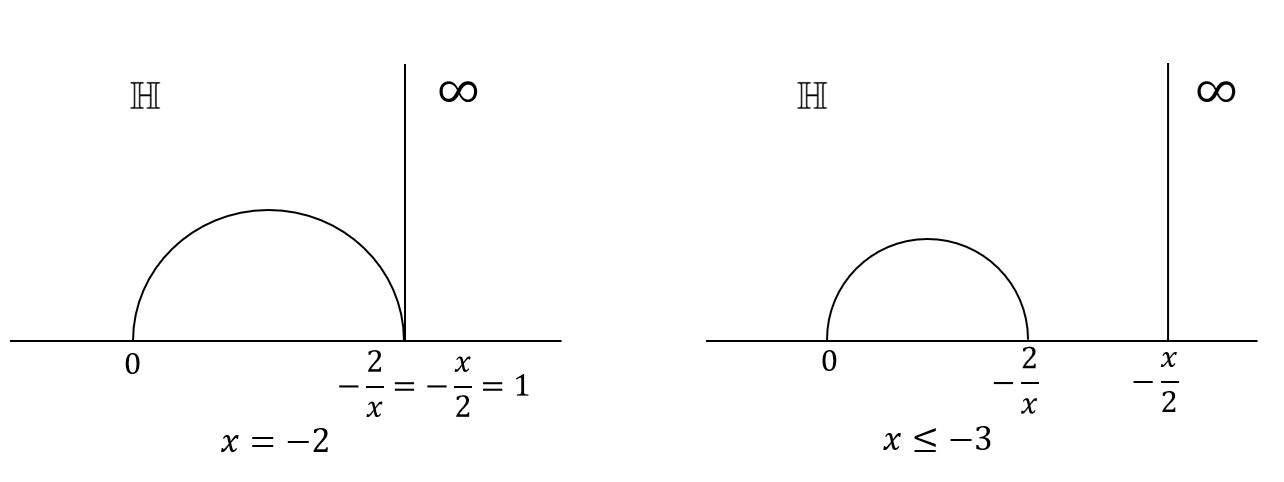}
\caption[Figure 6.4]{Fundamental domain in $\H$ of
$\Gamma^{(0),mat}$ for $x=-2$ and $x\leq -3$}
\label{Fig:6.4}
\end{figure}

Therefore $\mu(s_{e_1}^{(0),mat})$ and $\mu(s_{e_2}^{(0),mat})
\in\Isom(\H)$ are reflections along the hyperbolic line
$A(\infty,-\frac{x}{2})$ respectively $A(0,-\frac{2}{x})$.
As $x\leq -2$, we have $-\frac{2}{x}\leq -\frac{x}{2}$, so
$A(0,-\frac{2}{x})\cap A(-\frac{x}{2},\infty)=\emptyset$, see the
pictures in Figure \ref{Fig:6.4}. 

Theorem \ref{ta.2} (a) applies and shows that
$\langle \mu(s_{e_1}^{(0),mat}),\mu(s_{e_2}^{(0),mat})\rangle
\subset \Isom(\H)$ is a free Coxeter group with the two
given generators. 
Therefore also $\Gamma^{(0)}$ is a free Coxeter group with the 
two generators $s_{e_1}^{(0)}$ and $s_{e_2}^{(0)}$.

(c) The set $\Delta^{(0)}$ for $x=-2$: 
\begin{eqnarray*}
R^{(0)}&=& \{y_1e_1+y_2e_2\in H_\Z\,|\, 1=(y_1-y_2)^2\}\\
&=& (e_1+\Z f_1)\ \dot\cup\ (e_2+\Z f_1).
\end{eqnarray*}
For $m\in\Z$, $\varepsilon\in\{\pm 1\}$, 
\begin{eqnarray*}
s_{e_1}^{(0)}(\varepsilon e_1+mf_1)&=& -\varepsilon e_1+mf_1,\\
s_{e_2}^{(0)}(\varepsilon e_2+mf_1)&=& -\varepsilon e_2+mf_1,\\
s_{e_1}^{(0)}s_{e_2}^{(0)}(e_1+mf_1,e_2+mf_1)
&=&(e_1+(m+2)f_1,e_2+(m-2)f_1).
\end{eqnarray*}
This shows all claims on $\Delta^{(0)}$ in part (c).

(d) (i) Recall $\kappa_1+\kappa_2=-x$, $\kappa_1\kappa_2=1$, 
$0=\kappa_i^2+x\kappa_i+1$, $\kappa_2=\kappa_1^{-1}=-\kappa_1-x$,
$\kappa_1=-\kappa_2-x$. 

Because of $R^{(0)}=\{y_1e_1+y_2e_2\in H_\Z\,|\, 
1=y_1^2+xy_1y_2+y_2^2\}$ the map
\begin{eqnarray*}
u:R^{(0)}&\to&\{\textup{the units with norm }1\textup{ in }
\Z[\kappa_1]\}\\
y_1e_1+y_2e_2&\mapsto& y_1-\kappa_1y_2
\end{eqnarray*}
is well defined and a bijection. Because of Lemma \ref{tc.1} (a)
\begin{eqnarray*}
\{\textup{the units with norm }1\textup{ in }\Z[\kappa_1]\}
=\{\pm\kappa_1^l\,|\, l\in\Z\}.
\end{eqnarray*}
Now
\begin{eqnarray*}
s_{e_1}^{(0)}(u^{-1}(\varepsilon \kappa_1^l))
=u^{-1}(-\varepsilon\kappa_1^{-l})\quad\textup{and}\quad
s_{e_2}^{(0)}(u^{-1}(\varepsilon \kappa_1^l))
=u^{-1}(-\varepsilon\kappa_1^{2-l})
\end{eqnarray*}
follow from
\begin{eqnarray*}
s_{e_1}^{(0)}\uuuu{e}\begin{pmatrix}y_1\\y_2\end{pmatrix}
=\uuuu{e}\begin{pmatrix}-1&-x\\0&1\end{pmatrix}
\begin{pmatrix}y_1\\y_2\end{pmatrix}
=\uuuu{e}\begin{pmatrix}-y_1-xy_2\\y_2\end{pmatrix},\\
(-y_1-xy_2)-\kappa_1y_2 = -(y_1-\kappa_2y_2)
=-(y_1-\kappa_1y_2)^{-1},\\
s_{e_2}^{(0)}\uuuu{e}\begin{pmatrix}y_1\\y_2\end{pmatrix}
=\uuuu{e}\begin{pmatrix}1&0\\-x&-1\end{pmatrix}
\begin{pmatrix}y_1\\y_2\end{pmatrix}
=\uuuu{e}\begin{pmatrix}y_1\\-xy_1-y_2\end{pmatrix},\\
y_1-\kappa_1(-xy_1-y_2) = \kappa_1((\kappa_2+x)y_1+y_2)
=\kappa_1(-\kappa_1y_1+y_2)\\
=-\kappa_1^2(y_1-\kappa_2y_2)
=-\kappa_1^2(y_1-\kappa_1y_2)^{-1}.
\end{eqnarray*}
This shows $\Delta^{(0)}=R^{(0)}$ and that $\Delta^{(0)}$
splits into the two disjoint orbits 
$\Gamma^{(0)}\{e_1\}$ and $\Gamma^{(0)}\{e_2\}$. 

(ii) The formulas in part (i) show immediately 
$(s_{e_1}^{(0)}s_{e_2}^{(0)})^m(u^{-1}(\varepsilon\kappa_1^l))
=u^{-1}(\varepsilon\kappa_1^{l-2m})$ for $l,m\in\Z$,
$\varepsilon\in\{\pm 1\}$.
Together with $u(e_1)=1$ and $u(e_2)=-\kappa_1$ 
this implies the formulas in part (ii). 

(iv) This follows from the formulas for the action of
$\www{s}_{e_i}^{(0)}$ on $\www{R}^{(0)}$ and of
$s_{e_i}^{(0)}$ on $R^{(0)}$. 

(iii) First we consider the lower right branch of the hyperbola.
There $y_1-\kappa_1 y_2>0$ and $y_1-\kappa_2y_2>0$. We consider
$y_2$ as an implicit function in $y_1$. The equation 
\begin{eqnarray*}
1=y_1^2+xy_1y_2+y_2^2=(y_1-\kappa_1y_2)(y_1-\kappa_2y_2)
\end{eqnarray*}
implies
\begin{eqnarray*}
0&=&(y_1-\kappa_1y_2)(1-\kappa_2y_2')
+(y_1-\kappa_2y_2)(1-\kappa_1y_2')\\
&=& [(y_1-\kappa_1y_2)+(y_1-\kappa_2y_2)]
-[(y_1-\kappa_1y_2)\kappa_2+(y_1-\kappa_2y_2)\kappa_1]y_2',\\
\textup{so}&&y_2'>0\quad\textup{and}\quad 
(1-\kappa_1y_2')(1-\kappa_2y_2')<0,\\
0&=& (y_1-\kappa_1y_2)(-\kappa_2y_2'')
+2(1-\kappa_1y_2')(1-\kappa_2y_2')\\
&&+(y_1-\kappa_2y_2)(-\kappa_1y_2'')\\
&=&-[(y_1-\kappa_1y_2)\kappa_2+(y_1-\kappa_2y_2)\kappa_1]y_2''
+2(1-\kappa_1y_2')(1-\kappa_2y_2'),\\
\textup{so}&&y_2''<0.
\end{eqnarray*}
Therefore the lower right branch of the hyperbola is strictly
monotonously increasing and concave. 
The upper left branch is obtained from the lower right branch
by the reflection $H_\R\to H_\R,\ (y_1,y_2)\mapsto (y_2,y_1),$
along the diagonal. Therefore it is strictly monotonously 
increasing and convex.

By definition $s_{e_1}^{(0)}$ maps each horizontal line in $H_\R$
to itself, and $s_{e_2}^{(0)}$ maps each vertical line in $H_\R$
to itself. As they map $\Delta^{(0)}=R^{(0)}$ to itself,
this shows all statements in (iii). 

(e) Obviously $g_{1,2}\in O^{(0)}$ and 
$\{\pm \id,\pm g_{1,2}\}$ is a subgroup of $O^{(0)}$. 
Recall
\begin{eqnarray*}
H_\Z^\sharp \supset j^{(0)}(H_\Z)
=\langle j^{(0)}(e_1),j^{(0)}(e_2)\rangle\\
\textup{with}\quad 
j^{(0)}(e_1)(\uuuu{e})=(2,x),j^{(0)}(e_2)(\uuuu{e})=(x,2).
\end{eqnarray*}
Define $l\in H_\Z^\sharp$ with $l(\uuuu{e})=(1,0)$. Then
\begin{eqnarray*}
(l-l\circ(-\id)^{-1})(\uuuu{e})=2l(\uuuu{e})=(2,0),\\
\textup{so}\quad l-l\circ(-\id)^{-1}\notin j^{(0)}(H_\Z),
\textup{so}\quad -\id\notin\ker\tau^{(0)}.\\
(l-l\circ g_{1,2}^{-1})(\uuuu{e})=(1,-1),\\
\textup{so}\quad l-l\circ g_{1,2}^{-1}\notin j^{(0)}(H_\Z),
\textup{so}\quad g_{1,2}\notin\ker\tau^{(0)}.
\end{eqnarray*}
Denote for a moment by $\www{O}^{(0)}$ the subgroup of 
$O^{(0)}$ which is generated by $\{\pm\id,\pm g_{1,2}\}$
and $\Gamma^{(0)}$. We just saw
\begin{eqnarray*}
(\ker\tau^{(0)})\cap \www{O}^{(0)} =\Gamma^{(0)},\\
\textup{so}\quad 
\www{O}^{(0)}=\Gamma^{(0)}\rtimes \{\pm \id,\pm g_{1,2}\}.
\end{eqnarray*}
It remains to show that this subgroup is $O^{(0)}$. 

By the parts (c) and (d) (iii) the set $\Delta^{(0)}=R^{(0)}$
splits into two $\Gamma^{(0)}$ orbits,
$\Gamma^{(0)}\{e_1\}$ and $\Gamma^{(0)}\{e_2\}$.
The element $g_{1,2}$ interchanges them, so 
$\Delta^{(0)}=R^{(0)}$ is a single $\www{O}^{(0)}$ orbit.

Therefore each element of $O^{(0)}$ can be written as a product
of an element in $\www{O}^{(0)}$ and an element $g\in O^{(0)}$
with $g(e_1)=e_1$. It is sufficient to show $g\in\www{O}^{(0)}$.
Observe that the set $\{v\in H_\Z\,|\, I^{(0)}(v,v)>0\}$
consists of two components, with $e_1$ in one component
and $e_2$ in the other component. Because of $g(e_1)=e_1$,
$g(e_2)$ is in the same component as $e_2$. Now
$H_\Z=\Z e_1+\Z g(e_2)$ shows 
$g(e_2)\in\{e_2,s_{e_1}^{(0)}(-e_2)\}$ and 
$g\in\{\id,-s_{e_1}^{(0)}\}$. Therefore 
$\www{O}^{(0)}=O^{(0)}$. 
\hfill$\Box$

\begin{remarks}\label{t6.9}
(i) By part (iv) of Theorem \ref{t6.8} (d) the pairs
$(\Gamma^{(0)},\Delta^{(0)})$ with the action of $\Gamma^{(0)}$
on $\Delta^{(0)}$ are isomorphic for all $x\leq -2$.
This is interesting as in the case $x=-2$ the set $\Delta^{(0)}$
and this action could be written down in a very simple way.

The parts (i) and (iii) of Theorem \ref{t6.8} (d) offered two
ways to control the set $\Delta^{(0)}$ and this action also for
$x\leq -3$, a number theoretic way and a geometric way.
But both ways are less simple than $\Delta^{(0)}$ in the
case $x=-2$.

(ii) In the odd cases the situation will be partly similar, partly
different. The pairs $(\Gamma^{(1)},\Delta^{(1)})$ with the action
of $\Gamma^{(1)}$ on $\Delta^{(1)}$ are isomorphic for all
$x\leq -2$. In the case $x=-2$ the set $\Delta^{(1)}$ and this
action can be written down in a fairly simple way.
But we lack analoga of the parts (i) and (iii) in Theorem
\ref{t6.8} (d). We do not have a good control on the sets 
$\Delta^{(1)}$ for $x\leq -3$. 

(iii) For each $x\leq -2$ and $i\in\{1,2\}$ 
$\Stab_{\Gamma^{(0)}}(e_i)=\{\id\}$, so the map
$\Gamma^{(0)}\to \Gamma^{(0)}\{e_i\}$, 
$\gamma\mapsto \gamma(e_i)$,
is a bijection. The action of $\Gamma^{(0)}$ on 
$\Gamma^{(0)}\{e_i\}$ shows again immediately that $\Gamma^{(0)}$
is a free Coxeter group with generators $s_{e_1}^{(0)}$
and $s_{e_2}^{(0)}$. 
\end{remarks}

Now we come to the odd cases. As before we restrict to $x\in \Z_{<0}$.

\begin{theorem}\label{t6.10}
(a) We have 
\begin{eqnarray*}
O^{(1)}&\cong & SL_2(\Z),\\
\ker\tau^{(1)}&\cong& \Gamma(x):=\{A\in SL_2(\Z)\,|\, 
A\equiv E_2\mmod x\},
\end{eqnarray*}
\begin{eqnarray*}
s^{(1)}_{e_i}\uuuu{e}&=&\uuuu{e}\cdot s_{e_i}^{(1),mat}
\quad\textup{ with }\\
s^{(1),mat}_{e_1}&=&\begin{pmatrix}1&-x\\0&1\end{pmatrix},\quad
s^{(1),mat}_{e_2}=\begin{pmatrix}1&0\\x&1\end{pmatrix},\\
\Gamma^{(1)}&\cong& \Gamma^{(1),mat}:=
\langle s^{(1),mat}_{e_1},s^{(1),mat}_{e_2}\rangle 
\subset SL_2(\Z),
\end{eqnarray*}
The map from $\Delta^{(1)}$ to its image in 
$\widehat{\R}=\R\cup\{\infty\}$ under the composition
$C:\Delta^{(1)}\to\widehat{\R}$ of maps
\begin{eqnarray*}
\Delta^{(1)}&\to& (\R^*\Delta^{(1)})/\R^* 
=\{\textup{lines through vanishing cycles}\}\\
&\hookrightarrow& (H_\R-\{0\})/\R^*
=\{\textup{lines in }H_\R\}
\stackrel{\cong}{\longrightarrow}\widehat{\R},\\
&&\R^*\uuuu{e}\begin{pmatrix}y_1\\1\end{pmatrix}\mapsto y_1,\quad
\R^*e_1\mapsto\infty,
\end{eqnarray*}
is two-to-one.

(b) The case $x=-1$: $\Gamma^{(1),mat}=SL_2(\Z)$.
\begin{eqnarray*}
\Delta^{(1)}=\Gamma^{(1)}\{e_1\}
= \{y_1e_1+y_2e_2\in H_\Z\,|\, \gcd(y_1,y_2)=1\}
=H_\Z^{prim},
\end{eqnarray*}
where $H_\Z^{prim}$ denotes the set of primitive vectors in $H_\Z$.
The image $C(\Delta^{(1)})\subset\widehat{\R}$ is
$\widehat{\Q}=\Q\cup\{\infty\}$. 

(c) The case $x=-2$: 
$\Gamma^{(1)}\cong G^{free,2}$ 
is a free group with the two generators 
$s_{e_1}^{(1)}$ and $s_{e_2}^{(1)}$.
The (isomorphic) matrix group $\Gamma^{(1),mat}$ is 
\begin{eqnarray*}
\Gamma^{(1),mat}&=& \{\begin{pmatrix}a&b\\c&d\end{pmatrix}
\in SL_2(\Z)\,|\, a\equiv d\equiv 1(4),\ b\equiv c\equiv 0(2)\}.
\end{eqnarray*}
It is a subgroup of index 2 in the principal congruence subgroup
\begin{eqnarray*}
\Gamma(2)=\{\begin{pmatrix}a&b\\c&d\end{pmatrix}\,|\,
a\equiv d\equiv 1(2),\ b\equiv c\equiv 0(2)\}
\end{eqnarray*}
of $SL_2(\Z)$. 
The set $\Delta^{(1)}=\Gamma^{(1)}\{\pm e_1,\pm e_2\}$ is
$$\Delta^{(1)}= \{y_1e_1+y_2e_2\in H_\Z^{prim}\,|\, 
y_1+y_2\equiv 1(2)\}.$$
It splits into
the four disjoint orbits
\begin{eqnarray*}
\Gamma^{(1)}\{e_1\}&=& \{y_1e_1+y_2e_2\in H_\Z^{prim}\,|\, 
y_1\equiv 1(4),y_2\equiv 0(2)\}\\
\Gamma^{(1)}\{-e_1\}&=& \{y_1e_1+y_2e_2\in H_\Z^{prim}\,|\, 
y_1\equiv 3(4),y_2\equiv 0(2)\}\\
\Gamma^{(1)}\{e_2\}&=& \{y_1e_1+y_2e_2\in H_\Z^{prim}\,|\, 
y_1\equiv 0(2),y_2\equiv 1(4)\}\\
\Gamma^{(1)}\{-e_2\}&=& \{y_1e_1+y_2e_2\in H_\Z^{prim}\,|\,
y_1\equiv 0(2),y_2\equiv 3(4)\}. 
\end{eqnarray*}
The set $H_\Z^{prim}$ of primitive vectors is the disjoint union
of $\Delta^{(1)}$ and the set 
\begin{eqnarray*}
\{y_1e_1+y_2e_2\in H_\Z^{prim}\,|\,y_1\equiv y_2\equiv 1(2)\}.
\end{eqnarray*}
The image $C(\Delta^{(1)})\subset\widehat{\R}$ is
\begin{eqnarray*}
\{\infty\}\cup\{\frac{a}{b}\,|\, a\in\Z,b\in\N,
\gcd(a,b)=1,a\equiv 0(2)\textup{ or }b\equiv 0(2)\}
\subset\widehat{\Q},
\end{eqnarray*}
and is dense in $\widehat{\R}$.

(d) The cases $x\leq -3$. 

(i) $\Gamma^{(1)}\cong G^{free,2}$ 
is a free group with the two generators 
$s_{e_1}^{(1)}$ and $s_{e_2}^{(1)}$.

(ii) The matrix group $\Gamma^{(1),mat}$ is a Fuchsian group of
\index{Fuchsian group}
the second kind. It has infinite index in the group
\begin{eqnarray*}
\{\begin{pmatrix}a&b\\c&d\end{pmatrix}\in SL_2(\Z)
\,|\, a\equiv d\equiv 1(x^2),b\equiv c\equiv 0(x)\},
\end{eqnarray*}
which has finite index in $SL_2(\Z)$.

(iii) The image $C(\Delta^{(1)})\subset\widehat{\R}$ is a subset
of $\widehat{\Q}$ which is nowhere dense in $\widehat{\R}$.

(iv) Denote by $\www{s}_{e_i}^{(1)}$, $\www{\Gamma}^{(1)}$ and
$\www{\Delta}^{(1)}$ the objects for $x=-2$ and as before by
$s_{e_i}^{(1)}$, $\Gamma^{(1)}$ and $\Delta^{(1)}$ the objects for
an $x\leq -3$.

The map $\www{s}_{e_i}^{(1)}\mapsto s_{e_i}^{(1)}$ extends to a 
group isomorphism $\www{\Gamma}^{(1)}\to\Gamma^{(1)}$,
which maps the stabilizer $\langle \www{s}_{e_i}^{(1)}\rangle$
of $e_i$ in $\www{\Gamma}^{(1)}$ to the stabilizer 
$\langle s_{e_i}^{(1)}\rangle$ of $e_i$ in $\Gamma^{(1)}$.
The induced map
\begin{eqnarray*}
\www{\Delta}^{(1)}\to\Delta^{(1)},\ 
\www{\gamma}(\varepsilon e_i)\mapsto \gamma(\varepsilon e_i)
\quad\textup{for }\varepsilon\in\{\pm 1\},
\www{\gamma}\mapsto\gamma,
\end{eqnarray*}
is a bijection. The bijections
$\www{\Gamma}^{(1)}\to \Gamma^{(1)}$ and 
$\www{\Delta}^{(1)}\to \Delta^{(1)}$ are compatible with the
actions of $\www{\Gamma}^{(1)}$ on $\www{\Delta}^{(1)}$
and of $\Gamma^{(1)}$ on $\Delta^{(1)}$. 
Especially, $\Delta^{(1)}$ splits into the four disjoint orbits
$\Gamma^{(1)}\{e_1\}$, $\Gamma^{(1)}\{-e_1\}$, 
$\Gamma^{(1)}\{e_2\}$, $\Gamma^{(1)}\{-e_2\}$. 

(e) In the case $A_1^2$ $\Delta^{(0)}=\Delta^{(1)}=\{\pm e_1,\pm e_2\}$.
In all other rank 2 cases $\Delta^{(0)}\subsetneqq \Delta^{(1)}$.
\end{theorem}

{\bf Proof:}
(a) Because of $I^{(1)}(\uuuu{e}^t,\uuuu{e})
=\begin{pmatrix}0&x\\-x&0\end{pmatrix}$ we have
$O^{(1)}\cong SL_2(\Z)$. In order to see $\ker \tau^{(1)}\cong \Gamma(x)$,
consider the generators $l_1,l_2\in H_\Z^\sharp$ of $H_\Z^\sharp$
with $l_1(\uuuu{e})=(1,0)$ and $l_2(\uuuu{e})=(0,1)$. 
Observe first
$$j^{(1)}(e_1)(\uuuu{e})=(0,x),\ 
j^{(1)}(e_2)(\uuuu{e})=(-x,0),
\quad\textup{so }j^{(1)}(H_\Z)=x H_\Z^\sharp,$$
and second that $g\in O^{(1)}$ with 
$g^{-1}(\uuuu{e})=\uuuu{e}\cdot\begin{pmatrix}a&b\\c&d\end{pmatrix}$
satisfies 
\begin{eqnarray*}
(l_1-l_1\circ g^{-1})(\uuuu{e})&=&(1-a,b),\\
(l_2-l_2\circ g^{-1})(\uuuu{e})&=&(-c,1-d),
\end{eqnarray*}
so $g\in \ker \tau^{(1)}$ if and only if 
$\begin{pmatrix}a&b\\c&e\end{pmatrix}\equiv E_2\mmod x$. 

It remains to prove that the line 
$\R\cdot \delta\subset H_\R$
through a vanishing cycle $\delta\in\Delta^{(1)}$ intersects
$\Delta^{(1)}$ only in $\pm \delta$. To prove this we can
restrict to $\delta=e_i$. There it follows from the fact that
any matrix $A\in SL_2(\Z)$ with a zero in an entry $A_{ij}$
has entries $A_{i,j+1(2)}$, $A_{i+1(2),j}\in\{\pm 1\}$. 

(b) $\Gamma^{(1),mat}=SL_2(\Z)$ is well known.
The standard arguments for this are as follows.
The group $\mu(\Gamma^{(1),mat})\subset\Isom(\H)$ is generated
by 
\begin{eqnarray*}
\mu(s_{e_1}^{(1),mat})= \mu(\begin{pmatrix}1&1\\0&1\end{pmatrix})
&=&(z\mapsto z+1),\\
\mu(s_{e_1}^{(1),mat}s_{e_2}^{(1),mat}s_{e_1}^{(1),mat})
=\mu(\begin{pmatrix}0&1\\-1&0\end{pmatrix})&=&(z\mapsto -z^{-1}).
\end{eqnarray*}
One sees almost immediately that it acts transitively on 
$\widehat{\Q}$ and that the stabilizer 
$\langle \mu(s_{e_1}^{(1),mat})\rangle$ of $\infty$ in
$\mu(\Gamma^{(1),mat})$ coincides with the stabilizer of $\infty$
in $\mu(SL_2(\Z))$. Therefore
$\mu(\Gamma^{(1),mat})=\mu(SL_2(\Z))$.
But $-E_2\in \Gamma^{(1),mat}$ because of
$\begin{pmatrix}0&1\\-1&0\end{pmatrix}^2=-E_2$. 
Therefore $\Gamma^{(1),mat}=SL_2(\Z)$. 

The fact that $\mu(SL_2(\Z))$ acts transitively on $\widehat{\Q}$
shows $C(\Delta^{(1)})=\widehat{\Q}$. Together with
$-E_2\in SL_2(\Z)$ this shows
\begin{eqnarray*}
\Delta^{(1)}=\Gamma^{(1)}\{e_1\}=H_\Z^{prim}.
\end{eqnarray*}

\begin{figure}[H]
\includegraphics[width=0.8\textwidth]{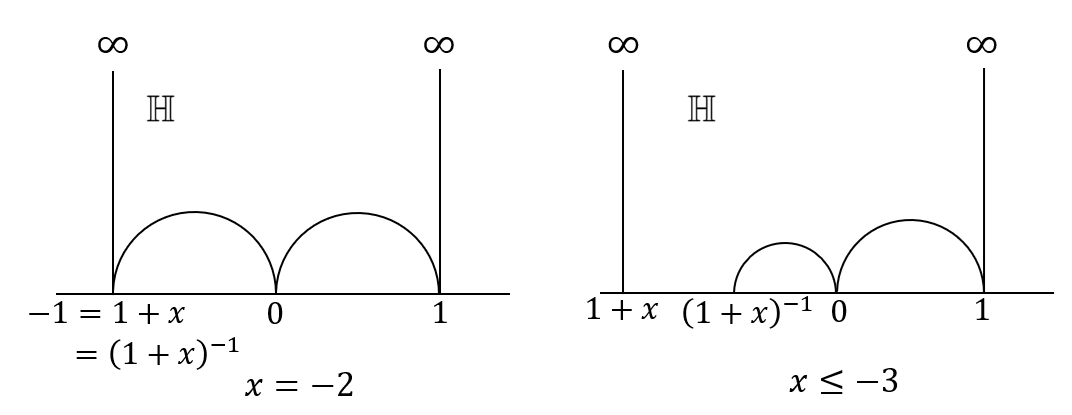}
\caption[Figure 6.5]{Fundamental domain in $\H$ of 
$\Gamma^{(1),mat}$ for $x=-2$ and $x\leq -3$}
\label{Fig:6.5}
\end{figure}

(c) and (d) The group $\Gamma^{(1)}$ for $x\leq -2$:
Recall the Remarks and Notations \ref{ta.1}. The elements
\begin{eqnarray*}
\mu(s_{e_1}^{(1),mat})=\mu(\begin{pmatrix}1&-x\\0&1\end{pmatrix})
&=&(z\mapsto z-x)\quad\textup{and}\\
\mu(s_{e_2}^{(1),mat})=\mu(\begin{pmatrix}1&0\\x&1\end{pmatrix})
&=&(z\mapsto \frac{z}{xz+1})
\end{eqnarray*}
of $\Isom(\H)$ 
are parabolic with fixed points $\infty$ respectively $0$ on
$\widehat{\R}$. Observe
\begin{eqnarray*}
\mu(s_{e_1}^{(1),mat})^{-1}(1)=1+x,\quad
\mu(s_{e_2}^{(1),mat})(1)=(1+x)^{-1},\\
(1+x)^{-1}\geq 1+x\quad\textup{for }x\leq -2.
\end{eqnarray*}
Therefore $\mu(s_{e_1}^{(1),mat})^{-1}(A(\infty,1))=A(\infty,1+x)$
and $\mu(s_{e_2}^{(1),mat})(A(0,1))=A(0,(1+x)^{-1})$ 
do not intersect. See Figure \ref{Fig:6.5}

Theorem \ref{ta.2} (c) applies and shows that 
$\mu(\Gamma^{(1),mat})$ is a free group with the two generators 
$\mu(s_{e_1}^{(1),mat})$ and $\mu(s_{e_2}^{(1),mat})$.
Therefore also $\Gamma^{(1)}$ is a free group with the two
generators $s_{e_1}^{(1)}$ and $s_{e_2}^{(1)}$. 
As therefore the map $\Gamma^{(1),mat}\to \mu(\Gamma^{(1),mat})$
is an isomorphism, $-E_2\notin\Gamma^{(1),mat}$ and
$-\id\notin \Gamma^{(1)}$. 

Theorem \ref{ta.2} (c) says also that the contractible open
set $\FF$ whose hyperbolic boundary consists of the four 
hyperbolic lines which were used above,
$A(1+x,\infty)$, $A(\infty,1)$, $A(1,0)$, $A(0,(1+x)^{-1})$,
is a fundamental domain for $\mu(\Gamma^{(1),mat})$.
If $x=-2$ then its euclidean boundary in $\widehat{\C}$ consists
of these four hyperbolic lines and the four points
$\infty$, $1$, $0$, $-1=1+x=(1+x)^{-1}$. 
If $x\leq -3$ then its euclidean boundary consists of these four
hyperbolic lines, the three points $\infty$, $1$, $0$, and
the interval $[1+x,(1+x)^{-1}]$.

(c) $\Gamma^{(1),mat}$ and $\Delta^{(1)}$ for $x=-2$:
The following facts together imply 
$\mu(\Gamma^{(1),mat})=\mu(\Gamma(2))$:
\begin{eqnarray*}
\Gamma^{(1),mat}\subset\Gamma(2),\quad
[SL_2(\Z):\Gamma(2)]=6,\quad -E_2\in \Gamma(2),\\
(\textup{hyperbolic area of the fundamental domain }\FF
\textup{ of }\mu(\Gamma^{(1),mat}))=2\pi,\\
(\textup{hyperbolic area of a fundamental domain of }
\mu(SL_2(\Z)))=\frac{\pi}{3}.
\end{eqnarray*}
Therefore either $\Gamma^{(1),mat}=\Gamma(2)$ or 
$\Gamma^{(1),mat}$ is a subgroup of index 2 in $\Gamma(2)$.
But $\Gamma^{(1),mat}$ is certainly a subgroup of the subgroup
\begin{eqnarray*}
\begin{pmatrix}a&b\\c&d\end{pmatrix}\in SL_2(\Z)\,|\, 
a\equiv d\equiv 1(4),b\equiv c\equiv 0(2)\}
\end{eqnarray*}
of $\Gamma(2)$ and does not contain $-E_2$. Therefore
$\Gamma^{(1),mat}$ coincides with this subgroup of $\Gamma(2)$
and has index 2 in $\Gamma(2)$. 

Therefore the orbits $\Gamma^{(1)}\{e_1\}$, 
$\Gamma^{(1)}\{-e_1\}$, $\Gamma^{(1)}\{e_2\}$, 
$\Gamma^{(1)}\{-e_2\}$ are contained in the right hand sides
of the equations in part (c) which describe them and are
disjoint. It remains to show equality.

We restrict to $\Gamma^{(1)}\{e_1\}$. The argument for
$\Gamma^{(1)}\{e_2\}$ is analogous, and the equations for
$\Gamma^{(1)}\{-e_1\}$ and $\Gamma^{(1)}\{-e_2\}$ follow
immediately. 

Suppose $y_1,y_3\in\Z$ with $y_1\equiv 1(4)$, $y_3\equiv  0(2)$,
$\gcd(y_1,y_3)=1$. We have to show 
$\uuuu{e}\begin{pmatrix}y_1\\y_3\end{pmatrix}
\in \Gamma^{(1)}\{e_1\}$. For that we have to find 
$y_2,y_4\in\Z$ with
$\begin{pmatrix}y_1&y_2\\y_3&y_4\end{pmatrix}\in \Gamma^{(1)}$,
so with $1=y_1y_4-y_2y_3$, $y_4\equiv 1(4)$, $y_2\equiv 0(2)$.

The condition $\gcd(y_1,y_3)$ implies existence of
$\www{y}_2,\www{y}_4\in\Z$ with $1=y_1\www{y}_4-\www{y}_2y_3$.

1st case, $\www{y_2}\equiv 0(2)$: Then 
$1=y_1\www{y}_4-\www{y}_2y_3$ shows $\www{y}_4\equiv 1(4)$,
and $(y_2,y_4)=(\www{y}_2,\www{y}_4)$ works.

2nd case, $\www{y}_2\equiv 1(2)$: Then
$(y_2,y_4)=(\www{y}_2+y_1,\www{y}_4+y_3)$ satisfies
$1=y_1y_4-y_2y_3$ and $y_2\equiv 0(2)$, so we are in the 1st case,
so $(y_2,y_4)$ works.

Therefore the orbit $\Gamma^{(1)}\{e_1\}$ is as claimed the set
$\{y_1e_1+y_2e_2\in H_\Z^{prim}\,|\, 
y_1\equiv 1(4),y_2\equiv 0(2)\}$. 

The statements on $H_\Z^{prim}$ and $C(\Delta^{(1)})$ are clear
now, too.

(d) $\Gamma^{(1),mat}$ and $\Delta^{(1)}$ for $x\leq -3$:
Part (i) was shown above.

(ii) The euclidean boundary of the fundamental domain $\FF$
above of $\mu(\Gamma^{(1),mat})$ contains the real interval
$[(1+x)^{-1},1+x]$. Therefore its hyperbolic area is $\infty$,
$\Gamma^{(1),mat}$ is a Fuchsian group of the second kind,
and the index of $\Gamma^{(1),mat}$ in $SL_2(\Z)$ is $\infty$
(see e.g. \cite[34.]{Fo51} \cite[\S 5.3, \S 8.1]{Be83}).

(iii) The set $C(\Delta^{(1)})\subset\widehat{\R}$ is the union 
of the $\mu(\Gamma^{(1),mat})$ orbits of $\infty$ and $0$ in
$\widehat{\R}$. Because $\Gamma^{(1),mat}$ is a Fuchsian group
of the second kind, these two orbits are nowhere dense in 
$\widehat{\R}$ (see e.g. \cite{Fo51} \cite{Be83}). 
For example, they contain no point of the open interval 
$((1+x)^{-1},(1+x))$ and of its $\mu(\Gamma^{(1),mat})$ orbit.

(iv) The groups $\www{\Gamma}^{(1)}$ and $\Gamma^{(1)}$ are
free groups with generators 
$\www{s}_{e_1}^{(1)}$, $\www{s}_{e_2}^{(1)}$ and
$s_{e_1}^{(1)}$, $s_{e_2}^{(1)}$. 
Therefore there is a unique isomorphism
$\www{\Gamma}^{(1)}\to\Gamma^{(1)}$ with
$\www{s}_{e_i}^{(1)}\mapsto s_{e_i}^{(1)}$.
The other statements follow immediately.

(e) The case $x=0$, so $A_1^2$: 
$$\Delta^{(0)}=\Delta^{(1)}=\{\pm e_1,\pm e_2\}.$$

The case $x=-1$, so $A_2$: 
$$\Delta^{(0)}=\{\pm e_1,\pm e_2,\pm (e_1+e_2)\}\subsetneqq
H_\Z^{prim}=\Delta^{(1)}.$$

The case $x=-2$, so $\P^1A_1$: 
\begin{eqnarray*}
\Delta^{(0)}&=&(e_1+\Z f_1)\, \dot\cup\, (e_2+\Z f_1)\\
&\subsetneqq& \{y_1e_1+y_2e_2\in H_\Z^{prim}\,|\, y_1+y_2\equiv 1(2)\}
=\Delta^{(1)}.
\end{eqnarray*}

The cases $x\leq -3$: Recall
$s_{e_1}^{(0)}s_{e_2}^{(0)}=-M=-s_{e_1}^{(1)}s_{e_2}^{(1)}$.
Therefore for $m\in\Z,\varepsilon\in\{\pm 1\}$ 
\begin{eqnarray*}
(s_{e_1}^{(0)}s_{e_2}^{(0)})^m(\uuuu{e})&=&
(-s_{e_1}^{(1)}s_{e_2}^{(1)})^m(\uuuu{e})\in(\Delta^{(1)})^2,\\
(s_{e_1}^{(0)}s_{e_2}^{(0)})^m s_{e_1}^{(0)}(e_1)&=&
-(-s_{e_1}^{(1)}s_{e_2}^{(1)})^m(e_1)\in \Delta^{(1)},\\
(s_{e_1}^{(0)}s_{e_2}^{(0)})^m s_{e_1}^{(0)}(e_2)&=&
(s_{e_1}^{(0)}s_{e_2}^{(0)})^{m+1}s_{e_2}^{(0)}(e_2)\in \Delta^{(1)},\\
&=&-(-s_{e_1}^{(1)}s_{e_2}^{(1)})^{m+1}(e_2)\in \Delta^{(1)}.
\end{eqnarray*}
This shows $\Delta^{(0)}\subset\Delta^{(1)}$. For example
$(s_{e_1}^{(1)})^{-1}(e_2)=xe_1+e_2$ satisfies
$$L(xe_1+e_2,xe_1+e_2)=\begin{pmatrix}x&1\end{pmatrix}
\begin{pmatrix}1&0\\x&1\end{pmatrix}\begin{pmatrix}x\\1\end{pmatrix}
=2x^2+1\neq 1,$$
so $\Delta^{(0)}\subsetneqq\Delta^{(1)}$. 
\hfill$\Box$

\section{The even rank 3 cases}\label{s6.3}

For $\uuuu{x}=(x_1,x_2,x_3)\in\Z^3$ consider the matrix
$S=S(\uuuu{x})=
\begin{pmatrix}1&x_1&x_2\\0&1&x_3\\0&0&1\end{pmatrix}
\in T^{uni}_3(\Z)$, and consider a unimodular bilinear lattice
$(H_\Z,L)$ with a triangular basis $\uuuu{e}=(e_1,e_2,e_3)$
with $L(\uuuu{e}^t,\uuuu{e})^t=S$.
In this section we will determine in all cases the even
monodromy group $\Gamma^{(0)}=\langle s_{e_1}^{(0)},
s_{e_2}^{(0)},s_{e_3}^{(0)}\rangle$ and in many, but not all,
cases the set 
$\Delta^{(0)}=\Gamma^{(0)}\{\pm e_1,\pm e_2,\pm e_3\}$
of even vanishing cycles.
The cases where we control $\Delta^{(0)}$ well contain all 
cases with $r(\uuuu{x})\in\{0,1,2,4\}$ (3 does not turn up). 

The group $\Br_3\ltimes \{\pm 1\}^3$ acts on the set
$\BB^{tri}$ of triangular bases of $(H_\Z,L)$, but this action
does not change $\Gamma^{(0)}$ and $\Delta^{(0)}$.
Therefore the analysis of the action of $\Br_3\ltimes\{\pm 1\}^3$
on $T^{uni}_3(\Z)$ in Theorem \ref{t4.6} allows to restrict
to the matrices $S(\uuuu{x})$ with $\uuuu{x}$ in the following
list:
\begin{eqnarray*}
S(\uuuu{x})\textup{ with }\uuuu{x}\in\Z^3_{\leq 0}
\textup{ and }r(\uuuu{x})>4,\\
S(A_1^3),\ S(\P^2),\ S(A_2A_1),\ S(A_3),\ S(\P^1A_1),\
S(\whh{A}_2),\ S(\HH_{1,2}),\\
S(-l,2,-l)\textup{ for }l\geq 3,\\
\left\{\begin{array}{r}
S(\uuuu{x})\textup{ with }\uuuu{x}\in\Z^3_{\geq 3}
\textup{ and }r(\uuuu{x})<0\textup{ and }\\
x_i\leq \frac{1}{2}x_jx_k
\textup{ for }\{i,j,k\}=\{1,2,3\}.\end{array}\right.
\end{eqnarray*}
The following of these matrices satisfy 
$\uuuu{x}\in\Z^3_{\leq 0}$:
\begin{eqnarray*}
S(\uuuu{x})\textup{ with }\uuuu{x}\in\Z^3_{\leq 0}
\textup{ and }r(\uuuu{x})>4,\\
S(A_1^3),\ S(A_2A_1),\ S(A_3),\ S(\P^1A_1),\ S(\whh{A}_2).
\end{eqnarray*}
These are all Coxeter matrices. Their even monodromy groups
$\Gamma^{(0)}$ are Coxeter groups and are well known. The cases with
$\uuuu{x}\in\{0,1,2\}^3$ are classical,
the extension to all $\uuuu{x}\in\Z^3_{\leq 0}$ has been done
by Vinberg \cite[Prop. 6, Thm. 1, Thm. 2, Prop. 17]{Vi71}.
If we write $(x_1,x_2,x_3)=(S_{12},S_{13},S_{23})$ then the
following holds \cite[5.3+5.4]{Hu90} \cite{Vi71} 
\cite[4.1+4.2]{BB05}:
All relations in $\Gamma^{(0)}=\langle s_{e_1}^{(0)},
s_{e_2}^{(0)},s_{e_3}^{(0)}\rangle$ are generated by the
relations
\begin{eqnarray}\label{6.1}
(s_{e_i}^{(0)})^2=\id && \textup{for }i\in\{1,2,3\},\\
(s_{e_i}^{(0)}s_{e_j}^{(0)})^2=\id 
&&\textup{for }\{i,j,k\}=\{1,2,3\}\textup{ with }s_{ij}=0
\label{6.2}\\
{}[\textup{equivalent: }&&s_{e_i}^{(0)}\textup{ and }
s_{e_j}^{(0)}\textup{ commute}],\nonumber\\
(s_{e_i}^{(0)}s_{e_j}^{(0)})^3=\id 
&&\textup{for }\{i,j,k\}=\{1,2,3\}\textup{ with }s_{ij}=-1,
\label{6.3}\\
\textup{no relation }&&
\textup{for }\{i,j,k\}=\{1,2,3\}\textup{ with }s_{ij}\leq -2.
\label{6.4}
\end{eqnarray}
Especially, $\Gamma^{(0)}\cong G^{fCox,3}$ is a free Coxeter
group with three generators if $\uuuu{x}\in\Z^3_{\leq -2}$.

In Theorem \ref{t6.11} we recover this result, we say more
about the Coxeter groups with $r\in\{0,1,2,4\}$, so the
cases $S(A_1^3),S(A_2A_1),S(A_3),S(\P^1A_1),S(\whh{A}_2)$,
and we treat also the other cases where $\Gamma^{(0)}$ is not
a Coxeter group. 

The only cases where $\Rad I^{(0)}\supsetneqq\{0\}$ are the
cases with $r(\uuuu{x})=4$, so the cases
$S(\P^1A_1)$, $S(\whh{A}_2)$, $S(\HH_{1,2})$ and 
$S(-l,2,-l)$ with $l\geq 3$. In these cases, we have the
exact sequence 
\begin{eqnarray}\{1\}\to \Gamma^{(0)}_u\to\Gamma^{(0)}\to
\Gamma^{(0)}_s\to\{1\}\label{6.5}
\end{eqnarray}
in Lemma \ref{t6.2} (d).

\begin{theorem}\label{t6.11}
(a) We have
\begin{eqnarray*}
s_{e_i}^{(0)}\uuuu{e}&=&\uuuu{e}\cdot s_{e_i}^{(0),mat}
\quad\textup{with}\quad 
s_{e_1}^{(0),mat}=
\begin{pmatrix}-1&-x_1&-x_2\\0&1&0\\0&0&1\end{pmatrix},\\ 
s_{e_2}^{(0),mat}&=&
\begin{pmatrix}1&0&0\\-x_1&-1&-x_3\\0&0&1\end{pmatrix},\ 
s_{e_3}^{(0),mat}=
\begin{pmatrix}1&0&0\\0&1&0\\-x_2&-x_3&-1\end{pmatrix},\\
\Gamma^{(0)}&\cong&\Gamma^{(0),mat} :=\langle
s_{e_1}^{(0),mat},s_{e_2}^{(0),mat},s_{e_3}^{(0),mat}\rangle
\subset GL_3(\Z),\\
R^{(0)}&=&\{y_1e_1+y_2e_2+y_3e_3\in H_\Z\,|\, \\
&&1=y_1^2+y_2^2+y_3^2+x_1y_1y_2+x_2y_1y_3+x_3y_2y_3\}.
\end{eqnarray*}

(b) In the cases $S(\uuuu{x})$ with $\uuuu{x}\in\Z^3_{\leq 0}$
and $r(\uuuu{x})>4$ and in the reducible cases 
$S(A_1^3),S(A_2A_1),S(\P^1A_1)$, all relations in $\Gamma^{(0)}$
are generated by the relations in \eqref{6.1}--\eqref{6.4}.
Especially
\begin{eqnarray*}
\Gamma^{(0)}(A_1^3)&\cong&\Gamma^{(0)}(A_1)\times 
\Gamma^{(0)}(A_1)\times \Gamma^{(0)}(A_1)
\cong (G^{fCox,1})^3\cong \{\pm 1\}^3,\\
\Gamma^{(0)}(A_2A_1)&\cong& \Gamma^{(0)}(A_2)\times
\Gamma^{(0)}(A_1)\cong D_6\times \{\pm 1\}\cong
S_3\times \{\pm 1\},\\
\Gamma^{(0)}(\P^1A_1)&\cong& \Gamma^{(0)}(\P^1)\times
\Gamma^{(0)}(A_1)\cong G^{fCox,2}\times\{\pm 1\}.
\end{eqnarray*}

(c) In the case $S(A_3)$ the group $\Gamma^{(0)}$ is the 
Weyl group \index{Weyl group} 
of the root system $A_3$, so 
$\Gamma^{(0)}=\ker\tau^{(0)}=O^{(0),*}\cong S_4$.

(d) In the case $S(\whh{A}_2)$ the group $\Gamma^{(0)}$ is the
Weyl group of the affine root system $\whh{A}_2$. 
More concretely, the following holds. 
\begin{eqnarray*}
\Rad I^{(0)}&=&\Z f_1\textup{ with }f_1=e_1+e_2+e_3,\\
\oooo{H_\Z}^{(0)}&=&\Z\oooo{e_1}^{(0)}\oplus \Z\oooo{e_2}^{(0)},\\
\Gamma^{(0)}_u&=&\bigl(\ker \tau^{(0)}\bigr)_u
=T(\oooo{j}^{(0)}(\oooo{H_\Z}^{(0)})\otimes\Z f_1)\\
&=&\langle T(\oooo{j}^{(0)}(\oooo{e_1}^{(0)})\otimes f_1), 
T(\oooo{j}^{(0)}(\oooo{e_2}^{(0)})\otimes f_1)
\rangle\cong\Z^2 \quad\textup{ with}\\
&&T(\oooo{j}^{(0)}(\oooo{e_1}^{(0)})\otimes f_1)(\uuuu{e}) 
=\uuuu{e}+f_1(2,-1,-1),\\
&&T(\oooo{j}^{(0)}(\oooo{e_2}^{(0)})\otimes f_1)(\uuuu{e}) 
=\uuuu{e}+f_1(-1,2,-1),\\
\Gamma^{(0)}_u&=&\bigl(\ker \tau^{(0)}\bigr)_u
\stackrel{1:3}{\subset} O^{(0),Rad}_u
=T(\oooo{H_\Z}^{(0),\sharp}\otimes \Z f_1),\\
\Gamma^{(0)}_s&=& (\ker\tau^{(0)})_s \cong\Gamma^{(0)}(A_2)
\cong D_6\cong S_3,\\
\Gamma^{(0)}_s&\stackrel{1:2}{\subset}&
O^{(0),Rad}_s =\Aut(\oooo{H_\Z}^{(0)},
\oooo{I}^{(0)})\cong D_{12},\\
\Gamma^{(0)}&=&\ker \tau^{(0)}=O^{(0),*}
\stackrel{1:6}{\subset} O^{(0),Rad}.
\end{eqnarray*}
The exact sequence \eqref{6.5} splits non-canonically with
$\Gamma^{(0)}_s\cong\langle s_{e_1}^{(0)},s_{e_2}^{(0)}\rangle
\subset\Gamma^{(0)}$ (for example). 

(e) The case $S(\HH_{1,2})$: The following holds. 
\begin{eqnarray*}
H_\Z&=&\Z f_3\oplus \Rad I^{(0)}\textup{ with }
f_3=e_1+e_2+e_3,\\
\oooo{H_\Z}^{(0)}&=& \Z\oooo{f_3}^{(0)},\\
\Rad I^{(0)}&=&\Z f_1\oplus \Z f_2\textup{ with } 
f_1=e_1+e_2,\ f_2=e_2+e_3,\\
\Gamma^{(0)}_u&=&(\ker \tau^{(0)})_u
=T(\oooo{j}^{(0)}(\oooo{H_\Z}^{(0)})\otimes\Rad I^{(0)})\\
&=&\langle T(\oooo{j}^{(0)}(\oooo{f_3}^{(0)})\otimes f_1),
T(\oooo{j}^{(0)}(\oooo{f_3}^{(0)})\otimes f_2)\rangle\cong \Z^2
\quad \textup{with }\\
&&T(\oooo{j}^{(0)}(\oooo{f_3}^{(0)})\otimes f_1)(f_1,f_2,f_3)=(f_1,f_2,f_3+2f_1),\\
&&T(\oooo{j}^{(0)}(\oooo{f_3}^{(0)})\otimes f_2)(f_1,f_2,f_3)=(f_1,f_2,f_3+2f_2),\\
\Gamma^{(0)}_u &\stackrel{1:4}{\subset}& O^{(0),\Rad}_u
=T(\oooo{H_\Z}^{(0),\sharp}\otimes \Rad I^{(0)}),\\
\Gamma^{(0)}_s&=&(\ker \tau^{(0)})_s =O^{(0),Rad}_s
\cong \Gamma^{(0)}(A_1)\cong\{\pm 1\},\\
\Gamma^{(0)}&=&\ker \tau^{(0)}
=O^{(0),*}\stackrel{1:4}{\subset} O^{(0),Rad}.
\end{eqnarray*}
The exact sequence \eqref{6.5} splits non-canonically with
$\Gamma^{(0)}_s\cong\langle -M\rangle\subset
\Gamma^{(0)}$ and $-M(f_1,f_2,f_3)=(f_1,f_2,-f_3)$.
Therefore
\begin{eqnarray*}
\Gamma^{(0)}=\{(f_1,f_2,f_3)\mapsto (f_1,f_2,\varepsilon f_3+
2\beta_1f_1+2\beta_2f_2)\,|\, \varepsilon\in\{\pm 1\},
\beta_1,\beta_2\in \Z\}.
\end{eqnarray*}

(f) The cases $S(-l,2,-l)$ with $l\geq 3$: The following holds. 
\begin{eqnarray*}
\Rad I^{(0)}&=&\Z f_1\textup{ with }f_1=e_1-e_3,\\  
\oooo{H_\Z}^{(0)}&=&\Z\oooo{e_1}^{(0)}\oplus\Z\oooo{e_2}^{(0)},\\
&&\oooo{I}^{(0)}((\oooo{e_1}^{(0)},\oooo{e_2}^{(0)})^t,
(\oooo{e_1}^{(0)},\oooo{e_2}^{(0)}))
=\begin{pmatrix}2&-l\\-l&2\end{pmatrix},\\
\Gamma^{(0)}_u&=&\langle T(\oooo{j}^{(0)}(\oooo{e_1}^{(0)})
\otimes f_1),
T(\oooo{j}^{(0)}(l\oooo{e_2}^{(0)})\otimes f_1)\rangle\cong\Z^2
\quad \textup{with }\\
&&T(\oooo{j}^{(0)}(\oooo{e_1}^{(0)})
\otimes f_1)(\uuuu{e})=\uuuu{e}+f_1(2,-l,2),\\
&&T(\oooo{j}^{(0)}(l\oooo{e_2}^{(0)})\otimes f_1)
(\uuuu{e})=\uuuu{e}+f_1(-l^2,2l,-l^2),\\
\Gamma^{(0)}_u&\stackrel{1:l}{\subset}&
(\ker\tau^{(0)})_u\stackrel{1:(l^2-4)}{\subset}
O^{(0),Rad}_u\cong\Z^2,\\
\Gamma^{(0)}_s &\cong& \Gamma^{(0)}(S(-l))\cong G^{fCox,2},\\
\Gamma^{(0)}_s &=& (\ker \tau^{(0)})_s\cap \ker\oooo{\sigma}
\stackrel{1:4}{\subset} O^{(0),Rad}_s,\\
\Gamma^{(0)}&\stackrel{1:l}{\subset}& O^{(0),*}
\stackrel{1:4(l^2-4)}{\subset} O^{(0),Rad}.
\end{eqnarray*}
The exact sequence \eqref{6.5} splits non-canonically with
$\Gamma^{(0)}_s\cong\langle s^{(0)}_{e_1},s^{(0)}_{e_2}\rangle
\subset\Gamma^{(0)}$ (for example).

(g) The case $S(\P^2)$ and the cases $S(\uuuu{x})$ with
$\uuuu{x}\in \Z^3_{\geq 3}$, $r(\uuuu{x})<0$ and
$x_i\leq \frac{1}{2}x_jx_k$ for $\{i,j,k\}=\{1,2,3\}$:
$\Gamma^{(0)}\cong G^{fCox,3}$ is a free Coxeter group with 
the three generators $s^{(0)}_{e_1},s^{(0)}_{e_2},s^{(0)}_{e_3}$.
\end{theorem}

{\bf Proof:}
(a) This follows from the definitions in Lemma \ref{t2.6} (a)
and in Definition \ref{t2.8}. Especially
\begin{eqnarray*}
R^{(0)}&=&\{\uuuu{e}\begin{pmatrix}y_1\\y_2\\y_3\end{pmatrix}
\in H_\Z\,|\, 2=I^{(0)}(
(\uuuu{e}\begin{pmatrix}y_1\\y_2\\y_3\end{pmatrix})^t,
\uuuu{e}\begin{pmatrix}y_1\\y_2\\y_3\end{pmatrix})\\
&=&\{\uuuu{e}\begin{pmatrix}y_1\\y_2\\y_3\end{pmatrix}
\in H_\Z\,|\, 2=(y_1\ y_2\ y_3)
\begin{pmatrix}2&x_1&x_2\\x_1&2&x_3\\x_2&x_3&2\end{pmatrix}
\begin{pmatrix}y_1\\y_2\\y_3\end{pmatrix}\}.
\end{eqnarray*}

(b) First we consider the cases $S(\uuuu{x})$ with 
$\uuuu{x}\in\Z^3_{\leq 0}$ and $r(\uuuu{x})>4$.
By Lemma \ref{t5.7} (b) $\sign I^{(0)}=(++-)$.
We will apply Theorem \ref{ta.4} with $I^{[0]}:=-I^{(0)}$,
which has signature $\sign I^{[0]}=(+--)$. 

The vectors $e_1,e_2,e_3$ are negative with respect to 
$I^{[0]}$. By Theorem \ref{ta.4} (c) (vi)
the reflection $s_{e_i}^{(0)}$ acts on the model $\KK/\R^*$
of the hyperbolic plane as reflection along the hyperbolic line
$((\R e_i)^\perp\cap\KK)/\R^*\subset \KK/\R^*$.
The corresponding three planes $(\R e_1)^\perp$, $(\R e_2)^\perp$
and $(\R e_3)^\perp$ in $H_\R$ intersect pairwise in the
following three lines
\begin{eqnarray*}
(\R e_1)^\perp\cap(\R e_2)^\perp=\R y^{[1]},&&  
y^{[1]}= (-2x_2+x_1x_3,-2x_3+x_1x_2,4-x_1^2),\\
(\R e_1)^\perp\cap(\R e_3)^\perp=\R y^{[2]},&&
y^{[2]}= (-2x_1+x_2x_3,4-x_2^2,-2x_3+x_1x_2),\\
(\R e_2)^\perp\cap(\R e_3)^\perp=\R y^{[3]},&&
y^{[3]}= (4-x_3^2,-2x_1+x_2x_3,-2x_2+x_1x_3).
\end{eqnarray*}
$\uuuu{x}\in\Z^3_{\leq 0}$ and $y^{[1]}=0$ would imply
$x_2=x_3=0$, $x_1=-2$, $r(\uuuu{x})=4$. But $r(\uuuu{x})>4$
by assumption. Therefore $y^{[1]}\neq 0$.
Analogously $y^{[2]}\neq 0$ and $y^{[3]}\neq 0$. 

One calculates
\begin{eqnarray*}
I^{[0]}(y^{[i]},y^{[i]})&=& 2(4-x_i^2)(r(\uuuu{x})-4)
\quad\textup{for }i\in\{1,2,3\}\\
&&\left\{\begin{array}{ll}
\leq 0&\textup{ for }x_i\leq -2,\\
>0&\textup{ for }x_i\in\{0,-1\}.\end{array}\right.
\end{eqnarray*}
Therefore two of the three hyperbolic lines
$((\R e_j)^\perp\cap\KK)/\R^*$ $(j\in\{1,2,3\})$
intersect in $\KK/\R^*$ if and only if the corresponding
$x_i$ is $0$ or $-1$. 

\medskip
{\bf Claim:}
{\it If $x_i\in\{0,-1\}$ then the angle between the two hyperbolic
lines at the intersection point $\R^* y^{[i]}\in\KK/\R^*$
is $\frac{\pi}{2}$ if $x_i=0$ and $\frac{\pi}{3}$ if $x_i=-1$.}

\medskip
We prove the claim in an indirect way. Observe in general
\begin{eqnarray}\nonumber
x_1=0\Rightarrow 
\bigl(s_{e_1}^{(0),mat}s_{e_2}^{(0),mat}\bigr)^2
&=&(\begin{pmatrix}-1&0&-x_2\\0&1&0\\0&0&1\end{pmatrix}
\begin{pmatrix}1&0&0\\0&-1&-x_3\\0&0&1\end{pmatrix})^2\\
&=&\begin{pmatrix}-1&0&-x_2\\0&-1&-x_3\\0&0&1\end{pmatrix}^2
=E_3,\label{6.6}\\ \nonumber
x_1=-1\Rightarrow 
\bigl(s_{e_1}^{(0),mat}s_{e_2}^{(0),mat}\bigr)^3
&=&(\begin{pmatrix}-1&1&-x_2\\0&1&0\\0&0&1\end{pmatrix}
\begin{pmatrix}1&0&0\\0&-1&-x_3\\0&0&1\end{pmatrix})^3\\
&=&\begin{pmatrix}-1&-1&-x_2-x_3\\0&-1&-x_3\\0&0&1\end{pmatrix}^3
=E_3,\label{6.7}
\end{eqnarray}
and analogously for $x_2$ and $x_3$.
Therefore the angle between the hyperbolic lines must be
$\frac{\pi}{2}$ if $x_i=0$ and $\frac{\pi}{3}$ or 
$\frac{2\pi}{3}$ if $x_i=-1$.
But in the case $x_1=x_2=x_3=-1$ the three intersection points
are the vertices of a hyperbolic triangle, so then the angles
are all $\frac{\pi}{3}$. Deforming $x_2$ and $x_3$ does not
change the angle at $\R^*y^{[1]}$, so it is $\frac{\pi}{3}$
if $x_1=-1$. This proves the Claim. \hfill$(\Box)$

\medskip
Now Theorem \ref{ta.2} (a) shows that in the group of 
automorphisms of $\KK/\R^*$ which is induced by $\Gamma^{(0)}$
all relations are generated by the relations in 
\eqref{6.1}--\eqref{6.4}. Therefore this holds also for 
$\Gamma^{(0)}$ itself.

Now we consider the three reducible cases
$S(A_1^3)$, $S(A_2A_1)$ and $S(\P^1A_1).$
Lemma \ref{t2.11} gives the first isomorphisms in part (b)
for $\Gamma^{(0)}(A_1^3)$, $\Gamma^{(0)}(A_2A_1)$ and 
$\Gamma^{(0)}(\P^1A_1)$. 
Lemma \ref{t2.12} (for $A_1$) and Theorem \ref{t6.8} (b)
and (c) give the second isomorphisms in part (b).
The isomorphisms show in these three cases that all 
relations in $\Gamma^{(0)}$ are generated by the relations
in \eqref{6.1}--\eqref{6.4}.

(c) It is classical that in the case of the $A_3$ root lattice
the monodromy group $\Gamma^{(0)}$ is the Weyl group and is
$\ker\tau^{(0)}\cong S_4$.

(d) The proof of Theorem \ref{t5.14} (b) (iii) shows
\begin{eqnarray*}
\Rad I^{(0)}&=&\ker\Phi_2(M)=\ker\Phi_1(M^{root})=\Z f_1,\\
\oooo{H_\Z}^{(0)}&=&\Z\oooo{e_1}^{(0)}\oplus\Z\oooo{e_2}^{(0)},\\
&& \oooo{I}^{(0)}((\oooo{e_1}^{(0)},\oooo{e_2}^{(0)})^t,
(\oooo{e_1}^{(0)},\oooo{e_2}^{(0)}))
=\begin{pmatrix}2&-1\\-1&2\end{pmatrix},
\end{eqnarray*}
so $(\oooo{H_\Z}^{(0)},\oooo{I}^{(0)})$ is an $A_2$ root lattice.
This was treated in Theorem \ref{t6.8} (b). We have
\begin{eqnarray*}
O^{(0),Rad}_s&=&\Aut(\oooo{H_\Z}^{(0)},
\oooo{I}^{(0)})\cong D_{12},\\
\oooo{e_3}^{(0)}&=&-\oooo{e_1}^{(0)}-\oooo{e_2}^{(0)},\\
R^{(0)}(\oooo{H_\Z}^{(0)},\oooo{I}^{(0)})&=&
\{\pm\oooo{e_1}^{(0)},\pm\oooo{e_2}^{(0)},\pm\oooo{e_3}^{(0)}\},\\
\oooo{\Gamma^{(0)}}=\Gamma^{(0)}_s
&=&\langle \oooo{s^{(0)}_{e_1}},\oooo{s^{(0)}_{e_2}}\rangle
\cong \Gamma^{(0)}(A_2)\cong D_6\cong S_3,\\
\Gamma^{(0)}_s&=& (\ker \tau^{(0)})_s\stackrel{1:2}{\subset} O^{(0),Rad}_s.
\end{eqnarray*}
Observe also
\begin{eqnarray*}
&&s^{(0)}_{e_2}s^{(0)}_{e_1}s^{(0)}_{e_2}s^{(0)}_{e_3}(\uuuu{e})\\
&=&\uuuu{e}
\begin{pmatrix}1&0&0\\1&-1&1\\0&0&1\end{pmatrix}
\begin{pmatrix}-1&1&1\\0&1&0\\0&0&1\end{pmatrix}
\begin{pmatrix}1&0&0\\1&-1&1\\0&0&1\end{pmatrix}
\begin{pmatrix}1&0&0\\0&1&0\\1&1&-1\end{pmatrix}\\
&=&\uuuu{e}+f_1(1,1,-2) = T(\oooo{j}^{(0)}(-\oooo{e_3}^{(0)})
\otimes f_1)(\uuuu{e}),
\end{eqnarray*}
so $T(j^{(0)}(-\oooo{e_3}^{(0)})\otimes f_1)\in \Gamma^{(0)}_u$.
Compare Lemma \ref{t6.2} (f) and recall that $\Gamma^{(0)}_s$
acts transitively on $\{\pm \oooo{e_1}^{(0)},
\pm\oooo{e_2}^{(0)},\pm\oooo{e_3}^{(0)}\}$. 
Therefore $T(\oooo{j}^{(0)}(\oooo{e_j}^{(0)})\otimes f_1)
\in\Gamma^{(0)}_u$ for $j\in\{1,2,3\}$
with
\begin{eqnarray*}
T(\oooo{j}^{(0)}(\oooo{e_1}^{(0)})\otimes f_1)(\uuuu{e})
&=&\uuuu{e}+f_1(2,-1,-1),\\
T(\oooo{j}^{(0)}(\oooo{e_2}^{(0)})\otimes f_1)(\uuuu{e})
&=&\uuuu{e}+f_1(-1,2,-1),
\end{eqnarray*}
so
\begin{eqnarray*}
\Gamma^{(0)}_u&=&T(\oooo{j}^{(0)}(\oooo{H_\Z}^{(0)})\otimes f_1)
=(\ker \tau^{(0)})_u \\
&\subset& T(\oooo{H_\Z}^{(0),\sharp}\otimes f_1)
=O^{(0),Rad}_u\cong\Z^2.
\end{eqnarray*}
$T(\oooo{H_\Z}^{(0),\sharp}\otimes f_1)$ is generated by
\begin{eqnarray*}
(\uuuu{e}\mapsto \uuuu{e}+f_1(1,-1,0))\quad\textup{and}\quad 
(\uuuu{e}\mapsto \uuuu{e}+f_1(0,1,-1)).
\end{eqnarray*}
Therefore $\Gamma^{(0)}_u\stackrel{1:3}{\subset}
O^{(0),Rad}_u$. 

Together $\Gamma^{(0)}_s=(\ker\tau^{(0)})_s
\stackrel{1:2}{\subset}\OO_s^{(0),Rad}$ and
$\Gamma^{(0)}_u=(\ker\tau^{(0)})_u
\stackrel{1:3}{\subset}\OO_u^{(0),Rad}$ show
\begin{eqnarray*}
\Gamma^{(0)}=\ker\tau^{(0)}\stackrel{1:6}{\subset}
O^{(0),Rad}
\end{eqnarray*}
\eqref{6.7} and $x_1=-1$ show
$\langle s^{(0)}_{e_1},s^{(0)}_{e_2}\rangle\cong D_6$, 
so $\Gamma^{(0)}_s\cong \langle s^{(0)}_{e_1},
s^{(0)}_{e_2}\rangle\subset \Gamma^{(0)}$, so the exact sequence
\eqref{6.5} splits non-canonically. 

(e) Recall from the proof of Theorem \ref{t5.14} (a) (i) that
\begin{eqnarray*}
\uuuu{f}=(f_1,f_2,f_3):=\uuuu{e}
\begin{pmatrix}1&0&1\\1&1&1\\0&1&1\end{pmatrix}
\end{eqnarray*}
is a $\Z$-basis of $H_\Z$ and 
\begin{eqnarray*}
\Rad I^{(0)}&=& \Z f_1\oplus \Z f_2,\\
H_\Z&=& \Z f_3\oplus \Rad I^{(0)},\\
\oooo{H_\Z}^{(0)}&=& \Z \oooo{f_3}^{(0)}.
\end{eqnarray*}
Also observe
\begin{eqnarray*}
s^{(0)}_{e_j}|_{\Rad I^{(0)}}&=&\id\quad\textup{ for }
i\in\{1,2,3\},\\
s^{(0)}_{e_1}(f_3)&=&-f_3+2f_2,\\ 
s^{(0)}_{e_2}(f_3)&=&-f_3+2f_1+2f_2,\\ 
s^{(0)}_{e_3}(f_3)&=&-f_3+2f_1,\\
\oooo{\Gamma^{(0)}}=\Gamma^{(0)}_s &=&\{\pm\id\}
=(\ker \tau^{(0)})_s=O^{(0),Rad}_s
\cong \Gamma^{(0)}(A_1)\cong\{\pm 1\}.
\end{eqnarray*}
Therefore
\begin{eqnarray*}
\Gamma^{(0)}_u\owns s^{(0)}_{e_1} s^{(0)}_{e_2}=
(\uuuu{f}\mapsto \uuuu{f}+2f_1(0,0,1))
=T(\oooo{j}^{(0)}(\oooo{f_3}^{(0)})\otimes f_1),\\
\Gamma^{(0)}_u\owns s^{(0)}_{e_3} s^{(0)}_{e_2}=
(\uuuu{f}\mapsto \uuuu{f}+2f_2(0,0,1))
=T(\oooo{j}^{(0)}(\oooo{f_3}^{(0)})\otimes f_2),
\end{eqnarray*}
so
\begin{eqnarray*}
\Gamma^{(0)}_u&=&(\ker \tau^{(0)})_u 
=T(\oooo{j}^{(0)}(\oooo{H_\Z}^{(0)})\otimes \Rad I^{(0)})\\
&\stackrel{1:4}{\subset}& O^{(0),Rad} 
=T(\oooo{H_\Z}^{(0),\sharp}\otimes \Rad I^{(0)})\\
&=&\langle (\uuuu{f}\mapsto \uuuu{f}+f_1(0,0,1)),
(\uuuu{f}\mapsto \uuuu{f}+f_2(0,0,1))\rangle\cong\Z^2.
\end{eqnarray*}
Together the statements on $\Gamma^{(0)}_s$ and $\Gamma^{(0)}_u$
imply 
\begin{eqnarray*}
\Gamma^{(0)}=\ker\tau^{(0)} \stackrel{1:4}{\subset}
O^{(0),Rad}.
\end{eqnarray*}
The exact sequence \eqref{6.5} splits non-canonically with
$\Gamma^{(0)}_s=\{\pm id\}\cong \langle -M\rangle
\subset\Gamma^{(0)}$ (for example).

(f) Recall from the proof of Theorem \ref{t5.14} (a) (ii) 
and (b) (iv) 
\begin{eqnarray*}
\Rad I^{(0)}&=& \Z f_1\quad\textup{with}\quad 
f_1=e_1-e_3,\quad\textup{so}\\
\oooo{H_\Z}^{(0)}&=& \Z\oooo{e_1}^{(0)}\oplus\Z 
\oooo{e_2}^{(0)},\\
&&\oooo{I}^{(0)}((\oooo{e_1}^{(0)},\oooo{e_2}^{(0)})^t,
(\oooo{e_1}^{(0)},\oooo{e_2}^{(0)}))
=\begin{pmatrix}2&-l\\-l&2\end{pmatrix}.
\end{eqnarray*}

Observe
\begin{eqnarray*}
s^{(0),mat}_{e_1}
=\begin{pmatrix}-1&l&-2\\0&1&0\\0&0&1\end{pmatrix},\ 
s^{(0),mat}_{e_2}
=\begin{pmatrix}1&0&0\\l&-1&l\\0&0&1\end{pmatrix},\\
s^{(0),mat}_{e_3}
=\begin{pmatrix}1&0&0\\0&1&0\\-2&l&-1\end{pmatrix},\\
\oooo{s_{e_1}^{(0)}}(\oooo{e_1}^{(0)},\oooo{e_2}^{(0)})
=\oooo{s_{e_3}^{(0)}}(\oooo{e_1}^{(0)},\oooo{e_2}^{(0)})
=(\oooo{e_1}^{(0)},\oooo{e_2}^{(0)})
\begin{pmatrix}-1&l\\0&1\end{pmatrix},\\
\oooo{s_{e_2}^{(0)}}(\oooo{e_1}^{(0)},\oooo{e_2}^{(0)})
=\begin{pmatrix}1&0\\l&-1\end{pmatrix}.
\end{eqnarray*}
Theorem \ref{t6.8} (d) shows
\begin{eqnarray*}
\oooo{\Gamma^{(0)}}=\Gamma^{(0)}_s
\cong \Gamma^{(0)}(S(-l))\cong G^{fCox,2}.
\end{eqnarray*}
Therefore with respect to the generators 
$\oooo{s_{e_1}^{(0)}},\oooo{s_{e_2}^{(0)}}$ and
$\oooo{s_{e_3}^{(0)}}$, all relations in 
$\Gamma^{(0)}_s$ are generated by the relations
\begin{eqnarray*}
(\oooo{s_{e_1}^{(0)}})^2=(\oooo{s_{e_2}^{(0)}})^2
=\oooo{s_{e_1}^{(0)}}\ \oooo{s_{e_3}^{(0)}}=\id.
\end{eqnarray*}
Therefore $\Gamma^{(0)}_u$ is generated by the set
$\{g s^{(0)}_{e_1}s^{(0)}_{e_3}g^{-1}\,|\, g\in\Gamma^{(0)}\}$
of conjugates of $s^{(0)}_{e_1}s^{(0)}_{e_3}$. Observe
\begin{eqnarray*}
s^{(0)}_{e_1}s^{(0)}_{e_3}(\uuuu{e})
&=&\uuuu{e}\begin{pmatrix}-1&l&-2\\0&1&0\\0&0&1\end{pmatrix}
\begin{pmatrix}1&0&0\\0&1&0\\-2&l&-1\end{pmatrix}\\
&=&\uuuu{e}\begin{pmatrix}3&-l&2\\0&1&0\\-2&l&-1\end{pmatrix}
=\uuuu{e}+f_1(2,-l,2),\\
\textup{so }s^{(0)}_{e_1}s^{(0)}_{e_3}
&=& T(\oooo{j}^{(0)}(\oooo{e_1}^{(0)})\otimes f_1),
\end{eqnarray*}
and recall Lemma \ref{t6.2} (f). 
The $\Z$-lattice generated by the $\Gamma^{(0)}$ orbit of $e_1$
is $\Z e_1\oplus \Z le_2\oplus \Z \gcd(2,l)e_3$. 
Therefore
\begin{eqnarray*}
\Gamma^{(0)}_u = \langle 
T(\oooo{j}^{(0)}(\oooo{e_1}^{(0)})\otimes f_1),
T(\oooo{j}^{(0)}(l\oooo{e_2}^{(0)})\otimes f_1)\rangle\cong\Z^2
\textup{ with}\\
T(\oooo{j}^{(0)}(l\oooo{e_2}^{(0)}))\otimes f_1)(\uuuu{e})
=\uuuu{e}+f_1(-l^2,2l,-l^2).
\end{eqnarray*}
Compare
\begin{eqnarray*}
(\ker \tau^{(0)})_u &=& 
T(\oooo{j}^{(0)}(\oooo{H_\Z}^{(0)})\otimes f_1)\\
&=&\langle 
T(\oooo{j}^{(0)}(\oooo{e_1}^{(0)})\otimes f_1),
T(\oooo{j}^{(0)}(\oooo{e_2}^{(0)})\otimes f_1)\rangle,\\
O^{(0),Rad} &=& 
T(\oooo{H_\Z}^{(0),\sharp}\otimes f_1)\\
&=&\langle (\uuuu{e}\mapsto \uuuu{e}+f_1(1,0,1)),
(\uuuu{e}\mapsto\uuuu{e}+f_1(0,1,0))\rangle.
\end{eqnarray*}
Therefore 
\begin{eqnarray*}
\Gamma^{(0)}_u\stackrel{1:l}{\subset} (\ker\tau^{(0)})_u
\stackrel{1:(l^2-4)}{\subset} O^{(0),Rad}_u\cong\Z^2.
\end{eqnarray*}
Theorem \ref{t6.8} (e) shows also
\begin{eqnarray*}
\Gamma^{(0)}_s =(\ker\tau^{(0)})_s\cap \ker\oooo{\sigma}
\stackrel{1:4}{\subset} O^{(0),Rad}_s.
\end{eqnarray*}
Therefore
\begin{eqnarray*}
\Gamma^{(0)}\stackrel{1:l}{\subset} O^{(0),*}
\stackrel{1:4(l^2-4)}{\subset} O^{(0),Rad}.
\end{eqnarray*}

(g) By Lemma \ref{t5.7} (b) $\sign I^{(0)}=(+--)$. 
We will apply Theorem \ref{ta.4} with 
$I^{[0]}=I^{(0)}$ and Theorem \ref{ta.2} (b).
The vectors $e_1,e_2$ and $e_3$ are positive.
By Theorem \ref{ta.4} (c) (vii) the reflection 
$s^{(0)}_{e_i}$ acts on the model $\KK/\R^*$
of the hyperbolic plane as an elliptic element of order $2$
with fixed point $\R^*e_i\in\KK/\R^*$.

Consider the three vectors $v_1,v_2,v_3\in H_\Z\subset H_\R$
\begin{eqnarray*}
v_1&:=& -x_3e_1+x_2e_2+x_1e_3,\\
v_2&:=&  x_3e_1-x_2e_2+x_1e_3,\\
v_3&:=&  x_3e_1+x_2e_2-x_1e_3,
\end{eqnarray*}
and observe
\begin{eqnarray*}
v_1+v_2= 2x_1e_3,\ v_1+v_3=2x_2e_2,\ v_2+v_3=2x_3e_1,\\
I^{(0)}(v_i,v_i)=r(\uuuu{x})\leq 0.
\end{eqnarray*}
The three planes $\R v_1\oplus \R v_2$, $\R v_1\oplus \R v_3$
and $\R v_2\oplus \R v_3$ contain the lines 
$\R e_3$, $\R e_2$ respectively $\R e_1$. Two of the three
planes intersect in one of the lines
$\R v_1$, $\R v_2$ and $\R v_3$,
and these three lines do not meet $\KK$.
Therefore the three hyperbolic lines
$((\R v_1\oplus\R v_2)\cap\KK)/\R^*$,
$((\R v_1\oplus\R v_3)\cap\KK)/\R^*$ and
$((\R v_2\oplus\R v_3)\cap\KK)/\R^*$ 
in $\KK/\R^*$ contain the points
$\R^* e_3$, $\R^* e_2$ respectively $\R^* e_1$ and do not meet.

Now Theorem \ref{ta.2} (b) shows that the group of 
automorphisms of $\KK/\R^*$ which is induced by $\Gamma^{(0)}$
is isomorphic to $G^{fCox,3}$. 
Therefore also $\Gamma^{(0)}$ itself is isomorphic to
$G^{fCox,3}$.\hfill$\Box$

\begin{remarks}\label{t6.12}
(i) In part (g) of Theorem \ref{t6.12} we have less informations
than in the other cases. We do not even know in which cases
in part (g) $\Gamma^{(0)}=\OO^{(0),*}$ respectively
$\Gamma^{(0)}\subsetneqq\OO^{(0),*}$ holds. 

(ii) In the case of $S(\P^2)$, the proof of Theorem \ref{t6.11}
(g) gave three hyperbolic lines in the model $\KK/\R^*$ which
form a degenerate hyperbolic triangle, so with vertices on the
euclidean boundary of the hyperbolic plane.
These vertices are the lines $\R^* v_1$, $\R^* v_2$, $\R^* v_3$,
which are isotropic in the case of $S(\P^2)$ because there
$I^{(0)}(v_i,v_i)=r(\uuuu{x})=0$. 
The reflections $s^{(0)}_{e_1}$, $s^{(0)}_{e_2}$, $s^{(0)}_{e_3}$
act as elliptic elements of order 2 with fixed points on these
hyperbolic lines. Therefore the degenerate hyperbolic triangle
is a fundamental domain of this action of $\Gamma^{(0)}$.

(iii) Milanov \cite[4.1]{Mi19} had a different point of view on
$\Gamma^{(0)}$ in the case of $S(\P^2)$.
He gave an isomorphism $\Gamma^{(0)}\cong U$ to a certain 
subgroup $U$ of index 3 in $PSL_2(\Z)$. First we describe $U$ in
(iv), then we present our way to see this isomorphism in (v).

(iv) The class in $PSL_2(\Z)$ of a matrix $A\in SL_2(\Z)$
is denoted by $[A]$. It is well known that there is an
isomorphism of the free product of $\Z/2\Z$ and $\Z/3\Z$
with $PSL_2(\Z)$,
\begin{eqnarray*}
\langle \alpha\,|\, \alpha^2=e\rangle * 
\langle \beta\,|\, \beta^3=e\rangle
\to PSL_2(\Z), \\
\alpha\mapsto
[\begin{pmatrix}0&-1\\1&0\end{pmatrix}],\ 
\beta\mapsto [\begin{pmatrix}-1&-1\\1&0\end{pmatrix}].
\end{eqnarray*}
Consider the character
\begin{eqnarray*}
\chi:\langle \alpha\,|\, \alpha^2=e\rangle * 
\langle \beta\,|\, \beta^3=e\rangle
\to \{1,e^{2\pi i/3},e^{2\pi i 2/3}\}, \ \alpha\mapsto 1,\ 
\beta\mapsto e^{2\pi i /3},
\end{eqnarray*}
and the corresponding character $\www{\chi}$ on $PSL_2(\Z)$.
Then
\begin{eqnarray*}
U&=&\ker\www{\chi}\stackrel{1:3}{\subset} PSL_2(\Z),\\
\ker\chi&=& \langle \alpha,\beta\alpha\beta^2,
\beta^2\alpha\beta\rangle\\
\textup{with}&&
\alpha\simeq [F_1],\beta\alpha\beta^2\simeq[F_2],
\beta^2\alpha\beta\simeq [F_3]\textup{ and}\\
F_1&=&\begin{pmatrix}0&-1\\1&0\end{pmatrix},\quad
F_2=\begin{pmatrix}-1&-2\\1&1\end{pmatrix},\quad
F_3=\begin{pmatrix}-1&-1\\2&1\end{pmatrix}.
\end{eqnarray*}
It is easy to see $\ker\chi\cong G^{fCox,3}$, with the three
generators $\alpha$, $\beta\alpha\beta^2$, $\beta^2\alpha\beta$.
It is well known and easy to see that
\begin{eqnarray*}
\langle F_1,F_2,F_3\rangle =\{\begin{pmatrix}a&b\\c&d\end{pmatrix}
\in SL_2(\Z)\,|\, a\equiv 
d\equiv 0\mmod 3\textup{ or }b\equiv c\mmod 3\}.
\end{eqnarray*}

The M\"obius transformations $\mu(F_1)$, $\mu(F_2)$, $\mu(F_3)$
are elliptic of order 2 with fixed points
$z_1=i$, $z_2=-1+i$, $z_3=-\frac{1}{2}+\frac{1}{2}i$.

The hyperbolic lines $\www{l_1}:=A(\infty,0)$
$\www{l_2}:=A(-1,\infty)$, $\www{l_3}:=A(0,-1)$
(notations from the Remarks and Notations \ref{ta.1} (iii)) 
form a degenerate hyperbolic triangle, and $\www{l_i}$ contains $z_i$.

(v) Consider the matrices
\begin{eqnarray*}
B&:=&\begin{pmatrix}z_1\oooo{z_1} & -z_2\oooo{z_2} & 
2z_3\oooo{z_3}\\
\Ree(z_1) & -\Ree(z_2) & 2\Ree(z_3) \\ 1&-1&2\end{pmatrix}
=\begin{pmatrix} 1&-2&1\\0&1&-1\\1&-1&2\end{pmatrix}
\quad\textup{and}\\
B^{-1}&=&\frac{1}{2}\begin{pmatrix}1&3&1\\-1&1&1\\-1&-1&1
\end{pmatrix}.
\end{eqnarray*}
One checks
\begin{eqnarray*}
B^t\begin{pmatrix}0&0&1\\0&-2&0\\1&0&0\end{pmatrix}B
&=&\begin{pmatrix}2&-3&3\\-3&2&-3\\3&-3&2\end{pmatrix}
=S(\P^2)+S(\P^2)^t,\\
B^{-1}\Theta(F_i)B&=& -s_{e_i}^{(0),mat}\quad
\textup{for }i\in\{1,2,3\}
\end{eqnarray*}
(see Theorem \ref{ta.4} (i) for $\Theta$). 
The $\Z$-basis $\uuuu{e}$ of $H_\Z$ and the $\R$-basis
$\uuuu{f}$ in Theorem \ref{ta.4} of $H_\R$
are related by $\uuuu{e}=\uuuu{f}\cdot B$. 
The tuple $(v_1,v_2,v_3)$ in the proof of 
Theorem \ref{t6.11} (g) is 
\begin{eqnarray*}
(v_1,v_2,v_3)&=&\uuuu{e}
\begin{pmatrix}3&-3&-3\\3&-3&3\\-3&-3&3\end{pmatrix}\\
&=&\uuuu{f}\cdot B\cdot 
\begin{pmatrix}3&-3&-3\\3&-3&3\\-3&-3&3\end{pmatrix}
=\uuuu{f}\cdot (-6)
\begin{pmatrix}1&0&1\\-1&0&0\\1&1&0\end{pmatrix}.
\end{eqnarray*}
Finally observe that $\vartheta:\H\to \KK/\R^*$ in
Theorem \ref{ta.4} extends to the euclidean boundary with
\begin{eqnarray*}
(\vartheta(-1),\vartheta(0),\vartheta(\infty))
=\R^*\cdot \uuuu{f}\cdot
\begin{pmatrix}1&0&1\\-1&0&0\\1&1&0\end{pmatrix}.
\end{eqnarray*}
So the points $-1,0,\infty$ are mapped to the points
$\R^* v_1$, $\R^* v_2$, $\R^* v_3$.
The groups $\Gamma^{(0),mat}=
\langle s^{(0),mat}_{e_i}\,|\, i\in\{1,2,3\}\rangle$ and
$\langle -s^{(0),mat}_{e_i}\,|\, i\in\{1,2,3\}\rangle$
are isomorphic because $\Gamma^{(0),mat}$ does not contain
$-E_3$ because else it would not be a free Coxeter group with
three generators. 

Therefore the group
$U=\langle [F_1],[F_2],[F_3]\rangle\subset PSL_2(\Z)$ 
is isomorphic to the group
$\langle B^{-1}\Theta(F_i)B\,|\, i\in\{1,2,3\}\rangle
=\langle -s_{e_i}^{(0),mat}\,|\, i\in\{1,2,3\}\rangle$
and to the groups $\Gamma^{(0),mat}$ and $\Gamma^{(0)}$. 
\end{remarks}

Now we turn to the study of the set $R^{(0)}$ of roots and the
subset $\Delta^{(0)}\subset R^{(0)}$ of vanishing cycles.
For the set $R^{(0)}$ Theorem \ref{t6.11} (a) gave the 
general formula
$R^{(0)}=\{y_1e_1+y_2e_2+y_3e_3\in H_\Z\,|\, 
1=Q_3(y_1,y_2,y_3)\}$ with the quadratic form
\begin{eqnarray*}
Q_3:\Z^3\to\Z,\quad (y_1,y_2,y_3)\mapsto 
y_1^2+y_2^2+y_3^2+x_1y_1y_2+x_2y_1y_3+x_3y_2y_3.
\end{eqnarray*}
It gave also a good control on $\Gamma^{(0)}$ for all cases
$S(\uuuu{x})$ with $\uuuu{x}\in\Z^3$.

With respect to $\Delta^{(0)}$ and $R^{(0)}$ we know less.
We have a good control on them for the cases with
$r(\uuuu{x})\in\{0,1,2,4\}$ and the reducible cases,
but not for all other cases.
Theorem \ref{t6.14} treats all cases except those
in the Remarks \ref{t6.13} (ii).

\begin{remarks}\label{t6.13}
(i) The cases $S(\HH_{1,2})$, $S(-l,2,-l)$ for $l\geq 3$
and the four cases $S(3,3,4)$, $S(4,4,4)$, $S(5,5,5)$, $S(4,4,8)$
(more precisely their $\Br_3\ltimes\{\pm 1\}^3$ orbits)
are the only cases in rank 3 where we know 
$\Delta^{(0)}\subsetneqq R^{(0)}$.

(ii) We do not know whether $\Delta^{(0)}=R^{(0)}$ or
$\Delta^{(0)}\subsetneqq R^{(0)}$ in the following cases:
\begin{list}{}{}
\item[(a)]
All cases $S(\uuuu{x})$ with $r(\uuuu{x})<0$ except the four
cases $S(3,3,4)$, $S(4,4,4)$, $S(5,5,5)$, $S(4,4,8)$.
With the action of $\Br_3\ltimes\{\pm 1\}^3$ and Theorem \ref{t4.6} (c) 
they can be reduced to the cases $S(\uuuu{x})$ with
$\uuuu{x}\in\Z^3_{\geq 3}$, $r(\uuuu{x})<0$, 
$x_i\leq\frac{1}{2}x_jx_k$ for $\{i,j,k\}=\{1,2,3\}$. 
\item[(b)]
The irreducible cases $S(\uuuu{x})$ with 
$\uuuu{x}\in \Z^3_{\leq 0}$, $r(\uuuu{x})>4$ and
$\uuuu{x}\notin \{0,-1,-2\}^3$.
\end{list}
\end{remarks}

\begin{theorem}\label{t6.14}
(a) Consider the reducible cases (these include $S(A_1^3)$,
$S(A_2A_1)$, $S(\P^1A_1)$). More precisely, suppose that
$\uuuu{x}=(x_1,0,0)$. Then the tuple $(H_\Z,L,\uuuu{e})$
splits into the two tuples
$(\Z e_1\oplus \Z e_2,L|_{\Z e_1\oplus \Z e_2},(e_1,e_2))$
and $(\Z e_3,L|_{\Z e_3},e_3)$ with sets
$\Delta^{(0)}_1=R^{(0)}_1\subset \Z e_1\oplus \Z e_2$
and $\Delta^{(0)}_2=R^{(0)}_2=\{\pm e_3\}\subset\Z e_3$ of
vanishing cycles and roots, and
\begin{eqnarray*}
\Delta^{(0)}=\Delta^{(0)}_1\ \dot\cup\ \Delta^{(0)}_2
=R^{(0)}=R^{(0)}_1\ \dot\cup\  R^{(0)}_2.
\end{eqnarray*}
$\Delta^{(0)}_1=R^{(0)}_1$ is given in Theorem \ref{t6.8}.

(b) Consider $S(\uuuu{x})$ with $\uuuu{x}\in\{0,-1,-2\}^3$
and $r(\uuuu{x})>4$. Then $\Delta^{(0)}=R^{(0)}$. 

(c) The case $S(A_3)$ is classical. There
\begin{eqnarray*}
\Delta^{(0)}=R^{(0)}=\{\pm e_1,\pm e_2,\pm e_3,
\pm (e_1+e_2),\pm (e_2+e_3),\pm (e_1+e_2+e_3)\}.
\end{eqnarray*}

(d) The case $S(\whh{A}_2)$: Recall $\Rad I^{(0)}=\Z f_1$ with
$f_1=e_1+e_2+e_3$. There
\begin{eqnarray*}
\Delta^{(0)}&=&R^{(0)}=\Gamma^{(0)}\{e_1\}\\
&=& (\pm e_1+\Z f_1)\ \dot\cup\  (\pm e_2+\Z f_1)\ \dot\cup\  
(\pm e_3+\Z f_1).
\end{eqnarray*}

(e) The case $S(\HH_{1,2})$: Recall 
$(f_1,f_2,f_3)=\uuuu{e}
\begin{pmatrix}1&0&1\\1&1&1\\0&1&1\end{pmatrix}$ and 
$\Rad I^{(0)}=\Z f_1\oplus \Z f_2$. The set of roots is
\begin{eqnarray*}
R^{(0)}= \pm e_1+\Rad I^{(0)} =
(e_1+\Rad I^{(0)})\ \dot\cup\  (-e_1+\Rad I^{(0)}),
\end{eqnarray*}
with
\begin{eqnarray*}
e_1+\Rad I^{(0)} = -e_2+\Rad I^{(0)} =
e_3+\Rad I^{(0)} = f_3+\Rad I^{(0)}.
\end{eqnarray*}
It splits into the four $\Gamma^{(0)}$ orbits
\begin{eqnarray*}
\Gamma^{(0)}\{e_1\}=\pm e_1+2\Rad I^{(0)},&&
\Gamma^{(0)}\{e_2\}=\pm e_2+2\Rad I^{(0)},\\
\Gamma^{(0)}\{e_3\}=\pm e_3+2\Rad I^{(0)},&&
\Gamma^{(0)}\{f_3\}=\pm f_3+2\Rad I^{(0)}.
\end{eqnarray*}
The set $\Delta^{(0)}$ of vanishing cycles consists of the
first three of these sets,
\begin{eqnarray*}
\Delta^{(0)}=\Gamma^{(0)}\{e_1\}\ \dot\cup\ \Gamma^{(0)}\{e_2\}
\ \dot\cup\ \Gamma^{(0)}\{e_3\},
\end{eqnarray*}
so $\Delta^{(0)}\subsetneqq R^{(0)}$.

(f) The cases $S(-l,2,-l)$ with $l\geq 3$: 
Recall $\Rad I^{(0)}=\Z f_1$ with $f_1=e_1-e_3$.
As the tuple $(\oooo{H_\Z}^{(0)},\oooo{I}^{(0)},
(\oooo{e_1}^{(0)},\oooo{e_2}^{(0)}))$ is isomorphic to the
corresponding tuple from the $2\times 2$ matrix
$S(-l)=\begin{pmatrix}1&-l\\0&1\end{pmatrix}$, 
its sets of roots and its set of even vanishing cycles
coincide because of Theorem \ref{t6.8}. 
These sets are called $R^{(0)}(S(-l))$. Then
\begin{eqnarray*}
R^{(0)} = 
\{\www{y}_1e_1+\www{y}_2e_2\in H_\Z^{(0)}\,|\, 
\www{y}_1\oooo{e_1}^{(0)}+\www{y}_2\oooo{e_2}^{(0)}\in 
R^{(0)}(S(-l))\} + \Rad I^{(0)}.
\end{eqnarray*}
The cases with $l$ even: $R^{(0)}$ splits into the following
$l+2$ $\Gamma^{(0)}$ orbits,
\begin{eqnarray*}
\Gamma^{(0)}\{e_1\},\quad \Gamma^{(0)}\{e_3\},\quad 
\Gamma^{(0)}\{e_2+mf_1\}\textup{ for }m\in\{0,1,...,l-1\}.
\end{eqnarray*}
The set $\Delta^{(0)}$ of vanishing cycles consists of the first
three $\Gamma^{(0)}$ orbits,
\begin{eqnarray*}
\Delta^{(0)}=\Gamma^{(0)}\{e_1\}\ \dot\cup\ \Gamma^{(0)}\{e_3\}
\ \dot\cup\ \Gamma^{(0)}\{e_2\}.
\end{eqnarray*}
The cases with $l$ odd: $R^{(0)}$ splits into the following
$l+1$ $\Gamma^{(0)}$ orbits,
\begin{eqnarray*}
\Gamma^{(0)}\{e_1\}=\Gamma^{(0)}\{e_3\},\quad 
\Gamma^{(0)}\{e_2+mf_1\}\textup{ for }m\in\{0,1,...,l-1\}.
\end{eqnarray*}
The set $\Delta^{(0)}$ of vanishing cycles consists of the first
two $\Gamma^{(0)}$ orbits,
\begin{eqnarray*}
\Delta^{(0)}=\Gamma^{(0)}\{e_1\}\ \dot\cup\ \Gamma^{(0)}\{e_2\}.
\end{eqnarray*}
In both cases $\Delta^{(0)}\subsetneqq R^{(0)}$.

(g) The case $S(\P^2)$. Then $\Delta^{(0)}=R^{(0)}$, and
$R^{(0)}$ splits into three $\Gamma^{(0)}$ orbits,
\begin{eqnarray*}
R^{(0)}&=& \Gamma^{(0)}\{e_1\}\ \dot\cup\ \Gamma^{(0)}\{e_2\}
\ \dot\cup\ \Gamma^{(0)}\{e_3\}\quad\textup{with}\\
\Gamma^{(0)}\{e_i\}&\subset& \pm e_i+3H_\Z\quad
\textup{for }i\in\{1,2,3\}
\end{eqnarray*}
(but we would like to have a better control on $R^{(0)}$).

(h) The cases $S(3,3,4)$, $S(4,4,4)$, $S(5,5,5)$ and $S(4,4,8)$.
Then $$\Delta^{(0)}\subsetneqq R^{(0)}.$$ 
\end{theorem}

{\bf Proof:}
(a) The splittings $\Delta^{(0)}=\Delta^{(0)}_1\ \dot\cup\ 
\Delta^{(0)}_2$ and 
$R^{(0)}=R^{(0)}_1\ \dot\cup\  R^{(0)}_2$ are part of Lemma
\ref{t2.11}. Lemma \ref{t2.12} gives for $A_1$
$\Delta^{(0)}_2=R^{(0)}_2=\{\pm e_3\}$.
Theorem \ref{t6.8} gives for any rank 2 case
$\Delta^{(0)}_1=R^{(0)}_1$.

(b) The cases where $(H_\Z,L,\uuuu{e})$ is reducible are
covered by part (a). In any irreducible case, the bilinear
lattice $(H_\Z,L,\uuuu{e})$ with triangular basis is 
\index{hyperbolic bilinear lattice}
{\it hyperbolic} in the sense of the definition before 
Theorem 3.12 in \cite{HK16}, because $I^{(0)}$ is 
indefinite, but the submatrices $(2)$, 
$\begin{pmatrix}2&x_1\\x_1&2\end{pmatrix}$, 
$\begin{pmatrix}2&x_2\\x_2&2\end{pmatrix}$,
$\begin{pmatrix}2&x_3\\x_3&2\end{pmatrix}$
of the matrix $I^{(0)}(\uuuu{e}^t,\uuuu{e})$ 
are positive definite or positive semidefinite.
Theorem 3.12 in \cite{HK16} applies and gives
$\Delta^{(0)}=R^{(0)}$.

(c) This is classical. It follows also with
\begin{eqnarray*}
Q_3(y_1,y_2,y_3)&=& 2(y_1^2+y_2^2+y_3^2-y_1y_2-y_2y_3)\\
&=&y_1^2+(y_1-y_2)^2+(y_2-y_3)^2+y_3^2
\end{eqnarray*}
and the transitivity of the action of $\Gamma^{(0)}$ on $R^{(0)}$. 

(d) The quotient lattice $(\oooo{H_\Z}^{(0)},\oooo{I}^{(0)})$
is an $A_2$ lattice with set of roots 
$\{\pm \oooo{e_1}^{(0)},\pm \oooo{e_2}^{(0)},\pm \oooo{e_3}^{(0)}
\}$. Therefore 
\begin{eqnarray*}
R^{(0)}=(\pm e_1+\Z f_1)\ \dot\cup\  (\pm e_2+\Z f_1)
\ \dot\cup\  (\pm e_3+\Z f_1).
\end{eqnarray*}
(One can prove this also using
$2Q_3=(y_1-y_2)^2+(y_1-y_3)^2+(y_2-y_3)^2$.) 
$\Gamma^{(0)}_s\cong D_6$ acts transitively on the set
$\{\pm \oooo{e_1}^{(0)},\pm \oooo{e_2}^{(0)},\pm \oooo{e_3}^{(0)}
\}$. The group $\Gamma^{(0)}_u\cong \Z^2$ contains the elements
\begin{eqnarray*}
(\uuuu{e}\mapsto \uuuu{e}+f_1(2,-1,-1))\quad\textup{and}\quad
(\uuuu{e}\mapsto \uuuu{e}+f_1(-1,2,-1)).
\end{eqnarray*}
Therefore it acts transitively on each of the six sets
$\varepsilon e_i+\Z f_1$ with $\varepsilon\in\{\pm 1\}$,
$i\in\{1,2,3\}$.
Thus $\Gamma^{(0)}$ acts transitively on $R^{(0)}$, so
$\Delta^{(0)}=R^{(0)}=\Gamma^{(0)}\{e_1\}$. 

(e) The quotient lattice $(\oooo{H_\Z}^{(0)},\oooo{I}^{(0)})$
is an $A_1$ lattice with set of roots $\{\pm \oooo{e_1}^{(0)}\}$.
Therefore
\begin{eqnarray*}
R^{(0)}=\pm e_1+\Rad I^{(0)} =(e_1+\Rad I^{(0)})\ \dot\cup\ 
(-e_1+\Rad I^{(0)}).
\end{eqnarray*}
(One can prove this also using $Q_3=(y_1-y_2+y_3)^2$.)
$s^{(0)}_{e_i}$ exchanges $e_i$ and $-e_i$, and
$s^{(0)}_{e_1}$ maps $f_3$ to $-f_3+2f_2$.
$\Gamma_u^{(0)}\cong \Z^2$ is generated by the elements
\begin{eqnarray*}
((f_1,f_2,f_3)\mapsto (f_1,f_2,f_3+2f_1)) = 
(\uuuu{e}\mapsto \uuuu{e}+f_1(2,-2,2)),\\
((f_1,f_2,f_3)\mapsto (f_1,f_2,f_3+2f_2)) = 
(\uuuu{e}\mapsto \uuuu{e}+f_2(2,-2,2)).
\end{eqnarray*}
Therefore $R^{(0)}$ splits into the four $\Gamma^{(0)}$ orbits
\begin{eqnarray*}
\Gamma^{(0)}\{e_1\} &=& \pm e_1+2\Rad I^{(0)},\quad
\Gamma^{(0)}\{e_2\} = \pm e_2+2\Rad I^{(0)},\\
\Gamma^{(0)}\{e_3\} &=& \pm e_3+2\Rad I^{(0)},\quad
\Gamma^{(0)}\{f_3\} = \pm f_3+2\Rad I^{(0)}.
\end{eqnarray*}
$\Delta^{(0)}$ consists of the first three of them.

(f) The set of roots of the quotient lattice
$(\oooo{H_\Z}^{(0)},\oooo{I}^{(0)})$ is called
$R^{(0)}(S(-l))$. Theorem \ref{t6.8} describes it. Therefore
\begin{eqnarray*}
R^{(0)}=\{\www{y}_1e_1+\www{y}_2e_2\in H_\Z\,|\, 
\www{y}_1\oooo{e_1}^{(0)}+\www{y}_2\oooo{e_2}^{(0)}
\in R^{(0)}(S(-l))\} + \Rad I^{(0)}.
\end{eqnarray*}
By Theorem \ref{t6.8} (d) (iv) and (c),
$R^{(0)}(S(-l))$ splits into the two 
$\Gamma^{(0)}_s$ orbits $\Gamma^{(0)}_s\{\oooo{e_1}^{(0)}\}$
and $\Gamma^{(0)}_s\{\oooo{e_2}^{(0)}\}$, and  the action
of $\Gamma^{(0)}_s$ on each of these two orbits is simply 
transitive. 
$\Gamma^{(0)}_u\cong\Z^2$ is generated by the elements
\begin{eqnarray*}
(\uuuu{e}\mapsto \uuuu{e}+f_1(2,-l,2))\textup{ and }
(\uuuu{e}\mapsto \uuuu{e}+f_1(-l^2,2l,-l^2)).
\end{eqnarray*}
Therefore for $m\in\{0,1,...,l-1\}$
\begin{eqnarray*}
\Gamma^{(0)}\{e_2+m f_1\}\cap (e_2+\Z f_1) =e_2+mf_1+\Z lf_1.
\end{eqnarray*}
If $l$ is odd then $1=\gcd(2,l^2)$ and
\begin{eqnarray*}
\Gamma^{(0)}\{e_1\}\supset e_1+\Z f_1\owns e_3=e_1-f_1,
\quad\textup{so }\Gamma^{(0)}\{e_1\}=\Gamma^{(0)}\{e_3\}.
\end{eqnarray*}
If $l$ is even then $2=\gcd(2,l^2)$ and
\begin{eqnarray*}
\Gamma^{(0)}\{e_1\}\cap (e_1+2\Z f_1) =e_1+\Z 2f_1
\not\owns e_3,
\quad\textup{so }\Gamma^{(0)}\{e_1\}\cap\Gamma^{(0)}\{e_3\}
=\emptyset.
\end{eqnarray*}
This shows all claims.

(g) The matrices $s^{(0),mat}_{e_i}\in M_{3\times 3}(\Z)$ with
$s^{(0)}_{e_i}(\uuuu{e})=\uuuu{e}\cdot s^{(0),mat}_{e_i}$ are
\begin{eqnarray*}
s^{(0),mat}_{e_1} = 
\begin{pmatrix} -1&3&-3\\0&1&0\\0&0&1\end{pmatrix},\\
s^{(0),mat}_{e_2} = 
\begin{pmatrix} 1&0&0\\3&-1&3\\0&0&1\end{pmatrix},\quad 
s^{(0),mat}_{e_3} = 
\begin{pmatrix} 1&0&0\\0&1&0\\-3&3&-1\end{pmatrix}.
\end{eqnarray*}
One sees $\Gamma^{(0)}\{e_i\}\subset \pm e_i+3H_\Z$ for
$i\in\{1,2,3\}$. Therefore $\Delta^{(0)}$ splits into three
$\Gamma^{(0)}$ orbits,
\begin{eqnarray*}
\Delta^{(0)} = \Gamma^{(0)}\{e_1\}\ \dot\cup\ \Gamma^{(0)}\{e_2\}
\ \dot\cup\ \Gamma^{(0)}\{e_3\}.
\end{eqnarray*}
It remains to show $\Delta^{(0)}=R^{(0)}$. 
Write $\www{\uuuu{e}}=(e_1,-e_2,e_3)$, so that
\begin{eqnarray*}
L(\www{\uuuu{e}}^t,\www{\uuuu{e}})^t=\www{S} =
\begin{pmatrix}1&3&3\\0&1&3\\0&0&1\end{pmatrix}\quad\textup{and}
\quad I^{(0)}(\www{\uuuu{e}}^t,\www{\uuuu{e}})^t=\www{S}+\www{S}^t
=\begin{pmatrix}2&3&3\\3&2&3\\3&3&2\end{pmatrix}.
\end{eqnarray*}
The quadratic form $\www{Q}_3:\Z^3\to\Z$ with
\begin{eqnarray*}
\www{Q}_3(\uuuu{y})&=& \frac{1}{2} 
I^{(0)}((\www{\uuuu{e}}
\begin{pmatrix}y_1\\y_2\\y_3\end{pmatrix})^t,
\www{\uuuu{e}}\begin{pmatrix}y_1\\y_2\\y_3\end{pmatrix})
=\frac{1}{2}(y_1\ y_2\ y_3)
\begin{pmatrix}2&3&3\\3&2&3\\3&3&2\end{pmatrix}
\begin{pmatrix}y_1\\y_2\\y_3\end{pmatrix}\\
&=&y_1^2+y_2^2+y_3^2+3(y_1y_2+y_1y_3+y_2y_3)
\end{eqnarray*}
can also be written in the following two ways which will be useful
below,
\begin{eqnarray}\label{6.8}
\www{Q}_3(\uuuu{y})&=& (y_1+y_2)(y_1+y_3)+(y_1+y_2)(y_2+y_3)
+(y_1+y_3)(y_2+y_3),\\
\www{Q}_3(\uuuu{y})&=& \frac{3}{2}(y_1+y_2+y_3)^2-
\frac{1}{2}(y_1^2+y_2^2+y_3^2).\label{6.9}
\end{eqnarray}
We have $R^{(0)}=\{y_1\www{e_1}+y_2\www{e_2}+y_3\www{e_3}
\in H_\Z\,|\, \www{Q}_3(\uuuu{y})=1\}.$
Define
\begin{eqnarray*}
\|a\| := \sqrt{y_1^2+y_2^2+y_3^2}\quad\textup{for}\quad 
a=\sum_{i=1}^3 y_ie_i\in H_\Z.
\end{eqnarray*}

\medskip
{\bf Claim:} {\it For any $a\in R^{(0)}-\{\pm e_1,\pm e_2,
\pm e_3\}$ an index $i\in\{1,2,3\}$ with
\begin{eqnarray*}
\| s^{(0)}_{e_i}(a)\| <\|a\|
\end{eqnarray*}
exists.}

\medskip
The Claim implies $\Delta^{(0)}=R^{(0)}$  because it says that
any $a\in R^{(0)}$ can be mapped by a suitable sequence of 
reflections in $\{s^{(0)}_{e_1},s^{(0)}_{e_2},s^{(0)}_{e_3}\}$
to $e_1$ or $e_2$ or $e_3$. It remains to prove the Claim.

{\bf Proof of the Claim:}
Suppose $a\in R^{(0)}-\{\pm e_1,\pm e_2,\pm e_3\}$ satisfies
$\| s^{(0)}_{e_i}(a)\|\geq a$ for any $i\in\{1,2,3\}$.
Write $a=y_1\www{e_1}+y_2\www{e_2}+y_3\www{e_3}$. 
For $j$ and $k$ with $\{i,j,k\}=\{1,2,3\}$
\begin{eqnarray*}
&&\| s^{(0)}_{e_i}(a)\| \geq \|a\| \iff 
\| s^{(0)}_{e_i}(a)\|^2 \geq \|a\|^2 \\
&\iff& (-y_i-3y_j-3y_k)^2 + y_j^2+y_k^2 \geq 
y_i^2+y_j^2+y_k^2\\
&\iff& 6y_i(y_j+y_k) + 9 (y_j+y_k)^2 \geq 0 \\
&\iff& \left\{\begin{array}{ll}
3(y_j+y_k)\geq -2y_i &\textup{ if }y_j+y_k>0,\\
3(y_j+y_k)\leq -2y_i &\textup{ if }y_j+y_k<0,\\
\textup{no condition} &\textup{ if }y_j+y_k=0.
\end{array}\right.
\end{eqnarray*}
$y_j+y_k=0$ is impossible because else by formula \eqref{6.8} 
\begin{eqnarray*}
1=\www{Q}_3(\uuuu{y})=
(y_i+y_j)(y_i+y_k)=(y_i+y_j)(y_i-y_j)=y_i^2-y_j^2,\\
\quad\textup{so }y_i=\pm 1,y_j=y_k=0,
\end{eqnarray*}
which is excluded by $a\in R^{(0)}-\{\pm e_1,\pm e_2,\pm e_3\}$.
Also $(y_1+y_2>0,y_1+y_3>0,y_2+y_3>0)$ and
$(y_1+y_2<0,y_1+y_3<0,y_2+y_3<0)$ are impossible because of 
$1=\www{Q}_3(\uuuu{y})$ and \eqref{6.8}.

We can suppose
\begin{eqnarray*}
y_1+y_2>0,\quad y_1+y_3>0,\quad y_2+y_3<0,\quad
y_1\geq y_2\geq y_3.
\end{eqnarray*}
Then
\begin{eqnarray*}
y_1>0,\quad y_3\in\Z\cap [-y_1+1,-1],\quad 
y_2\in [-y_3-1,y_3],\\
3(y_1+y_2)\geq -2y_3\quad\textup{ because of }y_1+y_2>0,\\
3(y_2+y_3)\leq -2y_1\quad\textup{ because of }y_2+y_3<0,
\end{eqnarray*}
so
\begin{eqnarray*}
y_1&\geq& 3(y_1+y_2+y_3)\geq y_3\geq -y_1+1>-y_1,\\
|y_1|&\geq& 3|y_1+y_2+y_3|,\\
y_1^2&\geq& 9(y_1+y_2+y_3)^2,\\
\www{Q}_3(\uuuu{y}) &\stackrel{\eqref{6.9}}{=}&
\frac{3}{2}(y_1+y_2+y_3)^2-\frac{1}{2}(y_1^2+y_2^2+y_3^2)\\
&\leq& \frac{1}{6}y_1^2-\frac{1}{2}(y_1^2+y_2^2+y_3^2)\leq 0,
\end{eqnarray*}
a contradiction. Therefore an $a\in R^{(0)}-\{\pm e_1,\pm e_2,
\pm e_3\}$ with $\| s^{(0)}_{e_i}(a)\| \geq \| a\|$ for each
$i\in\{1,2,3\}$ does not exist. The Claim is proved.

(h) By Theorem \ref{t6.11} (g) $\Gamma^{(0)}$ is a free Coxeter
group with generators $s_{e_1}^{(0)}$, $s_{e_2}^{(0)}$ and 
$s_{e_3}^{(0)}$. By Example \ref{t3.23} (iv) equality holds
in \eqref{3.3}, so
$$\BB^{dist}=\{\uuuu{v}\in (\Delta^{(0)})^3\,|\, 
s_{v_1}^{(0)}s_{v_2}^{(0)}s_{v_3}^{(0)}=-M\}$$
(see also Theorem \ref{t7.2} (a)). By Theorem \ref{t5.16}
(a)+(b)+(d)+(e) $Q\in G_\Z-G_\Z^{\BB}$. Lemma \ref{t3.22} (a)
and $QMQ^{-1}=M$ give
$$s_{Q(e_1)}^{(0)}s_{Q(e_2)}^{(0)}s_{Q(e_3)}^{(0)}
=Qs_{e_1}^{(0)}s_{e_2}^{(0)}s_{e_3}^{(0)}Q^{-1}
=Q(-M)Q^{-1}=-M.$$
If $Q(e_1),Q(e_2),Q(e_3)$ were all in $\Delta^{(0)}$ then
equality in \eqref{3.3} would imply 
$(Q(e_1),Q(e_2),Q(e_3))\in\BB^{dist}$ and $Q\in G_\Z^{\BB}$,
a contradiction. So $Q(e_1),Q(e_2),Q(e_3)$ are not all in 
$\Delta^{(0)}$. But of course they are in $R^{(0)}$. 
\hfill $\Box$

\begin{remarks}\label{t6.15}
(i) In view of Remark \ref{t6.13} (ii) it would be desirable to 
extend the proof of $\Delta^{(0)}=R^{(0)}$ in part (g) of Theorem 
\ref{t6.11} to other cases. The useful formulas \eqref{6.8} and 
\eqref{6.9} generalize as follows:
\begin{eqnarray*}
(x_1+x_2+x_3)Q_3(\uuuu{y}) = 
(x_1+x_2+x_3)(y_1^2+y_2^2+y_3^2)\\
-x_1x_2y_1^2-x_1x_3y_2^2-x_2x_3y_3^2
+ (x_1y_2+x_2y_3)(x_1y_1+x_3y_3)\\
+ (x_1y_2+x_2y_3)(x_2y_1+x_3y_2)
+ (x_1y_1+x_3y_3)(x_2y_1+x_3y_2).
\end{eqnarray*}
If $\uuuu{x}\in(\Z-\{0\})^3$ then
\begin{eqnarray*}
2 Q_3(\uuuu{y})&=& x_1x_2x_3(\frac{y_1}{x_3}+\frac{y_2}{x_2}
+\frac{y_3}{x_1})^2 
-(x_1x_2x_3-2x_3^2)\Bigl(\frac{y_1}{x_3}\Bigr)^2 \\
&-&(x_1x_2x_3-2x_2^2)\Bigl(\frac{y_2}{x_2}\Bigr)^2 
-(x_1x_2x_3-2x_1^2)\Bigl(\frac{y_3}{x_1}\Bigr)^2.
\end{eqnarray*}
Also the rephrasing in the proof of part (g) of the inequality
$\| s^{(0)}_{e_i}(a)\| \geq \| a\|$ generalizes naturally.
But the further arguments do not seem to generalize easily.

(ii) In view of Theorem \ref{t6.14}, we know the following for $n=3$: 
\begin{eqnarray*}
\Delta^{(0)}=R^{(0)}&&\textup{in the cases }S(\uuuu{x})\textup{ with }
\uuuu{x}\in\{0,-1,-2\}^3\textup{ with }r(\uuuu{x})>4,\\
&&\textup{in all reducible cases }\\
&&\textup{and in the cases }A_3,\whh{A}_2,\P^2.\\
\Delta^{(0)}\subsetneqq R^{(0)}&&\textup{in the cases }\HH_{1,2},
\ S(-l,2,-l)\textup{ with }l\geq 3,\\
&& S(3,3,3),\ S(4,4,4),\ S(5,5,5)\textup{ and }S(4,4,8).
\end{eqnarray*}
In the following cases with $n=3$
we do not know whether $\Delta^{(0)}=R^{(0)}$ or 
$\Delta^{(0)}\subsetneqq R^{(0)}$ holds: All cases $\uuuu{x}\in\Z^3$
with $r(\uuuu{x})<0$ except four cases 
and many cases $\uuuu{x}\in\Z^3$ with $r(\uuuu{x})>4$. 
\end{remarks}

\section{The odd rank 3 cases}\label{s6.4}

For $\uuuu{x}\in\Z^3$ consider the matrix 
$S=S(\uuuu{x})
=\begin{pmatrix}1&x_1&x_2\\0&1&x_3\\0&0&1\end{pmatrix}
\in T^{uni}_3(\Z)$, and consider a unimodular bilinear lattice
$(H_\Z,L)$ with a triangular basis $\uuuu{e}=(e_1,e_2,e_3)$
with $L(\uuuu{e}^t,\uuuu{e})^t=S$.

In this section we will determine in all cases the odd monodromy
group $\Gamma^{(1)}=\langle s_{e_1}^{(1)},s_{e_2}^{(1)},
s_{e_3}^{(1)}\rangle\subset O^{(1)}$ and in many, but not all
cases the set 
$\Delta^{(1)}=\Gamma^{(1)}\{\pm e_1,\pm e_2,\pm e_3\}$ of
odd vanishing cycles. 

Recall Remark \ref{t4.17}. The group $\Gamma^{(1)}$ and the set
$\Delta^{(1)}$ are determined by the triple 
$(H_\Z,I^{(1)},\uuuu{e})$, and here $I^{(1)}$ is needed only
up to the sign.
By Remark \ref{t4.17} and Lemma \ref{t4.18} we can restrict
to $\uuuu{x}$ in the union of the following three families.
It will be useful to split each of the first two families
into the three subfamilies on the right hand side.
\begin{eqnarray*}
&&(x_1,x_2,0)\textup{ with }x_1\geq x_2\geq 0,\quad
\left\{\begin{array}{l}
(x_1,0,0):\textup{ reducible cases,}\\
(1,1,0):\quad A_3\textup{ and }\whh{A}_2,\\
(x_1,x_2,0) \textup{ with }2\leq x_1\geq x_2>0,
\end{array}\right. \\\
&&(-l,2,-l)\textup{ with }l\geq 2,\quad
\left\{\begin{array}{l}
(-l,2,-l)\textup{ with }l\equiv 0(4),\\
(-l,2,-l)\textup{ with }l\equiv 2(4) 
\textup{ (this includes }\HH_{1,2}),\\
(-l,2,-l)\textup{ with }l\equiv 1(2),
\end{array}\right. \\
&&(x_1,x_2,x_3)\in\Z^3_{\geq 3}\textup{ with }
2x_i\leq x_jx_k\textup{ for }\{i,j,k\}=\{1,2,3\}\\
&&\hspace*{7cm}(\textup{this includes }\P^2).
\end{eqnarray*}

Recall
\begin{eqnarray*}
\uuuu{\www{x}}=(\www{x}_1,\www{x}_2,\www{x}_3):=
\gcd(x_1,x_2,x_3)^{-1}(x_1,x_2,x_3)\quad\textup{for }
\uuuu{x}\neq (0,0,0).
\end{eqnarray*}
Recall from section \ref{s5.3} the definition 
\begin{eqnarray*}
f_3:=-\www{x}_3e_1+\www{x}_2e_2-\www{x}_1e_3\in H_\Z^{prim}
\quad\textup{for }\uuuu{x}\neq (0,0,0)
\end{eqnarray*}
and the fact
\begin{eqnarray*}
\Rad I^{(1)}&=& \left\{\begin{array}{ll}
\Z f_3&\textup{ if }\uuuu{x}\neq (0,0,0),\\
H_\Z &\textup{ if }\uuuu{x}=(0,0,0).
\end{array}\right.
\end{eqnarray*}
Therefore in all cases except $\uuuu{x}=(0,0,0)$ the exact
sequence
\begin{eqnarray}\label{6.10}
\{1\}\to\Gamma^{(1)}_u\to\Gamma^{(1)}\to \Gamma^{(1)}_s\to\{1\}
\end{eqnarray}
in Lemma \ref{t6.2} (d) is interesting.

\begin{lemma}\label{t6.16}
Suppose $x_1\neq 0$ (this holds in the three families above
except for the case $\uuuu{x}=(0,0,0)$).

(a) The sublattice $\Z\oooo{e_1}^{(1)} + \Z \oooo{e_2}^{(1)}
\subset \oooo{H_\Z}^{(1)}$ has index $\www{x}_1$ in 
$\oooo{H_\Z}^{(1)}$. 

(b) $\oooo{I}^{(1)}$ is nondegenerate. 
For each $\Z$-basis $\uuuu{b}=(b_1,b_2)$ of $\oooo{H_\Z}^{(1)}$
\begin{eqnarray*}
\oooo{I}^{(1)}(\uuuu{b}^t,\uuuu{b})=\varepsilon\gcd(x_1,x_2,x_3)
\begin{pmatrix}0&-1\\1&0\end{pmatrix}
\end{eqnarray*}
for some $\varepsilon\in\{\pm 1\}$. Also 
$O^{(1),Rad}_s\cong SL_2(\Z)$.

(c) 
\begin{eqnarray*}
&&\oooo{e_1}^{(1)}\in \gcd(\www{x}_1,\www{x}_2)
\oooo{H_\Z}^{(1),prim},\quad 
\oooo{e_2}^{(1)}\in \gcd(\www{x}_1,\www{x}_3)
\oooo{H_\Z}^{(1),prim},\\
&&\oooo{e_3}^{(1)}\in \gcd(\www{x}_2,\www{x}_3)
\oooo{H_\Z}^{(1),prim}\textup{ if }
(\www{x}_2,\www{x}_3)\neq (0,0),\quad\textup{ else }
\oooo{e_3}^{(1)}=0.
\end{eqnarray*}
\end{lemma}

{\bf Proof:}
(a) 
\begin{eqnarray*}
\oooo{H_\Z}^{(1)}&=& \Z\oooo{e_1}^{(1)} + \Z\oooo{e_2}^{(1)}
+ \Z\oooo{e_3}^{(1)}\\
&=& \Z\oooo{e_1}^{(1)} + \Z\oooo{e_2}^{(1)} 
+ \Z\frac{1}{\www{x}_1}(-\www{x}_3\oooo{e_1}^{(1)}
+\www{x}_2\oooo{e_2}^{(1)})\\
&=& \Z\oooo{e_1}^{(1)} + \Z\oooo{e_2}^{(1)} 
+ \Z\frac{\xi}{\www{x}_1}h_2\quad\textup{with}\\
\xi&:=& \gcd(\www{x}_2,\www{x}_3),\quad 
h_2:=-\frac{\www{x}_3}{\xi}\oooo{e_1}^{(1)}
+\frac{\www{x}_2}{\xi}\oooo{e_2}^{(1)}.
\end{eqnarray*}
The element $h_2$ is in 
$(\Z \oooo{e_1}^{(1)}+\Z \oooo{e_2}^{(1)})^{prim}$. 
One can choose a second element 
$h_1\in \Z \oooo{e_1}^{(1)}+\Z \oooo{e_2}^{(1)}$
with $\Z \oooo{e_1}^{(1)}\oplus\Z \oooo{e_2}^{(1)}
=\Z h_1\oplus \Z h_2$. Then
\begin{eqnarray*}
\Z h_2+\Z\frac{\xi}{\www{x}_1}h_2=\Z \frac{1}{\www{x}_1}h_2
\end{eqnarray*}
because $\gcd(\www{x}_1,\xi)=1$. Therefore
\begin{eqnarray*}
\oooo{H_\Z}^{(1)} = 
(\Z h_1+\Z h_2)+\Z\frac{\xi}{\www{x}_1}h_2
=\Z h_1\oplus \Z\frac{1}{\www{x}_1}h_2
\stackrel{\www{x}_1:1}{\supset} \Z h_1\oplus \Z h_2.
\end{eqnarray*}

(b) $\oooo{I}^{(1)}(\oooo{e_1}^{(1)},\oooo{e_2}^{(1)})=x_1\neq0$.
With part (a) one sees
\begin{eqnarray*}
I^{(1)}(b_1,b_2)=\pm \frac{x_1}{\www{x}_1}=\pm 
\gcd(x_1,x_2,x_3).
\end{eqnarray*}
A rank two lattice with a nondegenerate skew-symmetric bilinear 
form has an automorphism group isomorphic to $SL_2(\Z)$.

(c) The proof of part (a) and 
$1=\gcd(\www{x}_3,\gcd(\www{x}_1,\www{x}_2))$ show 
\begin{eqnarray*}
\Q\oooo{e_1}^{(1)}\cap \oooo{H_\Z}^{(1)}
=\Z\oooo{e_1}^{(1)} + \Z 
\frac{-\www{x}_3}{\gcd(\www{x}_1,\www{x}_2)}\oooo{e_1}^{(1)}
=\Z \frac{1}{\gcd(\www{x}_1,\www{x}_2)}\oooo{e_1}^{(1)}.
\end{eqnarray*}
This shows $\oooo{e_1}^{(1)}\in 
\gcd(\www{x}_1,\www{x}_2) \oooo{H_\Z}^{(1),prim}$.
Analogously $\oooo{e_2}^{(1)}\in 
\gcd(\www{x}_1,\www{x}_3) \oooo{H_\Z}^{(1),prim}$.

If $(\www{x}_2,\www{x}_3)=(0,0)$ then $-\www{x}_1\oooo{e_3}^{(1)}
=\oooo{f_3}^{(1)}=0$, so $\oooo{e_3}^{(1)}=0$.
If $\www{x}_2\neq 0$ or $\www{x}_3\neq 0$, formulas as in the
proof of part (a) hold also for 
$\Z\oooo{e_1}^{(1)}+\Z\oooo{e_3}^{(1)}$ respectively
$\Z\oooo{e_2}^{(1)}+\Z\oooo{e_3}^{(1)}$. In both cases one
shows $\oooo{e_3}^{(1)}\in\gcd(\www{x}_2,\www{x}_3)
\oooo{H_\Z}^{(1),prim}$ as above. \hfill$\Box$

\bigskip

In Theorem \ref{t6.18} we will consider in many cases the three 
groups $\Gamma^{(1)}_u\subset (\ker\tau^{(1)})_u\subset
O^{(1),Rad}_u$. Their descriptions in Lemma \ref{t6.2} (e)
simplify because now $\Rad I^{(1)}=\Z f_3$ if
$\uuuu{x}\neq (0,0,0)$. In the cases $\uuuu{x}=(-l,2,-l)$ with
$l\equiv 2(4)$ also the larger group 
$$O^{(1),Rad}_{\pm}:=\{g\in O^{(1),Rad}\,|\, \oooo{g}=\pm\id\}
\stackrel{2:1}{\supset} O^{(1),Rad}_u$$ 
will be considered.
The following Lemma  fixes notations and gives a description of
$O^{(1),Rad}_{\pm}$ similar to the one for $O^{(1),Rad}_u$ 
in Lemma \ref{t6.2} (e). It makes also $O^{(1),Rad}_u$
and $(\ker\tau^{(1)})_u$ more explicit, and, under some
condition, $\Gamma^{(1)}_u$ and 
$\Gamma^{(1)}\cap O^{(1),Rad}_\pm$.

\begin{lemma}\label{t6.17}
Suppose $x_1\neq 0$. Denote  
\index{$t_\lambda^+,\ t_\lambda^-$}
\index{$\Hom_{0\textup{ or }2}(H_\Z,\Z)$} 
\begin{eqnarray*}
\Hom_0(H_\Z,\Z)&:=& \{\lambda:H_\Z\to \Z\,|\, \lambda
\textup{ is }\Z\textup{-linear},\lambda(f_3)=0\},\\
\Hom_2(H_\Z,\Z)&:=& \{\lambda:H_\Z\to \Z\,|\, \lambda
\textup{ is }\Z\textup{-linear},\lambda(f_3)=2\},\\
t^+_\lambda:H_\Z\to H_\Z&& \textup{with }
t^+_\lambda(a)=a+\lambda(a)f_3\quad\textup{for }
\lambda\in \Hom_0(H_\Z,\Z),\\
t^-_\lambda:H_\Z\to H_\Z&& \textup{with }
t^-_\lambda(a)=-a+\lambda(a)f_3\quad\textup{for }
\lambda\in \Hom_2(H_\Z,\Z).
\end{eqnarray*}

(a) Then $t^+_\lambda\in O^{(1),Rad}_u$ for 
$\lambda\in \Hom_0(H_\Z,\Z)$, 
$t^-_\lambda\in O^{(1),Rad}_{\pm}-O^{(1),\Rad}_u$ for 
$\lambda\in \Hom_2(H_\Z,\Z)$.
The maps
\begin{eqnarray*}
\Hom_0(H_\Z,\Z)&\to& O^{(1),Rad}_u,\quad \lambda\mapsto 
t^+_\lambda,\\
\Hom_2(H_\Z,\Z)&\to& O^{(1),Rad}_\pm-O^{(1),Rad}_u,\quad 
\lambda\mapsto t^-_\lambda,
\end{eqnarray*}
are bijections, and the first one is a group isomorphism.
For $\lambda_1,\lambda_2\in \Hom_0(H_\Z,\Z)$ and
$\lambda_3,\lambda_4\in \Hom_2(H_\Z,\Z)$
\begin{eqnarray*}
t^+_{\lambda_2}\circ t^+_{\lambda_1} 
&=& t^+_{\lambda_2+\lambda_1},\quad
t^-_{\lambda_3}\circ t^+_{\lambda_1} 
= t^-_{\lambda_3+\lambda_1},\\
t^+_{\lambda_1}\circ t^-_{\lambda_3} 
&=& t^-_{-\lambda_1+\lambda_3},\quad
t^-_{\lambda_4}\circ t^-_{\lambda_3} 
= t^+_{-\lambda_4+\lambda_3},\\
&&\textup{and especially }(t^-_{\lambda_3})^2=\id.
\end{eqnarray*}

(b) 
\begin{eqnarray*}
(\ker\tau^{(1)})_u=\{t^+_\lambda\,|\,\lambda(\uuuu{e})\in 
\langle (0,x_1,x_2),(-x_1,0,x_3),(x_2,x_3,0)\rangle_\Z\}.
\end{eqnarray*}

(c) If $\Gamma^{(1)}_u$ is the normal subgroup generated by
$t^+_{\lambda_1}$ for some $\lambda_1\in\Hom_0(H_\Z,\Z)$,
then $\Gamma^{(1)}_u=\{t^+_\lambda\,|\, \lambda(\uuuu{e})
\in L\}$ where $L\subset\Z^3$ is the smallest sublattice with
$\lambda_1(\uuuu{e})\in L$ and 
$L\cdot (s_{e_i}^{(1),mat})^{\pm 1}\subset L$ for $i\in\{1,2,3\}$.

(d) If $\Gamma^{(1)}\cap O^{(1),Rad}_\pm$ is the normal subgroup 
generated by $t^-_{\lambda_1}$ for some 
$\lambda_1\in\Hom_2(H_\Z,\Z)$,
then 
\begin{eqnarray*}
\Gamma^{(1)}_u&=&\{t^+_\lambda\,|\, \lambda(\uuuu{e})
\in L\}\quad\textup{and}\\
\Gamma^{(1)}\cap O^{(1),Rad}_\pm-\Gamma^{(1)}_u
&=&\{t^-_\lambda\,|\, \lambda(\uuuu{e})\in \lambda_1(\uuuu{e})
+L\}
\end{eqnarray*}
where $L\subset\Z^3$ is the smallest sublattice with
$\lambda_1(\uuuu{e})-\lambda_1(\uuuu{e})
\cdot (s_{e_i}^{(1),mat})^{\pm 1}\in L$ 
and $L\cdot (s_{e_i}^{(1),mat})^{\pm 1}\subset L$ 
for $i\in\{1,2,3\}$. 
\end{lemma}

{\bf Proof:}
(a) By definition $t^+_\lambda =T(\oooo{\lambda}\otimes f_3)$
where $\oooo{\lambda}\in \oooo{H_\Z}^{(1),\sharp}$ denotes
the element which is induced by $\lambda$.
The map $\Hom_0(H_\Z,\Z)\to O^{(1),Rad}_u$ is an isomorphism
by Lemma \ref{t6.2} (e).
The proofs of the other statements are similar or easy.

(b) The row vectors $(0,x_1,x_2),(-x_1,0,x_3),(-x_2,-x_3,0)$ 
are the rows of the matrix $I^{(1)}(\uuuu{e}^t,\uuuu{e})$.
Because of this and Lemma \ref{t6.2} (e) $(\ker\tau^{(1)})_u$
is as claimed. 

(c) This follows from
\begin{eqnarray*}
\Gamma^{(1)}_u&=& \langle g^{-1}\circ t^+_{\lambda_1}\circ g\,|\, 
g\in \Gamma^{(1)}\rangle,\\
g^{-1}\circ t^+_{\lambda}\circ g&=&t^+_{\lambda\circ g},\\
\textup{and }\lambda\circ (s_{e_i}^{(1)})^{\pm 1}(\uuuu{e})
&=& \lambda(\uuuu{e})\cdot (s_{e_i}^{(1),mat})^{\pm 1},
\end{eqnarray*}
for $\lambda\in\Hom_0(H_\Z,\Z)$. 

(d) Similar to the proof of part (c). \hfill$\Box$

\bigskip
The following theorem gives $\Gamma^{(1)}$ for $\uuuu{x}$
in one of the three families above and thus via Remark 
\ref{t4.17} and Lemma \ref{t4.18} in principle for all
$\uuuu{x}\in\Z^3$. Though recall that it is nontrivial to find
for a given $\uuuu{x}\in\Z^3$ an element in one of the three
families above which is in the 
$(G^{phi}\ltimes \www{G}^{sign})\rtimes\langle\gamma\rangle$
orbit of $\uuuu{x}$.

\begin{theorem}\label{t6.18}
(a) We have
\begin{eqnarray*}
s_{e_i}^{(1)}\uuuu{e}&=& \uuuu{e}\cdot s_{e_i}^{(1),mat}
\quad\textup{with}\quad 
s_{e_1}^{(1),mat}
=\begin{pmatrix}1&-x_1&-x_2\\0&1&0\\0&0&1\end{pmatrix},\\
s_{e_2}^{(1),mat}
&=&\begin{pmatrix}1&0&0\\x_1&1&-x_3\\0&0&1\end{pmatrix},\quad 
s_{e_3}^{(1),mat}
=\begin{pmatrix}1&0&0\\0&1&0\\x_2&x_3&1\end{pmatrix},\\
\Gamma^{(1)}&\cong& \Gamma^{(1),mat}
=\langle s_{e_1}^{(1),mat},s_{e_2}^{(1),mat},s_{e_3}^{(1),mat}
\rangle \subset SL_3(\Z).
\end{eqnarray*}

(b) In each reducible case $\uuuu{x}=(x_1,0,0)$
\begin{eqnarray*}
\Gamma^{(1)}\cong \Gamma^{(1)}(S(-x_1))\times \Gamma^{(1)}(A_1)
\cong \Gamma^{(1)}(S(-x_1))\times\{1\},
\end{eqnarray*}
and $\Gamma^{(1)}(S(-x_1))$ is given in Theorem \ref{t6.10},
with 
\begin{eqnarray*}
\Gamma^{(1)}(S(0))&\cong&\Gamma^{(1)}(A_1^2)\cong \{1\},\\
\Gamma^{(1)}(S(-1))&\cong&\Gamma^{(1)}(A_2)\cong SL_2(\Z),\\
\Gamma^{(1)}(S(-x_1))&\cong& G^{free,2}\quad\textup{ for }
x_1\geq 2.
\end{eqnarray*}
Also $\Gamma^{(1)}\cong \Gamma^{(1)}_s$ and 
$\Gamma^{(1)}_u=\{\id\}$. 

(c) The case $\uuuu{x}=(1,1,0)$: (This is the case of $A_3$
and $\whh{A}_2$.) 
\begin{eqnarray*}
\Rad I^{(1)}&=&\Z f_3\quad \textup{with}\quad f_3=e_2-e_3,\\
\oooo{H_\Z}^{(1)}&=& \Z \oooo{e_1}^{(1)}\oplus 
\Z\oooo{e_2}^{(1)},\\
\Gamma^{(1)}_u&=& (\ker\tau^{(1)})_u = O^{(1),Rad}_u
=\{t^+_\lambda\,|\, \lambda\in\langle\lambda_1,\lambda_2\rangle_\Z
\}\cong\Z^2, \\
&&\textup{with}\quad \lambda_1(\uuuu{e})=(1,0,0),\quad
\lambda_2(\uuuu{e})=(0,1,1),\\
\Gamma^{(1)}_s&=& (\ker \tau^{(1)})_s=O^{(1),Rad}_s
\cong SL_2(\Z),\\
\Gamma^{(1)}&=& \ker\tau^{(1)}=O^{(1),Rad}.
\end{eqnarray*}
The exact sequence \eqref{6.10} splits non-canonically with
$\Gamma^{(1)}_s\cong \langle s_{e_1}^{(1)},s_{e_2}^{(1)}\rangle
\subset\Gamma^{(1)}$ (for example). 

(d) The cases $\uuuu{x}=(x_1,x_2,0)$ with $2\leq x_1\geq x_2>0$:
Write $$x_{12}:=\gcd(x_1,x_2)=\frac{x_1}{\www{x}_1}
=\frac{x_2}{\www{x}_2}.$$ 
Then
\begin{eqnarray*}
\Rad I^{(1)}&=& \Z f_3\quad\textup{with}\quad
f_3=\www{x}_2e_2-\www{x}_1e_3,\\
\oooo{H_\Z}^{(1)}&=& \Z \oooo{e_1}^{(1)}\oplus 
\Z g_2\quad\textup{with }g_2:=\frac{1}{\www{x}_1}\oooo{e_2}^{(1)}
=\frac{1}{\www{x}_2}\oooo{e_3}^{(1)}\in \oooo{H_\Z}^{(1)},
\end{eqnarray*}
\begin{eqnarray*}
\Gamma^{(1)}_u&=&\{t^+_\lambda\,|\, \lambda\in\langle\lambda_1,
\lambda_2\rangle_\Z\}\cong\Z^2 \quad\textup{with}\\
&& \lambda_1(\uuuu{e})=x_{12}\www{x}_1\www{x}_2(1,0,0),\quad
\lambda_2(\uuuu{e})=x_1x_2(0,\www{x}_1,\www{x}_2),\\
(\ker\tau^{(1)})_u&=&\{t^+_\lambda\,|\, 
\lambda(\uuuu{e})\in\langle x_{12}(1,0,0),(0,x_1,x_2)
\rangle_\Z\},\\
O^{(1),Rad}_u&=&\{t^+_\lambda\,|\, \lambda(\uuuu{e})\in\langle
(1,0,0),(0,\www{x}_1,\www{x}_2)\rangle_\Z\},
\end{eqnarray*}
\begin{eqnarray*}
\Gamma^{(1)}_s&\cong& \Gamma^{(1)}(S(-x_{12}))\cong
\langle\begin{pmatrix}1&-x_{12}\\0&1\end{pmatrix},
\begin{pmatrix}1&0\\x_{12}&1\end{pmatrix}\rangle \\
&\cong& \left\{\begin{array}{ll}
G^{free,2}&\textup{ if }x_{12}>1\\
SL_2(\Z)&\textup{ if }x_{12}=1.
\end{array}\right. 
\end{eqnarray*}
This matrix group has finite index in $SL_2(\Z)$ if and only
if $x_{12}\in\{1,2\}$. 
The exact sequence \eqref{6.10} splits non-canonically.

(e) The cases $\uuuu{x}=(-l,2,-l)$ with $l\geq 2$ even:

(This includes the case $\uuuu{x}=(-2,2,-2)$ which is the case
of $\HH_{1,2}$.) 
\begin{eqnarray*}
\Rad I^{(1)}&=&\Z f_3\quad\textup{with}\quad 
f_3=\frac{l}{2}(e_1+e_3)+e_2,\\
\oooo{H_\Z}^{(1)}&=& \Z \oooo{e_1}^{(1)}\oplus 
\Z\oooo{e_3}^{(1)},\quad 
\oooo{e_2}^{(1)}
=-\frac{l}{2}(\oooo{e_1}^{(1)}+\oooo{e_3}^{(1)}),\\
\end{eqnarray*}
\begin{eqnarray*}
\langle s_{e_1}^{(1)},s_{e_3}^{(1)}\rangle
&\cong& \langle \oooo{s_{e_1}^{(1)}},\oooo{s_{e_3}^{(1)}}\rangle
\cong \langle \begin{pmatrix}1&-2\\0&1\end{pmatrix},
\begin{pmatrix}1&0\\2&1\end{pmatrix}\rangle
\stackrel{1:2}{\subset}\Gamma(2),\\
\langle s_{e_1}^{(1)},s_{e_3}^{(1)}\rangle
&\cong& G^{free,2},
\end{eqnarray*}
\begin{eqnarray*}
(\ker\tau^{(1)})_u&=&\{t^+_\lambda\,|\, \lambda(\uuuu{e})
\in\langle (-2,0,2),(-2,l,0)\rangle_\Z\},\\
O^{(1),Rad}_u&=&\{t^+_\lambda\,|\, \lambda(\uuuu{e})
\in\langle (-1,0,1),(-1,\frac{l}{2},0)\rangle_\Z\}.
\end{eqnarray*}

(i) The cases with $l\equiv 0(4)$: 
$\Gamma^{(1)}_s\cong 
\langle \oooo{s_{e_1}^{(1)}},\oooo{s_{e_3}^{(1)}}\rangle
\cong G^{free,2}$. 
The isomorphism 
$\Gamma^{(1)}_s\cong \langle s_{e_1}^{(1)},s_{e_3}^{(1)}\rangle
\subset\Gamma^{(1)}$ gives a splitting of the exact sequence
\eqref{6.10}. Here $-\id\notin\Gamma^{(1)}_s$. 
\begin{eqnarray*}
\Gamma^{(1)}_u&=&\{t^+_\lambda\,|\, \lambda\in\langle
\lambda_1,\lambda_2\rangle_\Z\}\cong\Z^2 \quad\textup{with}\\
&&\lambda_1(\uuuu{e})=(-l,0,l),\quad 
\lambda_2(\uuuu{e})=(2l,-l^2,0).
\end{eqnarray*}

(ii) The cases with $l\equiv 2(4)$: Here $-\id\in\Gamma^{(1)}_s$. 
\begin{eqnarray*}
\Gamma^{(1)}_s&\cong& 
\langle \oooo{s_{e_1}^{(1)}},\oooo{s_{e_3}^{(1)}},-\id\rangle
\cong \langle \begin{pmatrix}1&-2\\0&1\end{pmatrix},
\begin{pmatrix}1&0\\2&1\end{pmatrix},
\begin{pmatrix}-1&0\\0&-1\end{pmatrix}\rangle
=\Gamma(2)\\
&\cong& G^{free,2}\times\{\pm 1\}.
\end{eqnarray*}
The isomorphism $\Gamma^{(1)}_s/\{\pm \id\}\cong
\langle s_{e_1}^{(1)},s_{e_2}^{(1)}\rangle\subset\Gamma^{(1)}$
gives a splitting of the exact sequence
\begin{eqnarray*}
\{1\}\to \Gamma^{(1)}\cap O^{(1),Rad}_\pm\to\Gamma^{(1)}\to
\Gamma^{(1)}_s/\{\pm\id\}\to\{1\}.
\end{eqnarray*}
\begin{eqnarray*}
\Gamma^{(1)}_u&=&\{t^+_\lambda\,|\, \lambda\in\langle
2\lambda_1,\lambda_2\rangle_\Z\}\cong\Z^2 \quad\textup{with}\\
&&\lambda_1(\uuuu{e})=(-l,0,l),\quad 
\lambda_2(\uuuu{e})=(2l,-l^2,0),\\
\Gamma^{(1)}\cap O^{(1),Rad}_\pm -\Gamma^{(1)}_u
&=&\{t^-_\lambda\,|\, \lambda\in\lambda_3+\langle
2\lambda_1,\lambda_2\rangle_\Z\}\quad\textup{with}\\
&&\lambda_3(\uuuu{e})=(-l,2,l),\quad\lambda_3(f_3)=2. 
\end{eqnarray*}

(f) The cases $\uuuu{x}=(-l,2,-l)$ with $l\geq 3$ odd:
\begin{eqnarray*}
\Rad I^{(1)}&=& \Z f_3\quad\textup{with}\quad 
f_3=l(e_1+e_3)+2e_2,\\
\oooo{H_\Z}^{(1)}&=& \Z\oooo{e_1}^{(1)}\oplus 
\Z\oooo{g_2}^{(1)}\quad\textup{with}\quad 
g_2:=\frac{1}{2}(e_1+e_3)-\frac{l}{2}f_3\in H_\Z,\\
\www{\uuuu{e}}&:=& (e_1,g_2,f_3)
\quad\textup{is a }\Z\textup{-basis of }H_\Z.
\end{eqnarray*}
Consider
\begin{eqnarray*}
s_4&:=& \bigl(s_{e_3}^{(1)}s_{e_1}^{(1)})^{\frac{l^2-1}{4}}
s_{e_2}^{(1)}\in\Gamma^{(1)}.
\end{eqnarray*}
Then 
\begin{eqnarray*}
\Gamma^{(1)}_s&\cong& \langle \oooo{s_{e_1}^{(1)}},
\oooo{s_4}\rangle \cong 
\langle s_{e_1}^{(1)},s_4\rangle\cong SL_2(\Z),
\end{eqnarray*}
and the isomorphism $\Gamma^{(1)}_s\cong \langle s_{e_1}^{(1)},
s_4\rangle\subset \Gamma^{(1)}$ gives a splitting of the
exact sequence \eqref{6.10}.
\begin{eqnarray*}
\Gamma^{(1)}_u&=&\{t^+_\lambda\,|\, \lambda\in \langle\lambda_1,
\lambda_2\rangle_\Z\}\cong\Z^2 \quad\textup{with}\\
&&\lambda_1(\uuuu{e})=(-l,0,l),\quad 
\lambda_2(\uuuu{e})=(2l,-l^2,0),\\
(\ker\tau^{(1)})_u&=& O^{(1),Rad}_u= 
\{t^+_\lambda\,|\, \lambda(\uuuu{e})\in \langle (-1,0,1),
(2,-l,0)\rangle_\Z\}.
\end{eqnarray*}

(g) The cases $\uuuu{x}\in\Z^3_{\geq 3}$ with
$2x_i\leq x_jx_k$ for $\{i,j,k\}=\{1,2,3\}$: 

(This includes the case $\uuuu{x}=(3,3,3)$ which is the case
of $\P^2$.)
\begin{eqnarray*}
\Rad I^{(1)}&=& \Z f_3\quad\textup{with}\quad 
f_3=-\www{x}_3e_1+\www{x}_2e_2-\www{x}_1e_3,\\
\Gamma^{(1)}_u&=&\{\id\}\subsetneqq (\ker\tau^{(1)})_u\cong\Z^2,\\
\Gamma^{(1)}&\cong& \Gamma^{(1)}_s\cong G^{free,3},
\end{eqnarray*}
$\Gamma^{(1)}$ and $\Gamma^{(1)}_s$ are free groups
with the three generators $s_{e_1}^{(1)}$, $s_{e_2}^{(1)}$, 
$s_{e_3}^{(1)}$ respectively 
$\oooo{s_{e_1}^{(1)}}$, $\oooo{s_{e_2}^{(1)}}$, 
$\oooo{s_{e_3}^{(1)}}$. 
\end{theorem}

{\bf Proof:}
(a) This follows from the definitions in Lemma \ref{t2.6} (a)
and in Definition \ref{t2.8}.

(b) This follows from Lemma \ref{t2.11} and Lemma \ref{t2.12}.

(c)--(f) In (c)--(f) $(\ker\tau^{(1)})_u$ and 
$O^{(1),Rad}_u$ are calculated with Lemma \ref{t6.17} (a)
and (b). 

(c) The first statements $\Rad I^{(1)}=\Z f_3$
and $\oooo{H_\Z}^{(1)}=\Z\oooo{e_1}^{(1)}\oplus \Z\oooo{e_2}^{(1)}$
are known respectively obvious. 

Also $(e_1,e_2,f_3)$ is a $\Z$-basis of $H_\Z$. With respect
to this basis
\begin{eqnarray*}
s_{e_1}^{(1)}(e_1,e_2,f_3)=(e_1,e_2,f_3)
\begin{pmatrix}1&-1&0\\0&1&0\\0&0&1\end{pmatrix},\\
s_{e_2}^{(1)}(e_1,e_2,f_3)=(e_1,e_2,f_3)
\begin{pmatrix}1&0&0\\1&1&0\\0&0&1\end{pmatrix}.
\end{eqnarray*}
This shows 
$\langle s_{e_1}^{(1)},s_{e_2}^{(1)}\rangle\cong\langle
\oooo{s_{e_1}^{(1)}},\oooo{s_{e_2}^{(1)}}\rangle
\cong SL_2(\Z)$.
Together with Lemma \ref{t6.16} (b) and Lemma \ref{t6.2} (d)
we obtain
\begin{eqnarray*}
\Gamma^{(1)}_s= \langle
\oooo{s_{e_1}^{(1)}},\oooo{s_{e_2}^{(1)}}\rangle
=(\ker\tau^{(1)})_s=O^{(1),Rad}_s\cong SL_2(\Z),
\end{eqnarray*}
and that the exact sequence \eqref{6.10} splits non-canonically
with $\Gamma^{(1)}_s\cong 
\langle s_{e_1}^{(1)},s_{e_2}^{(1)}\rangle\subset\Gamma^{(1)}$. 

From the actions of $s_{e_1}^{(1)}$ and $s_{e_2}^{(1)}$ on 
$(e_1,e_2,f_3)$ and from
\begin{eqnarray*}
s_{e_3}^{(1)}((e_1,e_2,f_3))=(e_1,e_2,f_3)
\begin{pmatrix}1&0&0&\\1&1&0\\-1&0&1\end{pmatrix},
\end{eqnarray*}
one sees that the map
\begin{eqnarray*}
\Bigl( (e_1,e_2,f_3)\mapsto (e_1,e_2,f_3)\begin{pmatrix}
1&0&0\\0&1&0\\-1&0&1\end{pmatrix}\Bigr) =t^+_{-\lambda_1}
\end{eqnarray*}
is in $\Gamma^{(1)}_u$ and that
\begin{eqnarray*}
\Gamma^{(1)}=\langle s_{e_1}^{(1)},s_{e_2}^{(1)},
t^+_{\lambda_1}\rangle.
\end{eqnarray*}
Also 
\begin{eqnarray*}
(s_{e_1}^{(1)})^{-1}\circ t^+_{\lambda_1}\circ s_{e_1}^{(1)}
=t^+_{\lambda_1\circ s_{e_1}^{(1)}},\quad\textup{with}\quad 
t^+_{\lambda_1\circ s_{e_1}^{(1)}}(\uuuu{e})=(1,-1,-1),
\end{eqnarray*}
and therefore $t^+_{\lambda_2}= t^+_{\lambda_1}
-t^+_{\lambda_1\circ s_{e_1}^{(1)}} \in\Gamma^{(1)}_u$.
But $O^{(1),Rad}_u =\langle t^+_{\lambda_1},
t^+_{\lambda_2}\rangle_\Z$, so 
\begin{eqnarray*}
\Gamma^{(1)}_u=(\ker \tau^{(1)})_u = O^{(1),Rad}_u
=\langle t^+_{\lambda_1},t^+_{\lambda_2}\rangle_\Z
\end{eqnarray*}
In the diagram of exact sequences in Lemma \ref{t6.2} (d)
the inclusions in the second and fourth columns are bijections, 
so also the inclusions in the middle column, 
\begin{eqnarray*}
\Gamma^{(1)}=\ker\tau^{(1)}=O^{(1),Rad}.
\end{eqnarray*}

(d) $f_3=\www{x}_2e_2-\www{x}_1e_3$ implies
$\www{x}_2\oooo{e_2}^{(1)}=\www{x}_1\oooo{e_3}^{(1)}$. 
Lemma \ref{t6.16} (c) gives
$\oooo{e_2}^{(1)}\in \www{x}_1\oooo{H_\Z}^{(1),prim}$, so 
$g_2:=\frac{1}{\www{x}_1}\oooo{e_2}^{(1)}
=\frac{1}{\www{x}_2}\oooo{e_3}^{(1)}$ is in
$\oooo{H_\Z}^{(1),prim}$. Thus 
\begin{eqnarray*}
\oooo{H_\Z}^{(1)}=\Z\oooo{e_1}^{(1)}+\Z\oooo{e_2}^{(1)}
+\Z\oooo{e_3}^{(1)}=\Z\oooo{e_1}^{(1)}\oplus\Z g_2.
\end{eqnarray*}
First we consider $\Gamma^{(1)}_s$. Define 
$\uuuu{g}:=(\oooo{e_1}^{(1)},g_2)$. One sees
\begin{eqnarray*}
\oooo{s_{e_1}^{(1)}}\uuuu{g}=\uuuu{g}
\begin{pmatrix}1&-x_{12}\\0&1\end{pmatrix},\quad
\oooo{s_{e_2}^{(1)}}\uuuu{g}=\uuuu{g}
\begin{pmatrix}1&0\\x_1\www{x}_1&1\end{pmatrix},\quad
\oooo{s_{e_3}^{(1)}}\uuuu{g}=\uuuu{g}
\begin{pmatrix}1&0\\x_2\www{x}_2&1\end{pmatrix}.
\end{eqnarray*}

Choose $y_1,y_2\in\Z$ with $1=y_1\www{x}_1^2+y_2\www{x}_2^2$
and define
\begin{eqnarray*}
s_4:= (s_{e_2}^{(1)})^{y_1}(s_{e_3}^{(1)})^{y_2}\in\Gamma^{(1)}.
\end{eqnarray*}
Then
\begin{eqnarray*}
\oooo{s_4}\,\uuuu{g}=\uuuu{g}
\begin{pmatrix}1&0\\x_{12}&1\end{pmatrix},\quad
s_4\,\uuuu{e}=\uuuu{e}\cdot s_4^{mat}\textup{ with }
s_4^{mat}=\begin{pmatrix}1&0&0\\y_1x_1&1&0\\y_2x_2&0&1
\end{pmatrix}.
\end{eqnarray*}

$\oooo{s_{e_3}^{(1)}}$ and $\oooo{s_{e_2}^{(1)}}$ are powers
of $\oooo{s_4}$, so $\Gamma^{(1)}_s=\langle \oooo{s_{e_1}^{(1)}},
\oooo{s_4}\rangle$, so 
\begin{eqnarray*}
\Gamma^{(1)}_s\cong\langle 
\begin{pmatrix}1&-x_{12}\\0&1\end{pmatrix},
\begin{pmatrix}1&0\\x_{12}&1\end{pmatrix}\rangle
=\Gamma^{(1),mat}(S(-x_{12}))\subset SL_2(\Z).
\end{eqnarray*}
The matrix group $\Gamma^{(1),mat}(S(-x_{12}))$ was treated
in Theorem \ref{t6.10} (b)--(d): 
\begin{eqnarray*}
\Gamma^{(1),mat}(S(-x_{12}))\cong\langle 
\begin{pmatrix}1&-x_{12}\\0&1\end{pmatrix},
\begin{pmatrix}1&0\\x_{12}&1\end{pmatrix}\rangle
\left\{\begin{array}{ll}
\cong G^{free,2}&\textup{ if }x_{12}>1,\\
= SL_2(\Z)&\textup{ if }x_{12}=1.
\end{array}\right.
\end{eqnarray*}

{\bf The case $x_{12}>1$:}
Then $\Gamma^{(1)}_s=\langle \oooo{s_{e_1}^{(1)}},\oooo{s_4}
\rangle$ and $\langle s_{e_1}^{(1)},s_4\rangle\subset\Gamma^{(1)}$
are free groups with the two given generators. 
Then $\langle s_{e_1}^{(1)},s_4\rangle\subset\Gamma^{(1)}$
gives a splitting of the exact sequence \eqref{6.10}.
The generating relations in $\Gamma^{(1)}_s$ with respect to the
four elements $\oooo{s_{e_1}^{(1)}}$, $\oooo{s_{e_2}^{(1)}}$, 
$\oooo{s_{e_3}^{(1)}}$ and $\oooo{s_4}$ are
\begin{eqnarray}\label{6.11}
\oooo{s_{e_2}^{(1)}}=(\oooo{s_4})^{\www{x}_1^2},\quad 
\oooo{s_{e_3}^{(1)}}=(\oooo{s_4})^{\www{x}_2^2}.
\end{eqnarray}
Therefore $\Gamma^{(1)}_u$ is the normal subgroup of 
$\Gamma^{(1)}$ generated by the elements
$s_4^{\www{x}_1^2}(s_{e_2}^{(1)})^{-1}$ and
$s_4^{\www{x}_2^2}(s_{e_3}^{(1)})^{-1}$. 

{\bf The case $x_{12}=1$:}
Then $\Gamma^{(1)}_s=\langle \oooo{s_{e_1}^{(1)}},\oooo{s_4}
\rangle\cong SL_2(\Z)$ 

\medskip
{\bf Claim:} Also $\langle s_{e_1}^{(1)},s_4\rangle\cong
SL_2(\Z)$.

\medskip
{\sf Proof of the Claim:}
The generating relations in $SL_2(\Z)$ with respect to the 
generators 
$\oooo{s_{e_1}^{(1)}}^{mat}=\begin{pmatrix}1&-1\\0&1\end{pmatrix}$
and $\oooo{s_4}^{mat}=\begin{pmatrix}1&0\\1&1\end{pmatrix}$ are
\begin{eqnarray*}
\oooo{s_{e_1}^{(1)}}^{mat} \oooo{s_4}^{mat} 
\oooo{s_{e_1}^{(1)}}^{mat}
=\oooo{s_4}^{mat} \oooo{s_{e_1}^{(1)}}^{mat}
\oooo{s_4}^{mat}\quad\textup{and}\quad 
E_2=\bigl(\oooo{s_{e_1}^{(1)}}^{mat}\oooo{s_4}^{mat}\bigr)^6.
\end{eqnarray*}
One checks with calculations that they lift to $\Gamma^{(1),mat}$,
\begin{eqnarray*}
s_{e_1}^{(1),mat} s_4^{mat} s_{e_1}^{(1),mat}
=s_4^{mat} s_{e_1}^{(1),mat} s_4^{mat}\quad\textup{and}\quad 
E_3=\bigl(s_{e_1}^{(1),mat} s_4^{mat}\bigr)^6.
\hspace*{0.5cm}(\Box)
\end{eqnarray*}

\medskip
Because of the Claim, $\langle s_{e_1}^{(1)},s_4\rangle\subset
\Gamma^{(1)}$ gives a splitting of the exact sequence 
\eqref{6.10}.
The generating relations in $\Gamma^{(1)}_s$ with respect to
the four elements 
$\oooo{s_{e_1}^{(1)}}$, $\oooo{s_{e_2}^{(1)}}$, 
$\oooo{s_{e_3}^{(1)}}$ and $\oooo{s_4}$ are the relations
between $\oooo{s_{e_1}^{(1)}}$ and $\oooo{s_4}$ in the proof 
of the Claim and the relations in \eqref{6.11}. 
The relations in the proof of the Claim lift to $\Gamma^{(1)}$.
Therefore again $\Gamma^{(1)}_u$ is the normal subgroup of 
$\Gamma^{(1)}$ generated by the elements
$s_4^{\www{x}_1^2}(s_{e_2}^{(1)})^{-1}$ and 
$s_4^{\www{x}_2^2}(s_{e_3}^{(1)})^{-1}$.

{\bf Back to both cases $x_{12}\geq 1$ together:}
We have to determine these two elements of $\Gamma^{(1)}_u$.
The first one is given by the following calculation,
\begin{eqnarray*}
&&s_4^{\www{x}_1^2}(s_{e_2}^{(1)})^{-1}(\uuuu{e})\\
&=&\uuuu{e}
\begin{pmatrix}1&0&0\\y_1x_1\www{x}_1^2&1&0\\y_2x_2\www{x}_1^2
&0&1\end{pmatrix}
\begin{pmatrix}1&0&0\\-x_1&1&0\\0&0&1\end{pmatrix}
=\uuuu{e}\begin{pmatrix}1&0&0\\-x_1+y_1x_1\www{x}_1^2&1&0\\
y_2x_2\www{x}_1^2 &0&1\end{pmatrix}\\
&=&\uuuu{e}\begin{pmatrix}1&0&0\\-y_2x_1\www{x}_2^2&1&0\\
y_2x_2\www{x}_1^2&0&1\end{pmatrix}
=\uuuu{e}+f_3(-y_2x_{12}\www{x}_1\www{x}_2,0,0),\\
&&\textup{so }
s_4^{\www{x}_1^2}(s_{e_2}^{(1)})^{-1}(\uuuu{e})
=t^+_{-y_2\lambda_1}.
\end{eqnarray*}
A similar calculation gives
\begin{eqnarray*}
s_4^{\www{x}_2^2}(s_{e_3}^{(1)})^{-1}
&=& t^+_{y_1\lambda_1}.
\end{eqnarray*}
Observe $\gcd(y_1,y_2)=1$. 
Therefore $\Gamma^{(1)}_u$ is the normal subgroup of 
$\Gamma^{(1)}$ generated by $t^+_{\lambda_1}$. 
Lemma \ref{t6.17} (c) and the calculations 
\begin{eqnarray*}
\lambda_1\circ s_{e_1}^{(1)}(\uuuu{e})=x_{12}\www{x}_1\www{x}_2
(1,-x_1,-x_2),\\
\textup{so }
\langle \lambda_1,\lambda_1\circ s_{e_1}^{(1)}\rangle_\Z
=\langle \lambda_1,\lambda_2\rangle_\Z,\\
\lambda_1\circ (s_{e_i}^{(1)})^{\pm 1},
\lambda_2\circ (s_{e_i}^{(1)})^{\pm 1}\in 
\langle\lambda_1,\lambda_2\rangle_\Z,
\end{eqnarray*}
give 
\begin{eqnarray*}
\Gamma^{(1)}_u=\{t^+_\lambda\,|\, \lambda\in 
\langle\lambda_1,\lambda_2\rangle_\Z\}.
\end{eqnarray*}

(e) The first statements $\Rad I^{(1)}=\Z f_3$, 
$f_3=\frac{l}{2}(e_1+e_3)+e_2$ and 
$\oooo{H_\Z}^{(1)}=\Z \oooo{e_1}^{(1)}\oplus
\Z\oooo{e_3}^{(1)}$ are obvious, also
\begin{eqnarray*}
\oooo{s_{e_i}^{(1)}}(\oooo{e_1}^{(1)},\oooo{e_3}^{(1)})
=(\oooo{e_1}^{(1)},\oooo{e_3}^{(1)})\oooo{s_{e_i}^{(1)}}^{mat}
\quad\textup{with}\quad
\oooo{s_{e_1}^{(1)}}^{mat}
=\begin{pmatrix}1&-2\\0&1\end{pmatrix},\\
\oooo{s_{e_2}^{(1)}}^{mat}
=\begin{pmatrix}1+\frac{l^2}{2}&-\frac{l^2}{2}\\ 
\frac{l^2}{2}&1-\frac{l^2}{2}\end{pmatrix},\quad 
\oooo{s_{e_3}^{(1)}}^{mat}
=\begin{pmatrix}1&0\\2&1\end{pmatrix}.
\end{eqnarray*}
By Theorem \ref{t6.10} (c) 
$\langle \oooo{s_{e_1}^{(1)}}^{mat},
\oooo{s_{e_3}^{(1)}}^{mat}\rangle\cong G^{free,2}$ and
therefore
\begin{eqnarray*}
\langle s_{e_1}^{(1)},s_{e_3}^{(1)}\rangle\cong G^{free,2}.
\end{eqnarray*}
We have 
\begin{eqnarray*} 
\oooo{s_{e_2}^{(1)}}^{mat}&\equiv& 
\begin{pmatrix}1&0\\0&1\end{pmatrix}\mmod 4\quad
\textup{if }l\equiv 0(4),\\
\oooo{s_{e_2}^{(1)}}^{mat}&\equiv& 
\begin{pmatrix}3&2\\2&3\end{pmatrix}\mmod 4\quad
\textup{if }l\equiv 2(4).
\end{eqnarray*}
If $l\equiv 0(4)$ then by Theorem \ref{t6.10} (c)
$\oooo{s_{e_2}^{(1)}}^{mat}\in 
\langle \oooo{s_{e_1}^{(1)}}^{mat},
\oooo{s_{e_3}^{(1)}}^{mat}\rangle$, so then
$\Gamma^{(1)}_s=\langle \oooo{s_{e_1}^{(1)}},
\oooo{s_{e_3}^{(1)}}\rangle,$
and the isomorphism $\Gamma^{(1)}_s \cong 
\langle s_{e_1}^{(1)},s_{e_3}^{(1)}\rangle
\subset\Gamma^{(1)}$ gives a splitting of the exact
sequence \eqref{6.10}.

If $l\equiv 2(4)$ then by Theorem \ref{t6.10} (c)
$\oooo{s_{e_2}^{(1)}}^{mat}\notin 
\langle \oooo{s_{e_1}^{(1)}}^{mat},
\oooo{s_{e_3}^{(1)}}^{mat}\rangle$, so then
$\langle \oooo{s_{e_i}^{(1)}}^{mat}\,|\, 
i\in\{1,2,3\} \rangle=\Gamma^{(2)}
\cong G^{free,2}\times\{\pm 1\}$, $-\id\in\Gamma^{(1)}_s$,
and the isomorphism $\Gamma^{(1)}_s/\{\pm \id\} \cong 
\langle s_{e_1}^{(1)},s_{e_3}^{(1)}\rangle
\subset\Gamma^{(1)}$ gives a splitting of the exact
sequence in part (ii). 

\medskip
{\bf Claim:} $\Gamma^{(1)}_u$ if $l\equiv 0(4)$ and
$\Gamma^{(1)}\cap O^{(1),Rad}_\pm$ if $l\equiv 2(4)$
is the normal subgroup generated by
$s_4:=\bigl(s_{e_3}^{(1)}s_{e_1}^{(1)}\bigr)^{l^2/4}
s_{e_2}^{(1)}$.

\medskip
{\sf Proof of the Claim:}
Consider the $\Z$-basis $\www{\uuuu{e}}:=(e_1,e_1+e_3,f_3)$ 
of $H_\Z$. 
\begin{eqnarray*}
s_{e_1}^{(1)}\www{\uuuu{e}}&=&\www{\uuuu{e}}
\begin{pmatrix}1&-2&0\\0&1&0\\0&0&1\end{pmatrix},\  
s_{e_2}^{(1)}\www{\uuuu{e}}=\www{\uuuu{e}}
\begin{pmatrix}1&0&0\\\frac{l^2}{2}&1&0\\-l&0&1\end{pmatrix},\ 
s_{e_3}^{(1)}\www{\uuuu{e}}=\www{\uuuu{e}}
\begin{pmatrix}-1&-2&0\\2&3&0\\0&0&1\end{pmatrix},\\
s_4\, \www{\uuuu{e}}&=&\www{\uuuu{e}}
\begin{pmatrix}-1&0&0\\2&-1&0\\0&0&1\end{pmatrix}^{l^2/4}
\begin{pmatrix}1&0&0\\\frac{l^2}{2}&1&0\\-l&0&1\end{pmatrix}\\
&=&\www{\uuuu{e}}
\begin{pmatrix}-1&0&0\\0&-1&0\\0&0&1\end{pmatrix}^{l/2}
\begin{pmatrix}1&0&0\\-\frac{l^2}{2}&1&0\\0&0&1\end{pmatrix}
\begin{pmatrix}1&0&0\\\frac{l^2}{2}&1&0\\-l&0&1\end{pmatrix}\\
&=&\www{\uuuu{e}}
\begin{pmatrix}-1&0&0\\0&-1&0\\0&0&1\end{pmatrix}^{l/2}
\begin{pmatrix}1&0&0\\0&1&0\\-l&0&1\end{pmatrix}
=\www{\uuuu{e}}
\begin{pmatrix}(-1)^{l/2}&0&0\\0&(-1)^{l/2}&0\\-l&0&1\end{pmatrix},
\end{eqnarray*}
\begin{eqnarray}
s_4\, \uuuu{e}=(-1)^{l/2}\uuuu{e}+f_3(-l,1-(-1)^{l/2},l),\nonumber\\
s_4=\left\{\begin{array}{ll}t^+_{\lambda_1}\in\Gamma^{(1)}_u&
\textup{ if }l\equiv 0(4),\\
t^-_{\lambda_3}\in \Gamma^{(1)}\cap O^{(1),Rad}_\pm 
& \textup{ if }l\equiv 2(4).
\end{array}\right. \label{6.12}
\end{eqnarray}
The splittings above of the exact sequences give semidirect
products
\begin{eqnarray*}
\Gamma^{(1)}=\left\{\begin{array}{ll}
\Gamma^{(1)}_u\rtimes 
\langle s_{e_1}^{(1)},s_{e_3}^{(1)}\rangle 
& \textup{ if }l\equiv 0(4),\\
\Gamma^{(1)}\cap O^{(1),Rad}_\pm \rtimes 
\langle s_{e_1}^{(1)},s_{e_3}^{(1)}\rangle 
& \textup{ if }l\equiv 2(4).
\end{array}\right.
\end{eqnarray*}
$s_{e_2}^{(1)}$ turns up linearly in $s_4$.
Therefore $\Gamma^{(1)}=\langle s_{e_1}^{(1)},
s_{e_3}^{(1)},s_4\rangle$. Together these facts show
that $\Gamma^{(1)}_u$ respectively 
$\Gamma^{(1)}\cap O^{(1),Rad}_\pm $ is the normal 
subgroup generated by $s_4$. 
This finishes the proof of the Claim.
\hfill($\Box$)

\medskip
Now we can apply Lemma \ref{t6.17} (c) if $l\equiv 0(4)$
and Lemma \ref{t6.17} (d) if $l\equiv 2(4)$.
The following calculations show the claims on $\Gamma^{(1)}_u$
and (in the case $l\equiv 2(4)$) 
$\Gamma^{(1)}\cap O^{(1),Rad}_\pm$. 

{\bf The case $l\equiv 0(4)$:}
\begin{eqnarray*}
\lambda_1\circ s_{e_1}^{(1)}(\uuuu{e})=(-l,-l^2,3l),\\
\textup{so }\lambda_1\circ s_{e_1}^{(1)}
=3\lambda_1+\lambda_2,\\
\lambda_1\circ (s_{e_i}^{(1)})^{\pm 1},
\lambda_2\circ (s_{e_i}^{(1)})^{\pm 1}
\in\langle \lambda_1,\lambda_2\rangle_\Z.
\end{eqnarray*}

{\bf The case $l\equiv 2(4)$:}
\begin{eqnarray*}
\lambda_3\circ s_{e_1}^{(1)}(\uuuu{e})=(-l,2-l^2,3l),\\
\textup{so }\lambda_3\circ s_{e_1}^{(1)}
=\lambda_3+2\lambda_1+\lambda_2,\\
\lambda_3\circ (s_{e_2}^{(1)})^{-1}(\uuuu{e})=(l,2,-l),\\
\textup{so }\lambda_3\circ (s_{e_2}^{(1)})^{-1}
=\lambda_3-2\lambda_1,\\
\lambda_3\circ (s_{e_i}^{(1)})^{\pm 1}
\in\lambda_3+\langle 2\lambda_1,\lambda_2\rangle_\Z,\\
2\lambda_1\circ (s_{e_i}^{(1)})^{\pm 1},
\lambda_2\circ (s_{e_i}^{(1)})^{\pm 1}
\in\langle 2\lambda_1,\lambda_2\rangle_\Z.
\end{eqnarray*}

(f) $\Rad I^{(1)}=\Z f_3$ is known.
$e_2=-\frac{l}{2}(e_1+e_3)+\frac{1}{2}f_3$ 
shows $\oooo{H_\Z}^{(1)}=\Z \oooo{e_1}^{(1)}
\oplus \Z (\frac{1}{2}(\oooo{e_1}^{(1)}+ 
\oooo{e_3}^{(1)}))$. We have
\begin{eqnarray*}
\www{\uuuu{e}}=\uuuu{e}
\begin{pmatrix} 1&\frac{1-l^2}{2}&l\\
0&-l&2\\0&\frac{1-l^2}{2}&l \end{pmatrix},\quad 
\uuuu{e}=\www{\uuuu{e}}
\begin{pmatrix} 1&0&-1\\0&-l&2\\0&\frac{1-l^2}{2}&l
\end{pmatrix},
\end{eqnarray*}
so $\www{\uuuu{e}}$ is a $\Z$-basis of $H_\Z$. 

First we calculate $s_4$ with respect to the 
$\Q$-basis $\www{\uuuu{f}}
:=(e_1,\frac{1}{2}(e_1+e_3),f_3)$ of $H_\Q$,
\begin{eqnarray*}
s_4(\www{\uuuu{f}})&=&\www{\uuuu{f}}
\bigl( \begin{pmatrix}-1&-1&0\\4&3&0\\0&0&1\end{pmatrix}
\begin{pmatrix}1&-1&0\\0&1&0\\0&0&1\end{pmatrix}
\bigr)^{\frac{l^2-1}{4}}
\begin{pmatrix}1&0&0\\ l^2&1&0\\
-\frac{l}{2}&0&1\end{pmatrix}\\
&=& \www{\uuuu{f}}
\begin{pmatrix}-1&0&0\\4&-1&0\\0&0&1
\end{pmatrix}^{\frac{l^2-1}{4}}
\begin{pmatrix}1&0&0\\ l^2&1&0\\
-\frac{l}{2}&0&1\end{pmatrix}\\
\\
&=& \www{\uuuu{f}}
\begin{pmatrix}1&0&0\\1-l^2&1&0\\0&0&1\end{pmatrix}
\begin{pmatrix}1&0&0\\ l^2&1&0\\
-\frac{l}{2}&0&1\end{pmatrix}\
=\www{\uuuu{f}}
\begin{pmatrix}1&0&0\\1&1&0\\-\frac{l}{2}&0&1
\end{pmatrix}.
\end{eqnarray*}
Therefore
\begin{eqnarray}\label{6.13}
s_4(\www{\uuuu{e}})=\www{\uuuu{e}}
\begin{pmatrix}1&0&0\\1&1&0\\0&0&1\end{pmatrix},\quad
s_{e_1}^{(1)}(\www{\uuuu{e}})=\www{\uuuu{e}}
\begin{pmatrix}1&-1&0\\0&1&0\\0&0&1\end{pmatrix}.
\end{eqnarray}
Thus
\begin{eqnarray*}
\Gamma^{(1)}_s=\langle \oooo{s_{e_1}^{(1)}},\oooo{s_4}
\rangle \cong
\langle s_{e_1}^{(1)},s_4\rangle\subset\Gamma^{(1)},
\end{eqnarray*}
and this isomorphism gives a splitting of the exact
sequence \eqref{6.10}. 

\eqref{6.13} shows that the group
$\langle s_{e_1}^{(1)},s_4\rangle$ contains an element
$s_5$ with 
\begin{eqnarray*}
s_5(\www{\uuuu{e}})=\www{\uuuu{e}}
\begin{pmatrix}3&1&0\\-4&-1&0\\0&0&1\end{pmatrix},
\end{eqnarray*}
thus 
\begin{eqnarray*}
s_5s_{e_3}^{(1)}(\www{\uuuu{e}})&=&\www{\uuuu{e}}
\begin{pmatrix}3&1&0\\-4&-1&0\\0&0&1\end{pmatrix}
\begin{pmatrix}-1&-1&0\\4&3&0\\2l&l&1\end{pmatrix}
=\www{\uuuu{e}}
\begin{pmatrix}1&0&0\\0&1&0\\2l&l&1\end{pmatrix},\\
s_5s_{e_3}^{(1)}(\uuuu{e})&=&\uuuu{e}+
f_3(2l,-l^2,0),\\
s_5s_{e_3}^{(1)}&=& t^+_{\lambda_2}\quad\textup{with}
\quad \lambda_2(\uuuu{e})=(2l,-l^2,0).
\end{eqnarray*}
Also 
\begin{eqnarray*}
\Gamma^{(1)}=\langle s_{e_1}^{(1)},s_{e_2}^{(1)},
s_{e_3}^{(1)}\rangle
=\langle s_4,s_{e_1}^{(1)},s_{e_3}^{(1)}\rangle
=\langle s_4,s_{e_1}^{(1)},t^+_{\lambda_2}\rangle,
\end{eqnarray*}
so $\Gamma^{(1)}_u$ is the normal subgroup generated 
by $t^+_{\lambda_2}$. 

Now we can apply Lemma \ref{t6.17} (c).  
The following formulas show $\Gamma^{(1)}_u=\{t^+_\lambda\,|\, 
\lambda\in\langle \lambda_1,\lambda_2\rangle_\Z\}.$
\begin{eqnarray*}
(2l,-l^2,0) s_{e_1}^{(1),mat}=(2l,l^2,-4l)
=(2l,-l^2,0)+(0,2l^2,-4l),\\
2(2l,-l^2,0)+(0,2l^2,-4l)=(4l,0,-4l),\\
(2l,-l^2,0)s_{e_2}^{(1),mat}=(2l,-l^2,0)+(l^3,0,-l^3),\\
\gcd(4l,l^3)=l,\quad 
\textup{so }\lambda_1\in\Gamma^{(1)}_u,\\
\lambda_1\circ (s_{e_i}^{(1)})^{\pm 1},
\lambda_2\circ (s_{e_i}^{(1)})^{\pm 1}\in 
\langle\lambda_1,\lambda_2\rangle_\Z.
\end{eqnarray*}

(g) Because of the cyclic action of 
$\gamma\in (G^{phi}\ltimes \www{G}^{sign})\rtimes
\langle\gamma\rangle$ on $\Z^3$, we can suppose
$x_1=\max(x_1,x_2,x_3)$. With respect to the 
$\Q$-basis $(\oooo{e_1}^{(1)},\oooo{e_2}^{(1)})$ of 
$\oooo{H_\Q}^{(1)}$, $\oooo{s_{e_1}^{(1)}}$, 
$\oooo{s_{e_2}^{(1)}}$ and $\oooo{s_{e_3}^{(1)}}$
take the shape
$\oooo{s_{e_i}^{(1)}}(\oooo{e_1}^{(1)},\oooo{e_2}^{(1)})
=(\oooo{e_1}^{(1)},\oooo{e_2}^{(1)})B_i$ with
\begin{eqnarray*}
B_1=\begin{pmatrix}1&-x_1\\0&1\end{pmatrix},\quad
B_2=\begin{pmatrix}1&0\\x_1&1\end{pmatrix},\quad
B_3=\begin{pmatrix}1-\frac{x_3x_2}{x_1}&-\frac{x_3^2}{x_1}\\
\frac{x_2^2}{x_1}&1+\frac{x_3x_2}{x_1}\end{pmatrix}.
\end{eqnarray*}

We will show that the group 
$\langle \mu_1,\mu_2,\mu_3\rangle$
of M\"obius transformations $\mu_i:=\mu(B_i)$
is a free group with the three generators 
$\mu_1,\mu_2,\mu_3$. This implies first 
$\Gamma^{(1)}_s\cong G^{free,3}$ and then 
$\Gamma^{(1)}\cong \Gamma^{(1)}_s\cong G^{free,3}$
and $\Gamma^{(1)}_u=\{\id\}$.

We will apply Theorem \ref{ta.2} (c). 
$\mu_1,\mu_2$ and $\mu_3$ are parabolic,
\begin{eqnarray*}
\mu_1&=&\bigl(z\mapsto z-x_1\bigr)
\quad\textup{with fixed point }\infty,\\
\mu_2&=&\Bigl(z\mapsto \frac{z}{x_1z+1}\Bigr)
\quad\textup{with fixed point }0,\\
\mu_3&=&\Bigl(z\mapsto \frac{(1-\frac{x_3x_2}{x_1})z
-\frac{x_3^2}{x_1}}
{\frac{x_2^2}{x_1}z+(1+\frac{x_3x_2}{x_1})}
=\frac{(x_1-x_2x_3)z-x_3^2}{x_2^2z+(x_1+x_2x_3)}\Bigr)\\
&&\hspace*{2cm}\textup{with fixed point }-\frac{x_3}{x_2}.
\end{eqnarray*}
Consider 
\begin{eqnarray*}
r_1:=\mu_1(1)=1-x_1,\quad 
r_2:=\mu_2^{-1}(1)=\frac{1}{1-x_1},\\
r_3:=\mu_3(r_1)
=\frac{(x_1-x_2x_3)(1-x_1)-x_3^2}
{x_2^2(1-x_1)+(x_1+x_2x_3)}.
\end{eqnarray*}
It is sufficient to show the inequalities
\begin{eqnarray}\label{6.14}
-\infty < r_1<-\frac{x_3}{x_2}< r_3\leq r_2<0<1<\infty.
\end{eqnarray}

\begin{figure}[H]
\includegraphics[width=0.8\textwidth]{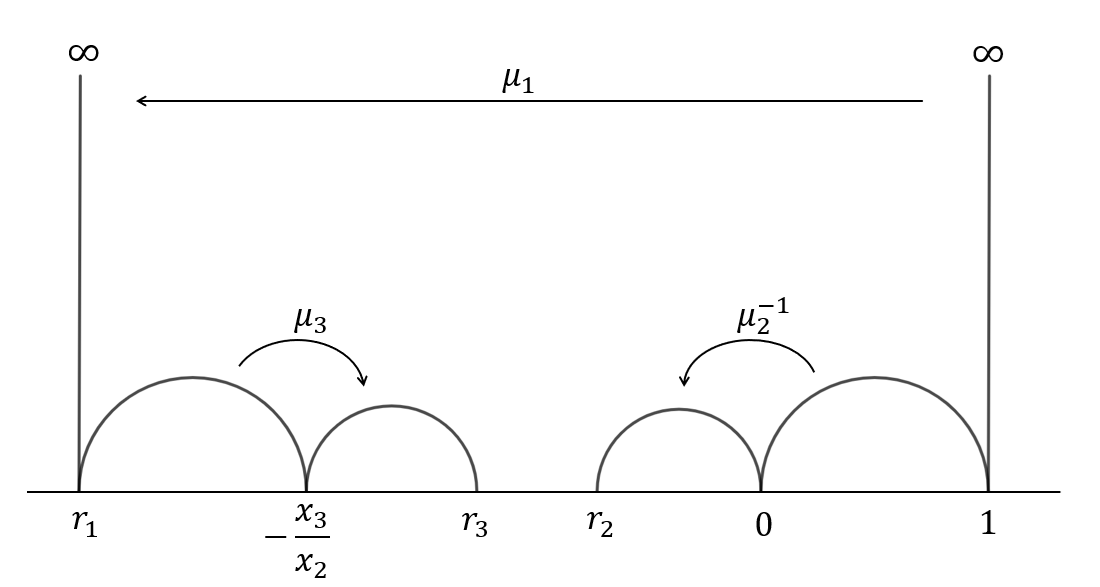}
\caption[Figure 6.6]{$G^{free,3}$ generated by three parabolic
M\"obius transformations, an application of Theorem 
\ref{ta.2} (c)}
\label{Fig:6.6}
\end{figure}

Then Theorem \ref{ta.2} (c) applies. Compare Figure 
\ref{Fig:6.6}.

We treat the case $\uuuu{x}=(3,3,3)$ separately and first. Then
\begin{eqnarray*}
(r_1,-\frac{x_3}{x_2},r_3,r_2)=(-2,-1,-\frac{1}{2},-\frac{1}{2}),
\end{eqnarray*}
so then \eqref{6.14} holds. 

From now on we suppose $x_1\geq 4$. 
The inequality $r_2<0$ is trivial. Consider the number
\begin{eqnarray*}
r_4:=-\frac{x_1+x_2x_3}{x_2^2}
=-\frac{x_1}{x_2^2}-\frac{x_3}{x_2}
<-\frac{x_3}{x_2},\quad\textup{with }
\mu_3(r_4)=\infty.
\end{eqnarray*}
We will show in this order the following claims:
\begin{list}{}{}
\item[(i)]
$r_1<r_4$, which implies $r_1<-\frac{x_3}{x_2}$.
\item[(ii)]
$r_3\in (-\frac{x_3}{x_2},\infty)$.
\item[(iii)]
The numerator $(x_2x_3-x_1)(x_1-1)-x_3^2$ 
of $r_3$ is positive.
\item[(iv)]
The denominator $x_2^2(1-x_1)+(x_1+x_2x_3)$ 
of $r_3$ is negative.
\item[(v)]
$r_3\leq r_2$.
\end{list}
Together (i)--(v) and $r_2<0$ show \eqref{6.14}.

Two of the three inequalities $2x_i\leq x_jx_k$ for 
$\{i,j,k\}=\{1,2,3\}$ in the assumption on $\uuuu{x}$ 
can be improved,
\begin{eqnarray*}
x_1x_2\geq 3x_1\geq 3x_3,\quad 
x_1x_3\geq 3x_1\geq 3x_2,\quad
\textup{we keep}\quad x_2x_3\geq 2x_1.
\end{eqnarray*}

(i) holds because
\begin{eqnarray*}
r_1<r_4\iff 1<x_1-\frac{x_1}{x_2^2}-\frac{x_3}{x_2}
\quad\textup{and}\\
x_1-\frac{x_1}{x_2^2}-\frac{x_3}{x_2}
\geq x_1-\frac{x_1}{9}-\frac{x_1}{3}=\frac{5x_1}{9}
\geq \frac{20}{9}.
\end{eqnarray*}

(ii) $r_1<r_4$ and $\mu_3(r_4)=\infty$ imply 
$r_3=\mu_3(r_1)\in(-\frac{x_3}{x_2},\infty)$.

(iii) holds because 
\begin{eqnarray*}
(x_2x_3-x_1)(x_1-1)-x_3^2>0\iff 
x_1-1>\frac{x_3^2}{x_2x_3-x_1}
=\frac{x_3}{x_2}\frac{1}{1-\frac{x_1}{x_2x_3}}\\
\textup{and }
\frac{x_3}{x_2}\frac{1}{1-\frac{x_1}{x_2x_3}}
\leq \frac{x_1}{3}\frac{1}{1-\frac{1}{2}}=\frac{2x_1}{3},\quad
x_1-1>\frac{2x_1}{3}\Leftarrow x_1\geq 4.
\end{eqnarray*}

(iv) holds because 
\begin{eqnarray*}
x_2^2(x_1-1)-(x_1+x_2x_3)>0\iff 
x_1-1>\frac{x_1+x_2x_3}{x_2^2}
=\frac{x_3}{x_2}(1+\frac{x_1}{x_2x_3})\\
\textup{and }
\frac{x_3}{x_2}(1+\frac{x_1}{x_2x_3})
\leq \frac{x_1}{3}(1+\frac{1}{2})=\frac{x_1}{2},\quad
x_1-1>\frac{x_1}{2}\Leftarrow x_1\geq 4.
\end{eqnarray*}

(ii)--(iv) show $r_3\in(-\frac{x_3}{x_2},0)$. Now
\begin{eqnarray*}
r_3\leq r_2&\iff &
\bigl[(x_2x_3-x_1)(x_1-1)-x_3^2\bigr](x_1-1)\\
&&-\bigl[x_2^2(x_1-1)-(x_1+x_2x_3)\bigr]\geq 0.
\end{eqnarray*}
The right hand side is
\begin{eqnarray*}
g(x_1,x_2,x_3):= (x_2x_3-x_1)(x_1-1)^2
-(x_2^2+x_3^2)(x_1-1)+(x_1+x_2x_3).
\end{eqnarray*}
This is symmetric and homogeneous of degree 2
in $x_2$ and $x_3$. 
\begin{eqnarray*}
\frac{\paa g}{\paa x_2}(\uuuu{x})
=(x_1-1)^2x_3-2(x_1-1)x_2+x_3,
\end{eqnarray*}
so for fixed $x_1$ and $x_3$ $g(\uuuu{x})$ 
takes its maximum in
\begin{eqnarray*}
x_2^0=\frac{(x_1-1)^2x_3+x_3}{2(x_1-1)}
=\frac{x_3(x_1-1)}{2}+\frac{x_3}{2(x_1-1)}
>\frac{3(x_1-1)}{2}>x_1,
\end{eqnarray*}
and is monotonously increasing left of $x_2^0$.

Because of the symmetry we can restrict to the cases
$x_1\geq x_2\geq x_3$. Then
$g(\uuuu{x})$ takes for fixed $x_1$ and $x_3$ its
minimum in $x_2=x_3$. 
\begin{eqnarray*}
\www{g}(x_1,x_3):= g(x_1,x_3,x_3)
=x_3^2(x_1-2)^2-(x_1-1)^2x_1+x_1.
\end{eqnarray*}
$\www{g}$ takes for fixed $x_1$ its minimum  at the
minimal possible $x_3$, which is
$x_3^0:=\max(3,\sqrt{2x_1})$ because of 
$2x_1\leq x_2x_3=x_3^2$. 
\begin{eqnarray*}
\www{g}(x_1,x_3^0)=\left\{\begin{array}{ll}
\www{g}(4,3)=4>0&\textup{ if }x_1=4,\\
\www{g}(x_1,\sqrt{2x_1})=x_1^3-6x_1^2+8x_1 & \\
\hspace*{1cm}=x_1(x_1-2)(x_1-4)>0 &\textup{ if }x_1\geq 5.
\end{array}\right.
\end{eqnarray*}
Therefore $r_3<r_2$ and \eqref{6.14} is proved.
\hfill$\Box$

\begin{remarks}\label{t6.19}
(i) In the proof of part (g) of Theorem \ref{t6.18}
the hyperbolic polygon $P$ whose relative boundary 
is the union of the six arcs
\begin{eqnarray*}
A(\infty,r_1),\ A(r_1,-\frac{x_3}{x_2}),\ 
A(-\frac{x_3}{x_2},r_3),\ A(r_2,0),\ A(0,1),\ 
A(1,\infty),
\end{eqnarray*}
is a fundamental domain for the action of the group
$\langle \mu_1,\mu_2,\mu_3\rangle$ on the upper half
plane $\H$.

(ii) In part (g) of Theorem \ref{t6.18} the case
$\uuuu{x}=(3,3,3)$ is especially interesting.
It is the only case within part (g) where $r_3=r_2$,
so the only case in part (g) where the hyperbolic 
polygon $P$ has finite hyperbolic area. 
In this case
\begin{eqnarray*}
\langle B_1,B_2,B_3\rangle =
\langle \begin{pmatrix}1&-3\\0&1\end{pmatrix},
\begin{pmatrix}1&0\\3&1\end{pmatrix},
\begin{pmatrix}-2&-3\\3&4\end{pmatrix}\rangle
=\Gamma(3).
\end{eqnarray*}
\end{remarks}

Now we turn to the study of the set $\Delta^{(1)}$
of odd vanishing cycles.
The shape of $\Delta^{(1)}$ and our knowledge on it
are very different for different $\uuuu{x}\in\Z^3$.

\begin{remarks}\label{t6.20}
(i) In the cases in the parts (c)--(f) of Theorem
\ref{t6.21} we will give $\Delta^{(1)}$ rather 
explicitly. There $\oooo{\Delta^{(1)}}$ is known,
and we can give a subset of $\Delta^{(1)}$ explicitly 
which maps bijectively to $\oooo{\Delta^{(1)}}$. 

Furthermore $\Gamma^{(1)}_u\cong\Z^2$, and for
$\varepsilon\in\{\pm 1\}$ and $i\in\{1,2,3\}$
there is a subset $F_{\varepsilon,i}\subset\Z f_3$ with
\begin{eqnarray*}
\Delta^{(1)}\cap (\varepsilon e_i+\Z f_3)
&=&\varepsilon e_i+ F_{\varepsilon,i}\quad\textup{and}\\
\Delta^{(1)}\cap (a+\Z f_3)
&=&a+ F_{\varepsilon,i}\quad\textup{for any }
a\in\Gamma^{(1)}\{\varepsilon e_i\}.
\end{eqnarray*}
In many, but not all cases $F_{\varepsilon,i}\subset\Z f_3$
is a sublattice.

(ii) On the contrary, in the cases in part (g) of
Theorem \ref{t6.21} we know rather little.
There $\Gamma^{(1)}_u=\{\id\}$, and remarkably the
projection $\Delta^{(1)}\to\oooo{\Delta^{(1)}}
\subset\oooo{H_\Z}^{(1)}$ is a bijection. 
In the case $\uuuu{x}=(3,3,3)$ 
we know $\oooo{\Delta^{(1)}}$. But the lift
$a\in\Delta^{(1)}$ of an element 
$\oooo{a}^{(1)}\in\oooo{\Delta^{(1)}}$ is difficult to
determine. See Lemma \ref{t6.26}. 
\end{remarks}

\begin{theorem}\label{t6.21}
(a) (Empty: (b)--(g) shall correspond to (b)--(g) in Theorem
\ref{t6.18}.)

(b) In each reducible case $\uuuu{x}=(x_1,0,0)$
\begin{eqnarray*}
\Delta^{(1)}&=&\Delta^{(1)}\cap (\Z e_1\oplus \Z e_2)
\ \dot\cup\ \{\pm e_3\}\\
\textup{with}&&
\Delta^{(1)}\cap(\Z e_1\oplus\Z e_2)\cong\Delta^{(1)}(S(-x_1)),
\end{eqnarray*}
and $\Delta^{(1)}(S(-x_1))$ is given in Theorem \ref{t6.10}.

(c) The case $\uuuu{x}=(1,1,0)$: 
\begin{eqnarray*}
\oooo{\Delta^{(1)}}&=&\oooo{H_\Z}^{(1),prim},\\
\Delta^{(1)} &=& (pr^{H,(1)})^{-1}(\oooo{\Delta^{(1)}})
=(\Z e_1\oplus\Z e_2)^{prim}+\Z f_3\subset H_\Z^{prim},\\
\Delta^{(1)}&=& \Gamma^{(1)}\{e_1\}.
\end{eqnarray*}
$\Delta^{(1)}$ and $\oooo{\Delta^{(1)}}$ consist each of 
one orbit. 

(d) The cases $\uuuu{x}=(x_1,x_2,0)$ with $2\leq x_1\geq x_2>0$:
Recall $x_{12}:=\gcd(x_1,x_2)$.

(i) The cases with $x_{12}=1$ and $x_2>1$: Then $x_1=\www{x}_1$, 
$x_2=\www{x}_2$ and $x_1>x_2>1$. Choose $y_1,y_2\in\Z$ with 
$1=y_1x_1^2+y_2x_2^2$. $\oooo{\Delta^{(1)}}$ consists of the 
three orbits
\begin{eqnarray*}
\Gamma^{(1)}_s\{\oooo{e_1}^{(1)}\}&=&\oooo{H_\Z}^{(1),prim},\\
\Gamma^{(1)}_s\{\oooo{e_2}^{(1)}\}
&=&x_1\oooo{H_\Z}^{(1),prim},\\
\Gamma^{(1)}_s\{\oooo{e_3}^{(1)}\}
&=&x_2\oooo{H_\Z}^{(1),prim}.
\end{eqnarray*}

$\Delta^{(1)}$ consists of five orbits:
\begin{eqnarray*}
\Gamma^{(1)}\{e_1\}=\Gamma^{(1)}\{\pm e_1\},\quad
\Gamma^{(1)}\{\varepsilon e_2\},\quad 
\Gamma^{(1)}\{\varepsilon e_3\}\quad\textup{with }
\varepsilon\in\{\pm 1\}.
\end{eqnarray*}
Here
\begin{eqnarray*}
&&\Delta^{(1)}\cap(e_1+\Z f_3)
= \Gamma^{(1)}\{e_1\}\cap (e_1+\Z f_3)\\
&&\hspace*{1cm}=\Gamma^{(1)}_u\{e_1\}
=e_1+\Z x_1x_2f_3,\\
&&\Delta^{(1)}\cap(\varepsilon e_2+\Z f_3)\\
&&\hspace*{1cm} = 
\Bigl(\Gamma^{(1)}\{\varepsilon e_2\}\cap (\varepsilon e_2+\Z f_3)
\Bigr) \ \dot\cup\ \Bigl( 
\Gamma^{(1)}\{-\varepsilon e_2\}\cap (\varepsilon e_2+\Z f_3)
\Bigr)\\
&&\hspace*{1cm} = \Bigl(\varepsilon e_2+\Z x_1^2x_2 f_3\Bigr)
\ \dot\cup\ 
\Bigl(\varepsilon e_2-\varepsilon 2y_2x_2 f_3+\Z x_1^2x_2f_3\Bigr)
,\\
&&\Delta^{(1)}\cap(\varepsilon e_3+\Z f_3)\\
&&\hspace*{1cm} = 
\Bigl(\Gamma^{(1)}\{\varepsilon e_3\}\cap (\varepsilon e_3+\Z f_3)
\Bigr)\ \dot\cup\ \Bigl(
\Gamma^{(1)}\{-\varepsilon e_3\}\cap (\varepsilon e_3+\Z f_3
\Bigr)\\
&&\hspace*{1cm} = \Bigl(\varepsilon e_3+\Z x_1x_2^2 f_3)\Bigr)
\ \dot\cup\ \Bigl(
\varepsilon e_3+\varepsilon 2y_1x_1 f_3+\Z x_1x_2^2f_3\Bigr).
\end{eqnarray*}

(ii) The cases with $x_2=1$: Then $x_1=\www{x}_1>x_2=\www{x}_2=x_{12}=1$.

$\oooo{\Delta^{(1)}}$ consists of the two orbits 
\begin{eqnarray*}
\Gamma^{(1)}_s\{\oooo{e_1}^{(1)}\}&=&
\Gamma^{(1)}_s\{\pm \oooo{e_1}^{(1)},\pm \oooo{e_3}^{(1)}\}
=\oooo{H_\Z}^{(1),prim},\\
\Gamma^{(1)}_s\{\oooo{e_2}^{(1)}\}
&=&\Gamma^{(1)}_s\{\pm \oooo{e_2}^{(1)}\}
=x_1\oooo{H_\Z}^{(1),prim}.
\end{eqnarray*}

$\Delta^{(1)}$ consists of three orbits,
\begin{eqnarray*}
\Gamma^{(1)}\{e_1\}=\Gamma^{(1)}\{\pm e_1,\pm e_3\},\quad
\Gamma^{(1)}\{\varepsilon e_2\},\quad \textup{with }
\varepsilon\in\{\pm 1\}.
\end{eqnarray*}
Here 
\begin{eqnarray*}
&&\Delta^{(1)}\cap(e_1+\Z f_3)
= \Gamma^{(1)}\{e_1\}\cap (e_1+\Z f_3)
=\Gamma^{(1)}_u\{e_1\}
=e_1+\Z x_1f_3,\\
&&\Delta^{(1)}\cap(\varepsilon e_2+\Z f_3)\\
&&\hspace*{1cm} = 
\Bigl(\Gamma^{(1)}\{\varepsilon e_2\}\cap (\varepsilon e_2+\Z f_3)
\Bigr) \ \dot\cup\ \Bigl( 
\Gamma^{(1)}\{-\varepsilon e_2\}\cap (\varepsilon e_2+\Z f_3)
\Bigr)\\
&&\hspace*{1cm} = \Bigl(\varepsilon e_2+\Z x_1^2 f_3\Bigr)
\ \dot\cup\ 
\Bigl(\varepsilon e_2-\varepsilon 2 f_3+\Z x_1^2f_3\Bigr).
\end{eqnarray*}

(iii) The cases with $x_{12}>1$ and $x_1>x_2>1$:
Then $\www{x}_1>\www{x}_2\geq 1$. 
Recall from Theorem \ref{t6.18} (s) 
$g_2=\frac{1}{\www{x}_1}\oooo{e_2}^{(1)}=
\frac{1}{\www{x}_2}\oooo{e_3}^{(1)}\in \oooo{H_\Z}^{(1)}$. 
Here $\oooo{\Delta^{(1)}}\subset\oooo{H_\Z}^{(1)}$ consists
of the six orbits (with $\varepsilon\in\{\pm 1\}$)
\begin{eqnarray*}
\Gamma^{(1)}_s\{\varepsilon \oooo{e_1}^{(1)}\}
&\subset&\oooo{H_\Z}^{(1),prim},\\
\Gamma^{(1)}_s\{\varepsilon \oooo{e_2}^{(1)}\}
&=&\www{x}_1\Gamma^{(1)}_s\{\varepsilon g_2\}
\subset \www{x}_1\oooo{H_\Z}^{(1),prim},\\
\Gamma^{(1)}_s\{\varepsilon \oooo{e_3}^{(1)}\}
&=&\www{x}_2\Gamma^{(1)}_s\{\varepsilon g_2\}
\subset \www{x}_2\oooo{H_\Z}^{(1),prim}
\end{eqnarray*}
Theorem \ref{t6.10} (c) and (d) describe the orbits
$\Gamma^{(1)}_s\{\oooo{e_1}^{(1)}\}$ and $\Gamma^{(1)}_s\{g_2\}$.
Also $\Delta^{(1)}$ consists of six orbits. 
\begin{eqnarray*}
\Delta^{(1)}\cap(\varepsilon e_1+\Z f_3)
&=& \varepsilon e_1 +\Z x_{12}\www{x}_1\www{x}_2f_3,\\
\Delta^{(1)}\cap(\varepsilon e_2+\Z f_3)
&=& \varepsilon e_2 +\Z x_1x_2\www{x}_1f_3,\\
\Delta^{(1)}\cap(\varepsilon e_3+\Z f_3)
&=& \varepsilon e_3 +\Z x_1x_2\www{x}_2f_3.
\end{eqnarray*}

(iv) The cases with $x_1=x_2\geq 2$: 
Then $x_{12}=x_1=x_2$, $\www{x}_1=\www{x}_2=1$,
$\oooo{e_2}^{(1)}=\oooo{e_3}^{(1)}=g_2$. Then
$\oooo{\Delta^{(1)}}\subset\oooo{H_\Z}^{(1)}$ consists
of the four orbits (with $\varepsilon\in\{\pm 1\}$)
\begin{eqnarray*}
\Gamma^{(1)}_s\{\varepsilon \oooo{e_1}^{(1)}\}
&\subset&\oooo{H_\Z}^{(1),prim},\\
\Gamma^{(1)}_s\{\varepsilon \oooo{e_2}^{(1)}\}
&=&\Gamma^{(1)}_s\{\varepsilon g_2\}\subset\oooo{H_\Z}^{(1),prim}.
\end{eqnarray*}
Theorem \ref{t6.10} (c) and (d) describe these orbits.
$\Delta^{(1)}$ consists of six orbits.
\begin{eqnarray*}
&&\Delta^{(1)}\cap(\varepsilon e_1+\Z f_3)
= \Gamma^{(1)}\{\varepsilon e_1\}\cap (\varepsilon e_1+\Z f_3)
= \Gamma^{(1)}_u\{\varepsilon e_1\}\\
&&\hspace*{1cm}= \varepsilon e_1 +\Z x_1f_3,\\
&&\Delta^{(1)}\cap(\varepsilon e_2+\Z f_3)
= \Delta^{(1)}\cap (\varepsilon e_3+\Z f_3)\\
&& \hspace*{1cm}= \Bigl(\Gamma^{(1)}\{\varepsilon e_2\}
\cap (\varepsilon e_2+\Z f_3)\Bigr)
\ \dot\cup\ 
\Bigl(\Gamma^{(1)}\{\varepsilon e_3\}
\cap (\varepsilon e_2+\Z f_3)\Bigr)\\
&& \hspace*{1cm}=
\Bigl(\varepsilon e_2 +\Z x_1^2f_3\Bigr)\ \dot\cup\ 
\Bigl(\varepsilon e_2 -\varepsilon f_3+\Z x_1^2f_3\Bigr).
\end{eqnarray*}

(e) The cases $\uuuu{x}=(-l,2,-l)$ with $l\geq 2$ even:
Consider $\varepsilon\in\{\pm 1\}$. 

(i) The cases with $l\equiv 0(4)$:
$\Delta^{(1)}$ consists of six orbits,
\begin{eqnarray*}
\Gamma^{(1)}\{\varepsilon e_1\}&=& 
\{y_1e_1+y_2e_3\in H_\Z^{prim}\,|\, 
y_1\equiv \varepsilon(4),y_2\equiv 0(2)\} + \Z lf_3,\\
\Gamma^{(1)}\{\varepsilon e_3\}&=& 
\{y_1e_1+y_2e_3\in H_\Z^{prim}\,|\, 
y_1\equiv 0(2),y_2\equiv \varepsilon(4)\} + \Z lf_3,\\
\Gamma^{(1)}\{\varepsilon e_2\}&=& 
\frac{l}{2}\{y_1e_1+y_2e_3\in H_\Z^{prim}\,|\, 
y_1\equiv 1(2),y_2\equiv 1(2)\}\\
&& + \varepsilon f_3 + \Z l^2f_3.
\end{eqnarray*}
$\oooo{\Delta^{(1)}}$ consists of five orbits,
$\Gamma_s^{(1)}\{\oooo{e_2}^{(1)}\}
=\Gamma_s^{(1)}\{-\oooo{e_2}^{(1)}\}$. 

(ii) The cases with $l\equiv 2(4)$:
$\Delta^{(1)}$ consists of six orbits,
\begin{eqnarray*}
&&\Gamma^{(1)}\{\varepsilon e_1\}= 
\Bigl(\{y_1e_1+y_2e_3\in H_\Z^{prim}\,|\, 
y_1\equiv \varepsilon(4),y_2\equiv 0(2)\} + \Z 2lf_3\Bigr)\\
&&\hspace*{0.5cm}
\dot\cup \Bigl(\{y_1e_1+y_2e_3\in H_\Z^{prim}\,|\, 
y_1\equiv -\varepsilon(4),y_2\equiv 0(2)\}+lf_3+\Z 2lf_3\Bigr),\\
&&\Gamma^{(1)}\{\varepsilon e_3\}= 
\Bigl(\{y_1e_1+y_2e_3\in H_\Z^{prim}\,|\, 
y_1\equiv 0(2),y_2\equiv \varepsilon(4)\} + \Z 2lf_3\Bigr)\\
&&\hspace*{0.5cm}
\dot\cup \Bigl(\{y_1e_1+y_2e_3\in H_\Z^{prim}\,|\, 
y_1\equiv 0(2),y_2\equiv -\varepsilon(4)\}+lf_3+\Z 2lf_3\Bigr),\\
&&\Gamma^{(1)}\{\varepsilon e_2\}= 
\frac{l}{2}\{y_1e_1+y_2e_3\in H_\Z^{prim}\,|\, 
y_1\equiv 1(2),y_2\equiv 1(2)\}\\
&&\hspace*{2cm}  + \varepsilon f_3 + \Z l^2f_3.
\end{eqnarray*}
Especially 
\begin{eqnarray*}
\Delta^{(1)}\cap (\varepsilon e_1+\Z f_3)
&=&\varepsilon e_1+\Z lf_3\\
\supsetneqq \Gamma^{(1)}\{\varepsilon e_1\}
\cap (\varepsilon e_1+\Z f_3)
&=&\varepsilon e_1+\Z 2l f_3,
\end{eqnarray*}
and similarly for $\varepsilon e_3$. 
$\oooo{\Delta^{(1)}}$ consists of three orbits,  
$\Gamma_s^{(1)}\{\oooo{e_i}^{(1)}\}
=\Gamma_s^{(1)}\{-\oooo{e_i}^{(1)}\}$ for 
$i\in\{1,2,3\}$.

(f) The cases $\uuuu{x}=(-l,2,-l)$ with $l\geq 3$ odd:
$\Delta^{(1)}$ consists of three orbits
(with $\varepsilon\in\{\pm 1\}$)
\begin{eqnarray*}
\Gamma^{(1)}\{e_1\} &=& \Gamma^{(1)}\{\pm e_1,\pm e_3\}
=(\Z e_1\oplus \Z g_2)^{prim}+\Z lf_3,\\
\Gamma^{(1)}\{\varepsilon e_2\}
&=&l(\Z e_1\oplus \Z g_2)^{prim}
+\varepsilon \frac{1-l^2}{2}f_3+\Z l^2 f_3.
\end{eqnarray*}
$\oooo{\Delta^{(1)}}$ consists of two orbits,
$\Gamma_s^{(1)}\{\oooo{e_2}^{(1)}\}
=\Gamma_s^{(1)}\{-\oooo{e_2}^{(1)}\}$.

(g) The cases $\uuuu{x}\in\Z^3_{\geq 3}$ with
$2x_i\leq x_jx_k$ for $\{i,j,k\}=\{1,2,3\}$:
$\Delta^{(1)}$ and $\oooo{\Delta^{(1)}}$ consist each of 
six orbits,
\begin{eqnarray*}
\Gamma^{(1)}\{\varepsilon e_i\}\quad\textup{respectively}\quad
\Gamma^{(1)}_s\{\varepsilon\oooo{e_i}^{(1)}\} 
\quad\textup{for }(\varepsilon,i)\in\{\pm 1\}\times\{1,2,3\}.
\end{eqnarray*}
The projection $\Delta^{(1)}\to\oooo{\Delta^{(1)}}$
is a bijection. 

(h) The case $\uuuu{x}=(3,3,3)$:
In this subcase of (g) the statements in (g) hold, and 
\begin{eqnarray*}
\Gamma^{(1)}\{\varepsilon e_i\}\subset \varepsilon e_i+3H_\Z,\quad
\Gamma^{(1)}_s\{\varepsilon\oooo{e_i}^{(1)}\} =
(\varepsilon \oooo{e_i}^{(1)}+3 \oooo{H_\Z}^{(1)})^{prim}
\end{eqnarray*}
for $\varepsilon\in\{\pm 1\}$, $i\in\{1,2,3\}$. 
\end{theorem}

{\bf Proof:}
(b) The splitting of $\Delta^{(1)}$ follows with Lemma \ref{t2.11}.
The first and second subset are treated by Theorem \ref{t6.10}
and Lemma \ref{t2.12}.

(c) $\oooo{\Delta^{(1)}}=\oooo{H_\Z}^{prim}$ follows with
$\Gamma^{(1)}_s\cong SL_2(\Z)$ and Theorem \ref{t6.10} (b). 
$\Delta^{(1)}$ is the full preimage in $H_\Z$, because
$\Gamma^{(1)}_u=O^{(1),Rad}_u=\{t^+_\lambda\,|\, 
\lambda(\uuuu{e})\in\langle (1,0,0),(0,1,1)\rangle_\Z\}$. 

(d) (i) All statements on $\oooo{\Delta^{(1)}}$ follow from
$\Gamma^{(1)}_s\cong SL_2(\Z)$ and 
$\oooo{e_2}^{(1)}=x_1g_2$, $\oooo{e_3}^{(1)}=x_2g_2$.

Recall the definition and the properties of $s_4\in\Gamma^{(1)}$
in the proof of Theorem \ref{t6.18}.

For $\Delta^{(1)}$ we have to use Lemma \ref{t6.7}.

\medskip
{\bf Claim 1:} For $a\in\{\pm e_1,\pm e_2,\pm e_3\}$ equality
in $\stackrel{(3)}{\supset}$ in Lemma \ref{t6.7} holds.

\medskip
{\sf Proof of Claim 1:}
Suppose $h\in \textup{Stab}_{\Gamma^{(1)}}(\oooo{a}^{(1)})$.
Then
\begin{eqnarray*}
\oooo{h}\in\textup{Stab}_{\Gamma^{(1)}_s}(\oooo{a}^{(1)})
=\left\{\begin{array}{ll}
\langle \oooo{s_{e_1}^{(1)}}\rangle
 & \textup{ if }a\in\{\pm e_1\},\\
\langle \oooo{s_4}\rangle 
 &\textup{ if }a\in\{\pm e_2,\pm e_3\},
\end{array}\right.
\end{eqnarray*}
so for a suitable $l\in\Z$ 
\begin{eqnarray*}
\oooo{h}(\oooo{s_{e_1}^{(1)}})^l=\id &\textup{and}& 
h(s_{e_1}^{(1)})^l\in \Gamma^{(1)}_u
\quad\textup{for }a\in\{\pm e_1\},\\
\oooo{h}(\oooo{s_4})^l=\id &\textup{and}& 
h(s_4)^l\in\Gamma^{(1)}_u
\quad\textup{for }a\in\{\pm e_2,\pm e_3\}.
\end{eqnarray*}
As $s_{e_1}^{(1)}\in\textup{Stab}_{\Gamma^{(1)}}
(\varepsilon e_1)$ and $s_4\in \textup{Stab}_{\Gamma^{(1)}}
(\varepsilon e_2)=\textup{Stab}_{\Gamma^{(1)}}
(\varepsilon e_3)$
$$h\in\Gamma^{(1)}_u\textup{Stab}_{\Gamma^{(1)}}(a).$$
\hfill($\Box$)

\medskip
By Lemma \ref{t6.7} (b) 
\begin{eqnarray}\label{6.15}
\Gamma^{(1)}\{\varepsilon e_i\}\cap (\varepsilon e_i+\Z f_3)
=\Gamma^{(1)}_u\{\varepsilon e_i\}.
\end{eqnarray}

\medskip
{\bf Claim 2:} $\bigl(s_{e_1}^{(1)}s_4 s_{e_1}^{(1)})^2
=t^-_{\lambda_3}$ with 
$\lambda_3\in\Hom_2(H_\Z,\Z)$, $\lambda_3(\uuuu{e})
=(0,2y_2x_2,-2y_1x_1)$. 

\medskip
{\sf Proof of Claim 2:}  This is a straightforward 
calculation with $s_{e_1}^{(1),mat}$ and 
$s_4^{mat}$.\hfill($\Box$)

\medskip
Therefore
\begin{eqnarray*}
\Gamma^{(1)}\{e_1\}&=&\Gamma^{(1)}\{-e_1\},\\
\Delta^{(1)}\cap (e_1+\Z f_3) 
&=&\Gamma^{(1)}\{e_1\}\cap (e_1+\Z f_3)
=\Gamma^{(1)}_u\{e_1\}\\
&=&e_1+\Z x_1x_2 f_3.
\end{eqnarray*}
Also 
\begin{eqnarray*}
\Gamma^{(1)}\{e_2\}\owns -e_2+2y_2x_2f_3,\quad 
\Gamma^{(1)}\{e_3\}\owns -e_3-2y_1x_1f_3.
\end{eqnarray*}
This together with \eqref{6.15} and the shape of $\Gamma^{(1)}_u$
shows the statements on 
$\Delta^{(1)}\cap (\varepsilon e_2+\Z f_3)$ and 
$\Delta^{(1)}\cap (\varepsilon e_3+\Z f_3)$.

(ii) All statements on $\oooo{\Delta^{(1)}}$ follow from
$\Gamma^{(1)}_s\cong SL_2(\Z)$ and $\oooo{e_2}^{(1)}=x_1g_2$
and $\oooo{e_3}^{(1)}=g_2$.

For $\Delta^{(1)}$ we have to use Lemma \ref{t6.7}. 
Claim 1 in the proof of part (i), its proof and the implication
\eqref{6.15} still hold. Now we can choose $(y_1,y_2)=(0,1)$,
so $s_4=s_{e_3}^{(1)}$. One calculates
\begin{eqnarray*}
s_{e_1}^{(1)}s_{e_3}^{(1)}s_{e_1}^{(1)}\, \uuuu{e}=\uuuu{e}
\begin{pmatrix}0&-x_1&-1\\0&1&0\\1&-x_1&0\end{pmatrix}.
\end{eqnarray*}
Therefore $\Gamma^{(1)}\{e_1\}=\Gamma^{(1)}\{\pm e_1,\pm e_3\}$.

Claim 2 in the proof of part (i) still holds. It gives
$(s_{e_1}^{(1)}s_{e_3}^{(1)}s_{e_1}^{(1)})^2=t^-_{\lambda_3}$
with $\lambda_3\in\Hom_2(H_\Z,\Z)$, 
$\lambda_3(\uuuu{e})=(0,2,0)$. 
Especially $t^-_{\lambda_3}(-e_2)=e_2-2f_3$. 

This fact, \eqref{6.15} and the shape of $\Gamma^{(1)}_u$
imply the statements on $\Delta^{(1)}\cap (e_1+\Z f_3)$
and $\Delta^{(1)}\cap (\varepsilon e_2+\Z f_3)$.

(iii) and (iv) 
By Theorem \ref{t6.18} (d) $\Gamma^{(1)}_s\cong 
\Gamma^{(1)}(S(-x_{12}))$. 
By Theorem \ref{t6.10} (c) and (d) 
$\Gamma^{(1)}_s\{\pm \oooo{e_1}^{(1)},\pm g_2\}$ consists
of the four orbits 
$\Gamma^{(1)}_s\{\varepsilon \oooo{e_1}^{(1)}\}$
and $\Gamma^{(1)}_s\{\varepsilon g_2\}$ with
$\varepsilon\in\{\pm 1\}$. 
Therefore if $x_1> x_2$ then $\oooo{\Delta^{(1)}}$
consists of the six orbits in (iii),
and if $x_1=x_2$ then $\oooo{\Delta^{(1)}}$ consists of
the four orbits $\Gamma^{(1)}_s\{\varepsilon \oooo{e_1}^{(1)}\}$
and $\Gamma^{(1)}_s\{\varepsilon \oooo{e_2}^{(1)}\}
=\Gamma^{(1)}_s\{\varepsilon \oooo{e_3}^{(1)}\}
=\Gamma^{(1)}_s\{\varepsilon g_2\}$.

For $\Delta^{(1)}$ we have to use Lemma \ref{t6.7}.
Claim 1 in the proof of part (i), its proof and the implication
\eqref{6.15} still hold.

In the case $x_1>x_2$, $\oooo{\Delta^{(1)}}$ consists of six
orbits. Part (a) of Lemma \ref{t6.7}, \eqref{6.15} and the
shape of $\Gamma^{(1)}_u$ imply the statements on
$\Delta^{(1)}\cap (\varepsilon e_i+\Z f_3)$ in part (iii).

In the case $x_1=x_2$, $\oooo{\Delta^{(1)}}$ consists of four
orbits. Then $e_3=e_2-f_3$, \eqref{6.15} and the shape of
$\Gamma^{(1)}_u$ imply the statements on 
$\Delta^{(1)}\cap (\varepsilon e_i+\Z f_3)$ in part (iv).

(e) Theorem \ref{t6.10} (c) and
\begin{eqnarray*}
s_{e_1}^{(1)}(e_1,e_3)=(e_1,e_3)
\begin{pmatrix}1&-2\\0&1\end{pmatrix},\quad 
s_{e_3}^{(1)}(e_1,e_3)=(e_1,e_3)
\begin{pmatrix}1&0\\2&1\end{pmatrix},
\end{eqnarray*}
imply
\begin{eqnarray*}
\langle s_{e_1}^{(1)},s_{e_3}^{(1)}\rangle
\{\varepsilon e_1\} &=& 
\{y_1e_1+y_2e_3\in H_\Z^{prim}\,|\, y_1\equiv\varepsilon (4),
y_2\equiv 0(2)\},\\
\langle s_{e_1}^{(1)},s_{e_3}^{(1)}\rangle
\{\varepsilon e_3\} &=& 
\{y_1e_1+y_2e_3\in H_\Z^{prim}\,|\, y_1\equiv 0(2)),
y_2\equiv \varepsilon(4)\}.
\end{eqnarray*}

In the case (i), the semidirect product 
$\Gamma^{(1)}=\Gamma^{(1)}_u\rtimes \langle s_{e_1}^{(1)},
s_{e_3}^{(1)}\rangle $ and the shape of $\Gamma^{(1)}_u$
show that $\Gamma^{(1)}\{\varepsilon e_1\}$ and 
$\Gamma^{(1)}\{\varepsilon e_3\}$ are as claimed.

In the case (ii), the semidirect product 
$\Gamma^{(1)}=\Gamma^{(1)}\cap O^{(1),Rad}_\pm 
\rtimes \langle s_{e_1}^{(1)},
s_{e_3}^{(1)}\rangle $ and the shape of $\Gamma^{(1)}
\cap O^{(1),Rad}_\pm$ 
show that $\Gamma^{(1)}\{\varepsilon e_1\}$ and 
$\Gamma^{(1)}\{\varepsilon e_3\}$ are as claimed.

The following fact was not mentioned in Theorem \ref{t6.10} (c):
\begin{eqnarray*}
\{y_1e_1+y_2e_3\in H_\Z ^{prim}\,|\, y_1\equiv y_2\equiv 1(2)\}
=\langle s_{e_1}^{(1)},s_{e_3}^{(1)}\rangle \{e_1+e_3\},
\end{eqnarray*}
so this set is a single 
$\langle s_{e_1}^{(1)},s_{e_3}^{(1)}\rangle$ orbit.
We skip its proof
(it contains the observation $s_{e_3}^{(1)}s_{e_1}^{(1)}
(e_1+e_3)=-e_1-e_3$).

This fact, the semidirect products above of $\Gamma^{(1)}$,
the shape of $\Gamma^{(1)}_u$ in case (i) and of 
$\Gamma^{(1)}\cap O^{(1),Rad}_\pm$ in case (ii),
and $e_2=-\frac{l}{2}(e_1+e_3)+f_3$ show in case (i) and case (ii)
that $\Gamma^{(1)}\{\varepsilon e_2\}$ is as claimed.

(f) Theorem \ref{t6.10} (b), $\www{\uuuu{e}}=(e_1, g_2,f_3)$,
\begin{eqnarray*}
s_{e_1}^{(1)}(\www{\uuuu{e}})&=&\www{\uuuu{e}}
\begin{pmatrix}1&-1&0\\0&1&0\\0&0&1\end{pmatrix},\quad
s_4(\www{\uuuu{e}})=\www{\uuuu{e}}
\begin{pmatrix}1&0&0\\1&1&0\\0&0&1\end{pmatrix},\\
e_3&=&-e_1+2 g_2+lf_3\quad\textup{and}\quad 
e_2=-l g_2+\frac{1-l^2}{2}f_3
\end{eqnarray*}
imply
\begin{eqnarray*}
\langle s_{e_1}^{(1)},s_4\rangle (e_1) &=& 
(\Z e_1\oplus \Z  g_2)^{prim},\\
\langle s_{e_1}^{(1)},s_4\rangle (e_3) &=& 
(\Z e_1\oplus \Z  g_2)^{prim}+lf_3,\\
\langle s_{e_1}^{(1)},s_4\rangle (e_2) &=& 
l(\Z e_1\oplus \Z  g_2)^{prim}+\frac{1-l^2}{2}f_3.
\end{eqnarray*}
The semidirect product 
$\Gamma^{(1)}=\Gamma^{(1)}_u\rtimes 
\langle s_{e_1}^{(1)},s_4\rangle$ and the shape of 
$\Gamma^{(1)}_u$ show  (with $\varepsilon\in\{\pm 1\}$)
\begin{eqnarray*}
\Gamma^{(1)}\{e_1\}&=& \Gamma^{(1)}\{\pm e_1,\pm e_3\}
=(\Z e_1\oplus \Z  g_2)^{prim}+\Z l f_3,\\
\Gamma^{(1)}\{\varepsilon e_2\}&=& 
l(\Z e_1\oplus \Z  g_2)^{prim}+\varepsilon\frac{1-l^2}{2}f_3
+ \Z l^2 f_3.
\end{eqnarray*}

(g) In the proof of part (g) of Theorem
\ref{t6.18} the hyperbolic polygon $P$ with the six arcs
in Remark \ref{t6.19} was used. It is a fundamental polygon
of the action of the group $\langle \mu_1,\mu_2,\mu_3\rangle$
on $\H$. Here $\mu_1,\mu_2,\mu_3$ are parabolic 
M\"obius transformations with fixed points $\infty,0$ and
$-\frac{x_3}{x_2}$. These fixed points are cusps of $P$.
This geometry implies
\begin{eqnarray*}
\textup{Stab}_{\langle \mu_1,\mu_2,\mu_3\rangle}(\infty)
&=&\langle\mu_1\rangle,\\
\textup{Stab}_{\langle \mu_1,\mu_2,\mu_3\rangle}(0)
&=&\langle\mu_2\rangle,\\
\textup{Stab}_{\langle \mu_1,\mu_2,\mu_3\rangle}(-\frac{x_3}{x_2})
&=&\langle\mu_3\rangle.
\end{eqnarray*}
As $\langle \mu_1,\mu_2,\mu_3\rangle\cong\Gamma^{(1)}_s
\bigl(\cong G^{free,3}\bigr)$ with 
$\mu_i\sim \oooo{s_{e_i}^{(1)}}$, this implies
\begin{eqnarray*}
\textup{Stab}_{\Gamma^{(1)}_s}(\{\pm \oooo{e_i}^{(1)}\}) 
= \langle \oooo{s_{e_i}^{(1)}}\rangle
=\textup{Stab}_{\Gamma^{(1)}_s}(\oooo{e_i}^{(1)}) 
\quad \textup{for}\quad i\in\{1,2,3\}.
\end{eqnarray*}
Especially 
$-\oooo{e_i}^{(1)}\notin \Gamma^{(1)}_s(\oooo{e_i}^{(1)})$, so the
orbits $\Gamma^{(1)}_s\{\oooo{e_i}^{(1)}\}$ and 
$\Gamma^{(1)}_s\{-\oooo{e_i}^{(1)}\}$ are disjoint. 

The cusps $\infty$, $0$ and $-\frac{x_3}{x_2}$ of the 
fundamental domain $P$ of the group 
$\langle \mu_1,\mu_2,\mu_3\rangle$ are in disjoint orbits
of $\langle \mu_1,\mu_2,\mu_3\rangle$.
Therefore the sets 
$\Gamma^{(1)}_s\{\pm \oooo{e_1}^{(1)}\}$, 
$\Gamma^{(1)}_s\{\pm \oooo{e_2}^{(1)}\}$ and 
$\Gamma^{(1)}_s\{\pm \oooo{e_3}^{(1)}\}$ are disjoint.
Therefore $\oooo{\Delta^{(1)}}$ consists of 
the six disjoint orbits 
$\Gamma^{(1)}_s\{\varepsilon \oooo{e_i}^{(1)}\}$ 
with $(\varepsilon,i)\in\{\pm 1\}\times \{1,2,3\}$. 
Therefore also $\Delta^{(1)}$ consists of the six orbits
$\Gamma^{(1)}\{\varepsilon e_i\}$ 
with $(\varepsilon,i)\in\{\pm 1\}\times \{1,2,3\}$. 

\medskip
{\bf Claim:} For $a\in\{\pm e_1,\pm e_2,\pm e_3\}$ equality in
$\stackrel{(1)}{\supset}$ and $\stackrel{(2)}{\supset}$
before Lemma \ref{t6.7} holds.

\medskip
{\bf Proof of the Claim:}
Equality in $\stackrel{(1)}{\supset}$ holds because of 
Lemma \ref{t6.7} (a) and because $\Delta^{(1)}$ and 
$\oooo{\Delta^{(1)}}$ each consist of six orbits. 

Equality in $\stackrel{(2)}{\supset}$ is by Lemma \ref{t6.7} (b)
equivalent to equality in $\stackrel{(3)}{\supset}$. 

Here $\Gamma^{(1)}\cong \Gamma^{(1)}_s(\cong G^{free,3})$,
$\Gamma^{(1)}_u=\{\id\}$, the lift $s_{e_i}^{(1)}\in \Gamma^{(1)}$
of $\oooo{s_{e_i}^{(1)}}\in \Gamma^{(1)}_s$ is in 
$\Stab_{\Gamma^{(1)}}(\varepsilon e_i)$, and therefore
\begin{eqnarray*}
\Stab_{\Gamma^{(1)}}(\varepsilon\oooo{e_i}^{(1)})
=\langle s_{e_i}^{(1)}\rangle 
=\Stab_{\Gamma^{(1)}}(\varepsilon e_i)
=\Gamma^{(1)}_u\cdot \Stab_{\Gamma^{(1)}}(\varepsilon e_i),
\end{eqnarray*}
so equality in $\stackrel{(3)}{\supset}$ holds.
The Claim is proved. \hfill ($\Box$)

\medskip
The Claim and $\Gamma^{(1)}_u=\{\id\}$ show
\begin{eqnarray*}
(\varepsilon e_i+\Z f_3)\cap \Delta^{(1)}
=\varepsilon e_i.
\end{eqnarray*}
Therefore the projection $\Delta^{(1)}\to\oooo{\Delta^{(1)}}$ 
is a bijection.

(h) In the case $\uuuu{x}=(3,3,3)$
\begin{eqnarray*}
\www{\uuuu{x}}=(1,1,1),\ f_3=-e_1+e_2-e_3,\ 
\oooo{e_3}^{(1)}=-\oooo{e_1}^{(1)}+\oooo{e_2}^{(1)},\\
\oooo{H_\Z}^{(1)}=\Z\oooo{e_1}^{(1)}\oplus 
\Z\oooo{e_2}^{(1)}.
\end{eqnarray*}
Recall that the matrices $B_1,B_2,B_3$ in the proof of Theorem
\ref{t6.18} (g) with 
$\oooo{s_{e_i}^{(1)}}(\oooo{e_1}^{(1)},\oooo{e_2}^{(1)})
=(\oooo{e_1}^{(1)},\oooo{e_2}^{(1)})B_i$ are here 
\begin{eqnarray*}
B_1=\begin{pmatrix}1&-3\\0&1\end{pmatrix},\ 
B_1=\begin{pmatrix}1&0\\3&1\end{pmatrix},\ 
B_3=\begin{pmatrix}-2&-3\\3&4\end{pmatrix},
\end{eqnarray*}
and generate $\Gamma(3)$. 
Because $s_{e_i}^{(1),mat}\equiv E_3\mmod 3$ and
$B_i\equiv E_2\mmod 3$,
\begin{eqnarray*}
\Gamma^{(1)}\{\varepsilon e_i\}\subset \varepsilon e_i +3H_\Z
\quad\textup{and}\quad
\Gamma^{(1)}_s\{\varepsilon \oooo{e_i}^{(1)}\}\subset 
\varepsilon\oooo{e_i}^{(1)}+3\oooo{H_\Z}^{(1)}.
\end{eqnarray*}

It remains to show
$\Gamma^{(1)}_s\{\oooo{e_i}^{(1)}\}=
(\oooo{e_i}^{(1)}+3\oooo{H_\Z}^{(1)})^{prim}$. 
This is equivalent to the following three statements:
\begin{list}{}{}
\item[(i)]
For $(a_1,a_3)\in\Z^2$ with $(a_1,a_3)\equiv (1,0)\mmod 3$ and 
$\gcd(a_1,a_3)=1$ a pair $(a_2,a_4)\in\Z^2$ with 
$\begin{pmatrix}a_1&a_2\\a_3&a_4\end{pmatrix}\in\Gamma(3)$
exists.
\item[(ii)]
For $(a_2,a_4)\in\Z^2$ with $(a_2,a_4)\equiv (0,1)\mmod 3$ and 
$\gcd(a_2,a_4)=1$ a pair $(a_1,a_3)\in\Z^2$ with 
$\begin{pmatrix}a_1&a_2\\a_3&a_4\end{pmatrix}\in\Gamma(3)$
exists.
\item[(iii)]
For $(b_1,b_2)\in\Z^2$ with $(b_1,b_2)\equiv (-1,1)\mmod 3$ and 
$\gcd(b_1,b_2)=1$ a matrix 
$\begin{pmatrix}a_1&a_2\\a_3&a_4\end{pmatrix}\in\Gamma(3)$ with 
$\begin{pmatrix}-a_1+a_2\\-a_3+a_4\end{pmatrix}
=\begin{pmatrix}b_1\\b_2\end{pmatrix}$ exists.
\end{list}
(i) is proved as follows. There exist $(\www{a_2},\www{a_4})\in\Z^2$
with $1=a_1\www{a_4}-a_3\www{a_2}$. Thus $\www{a_4}\equiv 1(3)$.
Let $\www{a_2}\equiv r(3)$ with $r\in\{0,1,2\}$. 
Choose $(a_2,a_4):=(\www{a_2}-ra_1,\www{a_4}-ra_3)$. 
The proofs of (ii) and (iii) are similar. 

The proof of Theorem \ref{t6.21} (h) is finished. 
\hfill$\Box$

\bigskip
In quite some cases $\Delta^{(0)}\subset \Delta^{(1)}$, but 
nevertheless in general $\Delta^{(0)}\not\subset\Delta^{(1)}$.
Corollary \ref{t6.22} gives some details.

\begin{corollary}\label{t6.22}
(a) $\Delta^{(0)}=\Delta^{(1)}$ holds only in the cases $A_1^n$, 
so the cases with $S=E_n$ for some $n\in\N$.
In all other cases $\Delta^{(1)}\not\subset R^{(0)}$. 

(b) $\Delta^{(0)}\subsetneqq\Delta^{(1)}$
in the cases $n=2$ except $A_1^2$, 
in the reducible cases with
$n=3$ except $A_1^3$ and in the case $A_3$.

(c) $\Delta^{(0)}\not\subset\Delta^{(1)}$ holds in the following
cases with $n=3$: $\whh{A}_2$, $\HH_{1,2}$, $S(-l,2,-l)$ with $l\geq 3$
and $\P^2$. 
\end{corollary}

{\bf Proof:}
(a) In the cases $A_1^n$ $\Delta^{(0)}=\Delta^{(1)}=
\{\pm e_1,...,\pm e_n\}$ by Lemma \ref{t2.12}.
In a case with $S\in T^{uni}_n(\Z)-\{E_n\}$ there is an entry
$S_{ij}\neq 0$ for some $i<j$, so $L(e_j,e_i)\neq 0$.
We can restrict to the rank 2 unimodular bilinear lattice
$(\Z e_i+\Z e_j,L|_{\Z e_i+\Z e_j},(e_i,e_j))$ with triangular
basis $(e_i,e_j)$. Part (e) of Theorem \ref{t6.10} shows that
it has odd vanishing cycles which are not roots.
They are also odd vanishing cycles of $(H_\Z,L,\uuuu{e})$,
so then $\Delta^{(1)}\not\subset R^{(0)}\supset \Delta^{(0)}$.

(b) For the cases $n=2$ see Theorem \ref{t6.10} (e).
The reducible cases with $n=3$ follow from the case $A_1$ and the
cases with $n=2$. In the case $A_3$ the twelve elements of $\Delta^{(0)}$
are given in Theorem \ref{t6.14} (c). The set $\Delta^{(1)}$ is by
Theorem \ref{t6.21} (c) 
\begin{eqnarray*}
(\pr^{H,(1)})^{-1}(\oooo{H_\Z}^{(1),prim})+\Z f_3.
\end{eqnarray*}
which contains $\Delta^{(0)}$ as a strict subset.

(c) The case $\whh{A}_2$: 
With $\uuuu{x}=(-1,-1,-1)$ we have $f_1=e_1+e_2+e_3$ and 
$f_3=e_1-e_2+e_3$. 
By Theorem \ref{t6.14} (d) 
\begin{eqnarray*}
\Delta^{(0)}=(\pm e_1+\Z f_1)\,\dot\cup\, (\pm e_2+\Z f_1)
\, \dot\cup\, (\pm e_3+\Z f_1).
\end{eqnarray*}
By Theorem \ref{t6.21} (c) 
\begin{eqnarray*}
\Delta^{(1)}= (\pr^{H,(1)})^{-1}(\oooo{H_\Z}^{(1),prim})+\Z f_3.
\end{eqnarray*}
Here for example for $m\in\Z-\{0,-1\}$
$$\Delta^{(0)}\owns e_2+mf_1=(2m+1)e_2+mf_3\notin\Delta^{(1)}.$$

The case $\HH_{1,2}$: With $\uuuu{x}=(-2,2,-2)$ we have
$\Rad I^{(0)}=\Z(e_1+e_2)\oplus \Z(e_2+e_3)$ and
$f_3=e_1+e_2+e_3$. By Theorem \ref{t6.14} (e)
$$\Delta^{(0)}=(\pm e_1+2\Rad I^{(0)})\,\dot\cup\, 
(\pm e_2+2\Rad I^{(0)})\,\dot\cup\, 
(\pm e_3+2\Rad I^{(0)}).$$
By Theorem \ref{t6.21} (e) (ii) 
$$\Delta^{(1)}\subset (\Z e_1+\Z e_3)^{prim}+ \Z f_3.$$
Here for example
\begin{eqnarray*}
\Delta^{(0)}&\owns& e_1+6(e_1+e_2)+4(e_2+e_3)\\
&=&(-3)e_1+(-6)e_3+10f_3\notin (\Z e_1+\Z e_3)^{prim}+ \Z f_3.
\end{eqnarray*}

This element is not contained in $\Delta^{(1)}$ 
because by Theorem \ref{t6.21} (e) 
$$\Delta^{(1)}\subset \Bigl((\Z e_1+\Z e_3)^{prim}+\Z f_3\Bigr)
\,\dot\cup\, 
\Bigl(\frac{l}{2}(\Z e_1+\Z e_3)^{prim}+\Z f_3\Bigr).$$

The cases $S(-l,2,-l)$: Recall
$\Rad I^{(0)}=\Z f_1$, 
$f_1=e_1-e_3$, 
\begin{eqnarray*}
&&f_3=\frac{1}{2}(le_1+2e_2+le_3)\quad\textup{if }l\geq 3\textup{ is even},\\
&&\left. \begin{array}{l}
f_3=le_1+2e_2+le_3\\
g_2=\frac{1}{2}(e_1+e_3)-\frac{l}{2}e_2\end{array}\right\} 
\textup{ if }l\geq 3\textup{ is odd.}
\end{eqnarray*}
Consider the element 
\begin{eqnarray*}
h_1(e_1)-la_1f_1&=&h_1(e_1-la_1f_1)\in H_\Z\quad\textup{with}\\
h_1&:=& (s_{e_1}^{(0)}s_{e_2}^{(0)})^3\in\Gamma^{(0)},\\
a_1&:=& \frac{1}{2}(l^5-4l^3+3l)
\in\Z
\end{eqnarray*}
By Theorem \ref{t6.11} (f) 
$T(\oooo{j}^{(0)}(\oooo{e_1}^{(0)})\otimes f_1)$ and 
$T(\oooo{j}^{(0)}(l\oooo{e_2}^{(0)})\otimes f_1)\in \Gamma_u^{(0)}$ with 
\begin{eqnarray*}
T(\oooo{j}^{(0)}(\oooo{e_1}^{(0)})\otimes f_1)(e_1)&=&e_1+2f_1,\\
T(\oooo{j}^{(0)}(l\oooo{e_2}^{(0)})\otimes f_1)(e_1)&=&e_1-l^2f_1,
\end{eqnarray*}
so
\begin{eqnarray*}
\Delta^{(0)}&\supset&\left\{\begin{array}{ll}
e_1+2\Z f_1 & \textup{ if }l\textup{ is even,}\\
e_1+\Z f_1 & \textup{ if }\textup{ is odd.} \end{array}\right.
\end{eqnarray*}
Therefore $h_1(e_1)-la_1f_1=h_1(e_1-la_1f_1)\in\Delta^{(0)}$. One calculates  
\begin{eqnarray*}
h_1(e_1)-la_1f_1&=& (s_{e_1}^{(0)}s_{e_2}^{(0)})^3 (e_1)-la_1f_1\\
&=&\uuuu{e}(\begin{pmatrix}-1&l&2\\0&1&0\\0&0&1\end{pmatrix}
\begin{pmatrix}1&0&0\\l&-1&l\\0&0&1\end{pmatrix})^2
\begin{pmatrix}1\\0\\0\end{pmatrix}-la_1f_1\\
&=&\uuuu{e}\begin{pmatrix}l^2-1&-l&l^2-2\\l&-1&l\\0&0&1\end{pmatrix}^3
\begin{pmatrix}1\\0\\0\end{pmatrix}-la_1f_1\\
&=&\uuuu{e}\begin{pmatrix}l^6-5l^4+6l^2-1\\l^5-4l^3+3l\\0\end{pmatrix}
-la_1f_1\\
&=&(-l^4+3l^2-1) e_1+
\left\{\begin{array}{ll}
2a_1f_3&\textup{ if }l\textup{ is even,}\\
a_1f_3&\textup{ if }l\textup{ is odd.}\end{array}\right.
\end{eqnarray*}
This element is not contained in $\Delta^{(1)}$ 
because $(-l^4+3l^2-1)\notin\{\pm 1,\pm \frac{l}{2},\pm l\}$ 
and because by Theorem \ref{t6.21} (e) and (f) 
\begin{eqnarray*}
\Delta^{(1)}&\subset& \Bigl((\Z e_1+\Z e_3)^{prim}+\Z f_3\Bigr)
\,\dot\cup\, 
\Bigl(\frac{l}{2}(\Z e_1+\Z e_3)^{prim}+\Z f_3\Bigr)\\
&& \textup{ if }
l\textup{ is even,}\\
\Delta^{(1)}&\subset& \Bigl((\Z e_1+\Z g_2)^{prim}+\Z f_3\Bigr)
\,\dot\cup\, 
\Bigl(l(\Z e_1+\Z g_2)^{prim}+\Z f_3\Bigr)\\
&& \textup{ if }l\textup{ is odd.}
\end{eqnarray*}

The case $\P^2$: With $\uuuu{x}=(-3,3,-3)$ we have
$f_3=e_1+e_2+e_3$. 
\begin{eqnarray*}
\Delta^{(1)}&\owns& s_{e_3}^{(1)}(s_{e_2}^{(1)})^{-2}(e_1)\\
&=&\uuuu{e}\begin{pmatrix}1&0&0\\0&1&0\\3&-3&1\end{pmatrix}
\begin{pmatrix}1&0&0\\3&1&-3\\0&0&1\end{pmatrix}
\begin{pmatrix}1&0&0\\3&1&-3\\0&0&1\end{pmatrix}
\begin{pmatrix}1\\0\\0\end{pmatrix}
=\uuuu{e}\begin{pmatrix}1\\6\\-15\end{pmatrix}.
\end{eqnarray*}
By Theorem \ref{t6.21} (g) the projection 
$\Delta^{(1)}\to\oooo{\Delta^{(1)}}$ is a bijection.
Therefore 
\begin{eqnarray*}
\Delta^{(1)}\not\owns (e_1+6e_2-15e_3)+9f_3
=10e_1+15e_2-6e_3.
\end{eqnarray*}
On the other hand
\begin{eqnarray*}
&&L(10e_1+15e_2-6e_3,10e_1+15e_2-6e_3)\\
&=&\begin{pmatrix}10&15&-6\end{pmatrix}
\begin{pmatrix}1&0&0\\-3&1&0\\3&-3&1\end{pmatrix}
\begin{pmatrix}10\\15\\-6\end{pmatrix} =1,
\end{eqnarray*}
so $10e_1+15e_2-6e_3\in R^{(0)}=\Delta^{(0)}$. \hfill$\Box$

\bigskip
Consider a triple $\uuuu{x}\in\Z^3$ and a corresponding
unimodular bilinear lattice $(H_\Z,L,\uuuu{e})$ with 
a triangular basis $\uuuu{e}$ with $L(\uuuu{e}^t,\uuuu{e})^t
=S(\uuuu{x})$. The Remarks \ref{t4.17} explained that the
tuple $(H_\Z,\pm I^{(1)},\Gamma^{(1)},\Delta^{(1)})$ depends
only on the $(G^{phi}\ltimes \www{G}^{sign})\rtimes
\langle\gamma\rangle$ orbit of $\uuuu{x}\in\Z^3$.
Lemma \ref{t4.18} gave at least one element of each orbit
of this group in $\Z^3$.
Theorem \ref{t6.18} and Theorem \ref{t6.21} gave detailed
information on the tuple $(H_\Z,\pm I^{(1)},\Gamma^{(1)},
\Delta^{(1)})$ for the elements in Lemma \ref{t4.18} (b)+(c)
and rather coarse information for the elements in
Lemma \ref{t4.18} (a). 

The next corollary uses this information to conclude that the 
$(G^{phi}\ltimes \www{G}^{sign})\rtimes\langle\gamma\rangle$ 
orbits of the elements in Lemma \ref{t4.18} (b)+(c) are
pairwise different and also different from the orbits of
the elements in Lemma \ref{t4.18} (a), because the 
corresponding tuples 
$(H_\Z,\pm I^{(1)},\Gamma^{(1)},\Delta^{(1)})$
are not isomorphic. As Theorem \ref{t6.18} and Theorem
\ref{t6.21} give only coarse information on the 
cases in Lemma \ref{t4.18} (a), also Corollary \ref{t6.23}
is vague about them.

The set of local minima in Lemma \ref{t4.18} (b)+(c) and
$(3,3,3)$ is called $\Lambda_1$, the set of local minima in
Lemma \ref{t4.18} (a) without $(3,3,3)$ is called $\Lambda_2$,
\begin{eqnarray*}
\Lambda_1&:=& \{(3,3,3)\}\cup\{(-l,2,-l)\,|\, l\geq 2\}\\
&&\cup\{(x_1,x_2,0)\,|\, x_1,x_2\in\Z_{\geq 0},x_1\geq x_2\},\\
\Lambda_2&:=&\{\uuuu{x}\in\Z^3_{\geq 3}\,|\, 2x_i\leq x_jx_k
\textup{ for }\{i,j,k\}=\{1,2,3\}\}-\{(3,3,3)\}.
\end{eqnarray*}

\begin{corollary}\label{t6.23}
Consider $\uuuu{x}$ and $\www{\uuuu{x}}\in \Lambda_1$ or 
$\uuuu{x}\in \Lambda_1$ and $\www{\uuuu{x}}\in \Lambda_2$.
Suppose $\uuuu{x}\neq \www{\uuuu{x}}$. 
Then the tuples 
$(H_\Z,\pm I^{(1)},\Gamma^{(1)},\Delta^{(1)})$ of
$\uuuu{x}$ and $\www{\uuuu{x}}$ are not isomorphic.
Consequently, the 
$(G^{phi}\ltimes \www{G}^{sign})\rtimes\langle\gamma\rangle$ 
orbits of $\uuuu{x}$ and $\www{\uuuu{x}}$ are disjoint.
\end{corollary}

{\bf Proof:}
In the following, (b), (c), (d), (d)(i), (d)(ii), (d)(iii), 
(d)(iv), (e), (e)(i), (e)(ii), (f), (g), (h)($\subset$ (g))
mean the corresponding families of cases in Theorem \ref{t6.21}.
Of course, (c) and (h) are single cases. 
We will first discuss how to separate the families by properties
of the isomorphism classes of the tuples 
$(H_\Z,\pm I^{(1)},\Gamma^{(1)},\Delta^{(1)})$, and then
how to separate the cases within one family.

The pair $(\Gamma^{(1)}_u,\Gamma^{(1)}_s)$ gives  the following
incomplete separation of families,
\begin{eqnarray*}
\begin{array}{l|c|r}
\Gamma^{(1)}_u\cong \ ? &\Gamma^{(1)}_s\cong \ ? & 
\textup{families}\\ \hline 
\{\id\}& G^{free,3} & (g) \\
\{\id\}&\not\cong G^{free,3} & (b) \\
\Z^2&  SL_2(\Z) & (c),(d)(i)+(ii),(f)\\
\Z^2&  SL_2(\Z)\times\{\pm 1\} & (e)(ii) \\
\Z^2&  G^{free,2} & (d)(iii)+(iv), (e)(i)
\end{array}
\end{eqnarray*}
The fundamental polygon $P$ in Remark \ref{t6.19} (ii) 
has finite area in the case (h) (i.e. 
$\uuuu{x}=(3,3,3)$) and infinite area in the other
cases in (g). So it separates the case (h) from the other
cases in (g).

The number of $\Br_3\ltimes\{\pm 1\}^3$ orbits
in $\oooo{\Delta^{(1)}}$ separates the families
(c),(d),(e)(i),(f) almost completely:
\begin{eqnarray*}
\begin{array}{l|c|c|c|c|c|c}
|\{\textup{orbits in }\oooo{\Delta^{(1)}}\}| & 
1 & 2 & 3 & 4 & 5 & 6 \\
\textup{families} & 
(c) & (d)(ii), (f) & (d)(i) & (d)(iv) & (e)(i) & (d)(iii)
\end{array}
\end{eqnarray*}

The separation of the family (d)(ii) from the family (f) is
more difficult and can be done as follows. In both families
of cases $\Delta^{(1)}$ consists of three orbits, and 
$\oooo{\Delta^{(1)}}$ consists of two orbits. The two
orbits $\Gamma^{(1)}\{e_2\}$ and $\Gamma^{(1)}\{-e_2\}$
unite to a single orbit
$\Gamma^{(1)}_s\{\oooo{e_2}^{(1)}\}
=\Gamma^{(1)}_s\{-\oooo{e_2}^{(1)}\}$. The set
\begin{eqnarray*}
\{n\in\N&|& \textup{there exists } a_1\in \Gamma^{(1)}\{e_2\}
\textup{ and }\varepsilon\in\{\pm 1\}\\
&&\textup{ with } a_1+\varepsilon nf_3\in \Gamma^{(1)}\{-e_2\}\}
\end{eqnarray*}
is well defined. Its minimum is $2$ in each case in (d)(ii)
because there $x_1>x_2=1$ and
\begin{eqnarray*}
\Gamma^{(1)}\{e_2\}\cap (e_2+\Z f_3)=e_2+\Z x_1^2f_3,\\
\Gamma^{(1)}\{-e_2\}\cap (e_2+\Z f_3)=e_2-2f_3+\Z x_1^2f_3.
\end{eqnarray*}
Its minimum is $1$ in each case in (f) because there
\begin{eqnarray*}
\Gamma^{(1)}\{e_2\}\cap (le_1+\Z f_3)
&=&le_1+\frac{1+l^2}{2}f_3+\Z l^2f_3,\\
\Gamma^{(1)}\{-e_2\}\cap (le_1+\Z f_3)
&=&le_1+\frac{-1+l^2}{2}+\Z l^2f_3.
\end{eqnarray*}

It remains to separate within each family (b), (d), (e) and (f)
the cases ((c) and (h) are single cases). 
The pair $(\oooo{H_\Z}^{(1)},\pm\oooo{I}^{(1)})$ and Lemma
\ref{t6.16} (b) allow to recover $\gcd(x_1,x_2,x_3)$ which
is as follows in these families,
\begin{eqnarray*}
\begin{array}{c|c|c|c|c}
\textup{family of cases} & (b) & (d) & (e) & (f) \\
\gcd(x_1,x_2,x_3) & x_1 & x_{12} & 2 & 1
\end{array}
\end{eqnarray*}
Within the family (b) this separates the cases.
For the family (d)
we need additionally the pair 
$(\www{x}_1,\www{x}_2)$ because 
$(x_1,x_2)=(x_{12}\www{x}_1,x_{12}\www{x}_2)$. 
The pair $(\www{x}_1,\www{x}_2)$ can be read off from
$\oooo{\Delta}^{(1)}\subset \oooo{H_\Z}^{(1)}$,
more precisely, from the relation of the 
$\Gamma^{(1)}_s$ orbits in $\oooo{\Delta^{(1)}}$
to the set $\oooo{H_\Z}^{(1),prim}\subset 
\oooo{H_\Z}^{(1)}$. 
In the family (e) one can read off $\frac{l}{2}$, and in the 
family (f) one can read off $l$ from the relation of the 
$\Gamma_s^{(1)}$ orbits in $\oooo{\Delta}^{(1)}$ to the set
$\oooo{H_\Z}^{(1),prim}\subset\oooo{H_\Z}^{(1)}$. 
\hfill$\Box$

\begin{remarks}\label{t6.24}
Let $(H_\Z,L)$ be a unimodular bilinear lattice of rank
$n\geq 2$, 
and let $\uuuu{e}$ be a triangular basis with matrix
$S=L(\uuuu{e}^t,\uuuu{e})^t\in T^{uni}_n(\Z)$.
Recall Theorem \ref{t3.7} (a).
If $S_{ij}\leq 0$ for $i<j$ then 
$(\Gamma^{(0)},\{s_{e_1}^{(0)},...,s_{e_n}^{(0)}\})$ is a 
Coxeter system, and the presentation in Definition
\ref{t3.15} of the Coxeter group $\Gamma^{(0)}$
is determined by $S$.
Especially $\Gamma^{(0)}\cong G^{fCox,n}$ if 
$S_{ij}\leq -2$ for $i<j$. 

One might hope for a similar easy control of $\Gamma^{(1)}$
if $S_{ij}\leq 0$ for $i<j$. 
In the cases with $n=2$ this works by Lemma \ref{t2.12}
and Theorem \ref{t6.10}:
\begin{eqnarray*}
\Gamma^{(1)}&\cong& \left\{\begin{array}{ll}
\{\id\}&\textup{ if }x=0,\\
SL_2(\Z)&\textup{ if }x=-1,\\ 
G^{free,2}&\textup{ if }x\leq -2.
\end{array}\right.
\end{eqnarray*}
But in the cases with $n=3$ this fails. 
The Remarks \ref{t4.17} show 
$\Gamma^{(1)}(S(\uuuu{x}))\cong \Gamma^{(1)}(S(-\uuuu{x}))$
for any $\uuuu{x}\in\Z^3$. The cases $S(\www{\uuuu{x}})$ with
$\www{\uuuu{x}}\in\Z^3_{\geq 0}$ lead by the action 
of $(G^{phi}\ltimes \www{G}^{sign})\rtimes\langle\gamma\rangle$
to all cases in Theorem \ref{t6.18}.

Especially, the cases $S(\www{\uuuu{x}})$ with 
$\www{\uuuu{x}}\in\Z^3_{\geq 2}$ contain the nice cases in
Theorem \ref{t6.18} (g) with $\Gamma^{(1)}\cong G^{free,3}$,
but also many other cases. Compare the family
$\{(3,3,l)\,|\, l\geq 2\}$ in the Examples \ref{t4.20} (iv)
or the case $S=S(2,2,2)\sim S(-2,-2,-2)\sim S(\HH_{1,2})$
with $\Gamma^{(1)}$ far from $G^{free,3}$.
\end{remarks}

\begin{remarks}\label{t6.25}
In the cases $\uuuu{x}\in\Z^3$ in Lemma \ref{t4.18} (a), so
$\uuuu{x}\in\Z^3_{\geq 2}$ with $2x_i\leq x_jx_k$ for 
$\{i,j,k\}=\{1,2,3\}$, Theorem \ref{t6.18} and Theorem 
\ref{t6.21} give rather coarse information,
\begin{eqnarray*}
\Gamma^{(1)}_u=\{\id\}\quad\textup{and}\quad 
\Gamma^{(1)}\cong\Gamma^{(1)}_s\cong G^{free,3}
\quad\textup{by Theorem \ref{t6.18}},\\
\Delta^{(1)}\to \oooo{\Delta^{(1)}}\quad\textup{is a bijection}
\qquad\textup{by Theorem \ref{t6.21}}.
\end{eqnarray*}
But it is nontrivial to determine the unique preimage in
$\Gamma^{(1),mat}$ of an element of $\Gamma^{(1)}_s$
and the unique preimage in $\Delta^{(1)}$ of an element of
$\oooo{\Delta^{(1)}}$. This holds especially for the case
$\uuuu{x}=(3,3,3)$ where $\Gamma^{(1)}_s\cong\Gamma(3)$
and $\oooo{\Delta^{(1)}}$ are known. Part (c) of the following
lemma gives for $\uuuu{x}=(3,3,3)$ 
at least an inductive way to determine the 
preimage in $\Gamma^{(1),mat}$ of a matrix in 
$\Gamma(3)\cong \Gamma^{(1)}_s$. 
\end{remarks}

\begin{lemma}\label{t6.26}
Consider the case $\uuuu{x}=(3,3,3)$. Denote by 
$L_{\P^2}:\Gamma(3)\to\Gamma^{(1),mat}$ the inverse of the
natural group isomorphism
\begin{eqnarray*}
\begin{CD}
\Gamma^{(1),mat} @>>> \Gamma^{(1)} @>>> \Gamma^{(1)}_s 
@>>> \Gamma(3),\\
s_{e_i}^{(1),mat} @>>> s_{e_i}^{(1)} @>>> 
\oooo{s_{e_i}^{(1)}} @>>> B_i \\
g^{mat} @<<< g @>>> \oooo{g} @>>> \oooo{g}^{mat}
\end{CD}
\end{eqnarray*}
with 
\begin{eqnarray*}
g(\uuuu{e})&=&\uuuu{e}\cdot g^{mat}\qquad\textup{and}\\ 
\oooo{g}(\oooo{e_1}^{(1)},\oooo{e_2}^{(1)})
&=&(\oooo{e_1}^{(1)},\oooo{e_2}^{(1)})\cdot \oooo{g}^{mat}
\end{eqnarray*}
for $g\in \Gamma^{(1)}$.  Define the subgroup of $SL_3(\Z)$
\index{$G^{(3,3,3)}\subset SL_3(\Z)$} 
\begin{eqnarray*}
G^{(3,3,3)}&:=& \{F\in SL_3(\Z)\,|\, F\equiv E_3\mmod 3, 
F\begin{pmatrix}-1\\1\\-1\end{pmatrix}=
\begin{pmatrix}-1\\1\\-1\end{pmatrix}\}.
\end{eqnarray*}
Define the map (st for standard) 
\begin{eqnarray*}
L_{st}:\Gamma(3)&\to& M_{3\times 3}(\Z),\quad 
\begin{pmatrix}a&b\\c&d\end{pmatrix}\mapsto 
\begin{pmatrix}a&b&1-a+b\\c&d&-1-c+d\\0&0&1\end{pmatrix},
\end{eqnarray*}
and the three matrices
\begin{eqnarray*}
K_1&:=& \begin{pmatrix}3&0&-3\\-3&0&3\\3&0&-3\end{pmatrix},\quad
K_2:= \begin{pmatrix}0&3&3\\0&-3&-3\\0&3&3\end{pmatrix},\\
K_3&:=&K_1+K_2= 
\begin{pmatrix}3&3&0\\-3&-3&0\\3&3&0\end{pmatrix}.
\end{eqnarray*}

(a) $L_{st}$ is an injective group homomorphism
$L_{st}:\Gamma(3)\to G^{(3,3,3)}$,
\begin{eqnarray*}
K_iK_j=0\quad\textup{for }i,j\in\{1,2,3\},\\
L_{st}(C)K_i=K_i\quad\textup{for }C\in\Gamma(3),\ 
i\in\{1,2,3\},\\
G^{(3,3,3)}=\{L_{st}(C)+\alpha K_1+\beta K_2\,|\, 
C\in \Gamma(3),\alpha,\beta\in\Z\}.
\end{eqnarray*}
The following sequence is an exact sequence of group 
homomorphisms,
\begin{eqnarray*}
\begin{CD}
\{1\} @>>> \Z^2 @>>> G^{(3,3,3)} @>>> \Gamma(3) @>>> \{1\}\\
 @. (\alpha,\beta) @>>> E_3+\alpha K_1+\beta K_2 @. @. \\
 @. @. L_{st}(C)+\alpha K_1+\beta K_2 @>>> C @. 
\end{CD}
\end{eqnarray*}
$L_{st}$ is a splitting of this exact sequence.

(b) $\Gamma^{(1),mat}\subset G^{(3,3,3)}$, and 
$L_{\P^2}:\Gamma(3)\to \Gamma^{(1),mat}$ is another splitting
of the exact sequence in part (b). It satisfies
\begin{eqnarray*}
L_{\P^2}(B_1)&=&L_{st}(B_1)=s_{e_1}^{(1),mat}
=\begin{pmatrix}1&-3&-3\\0&1&0\\0&0&1\end{pmatrix}
\quad\textup{for }B_1=\begin{pmatrix}1&-3\\0&1\end{pmatrix},\\
L_{\P^2}(B_2)&=&L_{st}(B_2)=s_{e_2}^{(1),mat}
=\begin{pmatrix}1&0&0\\3&1&-3\\0&0&1\end{pmatrix}
\quad\textup{for }B_2=\begin{pmatrix}1&0\\3&1\end{pmatrix},\\
L_{\P^2}(B_3)&=&L_{st}(B_3)+K_3=s_{e_3}^{(1),mat}
=\begin{pmatrix}1&0&0\\0&1&0\\3&3&1\end{pmatrix}
\quad\textup{for }B_3=\begin{pmatrix}-2&-3\\3&4\end{pmatrix},\\
L_{\P^2}(B_3^{-1})&=&L_{st}(B_3^{-1})-K_3.
\end{eqnarray*}

(c) An arbitrary element $C\in\Gamma(3)$ can be written in a 
unique way as a product 
\begin{eqnarray*}
C=C_1 B_3^{\varepsilon_1} C_2 B_3^{\varepsilon_2} C_3 ...
C_m B_3^{\varepsilon_m} C_{m+1}\\
\textup{with}\quad 
m\in\Z_{\geq 0},\ C_1,...,C_{m+1}\in 
\langle B_1^{\pm 1},B_2^{\pm 1}\rangle,\ 
\varepsilon_1,...,\varepsilon_m\in\{\pm 1\}.
\end{eqnarray*}
Then
\begin{eqnarray*}
L_{\P^2}(C)&=& L_{st}(C) + K_3\Bigl( 
\varepsilon_1 L_{st}(C_2B_3^{\varepsilon_2}C_3...C_m 
B_3^{\varepsilon_m}C_{m+1}) \Bigr. \\
&& \hspace*{2cm}+\varepsilon_2 
L_{st}(C_3B_3^{\varepsilon_3}C_4...C_m B_3^{\varepsilon_m}C_{m+1}) \\
&& \hspace*{2cm} \Bigl. + ... + 
\varepsilon_m L_{st}(C_{m+1})\Bigr).
\end{eqnarray*}
\end{lemma}

{\bf Proof:}
The parts (a) and (b) are easy.

(c) By part (b) $L_{\P^2}(C_j)=L_{st}(C_j)$ and 
$L_{\P^2}(B_3^{\varepsilon_j})=L_{st}(B_3^{\varepsilon_j})
+\varepsilon_j K_3$, so
\begin{eqnarray*}
L_{\P^2}(C)&=& L_{\P^2}(C_1)L_{\P^2}(B_3^{\varepsilon_1})
L_{\P^2}(C_2)...L_{\P^2}(C_m)L_{\P^2}(B_3^{\varepsilon_m})
L_{\P^2}(C_{m+1})\\
&=& L_{st}(C_1)(L_{st}(B_3^{\varepsilon_1})+\varepsilon_1K_3)
L_{st}(C_2)...\\
&& L_{st}(C_m)(L_{st}(B_3^{\varepsilon_m})+
\varepsilon_mK_3)L_{st}(C_{m+1}).
\end{eqnarray*}
Observe $K_3 L_{st}(\www{C})K_3=K_3K_3=0$ 
for $\www{C}\in\Gamma(3)$.
Therefore if one writes the product above as a sum of
$2^m$ terms, only the $1+m$ terms do not vanish in which $K_3$
turns up at most once. This leads to the claimed formula
for $L_{\P^2}(C)$.\hfill$\Box$

\chapter{Distinguished bases in the rank 2 and rank 3 cases}\label{s7}
\setcounter{equation}{0}
\setcounter{figure}{0}

In section \ref{s3.3} we introduced the set of distinguished bases
of a unimodular bilinear lattice $(H_\Z,L,\uuuu{e})$ with a triangular
basis. It is the orbit $\BB^{dist}=\Br_n\ltimes\{\pm 1\}^n(\uuuu{e})$ 
of $\uuuu{e}$ under the group $\Br_n\ltimes\{\pm 1\}^n$.
We also posed the question when this set can be characterized 
in an
easy way, more precisely, when the inclusions in \eqref{3.3} or
\eqref{3.4} are equalities,
\begin{eqnarray*}
\BB^{dist}\subset\{\uuuu{v}\in(\Delta^{(0)})^n\,|\, 
s_{v_1}^{(0)}...s_{v_n}^{(0)}=-M\},\hspace*{2cm}(3.3)\\
\BB^{dist}\subset\{\uuuu{v}\in(\Delta^{(1)})^n\,|\, 
s_{v_1}^{(1)}...s_{v_n}^{(1)}=M\}.\hspace*{2.4cm}(3.4)
\end{eqnarray*}
Theorem \ref{t3.2} (a) and (b) imply that \eqref{3.4} is an equality
if $\Gamma^{(1)}$ is a free group with generators 
$s_{e_1}^{(1)},...,s_{e_n}^{(1)}$ and that \eqref{3.3} is an equality
if $\Gamma^{(0)}$ is a free Coxeter group with generators 
$s_{e_1}^{(0)},...,s_{e_n}^{(0)}$, see the Examples \eqref{t3.23} (iv).
More generally, if $(\Gamma^{(0)},s_{e_1}^{(0)},...,s_{e_n}^{(0)})$
is a Coxeter system (Definition \ref{t3.5}) then by Theorem \ref{t3.6}
\eqref{3.3} is an equality, see the Examples \ref{t3.23} (v).
It is remarkable that the property $\sum_{i=1}^n \Z v_i=H_\Z$,
which each distinguished basis $\uuuu{v}\in \BB^{dist}$ satisfies
is not needed in the characterization in these cases.

These are positive results. In the sections \ref{s3.1}--\ref{s3.3}
we study systematically all cases of rank 2 and 3 and find also
negative results. 

In rank 2 in section \ref{s6.2} \eqref{3.3}
is always an equality, and \eqref{3.4} is an equality in all cases
except the case $A_1^2$. 

In the even rank 3 cases in section \ref{s7.2} \eqref{3.3} is 
in all cases except the case $\HH_{1,2}$ an equality. 
In the case $\HH_{1,2}$ the set on the right hand side contains
$\Br_3\ltimes\{\pm 1\}^3$ orbits of tuples $\uuuu{v}$ with
arbitrary large finite index $[H_\Z:\sum_{i=1}^3\Z v_i]$
and two orbits with index 1, $\BB^{dist}$ and one other orbit.

In the odd rank 3 cases in section \ref{s7.3} we understand
the set $B_1\cup B_2$ of triples $\uuuu{x}\in\Z^3$ such that 
\eqref{3.4} is an equality, and we also know a set $B_3\cup B_4$
of triples $\uuuu{x}\in\Z^3$ such that \eqref{3.4} becomes an
equality if one adds the condition $H_\Z=\sum_{i=1}^3 \Z v_i$.
But for $\uuuu{x}\in\Z^3-\cup_{j=1}^4B_j$, we know little.

Section \ref{s7.4} builds on section \ref{4.4} where for a 
unimodular bilinear lattice $(H_\Z,L,\uuuu{e})$ the stabilizer
$(\Br_3)_{\uuuu{x}/\{\pm 1\}^3}$ had been determined.
It determines the stabilizer $(\Br_3)_{\uuuu{e}/\{\pm 1\}^3}$. 
It uses the systematic results in chapter \ref{s5} on the group
$G_\Z$ and on the map $Z:(\Br_n\ltimes\{\pm 1\}^n)_S
\to G_\Z$ in the rank 3 cases. 

In the sections \ref{s4.3} and \ref{s4.4} 
the pseudo-graph $\GG(\uuuu{x})$ with
vertex set an orbit $\Br_3(\uuuu{x}/\{\pm 1\})$ and edge set
from generators of the group $G^{phi}\rtimes\langle\gamma\rangle$
had been crucial. In section \ref{s7.4} we introduce a variant
with the same vertex set, but different edge set, namely 
(now) oriented edges coming
from the elementary braids $\sigma_i^{\pm 1}$.
We also define the much larger $\sigma$-pseudo-graph
with vertex set a set $\BB^{dist}/\{\pm 1\}^3$ of distinguished
bases up to signs and oriented edges coming from the elementary braids 
$\sigma_i^{\pm 1}$. We consider especially the examples where
the set $\Br_3(\uuuu{x}/\{\pm 1\}^3)$ is finite.

\section{Distinguished bases in the rank 2 cases}
\label{s7.1}

In the rank 2 cases the inclusion \eqref{3.3} is always an equality,
and the inclusion \eqref{3.4} is almost always an equality,
namely in all cases except the case $A_1^2$. 

\begin{theorem}\label{t7.1}
Let $(H_\Z,L,\uuuu{e})$ be a unimodular bilinear lattice of rank
2 with a triangular basis $\uuuu{e}=(e_1,e_2)$ with matrix
$S=S(x)=\begin{pmatrix}1&x\\0&1\end{pmatrix}
=L(\uuuu{e}^t,\uuuu{e})^t$ with $x\in\Z$. Fix $k\in\{0,1\}$.

(a) The inclusion \eqref{3.3} respectively \eqref{3.4} 
in Remark \ref{t3.19} is an equality in all cases except
the odd case $A_1^2$, so the case $(k,x)=(1,0)$. In that case
the right hand side in \eqref{3.4} splits into the orbits
of the three pairs $(e_1,e_1)$, $(e_1,e_2)$, $(e_2,e_2)$. 

(b) The stabilizers in $\Br_2$ 
of $S/\{\pm 1\}^2$ and of $\uuuu{e}/\{\pm 1\}^2$ are 
\begin{eqnarray*}
(\Br_2)_{S/\{\pm 1\}^2}&=&\Br_2\quad\textup{and}\\
(\Br_2)_{\uuuu{e}/\{\pm 1\}^2}&=& \left\{\begin{array}{ll}
\langle \sigma_1^2\rangle & \textup{ if }x=0,\\
\langle \sigma_1^3\rangle & \textup{ if }x\in\{\pm 1\},\\
\{\id\} & \textup{ if }|x|\geq 2.
\end{array}\right.
\end{eqnarray*}
\end{theorem}

{\bf Proof:}
(a) The even and odd cases $A_1^2$: 
See the Examples \ref{t3.23} (iii). 

The cases with $|x|\geq 2$: Theorem \ref{t6.8} (c)+(d)
and Theorem \ref{t6.10} (c)+(d) show
\begin{eqnarray*}
\Gamma^{(k)}\cong \left\{\begin{array}{ll}
G^{fCox,2}&\textup{ with generators }s_{e_1}^{(0)},
s_{e_2}^{(0)}\textup{ if }k=0,\\
G^{free,2}&\textup{ with generators }s_{e_1}^{(1)},
s_{e_2}^{(1)}\textup{ if }k=1.
\end{array}\right.
\end{eqnarray*}
The Examples \ref{t3.23} (iv) apply and give
equality in \eqref{3.3} and \eqref{3.4}. 

The cases with $x=\pm 1$: We can restrict to the case $x=-1$.
The even case is a simple case of Example \ref{t3.23} (v)
(in the Remarks \ref{t7.2} we will offer an elementary proof
for the even case).

It remains to show equality in \eqref{3.4} in the odd case
$(k,x)=(1,-1)$. Consider $\uuuu{v}\in (\Delta^{(1)})^2$
with $s_{v_1}^{(1)}s_{v_2}^{(1)}=M$. 
Let $b:=I^{(1)}(v_1,v_2)\in\Z$. 
If $b=0$ then $v_2=\pm v_1$ and $M=(s_{v_1}^{(1)})^2$
would have an eigenvalue $1$, a contradiction.
Therefore $b\neq 0$ and $\Z v_1+\Z v_2$ has rank 2. Then
\begin{eqnarray*}
M(\uuuu{v})&=& s_{v_1}^{(1)}s_{v_2}^{(1)}(\uuuu{v})
= s_{v_1}^{(1)}(v_1+bv_2,v_2)\\
&=&(v_1+bv_2-b^2v_1,v_2-bv_1)
=\uuuu{v}\begin{pmatrix}1-b^2&-b\\b&1\end{pmatrix},
\end{eqnarray*}
$1=\tr M=(1-b^2)+1$, so $b=\pm 1$. 
By possibly changing the sign of $v_2$, we can suppose 
$b=-1=x$. Then 
\begin{eqnarray*}
I^{(1)}(\uuuu{v}^t,\uuuu{v})
=\begin{pmatrix}0&-1\\1&0\end{pmatrix}.
\end{eqnarray*}
Therefore $\uuuu{v}$ is a $\Z$-basis of $H_\Z$, and the
automorphism $g\in\Aut(H_\Z)$ with 
$(g(e_1),g(e_2))=(v_1,v_2)$ is in $O^{(1)}$.
By Lemma \ref{t3.22} (a) 
\begin{eqnarray*}
gMg^{-1}=g\circ((\pi_2\circ \pi_2^{(1)})(\uuuu{e}))\circ g^{-1}
=(\pi_2\circ\pi_2^{(1)})(\uuuu{v})=
s_{v_1}^{(1)}s_{v_2}^{(1)}= M,
\end{eqnarray*}
so $gMg^{-1}=M$, so $g\in G_\Z^{(1)}=G_\Z^M\cap O^{(1)}$. 
Theorem \ref{t5.5} can be applied and gives 
$\stackrel{(*)}{=}$, 
\begin{eqnarray*}
G_\Z^{(1)}\stackrel{(*)}{=}G_\Z\stackrel{(*)}{=} 
\{\pm (M^{root})^l\,|\, l\in\Z\}
\stackrel{(*)}{=}Z(\Br_2\ltimes \{\pm 1\}^2).
\end{eqnarray*}
Therefore $\uuuu{v}\in\BB^{dist}$. This shows equality in
\eqref{3.4}.

(b) Because of $\sigma_1\begin{pmatrix}1&x\\0&1\end{pmatrix}
=\begin{pmatrix}1&-x\\0&1\end{pmatrix}$,
the stabilizer $(\Br_2)_{S/\{\pm 1\}^2}$ is the whole group
$\Br_2=\langle \sigma_1\rangle$.  If $x=0$,
\begin{eqnarray*}
(e_1,e_2)\stackrel{\sigma_1}{\mapsto} (e_2,e_1)
\stackrel{\sigma_1}{\mapsto}(e_1,e_2),
\quad\textup{so }(\Br_2)_{\uuuu{e}/\{\pm 1\}^2}
=\langle \sigma_1^2\rangle.
\end{eqnarray*}
If $x=-1$, 
\begin{eqnarray*}
(e_1,e_2)\stackrel{\sigma_1}{\mapsto} (e_1+e_2,e_1)
\stackrel{\sigma_1}{\mapsto}(-e_2,e_1+e_2)
\stackrel{\sigma_1}{\mapsto}(e_1,-e_2),\\
\textup{so }(\Br_2)_{\uuuu{e}/\{\pm 1\}^2}
=\langle \sigma_1^3\rangle.
\end{eqnarray*}
If $|x|\geq 2$ Theorem \ref{t3.2} (a) or (b) and
$\Gamma^{(1)}\cong G^{free,2}$ or 
$\Gamma^{(0)}\cong G^{fCox,2}$ show
$(\Br_2)_{\uuuu{e}/\{\pm 1\}^2}=\{\id\}$.\hfill$\Box$

\begin{remarks}\label{t7.2}
(i) A direct elementary proof of equality in \eqref{3.3}
for the even case $A_2$, so $(k,x)=(0,-1)$, is instructive.
Recall from Theorem \ref{t6.8} (b) that
$\Delta^{(0)}=\{\pm e_1,\pm e_2,\pm (e_1+e_2)\}$. 
The map $\pi_2\circ\pi_2^{(0)}:(\Delta^{(0)})^2\to \Gamma^{(0)}$
has the three values $-M$, $M^2$ and $\id$ and the three fibers
\begin{eqnarray*}
(\pi_2\circ\pi_2^{(0)})^{-1}(-M)
&=& \{(\pm e_1,\pm e_2),(\pm (e_1+e_2),\pm e_1),
(\pm e_2,\pm (e_1+e_2)\}\\
&=&\BB^{dist},\\
(\pi_2\circ\pi_2^{(0)})^{-1}(M^2)
&=& \{(\pm e_2,\pm e_1),(\pm (e_1+e_2),\pm e_2),
(\pm e_1,\pm (e_1+e_2)\}\\
&=&
\Br_2\ltimes\{\pm 1\}^2(e_2,e_1),\\
(\pi_2\circ\pi_2^{(0)})^{-1}(\id)
&=& \{(\pm e_1,\pm e_1),(\pm e_2,\pm e_2),
(\pm (e_1+e_2),\pm (e_1+e_2)\}.
\end{eqnarray*}
This gives equality in \eqref{3.3} in the case $(k,x)=(0,-1)$.

(ii) Also in the cases $(k=0,x\leq -2)$ a direct elementary proof
of equality in \eqref{3.3} is instructive. 
Equality in \eqref{3.3} for $(k,x)=(0,-2)$ and Theorem 
\ref{t6.8} (d) (iv) imply equality in \eqref{3.3} for
$(k=0,x\leq -3)$. Therefore we restrict to the case
$(k,x)=(0,-2)$. Recall
\begin{eqnarray*}
\Rad I^{(0)}=\Z f_1\quad\textup{with}\quad f_1=e_1+e_2.
\end{eqnarray*}
By Theorem \ref{t6.8} (c)
\begin{eqnarray*}
\Delta^{(0)}=(e_1+\Z f_1)\ \dot\cup\ (-e_1+\Z f_1)
=(e_1+\Z f_1)\ \dot\cup\ (e_2+\Z f_1).
\end{eqnarray*}
One easily sees for $b_1,b_2\in\Z$
\begin{eqnarray*}
s_{e_1+b_1f_1}^{(0)}s_{e_2+b_2f_1}^{(0)}=-M
&\iff& b_1+b_2=0,
\end{eqnarray*}
thus
\begin{eqnarray*}
\{\uuuu{v}\in (\Delta^{(0)})^2\,|\, 
s_{v_1}^{(0)}s_{v_2}^{(0)}=-M\}
=\{\pm (e_1+bf_1),\pm (e_2-bf_1)\,|\, b\in\Z\}.
\end{eqnarray*}
This set is a single $\Br_2\ltimes \{\pm 1\}^2$ orbit
because of 
\begin{eqnarray*}
\delta_2\sigma_1(e_1+bf_1,e_2-bf_1)=(e_1+(b+1)f_1,e_2-(b+1)f_1).
\hspace*{1cm}\Box
\end{eqnarray*}
\end{remarks}

\section{Distinguished bases
in the even rank 3 cases}
\label{s7.2}

In the even cases with $n=3$ we have complete results on the
question when the inclusion in \eqref{3.3} is an equality.
It is one in all cases except the case $\HH_{1,2}$.

\begin{theorem}\label{t7.3}
Let $(H_\Z,L,\uuuu{e})$ be a unimodular bilinear lattice
of rank 3 with a triangular basis $\uuuu{e}$ with
matrix $S=S(\uuuu{x})=L(\uuuu{e}^t,\uuuu{e})^t\in T^{uni}_3(\Z)$.

(a) Suppose $S\notin (\Br_3\ltimes\{\pm 1\}^3)(S(\HH_{1,2}))$.
Then the inclusion in \eqref{3.3} is an equality.

(b) Suppose $S=S(\HH_{1,2})=S(-2,2,-2)$. Recall the basis
$(f_1,f_2,f_3)=\uuuu{e}
\begin{pmatrix}1&0&1\\1&1&1\\0&1&1\end{pmatrix}$ 
of $H_\Z$ with $\Rad I^{(1)}=\Z f_1\oplus \Z f_2$ and
$\Rad I^{(0)}=\Z f_3$. 
The set $\{\uuuu{v}\in (\Delta^{(0)})^3\,|\, 
(\pi_3\circ\pi_3^{(0)})(\uuuu{v})=-M\}$ splits into 
countably many orbits. The following list gives one
representative for each orbit,
\begin{eqnarray*}
(f_3-g_1,-f_3+g_1+g_2,f_3-g_2)\quad\textup{with}\quad 
\begin{pmatrix}g_1\\g_2\end{pmatrix}
=\begin{pmatrix}0&c_2\\c_1&c_3\end{pmatrix}
\begin{pmatrix}f_1\\f_2\end{pmatrix},\\
c_1\in\N\textup{ odd},\ c_2\in\Z\textup{ odd},
\ c_3\in\{0,1,...,|c_2|-1\}.
\end{eqnarray*}
The sublattice $\langle f_3-g_1,-f_3+g_1+g_2,f_3-g_2\rangle
=\langle f_3,g_1,g_2\rangle\subset H_\Z$ has finite index
$c_1\cdot |c_2|$ in $H_\Z$. It is $H_\Z$ in the following
two cases: 
\begin{eqnarray*}
\uuuu{e}=(f_3-f_2,-f_3+f_1+f_2,f_3-f_1),\\
\textup{so}\quad
(g_1,g_2,c_1,c_2,c_3)=(f_2,f_1,1,1,0),\\
(f_3+f_2,-f_3+f_1-f_2,f_3-f_1),\\
\textup{so}\quad
(g_1,g_2,c_1,c_2,c_3)=(-f_2,f_1,1,-1,0)
\end{eqnarray*}
(see also Example \ref{t3.23} (ii)) for the second case).
\end{theorem}

{\bf Proof:}
(a) We can replace $\uuuu{e}$ by an arbitrary element
$\www{\uuuu{e}}\in\BB^{dist}$. By Theorem \ref{t4.6}
the following cases exhaust all $\Br_3\ltimes\{\pm 1\}^3$ orbits
except that of $\HH_{1,2}$:
\begin{list}{}{}
\item[(A)]
$(H_\Z,L,\uuuu{e})$ is irreducible with 
$\uuuu{x}\in\Z^3_{\leq 0}$.
\item[(B)]
$r(\uuuu{x})\leq 0$ and $\uuuu{x}\neq (0,0,0)$.
\item[(C)]
$\uuuu{x}=(x_1,0,0)$ with $x_1\in\Z_{\leq 0}$, so 
$(H_\Z,L,\uuuu{e})$ is reducible (this includes the case
$A_1^3$).
\item[(D)]
$\uuuu{x}=(-l,2,-l)$ with $l\geq 3$.
\end{list}

The cases (A): $\Gamma^{(0)}$ is a Coxeter group 
by Theorem \ref{t3.7} (a). Theorem \ref{t3.7} (b) applies.

The cases (B): By Theorem \ref{t6.11} (g) $\Gamma^{(0)}$
is a free Coxeter group with generators 
$s_{e_1}^{(0)},s_{e_2}^{(0)},s_{e_3}^{(0)}$.
Theorem \ref{t3.7} (b) or Example \ref{t3.23} (iv) can be used.

The cases (C):
Consider a triple $\uuuu{v}\in (\Delta^{(0)})^3$ with 
$s_{v_1}^{(0)}s_{v_2}^{(0)}s_{v_3}^{(0)}=-M$. 
The set $\Delta^{(0)}$ splits into the subsets
$\Delta^{(0)}\cap (\Z e_1+\Z e_2)$ and $\{\pm e_3\}$.
Compare $-M|_{\Z e_3}=-\id|_{\Z e_3}$ with 
$s_{e_3}^{(0)}|_{\Z e_3}=-\id|_{\Z e_3}$ and 
$s_a^{(0)}|_{\Z e_3}=\id|_{\Z e_3}$ for 
$a\in \Delta^{(0)}\cap (\Z e_1+\Z e_2)$. 
All three $v_i\in\{\pm e_3\}$ is impossible because
$(s_{e_3}^{(0)})^3\neq -M$. Therefore there are $i,j,k$ with 
$\{i,j,k\}=\{1,2,3\}$, $i<j$, $v_i,v_j\in \Delta^{(0)}
\cap(\Z e_1+\Z e_2)$ and $v_k\in\{\pm e_3\}$.
The reflection $s_{v_k}^{(0)}$ acts trivially
on $\Z e_1+\Z e_2$ and commutes with 
$s_{v_i}^{(0)}$ and $s_{v_j}^{(0)}$. Therefore
\begin{eqnarray*}
s_{v_i}^{(0)}s_{v_j}^{(0)}s_{v_k}^{(0)}
=s_{v_1}^{(0)}s_{v_2}^{(0)}s_{v_3}^{(0)}
=(\pi_3\circ\pi_3^{(0})(\uuuu{v})=-M=
s_{e_1}^{(0)}s_{e_2}^{(0)}s_{e_3}^{(0)},\\
\textup{so}\quad s_{v_i}^{(0)}s_{v_j}^{(0)}
=s_{e_1}^{(0)}s_{e_2}^{(0)}.
\end{eqnarray*}
The reflections $s_{v_i}^{(0)}$, $s_{v_j}^{(0)}$,
$s_{e_1}^{(0)}$ and $s_{e_2}^{(0)}$ 
act trivially on $\Z e_3$. 
The inclusion in \eqref{3.3} is an equality because of
Theorem \ref{t7.1} (a) for the rank 2 cases.

The cases (D): The proof of these cases is prepared by
Lemma \ref{t7.4} and Lemma \ref{t7.5}.
The proof comes after the proof of Lemma \ref{t7.5}.

(b) The proof of part (b) is prepared by Lemma \ref{t7.6}
and comes after the proof of Lemma \ref{t7.6}.
\hfill ($\Box$)

\bigskip
The following lemma is related to  
$t^+_\lambda$ in Lemma \ref{t6.17}.
Recall also $j^{(k)}:H_\Z\to H_\Z^\sharp$, 
$a\mapsto I^{(k)}(a,.)$, in Definition \ref{t6.1}.

\begin{lemma}\label{t7.4}
Let $(H_\Z,L)$ be a unimodular bilinear lattice
of rank $n\in\N$. Fix $k\in\{0,1\}$. 
Suppose $\Rad I^{(k)}\neq\{0\}$ and choose an element
$f\in \Rad I^{(k)}-\{0\}$. Denote 
\index{$t_\lambda$}\index{$\Hom_{0,f}(H_\Z,\Z)$}
\begin{eqnarray*}
\Hom_{0,f}(H_\Z,\Z):=
\{\lambda:H_\Z\to\Z\,|\, \lambda \textup{ is }\Z\textup{-linear},
\lambda(f)=0\},\\
t_\lambda:H_\Z\to H_\Z\textup{ with }
t_\lambda(a)=a+\lambda(a)f\textup{ for }
\lambda\in \Hom_{0,f}(H_\Z,\Z).
\end{eqnarray*}
Then $t_\lambda\in O^{(k),Rad}_u$. 
The map 
\begin{eqnarray*}
\Hom_{0,f}(H_\Z,\Z)\to O^{(k),Rad}_u,\quad\lambda\mapsto
t_\lambda,
\end{eqnarray*}
is an injective group homomorphism.
For $b\in R^{(k)}$ (with $R^{(1)}=H_\Z$, see \ref{t3.9} (i)) 
and $a\in\Z$
\begin{eqnarray*}
s_{b+af}^{(k)} &=& s_b^{(k)}\circ t_{-aj^{(k)}(b)}
=t_{(-1)^kaj^{(k)}(b)}\circ s_b^{(k)},\\\
s_b^{(k)}\circ t_\lambda &=& 
t_{\lambda-(-1)^k\lambda(b)j^{(k)}(b)}\circ s_b^{(k)}.
\end{eqnarray*}
\end{lemma}

{\bf Proof:}
The proof is straightforward. We skip the details. \hfill$\Box$

\bigskip
The following lemma studies the Hurwitz action of $\Br_3$ on 
triples of reflections in $G^{fCox,2}$. It is related
to Theorem \ref{t3.2} (b).

\begin{lemma}\label{t7.5}
As in Definition \ref{t3.1}, $G^{fCox,2}$ denotes the 
free Coxeter group with two generators $z_1$ and $z_2$,
so generating relations are $z_1^2=z_2^2=1$. 

(a) Its set of reflections is
\begin{eqnarray*}
\Delta(G^{fCox,2})
=\bigcup_{i=1}^2 \{wz_iw^{-1}\,|\, w\in G^{fCox,2}\}
=\{(z_1z_2)^mz_1\,|\, m\in\Z\}.
\end{eqnarray*}
The complement of this set is
\begin{eqnarray*}
G^{fCox,2}-\Delta(G^{fCox,2})
=\{(z_1z_2)^m\,|\, m\in\Z\}.
\end{eqnarray*}
$\Delta(G^{fCox,2})$ respectively its complement consists
of the elements which can be written as words of odd
respectively even length in $z_1$ and $z_2$.

(b) The set
\begin{eqnarray*}
\{(w_1,w_2,w_3)\in (\Delta(G^{fCox,2})^3\,|\, 
w_1w_2w_3=z_1z_2z_1\}
\end{eqnarray*}
is the disjoint union of the following $\Br_3$ orbits:
\begin{eqnarray*}
\dot\bigcup_{m\in \Z_{\geq 0}}
\Br_3\bigl((z_1z_2z_1, (z_1z_2)^{1-m}z_1,(z_1z_2)^{1-m}z_1)\bigr).
\end{eqnarray*}
\end{lemma}

{\bf Proof:}
(a) Clear.

(b) The map
\begin{eqnarray*}
\{(w_1,w_2,w_3)\in(\Delta(G^{fCox,2}))^3\,|\, 
w_1w_2w_3=z_1z_2z_1\}&\to& M_{2\times 1}(\Z),\\
(w_1,w_2,w_3)&\mapsto& (m_1,m_2)^t\\ 
\textup{ with }
w_1w_2=(z_1z_2)^{m_1},\ w_2w_3=(z_1z_2)^{m_2},
\end{eqnarray*}
is a bijection because
\begin{eqnarray*}
z_1z_2z_1&=&(w_1w_2)w_2(w_2w_3),\quad\textup{so}\\
w_2&=&(w_1w_2)^{-1}z_1z_2z_1(w_2w_3)^{-1},\\
w_1&=&(w_1w_2)w_2,\\ 
w_3&=&w_2(w_2w_3),
\end{eqnarray*}
so a given column vector $(m_1,m_2)^t$ has a unique preimage.

The Hurwitz action of $\Br_3$ on the set on the left hand side
of the bijection above 
translates as follows to an action on $M_{2\times 1}(\Z)$.
\begin{eqnarray*}
\sigma_1(w_1,w_2,w_3)&=& (w_1w_2w_1,w_1,w_3),\\
w_1w_2w_1\cdot w_1&=&w_1w_2=(z_1z_2)^{m_1}, \\
w_1w_3&=&(w_1w_2)(w_2w_3)=(z_1z_2)^{m_1+m_2},\\
\sigma_1\begin{pmatrix}m_1\\m_2\end{pmatrix}
&=&\begin{pmatrix}m_1\\m_1+m_2\end{pmatrix}
=\begin{pmatrix}1&0\\1&1\end{pmatrix}
\begin{pmatrix}m_1\\m_2\end{pmatrix},\\
\sigma_2(w_1,w_2,w_3)&=& (w_1,w_2w_3w_2,w_2),\\
w_1\cdot w_2w_3w_2=(w_1w_2)(w_2w_3)^{-1}
&=&(z_1z_2)^{m_1-m_2},\\
w_2w_3w_2\cdot w_2&=&w_2w_3=(z_1z_2)^{m_2},\\
\sigma_2\begin{pmatrix}m_1\\m_2\end{pmatrix}
&=&\begin{pmatrix}m_1-m_2\\m_2\end{pmatrix}
=\begin{pmatrix}1&-1\\0&1\end{pmatrix}
\begin{pmatrix}m_1\\m_2\end{pmatrix}.
\end{eqnarray*}
So $\Br_3$ acts as multiplication with matrices in
$SL_2(\Z)$ from the left on $M_{2\times 1}(\Z)$. 
Each orbit has a unique element of the shape
$\begin{pmatrix}m\\0\end{pmatrix}$ with $m\in\Z_{\geq 0}$.
This element corresponds to 
\begin{eqnarray*}
(z_1z_2z_1,(z_1z_2)^{1-m}z_1,(z_1z_2)^{1-m}z_1).
\hspace*{1cm}\Box
\end{eqnarray*}

\bigskip
{\bf Proof} of Theorem \ref{t7.3} (a) in the cases (D),
$\uuuu{x}=(-l,2,-l)$ with $l\geq 3$:
Recall from Theorem \ref{t6.11} (f)
\begin{eqnarray*}
\Rad I^{(0)}&=&\Z f_1\quad\textup{with}\quad 
f_1=e_1-e_3,\\
\Gamma^{(0)}_s&\cong& G^{fCox,2}\quad\textup{with generators }
z_1=\oooo{s_{e_1}^{(0)}}=\oooo{s_{e_3}^{(0)}}, 
\ z_2=\oooo{s_{e_2}^{(0)}}.
\end{eqnarray*}
Suppose $\uuuu{v}\in (\Delta^{(0)})^3$ with 
$s_{v_1}^{(0)}s_{v_2}^{(0)}s_{v_3}^{(0)}=-M$.
We want to show $\uuuu{v}\in \BB^{dist}$ or equivalently
$(s_{v_1}^{(0)},s_{v_2}^{(0)},s_{v_3}^{(0)})\in \RR^{(0),dist}$.

First we look at the images in $\Gamma^{(0)}_s$:
$\oooo{s_{v_1}^{(0)}} \oooo{s_{v_2}^{(0)}} 
\oooo{s_{v_3}^{(0)}}=z_1z_2z_1$.
Because of Lemma \ref{t7.5} (b), we can make a suitable
braid group action and then suppose
\begin{eqnarray*}
(\oooo{s_{v_1}^{(0)}} ,\oooo{s_{v_2}^{(0)}}, 
\oooo{s_{v_3}^{(0)}})=(z_1z_2z_1,r,r)\textup{ with }
r=(z_1z_2)^{1-m}z_1\textup{ for some }m\in\Z_{\geq 0}.
\end{eqnarray*}
Write $\www{e_2}:=s_{e_1}^{(0)}(e_2)=e_2+le_1$ and observe
\begin{eqnarray*}
z_1z_2z_1=\oooo{s_{e_1}^{(0)}s_{e_2}^{(0)}s_{e_3}^{(0)}}
=\oooo{s_{\www{e_2}}^{(0)}}.
\end{eqnarray*}
After possibly changing the signs of $v_1$ and $v_3$, 
$\oooo{s_{v_1}^{(0)}}=\oooo{s_{\www{e_2}}^{(0)}}$ and
$\oooo{s_{v_2}^{(0)}}=r=\oooo{s_{v_3}^{(0)}}$ imply
\begin{eqnarray*}
v_1=\www{e_2}+a_1f_1\quad\textup{and}\quad 
v_3=v_2+a_2f_1\quad\textup{for some }a_1,a_2\in\Z.
\end{eqnarray*}
With Lemma \ref{t7.4} and 
$f_1=(f\textup{ in Lemma \ref{t7.4}})$ we calculate
\begin{eqnarray*}
-M&=& s_{e_1}^{(0)}s_{e_2}^{(0)}s_{e_3}^{(0)}
= s_{e_1}^{(0)}s_{e_2}^{(0)}s_{e_1-f_1}^{(0)}\\\
&=& s_{e_1}^{(0)}s_{e_2}^{(0)}s_{e_1}^{(0)}
t_{j^{(0)}(e_1)}
=s_{\www{e_2}}^{(0)}t_{j^{(0)}(e_1)},\\
-M&=& s_{v_1}^{(0)}s_{v_2}^{(0)}s_{v_3}^{(0)}
= s_{\www{e_2}+a_1f_1}^{(0)}s_{v_2}^{(0)}s_{v_2+a_2f_1}^{(0)}\\
&=& s_{\www{e_2}}^{(0)}t_{-a_1j^{(0)}(\www{e_2})}
s_{v_2}^{(0)}s_{v_2}^{(0)}t_{-a_2j^{(0)}(v_2)}\\
&=& s_{\www{e_2}}^{(0)}t_{-a_1j^{(0)}(\www{e_2})-a_2j^{(0)}(v_2)},
\end{eqnarray*}
so
\begin{eqnarray*}
j^{(0)}(e_1)=-a_1j^{(0)}(\www{e_2})-a_2j^{(0)}(v_2).
\end{eqnarray*}
Write 
\begin{eqnarray*}
\oooo{v_2}^{(0)}=b_1\oooo{e_1}^{(0)}+b_2\oooo{e_2}^{(0)}
\quad\textup{with }b_1,b_2\in\Z.
\end{eqnarray*} 
By Theorem \ref{t6.14} (f) the tuple 
$(\oooo{H_{\Z}}^{(0)},\oooo{I}^{(0)}, (\oooo{e_1}^{(0)},
\oooo{e_2}^{(0)}))$ is isomorphic to the corresponding tuple
from the $2\times 2$ matrix $S(-l)=
\begin{pmatrix}1&-l\\0&1\end{pmatrix}$. 
The set of roots of this tuple is called $R^{(0)}(S(-l))$.
It contains $\oooo{v_2}^{(0)}$, $\oooo{e_1}^{(0)}$,
$\oooo{e_1}^{(0)}$. By Theorem \ref{t6.8} (d)(i) the map
\begin{eqnarray*}
R^{(0)}(S(-l))&\to& \{\textup{units in }\Z[\kappa_1]
\textup{ with norm }1\}
=\{\pm \kappa_1^m\,|\, m\in\Z\}\\
y_1\oooo{e_1}^{(0)}+y_2\oooo{e_2}^{(0)} &\mapsto& 
y_1-\kappa_1y_2,
\end{eqnarray*}
is a bijection, where 
$\kappa_1=\frac{l}{2}+\frac{1}{2}\sqrt{l^2-4}$. 
The norm of $b_1-b_2\kappa_1$ is 
\begin{eqnarray*}
1=b_1^2-lb_1b_2+b_2^2.
\end{eqnarray*} 
Now
\begin{eqnarray*}
(2,-l)&=& j^{(0)}(e_1)(e_1,e_2)
=(-a_1j^{(0)}(\www{e_2})-a_2j^{(0)}(v_2))(e_1,e_2)\\
&=& -a_1((-l,2)+l(2,-l))-a_2(b_1(2,-l)+b_2(-l,2))\\
&=& (-a_1l-a_2b_1)(2,-l)+(-a_1-a_2b_2)(-l,2).
\end{eqnarray*}
so
\begin{eqnarray*}
a_1=-a_2b_2,\quad 1=-a_1l-a_2b_1=a_2(b_2l-b_1),\\
a_2=\pm 1\quad\textup{and}\quad b_1=-a_2+b_2l.
\end{eqnarray*}
Calculate 
\begin{eqnarray*}
0&=&-1+b_1^2-lb_1b_2+b_2^2=-1+(-a_2+b_2l)(-a_2)+b_2^2\\\
&=& b_2(b_2-a_2l).
\end{eqnarray*}
We obtain the four solutions
\begin{eqnarray*}
(a_1,a_2,b_1,b_2)&\in&\{(0,1,-1,0),(0,-1,1,0),\\
&&(-l,1,l^2-1,l),(-l,-1,1-l^2,-l)\}.
\end{eqnarray*}

In the case of the third solution 
\begin{eqnarray*}
\oooo{v_2}^{(0)}&=&(l^2-1)\oooo{e_1}^{(0)}+l\oooo{e_2}^{(0)}
=\oooo{s_{e_1}^{(0)}s_{e_2}^{(0)}(e_1)},\\
r&=&\oooo{s_{v_2}^{(0)}}
=\oooo{s_{s_{e_1}^{(0)}s_{e_2}^{(0)}(e_1)}^{(0)}}
=\oooo{s_{e_1}^{(0)}s_{e_2}^{(0)}s_{e_1}^{(0)}
s_{e_2}^{(0)}s_{e_1}^{(0)}}\\
&=& (z_1z_2)^2z_1=(z_1z_2)^{1-m}z_1\quad\textup{with }m=-1.
\end{eqnarray*}
As $m=-1$ is not in the set $\Z_{\geq 0}$, we can discard the
third solution. In fact, $(z_1z_2z_1,(z_1z_2)^2z_1,(z_1z_2)^2z_1)$
is in the $\Br_3$ orbit of $(z_1z_2z_1,z_1,z_1)$
because $\begin{pmatrix}-1\\0\end{pmatrix}$ is in the
$SL_2(\Z)$ orbit of $\begin{pmatrix}1\\0\end{pmatrix}$. 
We can discard also the fourth solution because its vector
$\oooo{v_2}^{(0)}$ differs from the vector $\oooo{v_2}^{(0)}$
in the third solution only by the sign.

Also the vector $\oooo{v_2}^{(0)}$ in the first solution differs
from the vector $\oooo{v_2}^{(0)}$ in the second solution 
only by the sign.

The second solution gives $\oooo{v_2}^{(0)}=\oooo{e_1}^{(0)}$
and thus for some $b_3\in\Z$
\begin{eqnarray*}
\uuuu{v}=(\www{e_2},e_1+b_3f_1,e_1+b_3f_1-f_1)
=(\www{e_2},e_1+b_3f_1,e_3+b_3f_1).
\end{eqnarray*}
The observation
\begin{eqnarray*}
\delta_2\sigma_2(\uuuu{v})
&=& \delta_2(\www{e_2},e_3+b_3f_1-2(e_1+b_3f_1),e_1+b_3f_1)\\
&=& (\www{e_2},e_1+(b_3+1)f_1,e_3+(b_3+1)f_1)
\end{eqnarray*}
shows $\uuuu{v}\in \Br_3\ltimes\{\pm 1\}^3(\www{e_2},e_1,e_3)$.
This orbit is $\BB^{dist}$ because 
$(\www{e_2},e_1,e_3)=\sigma_1(\uuuu{e})$. \hfill$\Box$ 

\bigskip
Lemma \ref{t7.6} states some facts which arise in the proof
of part (b) of Theorem \ref{t7.3} and which are worth 
to be formulated explicitly.

\begin{lemma}\label{t7.6}
Let $(H_\Z,L,\uuuu{e})$ be the unimodular bilinear lattice
of rank 3 with triangular basis $\uuuu{e}$ with matrix
$S=S(\HH_{1,2})=S(-2,2,-2)=L(\uuuu{e}^t,\uuuu{e})^t$.
Recall $\Rad I^{(0)}=\Z f_1+\Z f_2$ and
$R^{(0)}=\pm f_3+\Rad I^{(0)}$ (Theorem \ref{t6.14} (f)).

(a) For $g_1,g_2,g_3\in\Rad I^{(0)}$
\begin{eqnarray*}
s_{f_3-g_1}^{(0)}s_{f_3-g_2}^{(0)}s_{f_3-g_3}^{(0)}=-M
\iff g_2=g_1+g_3.
\end{eqnarray*}

(b) The map
\begin{eqnarray*}
\Phi:M_{2\times 2}(\Z)&\to& 
\{\uuuu{v}\in (R^{(0)})^3\,|\, (\pi_3\circ\pi_3^{(0)})(\uuuu{v})
=-M\}/\{\pm 1\}^3,\\
A &\mapsto& (f_3-g_1,f_3-g_1-g_3,f_3-g_3)/\{\pm 1\}^3\\
&&\textup{with }
\begin{pmatrix}g_1\\g_3\end{pmatrix}=
A \begin{pmatrix}f_1\\f_2\end{pmatrix},
\end{eqnarray*}
is a bijection. The action of $\Br_3$ on the right hand side
translates to the following action on the left hand side,
\begin{eqnarray*}
\sigma_1(A) = \begin{pmatrix}1&-1\\0&1\end{pmatrix} A,\quad
\sigma_2(A) = \begin{pmatrix}1&0\\1&1\end{pmatrix} A.
\end{eqnarray*}

(c) $\uuuu{v}\in (R^{(0)})^3$ with 
$(\pi_3\circ\pi_3^{(0)})(\uuuu{v})=-M$ satisfies either (i)
or (ii):
\begin{list}{}{}
\item[(i)]
There exists a permutation $\sigma\in S_3$ with
$v_i\in \Gamma^{(0)}\{e_{\sigma(i)}\}$
for $i\in\{1,2,3\}$.
\item[(ii)]
Either $v_1,v_2,v_3\in \Gamma^{(0)}\{f_3\}$ or
there exists a permutation $\sigma\in S_3$ and an $l\in\{1,2,3\}$
with $v_{\sigma(1)}\in \Gamma^{(0)}\{f_3\}$
and $v_{\sigma(2)},v_{\sigma(3)}\in \Gamma^{(0)}\{e_l\}$.
\end{list}
(i) holds if and only if $\Phi^{-1}(\uuuu{v}/\{\pm 1\}^3)$
has an odd determinant.

(d) Let $SL_2(\Z)$ act by multiplication from the left on
$\{A\in M_{2\times 2}(\Z)\,|\,\det A\textup{ is odd}\}$.
Each  orbit has a unique representative of the shape
\begin{eqnarray*}
\begin{pmatrix}0&c_2\\c_1&c_3\end{pmatrix}\quad
\textup{with}\quad c_1\in\N\textup{ odd}, 
c_2\in\Z\textup{ odd}, c_3\in\{0,1,...,|c_2|-1\}.
\end{eqnarray*}
\end{lemma}

{\bf Proof:} 
(a) For $g_1,g_2,g_3\in \Rad I^{(0)}$
\begin{eqnarray*}
s_{f_3-g_1}^{(0)}|_{\Rad I^{(0)}}=\id,\quad 
s_{f_3-g_1}^{(0)}(f_3+g_2)=-(f_3-g_2-2g_1),\\
s_{f_3-g_1}^{(0)}s_{f_3-g_2}^{(0)}s_{f_3-g_3}^{(0)}(f_3)
=-f_3+2(g_1-g_2+g_3).
\end{eqnarray*}
Compare $-M|_{\Rad I^{(0)}}=\id$, $-M(f_3)=-f_3$. 

(b) $\Phi$ is a bijection because of 
$R^{(0)}=\pm f_3+\Rad I^{(0)}$ and part (a).
The action of $\Br_3$ on the right hand side translates to the
claimed action on the left hand side because of the following,
\begin{eqnarray*}
\delta_1\sigma_1(f_3-g_1,f_3-g_1-g_3,f_3-g_3)
&=& (f_3-g_1+g_3,f_3-g_1,f_3-g_3),\\
\begin{pmatrix} g_1-g_3\\g_3\end{pmatrix}
&=& \begin{pmatrix} 1&-1\\0&1\end{pmatrix}
\begin{pmatrix} g_1\\g_3\end{pmatrix},\\
\delta_2\sigma_2(f_3-g_1,f_3-g_1-g_3,f_3-g_3)
&=& (f_3-g_1,f_3-2g_1-g_3,f_3-g_1-g_3),\\
\begin{pmatrix} g_1\\g_1+g_3\end{pmatrix}
&=& \begin{pmatrix} 1&0\\1&1\end{pmatrix}
\begin{pmatrix} g_1\\g_3\end{pmatrix}.
\end{eqnarray*}

(c) Recall that by Theorem \ref{t6.14} (e)
\begin{eqnarray*}
\Gamma^{(0)}\{e_i\}=\pm e_i+2\Rad I^{(0)},\quad
\Gamma^{(0)}\{f_3\}=\pm f_3+2\Rad I^{(0)},\\
(e_1,e_2,e_3)=(f_3-f_2,-f_3+f_1+f_2,f_3-f_1),\\
\Delta^{(0)}=\Gamma^{(0)}\{e_1\}\ \dot\cup\ 
\Gamma^{(0)}\{e_2\}\ \dot\cup\ \Gamma^{(0)}\{e_3\},\quad
R^{(0)}=\Delta^{(0)}\ \dot\cup\ \Gamma^{(0)}\{f_3\}.
\end{eqnarray*}
Observe that $g_1,g_2,g_3\in \Rad I^{(0)}$ with 
$g_2=g_1+g_3$ satisfy either (i)' or (ii)',
\begin{list}{}{}
\item[(i)']
$g_1,g_2,g_3\notin 2\Rad I^{(0)}$,
\item[(ii)']
There exists a permutation $\sigma\in S_3$ with 
$g_{\sigma(1)}\in 2\Rad I^{(0)}$ and 
$g_{\sigma(2)}-g_{\sigma(3)}\in 2\Rad I^{(0)}$.
\end{list}
$\uuuu{v}=(f_3-g_1,f_3-g_2,f_3-g_3)$ satisfies (i) if
(i)' holds, and it satisfies (ii) if (ii)' holds. 
If $\begin{pmatrix}g_1\\g_2\end{pmatrix}=A
\begin{pmatrix}f_1\\f_2\end{pmatrix}$ then (i)' holds
if and only if $\det A$ is odd.

(d) This is elementary. We skip the details. \hfill$\Box$

\bigskip
{\bf Proof of Theorem \ref{t7.3} (b):}
$\uuuu{v}\in (\Delta^{(0)})^3$ with 
$(\pi_3\circ\pi_3^{(0)})(\uuuu{v})=-M$ satisfies property (i)
in Lemma \ref{t7.6} (c) because it does not satisfy property 
(ii) in Lemma \ref{t7.6} (c). 
The parts (b) and (d) of
Lemma \ref{t7.6} show that the set
$(\Delta^{(0)})^3\cap (\pi_3\circ\pi_3^{(0)})^{-1}(-M)$
consists of countably many orbits. The parts (b) and (d)
of Lemma \ref{t7.6} also give the claimed representative 
in each orbit. The rest is obvious.
\hfill$\Box$

\section{Distinguished bases
in the odd rank 3 cases}
\label{s7.3}

Also in the odd cases with $n=3$ we have complete results on the 
question when the inclusion \eqref{3.4} is an equality.
It is one if and only if $\uuuu{x}\in B_1\cup B_2$ where
\index{$B_1,\ B_2,\ B_3,\ B_4,\ B_5$}
\begin{eqnarray*}
B_1&=& \{\uuuu{x}\in\Z^3-\{(0,0,0)\}\,| \\
&&\hspace*{1.5cm}
((G^{phi}\ltimes \www{G}^{sign})\rtimes\langle \gamma\rangle)
(\uuuu{x})\cap r^{-1}(\Z_{\leq 0})\neq\emptyset\},\\
B_2&=& \{\uuuu{x}\in\Z^3-\{(0,0,0)\}\,|\, 
S(\uuuu{x})\textup{ is reducible, i.e. there are }i,j,k\\
&& \hspace*{1cm}\textup{ with }
\{i,j,k\}=\{1,2,3\}\textup{ and }x_i\neq 0=x_j=x_k\},\\
B_3&:=& \{(0,0,0)\},\\
B_4&:=& \{\uuuu{x}\in\Z^3\,|\, S(\uuuu{x})\in 
(\Br_3\ltimes\{\pm 1\}^3)\Bigl(
\{S(A_3),S(\whh{A}_2),S(\HH_{1,2})\}\Bigr.\\
&&\Bigl.\hspace*{6cm}\cup 
\{S(-l,2,-l)\,|\, l\geq 3\}\Bigr)\},\\
B_5&:=& \Z^3-(B_1\cup B_2\cup B_3\cup B_4).
\end{eqnarray*}
$B_2$ is the set of $\uuuu{x}\neq (0,0,0)$ which give 
reducible cases.
$B_1$ contains $r^{-1}(\Z_{\leq 0})-\{(0,0,0)\}$, but is bigger.
$\uuuu{x}\in B_1$ if and only if the 
$(G^{phi}\ltimes \www{G}^{sign})\rtimes\langle \gamma\rangle$
orbit of $\uuuu{x}$ contains a triple $\www{\uuuu{x}}\in\Z^3$
as in Lemma \ref{t4.18} (a), so with 
$\www{\uuuu{x}}\in\Z^3_{\geq 3}$ and 
$2\www{x}_i\leq \www{x}_j\www{x}_k$ for $\{i,j,k\}=\{1,2,3\}$.
The Examples \ref{t4.20}  show that it is not so easy to
describe $B_1$ more explicitly.

In Theorem \ref{t7.7} we show 
$B_4\subset \Z^3-(B_1\cup B_2\cup B_3)$.
For $\uuuu{x}\in B_3\cup B_4$ 
the inclusion in \eqref{3.4} is not an equality, but
we can add the constraint
$\sum_{i=1}^3\Z v_i= H_\Z$ to \eqref{3.4} and obtain an equality.
For $\uuuu{x}\in B_5$ we do not know whether \eqref{3.4}
with the additional constraint $\sum_{i=1}^3\Z v_i= H_\Z$
becomes an equality.

\begin{theorem}\label{t7.7}
Let $(H_\Z,L,\uuuu{e})$ be a unimodular bilinear lattice of rank
$3$ with a triangular basis $\uuuu{e}$ with matrix
$L(\uuuu{e}^t,\uuuu{e})^t=S(\uuuu{x})\in T^{uni}_3(\Z)$
for some $\uuuu{x}\in \Z^3$.

(a) $(H_\Z,L,\uuuu{e})$ is reducible if and only if 
$\uuuu{x}\in B_2\cup B_3$. Then $\Gamma^{(1)}_u=\{\id \}$.

(b) The following conditions are equivalent:
\begin{list}{}{}
\item[(i)]
$\uuuu{x}\in B_1$.
\item[(ii)]
$\Gamma^{(1)}\cong G^{free,3}$.
\item[(iii)]
$(H_\Z,L,\uuuu{e})$ is irreducible and $\Gamma^{(1)}_u=\{\id\}$.
\end{list}

(c) $\Z^3=\dot\bigcup_{i\in\{1,2,3,4,5\}} B_i$.

(d) The inclusion in \eqref{3.4} is an equality 
$\iff \uuuu{x}\in B_1\cup B_2.$

(e) Consider $\uuuu{x}=(0,0,0)$. The set
$\{\uuuu{v}\in (\Delta^{(1)})^3\,|\, 
(\pi_3\circ\pi_3^{(1)})(\uuuu{v})=M\}$ is $(\Delta^{(1)})^3$. 
It consists of ten
$\Br_3\ltimes\{\pm 1\}^3$ orbits, the orbit $\BB^{dist}$
of $\uuuu{e}=(e_1,e_2,e_3)$ and the orbits of the nine 
triples
\begin{eqnarray*}
(e_1,e_1,e_1),(e_1,e_1,e_2),(e_1,e_2,e_2),(e_2,e_2,e_2),\\
(e_1,e_1,e_3),
(e_1,e_3,e_3),(e_3,e_3,e_3),(e_2,e_2,e_3),(e_2,e_3,e_3).
\end{eqnarray*}

(f) Consider $\uuuu{x}\in B_4\cup B_5$. 
Then $\Gamma^{(1)}_u\cong \Z^2$. The map
\begin{eqnarray*}
\Psi: \{\uuuu{v}\in (\Delta^{(1)})^3\,|\, 
(\pi_3\circ\pi_3^{(1)})(\uuuu{v})=M\} &\to& \N\cup\{\infty\},\\
\uuuu{v}&\mapsto& \Bigl(\textup{index of }\sum_{i=1}^3\Z v_i
\textup{ in }H_\Z\Bigr),
\end{eqnarray*}
has infinitely many values.
The set $\{\uuuu{v}\in (\Delta^{(1)})^3\,|\, 
(\pi_3\circ\pi_3^{(1)})(\uuuu{v})=M\}$ contains besides
$\BB^{dist}$ infinitely many $\Br_3\ltimes\{\pm 1\}^3$ orbits.

(g) For $\uuuu{x}\in B_3\cup B_4$ 
\begin{eqnarray*}
\BB^{dist}= \{\uuuu{v}\in (\Delta^{(1)})^3\,|\, 
(\pi_3\circ\pi_3^{(1)})(\uuuu{v})=M,\sum_{i=1}^3\Z v_i =H_\Z\}.
\end{eqnarray*}
\end{theorem}

{\bf Proof:}
(a) Compare Definition \ref{t2.10} in the case $n=3$.
For $\Gamma^{(1)}_u=\{\id\}$ see Theorem \ref{t6.18} (b).

(b) By the Remarks \ref{t4.17} the tuple 
$(H_\Z,\pm I^{(1)},\Gamma^{(1)},\Delta^{(1)})$ depends up to 
isomorphism only on the 
$(G^{phi}\ltimes \www{G}^{sign})\rtimes\langle \gamma\rangle$
orbit of $\uuuu{x}$. 
Lemma \ref{t4.18} gives representatives of each such orbit.
Theorem \ref{t6.18} studies their groups $\Gamma^{(1)}$.
Theorem \ref{t6.18} (b) treats $\uuuu{x}\in B_2\cup B_3$.
Theorem \ref{t6.18} (g) treats $\uuuu{x}\in B_1$.
Theorem \ref{t6.18} (c)--(f) treats 
$\uuuu{x}\in \Z^3-(B_1\cup B_2\cup B_3)$.

One sees 
\begin{eqnarray*}
\Gamma^{(1)}\cong G^{free,3}&\iff& \uuuu{x}\in B_1,\\
\Gamma^{(1)}_u=\{\id\}&\iff& \uuuu{x}\in B_1\cup B_2\cup B_3,\\
\Gamma^{(1)}_u\cong\Z^2&\iff& \uuuu{x}\in \Z^3 -(B_1\cup B_2
\cup B_3).
\end{eqnarray*}

(c) $B_1\cap B_2=\emptyset$, $B_2\cap B_4=\emptyset$ and
$(0,0,0)\notin B_1\cup B_2\cup B_4$ are clear.
$B_1\cap B_4=\emptyset$ follows from 
$\Gamma^{(1)}(\uuuu{x})\not\cong G^{free,3}$ for
$\uuuu{x}\in B_4$. 

(d) The parts (e) and (f) will give $\Longrightarrow$.
Here we show $\Longleftarrow$, first for $\uuuu{x}\in B_1$,
then for $\uuuu{x}\in B_2$.

Let $\uuuu{x}\in B_1$. Then by the Remarks \ref{t4.17} and
Theorem \ref{t6.18} (g) $\Gamma^{(1)}$ is a free group with
generators $s_{e_1}^{(1)}$, $s_{e_2}^{(1)}$ and $s_{e_3}^{(1)}$.
Example \ref{t3.23} (iv) applies. 

Let $\uuuu{x}\in B_2$. Because of the actions of $\gamma$
and $\www{G}^{sign}$ (here $G^{sign}$ is sufficient) on $B_2$
we can suppose $\uuuu{x}=(x_1,0,0)$ with $x_1\in \Z_{<0}$.
Then $e_3\in\Rad I^{(1)}$, $s_{e_3}^{(1)}=\id$,
$\Gamma^{(1)}=\langle s_{e_1}^{(1)},s_{e_2}^{(1)}\rangle$ and
by Theorem \ref{t6.21} (b)
$\Delta^{(1)}=\Delta^{(1)}\cap (\Z e_1+\Z e_2)\ \cup\ \{\pm e_3\}$.
The monodromy $M$ has the characteristic polynomial
$(t-1)(t^2-(2-r(\uuuu{x}))t+1)=(t-1)(t^2-(2-x_1^2)t+1)$, so 
three different eigenvalues.

Consider $\uuuu{v}\in (\Delta^{(1)})^3 $ with 
$s_{v_1}^{(1)}\circ s_{v_2}^{(1)}\circ s_{v_3}^{(1)}=M$.
Now $v_1,v_2,v_3\in\{\pm e_3\}$ is impossible because 
$M\neq \id$. Two of $v_1,v_2,v_3$ in $\{\pm e_3\}$
cannot be because then $M$ would have the eigenvalue $1$ with
multiplicity 3. 

\medskip
{\bf Claim:} All three $v_1,v_2,v_3\in\Delta^{(1)}\cap 
(\Z e_1+\Z e_2)$ is impossible.

\medskip
{\bf Proof of the Claim:}
Suppose $v_1,v_2,v_3\in \Delta^{(1)}\cap (\Z e_1+\Z e_2)$. 
First we consider a case with $x_1\leq -2$.
By Theorem \ref{t6.18} (b) and Theorem \ref{t6.10} (c)+(d) 
$\Gamma^{(1)}\cong G^{free,2}$ with generators 
$s_{e_1}^{(1)}$ and $s_{e_2}^{(1)}$. There is a unique group
homomorphism
\begin{eqnarray*}
\Gamma^{(1)}\to\{\pm 1\}\quad \textup{with}\quad
s_{e_1}^{(1)}\mapsto -1,\ s_{e_2}^{(1)}\mapsto -1.
\end{eqnarray*}
Each $s_{v_i}^{(1)}$ is conjugate to $s_{e_1}^{(1)}$ or
$s_{e_2}^{(1)}$ and thus has image $-1$. Also their product 
$s_{v_1}^{(1)}\circ s_{v_2}^{(1)}\circ s_{v_3}^{(1)}$ 
has image $-1$. But $M=s_{e_1}^{(1)}\circ s_{e_2}^{(1)}$ 
has image $1$, a contradiction.

Now consider the case $x_1=-1$. By Theorem \ref{t6.18} (b)
and Theorem \ref{t6.10} (a)+(b) $\Gamma^{(1)}\cong SL_2(\Z)$
with $s_{e_1}^{(1)}\sim \begin{pmatrix}1&1\\0&1\end{pmatrix}$
and $s_{e_2}^{(1)}\sim \begin{pmatrix}1&0\\-1&1\end{pmatrix}$.
It is well known that the group $SL_2(\Z)$ is isomorphic to
the group with the presentation
\begin{eqnarray*}
\langle x_1,x_2\,|\, x_1x_2x_1=x_2x_1x_2,\ 1=(x_1x_2)^6\rangle
\end{eqnarray*}
where $x_1\mapsto \begin{pmatrix}1&1\\0&1\end{pmatrix}$
and $x_2\mapsto \begin{pmatrix}1&0\\-1&1\end{pmatrix}$. 
The differences of the lengths of the words in 
$x_1^{\pm 1}$ and $x_2^{\pm 1}$ which are connected by these
relations are $3-3$ and $12-0$, so even. Therefore also in this
situation there is a unique group homomorphism
\begin{eqnarray*}
\Gamma^{(1)}\to\{\pm 1\}\quad \textup{with}\quad
s_{e_1}^{(1)}\to -1,\ s_{e_2}^{(1)}\to -1.
\end{eqnarray*}
The argument in the case 
$\Gamma^{(1)}\cong G^{free,2}$ goes here through, too. 
The Claim is proved. \hfill ($\Box$)

\medskip
Therefore a permutation $\sigma\in S_3$ with
$v_{\sigma(1)},v_{\sigma(2)}\in\Delta^{(1)}\cap 
(\Z e_1+\Z e_2)$, $\sigma(1)<\sigma(2)$ and 
$v_{\sigma(3)}\in\{\pm e_3\}$ exists. Then 
$s_{v_{\sigma(3)}}^{(1)}=\id$ and
\begin{eqnarray*}
s_{e_1}^{(1)}s_{e_2}^{(1)}=M=s_{v_1}^{(1)}s_{v_2}^{(1)}
s_{v_3}^{(1)}=s_{v_{\sigma(1)}}^{(1)}s_{v_{\sigma(2)}}^{(1)}.
\end{eqnarray*}
Because of Theorem \ref{t7.1} $(v_{\sigma(1)},v_{\sigma(2)})$
is in the $\Br_2\ltimes\{\pm 1\}^2$ orbit of $(e_1,e_2)$.
Therefore $\uuuu{v}$ is in 
$(\Br_3\ltimes\{\pm 1\}^3)(\uuuu{e})=\BB^{dist}$.

(e) See the Examples \ref{t3.23} (iii).

(f) $\Gamma^{(1)}_u\cong\Z^2$ for $\uuuu{x}\in B_4\cup B_5$
follows from Theorem \ref{t6.21} (c)--(f), the Remarks \ref{t4.17},
Lemma \ref{t4.18} and the definition of $B_4$ and $B_5$. 
The last statement in (f) follows from the middle statement
because the sublattice $\sum_{i=1}^3\Z v_i\subset H_\Z$ and
its index in $H_\Z$ are invariants of the 
$\Br_3\ltimes\{\pm 1\}^3$ orbit of $\uuuu{v}$. 

For the middle statement we consider
\begin{eqnarray*}
\uuuu{v}=(e_1+a(-\www{x}_3)f_3,e_2+a(\www{x}_2-\www{x}_1x_3)f_3,
e_3+a(-\www{x}_1)f_3)\quad\textup{with }a\in\Z.
\end{eqnarray*}
The next Lemma \ref{t7.8} implies
\begin{eqnarray*}
(\pi_3\circ\pi_3^{(1)})(\uuuu{v})&=&M,\\
\Bigl(\textup{index of }\sum_{i=1}^3\Z v_i\textup{ in }H_\Z\Bigr)
&=&\left|1+a\frac{r(\uuuu{x})}{\gcd(x_1,x_2,x_3)^2}\right|
\end{eqnarray*}
and that for a suitable $a_0\in\N$ and any $a\in\Z a_0$
$\uuuu{v}\in(\Delta^{(1)})^3$. As $r(\uuuu{x})\neq 0$
for $\uuuu{x}\in B_4\cup B_5$, this shows that the map
$\Psi$ has countably many values.

(g) Part (g) will be prepared by Lemma \ref{t7.10} and 
will be proved after the proof of Lemma \ref{t7.10}.
\hfill$\Box$

\begin{lemma}\label{t7.8}
Let $(H_\Z,L,\uuuu{e})$ be a unimodular bilinear lattice of rank 3
with a triangular basis $\uuuu{e}$ with matrix 
$L(\uuuu{e}^t,\uuuu{e})^t=S(\uuuu{x})$ for some
$\uuuu{x}\in \Z^3-\{(0,0,0)\}$. 
Recall $\Rad I^{(1)}=\Z f_3$ with
$f_3=-\www{x}_3e_1+\www{x}_2e_2-\www{x}_1e_3$ and
$(\www{x}_1,\www{x}_2,\www{x}_3)
=\gcd(x_1,x_2,x_3)^{-1}(x_1,x_2,x_3).$

(a) For $\uuuu{a}=(a_1,a_2,a_3)\in\Z^3$
\begin{eqnarray*}
&&s_{e_1+a_1f_3}^{(1)}\circ s_{e_2+a_2f_3}^{(1)}\circ 
s_{e_3+a_3f_3}^{(1)} =M\\
&\iff& (a_1,a_2,a_3)\in \Z(-\www{x}_3,\www{x}_2-\www{x}_1x_3,
-\www{x}_1).
\end{eqnarray*}

(b) For $\uuuu{a}=(a_1,a_2,a_3)=a
(-\www{x}_3,\www{x}_2-\www{x}_1x_3,-\www{x}_1)$ with $a\in\Z$,
the index of $\sum_{i=1}^3\Z(e_i+a_if_3)$ in $H_\Z$ is
$|1+a\frac{r(\uuuu{x})}{\gcd(x_1,x_2,x_3)^2}|$.

(c) If $\Gamma^{(1)}_u\cong\Z^2$ then there is a number
$a_0\in\N$ with
\begin{eqnarray*}
(e_1+a(-\www{x}_3)f_3,e_2+a(\www{x}_2-\www{x}_1x_3)f_3,
e_3+a(-\www{x}_1)f_3)\in (\Delta^{(1)})^3
\textup{ for }a\in \Z a_0.
\end{eqnarray*}
\end{lemma}

{\bf Proof:}
(a) With Lemma \ref{t7.4} and 
$f_3=(f\textup{ in Lemma \ref{t7.4}})$ one calculates
\begin{eqnarray*}
&&s_{e_1+a_1f_3}^{(1)}\circ s_{e_2+a_2f_3}^{(1)}\circ 
s_{e_3+a_3f_3}^{(1)} \circ M^{-1}\\
&=& t_{-a_1j^{(1)}(e_1)}\circ s_{e_1}^{(1)}
\circ t_{-a_2j^{(1)}(e_2)}\circ s_{e_2}^{(1)}
\circ t_{-a_3j^{(1)}(e_3)}\circ s_{e_3}^{(1)}\circ M^{-1}\\
&=& 
t_{-a_1j^{(1)}(e_1)}\circ 
t_{-a_2j^{(1)}(e_2)-a_2j^{(1)}(e_2)(e_1)j^{(1)}(e_1)}\circ
s_{e_1}^{(1)}\\
&&\circ \ 
t_{-a_3j^{(1)}(e_3)-a_3j^{(1)}(e_3)(e_2)j^{(1)}(e_2)}\circ 
s_{e_2}^{(1)}\circ s_{e_3}^{(1)}\circ M^{-1}\\
&=&
t_{-A}
\end{eqnarray*}
with
\begin{eqnarray*}
A&=& a_1j^{(1)}(e_1)+a_2j^{(1)}(e_2)
+a_2I^{(1)}(e_2,e_1)j^{(1)}(e_1)\\
&&+a_3j^{(1)}(e_3)+a_3I^{(1)}(e_3,e_2)j^{(1)}(e_2)
+a_3I^{(1)}(e_3,e_1)j^{(1)}(e_1)\\
&&+a_3I^{(1)}(e_3,e_2)I^{(1)}(e_2,e_1)j^{(1)}(e_1)\\
&=& j^{(1)}\Bigl((a_1-a_2x_1-a_3x_2+a_3x_1x_3)e_1
+(a_2-a_3x_3)e_2+a_3e_3\Bigr).
\end{eqnarray*}
$t_{-A}=\id$ holds if and only if $A=0$, so if and only if
\begin{eqnarray*}
(a_1-a_2x_1-a_3x_2+a_3x_1x_3)e_1+(a_2-a_3x_3)e_2+a_3e_3
\in \Rad I^{(1)}=\Z f_3.
\end{eqnarray*}
The ansatz that it is 
$af_3=a(-\www{x}_3e_1+\www{x}_2e_2-\www{x}_1e_3)$ with $a\in\Z$ 
gives
\begin{eqnarray*}
-a\www{x}_1=a_3,\ a\www{x}_2=a_2-a_3x_3,\ 
-a\www{x}_3=a_1-a_2x_1-a_3x_2+a_3x_1x_3,\\
\textup{so}\quad (a_1,a_2,a_3)=a(-\www{x}_3,\www{x}_2-\www{x}_1x_3,
-\www{x}_1).
\end{eqnarray*}

(b) Write $\uuuu{a}=\www{a}\gcd(x_1,x_2,x_3)(-x_3,x_2-x_1x_3,-x_1)$ with
$\www{a}=\gcd(x_1,x_2,x_3)^{-2}a\in \gcd(x_1,x_2,x_3)^{-2}\Z$.
Then
\begin{eqnarray*}
&&(e_1+a_1f_3,e_2+a_2f_3,e_3+a_3f_3)\\
&=&\uuuu{e}\begin{pmatrix} 
1+\www{a}(-x_3)(-x_3) & \www{a}(x_2-x_1x_3)(-x_3) & 
\www{a}(-x_1)(-x_3) \\
\www{a}(-x_3)x_2 & 1+\www{a}(x_2-x_1x_3)x_2 & 
\www{a}(-x_1)x_2 \\
\www{a}(-x_3)(-x_1) & \www{a}(x_2-x_1x_3)(-x_1) & 
1+\www{a}(-x_1)(-x_1)\end{pmatrix}.
\end{eqnarray*}
The determinant of this matrix is $1+\www{a}r(\uuuu{x})$.
The index of the lattice $\sum_{i=1}^3\Z(e_i+a_if_3)$ in
$H_\Z$ is the absolute value of this determinant.

(c) Suppose $\Gamma^{(1)}_u\cong\Z^2$.
Compare Lemma \ref{t6.17}. The set
\begin{eqnarray*}
\Lambda:=\{\lambda\in \Hom_0(H_\Z,\Z)\,|\, 
t_\lambda^+\in \Gamma^{(1)}_u\}
\end{eqnarray*}
is a sublattice of rank 2 in the lattice $\Hom_0(H_\Z,\Z)$
of rank 2. For $i\in\{1,2,3\}$ 
\begin{eqnarray*}
\Gamma^{(1)}_u\{e_i\}=\{e_i+\lambda(e_i)f_3\,|\, 
\lambda\in\Lambda\}
\subset (e_i+\Z f_3)\cap\Delta^{(1)}.
\end{eqnarray*}
The triple $(H_\Z,L,\uuuu{e})$ is irreducible because of
$\Gamma_u^{(1)}\cong\Z^2$ and Theorem \ref{t7.7} (a). 
Therefore it is not reducible with a summand of type $A_1$,
and thus $\{e_1,e_2,e_3\}\cap\Rad I^{(1)}=\emptyset$.
Because of this and because $\Lambda$ has finite index in
$\Hom_0(H_\Z,\Z)$, there is a number $b_i\in\N$ with
\begin{eqnarray*}
\Z b_i= \{b\in\Z\,|\, e_i+bf_3\in\Gamma^{(1)}_u\{e_i\}\}.
\end{eqnarray*}
For each $\uuuu{a}=(a_1,a_2,a_3)$ with 
$a_i\in \Z b_i$ $e_i+a_if_3\in\Delta^{(1)}$. 
Any number $a_0\in\N$ (for example the smallest one) with
\begin{eqnarray*}
a_0\www{x}_3\in \Z b_1,\ a_0(\www{x}_2-\www{x}_1x_3)\in\Z b_2,\ 
a_0\www{x}_1\in \Z b_3
\end{eqnarray*}
works. \hfill$\Box$

\begin{remarks}\label{t7.9}
If $\uuuu{x}\in B_1\cup B_2$ then the map 
$\Delta^{(1)}\to \oooo{\Delta^{(1)}}$ is a bijection by
Theorem \ref{t6.18} (b)+(g). Therefore then
\begin{eqnarray*}
\uuuu{v}=(e_1+a(-\www{x}_3)f_3,e_2+a(\www{x}_2-\www{x}_1x_3)f_3,
e_3+a(-\www{x}_1)f_3)\in H_\Z^3
\end{eqnarray*}
for $a\in\Z-\{0\}$ satisfies 
$(\pi_3\circ\pi_3^{(1)})(\uuuu{v})=M$,
but $\uuuu{v}\notin (\Delta^{(1)})^3$. This fits to 
Theorem \ref{t7.7} (d).
\end{remarks}

\begin{lemma}\label{t7.10}
The $\Br_3\ltimes\{\pm 1\}^3$ orbits in $r^{-1}(4)\subset\Z^3$
are classified in Theorem \ref{t4.6} (e). They are separated
by the isomorphism classes of the pairs 
$(\oooo{H_\Z}^{(1)},\oooo{M})$ for corresponding unimodular
bilinear lattices $(H_\Z,L,\uuuu{e})$ with triangular bases
$\uuuu{e}$ with $L(\uuuu{e}^t,\uuuu{e})^t=S(\uuuu{x})$
and $r(\uuuu{x})=4$. More precisely, $\oooo{H_\Z}^{(1)}$
has a $\Z$-basis $(c_1,c_2)$ with 
$\oooo{M}(c_1,c_2)=(c_1,c_2)
\begin{pmatrix}-1&\gamma\\0&-1\end{pmatrix}$ with a unique
$\gamma\in\Z_{\geq 0}$, which is as follows:
\begin{eqnarray*}
\begin{array}{c|c|c|c|c|c}
S(\uuuu{x}) & S(\HH_{1,2}) & S(\P^1A_1) & S(\whh{A}_2) & 
S(-l,2,-l) & S(-l,2,-l)\\
 & & & & \textup{ with }l\equiv 0(2) & 
 \textup{ with }l\equiv 1(2) \\ \hline 
\gamma & 0 & 2 & 3 & \frac{l^2}{2}-2 & l^2-4\end{array}
\end{eqnarray*}
The numbers $\gamma$ in this table are pairwise different.
\end{lemma}

{\bf Proof:}
$S(\HH_{1,2})$: See Theorem \ref{t5.14} (a) (i).

$S(\P^1A_1)$: By Theorem \ref{t5.13}
$(\oooo{H_\Z}^{(1)},\oooo{M})\cong (H_{\Z,1},M_1)$ which comes
from $S(\P^1)=\begin{pmatrix}1&-2\\0&1\end{pmatrix}$ with
$S(\P^1)^{-1}S(\P^1)^t=\begin{pmatrix}-3&2\\-2&1\end{pmatrix}$.
This monodromy matrix is conjugate to 
$\begin{pmatrix}-1&2\\0&-1\end{pmatrix}$ with respect to
$GL_2(\Z)$.

$S(\whh{A}_2)$: Compare Theorem \ref{t5.14} (b) (iii) and its
proof:
\begin{eqnarray*}
f_3&=&e_1-e_2+e_3,\\
\oooo{H_\Z}^{(1)}&=&\Z \oooo{e_1}^{(1)}+\Z \oooo{e_2}^{(1)},\\
M\uuuu{e}&=&\uuuu{e}
\begin{pmatrix}-2&-1&2\\-2&0&1\\-1&-1&1\end{pmatrix},\\
\oooo{M}(\oooo{e_1}^{(1)},\oooo{e_2}^{(1)})
&=&(\oooo{e_1}^{(1)},\oooo{e_2}^{(1)})
\begin{pmatrix}-1&0\\-3&-1\end{pmatrix}.
\end{eqnarray*}
This monodromy matrix is conjugate to
$\begin{pmatrix}-1&3\\0&-1\end{pmatrix}$ with respect to
$GL_2(\Z)$. 

$S(-l,2,-l)$ with $l\geq 4,l\equiv 0(2)$: 
See Theorem \ref{t5.14} (a) (ii). Here
$(\oooo{H_\Z}^{(1)},\oooo{M})\cong (H_{\Z,1},M_1)$. 

$S(-l,2,-l)$ with $l\geq 3,l\equiv 1(2)$: 
Compare Theorem \ref{t5.14} (b) (iv) and its proof.
Define elements $a_1,a_2\in H_\Z$,
\begin{eqnarray*}
a_1&:=& \frac{l+1}{2}e_1+e_2+\frac{l+1}{2}e_3
= \frac{1}{2}f_1+\frac{1}{2}f_3,\\
a_2&:=& -e_1
= \frac{1}{2}\www{f}_2-\frac{l^2}{4}f_1-\frac{l}{4}f_3.
\end{eqnarray*}
The triple $(a_1,a_2,f_3)$ is a $\Z$-basis of $H_\Z$. 
The equality in the proof of Theorem \ref{t5.14} (b) (iv),
$$M(f_1,\www{f}_2)=(f_1,\www{f}_2)
\begin{pmatrix}-1&l^2-4\\0&-1\end{pmatrix},$$
implies
$$M(\oooo{a_1}^{(1)},\oooo{a_2}^{(1)})
=(\oooo{a_2}^{(1)},\oooo{a_2}^{(1)})
\begin{pmatrix}-1&l^2-4\\0&-1\end{pmatrix}.\hspace*{1cm}\Box$$

{\bf Proof of Theorem \ref{t7.7} (g):}
The case $\uuuu{}(0,0,0)$ is treated first and separately.
Compare part (e). Of the ten triples listed there, only
the triple $\uuuu{v}=(e_1,e_2,e_3)$ satisfies
$\sum_{i=1}^3\Z v_i =H_\Z$. This shows part (g) in the 
case $\uuuu{x}=(0,0,0)$.

Now consider $\uuuu{x}\in B_4$. We can suppose
\begin{eqnarray*}
\uuuu{x}\in\{(-1,0,-1),(-1,-1,-1),(-2,2,-2)\}
\cup\{(-l,2,-l)\,|\, l\geq 3\},
\end{eqnarray*}
which are the cases $S(\uuuu{x})\in \{S(A_3),S(\whh{A}_2),
S(\HH_{1,2})\}\cup\{S(-l,2,-l)\,|\, l\geq 3\}$. 
Consider $\uuuu{v}\in (\Delta^{(1)})^3$ with
$(\pi_3\circ\pi_3^{(1)})(\uuuu{v})=M$ and 
$\uuuu{v}$ a $\Z$-basis of $H_\Z$. We want to show 
$\uuuu{v}\in \BB^{dist}$. We have
\begin{eqnarray*}
I^{(1)}(\uuuu{v}^t,\uuuu{v})=\begin{pmatrix}
0&y_1&y_2\\-y_1&0&y_3\\-y_2&-y_3&0\end{pmatrix}
=S(\uuuu{y})-S(\uuuu{y})^t\quad\textup{for some }
\uuuu{y}\in\Z^3.
\end{eqnarray*}
Define a new Seifert form $\www{L}:H_\Z\times H_\Z\to \Z$
by $\www{L}(\uuuu{v}^t,\uuuu{v})^t=S(\uuuu{y})$
(only at the end of the proof it will turn out that
$\www{L}=L$). Then
\begin{eqnarray*}
\www{L}^t-\www{L}=I^{(1)}=L^t-L.
\end{eqnarray*}
$\uuuu{v}$ is a triangular basis with respect to
$(H_\Z,\www{L})$. By Theorem \ref{t2.7} for 
$(H_\Z,\www{L})$ (alternatively, one can calculate the product
of the matrices of $s_{v_1}^{(1)}$, $s_{v_2}^{(1)}$ and
$s_{v_3}^{(1)}$ with respect to $\uuuu{v}$) 
\begin{eqnarray*}
M\uuuu{v}=(\pi_3\circ\pi_3^{(1)})(\uuuu{v})=
(s_{v_1}^{(1)}\circ s_{v_2}^{(1)}\circ s_{v_3}^{(1)})
(\uuuu{v})=\uuuu{v}S(\uuuu{y})^{-1}S(\uuuu{y})^t.
\end{eqnarray*}
Then
\begin{eqnarray*}
3-r(\uuuu{x})=\tr(M)=\tr(S(\uuuu{y})^{-1}S(\uuuu{y})^t)
=3-r(\uuuu{y}),
\end{eqnarray*}
so $r(\uuuu{x})=r(\uuuu{y})$. 
In the case of $A_3$, $r^{-1}(2)$ is a unique 
$\Br_3\ltimes\{\pm 1\}^3$ orbit.
In the cases of $\whh{A}_2$, $\HH_{1,2}$ and 
$\uuuu{x}\in \{(-l,2,-l)\,|\, l\geq 3\}$, 
Lemma \ref{t7.10} and 
$M\uuuu{v}=\uuuu{v}S(\uuuu{y})^{-1}S(\uuuu{y})^t$
show that $\uuuu{y}$ is in the same $\Br_3\ltimes\{\pm 1\}^3$
orbit as $\uuuu{x}$. Therefore in any case there is an element
of $\Br_3\ltimes\{\pm 1\}^3$ which maps $\uuuu{v}$ to
a $\Z$-basis $\uuuu{w}$ of $H_\Z$ with
\begin{eqnarray*}
\www{L}(\uuuu{w}^t,\uuuu{w})^t=S(\uuuu{x}).
\end{eqnarray*}
Then
\begin{eqnarray*}
I^{(1)}(\uuuu{w}^t,\uuuu{w})=S(\uuuu{x})-S(\uuuu{x})^t
=I^{(1)}(\uuuu{e}^t,\uuuu{e}).
\end{eqnarray*}
Define an automorphism $g\in O^{(1)}$ by $g(\uuuu{e})=\uuuu{w}$.
Because of
\begin{eqnarray*}
M&=& s_{v_1}^{(1)}s_{v_2}^{(1)}s_{v_2}^{(1)}
= s_{w_1}^{(1)}s_{w_2}^{(1)}s_{w_3}^{(1)}\\
&=& s_{g(e_1)}^{(1)}s_{g(e_2)}^{(1)}s_{g(e_3)}^{(1)}
=gs_{e_1}^{(1)}s_{e_2}^{(1)}s_{e_3}^{(1)}g^{-1}
=gMg^{-1}
\end{eqnarray*}
$g$ is in $G_\Z^{(1)}$. But for the considered cases of
$\uuuu{x}$
\begin{eqnarray*}
G_\Z^{(1)}\stackrel{\textup{Theorem \ref{t5.14}}}{=}
G_\Z \stackrel{\textup{Theorem \ref{t3.28}}}{=}
Z(\Br_3\ltimes \{\pm 1\}^3).
\end{eqnarray*}
Therefore there is an element of $\Br_3\ltimes \{\pm 1\}^3$
which maps $\uuuu{e}$ to $\uuuu{w}$.
Altogether $\uuuu{v}\in \Br_3\ltimes\{\pm 1\}^3(\uuuu{e})
=\BB^{dist}$.
(Now also $L(\uuuu{v}^t,\uuuu{v})^t\in T^{uni}_3(\Z)$ and
thus $L=\www{L}$ are clear.)\hfill$\Box$

\section[The stabilizers of distinguished bases]
{The stabilizers of distinguished bases in the rank 3 cases}
\label{s7.4}

Let $(H_\Z,L,\uuuu{e})$ be a unimodular bilinear lattice
of rank $3$ with a triangular basis $\uuuu{e}$ with
$L(\uuuu{e}^t,\uuuu{e})^t=S(\uuuu{x})\in T^{uni}_3(\Z)$
for some $x\in\Z^3$. We are interested in the stabilizer
$(\Br_3)_{\uuuu{e}/\{\pm 1\}^3}$. The surjective map
\begin{eqnarray*}
\BB^{dist}=(\Br_3\ltimes\{\pm 1\}^3)(\uuuu{e})&\to&
(\Br_3\ltimes\{\pm 1\}^3)(\uuuu{x})\\
\www{\uuuu{e}}&\mapsto& \www{\uuuu{x}}\ \textup{ with }\ 
L(\www{\uuuu{e}}^t,\www{\uuuu{e}})^t=S(\www{\uuuu{x}}),
\end{eqnarray*}
is $\Br_3\ltimes\{\pm 1\}^3$ equivariant.
By Theorem \ref{t4.13} (a) 
$$\Z^3=\dot\bigcup_{\uuuu{x}\in \bigcup_{i=1}^{24}C_i}
(\Br_3\ltimes\{\pm 1\}^3)(\uuuu{x}).$$ 
Therefore we can and will restrict to 
$\uuuu{x}\in \bigcup_{i=1}^{24}C_i$.

The stabilizer $(\Br_3)_{\uuuu{e}/\{\pm 1\}^3}$ is
by Lemma \ref{t3.25} (e) the kernel of the group 
antihomomorphism
\begin{eqnarray*}
\oooo{Z}:(\Br_3)_{\uuuu{x}/\{\pm 1\}^3}\to 
G_\Z^{\BB}/Z((\{\pm 1\}^3)_{\uuuu{x}}).
\end{eqnarray*}
Here this simplifies to 
\begin{eqnarray*}
\oooo{Z}:(\Br_3)_{\uuuu{x}/\{\pm 1\}^3}\to 
G_\Z/Z((\{\pm 1\}^3)_{\uuuu{x}})
\end{eqnarray*}
because $G_\Z^{\BB}=G_\Z$ in the reducible cases
(and also in most irreducible cases) and
$Z((\{\pm 1\}^3)_{\uuuu{x}})=\{\pm \id\}$, which is
a normal subgroup of $G_\Z$, in the irreducible
cases by Lemma \ref{t3.25} (f). 

Theorem \ref{t4.16} gives the stabilizer
$(\Br_3)_{\uuuu{x}/\{\pm 1\}^3}$ in all cases.
The following Theorem \ref{t7.11} gives the stabilizer
$(\Br_3)_{\uuuu{e}/\{\pm 1\}^3}$ in all cases.

\begin{theorem}\label{t7.11}
Consider a local minimum $\uuuu{x}\in C_i
\subset\Z^3$ for some $i\in\{1,...,24\}$ 
and the pseudo-graph $\GG_j$ with 
$\GG_j=\GG(\uuuu{x})$. In the following table,
the entry in the fourth column and in the line of $C_i$
is the stabilizer $(\Br_3)_{\uuuu{e}/\{\pm 1\}^3}$. 
The first, second and third column are copied from the
table in Theorem \ref{t4.13}. 
\begin{eqnarray*}
\begin{array}{l|l|l|l}
 & \textup{sets} & 
(\Br_3)_{\uuuu{x}/\{\pm 1\}^3} 
& (\Br_3)_{\uuuu{e}/\{\pm 1\}^3}\\ \hline 
\GG_1 & C_1\ (A_1^3)
& \Br_3 & \Br_3^{pure}\\
\GG_1 & C_2\ (\HH_{1,2})  
& \Br_3 & \langle (\sigma^{mon})^2 \rangle \\
\GG_2 & C_3\ (A_2A_1)
& \langle \sigma_1,\sigma_2^2\rangle & 
\langle \sigma_2^2,(\sigma^{mon})^{-1}\sigma_1^2,
\sigma^{mon}\sigma_1\rangle \\
 & & & = \langle \sigma_2^2, \sigma_1\sigma_2^2\sigma_1^{-1},
 \sigma_1^3\rangle \\
\GG_2 & C_4\ (\P^1A_1),C_5
& \langle \sigma_1,\sigma_2^2\rangle & 
\langle \sigma_2^2,(\sigma^{mon})^{-1}\sigma_1^2\rangle \\
 & & & =\langle \sigma_2^2,\sigma_1\sigma_2^2\sigma_1^{-1}
\rangle \\ 
\GG_3 & C_6\ (A_3) 
& \langle\sigma_1\sigma_2,\sigma_1^3\rangle & 
\langle (\sigma_1\sigma_2)^4,\sigma_1^3\rangle \\
\GG_4 & C_7 \ (\widehat{A}_2) 
& \langle\sigma_2\sigma_1,\sigma_1^3\rangle & 
\langle \sigma_1^3,\sigma_2^3,\sigma_2\sigma_1^3\sigma_2^{-1}
\rangle \\
\GG_5 & C_8,\ C_9\ ((-l,2,-l))
&\langle\sigma^{mon},\sigma_1^{-1}\sigma_2^{-1}\sigma_1\rangle
 & \langle 
(\sigma^{mon})^2\sigma_1^{-1}\sigma_2^{l^2-4}\sigma_1\rangle \\
\GG_6 & C_{10}\ (\P^2),\ C_{11},\ C_{12} 
& \langle\sigma_2\sigma_1\rangle & 
\langle \id\rangle \\
\GG_7 & C_{13}\ (\textup{e.g. }(4,4,8)) 
& \langle \sigma_2\sigma_1^2\rangle & 
\langle\id \rangle \\
\GG_8 & C_{14}\ (\textup{e.g. }(3,4,6)) 
& \langle \sigma^{mon}\rangle & 
\langle \id\rangle \\
\GG_9 & C_{15},\ C_{16},\ C_{23},\ C_{24} 
&\langle \sigma^{mon}\rangle & 
\langle \id\rangle \\
\GG_{10} & C_{17}\ (\textup{e.g. }(-2,-2,0)) 
& \langle \sigma^{mon},\sigma_2\rangle & 
\langle \sigma_2^2\rangle \\
\GG_{11} & C_{18}\ (\textup{e.g. }(-3,-2,0)) 
& \langle \sigma^{mon},\sigma_2^2\rangle & 
\langle \sigma_2^2\rangle \\
\GG_{12} & C_{19}\ (\textup{e.g. }(-2,-1,0)) 
& \langle \sigma^{mon},\sigma_2^2,
\sigma_2\sigma_1^3\sigma_2^{-1}\rangle & 
\langle \sigma_2^2,\sigma_2\sigma_1^3\sigma_2^{-1}\rangle \\
\GG_{13} & C_{20}\ (\textup{e.g. }(-2,-1,-1)) 
& \langle \sigma^{mon},\sigma_2^3,
\sigma_2\sigma_1^3\sigma_2^{-1}\rangle & 
\langle \sigma_2^3,\sigma_2\sigma_1^3\sigma_2^{-1}\rangle \\
\GG_{14} & C_{21},\ C_{22}
& \langle \sigma^{mon},\sigma_2^3\rangle & 
\langle \sigma_2^3\rangle 
\end{array}
\end{eqnarray*}
\end{theorem}

{\bf Proof:}
{\bf The reducible case $\GG_1\, \&\, C_1\, (A_1^3)$:} 
Here $\uuuu{x}=(0,0,0)$ and 
$(\Br_3)_{\uuuu{x}/\{\pm 1\}^3}=\Br_3$. 
$$\BB^{dist}=\{(\varepsilon_1 e_{\sigma(1)},
\varepsilon_2 e_{\sigma(2)},\varepsilon_3 e_{\sigma(3)})\,|\, 
\varepsilon_1,\varepsilon_2,\varepsilon_3\in\{\pm 1\},
\sigma\in S_3\},$$
and $\Br_3$ acts by permutation of the entries of triples
on $\BB^{dist}$. 
Therefore $(\Br_3)_{\uuuu{e}/\{\pm 1\}^3}$ is the kernel of
the natural group homomorphism $\Br_3\to S_3$, 
so it is the subgroup $\Br_3^{pure}$ of pure braids
(see Remark \ref{t8.2} (vii) for this group). 

{\bf The case $\GG_1\, \&\, C_2\, (\HH_{1,2})$:}
Here $\uuuu{x}=(2,2,2)\sim(-2,2,-2)$ and 
$(\Br_3)_{\uuuu{x}/\{\pm 1\}^3}=\Br_3$. 
Recall the case $\HH_{1,2}$
in the proof of Theorem \ref{t5.14},
recall the $\Z$-basis $\www{f}$ of 
$H_\Z=H_{\Z,1}\oplus H_{\Z,2}$, and recall
$$G_\Z=\{g\in\Aut(H_\Z,1)\,|\, \det g=1\}\times
\Aut(H_{\Z,2})\cong SL_2(\Z)\times \{\pm 1\}.$$
We found in the proof of Theorem \ref{t5.14} (c)
\begin{eqnarray*}
Z(\delta_2\sigma_1)=(\www{f}\mapsto \www{f}
\begin{pmatrix}1&-1&0\\0&1&0\\0&0&1\end{pmatrix}),\quad
Z(\delta_3\sigma_2)=(\www{f}\mapsto \www{f}
\begin{pmatrix}1&0&0\\1&1&0\\0&0&1\end{pmatrix}).
\end{eqnarray*}
The group antihomomorphism
$\oooo{Z}:(\Br_3)_{\uuuu{x}/\{\pm 1\}^3}=\Br_3\to 
G_\Z/\{\pm \id\}\cong SL_2(\Z)$
is surjective with $\oooo{Z}(\sigma_1)\equiv A_1$ and
$\oooo{Z}(\sigma_2)\equiv A_2$.
It almost coincides with the group homomorphism
$\Br_3\to SL_2(\Z)$ in Remark \ref{t4.15} (i).
It has the same kernel $\langle (\sigma^{mon})^2\rangle$.

{\bf The reducible cases 
$\GG_2\, \&\, C_3\, (A_2A_1), C_4\, (\P^1A_1),C_5$:}
Here $\uuuu{x}=(x_1,0,0)$ with $x_1\leq -1$ and 
$(\Br_3)_{\uuuu{x}/\{\pm 1\}^3}=\langle\sigma_1,\sigma_2^2
\rangle$. The quotient group 
$(\Br_3)_{\uuuu{x}/\{\pm 1\}^3}/(\Br_3)_{\uuuu{e}/\{\pm 1\}^3}$
is by Theorem \ref{t3.28} (c)  and Lemma \ref{t3.25} (e) 
isomorphic to the quotient group 
$G_\Z/Z((\{\pm 1\}^3)_{\uuuu{x}})$. 
Here $Z((\pm 1\}^3)_{\uuuu{x}})
=\langle (-1,-1,-1),(-1,-1,1)\rangle$ with 
\begin{eqnarray*}
Z((-1,-1,-1))&=&-\id,\\
Z((-1,-1,1))&=&(\uuuu{e}\mapsto (-e_1,-e_2,e_3))=Q.
\end{eqnarray*}
Define 
$$M^{root}:=Z(\delta_2\sigma_1)=
(\uuuu{e}\mapsto \uuuu{e}\begin{pmatrix}-x_1&-1&0\\
1&0&0\\0&0&1\end{pmatrix}),$$
and recall from Theorem \ref{t5.13} and Theorem \ref{t5.5} 
\begin{eqnarray*}
G_\Z=\{\pm (M^{root})^l\,|\, l\in\Z\}\times
\{\id,Q\}.
\end{eqnarray*}
Therefore
$$(\Br_3)_{\uuuu{x}/\{\pm 1\}^3}
/(\Br_3)_{\uuuu{e}/\{\pm 1\}^3}
\cong G_\Z/Z((\{\pm 1\}^3)_{\uuuu{x}})\cong
\{(M^{root})^l\,|\, l\in\Z\}.$$
In the case $C_3\, (A_2A_1)$ $x_1=-1$ and $M^{root}$
has order three, so the quotient group 
$(\Br_3)_{\uuuu{x}/\{\pm 1\}^3}/(\Br_3)_{\uuuu{e}/\{\pm 1\}^3}$
is cyclic of order three with generator the class $[\sigma_1]$
of $\sigma_1$.
In the cases $C_4$ and $C_5$ $x_1\leq -2$ and $M^{root}$
has infinite order, so the quotient group
$(\Br_3)_{\uuuu{x}/\{\pm 1\}^3}/(\Br_3)_{\uuuu{e}/\{\pm 1\}^3}$
is cyclic of infinite order with generator the class 
$[\sigma_1]$ of $\sigma_1$.

The cases $C_4\, (\P^1A_1), C_5$:
Theorem \ref{t7.1} (b) can be applied to the subbasis
$(e_2,e_3)$ with $x_3=0$. It shows
$\sigma_2^2\in (\Br_3)_{\uuuu{e}/\{\pm 1\}^3}$.
Therefore $(\Br_3)_{\uuuu{e}/\{\pm 1\}^3}$ contains the
normal closure of $\sigma_2^2$ in 
$(\Br_3)_{\uuuu{x}/\{\pm 1\}^3}
=\langle\sigma_1,\sigma_2^2\rangle$. 
This normal subgroup is obviously
$$\langle \sigma_1^l\sigma_2^2\sigma_1^{-l}\,|\,l\in\Z\rangle.$$
It can also be written with two generators, namely it is
$$\langle \sigma_2^2,\sigma_1\sigma_2^2\sigma_1^{-1}\rangle 
= \langle \sigma_2^2,(\sigma^{mon})^{-1}\sigma_1^2\rangle.$$
The equality of left and right side follows from
$$\sigma_2^2\cdot\sigma_1\sigma_2^2\sigma_1^{-1}
=\sigma_2^2\sigma_1\sigma_2^2\sigma_1\cdot \sigma_1^{-2}
=\sigma^{mon}\sigma_1^{-2}.$$
The equality of this group with 
$\langle \sigma_1^l\sigma_2^2\sigma_1^{-l}\,|\,l\in\Z\rangle$
follows from the fact that $\sigma^{mon}$ is in the center
of $\Br_3$. 
The quotient group $\langle\sigma_1,\sigma_2^2\rangle/
\langle \sigma_1^l\sigma_2^2\sigma_1^{-l}\,|\, l\in\Z\rangle$
is cyclic of infinite order with generator the class
$[\sigma_1]$ of $\sigma_1$. Therefore
\begin{eqnarray*}
(\Br_3)_{\uuuu{e}/\{\pm 1\}^3}
&=& \langle \sigma_2^2,(\sigma^{mon})^{-1}\sigma_1^2\rangle
=\langle \sigma_2^2,\sigma_1\sigma_2^2\sigma_1^{-1}\rangle \\
&=& 
\langle \sigma_1^l\sigma_2^2\sigma_1^{-l}\,|\, l\in\Z\rangle\\
&=& (\textup{the normal closure of }
\sigma_2^2\textup{ in }\langle\sigma_1,\sigma_2^2\rangle).
\end{eqnarray*}

The case $C_3\, (A_2A_1)$: 
Theorem \ref{t7.1} (b) can be applied to the subbasis
$(e_1,e_2)$ with $x_1=-1$ and to the subbasis $(e_2,e_3)$
with $x_3=0$. It shows $\sigma_1^3$ and $\sigma_2^2
\in (\Br_3)_{\uuuu{e}/\{\pm 1\}^3}$. Therefore
$(\Br_3)_{\uuuu{e}/\{\pm 1\}^3}$ contains the normal 
closure of $\sigma_1^3$ and $\sigma_2^2$ in
$(Br_3)_{\uuuu{x}/\{\pm 1\}^3}=\langle \sigma_1,\sigma_2^2
\rangle$. 
The quotient group 
$$\langle\sigma_1,\sigma_2^2\rangle /
(\textup{the normal closure of }\sigma_1^3\textup{ and }
\sigma_2^2\textup{ in }\langle \sigma_1,\sigma_2^2\rangle)$$
is cyclic of order three with generator the class
$[\sigma_1]$ of $\sigma_1$. Therefore
\begin{eqnarray*}
(\Br_3)_{\uuuu{e}/\{\pm 1\}^3}=
(\textup{the normal closure of }\sigma_1^3
\textup{ and }\sigma_2^2\textup{ in }
\langle\sigma_1,\sigma_2^2\rangle).
\end{eqnarray*}
It coincides with the subgroup generated by $\sigma_1^3$ and
by the normal closure $\langle \sigma_2^2,(\sigma^{mon})^{-1}
\sigma_1^2\rangle$ of $\sigma_2^2$ in $\langle\sigma_1,
\sigma_2^2\rangle$. Therefore 
\begin{eqnarray*}
(\Br_3)_{\uuuu{e}/\{\pm 1\}^3}&=& 
\langle\sigma_2^2,(\sigma^{mon})^{-1}\sigma_1^2,\sigma_1^3
\rangle
=\langle \sigma_2^2,(\sigma^{mon})^{-1}\sigma_1^2,
\sigma^{mon}\sigma_1\rangle.
\end{eqnarray*}

{\bf The case $\GG_3\, \&\, C_6\, (A_3)$:}
Here $\uuuu{x}=(-1,0,-1)$ and 
$(\Br_3)_{\uuuu{x}/\{\pm 1\}^3}=\langle \sigma_1\sigma_2,
\sigma_1^3\rangle$. 
By Theorem \ref{t5.14} (b) 
$G_\Z=\{\pm M^l\,|\, l\in\{0,1,2,3\}\}$,
and $M$ has order four. By Theorem \ref{t3.28} and
Lemma \ref{t3.25} (f), the antihomomorphism
$$\oooo{Z}:(\Br_3)_{\uuuu{x}/\{\pm 1\}^3} / 
(\Br_3)_{\uuuu{e}/\{\pm 1\}^3}\to G_\Z/\{\pm \id\}$$
is an antiisomorphism. Therefore the quotient group
$(\Br_3)_{\uuuu{x}/\{\pm 1\}^3} / 
(\Br_3)_{\uuuu{e}/\{\pm 1\}^3}$
is cyclic of order four. 

Theorem \ref{t7.1} (b) can be applied to the subbasis
$(e_1,e_2)$ with $x_1=-1$. It shows 
$\sigma_1^3\in (\Br_3)_{\uuuu{e}/\{\pm 1\}^3}$.
Observe also
\begin{eqnarray*}
\delta_1\sigma_1\sigma_2(\uuuu{e})
&=& \delta_1\sigma_1(e_1,e_3+e_2,e_2)
=\delta_1(e_1+e_2+e_3,e_1,e_2)\\
&=& (-e_1-e_2-e_3,e_1,e_2)=-M^{-1}(\uuuu{e}),\\
\textup{so }
Z(\delta_1\sigma_1\sigma_2)&=& -M^{-1}.
\end{eqnarray*}
$M$ and $-M^{-1}$ have order four. Therefore
$(\sigma_1\sigma_2)^4\in (\Br_3)_{\uuuu{e}/\{\pm 1\}^3}$.
Thus 
$(\Br_3)_{\uuuu{e}/\{\pm 1\}^3}
\supset\langle (\sigma_1\sigma_2)^4,\sigma_1^3\rangle.$

We will show first that 
$\langle (\sigma_1\sigma_2)^4,\sigma_1^3\rangle$
is a normal subgroup of 
$\langle \sigma_1\sigma_2,\sigma_1^3\rangle$ and then that
the quotient group is cyclic of order four. This will imply
$(\Br_3)_{\uuuu{e}/\{\pm 1\}^3}
=\langle (\sigma_1\sigma_2)^4,\sigma_1^3\rangle.$

Recall that $\sigma^{mon}=(\sigma_2\sigma_1)^3
=(\sigma_1\sigma_2)^3$ generates the center of $\Br_3$.
Therefore 
\begin{eqnarray*}
(\sigma_1\sigma_2)^l\sigma_1^3(\sigma_1\sigma_2)^{-l}
=(\sigma_1\sigma_2)^{4l}\sigma_1^3(\sigma_1\sigma_2)^{-4l}
\in\langle (\sigma_1\sigma_2)^4,\sigma_1^3\rangle
\textup{ for any }l\in\Z.
\end{eqnarray*}
Thus $\langle (\sigma_1\sigma_2)^4,\sigma_1^3\rangle$
is a normal subgroup of 
$\langle \sigma_1\sigma_2,\sigma_1^3\rangle$.
This also shows that the quotient group
$\langle\sigma_1\sigma_1,\sigma_1^3\rangle /
\langle (\sigma_1\sigma_2)^4,\sigma_1^3\rangle$ is cyclic
of order four. Therefore
$\langle (\sigma_1\sigma_2)^4,\sigma_1^3\rangle
=(\Br_3)_{\uuuu{e}/\{\pm 1\}^3}$.

{\bf The case $\GG_4\, \&\,  C_7\, (\whh{A}_2)$:}
Here $\uuuu{x}=(-1,-1,-1)$ and 
$(\Br_3)_{\uuuu{x}/\{\pm 1\}^3}=\langle \sigma_2\sigma_1,
\sigma_1^3\rangle$. 
By Theorem \ref{t5.14} (b) 
$G_\Z=\{\pm (M^{root})^l\,|\, l\in\Z\}$, and $M^{root}$
has infinite order. 
By Theorem \ref{t3.28} and Lemma \ref{t3.25} (f),
the antihomomorphism
$$\oooo{Z}:(\Br_3)_{\uuuu{x}/\{\pm 1\}^3} / 
(\Br_3)_{\uuuu{e}/\{\pm 1\}^3}\to G_\Z/\{\pm \id\}$$
is an antiisomorphism. By Theorem \ref{t3.26} (c) and
Theorem \ref{t5.14} (b) $Z(\delta_3\sigma_2\sigma_1)=M^{root}$.
Therefore the quotient group
$(\Br_3)_{\uuuu{x}/\{\pm 1\}^3} / 
(\Br_3)_{\uuuu{e}/\{\pm 1\}^3}$
is cyclic of infinite order with generator the class
$[\sigma_2\sigma_1]$ of $\sigma_2\sigma_1$. 

Theorem \ref{t7.1} (b) can be applied to the subbasis
$(e_1,e_2)$ with $x_1=-1$. It shows 
$\sigma_1^3\in (\Br_3)_{\uuuu{e}/\{\pm 1\}^3}$.
Therefore $(\Br_3)_{\uuuu{e}/\{\pm 1\}^3}$
contains the normal closure of $\sigma_1^3$ in
$\langle \sigma_2\sigma_1,\sigma_1^3\rangle$.
We will first determine this normal closure and then show 
that it equals $(\Br_3)_{\uuuu{e}/\{\pm 1\}^3}$.

As $\sigma^{mon}=(\sigma_2\sigma_1)^3$ generates the center
of $\Br_3$, 
\begin{eqnarray*}
(\sigma_2\sigma_1)^{\varepsilon+3l}\sigma_1^3
(\sigma_2\sigma_1)^{-\varepsilon-3l}
=(\sigma_2\sigma_1)^{\varepsilon}\sigma_1^3
(\sigma_2\sigma_1)^{-\varepsilon}
\quad\textup{for }\varepsilon\in\{0;\pm 1\},l\in\Z.
\end{eqnarray*}
One sees
\begin{eqnarray*}
(\sigma_2\sigma_1)\sigma_1(\sigma_2\sigma_1)^{-1}
&=&\sigma_2\sigma_1\sigma_2^{-1},\quad\textup{so}\\
(\sigma_2\sigma_1)\sigma_1^3(\sigma_2\sigma_1)^{-1}
&=&\sigma_2\sigma_1^3\sigma_2^{-1},\\
(\sigma_2\sigma_1)^{-1}\sigma_1(\sigma_2\sigma_1)
&=&\sigma_1^{-1}(\sigma_2^{-1}\sigma_1\sigma_2)\sigma_1
\stackrel{\eqref{4.14}}{=}
\sigma_1^{-1}(\sigma_1\sigma_2\sigma_1^{-1})\sigma_1
=\sigma_2, \quad\textup{so}\\
(\sigma_2\sigma_1)^{-1}\sigma_1^3(\sigma_2\sigma_1)
&=&\sigma_2^3.
\end{eqnarray*}
Therefore the normal closure of $\sigma_1^3$ in
$\langle \sigma_2\sigma_1,\sigma_1^3\rangle$ is
$\langle \sigma_1^3,\sigma_2^3,\sigma_2\sigma_1^3\sigma_2^{-1}
\rangle$. The quotient group is an infinite cyclic
group with generator the class $[\sigma_2\sigma_1]$
of $\sigma_2\sigma_1$. Therefore
\begin{eqnarray*}
(\Br_3)_{\uuuu{e}/\{\pm 1\}^3}&=&
\langle \sigma_1^3,\sigma_2^3,\sigma_2\sigma_1^3\sigma_2^{-1}
\rangle\\
&=& (\textup{the normal closure of }\sigma_1^3\textup{ in }
\langle \sigma_2\sigma_1,\sigma_1^3\rangle).
\end{eqnarray*}

{\bf The cases $\GG_5\, \&\, C_8$:}
Here $\uuuu{x}=(-l,2,-l)$ with $l\geq 3$ odd and 
$(\Br_3)_{\uuuu{x}/\{\pm 1\}^3}=\langle \sigma^{mon},
\sigma_1^{-1}\sigma_2^{-1}\sigma_1\rangle$. 

By Theorem \ref{t5.14} (b) 
$G_\Z=\{\pm (M^{root})^l\,|\,l\in\Z\}$
with $M^{root}$ as in Theorem \ref{t5.14} (b) with
$(M^{root})^{l^2-4}=-M$. Because $l$ is odd, the cyclic group
$G_\Z/\{\pm \id\}$ with generator $[M^{root}]$ can also be 
written as
$$G_\Z/\{\pm \id\}
=\langle [M],[(M^{root})^2]\rangle
\quad\textup{with }[M]=[M^{root}]^{l^2-4}.$$

By Theorem \ref{t3.28} and Lemma \ref{t3.25} (f), the
antihomomorphism 
$$\oooo{Z}:(\Br_3)_{\uuuu{x}/\{\pm 1\}^3}\to G_\Z/\{\pm \id\}$$
is surjective with kernel $(\Br_3)_{\uuuu{e}/\{\pm 1\}^3}$.
By the proof of Theorem \ref{t5.14} (b)(iv)
$[M]=\oooo{Z}(\sigma^{mon})$ and
$[(M^{root})^2]=\oooo{Z}(\sigma_1^{-1}\sigma_2^{-1}\sigma_1)$.
The single relation between $[M]$ and $[(M^{root})^2]$ is
$[\id]=[M]^2([(M^{root})^2])^{4-l^2}$. Therefore the 
kernel $(Br_3)_{\uuuu{e}/\{\pm 1\}^3}$ of the group
antihomomorphism $\oooo{Z}$ above is
generated by 
$(\sigma^{mon})^2(\sigma_1^{-1}\sigma_2^{-1}\sigma_1)^{4-l^2}
=(\sigma^{mon})^2\sigma_1^{-1}\sigma_2^{l^2-4}\sigma_1$.

{\bf The cases $\GG_6\, \&\, C_9$:}
Here $\uuuu{x}=(-l,2,-l)$ with $l\geq 4$ even and 
$(\Br_3)_{\uuuu{x}/\{\pm 1\}^3}=\langle \sigma^{mon},
\sigma_1^{-1}\sigma_2^{-1}\sigma_1\rangle$. 

By Theorem \ref{t3.28} and Lemma \ref{t3.25} (f), the
antihomomorphism 
$$\oooo{Z}:(\Br_3)_{\uuuu{x}/\{\pm 1\}^3}\to G_\Z/\{\pm \id\}$$
is surjective with kernel $(\Br_3)_{\uuuu{e}/\{\pm 1\}^3}$.

By Theorem \ref{t5.14} (a) and the proof of Theorem
\ref{t5.14} (c)
\begin{eqnarray*}
G_\Z&=&\langle -\id,\www{M},Q\rangle\\
\textup{with}\quad 
\www{M}&=&Z(\delta_3\sigma_1^{-1}\sigma_2^{-1}\sigma_1),\
Q=Z(\sigma^{mon})
Z(\delta_3\sigma_1^{-1}\sigma_2^{-1}\sigma_1)^{2-l^2/2},
\end{eqnarray*}
$\www{M}$ and $Q$ commute, $\www{M}$ has infinite order, 
$Q$ has order two.
Therefore the kernel $(\Br_3)_{\uuuu{e}/\{\pm 1\}^3}$
of the group antihomomorphism $\oooo{Z}$ above is generated by
$(\sigma^{mon}(\sigma_1^{-1}\sigma_2^{-1}\sigma_1)^{2-l^2/2})^2
=(\sigma^{mon})^2\sigma_1^{-1}\sigma_2^{l^2-4}\sigma_1$.

{\bf The cases $\GG_6\, \&\, C_{10}\, (\P^2),C_{11},C_{12}$:}
Here $(\Br_3)_{\uuuu{x}/\{\pm 1\}^3}
=\langle \sigma_2\sigma_1\rangle$.
Here $Z(\sigma_2\sigma_1)=M^{root}$ is a third root of 
the monondromy. The monodromy $M$ and $M^{root}$ have infinite
order. Therefore the kernel $(\Br_3)_{\uuuu{e}/\{\pm 1\}^3}$
of the group antihomomorphism
$$\oooo{Z}:(\Br_3)_{\uuuu{x}/\{\pm 1\}^3}
=\langle \sigma_2\sigma_1\rangle\to G_\Z/\{\pm \id\}$$
is $\langle \id\rangle$.

{\bf The cases $\GG_7\, \&\, C_{13}\, (\textup{e.g. }(4,4,8))$:}
Here $(\Br_3)_{\uuuu{x}/\{\pm 1\}^3}
=\langle \sigma_2\sigma_1^2\rangle$.
Here $Z(\sigma_2\sigma_1)=M^{root}$ is a root of 
the monondromy. The monodromy $M$ and $M^{root}$ have infinite
order. Therefore the kernel $(\Br_3)_{\uuuu{e}/\{\pm 1\}^3}$
of the group antihomomorphism
$$\oooo{Z}:(\Br_3)_{\uuuu{x}/\{\pm 1\}^3}
=\langle \sigma_2\sigma_1^2\rangle\to G_\Z/\{\pm \id\}$$
is $\langle \id\rangle$.

{\bf The cases $\GG_8\, \&\, C_{14}, \GG_9\, \&\,  C_{15},C_{16},C_{23},C_{24}$:}
Here $(\Br_3)_{\uuuu{x}/\{\pm 1\}^3}
=\langle \sigma^{mon}\rangle$.
The monodromy $M=Z(\sigma^{mon})$ has infinite order. 
Therefore the kernel $(\Br_3)_{\uuuu{e}/\{\pm 1\}^3}$
of the group antihomomorphism
$$\oooo{Z}:(\Br_3)_{\uuuu{x}/\{\pm 1\}^3}
=\langle \sigma^{mon}\rangle\to G_\Z/\{\pm \id\}$$
is $\langle \id\rangle$.

{\bf The cases $\GG_{10}\, \&\, C_{17}\, (\textup{e.g. }
(-2,-2,0))$:}
Here $(\Br_3)_{\uuuu{x}/\{\pm 1\}^3}
=\langle \sigma^{mon},\sigma_2\rangle$.
Recall from Theorem \ref{t5.16} (c) that
\begin{eqnarray*}
G_\Z=\{\id,Q\}\times \{\pm M^l\,|\, l\in\Z\}
=\langle -\id,Q,M\rangle,
\end{eqnarray*}
$Q$ and $M$ commute, $Q$ has order two, $M$ has infinite order,
$-Q=Z(\sigma_2)$ (see the proof of Theorem \ref{5.17} (e)), 
$M=Z(\sigma^{mon})$. 
Therefore the kernel $(\Br_3)_{\uuuu{e}/\{\pm 1\}^3}$
of the group antihomomorphism
$$\oooo{Z}:(\Br_3)_{\uuuu{x}/\{\pm 1\}^3}
=\langle \sigma^{mon},\sigma_2\rangle\to G_\Z/\{\pm \id\}$$
is $\langle \sigma_2^2\rangle$.

{\bf The cases $\GG_{11}\, \&\, C_{18}, \GG_{12}\, \&\, C_{19},
\GG_{13}\, \&\, C_{20},\GG_{14}\, \&\, C_{21},C_{22}$:}
By Theorem \ref{t4.16} in all cases 
$(\Br_3)_{\uuuu{x}/\{\pm 1\}^3}$ is generated by $\sigma^{mon}$
and some other generators. We claim that the other generators
are all in $(\Br_3)_{\uuuu{e}/\{\pm 1\}^3}$. 
Application of Theorem \ref{t7.1} (b) to the subbasis
$(e_2,e_3)$ shows this for the following generators:\\
$\sigma_2^2$ in the cases $\GG_{11}\& C_{18}$ and 
$G_{12}\& C_{19}$ because there $x_3=0$;\\
$\sigma_2^3$ in the cases $\GG_{13}\& C_{20}$ and 
$G_{14}\& C_{21},C_{22}$ because there $x_3=-1$.

$\sigma_2^{-1}$ maps $\uuuu{e}$ to the basis
$(e_1,e_3,s^{(0)}_{e_3}(e_2))$. Therefore application of
Theorem \ref{t7.1} (b) to the subbasis $(e_1,e_3)$ 
with $x_2=-1$ in the cases $\GG_{12}\& C_{19}$
and $\GG_{13}\& C_{20}$ shows 
$\sigma_2\sigma_1^3\sigma_2^{-1}
\in (\Br_3)_{\uuuu{e}/\{\pm 1\}^3}$. The claim 
$\langle \textup{other generators}\rangle
\subset (\Br_3)_{\uuuu{e}/\{\pm 1\}^3}$ is proved. 

The monodromy $M=Z(\sigma^{mon})$ has infinite order. 
Therefore the kernel $(\Br_3)_{\uuuu{e}/\{\pm 1\}^3}$
of the group antihomomorphism
$$\oooo{Z}:(\Br_3)_{\uuuu{x}/\{\pm 1\}^3}
=\langle \sigma^{mon},\textup{other generators}
\rangle\to G_\Z/\{\pm \id\}$$
is $\langle \textup{other generators}\rangle$. \hfill$\Box$

\begin{remarks}\label{t7.12}
(i) In the cases $\GG_6\,\&\, C_{10},C_{11},C_{12}$, 
$\GG_7\,\&\, C_{13}$, $\GG_8\,\&\, C_{14}$ and 
$\GG_9\,\&\, C_{15},C_{16},C_{23},C_{24}$, the even monodromy
group $\Gamma^{(0)}$ is a free Coxeter group with three
generators by Theorem \ref{t6.11} (b) and (g).
Example \ref{t3.23} (iv), which builds on Theorem \ref{t3.2},
shows $(\Br_3)_{\uuuu{e}/\{\pm 1\}^3}=\langle\id\rangle$.

Using this fact, one derives also 
$(\Br_3)_{\uuuu{x}/\{\pm 1\}^3}$ in the following way. 
In all cases $(H_\Z,L,\uuuu{e})$ is irreducible.
The group antihomomorphism
$$\oooo{Z}:(\Br_3)_{\uuuu{x}/\{\pm 1\}^3}\to G_\Z/\{\pm \id\}$$
is injective because the kernel is
$(\Br_3)_{\uuuu{e}/\{\pm 1\}^3}=\langle\id\rangle$.
By Theorem \ref{t3.28} $\oooo{Z}$ is surjective in almost
all cases. The proof of Theorem \ref{t3.28} provides
in all cases preimages of generators of
$\Imm(\oooo{Z})\subset  G_\Z/\{\pm \id\}$.
These preimages generate $(\Br_3)_{\uuuu{x}/\{\pm 1\}^3}$.

So, the arguments here on one side and the Theorems \ref{t4.13},
\ref{t4.16} and \ref{t7.11} on the other side 
offer two independent ways to derive the stabilizers 
$(\Br_3)_{\uuuu{x}/\{\pm 1\}^3}$ and 
$(\Br_3)_{\uuuu{e}/\{\pm 1\}^3}$ in the considered cases.

(ii) But the arguments in (i) cannot easily be adapted to the
other cases. In the cases $\GG_{10}\,\&\, C_{17}$, 
$\GG_{11}\,\&\,  C_{18}$, $\GG_{12}\,\&\,  C_{19}$, $\GG_{13}\,\&\,  C_{20}$
and $\GG_{14}\,\&\,  C_{21},C_{22}$, the even monodromy group
$\Gamma^{(0)}$ is a non-free Coxeter group. 

Theorem \ref{t3.2} (b) generalizes in Theorem \ref{t3.7} (b)
to a statement on the size of $\BB^{dist}$, but not to a 
statement on the stabilizer $(\Br_3)_{\uuuu{e}/\{\pm 1\}^3}$.

The cases $\GG_1$, $\GG_2$, $\GG_3$, $\GG_4$ and $\GG_5$
are with the exception of the reducible cases $C_5$
and the case $C_{10}$ the cases with $r(\uuuu{x})\in\{0,1,3,4\}$. 
For them it looks possible, but difficult, to generalize the
arguments in (i). 

The conceptual derivation of the stabilizer groups for all
cases with the Theorems \ref{t4.13}, \ref{t4.16} and
\ref{t7.11} is more elegant.
\end{remarks}

\begin{remarks}\label{t7.13}
(i) The table in Theorem \ref{t7.11} describes the stabilizer
$(\Br_3)_{\uuuu{e}/\{\pm 1\}^3}$ by generators, except for
the case $\GG_1\,\&\,  C_1(A_1^3)$ where 
$(\Br_3)_{\uuuu{e}/\{\pm 1\}^3}=\Br_3^{pure}$. In fact
\begin{eqnarray*}
\Br_3^{pure}&=& \langle \sigma_1^2,\sigma_2^2,
\sigma_2\sigma_1^2\sigma_2^{-1},
\sigma_2^{-1}\sigma_1^2\sigma_2\rangle\\
&=& (\textup{the normal closure of }\sigma_1^2\textup{ in }
\Br_3),
\end{eqnarray*}
because 
$\sigma_2\sigma_1^2\sigma_2^{-1}
=\sigma_1^{-1}\sigma_1^2\sigma_2$, 
$\sigma_2^{-1}\sigma_1^2\sigma_2
=\sigma_1\sigma_2^2\sigma_1^{-1}$ by \eqref{4.14}.

(ii) In some cases the proof of Theorem \ref{t7.11} provides
elements so that $(\Br_3)_{\uuuu{e}/\{\pm 1\}^3}$ is the
normal closure of these elements in 
$(\Br_3)_{\uuuu{x}/\{\pm 1\}^3}$:
\begin{eqnarray*}
\begin{array}{lcl}
 & \textup{elements} & (\Br_3)_{\uuuu{x}/\{\pm 1\}^3} \\
\hline 
\GG_2\,\&\,  C_3\, (A_2A_1) & \sigma_1^3,\sigma_2^2 & 
\langle \sigma_1,\sigma_2^2\rangle \\
\GG_2\,\&\,  C_4\, (\P^1A_1),C_5 & \sigma_2^2 &  
\langle \sigma_1,\sigma_2^2\rangle \\
\GG_4\,\&\,  C_7\, (\whh{A}_2) & \sigma_1^3 & 
\langle \sigma_2\sigma_1,\sigma_1^3\rangle
\end{array}
\end{eqnarray*}

(iii) In the case $\GG_3\,\&\,  C_6\, (A_3)$, the stabilizer
$(\Br_3)_{\uuuu{e}/\{\pm 1\}^3}
=\langle (\sigma_1\sigma_2)^4,\sigma_1^3\rangle$ was
determined already in \cite[Satz 7.3]{Yu90}.
\end{remarks}

The pseudo-graph $\GG(\uuuu{x})$ for 
$\uuuu{x}\in\bigcup_{i=1}^{24}C_i$ with vertex set
$\VV=\Br_3(\uuuu{x}/\{\pm 1\}^3)$ in Definition \ref{t4.9} (f),
Lemma \ref{t4.10} and the Examples \ref{t4.11} had been very
useful. All except two edges came from the generators 
$\varphi_1,\varphi_2,\varphi_3$ of the free Coxeter
group $G^{phi}$, and two edges came from $\gamma(v_0)$ and
$\gamma^{-1}(v_0)$. 
An a priori more natural choice of edges comes from the 
elementary braids $\sigma_1$ and $\sigma_2$. It is less useful,
but also interesting.

\begin{definition}\label{t7.14}
Let $\VV$ be a non-empty finite or countably infinite set
on which $\Br_3$ acts. The triple
$\GG_\sigma(\VV):=(\VV,\EE_1,\EE_2)$ with 
$\EE_1:=\{(v,\sigma_1(v))\,|\, v\in\VV\}$ and 
$\EE_2:=\{(v,\sigma_2(v))\,|\, v\in \VV\}$ is called
{\it $\sigma$-pseudo-graph of $\VV$}. 
\index{$\sigma$-pseudo-graph}
Here $\EE_1$ and $\EE_2$ are two families of directed edges.
A {\it loop} in $\EE_i$ is an edge $(v,\sigma_i(v))=(v,v)$. 
\end{definition}

\begin{remarks}\label{t7.15}
(i) In a picture of a $\sigma$-pseudo-graph, edges in $\EE_1$
and in $\EE_2$ are denoted as follows.

\includegraphics[height=0.03\textheight]{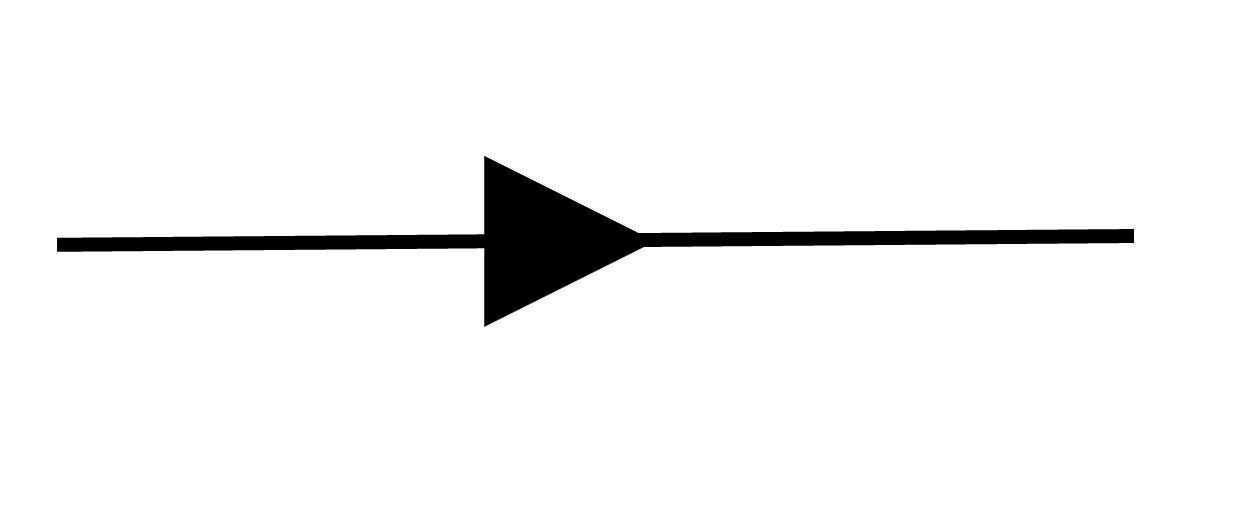} 
\quad an edge in $\EE_1$,

\includegraphics[height=0.03\textheight]{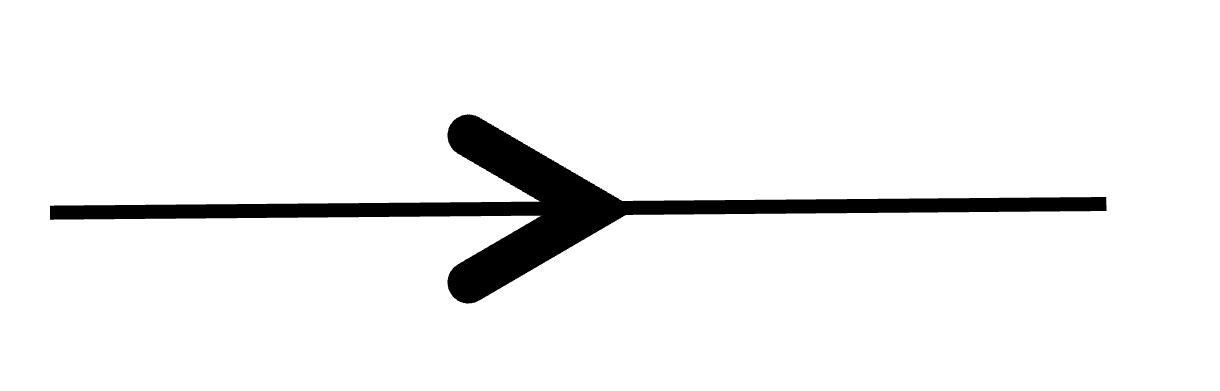} 
\quad an edge in $\EE_2$.

(ii) Consider a $\sigma$-pseudo-graph $\GG_\sigma(\VV)$.
Because $\sigma_1:\VV\to\VV$ and $\sigma_2:\VV\to\VV$
are bijections, each vertex $v\in\VV$ is starting point of 
one edge in $\EE_1$ and one edge in $\EE_2$ and end point of
one edge in $\EE_1$ and one edge in $\EE_2$. 
The $\sigma$-pseudo-graph is connected if and only if 
$\VV$ is a single $\Br_3$ orbit. 

(iii) Let $(H_\Z,L,\uuuu{e})$ be a unimodular bilinear lattice
with a triangular basis $\uuuu{e}$ with 
$L(\uuuu{e}^t,\uuuu{e})^t=S(\uuuu{x})$ for some
$\uuuu{x}\in\Z^3$. Two $\sigma$-pseudo-graphs are associated
to it, $\GG_\sigma(\BB^{dist}/\{\pm 1\}^3)$ and 
$\GG_\sigma(\Br_3(\uuuu{x}/\{\pm 1\}^3))$. The natural map
\begin{eqnarray*}
\BB^{dist}/\{\pm 1\}^3&\to& \Br_3(\uuuu{x}/\{\pm 1\}^3)\\
\www{\uuuu{e}}/\{\pm 1\}^3&\mapsto& 
\www{\uuuu{x}}/\{\pm 1\}^3\textup{ with }
L(\www{\uuuu{e}}^t,\www{\uuuu{e}})^t=S(\www{\uuuu{x}}),
\end{eqnarray*}
is $\Br_3$ equivariant and surjective. It induces a 
{\it covering} 
\begin{eqnarray*}
\GG_\sigma(\BB^{dist}/\{\pm 1\}^3)&\to 
\GG_\sigma(\Br_3(\uuuu{x}/\{\pm 1\}^3))
\end{eqnarray*}
of $\sigma$-pseudo-graphs. This is even a {\it normal covering}
with group of deck transformations 
$G_\Z^{\BB}/Z((\{\pm 1\}^3)_{\uuuu{x}})$ where
$G_\Z^{\BB}=Z((\Br_3\ltimes\{\pm 1\}^3)_{\uuuu{x}})\subset G_\Z$
is as in Lemma \ref{t3.25} (e). 
Now we explain what this means and why it holds.

The group $G_\Z^{\BB}$ acts transitively on the fiber over 
$\uuuu{x}$ of the map 
$\BB^{dist}\to (\Br_3\ltimes\{\pm 1\}^3)(\uuuu{x})$.
By Lemma \ref{t3.22} (a) the action of this group 
$G_\Z^{\BB}$ and the action of the group $\Br_3\ltimes\{\pm 1\}^3$ 
on $\BB^{dist}$ commute, so that $G_\Z^{\BB}$
acts transitively on each fiber of the map
$\BB^{dist}\to (\Br_3\ltimes\{\pm 1\}^3)(\uuuu{x})$.
Therefore the group $G_\Z^{\BB}/Z((\{\pm 1\}^3)_{\uuuu{x}})$ 
acts simply transitively on each fiber of the covering
$\GG_\sigma(\BB^{dist}/\{\pm 1\}^3)\to 
\GG_\sigma(\Br_3(\uuuu{x}/\{\pm 1\}^3))$
and is a group of automorphisms of the $\sigma$-pseudo-graph
$\GG_\sigma(\BB^{dist}/\{\pm 1\}^3)$. The quotient by this
group is the $\sigma$-pseudo-graph 
$\GG_\sigma(\Br_3(\uuuu{x}/\{\pm 1\}^3))$.
These statements are the meaning of the 
{\it normal covering}
$\GG_\sigma(\BB^{dist}/\{\pm 1\}^3)\to 
\GG_\sigma(\Br_3(\uuuu{x}/\{\pm 1\}^3))$. 

(iv) In part (iii) the $\sigma$-pseudo-graph 
$\GG_\sigma(\BB^{dist}/\{\pm 1\}^3)$ contains no loops,
and for any $v\in\BB^{dist}/\{\pm 1\}^3$ 
$\sigma_1(v) \neq \sigma_2(v)$.
The $\sigma$-pseudo-graph 
$\GG_\sigma(\Br_3(\uuuu{x}/\{\pm 1\}^3))$
contains loops in a few cases. It contains a vertex
$v$ with $\sigma_1(v)=\sigma_2(v)$ only in the cases
$\uuuu{x}=(0,0,0)$ ($A_1^3$) and 
$\uuuu{x}=(2,2,2)$ ($\HH_{1,2}$) where 
$\Br_3(\uuuu{x}/\{\pm 1\}^3)$ has only one vertex anyway. 
\end{remarks}

\begin{figure}
\parbox{2cm}{\includegraphics[height=0.15\textheight]
{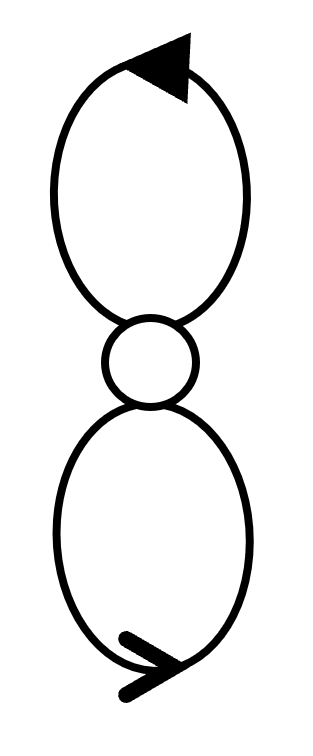}}
$\uuuu{x}\in C_1\cup C_2=\{(0,0,0),(2,2,2)\},\, A_1^3,\HH_{1,2}$\\ 
\includegraphics[height=0.15\textheight]{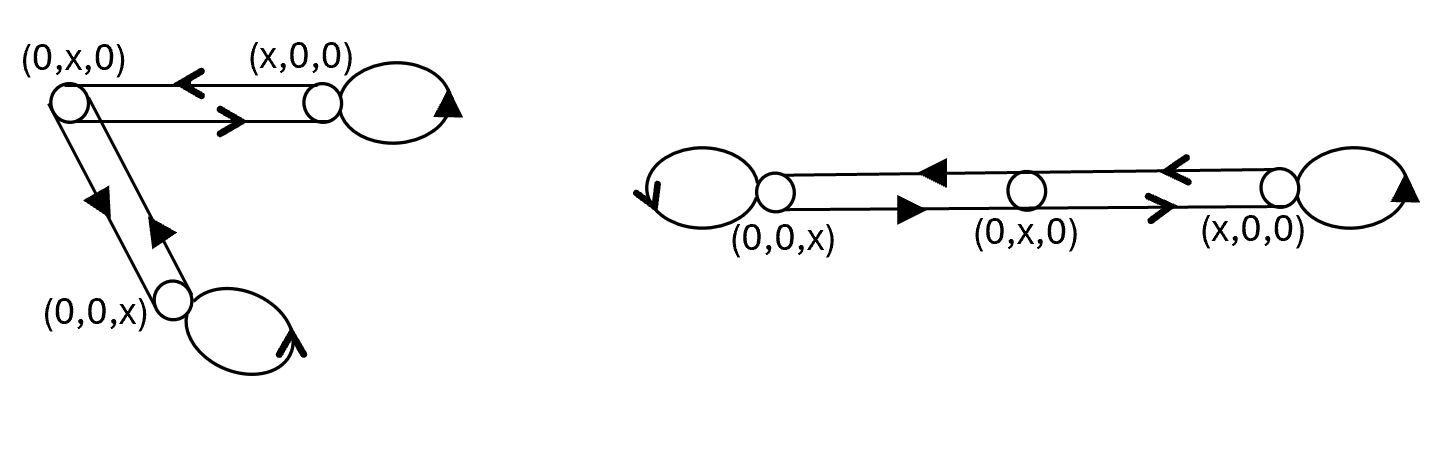}\\
Two equivalent pictures for the cases 
$\uuuu{x}=(x,0,0)\in C_3\cup C_4\cup C_5=\{(\www{x},0,0)\,|\, \www{x}<0\}$\\
($A_2A_1,\P^1A_1$, other reducible cases without $A_1^3$)\\
\parbox{6cm}{\includegraphics[height=0.2\textheight]{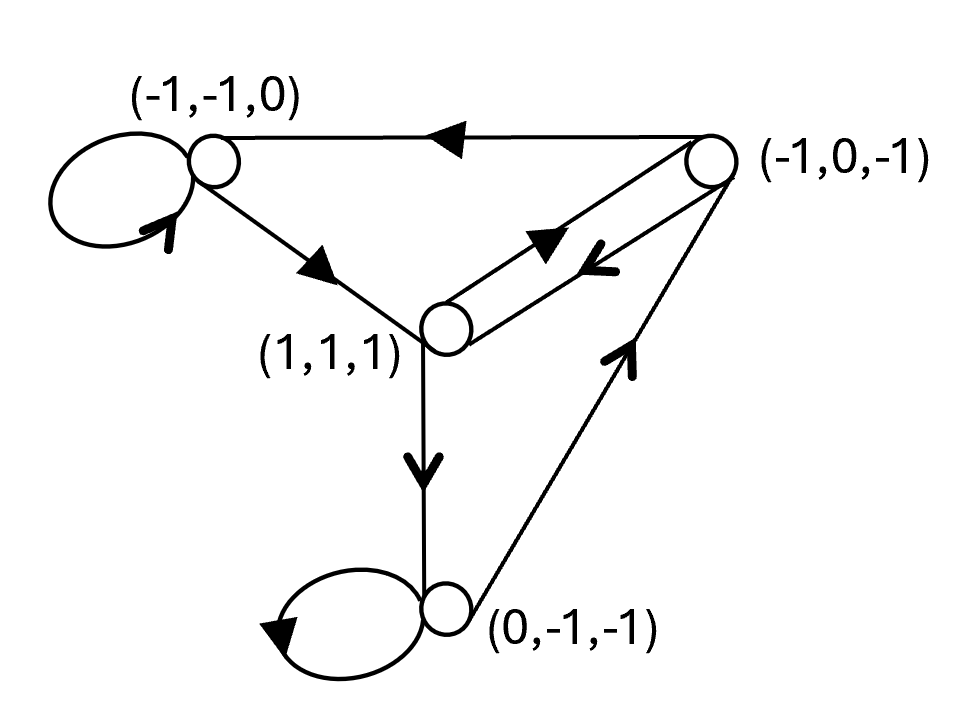}}
$\uuuu{x}=(-1,0,-1)\in C_6,\ A_3$\\
\parbox{6cm}{\includegraphics[height=0.2\textheight]{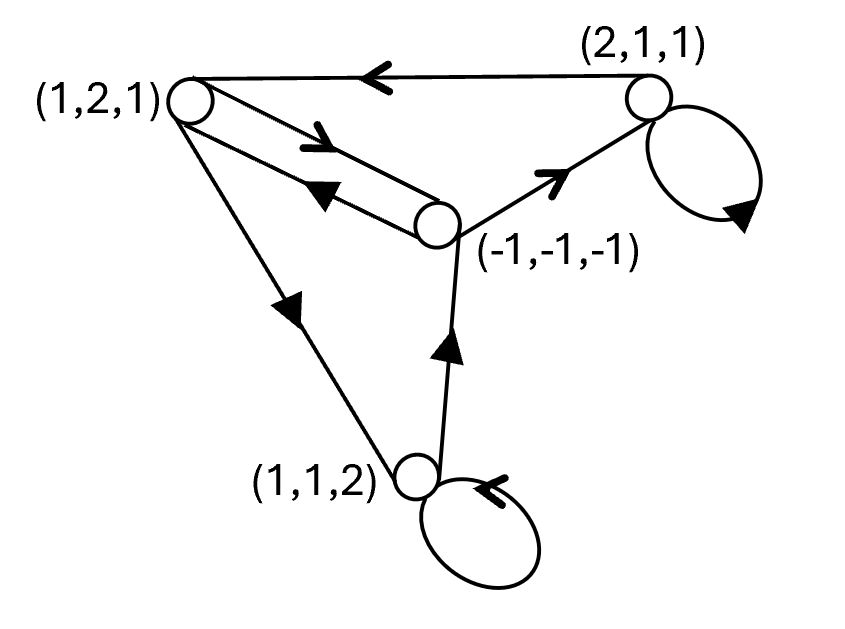}}
$\uuuu{x}=(-1,-1,-1)\in C_7,\, \whh{A}_2$
\caption[Figure 7.1]{Examples \ref{t7.16} (i): The $\sigma$-pseudo-graphs for
the finite $\Br_3$ orbits in $\Z^3/\{\pm 1\}^3$}
\label{Fig:7.1}
\end{figure}

\begin{examples}\label{t7.16}
(i) By Theorem \ref{t4.13} (a) $\Z^3/\{\pm 1\}^3$
consists of the $\Br_3$ orbits $\Br_3(\uuuu{x}/\{\pm 1\}^3)$
for $\uuuu{x}\in \bigcup_{i=1}^{24}C_i$. 
Precisely for $\uuuu{x}\in\bigcup_{i=1}^7C_i$ such an orbit 
is finite. This led to the four pseudo-graphs 
$\GG_1,\GG_2,\GG_3,\GG_4$ in the Examples \ref{t4.11}. 
The four corresponding 
$\sigma$-pseudo-graphs $\GG_\sigma(\Br_3(\uuuu{x}/\{\pm 1\}^3))$ 
are listed in Figure \ref{Fig:7.1}. 
A vertex $\www{\uuuu{x}}/\{\pm 1\}^3\in\Z^3/\{\pm 1\}^3$
is denoted by a representative 
$\www{\uuuu{x}}\in\Z^3_{>0}\cup\Z^3_{\leq 0}$. 
The vertices are positioned at the same places as in the
pictures in the Examples \ref{t4.11} for 
$\GG_1,\GG_2,\GG_3,\GG_4$.

(ii) The case $\HH_{1,2}$, $\uuuu{x}=(2,2,2)$:
Here $\Br_3(\uuuu{x}/\{\pm 1\}^3)=\{\uuuu{x}/\{\pm 1\}^3\}$
has only one vertex, but the group
$$G_\Z^{\BB}/Z((\{\pm 1\}^3)_{\uuuu{x}}) =G_\Z/\{\pm \id\}
\cong SL_2(\Z)$$ 
is big. There is a natural bijection
$\BB^{dist}/\{\pm 1\}^3\to SL_2(\Z)$,
and the elementary braids $\sigma_1$ and $\sigma_2$ act
by multiplication from the left with the matrices 
$A_1=\begin{pmatrix} 1 &-1\\0&1\end{pmatrix}$ and
$A_2=\begin{pmatrix} 1&0\\1&1\end{pmatrix}$ on $SL_2(\Z)$.
This gives a clear description of the $\sigma$-pseudo-graph
$\GG_\sigma(\BB^{dist}/\{\pm 1\}^3)$. We do not attempt
a picture.

(iii) The reducible case $A_1^3$, $\uuuu{x}=(0,0,0)$:
Also here $\Br_3(\uuuu{x}/\{\pm 1\}^3)=\{\uuuu{x}/\{\pm 1\}^3\}$
has only one vertex. The group
$$ G_\Z^{\BB}/Z((\{\pm 1\}^3)_{\uuuu{x}})
=G_\Z/\{\pm 1\}^3\cong O_3(\Z)/\{\pm 1\}^3\cong S_3$$
has six elements. 
Therefore $\BB^{dist}/\{\pm 1\}^3$ has six elements,
and $\sigma_1$ and $\sigma_2$ act as involutions.
The right hand side of the first line in Figure 
\ref{Fig:7.3} gives the $\sigma$-pseudo-graph
$\GG_\sigma(\BB^{dist}/\{\pm 1\}^3)$. 
Part (iv) offers a different description which
applies also to $A_1^3$ if one sees it as $A_1^2A_1$.

(iv) The reducible cases, 
$\uuuu{x}\in \bigcup_{i\in\{3,4,5\}}C_i$
($A_2A_1$, $\P^1A_1$, other reducible cases):
Here $(H_\Z,L,\uuuu{e})=(H_{\Z,1},L_1,(e_1,e_2))\oplus
(H_{\Z,2},L_2,e_3)$ with $H_{\Z,1}=\Z e_1\oplus \Z e_2$
and $H_{\Z,2}=\Z e_3$. 

The group of deck transformations of the normal covering
$\GG_\sigma(\BB^{dist}/\{\pm 1\}^3)\to
\GG_\sigma(\Br_3(\uuuu{x}/\{\pm 1\}^3))$
is 
$$ G_\Z^{\BB}/Z((\{\pm 1\}^3)_{\uuuu{x}})
=G_\Z/\{\pm \id,\pm Q\}\cong \{(M^{root})^l\,|\, l\in\Z\}.$$
Here $M^{root}$ has
order 3 in the case $A_2A_1$ and infinite order in the
other cases. Therefore the $\sigma$-pseudo-graph 
$\GG_\sigma(\BB^{dist}/\{\pm 1\}^3)$
can be obtained by a triple or infinite covering of the
$\sigma$-pseudo-graph $\GG_\sigma(\Br_3(\uuuu{x}/\{\pm 1\}^3)$.

\begin{figure}[H]
\includegraphics[width=0.6\textwidth]{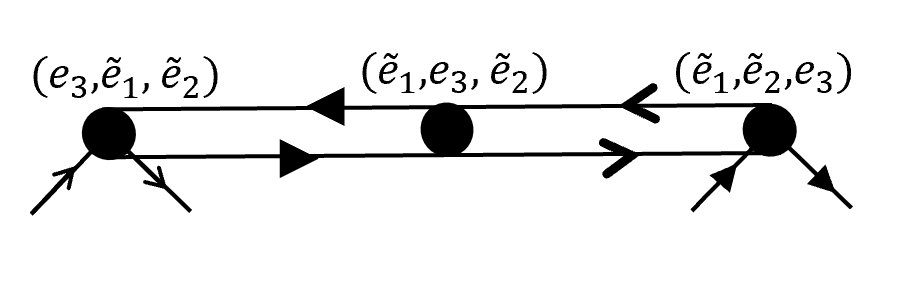}
\caption[Figure 7.2]{In the reducible cases (without $A_1^3$)
one sheet of the covering
$\GG_\sigma(\BB^{dist}/\{\pm 1\}^3)\to 
\GG_\sigma(\Br_3(\uuuu{x}/\{\pm 1\}^3))$}
\label{Fig:7.2}
\end{figure}

More concretely, the type of the covering is determined by the 
$\Br_1$ orbit of distinguished bases up to signs of 
$(H_{\Z,1},L_1,(e_1,e_2))$. 
One such distinguished basis modulo signs
$(\www{e}_1,\www{e}_2)/\{\pm 1\}^2$ gives rise to one sheet
in the covering 
$\GG_\sigma(\BB^{dist}/\{\pm 1\}^3)\to 
\GG_\sigma(\Br_3(\uuuu{x}/\{\pm 1\}^3))$.
Figure \ref{Fig:7.2} shows the part of a $\sigma$-pseudo-graph
which corresponds to one such sheet. 

\begin{figure}
\includegraphics[height=0.2\textheight]{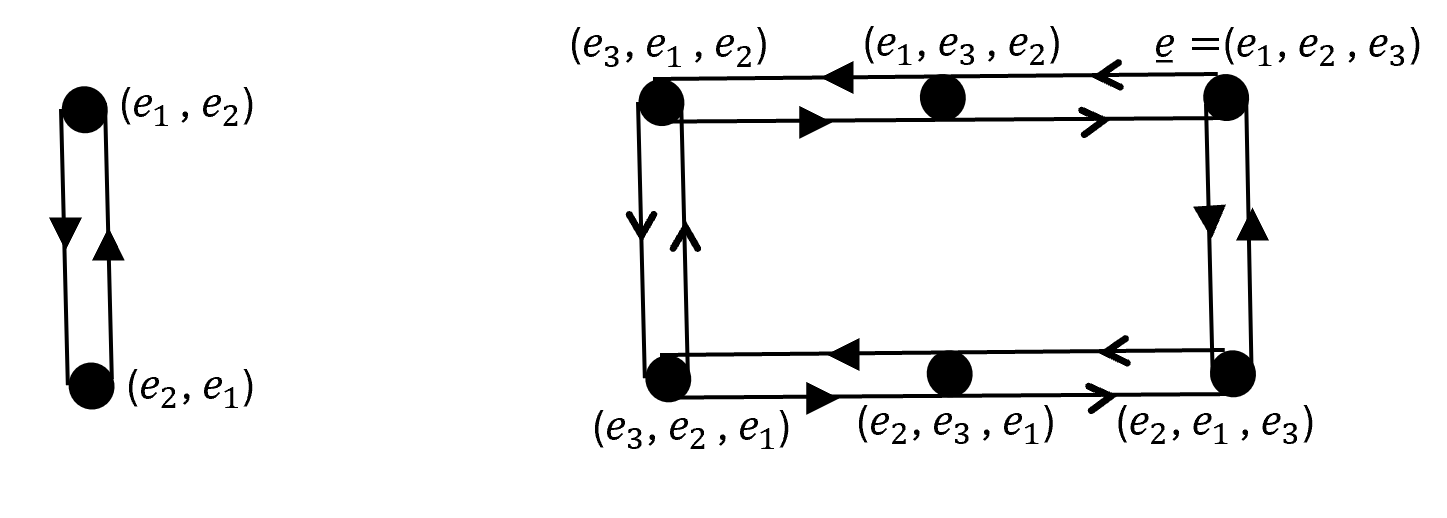}\\
$A_1^2$ \hspace*{7cm} $A_1^3$ \\
\hspace*{-0.5cm}\includegraphics[height=0.2\textheight]{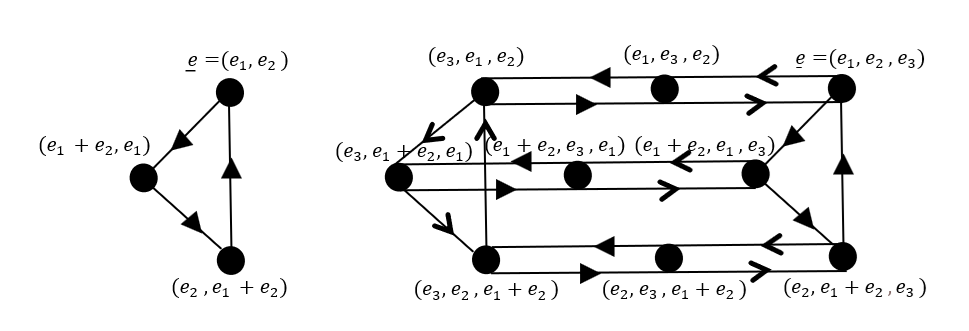}\\
$A_2$ \hspace*{7cm}$A_2A_1$ \\
\hspace*{1cm}\includegraphics[height=0.3\textheight]
{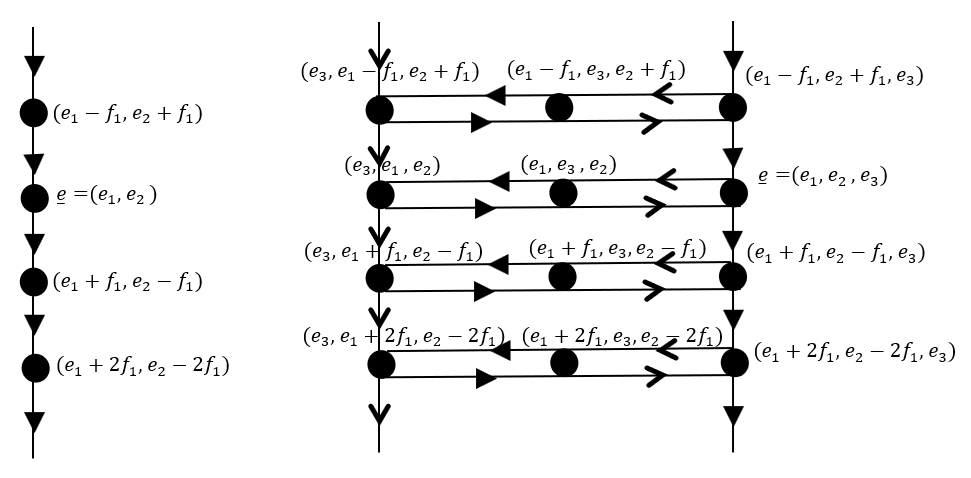}
$\P^1$, $f_1=e_1+e_2$ \hspace*{5cm} $\P^1A_1$ \\
\caption[Figure 7.3]{Examples \ref{t7.16} (iii): The $\sigma$-pseudo-graphs for
distinguished bases modulo signs in the reducible cases}
\label{Fig:7.3}
\end{figure}

The six pictures in Figure \ref{Fig:7.3} 
show on the left hand side analogous
$\sigma_1$-pseudo-graphs for the distinguished bases
modulo signs of the rank 2 cases
$A_1^2$, $A_2$ and $\P^1$, and on the right hand side the
$\sigma$-pseudo-graphs 
$\GG_\sigma(\BB^{dist}/\{\pm 1\}^3)$ 
for $A_1^3$, $A_2A_1$ and $\P^1A_1$
(respectively only a part of the $\sigma$-pseudo-graph
in the case of $\P^1A_1)$. 
The $\sigma$-pseudo-graph 
$\GG_\sigma(\BB^{dist}/\{\pm 1\}^3)$ for $\uuuu{x}=(x,0,0)$
with $x<-2$ looks the same as the one for $\P^1A_1$, though
of course the distinguished bases are different.

(v) The case $A_3$, $\uuuu{x}=(-1,0,-1)$: 
The $\sigma$-pseudo-graph 
$\GG_\sigma(\BB^{dist}/\{\pm 1\}^3)$
was first given in \cite[page 40, Figur 6]{Yu90}. 
We recall and explain it in our words.
The group of deck transformations of the normal covering
$\GG_\sigma(\BB^{dist}/\{\pm 1\}^3)\to
\GG_\sigma(\Br_3(\uuuu{x}/\{\pm 1\}^3))$
is 
$$ G_\Z^{\BB}/Z((\{\pm 1\}^3)_{\uuuu{x}})
=G_\Z/\{\pm \id\}\cong \{M^l\,|\, l\in\{0,1,2,3\}\}.$$
Here the monodromy $M$ acts in the natural way,
$$M((\www{e}_1,\www{e}_2,\www{e}_3)/\{\pm 1\}^3)
=(M(\www{e}_1),M(\www{e}_2),M(\www{e}_3))/\{\pm 1\}^3),$$
on $\BB^{dist}/\{\pm 1\}^3$ and has order four,
$M^4=\id$. Here $M$ and its powers are
\begin{eqnarray*}
M(\uuuu{e})=\uuuu{e}
\begin{pmatrix}0&0&1\\ -1&0&1\\0&-1&1\end{pmatrix},
M^2(\uuuu{e})=\uuuu{e}
\begin{pmatrix}0&-1&1\\0&-1&0\\1&-1&0\end{pmatrix}, \\
M^3(\uuuu{e})=\uuuu{e}
\begin{pmatrix}1&-1&0\\1&0&-1\\1&0&0\end{pmatrix}.
\end{eqnarray*}

Because of the shape of 
$\GG_\sigma(\Br_3(\uuuu{x}/\{\pm 1\}^3))$,
$$b^{1,0}:=\uuuu{e}/\{\pm 1\}^3,\quad
b^{2,0}:=\sigma_1 b^{1,0},\quad
b^{3,0}:=\sigma_1^2 b^{1,0},\quad
b^{4,0}:=\sigma_2^{-1} b^{1,0}$$
form one sheet of the fourfold cyclic covering
$\GG_\sigma(\BB^{dist}/\{\pm 1\}^3)\to
\GG_\sigma(\Br_3(\uuuu{x}/\{\pm 1\}^3))$.
Define
$$b^{i,l}:=M^l(b^{i,0})\quad \textup{for}\quad
l\in\{1,2,3\}.$$
Then
$$\BB^{dist}/\{\pm 1\}^3
=\{b^{i,l}\,|\, i\in\{1,2,3,4\},l\in\{0,1,2,3\}\}$$
has sixteen elements. We claim for $l\in\{0,1,2,3\}$
\begin{eqnarray*}
\begin{array}{ll}
\sigma_1 b^{1,l}=b^{2,l},& 
\sigma_2 b^{1,l}=b^{3,l+3(\mmod 4)},\\
\sigma_1 b^{2,l}=b^{3,l},&
\sigma_2 b^{2,l}=b^{2,l+2(\mmod 4)},\\
\sigma_1 b^{3,l}=b^{1,l},&
\sigma_2 b^{3,l}=b^{4,l+1(\mmod 4)},\\
\sigma_1 b^{4,l}=b^{4,l+2(\mmod 4)},&
\sigma_2 b^{4,l}=b^{1,l}.
\end{array}
\end{eqnarray*}
It is sufficient to prove the claim for $l=0$.
The equations
$\sigma_1 b^{1,0}=b^{2,0}$,  
$\sigma_1 b^{2,0}=b^{3,0}$,
$\sigma_2 b^{4,0}=b^{1,0}$ 
follow from the definitions of 
$b^{2,0}$, $b^{3,0}$, $b^{4,0}$.
The inclusion $\sigma_1^3\in (\Br_3)_{\uuuu{e}/\{\pm 1\}^3})$
gives $\sigma_1 b^{3,0}=b^{1,0}$. It remains to show
$$\sigma_1 b^{4,0}=b^{4,2},\\
\sigma_2 b^{1,0}=b^{3,3},\ 
\sigma_2 b^{2,0}=b^{2,2},\
\sigma_2 b^{3,0}=b^{4,1}.$$
One sees
\begin{eqnarray*}
\begin{array}{llll}
i & b^{i,0} & \www{\uuuu{e}} & \www{\uuuu{x}}
\textup{ with }L(\www{\uuuu{e}}^t,\www{\uuuu{e}})
=S(\www{\uuuu{x}})\\ \hline 
1 & \uuuu{e}/\{\pm \}^3 & \uuuu{e} & (-1,0,-1)\\
2 & \sigma_1(\uuuu{e})/\{\pm 1\}^3 & 
\sigma_1(\uuuu{e})=(e_1+e_2,e_1,e_3) & (1,-1,0) \\
3 & \sigma_1^2(\uuuu{e})/\{\pm 1\}^3 & 
\sigma_1^2(\uuuu{e})=(-e_2,e_1+e_2,e_3)
& (-1,1,-1)\\
4 & \sigma_2^{-1}(\uuuu{e})/\{\pm 1\}^3 & 
\sigma_2^{-1}(\uuuu{e})=(e_1,e_3,e_2+e_3) & (0,-1,1)
\end{array}
\end{eqnarray*}
\begin{eqnarray*}
\sigma_1 b^{4,0}&=& (e_3,e_1,e_2+e_3)/\{\pm 1\}^3
=M^2b^{4,0}=b^{4,2},\\
\sigma_2 b^{1,0}&=& (e_1,e_2+e_3,e_2)/\{\pm 1\}^3
=M^3b^{3,0}=b^{3,3},\\
\sigma_2 b^{2,0}&=& (e_1+e_2,e_3,e_1)/\{\pm 1\}^3
=M^2b^{2,0}=b^{2,2}.
\end{eqnarray*}
$\sigma_2b^{1,0}=b^{3,3}$, $\sigma_2b^{4,0}=b^{1,0}$
and $\sigma_2^3\in (\Br_3)_{\uuuu{e}/\{\pm 1\}^3}$ show
$\sigma_2b^{3,3}=b^{4,0}$.
This implies $\sigma_2b^{3,0}=b^{4,1}$.  
The claim is proved. The $\sigma$-pseudo-graph 
$\GG_\sigma(\BB^{dist}/\{\pm 1\}^3)$ is given in Figure \ref{Fig:7.4}.

\begin{figure}
\includegraphics[width=0.7\textwidth]{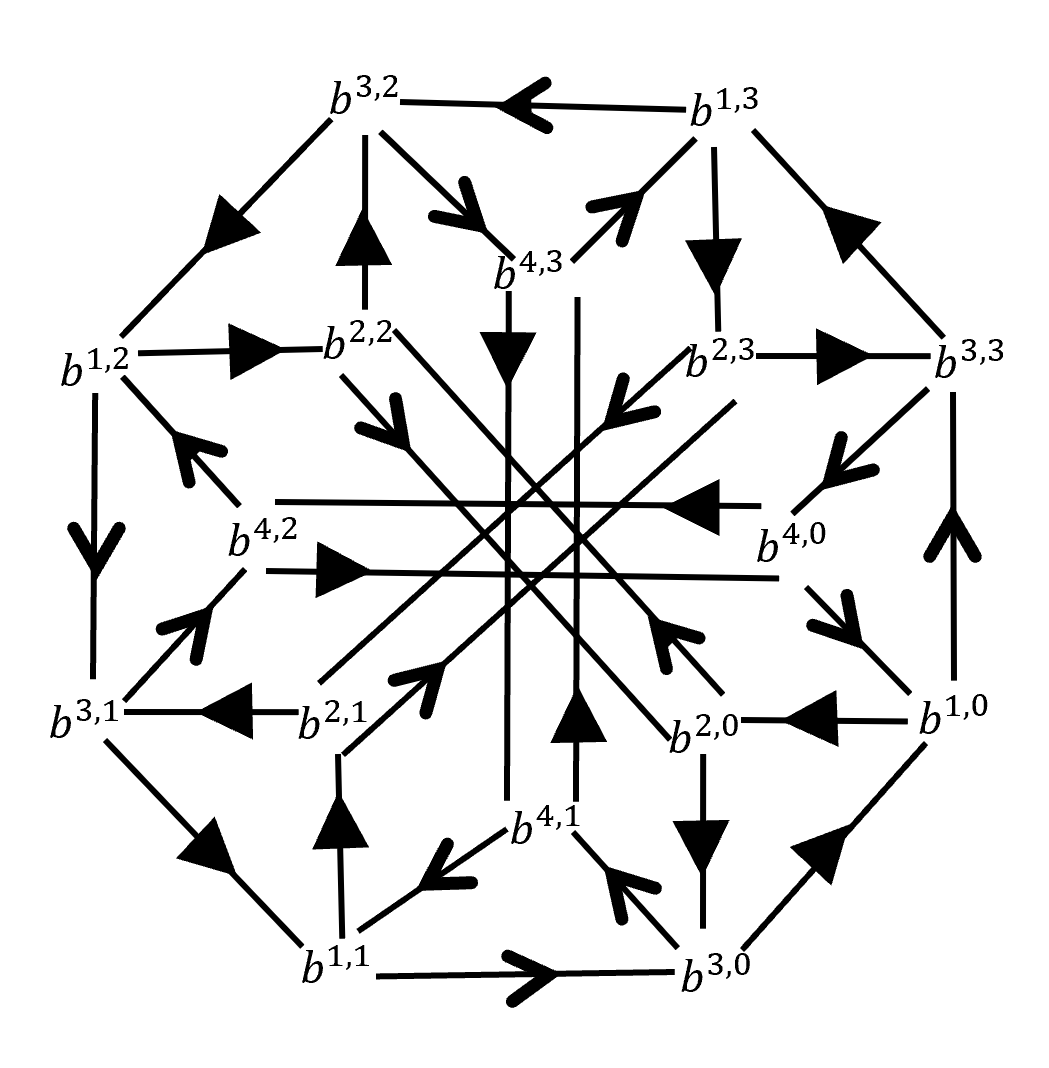}
\caption[Figure 7.4]{Example \ref{t7.16} (iv): The $\sigma$-pseudo-graph for
distinguished bases modulo signs in the case $A_3$}
\label{Fig:7.4}
\end{figure}

(vi) The case $\whh{A}_2$, $\uuuu{x}=(-1,-1,-1)$: 
The group of deck transformations of the normal covering
$\GG_\sigma(\BB^{dist}/\{\pm 1\}^3)\to
\GG_\sigma(\Br_3(\uuuu{x}/\{\pm 1\}^3))$
is 
$$ G_\Z^{\BB}/Z((\{\pm 1\}^3)_{\uuuu{x}})
=G_\Z/\{\pm \id\}\cong \{(M^{root})^l\,|\, l\in\Z\}.$$
Here $M^{root}$ acts in the natural way
on $\BB^{dist}/\{\pm 1\}^3$.
$M^{root}$ has infinite order and satisfies $(M^{root})^3=-M$.
Recall $f_1=e_1+e_2+e_3$, $\Z f_1=\Rad I^{(0)}$, 
\begin{eqnarray*}
M^{root}(\uuuu{e})&=&\uuuu{e}
\begin{pmatrix}1&1&-1\\1&0&0\\0&1&0\end{pmatrix},\quad
M^{root}(f_1)=f_1,\\
(M^{root})^2(\uuuu{e})&=&\uuuu{e}+f_1(1,0,-1).
\end{eqnarray*}

Because of the shape of 
$\GG_\sigma(\Br_3(\uuuu{x}/\{\pm 1\}^3))$,
$$b^{1,0}:=\uuuu{e}/\{\pm 1\}^3,\quad
b^{2,0}:=\sigma_1 b^{1,0},\quad
b^{3,0}:=\sigma_1^2 b^{1,0},\quad
b^{4,0}:=\sigma_2 b^{1,0}$$
form one sheet of the infinite cyclic covering
$\GG_\sigma(\BB^{dist}/\{\pm 1\}^3)\to
\GG_\sigma(\Br_3(\uuuu{x}/\{\pm 1\}^3))$.
Define
$$b^{i,l}:=(M^{root})^l(b^{i,0})\quad \textup{for}\quad
i\in\{1,2,3,4\},\ l\in\Z-\{0\}.$$
Then
$$\BB^{dist}/\{\pm 1\}^3
=\{b^{i,l}\,|\, i\in\{1,2,3,4\},l\in\Z\}\}.$$
We claim for $l\in\Z$
\begin{eqnarray*}
\begin{array}{ll}
\sigma_1 b^{1,l}=b^{2,l},&
\sigma_2 b^{1,l}=b^{4,l},\\
\sigma_1 b^{2,l}=b^{3,l},&
\sigma_2 b^{2,l}=b^{1,l+1},\\
\sigma_1 b^{3,l}=b^{1,l},&
\sigma_2 b^{3,l}=b^{3,l+2},\\
\sigma_1 b^{4,l}=b^{4,l+2},&
\sigma_2 b^{4,l}=b^{2,l-1}.
\end{array}
\end{eqnarray*}
It is sufficient to prove the claim for $l=0$.
The equations
$\sigma_1 b^{1,0}=b^{2,0}$,  
$\sigma_1 b^{2,0}=b^{3,0}$,
$\sigma_2 b^{1,0}=b^{4,0}$ 
follow from the definitions of 
$b^{2,0}$, $b^{3,0}$, $b^{4,0}$.
The inclusion $\sigma_1^3\in (\Br_3)_{\uuuu{e}/\{\pm 1\}^3}$
gives $\sigma_1 b^{3,0}=b^{1,0}$. It remains to show
$$\sigma_1 b^{4,0}=b^{4,2},\\
\sigma_2 b^{2,0}=b^{1,1},\ 
\sigma_2 b^{3,0}=b^{3,2},\
\sigma_2 b^{4,0}=b^{2,-1}.$$
One sees 
\begin{eqnarray*}
\begin{array}{llll}
i & b^{i,0} & \www{\uuuu{e}} & \www{\uuuu{x}}
\textup{ with }L(\www{\uuuu{e}}^t,\www{\uuuu{e}})
=S(\www{\uuuu{x}})\\ \hline 
1 & \uuuu{e}/\{\pm \}^3 & \uuuu{e} & (-1,-1,-1)\\
2 & \sigma_1(\uuuu{e})/\{\pm 1\}^3 & 
\sigma_1(\uuuu{e})=(e_1+e_2,e_1,e_3) & (1,-2,-1) \\
3 & \sigma_1^2(\uuuu{e})/\{\pm 1\}^3  & 
\sigma_1^2(\uuuu{e})=(-e_2,e_1+e_2,e_3)
& (-1,1,-2)\\
4 & \sigma_2(\uuuu{e})/\{\pm 1\}^3 & 
\sigma_2(\uuuu{e})=(e_1,e_2+e_3,e_2) & (-2,-1,1)
\end{array}
\end{eqnarray*}
\begin{eqnarray*}
\sigma_1 b^{4,0}&=& (2e_1+e_2+e_3,e_1,e_2)/\{\pm 1\}^3\\
&=&(e_1+f_1,e_2+e_3-f_1,e_2)/ \{\pm 1\}^3
=(M^{root})^2b^{4,0}=b^{4,2},\\
\sigma_2 b^{2,0}&=& (e_1+e_2,e_1+e_3,e_1)/\{\pm 1\}^3
=M^{root}b^{1,0}=b^{1,1},\\
\sigma_2 b^{3,0}&=& (-e_2,2e_1+2e_2+e_3,e_1+e_2)/\{\pm 1\}^3\\
&=&(-e_2,e_1+e_2+f_1,-e_3+f_1)/\{\pm 1\}^3
=M^2b^{3,0}=b^{3,2}.
\end{eqnarray*}
$\sigma_2b^{2,0}=b^{1,1}$ implies $\sigma_2b^{2,-1}=b^{1,0}$.
This, $\sigma_2b^{1,0}=b^{4,0}$ and 
$\sigma_2^3\in (\Br_3)_{\uuuu{e}/\{\pm 1\}^3}$ show
$\sigma_2b^{4,0}=b^{2,-1}$. The claim is proved.
A part of the $\sigma$-pseudo-graph 
$\GG_\sigma(\BB^{dist}/\{\pm 1\}^3)$ is given in Figure
\ref{Fig:7.5}.
\begin{figure}
\includegraphics[width=1.0\textwidth]{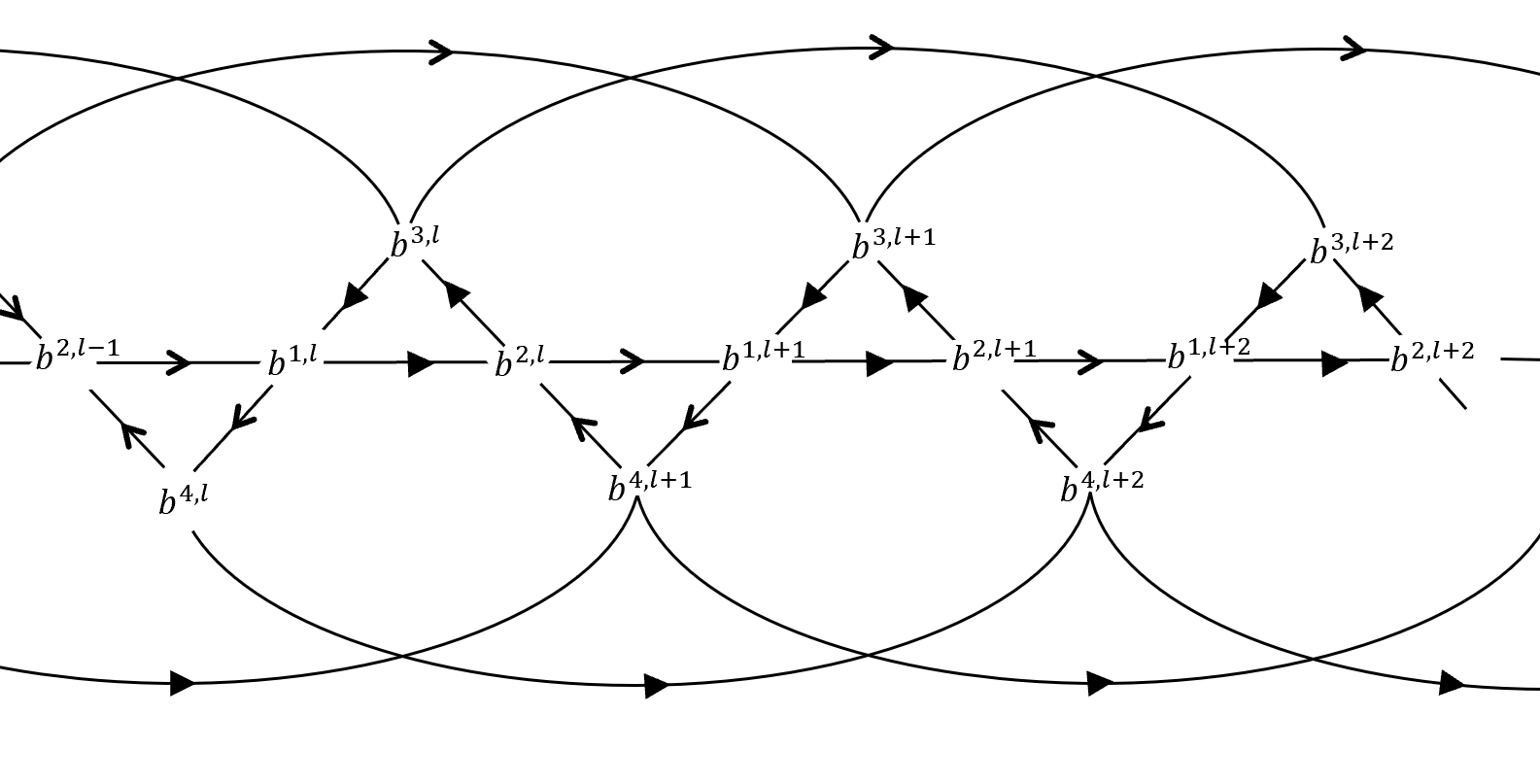}
\caption[Figure 7.5]{Example \ref{t7.16} (v): 
A part of the $\sigma$-pseudo-graph for
distinguished bases modulo signs in the case $\whh{A}_2$}
\label{Fig:7.5}
\end{figure}

\end{examples}

\chapter[Manifolds induced by braid group orbits]
{Manifolds induced by the orbit of distinguished
bases and the orbit of distinguished matrices}\label{s8}
\setcounter{equation}{0}
\setcounter{figure}{0}

A single matrix $S\in T^{uni}_n(\Z)$ induces a unimodular
bilinear lattice $(H_\Z,L,\uuuu{e})$ with triangular basis
$\uuuu{e}$. This triple is unique up to isomorphism. 
The chapters \ref{s2}--\ref{s7} studied this structure 
and the associated structures $\Gamma^{(k)}$, 
$\Delta^{(k)}$ and $\BB^{dist}$ in
general, and also systematically in the cases of rank 2 and 3.
These chapters are all of an algebraic/combinatorial type.

Chapter \ref{s8} comes to geometry. 
It will present two complex $n$-dimensional
manifolds $C_n^{\uuuu{e}/\{\pm 1\}^n}$ and $C_n^{S/\{\pm 1\}^n}$
which are induced by $S$ respectively by 
$(H_\Z,L,\uuuu{e})$. 
If $n\geq 2$ they are certain coverings of the configuration
space of the braid group $\Br_n$. 
They can also be seen as the results of gluing many 
{\it Stokes regions}, one for each element in
$\BB^{dist}/\{\pm 1\}^n$ respectively in
$\SSS^{dist}/\{\pm 1\}^n$. 
Their definition is not difficult, but worth to be 
discussed in detail. This is done in section \ref{s8.1}. 

Many of them are important in algebraic geometry and singularity 
theory. See chapter \ref{s10} for statements on the cases 
in singularity theory. 
They all are carriers of much richer structures.
One version of these richer structures is presented in
section \ref{s8.5}, certain {\it $\Z$-lattice bundles}
on them. Another rather elementary differential geometric
structure on them is given in section \ref{s8.2}.

Section \ref{s8.2} recalls the notion of an F-manifold
and states that $C_n^{\uuuu{e}/\{\pm 1\}^n}$ and
$C_n^{S/\{\pm 1\}^n}$ are semisimple F-manifolds with
Euler field and empty Maxwell stratum. 
As such, often they have partical compactifications to which
the F-manifold structure extends.

Section \ref{s8.3} considers first the case when 
$(H_\Z,L,\uuuu{e})$ is reducible. It states that then
$C_n^{\uuuu{e}/\{\pm 1\}^n}$ embeds into the product of the
corresponding manifolds for the summands of 
$(H_\Z,L,\uuuu{e})$. Afterwards it discusses the rank 2 cases.

Section \ref{s8.4} recalls a notion which was developed
within singularity theory, the {\it distinguished systems of
$n$ paths}. There are a priori two actions of the braid group
$\Br_n$ on sets of (homotopy classes) of distinguished systems
of paths. They are discussed and compared in section \ref{s8.4}.

They are needed in section \ref{s8.5}. It constructs
from $(H_\Z,L,\uuuu{e})$ natural families of $\Z$-lattices
structures over the manifolds $C_n^{univ}$ and
$C_n^{\uuuu{e}/\{\pm 1\}^n}$. Before, it constructs a single
$\Z$-lattice structure on $\C-\{u_1,...,u_n\}$ from
$(H_\Z,L,\uuuu{e})$ and the additional choices of 
a pair $(\uuuu{u},r)\in C_n^{pure}\times\R$ with $r$ big enough
and of a distinguished system of paths with starting point $r$
and endpoints in $\{u_1,...,u_n\}$. 

The $\Z$-lattice bundle over $C_n^{\uuuu{e}/\{\pm 1\}^n}$
leads in general (with an additional choice) 
to a {\it Dubrovin-Frobenius manifold structure} 
on a Zariski open subset of $C_n^{\uuuu{e}/\{\pm 1\}^n}$. 
Remark \ref{t8.9} (ii) offers one definition of a 
Dubrovin-Frobenius manifold and references.
We do not go more into this big story.

\section[Natural coverings of the configuration space]
{Natural coverings of the configuration space
of the braid group $\Br_n$}
\label{s8.1}

First we fix some notations and recall some facts from the 
theory of covering spaces and some facts around the
configuration space of the braid group $\Br_n$. After this  
the Remarks \ref{t8.4} start the discussion of the manifolds
$C_n^{\uuuu{e}/\{\pm 1\}^n}$ and $C_n^{S/\{\pm 1\}^n}$.

\begin{definition}\label{t8.1}
Let $n\in \Z_{\geq 2}$. Define
\index{$C_n^{conf},\ C_n^{pure},\ C_n^{univ}$}
\index{$D_n^{conf},\ D_n^{pure},\ D_n^{univ}$} 
\index{$b_n^{conf},\ b_n^{pure},\ b_n^{univ}$}
\begin{eqnarray*}
D_n^{pure}&:=& \bigcup_{1\leq i<j\leq n}
\{\uuuu{u}\in\C^n\,|\, u_i=u_j\}\subset \C^n,\\
&&\textup{the union of the partial diagonals in }\C^n,\\
C_n^{pure}&:=& \C^n-D_n^{pure},\\
C_n^{conf}&:=& C_n^{pure}/S_n \subset \C^n/S_n\\
D_n^{conf}&:=& D_n^{pure}/S_n \subset \C^n/S_n,\quad 
\textup{so }\C^n/S_n=C_n^{conf}\, \dot\cup\, 
D_n^{conf},\\
C_n^{univ}&:=&(\textup{the universal covering of }
C_n^{pure}\textup{ and }C_n^{conf}),\\
b_n^{pure}&:=& (i,2i,...,ni)\in C_n^{pure},\\
b_n^{conf}&:=&[b_n^{pure}]=b_n^{pure}/S_n,\\
b_n^{univ}&:=&(\textup{the preimage of }b_n^{pure}
\textup{ in }C_n^{univ}
\textup{ which corresponds}\\
&&\textup{to the trivial path in }C_n^{pure}).
\end{eqnarray*}
\index{universal covering}
$C_n^{conf}$ is called {\it configuration space} of the braid
group $\Br_n$. The three covering maps are denoted as follows,
\index{$\pr_n^{p,c},\ \pr_n^{u,p},\ \pr_n^{u,c}$}
\begin{eqnarray*}
\pr_n^{p,c}:C_n^{pure}\to C_n^{conf},&&\\
\pr_n^{u,p}:C_n^{univ}\to C_n^{pure},&&\\
\pr_n^{u,c}=\pr_n^{p,c}\circ\pr_n^{u,p}:C_n^{univ}\to C_n^{conf}.&&
\end{eqnarray*}
\end{definition}

\begin{remarks}\label{t8.2}
(i) The permutation group $S_n$ acts from the left on
$\C^n$ by the following action for $\sigma\in S_n$,
\begin{eqnarray*}
\sigma(\uuuu{u})&=& \sigma((u_1,u_2,...,u_n))
:=(u_{\sigma^{-1}(1)},u_{\sigma^{-1}(2)},...,
u_{\sigma^{-1}(n)}).
\end{eqnarray*}
So the entry of $\uuuu{u}$ at the $\sigma^{-1}(j)$-th place
is put to the $j$-th place. If $\tau\in S_n$, then
$\sigma(\tau(\uuuu{u}))$ has at the $j$-th place the
entry $u_{\tau^{-1}(\sigma^{-1}(j))}$ of $\tau(\uuuu{u})$
at the $\sigma^{-1}(j)$-th place, so the entry
$u_{(\sigma\circ\tau)^{-1}(j)}$. Therefore this is an 
action from the left of $S_n$ on $\C^n$. 
Although this is an action from the left, we write in the
quotient $\C^n/S_n$ the group on the right.

(ii) The quotient $\C^n/S_n$ is again isomorphic to $\C^n$.
It can be identified with the set
\begin{eqnarray*}
\C[x]_n:=\{f(x)=x^n+\sum_{j=1}^nf_jx^{n-j}\in\C[x]\,|\, 
f=(f_1,...,f_n)\in\C^n\}
\end{eqnarray*}
of complex unitary polynomials of degree $n$. The identification
is given by the map 
$[\uuuu{u}]=[(u_1,...,u_n)]\mapsto \prod_{i=1}^n(x-u_i)$. 

The covering $C_n^{pure}\to C_n^{conf}$ extends to a branched
covering $\C^n\to \C^n/S_n$. It maps each partial diagonal
$\{\uuuu{u}\in\C^n\,|\, u_i=u_j\}$ for $1\leq i<j\leq n$
to the same irreducible algebraic hypersurface
$D_n^{conf}=D_n^{pure}/S_n\subset \C^n/S_n$. 
This hypersurface $D_n^{conf}$ and its image 
$$\C[x]_n^{mult}:=\{f(x)\in\C[x]_n\,|\, f(x) 
\textup{ has multiple roots}\}$$
\index{$\C[x]_n\, \C[x]_n^{mult},\ \C[x]_n^{reg}$} 
in $\C[x]_n$ are called {\it discriminants}. \index{discriminant}
The configuration space $C_n^{conf}$ is identified with the
complement
$$\C[x]_n^{reg}:=\{f(x)\in\C[x]_n\,|\, f(x) 
\textup{ has no multiple roots}\}.$$

(iii) The braid group $\Br_n$ is isomorphic to the fundamental
\index{braid group}
group $\pi_1(C_n^{conf},b_n^{conf})$
\index{fundamental group}
(e.g. \cite[Theorem 1.8]{Bi75}, \cite[1.4.3]{KT08}).
A closed loop $\gamma:[0,1]\to C_n^{conf}$ with 
$\gamma(0)=\gamma(1)=b_n^{conf}$ can be visualized by a braid
with $n$ strings. Figure \ref{Fig:8.1} shows five braids with three
strings which represent the braids
$\sigma_1,\sigma_1^{-1},\sigma_2,\sigma_2^{-1},
\sigma_2^{-1}\sigma_1$. The interval $[0,1]$ is going top down.
The starting points and the end points of the strings are
$i,2i,...,ni\in\C$ where $\C$ 
is the complex line in which the entries of the loop live.
Recall that in a fundamental group a product of elements
is gone through from left to right. Therefore in the braid
$\sigma_2^{-1}\sigma_1$ in Figure \ref{Fig:8.1}  
$\sigma_2^{-1}$ comes first and $\sigma_1$ comes second.

\begin{figure}
\includegraphics[width=0.9\textwidth]{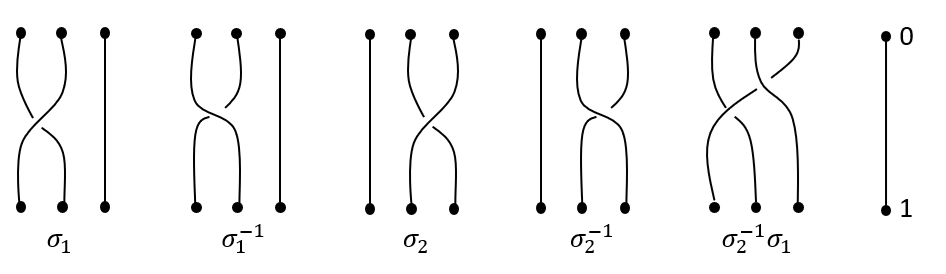}
\caption[Figure 8.1]{Five braids}
\label{Fig:8.1}
\end{figure}

The relations 
\begin{eqnarray*}
\sigma_j\sigma_{j+1}\sigma_j&=&\sigma_{j+1}\sigma_j\sigma_{j+1}
\quad\textup{for }j\in\{,...,n-1\}\\ \textup{and}\qquad  
\sigma_j\sigma_k&=&\sigma_k\sigma_j\quad\textup{for }|j-k|>1
\end{eqnarray*}
are now obvious. It is true, but not obvious that they are 
generating relations in $\Br_n$ for the generators 
$\sigma_1,...,\sigma_{n-1}$. 

(iv) $C_n^{conf}$ is a connected complex manifold and especially
locally simply connected. Therefore the main results of the
theory of covering spaces apply. The {\it Galois correspondence}
\index{covering}\index{connected covering}\index{Galois correspondence}
for connected coverings tells the following 
\cite[Theorem 1.38]{Ha01}.

Each connected covering 
$p:(Y,y_0)\to (C_n^{conf},b_n^{conf})$
with base point $y_0\in Y$ with $p(y_0)=b_n^{conf}$ gives rise
to an injective group homomorphism
$p_*:\pi_1(Y,y_0)\to \pi_1(C_n^{conf},b_n^{conf})$.
Vice versa, for each subgroup $U\subset \Br_n$, there is such
a connected covering $p:(Y,y_0)\to(C_n^{conf},b_n^{conf})$
with $p_*(\pi_1(Y,y_0))=U$, and it is unique up to 
isomorphism of connected coverings with base points.
This gives a 1--1 correspondence between connected coverings
of $C_n^{conf}$ with base points up to isomorphism and
subgroups of $\Br_n$. If one forgets the base points, one
obtains a 1--1 correspondence between connected coverings and
conjugacy classes of subgroups of $\Br_n$. 

(v) The group of deck transformations 
\index{deck transformation}
of a connected covering
$p:(Y,y_0)\to(C_n^{conf},b_n^{conf})$ is the group of
automorphisms of $Y$ which map each fiber of $p$ to itself.
A deck transformation $f:Y\to Y$ is determined by the one
value $f(y_0)\in p^{-1}(b_n^{conf})$. The group of 
deck transformations is isomorphic to the quotient
$$N(p_*(\pi_1(Y,y_0)))/p_*(\pi_1(Y,y_0)),$$
where $N(p_*(\pi_1(Y,y_0)))\subset \Br_n$ is the largest 
subgroup of $\Br_n$ in which $p_*(\pi_1(Y,y_0))$ is a normal
subgroup. 

The covering is normal if $N(p_*(\pi_1(Y,y_0)))=\Br_n$.
This holds if and only if the group of deck transformations
acts transitively on the fiber $p^{-1}(b_n^{conf})$
(and then transitively on each fiber of $p$).

For a given subgroup $U\subset \Br_n$, the covering 
$p:(Y,y_0)\to (C_n^{conf},b_n^{conf})$ can be constructed
as the quotient $C_n^{univ}/U$ where $U$ acts from the
left as subgroup of the group $\Br_n$, which acts as group of
deck transformations of $C_n^{univ}$. 

(vi) The universal covering \index{universal covering}
$\pr_n^{u,c}:C_n^{univ}\to C_n^{conf}$ has trivial 
fundamental group $\{\id\}\subset \Br_n$ which is normal in
$\Br_n$, so it is a normal covering with $\Br_n$ as group
of deck transformations. One can identify $C_n^{univ}$ as a set
with the following set,
\begin{eqnarray*}
\bigcup_{b\in C_n^{conf}}
\{\textup{homotopy classes of paths in }C_n^{conf}
\textup{ from }b_n^{conf}\textup{ to }b\}.
\end{eqnarray*}
$\Br_n$ acts on this set from the left as follows:
a braid $\alpha$ maps a homotopy class $\beta$ to the
homotopy class $\alpha\beta$, where one first goes along a 
loop which represents $\alpha$ and then along a path which
represents $\beta$. Remarkably, this gives an action {\it from
the left} of $\Br_n$ on $C_n^{univ}$, although in $\Br_n$
in a product of two loops one walks first along the left loop,
whereas in a product of two automorphisms of $C_n^{univ}$ the
right automorphism is carried out first.

(vii) The covering $\pr_n^{p,c}:C_n^{pure}\to C_n^{conf}$
is normal of degree $n!$ with group of deck transformations
isomorphic to $S_n$. The fundamental group 
$\pi_1(C_n^{pure},b_n^{pure})$ embeds into $\Br_n$ as group
$\Br_n^{pure}$ of {\it pure braids}. 
\index{pure braid}
A pure braid is
represented by a closed loop in $C_n^{pure}$ which starts
and ends at $b_n^{pure}$. So in a pure braid with $n$ strings,
each string ends at the point at which it starts.
There is an exact sequence
$$\{1\}\to \Br_n^{pure}\to\Br_n\to S_n\to\{1\}$$
of group homomorphisms. 
$\Br_n^{pure}$ is the kernel of the map
$\Br_n\to S_n$, so normal in $\Br_n$, so
$N(\Br_n^{pure})=\Br_n$. 
Here given a braid $\alpha\in \Br_n$,
that permutation $\sigma\in S_n$ is associated to it, 
such that for a closed loop in $C_n^{conf}$ which represents
$\alpha$, the lift to $C_n^{pure}$ which starts at 
$b_n^{pure}=(i,2i,...,ni)$ ends at 
$(\sigma^{-1}(1)i,\sigma^{-1}(2)i,...,\sigma^{-1}(n)i)$. 
This fits to the action of $S_n$ on $\C^n$ in part (i) of
these remarks. 

(viii) Arnold \cite{Ar68} showed that the universal covering
$C_n^{univ}$ is homeomorphic to an open ball in $\R^{2n}$.
Kaliman \cite{Ka75} \cite{Ka93} showed that $C_n^{univ}$ is
biholomorphic to $\C^2\times M_{n-2}^{univ}$ where 
$M_{n-2}^{univ}$ is a bounded domain in $\C^{n-2}$ and a 
Teichm\"uller space. 

We will need explicitly only the cases $C_2^{univ}\cong\C^2$
and $C_3^{univ}\cong \C^2\times\H$. They are treated together
with the actions of $\Br_2$ respectively $\Br_3$ as groups
of deck transformations in Theorem \ref{t8.12} ($n=2$) and
in Theorem \ref{t9.1} and Remark \ref{t9.2} ($n=3$). 

(ix) As $C_n^{conf}$ is a connected complex manifold,
each connected covering is a connected complex manifold.
As $C_n^{conf}$ is an affine algebraic manifold, each finite
connected covering is an affine algebraic manifold.
This follows from Grothendieck's version of the Riemann
existence theorem, relating especially finite etale coverings of 
an affine algebraic manifold 
\index{affine algebraic manifold}
to finite coverings of the underlying complex manifold
\cite[XII Th\'eor\`eme 5.1]{Gr61}.
\end{remarks}

\begin{definition}\label{t8.3}
Fix $n\in\Z_{\ge 2}$.
Let $X\neq\emptyset$ be a set on which $\Br_n$ acts
transitively from the left, and let $x_0\in X$. 
The covering of $C_n^{conf}$ with fundamental group the
stabilizer $(\Br_n)_{x_0}$ and with a base point 
$b_n^{x_0}\in C_n^{x_0}$ is denoted by 
$$\pr_n^{x_0,c}:(C_n^{x_0},b_n^{x_0})\to(C_n^{conf},b_n^{conf}).
$$ 
\index{$C_n^{x_0}$}\index{$b_n^{x_0}$}
It is isomorphic to the quotient manifold 
$C_n^{univ}/(\Br_n)_{x_0}$ (in the quotient, we write the
group on the right although the action on $C_n^{univ}$ is
from the left). 
\end{definition}

\begin{remarks}\label{t8.4}
Let $(H_\Z,L,\uuuu{e})$ be a unimodular bilinear lattice
of rank $n\geq 2$ 
with a triangular basis $\uuuu{e}$ and matrix
$S=L(\uuuu{e}^t,\uuuu{e})^t\in T^{uni}_n(\Z)$. 
Consider the coverings
\begin{eqnarray*}
\pr_n^{\uuuu{e},c}:(C_n^{\uuuu{e}/\{\pm 1\}^n},
b_n^{\uuuu{e}/\{\pm 1\}^n})\to (C_n^{conf},b_n^{conf}),\\
\pr_n^{S,c}:(C_n^{S/\{\pm 1\}^n},
b_n^{S/\{\pm 1\}^n})\to (C_n^{conf},b_n^{conf}).
\end{eqnarray*}
\index{$C_n^{\uuuu{e}/\{\pm 1\}^n},\ C_n^{S/\{\pm 1\}^n}$}
\index{$b_n^{\uuuu{e}/\{\pm 1\}^n},\ b_n^{S/\{\pm 1\}^n}$}
We claim that especially the covering 
$C_n^{\uuuu{e}/\{\pm 1\}^n}$ is important.

(i) They are the following quotients of the universal covering,
\begin{eqnarray*}
C_n^{\uuuu{e}/\{\pm 1\}^n}&\cong & 
C_n^{univ}/(\Br_n)_{\uuuu{e}/\{\pm 1\}^n},\quad
\textup{with }b_n^{\uuuu{e}/\{\pm 1\}^n}=[b_n^{univ}],\\
C_n^{S/\{\pm 1\}^n}&\cong & 
C_n^{univ}/(\Br_n)_{S/\{\pm 1\}^n},\quad
\textup{with }b_n^{S/\{\pm 1\}^n}=[b_n^{univ}].
\end{eqnarray*}

(ii) Because of Lemma \ref{t3.25} (e)
$(\Br_n)_{\uuuu{e}/\{\pm 1\}^n}$ is a normal subgroup of
$(\Br_n)_{S/\{\pm 1\}^n}$. Therefore there is a normal
covering $C_n^{\uuuu{e}/\{\pm 1\}^n}\to C_n^{S/\{\pm 1\}^n}$
with group of deck transformations the quotient group
$(\Br_n)_{S/\{\pm 1\}^n}/(\Br_n)_{\uuuu{e}/\{\pm 1\}^n}$, 
and the covering space $C_n^{S/\{\pm 1\}^n}$ is canonically
isomorphic to the quotient 
$$C_n^{\uuuu{e}/\{\pm 1\}^n} / \frac{(\Br_n)_{S/\{\pm 1\}^n}}
{(\Br_n)_{\uuuu{e}/\{\pm 1\}^n}}$$
(we write the group on the right although it acts from the
left). 

(iii) Recall from Lemma \ref{t3.25} (e) that the quotient
group $(\Br_n)_{S/\{\pm 1\}^n}/(\Br_n)_{\uuuu{e}/\{\pm 1\}^n}$
is isomorphic to the group $G_\Z^{\BB}/Z((\{\pm 1\}^n)_S)$.
Often $G_\Z^{\BB}=G_\Z$.

Recall from Lemma \ref{t3.25} (f) that in the case of
$(H_\Z,L,\uuuu{e})$ irreducible 
$Z((\{\pm 1\}^n)_S)=\{\pm \id\}\subset G_\Z$, so
the quotient group above is then isomorphic to
$G_\Z^{\BB}/\{\pm \id\}\subset G_\Z/\{\pm \id\}$. 
Often equality holds. 

(iv) If $I^{(0)}$ is positive definite, then
$\www{S}+\www{S}^t$ for $\www{S}\in \SSS^{dist}$ is
positive definite with $(\www{S}+\www{S}^t)_{jj}=2$, so
$\www{S}_{ij}\in\{0,\pm 1\}$ for $i<j$. 
Then the sets $R^{(0)}$, $\BB^{dist}$ and $\SSS^{dist}$ are 
finite, the covering $C_n^{\uuuu{e}/\{\pm 1\}^n}\to
C_n^{conf}$ is finite of degree $|\BB^{dist}/\{\pm 1\}^n|$,
the covering $C_n^{S/\{\pm 1\}^n}\to C_n^{conf}$ is
finite of degree $|\SSS^{dist}/\{\pm 1\}^n|$, and the
manifolds $C_n^{\uuuu{e}/\{\pm 1\}^n}$ and 
$C_n^{S/\{\pm 1\}^n}$ are affine algebraic. 

(v) If $I^{(0)}$ is positive semidefinite, then
$\www{S}+\www{S}^t$ for $\www{S}\in \SSS^{dist}$ is
positive semidefinite with $(\www{S}+\www{S}^t)_{jj}=2$, so
$\www{S}_{ij}\in\{0,\pm 1,\pm 2\}$ for $i<j$. 
Then the set $\SSS^{dist}$ is finite, 
the covering $C_n^{S/\{\pm 1\}^n}\to C_n^{conf}$ is
finite of degree $|\SSS^{dist}/\{\pm 1\}^n|$, and the
manifold $C_n^{S/\{\pm 1\}^n}$ is affine algebraic. 

(vi) Often the covering space $C_n^{\uuuu{e}/\{\pm 1\}^n}$
has a natural smooth partial compactification
\index{partial compactification}
$\oooo{C_n^{\uuuu{e}/\{\pm 1\}^n}}^p$, and the covering map
$\pr_n^{\uuuu{e},c}:C_n^{\uuuu{e}/\{\pm 1\}^n}\to 
C_n^{conf}$ extends to a holomorphic map
$\oooo{C_n^{\uuuu{e}/\{\pm 1\}^n}}^p\to \C^n/S_n$. 

(vii) If $H_\Z$ has rank $n=1$, we have $S=(1)$, and we set
\begin{eqnarray*}
C_1^{univ}=C_1^{\uuuu{e_1}/\{\pm 1\}}=C_1^{S/\{\pm 1\}}
=C_1^{pure}=C_1^{conf}=\C,\\
\Br_1=(\Br_1)_{\uuuu{e_1}/\{\pm 1\}}
=(\Br_1)_{S/\{\pm 1\}}=\{1\}.
\end{eqnarray*}
Then the statements in the parts (i)--(iv) still hold,
though in a rather trivial way. The coordinate on $\C$
is here called $u_1$. 
\end{remarks}

The covering $C_n^{x_0}\to C_n^{conf}$ in Definition \ref{t8.3}
can and will be described in Lemma \ref{t8.6} in a more
explicit way, by gluing manifolds with boundaries, one for each
element of the set $X$. The manifolds with boundaries are all
copies of $\oooo{F_n}$, which is given in Definition \ref{t8.5}.

\begin{definition}\label{t8.5}
Fix $n\in\Z_{\ge 2}$.
Define for $j\in\{1,...,n-1\}$
\begin{eqnarray*}
W_n^{j,+}&:=&\{\uuuu{u}\in C_n^{pure}\,|\, 
\Imm(u_1)\leq ...\leq \Imm(u_n),\\  
&& \hspace*{2.5cm} \Imm(u_j)=\Imm(u_{j+1}),\ 
\Ree(u_j)<\Ree(u_{j+1})\},\\
W_n^{j,-}&:=&\{\uuuu{u}\in C_n^{pure}\,|\, 
\Imm(u_1)\leq ...\leq \Imm(u_n),\\  
&& \hspace*{2.5cm} \Imm(u_j)=\Imm(u_{j+1}),\ 
\Ree(u_j)>\Ree(u_{j+1})\},\\
W_n&:=& \bigcup_{j=1}^{n-1}\bigl(W_n^{j,+}\cup W_n^{j,-}\bigr),
\\
F_n&:=& \{\uuuu{u}\in C_n^{pure}\,|\, 
\Imm(u_1)<...<\Imm(u_n)\}.
\end{eqnarray*}
\index{$W_n^+,\ W_n^-,\ W_n$}\index{$F_n$} 
\end{definition}

\begin{lemma}\label{t8.6}
Fix $n\in\Z_{\geq 2}$. 

(a) The set $W_n^{j,\varepsilon}\subset C_n^{pure}$ 
for $\varepsilon\in\{\pm 1\}$ is called {\sf wall} or
\index{wall}
{\sf $\sigma_j^\varepsilon$-wall}. 
It is closed in $C_n^{pure}$, 
it has real codimension 1 in $C_n^{pure}$, 
it has boundary in $C_n^{pure}$ if $n\geq 3$,
it is (outside its boundary if $n\geq 3$)
smooth and affine linear, it is contractible. 
If $n\geq 4$ its boundary in $C_n^{pure}$ is only
piecewise smooth and affine linear.

The set $F_n$ is an open convex polyhedron in 
$\C^n\supset C_n^{pure}$ and therefore contractible.
Its boundary consists of the union $W_n$ of the 
$2(n-1)$ walls. Two walls $W_n^{j_1,\varepsilon_1}$ and
$W_n^{j_2,\varepsilon_2}$ with $j_1\neq j_2$ intersect
in real codimension 1. Two walls $W_n^{j,+}$ and $W_n^{j,-}$
do not intersect.

(b) The covering $\pr_n^{p,c}:C_n^{pure}\to C_n^{conf}$ 
restricts to a homeomorphism from $F_n$ to its image in
$C_n^{conf}$. The complement of $\pr_n^{p,c}(F_n)$ in
$C_n^{conf}$ is $\pr_n^{p,c}(W_n)$. It has real codimension 1
and is the boundary of $\pr_n^{p,c}(F_n)$. 
The covering $\pr_n^{p,c}$ maps $W_n^{j,+}$ and $W_n^{j,-}$
to the same set. The covering $\pr_n^{p,c}$ is not injective
on $W_n^{j,\varepsilon}$ because it maps the disjoint sets
$W_n^{j,\varepsilon}\cap W_n^{j_1,+}$ and 
$W_n^{j,\varepsilon}\cap W_n^{j_1,-}$ for $j_1\neq j$ to
the same set.

(c) Let $X\neq \emptyset$ be a set on which $\Br_n$ acts
transitively from the left, and let $x_0\in X$. 
The following are basic general facts from covering theory.
The group $\Br_n$ acts from the left on the fiber
$(\pr_n^{x_0,c})^{-1}(b_n^{conf})$ in the following way. 
For $\alpha\in Br_n=\pi_1(C_n^{conf},b_n^{conf})$ 
and $b\in (\pr_n^{x_0,c})^{-1}(b_n^{conf})$ let $\www{\alpha}$
be the lift to $C_n^{x_0}$ with end point $\www{\alpha}(1)=b$
of a loop in $C_n^{conf}$ which represents $\alpha$.
Then $\alpha(b):=\www{\alpha}(0)$. 

Write $U:=(\Br_n)_{x_0}$. This action induces the two bijections
\begin{eqnarray*}
\begin{array}{ccccc}
X&\stackrel{1:1}{\longleftarrow} &\Br_n/U & 
\stackrel{1:1}{\longrightarrow}  & 
(\pr_n^{x_0,c})^{-1}(b_n^{conf})\\ 
x=\alpha(x_0) & \longmapsfrom & \alpha U & 
\longmapsto & \alpha(b_n^{x_0})
\end{array}
\end{eqnarray*}
For $x=\alpha(x_0)$ write $b_n^{x_0,x}:=\alpha(b_n^{x_0})
\in (\pr_n^{x_0,c})^{-1}(b_n^{conf})$
(especially $b_n^{x_0}=b_n^{x_0,x_0}$). 
So there is a natural bijection between $X$ and the fiber
in $C_n^{x_0}$ over $b_n^{conf}$, and $x$ corresponds to
$b_n^{x_0,x}$. 

(d) Let $X$ and $x_0$ be as in part (c). 
The covering $\pr_n^{x_0,c}:(C_n^{x_0},b_n^{x_0})\to 
(C_n^{conf},b_n^{conf})$ restricts over the open subset
$\pr_n^{p,c}(F_n)$ of $C_n^{conf}$ to an even covering, that
means that $(\pr_n^{x_0,c})^{-1}(\pr_n^{p,c}(F_n))$ is a union
of disjoint open subsets of $C_n^{x_0}$ each of which is 
mapped homeomorphically to $\pr_n^{p,c}(F_n)$. 

The components of $(\pr_n^{x_0,c})^{-1}(\pr_n^{p,c}(F_n))$
are called {\sf Stokes regions}. \index{Stokes region}
Each contains one element of the fiber 
$(\pr_n^{x_0,c})^{-1}(b_n^{conf})$. For $x\in X$
let $F_n^{x_0,x}$ be the component which contains the
point $b_n^{x_0,x}$. 

The complement 
$C_n^{x_0}-(\pr_n^{x_0,c})^{-1}(\pr_n^{p,c}(F_n))$ is the 
preimage of the walls $\pr_n^{p,c}(W_n)$ in $C_n^{conf}$,
it has everywhere real codimension 1. It is the boundary 
of $(\pr_n^{x_0,c})^{-1}(\pr_n^{p,c}(F_n))$.

(e) In $C_n^{conf}$ the images of the walls $W_n^{j,+}$ and
$W_n^{j,-}$ coincide, but this image is a smooth oriented
real hypersurface and has two sides, the
$W_n^{j,+}$-side and the $W_n^{j,-}$-side. 
The braid $\sigma_j$ has a representative by a loop which goes
once through this hypersurface, from the $W_n^{j,+}$-side
to the $W_n^{j,-}$-side. 

(f) Let $X$ and $x_0$ be as in the parts (c) and (d). 
The covering space $C_n^{x_0}$ can be obtained by glueing
copies of the closure $\oooo{F_n}$ in $C_n^{pure}$ in the
following way along their boundaries. 
More precisely, it is isomorphic to the quotient of the
disjoint union $\bigcup_{x\in X}(\oooo{F_n}\times \{x\})$
by the equivalence relation which is generated by the following
equivalences.
\begin{eqnarray*}
(\uuuu{u},x)&\sim& ((u_1,...,u_{j-1},u_{j+1},u_j,u_{j+2},...,
u_n)),\sigma_j^{-1}(x)) \quad \textup{for }x\in X\\
&& \textup{ and }\uuuu{u}\in W_n^{j,+}
\textup{ (and thus }(u_1,...,u_{j+1},u_j,...,u_n)\in W_n^{j,-}.
\end{eqnarray*}
So the walls $W_n^{j,+}\times\{x\}$ and 
$W_n^{j,-}\times \{\sigma_j^{-1}(x)\}$ are identified. 
The isomorphism maps $F_n\times\{x\}$ to $F_n^{x_0,x}$. 
\end{lemma}

{\bf Proof:}
(a) Trivial.

(b) Trivial.

(c) The bijection $\Br_n/U\to X$, $\alpha U\to \alpha(x_0)$,
is elementary group theory. The group action from the left
of $\Br_n$ on the fiber $(\pr_n^{x_0,c})^{-1}(b_n^{conf})$
is described for example in \cite[page 65, section  1.3]{Ha01}.
Here the image $\alpha(b)$ of the {\it end} point 
$b=\www{\alpha}(1)$ under $\alpha$ is the {\it starting} point 
$\www{\alpha}(0)$, because composition of paths is written
from left to right, whereas the action of 
$\Br_n$ on the fiber $(\pr_n^{x_0,c})^{-1}(b_n^{conf})$ 
shall be an action from the left. 

The stabilizer of the base point $b_n^{x_0}$ is by construction
of the covering the group $U=(\Br_n)_{x_0}$.
This gives the second bijection 
$\Br_n/U\to (\pr_n^{x_0,c})^{-1}(b_n^{conf})$,
$\alpha U\mapsto \alpha(b_n^{x_0})$. 

(d) This follows from the shape of the subset 
$\pr_n^{p,c}(F_n)$ in $C_n^{conf}$, from the definition of
a covering of $C_n^{conf}$ and from part (c).

\begin{figure}
\includegraphics[width=1.0\textwidth]{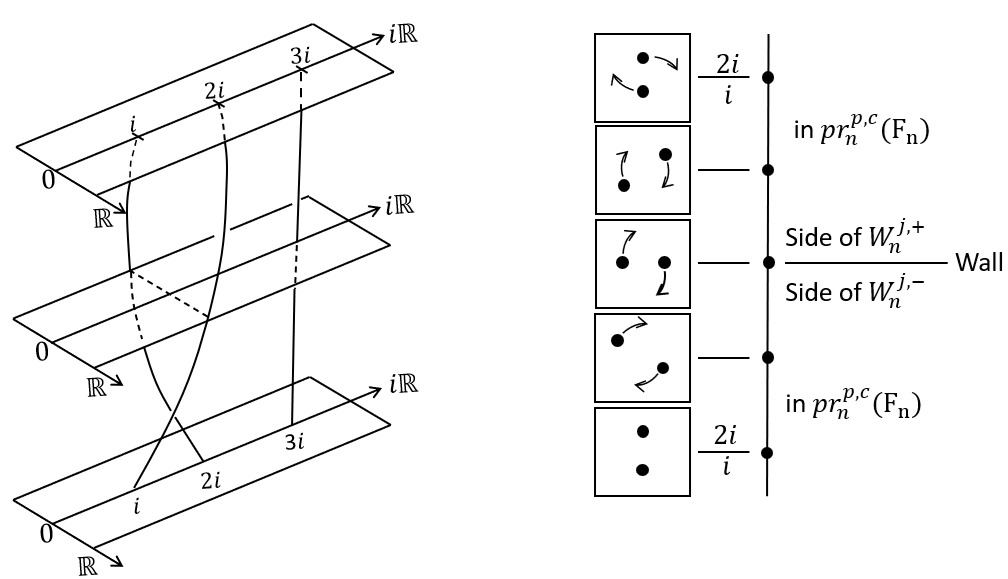}
\caption[Figure 8.2]{Crossing a wall with an elementary braid}
\label{Fig:8.2}
\end{figure}

(e) Compare the description of the elementary braid $\sigma_j$
in Remark \ref{t8.2} (iii) and the two pictures in Figure
\ref{Fig:8.2}.

The path $[0,1]\to C_n^{pure}$, $t\mapsto \uuuu{u}(t)$, 
which represents $\sigma_j$ is chosen such that
$u_1(t),...,u_{j-1}(t),u_{j+2}(t),...,u_n(t)$ are constant
and $u_j(t)$ and $u_{j+1}(t)$ have a constant center
$\frac{1}{2}(u_j(t)+u_{j+1}(t))$ and turn around halfway
clockwise.

(f) This follows from the parts (c), (d) and (e).
\hfill$\Box$

\begin{remarks}\label{t8.7}
Let $X$ and $x_0$ be as in the parts (c), (d) and (f) in Lemma
\ref{t8.6}. 

(i) The union
\begin{eqnarray*}
W_n^{sing}:= \bigcup_{j_1,j_2\in\{1,...,n-1\},j_1\neq j_2}
\bigcup_{\varepsilon_1,\varepsilon_2\in\{\pm 1\}}
W_n^{j_1,\varepsilon_1}\cap W_n^{j_2,\varepsilon_2}
\subset W_n\subset C_n^{pure}
\end{eqnarray*}
has real codimension 2 in $C_n^{pure}$. Also the image
$\pr_n^{p,c}(W_n^{sing})\subset C_n^{conf}$ has real
codimension 2 in $C_n^{conf}$.

(ii) Choose a sequence 
$(\sigma_{j_a}^{\varepsilon_a})_{a\in\{1,...,l\}}$ 
of elementary braids for some $l\in\N$, $j_a\in\{1,...,n-1\}$,
$\varepsilon_a\in\{\pm 1\}$. Choose a closed path $p$ in
$C_n^{conf}$ with starting point and end point $b_n^{conf}$
which avoids $\pr_n^{p,c}(W_n^{sing})$ and which 
passes $l$ times through $\pr_n^{p,c}(W_n)$,
first from the $W_n^{j_1,\varepsilon_1}$-side through 
$\pr_n^{p,c}(W_n^{j_1,\varepsilon_1})$,
then from the $W_n^{j_2,\varepsilon_2}$-side through
$\pr_n^{p,c}(W_n^{j_2,\varepsilon_2})$, and so on.
The path $p$ represents the braid
$\alpha=\sigma_{j_1}^{\varepsilon_1}\sigma_{j_2}^{\varepsilon_2}
...\sigma_{j_l}^{\varepsilon_l}$. 
The lift of $p$ to $C_n^{x_0}$ which starts at $b_n^{x_0,x_0}$
starts in the Stokes region $F_n^{x_0,x_0}$ and goes through
the wall $W_n^{j_1,\varepsilon_1}$ of it to the Stokes region
$F_n^{x_0,\sigma_{j_1}^{-\varepsilon_1}(x_0)}$,
then through the wall $W_n^{j_2,\varepsilon_2}$ of this
Stokes region to the Stokes region
$F_n^{x_0,\sigma_{j_2}^{-\varepsilon_2}
\sigma_{j_1}^{-\varepsilon_1}(x_0)}$, and so on. It ends in the
Stokes region $F_n^{x_0,\alpha^{-1}(x_0)}$ at the
point $b_n^{x_0,\alpha^{-1}(x_0)}$. 

(iii) A homotopy from a path as in part (ii)
for $\sigma_j\sigma_{j+1}\sigma_j$ to a path as in part (ii)
for $\sigma_{j+1}\sigma_j\sigma_{j+1}$ contains at least
one path which contains a point in 
$\pr_n^{p,c}(W_n^{j,+}\cap W_n^{j+1,+})$. A preimage in
$C_n^{pure}$ of this point satisfies
$\Imm(u_j)=\Imm(u_{j+1})=\Imm(u_{j+2})$. 

(iv) If the covering $\pr_n^{x_0,c}:C_n^{x_0}\to C_n^{conf}$
is normal, the deck transformation group of it acts
transitively on the set of Stokes regions 
$\{F_n^{x_0,x}\,|\, x\in X\}$. Then $F_n^{x_0,x_0}$ is
a fundamental domain in $C_n^{x_0}$ for the action of
the deck transformation group. 
\end{remarks}

\section{Semisimple F-manifolds}
\label{s8.2}

Each covering space of $C_n^{conf}$ is a semisimple
F-manifold with Euler field and with empty Maxwell stratum.
We explain these notions here.

\begin{definition}\label{t8.8}
\cite{HM98}\cite[Definition 2.8]{He02}

(a) A (holomorphic) 
\index{F-manifold}
{\it F-manifold} is a tuple $({\MM},\circ,e)$
\index{$(M,\circ,e,E)$} 
where ${\MM}$ is a complex manifold, $\circ$ is a commutative
and associative multiplication on the holomorphic tangent bundle
\index{multiplication on the holomorphic tangent bundle}
$T{\MM}$ of ${\MM}$ which satisfies the integrability condition
\begin{eqnarray*}
\Lie_{X\circ Y}&=& Y\circ \Lie_X(\circ)+ X\circ\Lie_Y(\circ)
\end{eqnarray*}
for local holomorphic vector fields $X$ and $Y$, 
and $e$ is a global holomorphic vector field with
$e\,\circ=\id$. It is called {\it unit field}.
\index{unit field}

(b) An {\it Euler field} \index{Euler field}
$E$ in an F-manifold $({\MM},\circ,e)$
is a global holomorphic vector field with $\Lie_E(\circ)=\circ$.

(c) An F-manifold is semisimple 
\index{semisimple F-manifold}
at a point $t\in {\MM}$ if the 
algebra $(T_t{\MM},\circ_t,e_t)$ is semisimple, that means,
it splits canonically into a sum of 1-dimensional
algebras, so $T_t{\MM}=\bigoplus_{j=1}^n\C\cdot e_{t,j}$ with
$e_{t,i}\circ e_{t,j}=\delta_{ij}e_{t,i}$
(and automatically $e=\sum_{j=1}^n e_j$). 
\index{semisimple multiplication}
\end{definition}

\begin{remarks}\label{t8.9}
(i) The theory of F-manifolds is developed in 
\cite[chapters 1--5]{He02}.

(ii) The condition to be semisimple is an open condition.
If an F-manifold is semisimple at one point, it is generically
semisimple. More precisely, then the set $\KK_3\subset {\MM}$
of points where it is not semisimple is either empty
or an analytic hypersurface \cite[Proposition 2.6]{He02}. 
The set $\KK_3$ is called {\it caustic}. 
\index{caustic}\index{$\KK_3$: caustic}

(iii) A complex manifold with a semisimple (commutative and
associative) multiplication on the holomorphic tangent bundle
has locally a basis $e_1,...,e_n$ of holomorphic vector fields
which are {\it the partial units}, so 
$e_i\circ e_j=\delta_{ij}e_i$. It is unique up to indexing.
Here the integrability condition of an F-manifold says
$[e_i,e_j]=0$, so the partial units are locally coordinate
vector fields \cite[Theorem 3.2]{He02}.
The corresponding coordinates are unique up to addition of
constants.

The choice of an Euler field $E$ makes them unique, because
the eigenvalues $u_1,...,u_n$ of the endomorphism 
$E\,\circ$ on the holomorphic tangent bundle are such local
coordinates with $e_1=\frac{\paa}{\paa u_1},...,
e_n=\frac{\paa}{\paa u_n}$ and $E=\sum_{j=1}^n u_je_j$.
These local coordinates are unique up to indexing in the
case of a semisimple F-manifold with Euler field and are
called {\it canonical coordinates}. \index{canonical coordinates}

(iv) Let $({\MM},\circ,e,E)$ be a generically semisimple F-manifold
with Euler field. Besides the caustic $\KK_3$ there is often
a second analytic hypersurface $\KK_2$, the {\it Maxwell 
stratum}. \index{Maxwell stratum}\index{$\KK_2$: Maxwell stratum} 
It is the closure of the set of points $t\in {\MM}$
such that $(T_t{\MM},\circ_t,e_t)$ is semisimple
(so locally canonical coordinates exist), but at least two
eigenvalues of $E\, \circ$ coincide.
Usually $\KK_3$ and $\KK_2$ intersect. 
The map
\begin{eqnarray*}
\LL:{\MM}\to \C[x]_n,\quad t\mapsto
\textup{ (the characteristic polynomial of }E_t\,\circ),
\end{eqnarray*}
is holomorphic. In singularity theory it is called
{\it Lyashko-Looijenga map}. 
\index{Lyashko-Looijenga map}\index{$\LL$}
$\KK_3\cup\KK_2$ is mapped to the 
discriminant $\C[x]_n^{mult}$,
the complement ${\MM}-(\KK_3\cup\KK_2)$ is mapped locally
biholomorphically to the complement $\C[x]_n^{reg}\cong
C_n^{conf}$. Near a generic point of $\KK_2$
there are local coordinates such that two of them take the
same value at points of $\KK_2$. Therefore $\LL$ is
a branched covering of degree two near a generic point
of $\KK_2$. The index 3 in $\KK_3$ stems from the fact that
in a family of important cases $\LL$ is a branched covering
of order 3 near generic points of $\KK_3$.

(v) Vice versa, a locally biholomorphic map
$\LL:{\MM}\to\C[x]_n^{reg}\cong C_n^{conf}$ induces on ${\MM}$ the structure
of a semisimple F-manifold with Euler field and empty
Maxwell stratum. 

(vi) Define \index{$\C[x]_n^{mult,reg}$}
\begin{eqnarray*}
\C[x]_n^{mult,reg}:=\{f\in \C[x]^{mult}\,|\, 
f\textup{ has precisely }n-1\textup{ different roots}\}.
\end{eqnarray*}
It is an open part of $\C[x]_n^{mult}$. The complement is
a complex hypersurface in $\C[x]_n^{mult}$ (which itself is a
complex hypersurface in $\C[x]_n$).

Consider the following situation:
A locally biholomorphic map $\LL:{\MM}\to \C[x]_n^{reg}$,
a point $p\in \C[x]_n^{mult,reg}$, a neighborhood
$U\subset \C[x]_n$ of $p$ with 
$U\cap \C[x]_n^{mult}\subset\C[x]_n^{mult,reg}$,
$U':=U-\C[x]_n^{mult}$, a component $\www{U'}$ of
$\LL^{-1}(U')$ such that the restriction
$\LL:\www{U'}\to U'$ is an $m$-fold covering  for some
$m\in\N$. Then $\www{U'}$ 
has a partial compactification $\www{U}$
\index{partial compactification} 
such that $\LL$ extends to a branched covering 
$\www{\LL}:\www{U}\to U$, which is branched of order $m$ 
over $U-U'=U\cap \C[x]_n^{mult}$. 
If $m\geq 2$ the F-manifold structure extends from
$\www{U'}$ to $\www{U}$.
If $m=2$ $\www{U}-\www{U'}$ is the Maxwell stratum in $\www{U}$.
If $m\geq 3$ $\www{U}-\www{U'}$ is the caustic in $\www{U}$.
See the Examples \ref{t8.10}.

(vii) A semisimple F-manifold with Euler field is locally
a rather trivial object. More interesting is a generically
semisimple F-manifold near points of the caustic $\KK_3$.

(viii) The notion of a {\it Dubrovin-Frobenius manifold}
\index{Dubrovin-Frobenius manifold}\index{Frobenius manifold}
is much richer. It is an F-manifold $({\MM},\circ,e,E)$ with
Euler field together with a flat holomorphic metric $g$
on the holomorphic tangent bundle such that $e$ is flat
(with respect to the Levi-Civita connection of $g$),
$g$ is {\it multiplication invariant}, that means
\begin{eqnarray*}
g(X\circ Y,Z)=g(X,Y\circ Z)
\end{eqnarray*}
for local holomorphic vector fields $X,Y,Z$, 
and $$\Lie_E(g)=(2-d)g \quad\textup{for some}\quad d\in\C.$$ 
It was first defined by Dubrovin \cite{Du92}\cite{Du96}.
To see the equivalence of the definition above with the
definition in \cite{Du92}\cite{Du96} one needs 
\cite[Theorem 2.15 and Lemma 2.16]{He02}.

Dubrovin-Frobenius manifolds arise in singularity theory,
in quantum cohomology, in mirror symmetry and in integrable
systems. They are all highly transcendental and interesting
objects. The construction in singularity theory of a metric
$g$ which makes a semisimple F-manifold with Euler field
into a Dubrovin-Frobenius manifold uses additional 
structure above the F-manifold. One version of this
additional structure is presented in section \ref{s8.5}.
\end{remarks}

\begin{examples}\label{t8.10}
(i) Of course, $\C^n$ with coordinates
$u_1,...,u_n$ is a semisimple F-manifold with
partial units $e_i=\frac{\paa}{\paa u_i}$ and Euler
field $E=\sum_{i=1}^n u_ie_i$. It is called of type $A_1^n$.
The Maxwell stratum is the union $D_n^{pure}$ of partial
diagonals.

(ii) For $m\in\Z_{\geq 2}$ the manifold ${\MM}=\C^2$ with
coordinates $z_1,z_2$ becomes a generically semisimple
F-manifold with Euler field in the following way
\cite[Theorem 4.7]{He02}. 
It is called of type $I_2(m)$. \index{$I_2(m)$} 
Denote by $\paa_1:=\frac{\paa}{\paa z_1}$
and $\paa_2:=\frac{\paa}{\paa z_2}$ the coordinate
vector fields and define the unit field, 
the multiplication $\circ$ 
and the Euler field $E$ as follows,
\begin{eqnarray*}
e&:=& \paa_1,\quad\textup{ so }\paa_1\,\circ =\id
\textup{ on the holomorphic tangent bundle},\\
\paa_2\circ\paa_2&:=& z_2^{m-2}\paa_1,\\
E&:=& z_1\paa_1+\frac{2}{m}z_2\paa_2.
\end{eqnarray*}
If $m=2$ it is everywhere semisimple, and the Maxwell
stratum is $\KK_2=\{(z_1,z_2)\in {\MM}\,|\, z_2=0\}$. 
If $m\geq 3$ it is generically semisimple, and the caustic is
$\KK_3=\{(z_1,z_2)\in {\MM}\,|\, z_2=0\}$. 

For $m=2$ there are global canonical coordinates, 
for even $m\geq 3$ there are global canonical coordinates 
on ${\MM}-\KK_3$, for odd $m\geq 3$ there are local canonical
coordinates on ${\MM}-\KK_3$, which in fact are globally 
2-valued. In all cases they are
\begin{eqnarray*}
u_{1/2}=z_1\pm \frac{2}{m}z_2^{m/2}.
\end{eqnarray*}
The partial units are
\begin{eqnarray*}
e_{1/2}=\frac{1}{2}e\pm \frac{1}{2} z_2^{1-\frac{m}{2}}
\partial_2.
\end{eqnarray*}
The Lyashko-Looijenga map
\begin{eqnarray*}
\LL:{\MM}\to \C[x]_2,&& (z_1,z_2)\mapsto (-u_1-u_2,u_1u_2)=
(-2z_1,z_1^2-\frac{4}{m^2}z_2^m)\\
&& \cong x^2+(-u_1-u_2)x+u_1u_2=(x-u_1)(x-u_2),
\end{eqnarray*}
restricts to an $m$-fold covering
$\LL:{\MM}-\{z_2=0\}\to \C[x]^{reg}$, 
it maps $\{z_2=0\}$ to $\C[x]^{mult}$, and it is branched
of order $m$ along $\{z_2=0\}$. 

(iii) The type $I_2(2)$ is the type $A_1^2$, the type
$I_2(3)$ is also called type $A_2$. 

(iv) Theorem 4.7 in \cite{He02} says also that each 
2-dimensional irreducible germ of a generically semisimple 
F-manifold with Euler field is isomorphic to the germ
$(({\MM},t^0),\circ,e,E)$ in the F-manifold ${\MM}=\C^2$ of type
$I_2(m)$ for some $t^0=(0,z_2^0)\in\KK_3$
and some $m\in\Z_{\geq 3}$. 

(v) If $({\MM}^{(i)},\circ^{(i)},e^{(i)},E^{(i)})$ for 
$i\in\{1,2\}$ are F-manifolds with Euler fields, their product
becomes an F-manifold
$({\MM}^{(1)}\times {\MM}^{(2)},\circ^{(1)}\times \circ^{(2)},
e^{(1)}+e^{(2)},E^{(1)}+E^{(2)})$ in the natural way,
with Maxwell stratum 
$\KK_2^{(1)}\times {\MM}^{(2)}\, \cup\, {\MM}^{(1)}\times \KK_2^{(2)}$
and caustik 
$\KK_3^{(1)}\times {\MM}^{(2)}\, \cup\, {\MM}^{(1)}\times \KK_3^{(2)}$.
For example, one obtains for $m\geq 2$ and $n\geq 3$ the
generically semisimple F-manifold of type 
$I_2(m)A_1^{n-2}$ with Euler field on $\C^2\times\C^{n-2}$.
\end{examples}

\section{Reducible cases and rank 2 cases}
\label{s8.3}

Let $(H_\Z,L,\uuuu{e})$ be a unimodular bilinear lattice
of rank $n\geq 2$ with a triangular basis $\uuuu{e}$
and matrix $S=L(\uuuu{e}^t,\uuuu{e})^t\in T^{uni}_n(\Z)$.
First we state a result on $C_n^{\uuuu{e}/\{\pm 1\}^n}$ in the
case when $(H_\Z,L,\uuuu{e})$ is reducible.
Then we treat the rank 2 cases.

\begin{theorem}\label{t8.11}
Let $(H_\Z,L,\uuuu{e})= (H_{\Z,1},L_1,\uuuu{e}^{(1)})
\oplus (H_{\Z,2},L_2,\uuuu{e}^{(2)})$ be reducible
with $\uuuu{e}^{(1)}=(e_1,e_2,...,e_{n_1})$ and
$\uuuu{e}^{(2)}=(e_{n_1+1},...,e_n)$, so 
$S_{ij}=0$ for $i\leq n_1<j$. Write $n_2:=n-n_1$. 

There is a natural embedding
\begin{eqnarray*}
C_n^{\uuuu{e}/\{\pm 1\}^n}\hookrightarrow
C_{n_1}^{\uuuu{e}^{(1)}/\{\pm 1\}^{n_1}}\times
C_{n_2}^{\uuuu{e}^{(2)}/\{\pm 1\}^{n_2}}.
\end{eqnarray*}
It is an embedding of semisimple F-manifolds with Euler fields.
The complement is a complex hypersurface.
This hypersurface is the Maxwell stratum of the semisimple
F-manifold with Euler field on the right hand side.
\end{theorem}

{\bf Proof:}
In order to write down the embedding concretely, we will make
use of the construction of the three manifolds by the glueing
of Stokes regions in Lemma \ref{t8.6} (f). Here
\begin{eqnarray*}
x_0&:=&\uuuu{e}/\{\pm 1\}^n\in 
X:=\Br_n(x_0)=\BB^{dist}/\{\pm 1\}^n,\\
x_0^{(1)}&:=& \uuuu{e}^{(1)}/\{\pm 1\}^{n_1}\in 
X^{(1)}:=\Br_{n_1}(x_0^{(1)})=\BB^{dist,(1)}/\{\pm 1\}^{n_1},\\
x_0^{(2)}&:=& \uuuu{e}^{(2)}/\{\pm 1\}^{n_2}\in 
X^{(2)}:=\Br_{n_2}(x_0^{(2)})=\BB^{dist,(2)}/\{\pm 1\}^{n_2},
\end{eqnarray*}
\begin{eqnarray*}
C_n^{x_0}&=& \Bigl(\bigcup_{x\in X}\oooo{F_n}
\times \{x\}\Bigr)_\sim,\\
C_{n_1}^{x_0^{(1)}}&=& \Bigl(\bigcup_{x^{(1)}\in X^{(1)}}
\oooo{F_{n_1}}\times \{x^{(1)}\}\Bigr)_{\sim^{(1)}},\\
C_{n_2}^{x_0^{(2)}}&=& \Bigl(\bigcup_{x^{(2)}\in X^{(2)}}
\oooo{F_{n_2}}\times \{x^{(2)}\}\Bigr)_{\sim^{(2)}}.
\end{eqnarray*}

For any $x=(x_1,...,x_n)\in X$ there are unique indices
$i_1,...,i_{n_1}$ and $j_1,...,j_{n_2}$ with
\begin{eqnarray*}
1\leq i_1<...<i_{n_1}\leq n,\quad  
1\leq j_1<...<j_{n_2}\leq n,\\
\{i_1,...,i_{n_1},j_1,...,j_{n_2}\}
=\{1,2,...,n\},\\
x^{(1)}:=(x_{i_1},x_{i_2},...,x_{i_{n_1}})\in X^{(1)},\
x^{(2)}:=(x_{j_1},x_{j_2},...,x_{j_{n_1}})\in X^{(2)}.
\end{eqnarray*}
In words, the distinguished basis $x$ up to signs for
$(H_\Z,L,\uuuu{e})$ is a shuffling of the distinguished basis
$x^{(1)}$ up to signs for $(H_{\Z,1},L_1,\uuuu{e}^{(1)})$
and the distinguished basis $x^{(2)}$ up to signs for 
$(H_{\Z,2},L_2,\uuuu{e}^{(2)})$. 
Here $x^{(1)}$ and $x^{(2)}$ are unique.

A point $(\uuuu{u},x)\in \oooo{F_n}\times\{x\}$ is mapped
to the point 
\begin{eqnarray*}
(\uuuu{u}^{(1)},x^{(1)};\uuuu{u}^{(2)},x^{(2)})
&:=& ((u_{i_1},u_{i_2},...,u_{i_{n_1}}),x^{(1)};
(u_{j_1},u_{j_2},...,u_{j_{n_2}}),x^{(2)})\\
&\in& (\oooo{F_{n_1}}\times\{x^{(1)}\})\times
(\oooo{F_{n_2}}\times \{x^{(2)}\}).
\end{eqnarray*}

The following observations (1), (2) and (3) are crucial.

(1) Given $x^{(1)}\in X^{(1)}$ and $x^{(2)}\in X^{(2)}$, there
are $\begin{pmatrix}n\\n_1\end{pmatrix}$ distinguished
bases $y\in X$ up to signs with $y^{(1)}=x^{(1)}$ and
$y^{(2)}=x^{(2)}$. They are all the tuples
obtained by shuffling of $x^{(1)}$ and $x^{(2)}$.

(2) For $x^{(1)}\in X^{(1)}$ and $x^{(2)}\in X^{(2)}$ the map 
\begin{eqnarray*}
\dot\bigcup_{y\in X:\, y^{(1)}=x^{(1)},y^{(2)}=x^{(2)}}
\oooo{F_n}\times\{y\}\to 
(\oooo{F_{n_1}}\times\{x^{(1)}\})\times
(\oooo{F_{n_2}}\times\{x^{(2)}\})
\end{eqnarray*}
is almost surjective and almost injective.
More precisely, the following two points (i) and (ii) hold.

(i) The map 
\begin{eqnarray*}
\left(\dot\bigcup_{y\in X:\, y^{(1)}=x^{(1)},y^{(2)}=x^{(2)}}
\oooo{F_n}\times\{y\}\right)_\sim\to 
(\oooo{F_{n_1}}\times\{x^{(1)}\})\times
(\oooo{F_{n_2}}\times\{x^{(2)}\})
\end{eqnarray*}
is well defined and injective because of the following.
The restriction to 
$\dot\bigcup_{y\in X:\, y^{(1)}=x^{(1)},y^{(2)}=x^{(2)}}
F_n\times\{y\}$ is injective, anyway. 
Consider $\uuuu{u}\in W_n^{j,+}$ and $y\in X$ with 
$x^{(1)}=y^{(1)}=(\sigma_j^{-1}y)^{(1)}$,
$x^{(2)}=y^{(2)}=(\sigma_j^{-1}y)^{(2)}$.
Then $\Imm(u_j)=\Imm(u_{j+1})$, and 
$(\uuuu{u},x)$ and 
$((u_1,...,u_{j+1},u_j,...,u_n),\sigma_j^{-1}(x))$ 
have the same image.

(ii) Consider any two points $\uuuu{u}^{(1)}\in \oooo{F_{n_1}}$
and $\uuuu{u}^{(2)}\in \oooo{F_{n_2}}$ such that the 
entries of $\uuuu{u}^{(1)}$ are pairwise different
from the entries of $\uuuu{u}^{(2)}$.
There are indices $\www{i}_1,...,\www{i}_{n_1}$ and $\www{j}_1,...,\www{j}_{n_2}$ with
\begin{eqnarray*}
1\leq \www{i}_1<...<\www{i}_{n_1}\leq n,\quad 
1\leq \www{j}_1<...<\www{j}_{n_2}\leq n,\\
\{\www{i}_1,...,\www{i}_{n_1},\www{j}_1,...,\www{j}_{n_2}\}=\{1,...,n\},
\end{eqnarray*}
such that the tuple $\uuuu{u}\in\C^n$ which is defined by 
\begin{eqnarray*}
\uuuu{u}^{(1)}=(u_{\www{i}_1},...,u_{\www{i}_{n_1}}),\quad 
\uuuu{u}^{(2)}=(u_{\www{j}_1},...,u_{\www{j}_{n_2}}),
\end{eqnarray*}
satisfies
\begin{eqnarray*}
\Imm(u_1)\leq \Imm(u_2)\leq ...\leq \Imm(u_n).
\end{eqnarray*}
The tuple $y=(y_1,...,y_n)$ which is defined by
\begin{eqnarray*}
x^{(1)}=(y_{\www{i}_1},y_{\www{i}_2},...,y_{\www{i}_{n_1}}),\quad
x^{(2)}=(y_{\www{j}_1},y_{\www{j}_2},...,y_{\www{j}_{n_2}}),
\end{eqnarray*}
is in $X$ and is obtained by shuffling of $x^{(1)}$
and $x^{(2)}$. 
The pair $(\uuuu{u},y)$ is mapped to the tuple
$(\uuuu{u}^{(1)},x^{(1)};\uuuu{u}^{(2)},x^{(2)})$.
Therefore the only points 
$(\uuuu{u}^{(1)},x^{(1)};\uuuu{u}^{(2)},x^{(2)})
\in (\oooo{F_{n_1}}\times \{x^{(1)}\})\times 
(\oooo{F_{n_2}}\times\{x^{(2)}\})$
which are not met by the map in part (i) are those points
where some entry of $\uuuu{u}^{(1)}$ coincides 
with some entry of $\uuuu{u}^{(2)}$.
In the semisimple F-manifold 
$(F_{n_1}\times \{x^{(1)}\})\times (F_{n_2}\times\{x^{(2)}\})$
with Euler field these points form the Maxwell stratum.

(3) (i) The  map
\begin{eqnarray*}
\dot\bigcup_{x\in X}\oooo{F_n}\times\{x\}\to
\Bigl(\dot\bigcup_{x^{(1)}\in X^{(1)}}
\oooo{F_{n_1}}\times\{x^{(1)}\}\Bigr)\times
\Bigl(\dot\bigcup_{x^{(2)}\in X^{(2)}}
\oooo{F_{n_2}}\times\{x^{(2)}\}\Bigr)
\end{eqnarray*}
is well defined (and almost surjective and almost injective).
This follows from (2)(i) and (ii).

(ii) {\bf Claim:} 
The map in (3)(i) induces a well defined injective map
on the quotients
\begin{eqnarray*}
&&\Bigl(\dot\bigcup_{x\in X}\oooo{F_n}\times\{x\}\Bigr)_{\sim}
\to\\
&&\Bigl(\dot\bigcup_{x^{(1)}\in X^{(1)}}
\oooo{F_{n_1}}\times\{x^{(1)}\}\Bigr)_{\sim^{(1)}}\times
\Bigl(\dot\bigcup_{x^{(2)}\in X^{(2)}}
\oooo{F_{n_2}}\times\{x^{(2)}\}\Bigr)_{\sim^{(2)}}.
\end{eqnarray*}

Proof of the Claim:
The restriction to $\dot\bigcup_{x\in X}F_n\times\{x\}$
is well defined and injective, anyway. 
Consider $\uuuu{u}\in W_n^{j,+}$ and $x\in X$.
If $(\sigma_j^{-1}(x))^{(1)}=x^{(1)}$ and 
$(\sigma_j^{-1}(x))^{(2)}=x^{(2)}$, we are in the case
treated in (2)(ii). If not, then either
$x_j$ and $x_{j+1}$ are both in $H_{\Z,1}/\{\pm 1\}^{n_1}$ 
or $x_j$ and $x_{j+1}$ are both in $H_{\Z,2}/\{\pm 1\}^{n_2}$. 
In the first case 
$(\uuuu{u}^{(1)},x^{(1)})\sim^{(1)}
((u_1,...,u_{j+1},u_j,...,u_n)^{(1)},
(\sigma_j^{-1}(x))^{(1)})$,
in the second case
$(\uuuu{u}^{(2)},x^{(2)})\sim^{(2)}
((u_1,...,u_{j+1},u_j,...,u_n)^{(2)},
(\sigma_j^{-1}(x))^{(2)})$.
This proves the Claim in (3)(ii).

The observations (1)--(3) show that there is a natural open
embedding 
\begin{eqnarray*}
C_n^{x_0}\hookrightarrow C_{n_1}^{x_0^{(1)}} 
\times C_{n_2}^{x_0^{(2)}}
\end{eqnarray*}
and that the complement is a hypersurface, which is the
Maxwell stratum of 
$C_{n_1}^{x_0^{(1)}} \times C_{n_2}^{x_0^{(2)}}$
as a semisimple F-manifold with Euler field.
\hfill$\Box$ 

\bigskip

Now we come to the rank 2 cases.
The parts (a)--(c) of the following theorem treat
$C_2^{pure}$, $C_2^{conf}$ and $C_2^{univ}$ as
semisimple F-manifolds with Euler fields,
part (d) considers $C_2^{\uuuu{e}/\{\pm 1\}^2}$ for
a unimodular bilinear lattice of rank 2 with a triangular basis.

\begin{theorem}\label{t8.12}
(a) $C_2^{pure}$ is isomorphic to the semisimple F-manifold
with Euler field $I_2(2)$ in the Examples \ref{t8.10},
minus its Maxwell stratum. In formulas:
\begin{eqnarray*}
C_2^{pure}=\{(u_1,u_2)\in\C^2\,|\, u_1\neq u_2\}
&\stackrel{\cong}{\longrightarrow}&
\C\times \C^*,\\
(u_1,u_2)&\mapsto& (z_1,z_2)
=(\frac{u_1+u_2}{2},\frac{u_1-u_2}{2}).
\end{eqnarray*}
Therefore $C_2^{pure}$ extends to an F-manifold on $\C^2$
with Maxwell stratum $\{(u_1,u_2)\in \C^2\,|\, u_1=u_2\}$.
Multiplication $\circ$, unit field $e$, Euler field $E$
and partial units $e_1$ and $e_2$ are explicitly as follows, 
\begin{eqnarray*}
u_{1/2}=z_1\pm z_2,\quad e_i=\frac{\paa}{\paa u_i},
\quad e_i\circ e_j=\delta_{ij}e_i,\\
e=e_1+e_2=\frac{\paa}{\paa z_1},
\quad\textup{ so }e\,\circ =\id,\\
\frac{\paa}{\paa z_2}=e_1-e_2,\quad 
\frac{\paa}{\paa z_2}\circ \frac{\paa}{\paa z_2}=e,\\
E=u_1e_1+u_2e_2=z_1\frac{\paa}{\paa z_1}
+z_2\frac{\paa}{\paa z_2}.
\end{eqnarray*}

(b) The covering $C_2^{pure}\to C_2^{conf}$ has degree 2.
In the coordinates $(z_1,z_2)$ on $C_2^{pure}$ it extends to
the branched covering
\begin{eqnarray*}
\C^2\to \C^2,\quad (z_1,z_2)
\mapsto (\www{z}_1,\www{z}_2) = (z_1,z_2^2).
\end{eqnarray*}
On $C_2^{conf}\cong \C\times\C^*$ with coordinates
$(\www{z}_1,\www{z}_2)$, unit field $e$, Euler field $E$, 
multiplication $\circ$ and partial units $e_1$ and $e_2$ 
are as follows,
\begin{eqnarray*}
e=e_1+e_2=\frac{\paa}{\paa \www{z}_1},\quad 
E=\www{z}_1e
+2\www{z}_2\frac{\paa}{\paa \www{z}_2},\\
\frac{\paa}{\paa\www{z}_2}\circ \frac{\paa}{\paa\www{z}_2}
=\frac{1}{4\www{z}_2}e,\\
u_{1/2}=\www{z}_1\pm\sqrt{\www{z}_2},\quad 
e_{1/2}=\frac{1}{2}e
\pm \sqrt{\www{z}_2}\frac{\paa}{\paa \www{z}_2}.
\end{eqnarray*}
The canonical coordinates and the partial units are 2-valued.
The multiplication does not extend to $\C\times\{0\}$. 

(c) The universal covering $C_2^{univ}\to C_2^{pure}$ can be
written as follows,
\begin{eqnarray*}
C_2^{univ}\cong \C\times\C&\to& \C\times\C^*\cong C_2^{pure},\\
(z_1,\zeta)&\mapsto& (z_1,\exp(\zeta))=(z_1,z_2).
\end{eqnarray*}
In the semisimple F-manifold $C_2^{univ}\cong\C^2$
with coordinates $(z_1,\zeta)$, unit field $e$, Euler field $E$, 
multiplication $\circ$ and partial units $e_1$
and $e_2$ are as follows,
\begin{eqnarray*}
e=e_1+e_2=\frac{\paa}{\paa z_1},
\quad E=z_1e +\frac{\paa}{\paa \zeta},\\
\frac{\paa}{\paa\zeta}\circ\frac{\paa}{\paa\zeta}
=\exp(2\zeta)e,\\
u_{1/2}=z_1\pm \exp(\zeta),\quad 
e_{1/2}=\frac{1}{2}e\pm 
\frac{1}{2}\exp(-\zeta)\frac{\paa}{\paa \zeta},\\
\end{eqnarray*}

(d) Let $(H_\Z,L,\uuuu{e})$ be a unimodular bilinear lattice
of rank 2 with triangular basis $\uuuu{e}$ and matrix
$S=L(\uuuu{e}^t,\uuuu{e})^t
=\begin{pmatrix}1&x\\0&1\end{pmatrix}\in T^{uni}_2(\Z)$
for some $x\in\Z_{\leq 0}$. 
Then 
$$C_2^{S/\{\pm 1\}^2}=C_2^{conf},$$
and 
\begin{eqnarray*}
C_2^{\uuuu{e}/\{\pm 1\}^2}
=\left\{\begin{array}{ll}
C_2^{pure}&\textup{ if }x=0,\quad \textup{the case }A_1^2,\\
C_2^{A_2}&\textup{ if }x=-1,\quad \textup{the case }A_2,\\
C_2^{univ}&\textup{ if }x\leq -2,\quad\textup{the cases }\P^1
\textup{ and beyond}.\end{array}\right. 
\end{eqnarray*}
Here $C_2^{A_2}\cong\C\times \C^*$ with coordinates $(z_1,z_2)$
is the semisimple F-manifold of type $A_2$ in the Examples
\ref{t8.10} (ii), minus its Maxwell stratum $\C\times\{0\}$. 
Explicitly, unit field $e$, Euler field $E$, multiplication 
$\circ$ and partial units are as follows,
\begin{eqnarray*}
e=e_1+e_2=\frac{\paa}{\paa z_1},\quad 
E=z_1e
+\frac{2}{3}z_2\frac{\paa}{\paa z_2},\\
\frac{\paa}{\paa z_2}\circ \frac{\paa}{\paa z_2}=z_2e,\\
u_{1/2}=z_1\pm\frac{2}{3}z_2^{3/2},\quad 
e_{1/2}=\frac{1}{2}e
\pm \frac{1}{2}z_2^{-1/2}\frac{\paa}{\paa z_2}.
\end{eqnarray*}
The canonical coordinates and the partial units are 2-valued.
The semisimple F-manifold $C_2^{A_2}$ extends to $\C\times\{0\}$
with $\C\times\{0\}$ as caustic.
\end{theorem}

{\bf Proof:}
(a) Trivial.

(b) Locally on $C_2^{pure}$ $(\www{z}_1,\www{z}_2)=(z_1,z_2^2)$
is a coordinate system with $\ddd \www{z}_2=2z_2\ddd z_2$
and $\frac{\paa}{\paa \www{z}_2}
=\frac{1}{2z_2}\frac{\paa}{\paa z_2}$. One applies part (a).

(c) Locally on $C_2^{univ}\cong \C\times\C$ 
$(z_1,z_2)=(z_1,e^\zeta)$ is a coordinate system with
$\ddd \zeta=e^{-\zeta}\ddd (e^\zeta)=z_2^{-1}\ddd z_2$
and $\frac{\paa}{\paa \zeta}=z_2\frac{\paa}{\paa z_2}$. 
One applies part (a). 

(d) Recall from Theorem \ref{t7.1} (b)
\begin{eqnarray*}
(\Br_2)_{S/\{\pm 1\}^2}&=& \Br_2 \quad\textup{for any }x,\\
(\Br_2)_{\uuuu{e}/\{\pm 1\}^2} &=& 
\left\{\begin{array}{ll}
\langle\sigma_1^2\rangle 
&\textup{ for }x=0,\textup{ the case } A_1^2,\\
\langle\sigma_1^3\rangle
&\textup{ for }x=-1,\textup{ the case }A_2,\\
\{\id\}
&\textup{ for }x\leq -2,\textup{ the cases }\P^1
\textup{ and beyond.}\end{array}\right.
\end{eqnarray*}
Therefore 
$$C_2^{S/\{\pm 1\}^2}=C_2^{conf},$$
and $C_2^{\uuuu{e}/\{\pm 1\}^2}$ is a cyclic covering of
$C_2^{conf}$ of degree 2 for $x=0$, 
of degree 3 for $x=-1$ and of infinite
degree for $x\leq -2$. With the parts (a)--(c) and the case of
type $A_2$ in Example \ref{t8.10} (ii) this
shows part (d).\hfill$\Box$

\section[Distinguished systems of paths]
{Distinguished systems of paths, two braid group actions}
\label{s8.4}

For the $\Z$-lattice bundles in section \ref{s8.5},
{\it distinguished systems of paths} are needed.
This is a notion which was developed within the theory of
isolated hypersurfac singularities. 
A good source is \cite[5.2 and 5.7]{Eb01}.
Definition \ref{t8.13} and Theorem \ref{t8.14}
define them and describe one action of the braid group
$\Br_n$ on the set of distinguished systems of $n$ paths
with fixed starting point and fixed set of endpoints.

Though there is also an action of the groupoid 
of homotopy classes of paths in $C_n^{pure}$ on the
set of distinguished systems of $n$ paths with fixed
starting point, but variable set of endpoints
(Remarks \ref{t8.15} and Definition/Lemma \ref{t8.16}).
It leads to a second action of the braid group $\Br_n$
on the set of distinguished systems of $n$ paths
with fixed starting point and fixed set of endpoints
(Remark \ref{t8.17} (ii)). 
This second action is not discussed in \cite{Eb01}.
Lemma \ref{t8.18} clarifies the relationship with the
first action. The two actions commute, and one can be 
expressed through the other. Example \ref{t8.19}
gives an example.

\begin{definition}\label{t8.13}
Let $\uuuu{u}=(u_1,u_2,...,u_n)\in C_n^{pure}$ and 
$r\in\R$ with $r>\max(\Ree u_1,...,\Ree u_n)$.

(a) A {\it $(\uuuu{u},r)$-path} $\gamma$ 
\index{$(\uuuu{u},r)$-path}
is a continuous
embedding $\gamma:[0,1]\to \C$ with $\gamma(0)=r$,
$\gamma(1)\in\{u_1,...,u_n\}$ and $\gamma((0,1))
\subset\{z\in\C\,|\, \Ree z<r\}-\{u_1,...,u_n\}$. 

(b) \cite[5.2]{Eb01} A {\it distinguished system of paths}
\index{distinguished system of paths}
is a tuple $(\gamma_1,...,\gamma_n)$ of $(\uuuu{u},r)$-paths
together with a permutation $\sigma\in S_n$ such that
$\gamma_j(1)=u_{\sigma(j)}$, 
$\gamma_i((0,1))\cap \gamma_j((0,1))=\emptyset$ for $i\neq j$,
and $\gamma_1,...,\gamma_n$ leave their common starting point
$r=\gamma_1(0)=...=\gamma_n(0)$ in clockwise order.

(c) \cite[5.7]{Eb01} Two distinguished systems of 
$(\uuuu{u},r)$-paths $(\gamma_1,...,\gamma_n;\sigma)$ and
$(\www{\gamma}_1,...,\www{\gamma}_n;\sigma)$ are
\index{homotopic}
{\it homotopic}, if continuous maps $\Gamma_1,...,\Gamma_n
:[0,1]\times [0,1]\to\C$ exist such that 
$(\Gamma_1(.,s),...,\Gamma_n(.,s);\sigma)$ is for each
$s\in[0,1]$ a distinguished system of $(\uuuu{u},r)$-paths and 
$(\Gamma_1(.,0),...,\Gamma_n(.,0))=(\gamma_1,...,\gamma_n)$,
$(\Gamma_1(.,1),...,\Gamma_n(.,1))=(\www{\gamma}_1,...,
\www{\gamma}_n)$. 
Denote by $\PP(\uuuu{u},r)$ the set of all homotopy classes of
distinguished systems of $(\uuuu{u},r)$-paths.

(d) A {\it standard system of paths} 
\index{standard system of paths}
is a distinguished system
of paths $(\gamma_1,...,\gamma_n)$ with the following properties
(1)--(4).\\
(1) $\sigma\in S_n$ is determined by 
\begin{list}{}{}
\item[(i)]
$\Imm u_{\sigma(1)}\leq ...\leq \Imm u_{\sigma(n)}$.
\item[(ii)]
If $\Imm u_{\sigma(j)}=\Imm u_{\sigma(j+1)}$ then
$\Ree u_{\sigma(j)}<\Ree u_{\sigma(j+1)}$.
\end{list}
(2) Choose $r_1\in \R$ with $\max(\Ree u_1,...,\Ree u_n)
< r_1 < r$. 
Consider $j\in\{1,...,n\}$ and $n_1\in\N$ with
\begin{eqnarray*}
\Imm u_{\sigma(j-n_1)}<\Imm u_{\sigma(j-n_1+1)}
=...=\Imm u_{\sigma(j)}<\Imm u_{\sigma(j+1)}.
\end{eqnarray*}
(3) Choose $\varepsilon_j>0$ so small that
$n_1\varepsilon_j < \Imm u_{\sigma(j)}-\Imm u_{\sigma(j-n_1)}$.
\\
(4) Choose the path $\gamma_{j-k}$ for $k\in\{0,1,...,n_1-1\}$
as follows:\\
First go straight from $r$ to 
$r_1+i(\Imm u_{\sigma(j)}-k\varepsilon_j)$.\\
Then go straight (horizontally to the left) to
$\Ree(u_{\sigma(j-k)})+i(\Imm u_{\sigma(j)}-k\varepsilon_j)$.\\
Then go straight (vertically upwards) to $u_{\sigma(j-k)}$. 
\end{definition}

\begin{figure}
\includegraphics[width=1.0\textwidth]{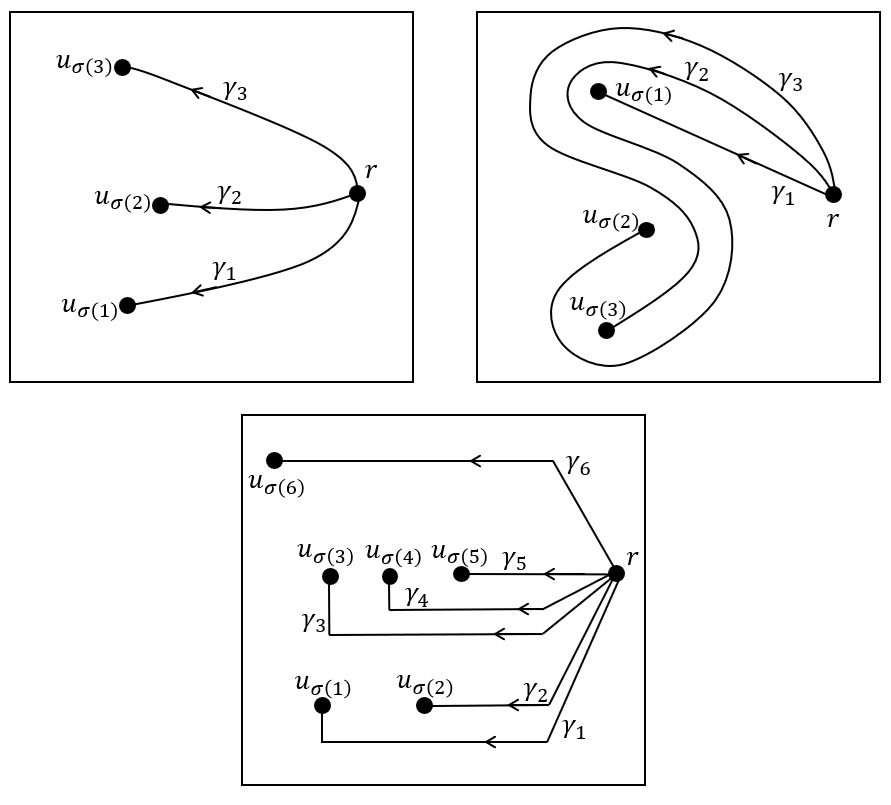}
\caption[Figure 8.3]{Three distinguished systems of paths}
\label{Fig:8.3}
\end{figure}

Obviously any two standard systems of paths are homotopic.
Figure \ref{Fig:8.3} shows three distinguished
systems of paths. The lower one is a standard system of paths,
the upper left one is homotopic to a standard system of paths,
the upper right one not.

\begin{theorem}\label{t8.14}
\cite[5.7]{Eb01}
Let $\uuuu{u}=(u_1,...,u_n)\in C_n^{pure}$ and 
let $r\in\R$ with $r>\max(\Ree u_1,...,\Ree u_n)$. 

(a) (Definition) Define actions of the elementary braids
$\sigma_j$ and $\sigma_j^{-1}$ in $\Br_n$ on the set
$\PP(\uuuu{u},r)$ as follows.
$\sigma_j$ maps the homotopy class of a distinguished system
of paths $(\gamma_1,...,\gamma_n)$ to the homotopy class
of a distinguished system of paths $(\gamma_1',...,\gamma_n')$
where
\begin{eqnarray*}
\gamma_i'&:=&\gamma_i\textup{ for }i\in\{1,...,n\}-\{j,j+1\},\\
\gamma_{j+1}'&:=&\gamma_j,\\
\gamma_j'&:& \textup{\it\  go from }r\textup{\it\ along }\gamma_j
\textup{\it\ almost to }u_{\sigma(j)},\\
&&\textup{\it\  turn around }u_{\sigma(j)}
\textup{\it\  once clockwise},\\ 
&&\textup{\it\ go along }\gamma_j^{-1}\textup{\it\ back to }r,\\
&& \textup{\it\  go from }r\textup{\it\   along }\gamma_{j+1}
\textup{\it\  to }u_{\sigma(j+1)},\\
\sigma'&:=& \sigma\circ (j\, j+1).
\end{eqnarray*}
$\sigma_j^{-1}$ maps the homotopy class of a distinguished system
of paths $(\gamma_1,...,\gamma_n)$ to the homotopy class
of a distinguished system of paths 
$(\gamma_1'',...,\gamma_n'')$
where
\begin{eqnarray*}
\gamma_i''&:=&\gamma_i\textup{ for }i\in\{1,...,n\}-\{j,j+1\},\\
\gamma_j''&:=&\gamma_{j+1},\\
\gamma_{j+1}''&:& \textup{\it\  go from }r\textup{\it\  along }
\gamma_{j+1}\textup{\it\  almost to }u_{\sigma(j+1)},\\
&&\textup{\it\  turn around }u_{\sigma(j+1)}
\textup{\it\  once counter clockwise,}\\
&&\textup{\it\ go along }\gamma_{j+1}^{-1}\textup{\it\  back to }r,\\
&& \textup{\it\  go from }r\textup{\it\  along }\gamma_j
\textup{\it\  to }u_{\sigma(j)},\\
\sigma''&:=& \sigma\circ (j\, j+1).
\end{eqnarray*}
Figure \ref{Fig:8.4} shows the old and new paths under both operations.

\begin{figure}
\includegraphics[width=0.9\textwidth]{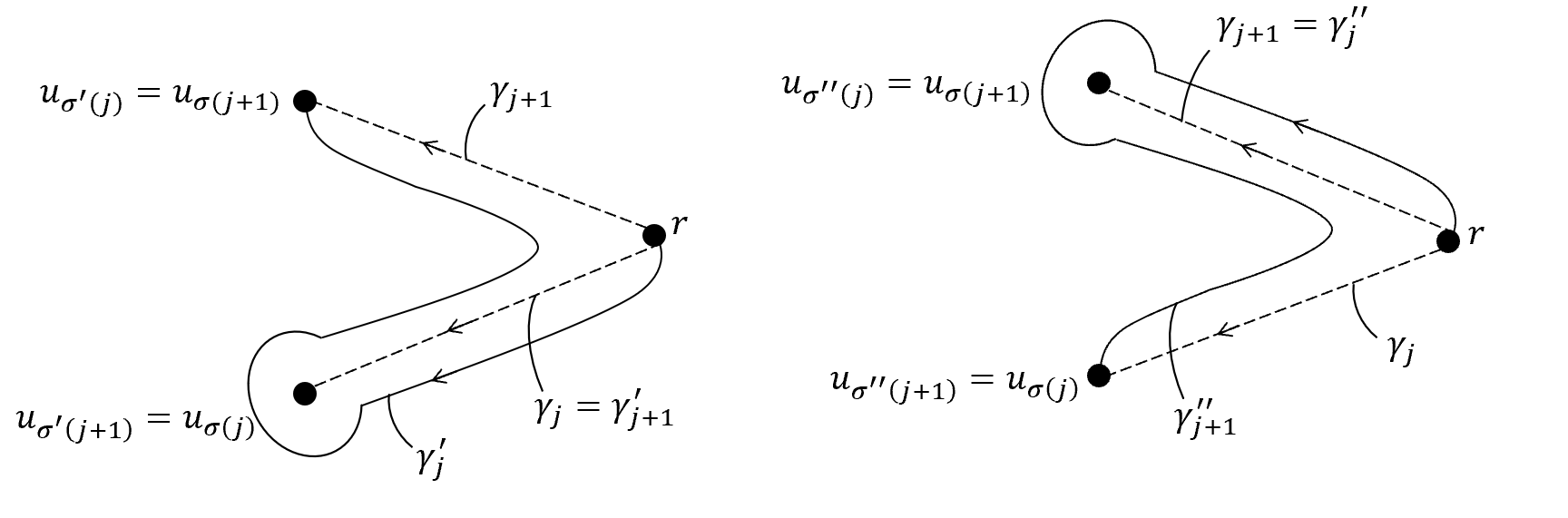}
\caption[Figure 8.4]{Actions of $\sigma_j$ and $\sigma_j^{-1}$ 
on distinguished systems of paths}
\label{Fig:8.4}
\end{figure}

(b) The actions of $\sigma_j$ and $\sigma_j^{-1}$ extend to
an action of the braid group $\Br_n$ from the left on the set 
$\PP(\uuuu{u},r)$ of homotopy classes of distinguished systems
of $(\uuuu{u},r)$-paths. 

(c) The action of $\Br_n$ on the set $\PP(\uuuu{u},r)$ is
transitive.
\end{theorem}

\begin{figure}
\includegraphics[width=1.0\textwidth]{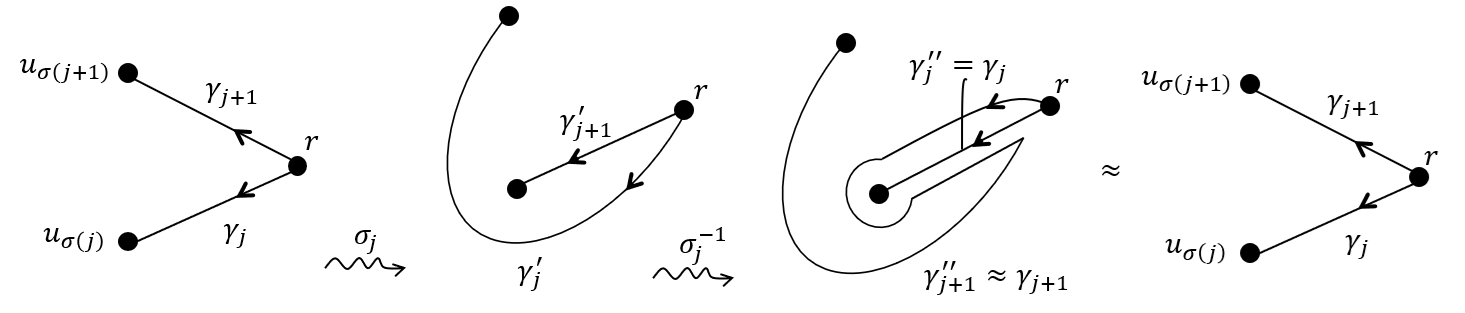}
\caption[Figure 8.5]{The action of $\sigma_j^{-1}\sigma_j$ on distinguished
systems of paths is trivial}
\label{Fig:8.5}
\end{figure}

The sequence of pictures in Figure \ref{Fig:8.5} illustrates
that the action of $\sigma_j^{-1}\sigma_j$ on 
$\PP(\uuuu{u},r)$ is trivial. 
The pictures 5.27 and 5.28 in \cite{Eb01} illustrate that 
the braids $\sigma_j\sigma_{j+1}\sigma_j$ and 
$\sigma_{j+1}\sigma_j\sigma_{j+1}$ act in the same way on 
$\PP(\uuuu{u},r)$. 

In the case of $\uuuu{u}=b_n^{pure}=(i,2i,...,ni)$ 
and $r=1$ ($r=1$ is here an example, any $r>0$ works here) there
is also an action of $\Br_n$ from the right on the set
$\PP(b_n^{pure},1)$. It comes from an action of the
fundamental groupoid of $C_n^{pure}$ on the union 
$\bigcup_{(\uuuu{u},r)}\PP(\uuuu{u},r)$.
This action and the compatibility with the action of $\Br_n$ 
from the left on $\PP(b_n^{conf},1)$ in Theorem \ref{t8.14} 
are not discussed in \cite{Eb01}.
They are given in the points \ref{t8.15} to \ref{t8.19}.

\begin{remark}\label{t8.15}
\cite[1.7]{Sp66}
For $\uuuu{u}^0$ and $\uuuu{u}^1\in C_n^{pure}$ denote by
$H(\uuuu{u}^0,\uuuu{u}^1)$ the set of homotopy classes
\index{homotopy class of a path}
of (continuous maps =) paths $\alpha:[0,1]\to C_n^{pure}$
with $\alpha(0)=\uuuu{u}^0$ and $\alpha(1)=\uuuu{u}^1$.
The tuple $(H(\uuuu{u}^0,\uuuu{u}^1))_{\uuuu{u}^0,\uuuu{u}^1
\in C_n^{pure}}$ is a {\it groupoid} \index{groupoid} with
$[\alpha_0][\alpha_1]=[\alpha_0\alpha_1]\in H(\uuuu{u}^0,
\uuuu{u}^2)$ for 
$[\alpha_0]\in H(\uuuu{u}^0,\uuuu{u}^1)$,
$[\alpha_1]\in H(\uuuu{u}^1,\uuuu{u}^2)$.
The multiplication is associative.
$[\alpha][\alpha^{-1}]=[\textup{trivial path}]\in 
H(\uuuu{u}^0,\uuuu{u}^0)$ for 
$[\alpha]\in H(\uuuu{u}^0,\uuuu{u}^1)$.
See \cite[1.7]{Sp66} for more details on this groupoid. 
\end{remark}

\begin{definition/lemma}\label{t8.16}
(a) Let $\alpha:[0,1]\to C_n^{pure}$ be a path with 
$[\alpha]\in H(\uuuu{u}^0,\uuuu{u}^1)$.
We will construct a natural bijection 
$\Phi(\alpha):\PP(\uuuu{u}^0,r_0)\to\PP(\uuuu{u}^1,r_1)$
(in the notation $\Phi(\alpha)$ we should include
$r_0$ and $r_1$, but for simplicity we drop them).

Choose $r\in\R$ with 
$$r\geq \max(r_0,r_1),\quad 
r>\max_{i\{1,...,n\}}\max_{t\in [0,1]}\Ree\alpha_i(t).$$ 
Because there are natural 
bijections $\PP(\uuuu{u}^0,r_0)\to \PP(\uuuu{u}^0,r)$
and $\PP(\uuuu{u}^1,r_1)\to\PP(\uuuu{u}^1,r)$, we can
replace $r_0$ and $r_1$ by $r$.
Choose a continuous family in $s\in[0,1]$ of homeomorphisms
$\Psi_s:\C\to\C$ with the properties
\begin{eqnarray*}
&& \Psi_s|_{\{z\in\C\,|\, \Ree z\geq r\}}=\id,\\
&& \Psi_s(\{u_1^0,...,u_n^0\})=\{\alpha_1(s),...,\alpha_n(s)\}.
\end{eqnarray*}

For a distinguished system of paths
$(\uuuu{\gamma};\sigma)=(\gamma_1,...,\gamma_n;\sigma)$ with 
$[(\uuuu{\gamma};\sigma)]\in\PP(\uuuu{u}^0,r)$,
the composition $(\Psi_s\circ\uuuu{\gamma};\sigma)
=(\Psi_s\circ\gamma_1,...,\Psi_s\circ\gamma_n;\sigma)$
is a distinguished system of paths with 
$[(\Psi_s\circ\uuuu{\gamma};\sigma)]\in\PP(\alpha(s),r)$.
The class $[(\Psi_s\circ\uuuu{\gamma};\sigma)]$
depends only on the class $[(\uuuu{\gamma};\sigma)]$
and on the path $\alpha$. Define
\begin{eqnarray*}
\Phi(\alpha)([(\uuuu{\gamma};\sigma)])
:= [(\Psi_1\circ\uuuu{\gamma};\sigma)]\in 
\PP(\uuuu{u}^1,r).
\end{eqnarray*}

(b) The map $\Phi(\alpha):\PP(\uuuu{u}^0,r)\to
\PP(\uuuu{u}^1,r)$ depends only on the homotopy class
$[\alpha]\in H(\uuuu{u}^0,\uuuu{u}^1)$.

(c) If $[\alpha_1]\in H(\uuuu{u}^0,\uuuu{u}^1)$ and
$[\alpha_2]\in H(\uuuu{u}^1,\uuuu{u}^2)$ then
$[\alpha_1\alpha_2]\in H(\uuuu{u}^0,\uuuu{u}^2)$ and
$$\Phi(\alpha_1\alpha_2)=\Phi(\alpha_2)\Phi(\alpha_1).$$
The groupoid $(H(\uuuu{u}^0,\uuuu{u}^1))_{\uuuu{u}^0,
\uuuu{u}^1\in C_n^{pure}}$ acts from the right
on the tuple $(\PP(\uuuu{u},r))_{\uuuu{u}\in C_n^{pure},
r>\max(\Ree u_1,...,\Ree u_n)}$.
\end{definition/lemma}

\begin{remarks}\label{t8.17}
(i) Let $\uuuu{u}\in C_n^{pure}$ and $r>\max(\Ree u_1,...,
\Ree u_n)$. Let $\tau\in S_n$. Then $\tau$ as deck transformation
of the covering $\pr_n^{p,c}:C_n^{pure}\to C_n^{conf}$ maps
$\uuuu{u}=(u_1,...,u_n)$ to
$\tau(\uuuu{u})=(u_{\tau^{-1}(1)},...,u_{\tau^{-1}(n)})$,
see the Remarks \ref{t8.2} (i) and (vii).
A distinguished basis $(\uuuu{\gamma};\sigma)$ for 
$\uuuu{u}$ and $r$ is almost the same as the distinguished basis
$(\uuuu{\gamma},\tau\sigma)$ for $\tau(\uuuu{u})$ and $r$.
They just differ by the indexing of the entries of $\uuuu{u}$.
This gives a natural bijection 
$\PP(\uuuu{u},r)\to \PP(\tau(\uuuu{u}),r)$,
$[(\uuuu{\gamma};\sigma)]\mapsto [(\uuuu{\gamma};\tau\sigma)]$.

(ii) A closed path $\beta$ in $C_n^{conf}$ with starting point
and end point $b_n^{conf}$ represents a braid $[\beta]\in \Br_n$.
It lifts to a path $\www{\beta}$ in $C_n^{pure}$ from
$b_n^{pure}=(i,2i,...,ni)$ to 
$\sigma(b_n^{pure})=(\sigma^{-1}(1)i,\sigma^{-1}(2)i,...,
\sigma^{-1}(n)i)$ 
for some $\sigma\in S_n$ with homotopy class
$[\www{\beta}]\in H(b_n^{pure},\sigma(b_n^{pure})]$.
Now $\Phi(\www{\beta})$ maps $\PP(b_n^{pure},1)$
to $\PP(\sigma(b_n^{pure}),1)$ which can be identified
with $\PP(b_n^{pure},1)$ because of (i).
Therefore $[\beta]$ acts on $\PP(b_n^{pure},1)$. 
We obtain an action of $\Br_n$ on $\PP(b_n^{pure},1)$,
which is an action from the right because of Lemma \ref{t8.16}
(c).
\end{remarks}

Remark \ref{t8.17} (ii) gives an intuitive geometric action
of $\Br_n$ on $\PP(b_n^{pure},1)$ from the right.
It comes from moving the endpoints of a distinguished system
of paths along a braid and deforming the distinguished system
of paths accordingly.
Lemma \ref{t8.18} connects it to the more algebraic action of
$\Br_n$ on $\PP(b_n^{pure},1)$ from the left in Theorem
\ref{t8.14} (a)+(b). 
Lemma \ref{t8.18} shows first that these actions commute and
second that one action can be expressed through the other
action.

\begin{lemma}\label{t8.18}
The action of $\Br_n$ on $\PP(b_n^{pure},1)$ from the left
in Theorem \ref{t8.14} and the action of $\Br_n$ on
$\PP(b_n^{pure},1)$ from the right in Remark \ref{t8.17} (ii)
commute,
\begin{eqnarray*}
\beta_1([(\uuuu{\gamma};\sigma)].\beta_2)
= (\beta_1[(\uuuu{\gamma};\sigma)]).\beta_2
\quad\textup{for }\beta_1,\beta_2\in\Br_n,
[(\uuuu{\gamma};\sigma)]\in\PP(b_n^{pure},1),
\end{eqnarray*}
and satisfy
\begin{eqnarray*}
\beta[(\uuuu{\gamma};\sigma)] = [(\uuuu{\gamma};\sigma)].\beta
\quad\textup{for }\beta\in\Br_n,[(\uuuu{\gamma};\sigma)]
\in \PP(b_n^{pure},1).
\end{eqnarray*}
\end{lemma}

{\bf Proof:}
The action on the right comes from a continuous family of 
homeomorphisms as in Lemma \ref{t8.16} (a) which deforms
the distinguished system of paths. The action on the left
comes from a construction of new paths which follow closely
the old paths. The new paths are deformed together with the
old paths. A deformation of the old and new paths
does not change this construction of new paths from old paths.
Therefore the actions commute.

\begin{figure}
\includegraphics[width=0.7\textwidth]{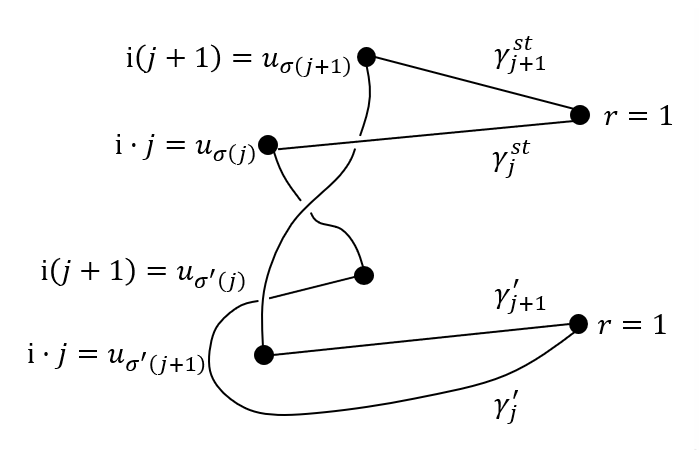}
\caption[Figure 8.6]{Moving along a braid in $C_n^{pure}$
changes a distinguished system of paths}
\label{Fig:8.6}
\end{figure}

The action of $\Br_n$ from the left on $\PP(b_n^{pure},1)$
is transitive by Theorem \ref{t8.14} (c).
Therefore it is sufficient to show for an elementary braid
$\sigma_j$ and a standard system of paths
$(\uuuu{\gamma}^{st};\sigma)$ with $[(\uuuu{\gamma}^{st};\sigma)]
\in \PP(b_n^{conf},1)$ 
$$\sigma_j[(\uuuu{\gamma}^{st};\sigma)] =
[(\uuuu{\gamma}^{st};\sigma)].\sigma_j.$$
This is obvious respectively shown in the picture in Figure
\ref{Fig:8.6}. \hfill$\Box$

\begin{example}\label{t8.19}
Figure \ref{Fig:8.7} shows for $n=3$ on the left the action of 
$\sigma_2^{-1}\sigma_1$ on $[(\uuuu{\gamma}^{st};\sigma)]$
from the left and on the right the action of 
$\sigma_2^{-1}\sigma_1$ on $[(\uuuu{\gamma}^{st};\sigma)]$
from the right.

\begin{figure}
\includegraphics[width=1.0\textwidth]{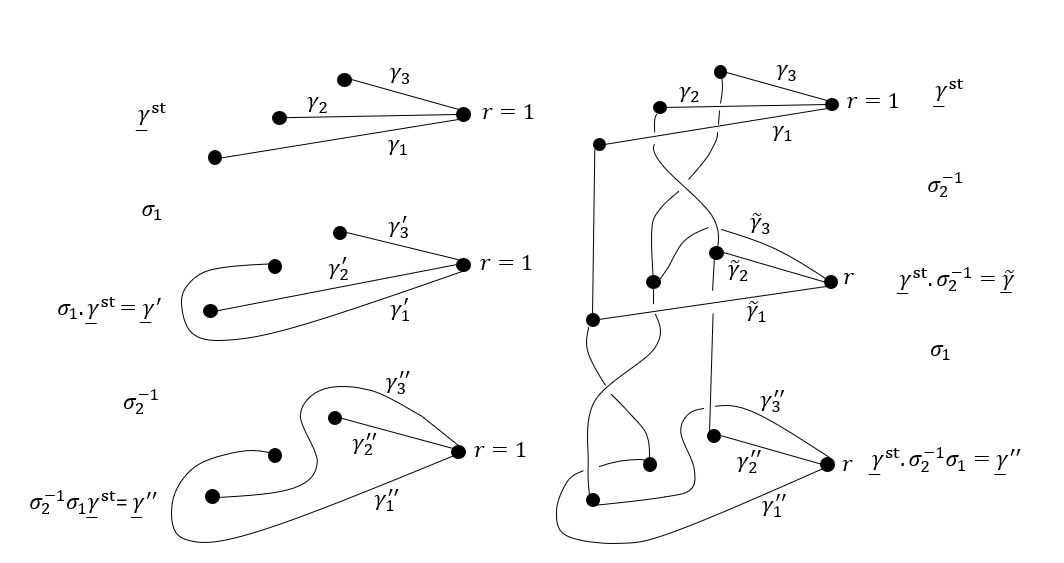}
\caption[Figure 8.7]{Actions of $\sigma_2^{-1}\sigma_1$ from
the left and from the right: same result}
\label{Fig:8.7}
\end{figure}
\end{example}

\section[Two families of $\Z$-lattice structures over
$C_n^{\uuuu{e}/\{\pm 1\}^n}$]{Two families of $\Z$-lattice structures over the 
manifold $C_n^{\uuuu{e}/\{\pm 1\}^n}$}
\label{s8.5}

Let $(H_\Z,L,\uuuu{e})$ be a unimodular bilinear lattice of
rank $n\in\N$ with a triangular basis $\uuuu{e}$ and 
matrix $S=L(\uuuu{e}^t,\uuuu{e})^t\in T^{uni}_n(\Z)$. 
The manifold $C_n^{\uuuu{e}/\{\pm 1\}^n}$ is not just 
a semisimple F-manifold with Euler field.
It is the carrier of much richer structures, which we call
{\it $\Z$-lattice structures} and which we explain in this
section.

\begin{remark}\label{t8.20}
A $\Z$-lattice structure leads often in several steps
to a strong enrichment of the F-manifold structure, namely a 
Dubrovin-Frobenius manifold. We will not go through these
steps and not treat Dubrovin-Frobenius manifolds here,
but we will say a few words in this remark.
Details on the steps are given in \cite{Sa02} \cite{He02} 
\cite{He03}.

The $\Z$-lattice structure can be extended to a holomorphic
vector bundle on $\C\times C_n^{\uuuu{e}/\{\pm 1\}^n}$
with a meromorphic connection with a logarithmic pole
along a certain discriminant which will be defined below
and a $\Z$-lattice bundle on the complement.
A Fourier-Laplace transformation of this bundle
along the first factor $\C$ of the base space leads to another
holomorphic vector bundle 
on $\C\times C_n^{\uuuu{e}/\{\pm 1\}^n}$ 
with a meromorphic connection
with a pole of Poincar\'e rank 1 along the hypersurface
$\{0\}\times C_n^{\uuuu{e}/\{\pm 1\}^n}$,
a $\Z$-lattice bundle on the complement and a certain pairing.
In the notation of \cite{He03} it is a {\it TEZP-structure}
(Twistor Extension $\Z$-lattice Pairing).
From this and from some additional choice one can usually
construct a Dubrovin-Frobenius manifold structure
on $C_n^{\uuuu{e}/\{\pm 1\}^n}$ or at least on the complement
of a hypersurface in it. 
See \cite{Sa02} \cite{He02} \cite{He03} for this construction.
Dubrovin-Frobenius manifolds were first defined by Dubrovin
\cite{Du92}\cite{Du96}. \index{Dubrovin-Frobenius manifold}
\index{Frobenius manifold}
A definition is given above in Remark \ref{t8.9} (viii). 
Here we will restrict to explain two natural families of 
$\Z$-lattice structures over $C_n^{\uuuu{e}/\{\pm 1\}^n}$. 
\end{remark}

Definition/Lemma \ref{t8.21} (c) presents the notion of a single 
$\Z$-lattice structure. The construction of natural
families of $\Z$-lattice structures over $C_n^{univ}$
and $C_n^{\uuuu{e}/\{\pm 1\}^n}$ comes 
in Definition \ref{t8.22} (c) and Theorem \ref{t8.23} (d).
The unimodular lattice $(H_\Z,L,\uuuu{e})$ with distinguished 
basis $\uuuu{e}$ induces such families. 
The construction uses distinguished systems of paths. 
Theorem \ref{t8.23} (a)--(c) 
discusses how $\Z$-lattice structures over different points
in $C_n^{univ}$ in such a family of $\Z$-lattice structures 
are related by the actions of braids.

\begin{definition/lemma}\label{t8.21}
Let $(H_\Z,L,\uuuu{e})$ and $S$ be as above.

(a) (Definition) The {\sf discriminant} 
\index{discriminant}\index{$D_{1,n}^{univ}$} 
$D_{1,n}^{univ}\subset \C\times C_n^{univ}$ 
is the smooth complex hypersurface
$$\{(\tau,b)\in \C\times C_n^{univ}\,|\, \LL(b)(\tau)=0\}.$$
Recall the map $\LL:C_n^{univ}\to \C[x]_n^{reg}$ in 
Remark \ref{t8.9} (iv). The polynomial 
$\LL(b)\in \C[x]_n^{reg}\cong C_n^{conf}$
corresponds to $\pr_n^{u,c}(b)\in C_n^{conf}$. 

(Trivial lemma) The projection 
$\pr_2:\C\times C_n^{univ}\to C_n^{univ}$,
$(\tau,b)\mapsto b$, restricts to a covering 
$D_{1,n}^{univ}\to C_n^{univ}$ of degree $n$. 

(b) (Trivial lemma) The group $\Br_n$ of deck transformations
of the universal covering $\pr_n^{u,c}:C_n^{univ}\to C_n^{conf}$
extends with $\id_\C$ on the first factor of 
$\C\times C_n^{univ}$ to a group of automorphisms of
$\C\times C_n^{univ}$. It leaves the discriminant
$D_{1,n}^{univ}$ invariant. 

(Definition) The image of $D_{1,n}^{univ}$ under the group
$\id_\C\times(\Br_n)_{\uuuu{e}/\{\pm 1\}^n}$ is the 
discriminant 
\index{$D_{1,n}^{\uuuu{e}/\{\pm 1\}^n}$}
$D_{1,n}^{\uuuu{e}/\{\pm 1\}^n}
\subset \C\times C_n^{\uuuu{e}/\{\pm 1\}^n}$,
$$D_{1,n}^{\uuuu{e}/\{\pm 1\}^n}
=\{(\tau,b)\in \C\times C_n^{\uuuu{e}/\{\pm 1\}^n}\,|\, 
\LL(b)(\tau)=0\},$$
where now $\LL$ is the Lyashko-Looijenga map 
$\LL:C_n^{\uuuu{e}/\{\pm 1\}^n}\to \C[x]_n^{reg}.$

(c) (Definition) A $\Z$-lattice structure means an (automatically
flat) $\Z$-lattice bundle over $\C-\{n\textup{ points}\}$.

(Remark) Therefore a $\Z$-lattice bundle over 
$\C\times C_n^{univ}-D_{1,n}^{univ}$ is considered as a 
({\sf flat} or {\sf isomonodromic}) family of 
$\Z$-lattice structures over $C_n^{univ}$.  
A $\Z$-lattice bundle over 
$\C\times C_n^{\uuuu{e}/\{\pm 1\}^n} -
D_{1,n}^{\uuuu{e}/\{\pm 1\}^n}$ is considered as a 
({\sf flat} or {\sf isomonodromic}) family of 
$\Z$-lattice structures over $C_n^{\uuuu{e}/\{\pm 1\}^n}$.

(Remark) The $\Z$-lattice structures which we will construct
will come equipped with an (automatically flat) even or odd
intersection form on each fiber.
\end{definition/lemma}

Recall that a unimodular bilinear lattice $(H_\Z,L)$ with
triangular basis $\uuuu{e}$ and matrix 
$S=L(\uuuu{e}^t,\uuuu{e})^t\in T^{uni}_n(\Z)$ was fixed
at the beginning of this section. 
Together with the choice of a pair $(\uuuu{u},r)\in 
C_n^{pure}\times \R$ with $r>\max(\Ree u_1,...,\Ree u_n)$
and the choice of a homotopy class $[(\uuuu{\gamma};\sigma)]\in 
\PP(\uuuu{u},r)$ of a distinguished system of paths
$(\uuuu{\gamma};\sigma)$, it induces a $\Z$-lattice
structure on $\C-\{u_1,...,u_n)$. It also induces
(without any additional choice) a natural 
$\Z$-lattice bundle over $C_n^{univ}$. 
The two constructions are given in Definition \ref{t8.22}
(b) and (c). 

\begin{definition}\label{t8.22}
(a) Let $\uuuu{u}\in C_n^{pure}$ and $r\in\R$ with 
$r>\max(\Ree u_1,...,\Ree u_n)$.
Let $\gamma:[0,1]\to \C$ be a $(\uuuu{u},r)$-path
(Definition \ref{t8.13} (a)).
An {\it associated  loop} $\delta:[0,1]\to\C$
\index{associated loop} 
is a closed path which starts and ends at $r$,
which first goes along $\gamma$ almost to $\gamma(1)$,
which then moves on a small circle once counterclockwise 
around $\gamma(1)$  and which finally goes along 
$\gamma^{-1}$ back to $r$.

Figure \ref{Fig:8.8} illustrates in a picture three loops
which are associated to three paths.

\begin{figure}
\includegraphics[width=0.5\textwidth]{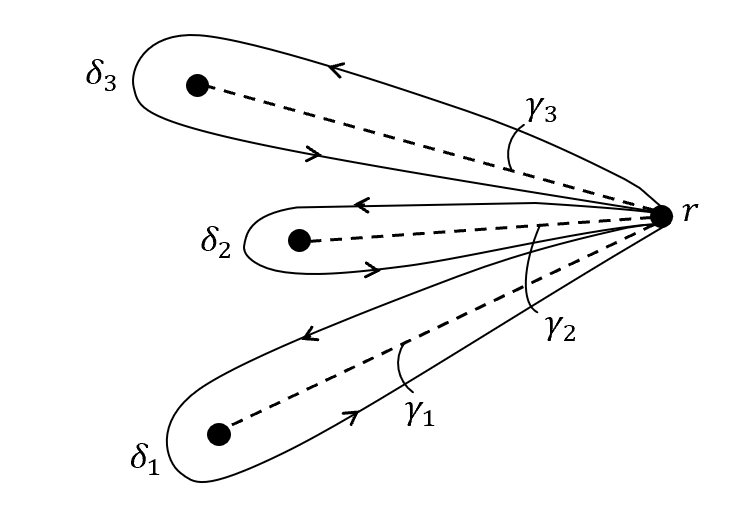}
\caption[Figure 8.8]{Three loops associated to three paths}
\label{Fig:8.8}
\end{figure}

(b) Let $\uuuu{u}\in C_n^{pure}$ and $r\in\R$ with 
$r>\max(\Ree u_1,...,\Ree u_n)$.
Let $(\uuuu{\gamma};\sigma)$ be a distingushed system of
paths with $[(\uuuu{\gamma};\sigma)]\in 
\PP(\uuuu{u},r)$.
Let $\uuuu{f}\in\BB^{dist}$ be a distinguished basis of $H_\Z$.
Let $k\in\{0,1\}$. 
Then $V_\Z^{(k)}((\uuuu{u},r),[(\uuuu{\gamma};\sigma)],
\uuuu{f})$ is the $\Z$-lattice bundle on $\C-\{u_1,...,u_n\}$
which is determined by the following conditions (i) and (ii):
\begin{list}{}{}
\item[(i)]
The restriction of the bundle to 
$\{z\in\C\,|\, \Ree z\geq r\}$ is trivial, and each fiber over
a point in this set is identified with $H_\Z$.
\item[(ii)]
Let $\delta_1,...,\delta_n$ be loops which are associated to
$\gamma_1,...,\gamma_n$. The monodromy along the loop $\delta_j$
is with the identification of the fiber 
$V_\Z^{(k)}((\uuuu{u},r),[(\uuuu{\gamma};\sigma)],\uuuu{f})_r$ 
with $H_\Z$ the automorphism
$s^{(k)}_{f_j}:H_\Z\to H_\Z$.
\end{list}

It is well known that (i) and (ii) determine uniquely
a $\Z$-lattice bundle on $\C-\{u_1,...,u_n\}$.
This bundle 
\index{$V_\Z^{(k)}((\uuuu{u},r),[(\uuuu{\gamma};\sigma)],\uuuu{f})$} 
$V_\Z^{(k)}((\uuuu{u},r),[(\uuuu{\gamma};\sigma)],
\uuuu{f})$ is a $\Z$-lattice structure in the sense of 
Definition \ref{t8.21} (c).

(c) Let $k\in\{0,1\}$. $V_\Z^{(k),univ}$ 
\index{$V_\Z^{(k),univ}$}is the 
$\Z$-lattice bundle on $\C\times C_n^{univ}-D_{1,n}^{univ}$
which is determined by the following conditions (i) and (ii):
\begin{list}{}{}
\item[(i)]
The restriction of the bundle to 
$\{(z,b)\in\C\times C_n^{univ}\,|\, 
\Ree z>\max(\Ree u\,|\, (u,b)\in D_{1,n}^{univ})\}$
is trivial, and each fiber over a point in this set 
is identified with $H_\Z$.
\item[(ii)]
The restriction of $V_\Z^{(k),univ}$ to 
$(\C-\{i,2i,...,ni\})\times \{b_n^{univ}\}$ is equal to
the $\Z$-lattice bundle 
$V_\Z^{(k)}((b_n^{pure},1),[(\uuuu{\gamma}^{st};\id)],
\uuuu{e})$. Here $(\uuuu{\gamma}^{st},\id)$ is a standard 
system of paths which start at $r=1$
and end at the entries of $b_n^{pure}=(i,2i,...,ni)$. 
\end{list}

Because $\pr_2:\C\times C_n^{univ}$ restricts to a trivial
$n$-sheeted covering $\pr_2:D_{1,n}^{univ}\to C_n^{univ}$
over the simply connected manifold $C_n^{univ}$,
(i) and (ii) determine uniquely a $\Z$-lattice bundle
on $\C\times C_n^{univ}-D_{1,n}^{univ}$. 
The bundle $V_\Z^{(k),univ}$ is a flat family of 
$\Z$-lattice structures over the manifold $C_n^{univ}$
in the sense of Definition \ref{t8.21} (c).
\end{definition}

The next theorem is central in this section.
It will construct in part (d) a $\Z$-lattice bundle over
$\C\times C_n^{\uuuu{e}/\{\pm 1\}^n}-
D_{1,n}^{\uuuu{e}/\{\pm 1\}^n}$.
This is a flat family of $\Z$-lattice structures over
$C_n^{\uuuu{e}/\{\pm 1\}^n}$ in the sense of Definition
\ref{t8.21} (c). 
The parts (a), (b) and (c) prepare this.

Part (a) starts with a $\Z$-lattice bundle 
$V_\Z^{(k)}((\uuuu{u},r),[(\uuuu{\gamma};\sigma)],\uuuu{f})$
on $\C-\{u_1,...,u_n\}$ as in Definition \ref{t8.22} (b).
It answers the question how one has to change $\uuuu{f}$
in order to obtain the same $\Z$-lattice bundle if one changes
the distinguished system of paths from 
$[(\uuuu{\gamma};\sigma)]$ to $\beta[(\uuuu{\gamma};\sigma)]$
with some braid $\beta\in\Br_n$.

Part (b) starts with the $\Z$-lattice bundle
$V_\Z^{(k),univ}$ on $\C\times C_n^{univ}-D_{1,n}^{univ}$
as in Definition \ref{t8.22} (c).
It considers for a point $\uuuu{\www{u}}\in C_n^{univ}$ 
and its image $\uuuu{u}^1=\pr^{u,p}(\uuuu{\www{u}})
\in C_n^{pure}$ the restriction of $V_\Z^{(k),univ}$ to a
$\Z$-lattice bundle on 
$(\C-\{u_1^1,...,u_n^1\})\times\{\uuuu{\www{u}}\}$.
It studies how this restriction is obtained from the 
restriction to $(\C-\{i,2i,...,ni\})\times\{b_n^{univ}\}$
if one moves from $b_n^{univ}$ to $\uuuu{\www{u}}$.
The groupoid action in Lemma \ref{t8.18} is used.

Part (c) specializes the situation in part (b) to the case
where $\pr^{u,c}(\uuuu{\www{u}})=b_n^{conf}$. Then a path in
$C_n^{univ}$ from $b_n^{univ}$ to $\uuuu{\www{u}}$ corresponds
to a braid $\beta\in \Br_n$. Part (c) shows with the help of
the parts (a) and (b) and Lemma \ref{t8.18} that then the 
restriction of $V_\Z^{(k),univ}$ to
$(\C-\{i,2i,...,ni\})\times\{\uuuu{\www{u}}\}$ is simply
$V_\Z^{(k)}((b_n^{pure},1),[(\uuuu{\gamma^{st}};\id)],
\beta^{-1}(\uuuu{e}))$.

\begin{theorem}\label{t8.23}
(a) Let $(\uuuu{\gamma};\sigma)$ be a 
distinguished system of paths
with $[(\uuuu{\gamma};\sigma)]\in\PP(b_n^{pure},1)$.
Let $\uuuu{f}$ be a distinguished basis of $H_\Z$.
Let $k\in\{0,1\}$. Let $\beta\in \Br_n$. Then
\begin{eqnarray*}
V_\Z^{(k)}((b_n^{pure},1),[(\uuuu{\gamma};\sigma)],
\uuuu{f})
=V_\Z^{(k)}((b_n^{pure},1),\beta[(\uuuu{\gamma};\sigma)],
\beta(\uuuu{f})).
\end{eqnarray*}

(b) Let $\alpha:[0,1]\to C_n^{pure}$ be a path with
$[\alpha]\in H(b_n^{pure},\uuuu{u}^1)$ for some
$\uuuu{u}^1\in C_n^{pure}$. 
Let $\www{\alpha}:[0,1]\to C_n^{univ}$ be the lift to
$C_n^{univ}$ of $\alpha$ with starting point 
$\www{\alpha}(0)=b_n^{univ}$. Let $k\in\{0,1\}$. 
The restriction of $V_\Z^{(k),univ}$ to
$(\C-\{u_1^1,...,u_n^1\})\times\{\www{\alpha}(1)\}$ is
equal to 
$V_\Z^{(k)}((\uuuu{u}^1,r),
\Phi(\alpha)[(\uuuu{\gamma}^{st};\id)],\uuuu{e})$
for some $r\in \R_{\geq 1}$ with 
$r>\max(\Ree u_1^1,...,\Ree u_n^1)$. 

(c) Let $\beta\in \Br_n$. 
Let $\www{\beta}:[0,1]\to C_n^{univ}$ be the lift to $C_n^{univ}$
with $\www{\beta}(0)=b_n^{univ}$ of a loop in $C_n^{conf}$
which represents $\beta$. Let $k\in\{0,1\}$. 
Then the restriction of $V_\Z^{(k),univ}$ to
$(\C-\{i,2i,...,ni\})\times \{\www{\beta}(1)\}$
is equal to $V_\Z^{(k)}((b_n^{pure},1),
[(\uuuu{\gamma}^{st},\id)],\beta^{-1}(\uuuu{e}))$.

(d) For $\beta\in (\Br_n)_{\uuuu{e}/\{\pm 1\}^n}$ denote
the deck transformation of $C_n^{univ}$ which it induces
by $\beta^{deck}:C_n^{univ}\to C_n^{univ}$.
It extends to an automorphism 
$\beta^{deck,V,(k)}$ of $V_\Z^{(k),univ}$ which is
$(\id_{H_\Z},\id_\C\times\beta^{deck})$ over 
$\{(z,b)\in\C\times C_n^{univ}\,|\, \Ree z> 
\max(\Ree u\,|\,(u,b)\in D_{1,n}^{univ})\}$. 
Therefore the quotient 
$V_\Z^{(k),univ}/ (\Br_n)_{\uuuu{e}/\{\pm 1\}^n}^{deck,V,(k)}$
is a $\Z$-lattice bundle over 
$\C\times C_n^{\uuuu{e}/\{\pm 1\}^n}
-D_{1,n}^{\uuuu{e}/\{\pm 1\}^n}$
whose restriction to
$$(\{(z,b)\in \C\times C_n^{\uuuu{e}/\{\pm 1\}^n}\,|\,
\Ree z>\max(\Ree u\,|\, (u,b)\in 
D_{1,n}^{\uuuu{e}/\{\pm 1\}^n})\}$$ 
is trivial. We call this bundle
\index{$V_\Z^{(k),\uuuu{e}/\{\pm 1\}^n}$}
$V_\Z^{(k),\uuuu{e}/\{\pm 1\}^n}$.
It is a family of $\Z$-lattice structures over 
$C_n^{\uuuu{e}/\{\pm 1\}^n}$. 
\end{theorem}

{\bf Proof:}
(a) Let $\delta_1,...,\delta_n$ be loops associated to
$\gamma_1,...,\gamma_n$. By definition, in the bundle
$V_\Z^{(k)}((b_n^{pure}),[(\uuuu{\gamma};\sigma)],\uuuu{f})$
the local monodromy along the loop $\delta_j$
is given by $s^{(k)}_{f_j}$. 
Write $[(\uuuu{\gamma}';\sigma')]
:= \beta[(\uuuu{\gamma};\sigma)]$ and 
$\uuuu{f}':=\beta(\uuuu{f})$. 

We have to show that in the bundle 
$V_\Z^{(k)}((b_n^{pure}),[(\uuuu{\gamma};\sigma)],\uuuu{f})$
the local monodromy around the loop $\delta_j'$
is given by $s^{(k)}_{f_j'}$. 

First we consider $\beta=\sigma_j$. Then
\begin{eqnarray*}
&&\gamma_i'=\gamma_i\textup{ for }i\in\{1,...,\}-\{j,j+1\},\\
&& \gamma_{j+1}'=\gamma_j,\\
&&\gamma_j'\textup{ is homotopic to }\delta_j^{-1}\gamma_{j+1},\\
\textup{so }&&\delta_i'=\delta_i\textup{ for }
i\in\{1,...,n\}-\{j,j+1\},\\
&&\delta_{j+1}'=\delta_j,\\
&&\delta_j'\textup{ is homotopic to }\delta_j^{-1}\delta_{j+1}
\delta_j.
\end{eqnarray*}
Also 
\begin{eqnarray*}
&&f_i'=f_i\textup{ for }i\in\{1,...,n\}-\{j,j+1\},\\
&&f_{j+1}'=f_j,\\
&&f_j'=s_{f_j}^{(k)}(f_{j+1}).
\end{eqnarray*}
Because of 
\begin{eqnarray*}
s^{(k)}_{f_j'} = s^{(k)}_{s^{(k)}_{f_j}(f_{j+1})}
=s^{(k)}_{f_j} s^{(k)}_{f_{j+1}} (s^{(k)}_{f_j})^{-1}
\end{eqnarray*}
the local monodromy along $\delta_i'$ is indeed given
by $s^{(k)}_{f_i'}$ for $i\in\{1,...,n\}$ 
(recall that a composition of paths is read from the left,
a composition of automorphisms is read from the right).
The case $\beta=\sigma_j^{-1}$ is treated analogously.
The general case $\beta\in \Br_n$ follows.

(b) In the restriction of $V_\Z^{(k),univ}$ to
$(\C-\{i,2i,...,ni\})\times \{b_n^{univ}\}$
the local monodromies along loops associated to a standard
system of paths $(\uuuu{\gamma}^{st};\id)$ are
$s^{(k)}_{e_1},...,s^{(k)}_{e_n}$. 

Moving in the base space $C_n^{univ}$ from $b_n^{univ}$ to
$\www{\alpha}(1)$, the standard system of paths is deformed
to a representative of $\Phi(\alpha)[(\uuuu{\gamma}^{st};\id)]$,
associated loops are deformed accordingly, and the local
monodromies along the deformed associated loops are still
$s^{(k)}_{e_1},...,s^{(k)}_{e_n}$.
This shows the claim.

(c) One has to apply part (b), Remark \ref{t8.17} (ii), 
Lemma \ref{t8.18} and part (a), as follows.
Let $\oooo{\beta}:[0,1]\to C_n^{pure}$ be the lift of a 
representative of $\beta$ with $\oooo{\beta}(0)=b_n^{pure}$.
Then $\oooo{\beta}(1)=\sigma(b_n^{pure})$ for some
$\sigma\in S_n$. We have the following equalities
\begin{eqnarray*}
&&\textup{the restriction of }V_\Z^{(k),univ}\textup{ to }
(\C-\{i,2i,...,ni\})\times \{\www{\beta}(1)\}\\
&=& 
V_\Z^{(k)}((\sigma(b_n^{pure}),1),
\Phi(\oooo{\beta})[(\uuuu{\gamma}^{st};\id)],\uuuu{e})
\quad(\textup{by part (b)})\\
&=&
V_\Z^{(k)}((b_n^{pure},1),[(\uuuu{\gamma}^{st};\id)].\beta,
\uuuu{e}) \quad(\textup{by Remark \ref{t8.17} (ii)})\\
&=&
V_\Z^{(k)}((b_n^{pure},1),\beta[(\uuuu{\gamma}^{st};\id)],
\uuuu{e})\quad(\textup{by Lemma \ref{t8.18}})\\
&=&
V_\Z^{(k)}((b_n^{pure},1),[(\uuuu{\gamma}^{st};\id)],
\beta^{-1}(\uuuu{e}))\quad(\textup{by part (a)}).
\end{eqnarray*}

(d) For $\beta\in (\Br_n)_{\uuuu{e}/\{\pm 1\}^n }$
$\beta^{-1}(\uuuu{e})$ coincides with $\uuuu{e}$ up to signs,
$\beta^{-1}(\uuuu{e}/\{\pm 1\}^n)=\uuuu{e}/\{\pm 1\}^n$.
Because of $s^{(k)}_{-e_i}=s^{(k)}_{e_i}$ and part (c),
the restriction of $V_\Z^{(k),univ}$ to
$(\C-\{i,2i,...,ni\})\times\{\www{\beta}(1)\}$
is canonically isomorphicm to the restriction of
$V_\Z^{(k),univ}$ to 
$(\C-\{i,2i,...,ni\})\times \{b^{univ}\}$.
Therefore $\beta^{deck}$ extends to an automorphism
of $V_\Z^{(k),univ}$ as claimed.
The quotient bundle $V_\Z^{(k),\uuuu{e}/\{\pm 1\}^n}$
is well defined.\hfill$\Box$

\bigskip
Part (a) of the next lemma says that in all 
$\Z$-lattice structures and $\Z$-lattice bundles which are
constructed from $(H_\Z,L,\uuuu{e})$, the bilinear form $I^{(k)}$
becomes a flat bilinear form on the bundle. 
Part (b) is a nice observation of K. Saito \cite[Lemma 2]{Sa82}.
It says that $I^{(k)}$ is essentially the unique flat bilinear
form on these bundles.

\begin{lemma}\label{t8.24} 
Let $k\in\{0,1\}$.

(a) Each of the following $\Z$-lattice bundles comes with
a $(-1)^k$-symmetric bilinear form on its fibers, which is
induced by $I^{(k)}$ on $H_\Z$ and which is also called
$I^{(k)}$: \\
the bundle 
$V_\Z^{(k)}((\uuuu{u},r),[(\uuuu{\gamma};\sigma)],
\uuuu{f})$ over $\C-\{u_1,...,u_n\}$ 
in Definition \ref{t8.22} (b),\\
the bundle
$V_\Z^{(k),univ}$ over $\C\times C_n^{univ}-D_{1,n}^{univ}$
in Definition \ref{t8.22} (c),\\
the bundle $V_\Z^{(k),\uuuu{e}/\{\pm 1\}^n}$ over 
$\C\times C_n^{\uuuu{e}/\{\pm 1\}^n}
-D_{1,n}^{\uuuu{e}/\{\pm 1\}^n}$ in Theorem \ref{t8.23} (d).

(b) \cite[Lemma 2]{Sa82} Suppose that $(H_\Z,L,\uuuu{e})$ is
irreducible and $(n,k)\neq (1,1)$. In each of the  bundles
in part (a), the only flat bilinear forms on its fibers
are the bilinear bilinear forms $a\cdot I^{(k)}$
with $a\in\Z$, so the scalar multiples of $I^{(k)}$.
\end{lemma}

{\bf Proof:}
(a) This follows in all cases from $s^{(k)}_{f_i}\in\OO^{(k)}$,
which says that $s^{(k)}_{f_i}:H_\Z\to H_\Z$ respects 
$I^{(k)}$. 

(b) In each of the bundles in part (a), the monodromy group,
which is the image of the group antihomomorphism
\begin{eqnarray*}
\pi_1(\textup{base space},\textup{base point})\to O^{(k)}
\end{eqnarray*}
is $\Gamma^{(k)}$. 
This is clear for the bundle 
$V_\Z^{(k)}((\uuuu{u},r),[(\uuuu{\gamma};\sigma)],\uuuu{f})$
in Definition \ref{t8.22} (b) and for $V_\Z^{(k),univ}$.
It follows for the bundle $V_\Z^{\uuuu{e}/\{\pm 1\}^n}$
with the fact that the restriction of this bundle to
$\{(z,b)\in\C\times C_n^{\uuuu{e}/\{\pm 1\}^n}\,|\, 
\Ree z>\max(\Ree u\,|\, (u,b)\in D_{1,n}^{\uuuu{e}/\{\pm 1\}^n})
\}$ is trivial
(without this fact, one might have additional
{\it transversal} monodromy).

The proof of Lemma 2 in \cite{Sa82} shows that the only
$\Gamma^{(k)}$ invariant bilinear forms on $H_\Z$ are the
forms $a\cdot I^{(k)}$ with $a\in \Z$.
\hfill$\Box$

\begin{example}\label{t8.25}
The case $A_2$, $(H_\Z,L,\uuuu{e})$ of rank 2 with the matrix
$S=L(\uuuu{e}^t,\uuuu{e})^t
=\begin{pmatrix}1&-1\\0&1\end{pmatrix}$.
We saw in Theorem \ref{t8.12} 
\begin{eqnarray*}
C_2^{\uuuu{e}/\{\pm 1\}^2}=C_2^{A_2}\cong \C\times\C^*
\quad\textup{with coordinates }(z_1,z_2)\\
\textup{with }u_{1/2}=z_1\pm \frac{2}{3}z_2^{3/2},\quad
z_1=\frac{1}{2}(u_1+u_2),
z_2=(\frac{3}{4}(u_1-u_2))^{2/3}.
\end{eqnarray*}
The base point $b_2^{\uuuu{e}/\{\pm 1\}^2}\in 
C_2^{\uuuu{e}/\{\pm 1\}^2}$ is a point with $(u_1,u_2)=(i,2i)$.
We choose $b_2^{\uuuu{e}/\{\pm 1\}^2}$ with 
$(z_1^0,z_2^0):=(\frac{3}{2}i,
(\frac{3}{4})^{2/3}e^{2\pi i/6})$ in the
coordinates $(z_1,z_2)$. Recall
\begin{eqnarray*}
(e_1,e_2)\stackrel{\sigma_1}{\longmapsto}
(e_1+e_2,e_1)\stackrel{\sigma_1}{\longmapsto}
(-e_2,e_1+e_2)\stackrel{\sigma_1}{\longmapsto}
(e_1,-e_2).
\end{eqnarray*}
The set $\BB^{dist}/\{\pm 1\}^2$ has the three elements
\begin{eqnarray*}
x_0&=&(e_1,e_2)/\{\pm 1\}^2=\sigma_1^{-1}x_2,\\
x_1&=&(-e_2,e_1+e_2)/\{\pm 1\}^2 =\sigma_1^{-1}x_0,\\
x_2&=&(e_1+e_2,e_1)/\{\pm 1\}^2 = \sigma_1^{-1}x_1.
\end{eqnarray*}
The walls in $C_2^{\uuuu{e}/\{\pm 1\}^2}$ are the
real hypersurfaces of points $(z_1,z_2)$ with 
$\Imm u_1=\Imm u_2$. They are in the coordinates $(z_1,z_2)$
\begin{eqnarray*}
\C\times\{z_2\in\C^*\,|\, \Imm z_2^{3/2}=0\}
=\C\times (\R_{>0}\, \dot\cup\,  \R_{>0}\cdot e^{2\pi i /3}
\, \dot\cup\, \R_{>0}\cdot e^{2\pi i 2/3}).
\end{eqnarray*}
Between them are the three Stokes regions 
$F_2^{x_0,x_0},F_2^{x_0,x_1},F_2^{x_0,x_2}$.

The picture in Figure \ref{Fig:8.9} shows in the middle a part of
the slice $\{z_1^0\}\times \C \subset 
C_2^{\uuuu{e}/\{\pm 1\}^2}$ and the intersection of this part
with the three walls and the three Stokes regions.
The three lines with arrows are the three lifts of 
$\sigma_1\in \pi_1(C_2^{conf},b_2^{conf})$ to paths in
$C_2^{\uuuu{e}/\{\pm 1\}^2}$.
In the outer part, the picture shows for different values of
$(z_1^0,z_2)$ standard systems of paths
(recall $u_{1/2}=z_1\pm \frac{2}{3}z_2^{3/2}$).
The small circle between the points $u_1$ and $u_2$ means the
value 0. 

We have chosen the indexing of $u_1$ and $u_2$ in each Stokes
region so that $\Imm u_1<\Imm u_2$. 
Therefore there is a discontinuity of indexing along each wall.

\begin{figure}
\includegraphics[width=0.9\textwidth]{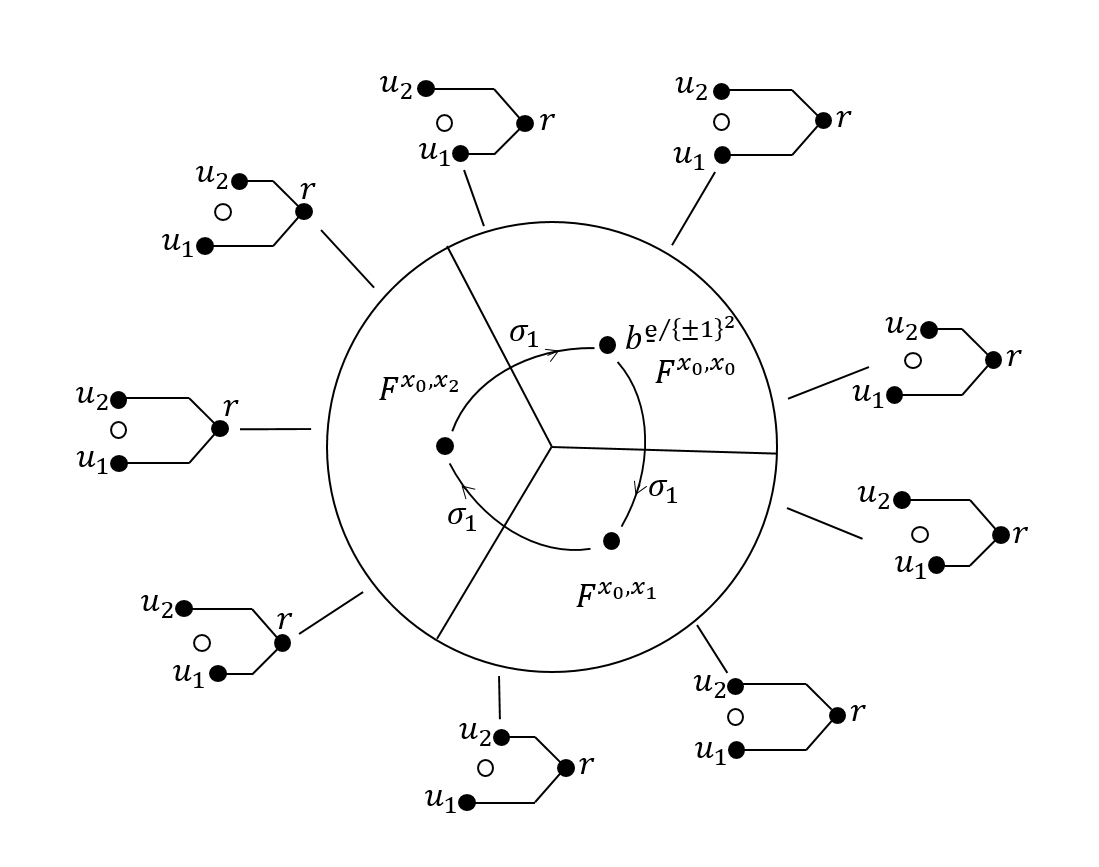}
\caption[Figure 8.9]{Moving with powers of $\sigma_1$ through the 
three Stokes regions of $C_2^{A_2}$}
\label{Fig:8.9}
\end{figure}

\end{example}

\begin{remarks}\label{t8.26}
(i) Part (c) and part (a) of Theorem \ref{t8.23} give the
following.
Let $k\in\{0,1\}$ and let $\beta\in \Br_n$. 
Let $\www{\beta}:[0,1]\to C_n^{\uuuu{e}/\{\pm 1\}^n}$
be the lift to $C_n^{\uuuu{e}/\{\pm 1\}^n}$ with 
$\www{\beta}(0)=b_n^{\uuuu{e}/\{\pm 1\}^n}$ of a loop in
$C_n^{conf}$ which represents $\beta$. 
Then the restriction of $V_\Z^{(k),\uuuu{e}/\{\pm 1\}^n}$
to $(\C-\{i,2i,...,ni\})\times \{\www{\beta}(1)\}$
is equal to $V_\Z^{(k)}((b_n^{pure},1),
[(\uuuu{\gamma}^{st},\id)],\beta^{-1}(\uuuu{e}))$.

(ii) This fits well to the indexing of the Stokes regions
in Lemma \ref{t8.6} (d) in the case
$X=\BB^{dist}/\{\pm 1\}^n$.
Then $x_0:=\uuuu{e}/\{\pm 1\}^n$.
The Stokes region which contains $\www{\beta}(1)$ is
$F_n^{x_0,x}$ with $x=\beta^{-1}(\uuuu{e})/\{\pm 1\}^n\in X$.
See also Remark \ref{t8.7} (ii). 

(iii) One can modify the construction of the bundle
$V_\Z^{(k),\uuuu{e}/\{\pm 1\}^n}$ in Theorem \ref{t8.23} (d).
Suppose that a group antihomomorphism
\begin{eqnarray*}
\pi_1(C_n^{\uuuu{e}/\{\pm 1\}^n},b_n^{\uuuu{e}/\{\pm 1\}^n})
\to (\{\pm 1\}^n)_S\subset G_\Z
\end{eqnarray*}
is given. Here an element of $\{\pm 1\}^n$ means an automorphism
of $H_\Z$ which changes some signs in the basis $\uuuu{e}$.
Then $(\{\pm 1\}^n)_S\subset G_\Z$.

One can twist the automorphism $\beta^{deck,V,(k)}$
in Theorem \ref{t8.23} (d) with the correct element of
$(\{\pm 1\}^n)_S$ and obtain a twisted quotient bundle
over $\C\times C_n^{\uuuu{e}/\{\pm 1\}^n}-
D_{1,n}^{\uuuu{e}/\{\pm 1\}^n}$. 

But one looses the property of $V_\Z^{(k),\uuuu{e}/\{\pm 1\}^n}$
to be a trivial bundle over the open set
\begin{eqnarray*}
\{(z,b)\in\C\times C_n^{\uuuu{e}/\{\pm 1\}^n}\,|\, 
\Ree z>\max(\Ree u\,|\, (u,b)\in 
D_{1,n}^{\uuuu{e}/\{\pm 1\}^n})\}.
\end{eqnarray*}
This property is needed for a possible extension of such a 
bundle over a partial compactification of 
$C_n^{\uuuu{e}/\{\pm 1\}^n}$. 
Therefore the bundle $V_\Z^{(k),\uuuu{e}/\{\pm 1\}^n}$
is more natural than the twisted bundles.

(iv) Consider a bundle $V_\Z^{(k)}((\uuuu{u},r),
[(\uuuu{\gamma};\sigma)],\uuuu{f})$ as in Definition \ref{t8.22}
(a). In most cases one can recover $\uuuu{f}/\{\pm 1\}^n$
from the bundle and from $[(\uuuu{\gamma};\sigma)]$, 
namely in all cases with $k=0$ and in all cases with $k=1$ 
where is 
$(H_\Z,L,\uuuu{e})$ is not reducible with at least two 
summands of type $A_1$, so especially in all irreducible cases.

We explain this. In any case, one recovers from the bundle
and the distinguished system of paths $(\uuuu{\gamma};\sigma)$
the local monodromies $s^{(k)}_{f_1},...,s^{(k)}_{f_n}$
along the associated loops.
With Lemma \ref{t3.15} (b) and (c) one obtains from the tuple
$(s_{f_1}^{(k)},...,s_{f_n}^{(k)})$ of reflections or 
transvections the distinguished basis $\uuuu{f}/\{\pm 1\}^n$
up to signs. 
\end{remarks}

\chapter[The manifolds in the rank 3 cases]
{The manifolds in the rank 3 cases}\label{s9}
\setcounter{equation}{0}
\setcounter{figure}{0}

The manifolds $C_n^{\uuuu{e}/\{\pm 1\}^n}$ and 
$C_n^{S/\{\pm 1\}^n}$ from chapter \ref{s8} had been
discussed in the rank $2$ cases in section \ref{s8.3}.
There only three different manifolds 
$C_2^{\uuuu{e}/\{\pm 1\}^n}$ arose, namely $C_2^{pure}$, 
$C_2^{A_2}$ and $C_2^{univ}$, and only one manifold
$C_2^{S/\{\pm 1\}^2}=C_2^{conf}$.
In rank 3 the classification is already much richer.
The rank 3 cases are treated in this chapter.

Section \ref{s9.1} makes the deck transformations of
$C_3^{univ}$ as universal covering of $C_3^{conf}$ and 
of $C_3^{pure}$ explicit. 
For that it introduces new coordinates on $C_3^{pure}$. 
In fact, the new coordinates come from the restriction to
$C_3^{pure}$ of a blowing up of the complex line
$\{u_1=u_2=u_3\}$ in $\C^3\supset C_3^{pure}$. 
For the deck transformations the Schwarzian triangle function
$T:\H\to \C-\{0,1\}\cong \H/\oooo{\Gamma(2)}$ and a 
lift $\kappa:\H\to\C$ with $T$ of the logarithm
$\ln:\C-\{0\}\dashrightarrow \C$ are needed. Both are treated
in Appendix \ref{sc}. 

Theorem \ref{t9.3} in section \ref{s9.2} gives the main result,
a description of all possible manifolds
$C_3^{\uuuu{e}/\{\pm 1\}^3}$ and $C_3^{S/\{\pm 1\}^3}$
in rank 3. 
It starts with a unimodular bilinear lattice $(H_\Z,L,\uuuu{e})$
with a triangular basis $\uuuu{e}$ and matrix
$S=S(\uuuu{x})=L(\uuuu{e}^t,\uuuu{e})^t\in T^{uni}_3(\Z)$
for some $\uuuu{x}\in\Z^3$. 
Thanks to Theorem \ref{t4.14} we can and will restrict to 
a local minimum $\uuuu{x}\in C_i$ for some $i\in\{1,2,...,24\}$.
Theorem \ref{t7.11} lists all 16 possible pairs
$((\Br_3)_{\uuuu{e}/\{\pm 1\}^3},
(\Br_3)_{S/\{\pm 1\}^3})$ of stabilizers, the case
$S(-l,2,-l)$ with parameter $l\in \Z_{\geq 3}$. 

Theorem \ref{t9.3} gives in all 16 cases the manifolds
$C_3^{\uuuu{e}/\{\pm 1\}^3}$ and $C_3^{S/\{\pm 1\}^3}$.
Though this is coarse information.
The covering map of the normal covering
$C_3^{\uuuu{e}/\{\pm 1\}^3}\to C_3^{S/\{\pm 1\}^3}$
is also important.
But its description requires the whole discussion in the proof,
and therefore it is given only in the proof of Theorem
\ref{t9.3}. 

Also the F-manifold structure and the Euler field are important,
as well as possible partial compactifications of 
$C_3^{\uuuu{e}/\{\pm 1\}^3}$ such that the F-manifold structure
extends. They are treated in section \ref{s9.3} 
in Lemma \ref{t9.7} and Corollary \ref{t9.8}.
The cases $A_3$ and $A_2A_1$ have partial compactifications
which are well known from singularity theory. Seeing them
here is a bit involved, especially in the case $A_3$.

\section{The deck transformations on $C_3^{univ}$}\label{s9.1}

Theorem \ref{t9.1} first gives new coordinates on
$C_3^{pure}$. They are better suited to treat
the deck transformations of $C_3^{univ}$ as universal covering 
$\pr_3^{u,c}:C_3^{univ}\to C_3^{conf}$ and as universal covering
$\pr_3^{u,p}:C_3^{univ}\to C_3^{pure}$. 
It writes the F-manifold structure in these new coordinates.
It gives in these coordinates explicitly the six deck 
transformations of $C_3^{pure}$ as normal covering of
$C_3^{conf}$. 
Finally, it gives explicitly several deck transformations
of $C_3^{univ}$ as universal covering of $C_3^{conf}$. 
Here the Schwarzian triangle function
$T:\H\to\C-\{0,1\}$ and the lift $\kappa:\H\to\C$ with $T$
of the logarithm $\ln:\C-\{0\}\dashrightarrow\C$ 
from Appendix \ref{sc} are used.

\begin{theorem}\label{t9.1}
(a) The map
\begin{eqnarray*}
f:C_3^{pure}&\to& \C\times \C^*\times (\C-\{0,1\}),\\
(u_1,u_2,u_3) &\mapsto&  
(z_1,z_2,z_3)=(u_1+u_2+u_3,u_2-u_1,\frac{u_3-u_1}{u_2-u_1}),\\
\end{eqnarray*}
is an isomorphism of complex manifolds. The inverse map is
\begin{eqnarray*}
f^{-1}:\C\times\C^*\times (\C-\{0,1\}) \to C_3^{pure}\\
(z_1,z_2,z_3) \mapsto  (u_1,u_2,u_3)= 
\frac{1}{3}(z_1-z_2(1+z_3),\\
z_1+z_2(2-z_3),z_1+z_2(-1+2z_3)).
\end{eqnarray*}
Let $\paa_j:=\frac{\paa}{\paa z_j}$ be the coordinate vector
fields of the coordinates $(z_1,z_2,z_3)$. The base changes
between them and the partial units 
$e_i=\frac{\paa}{\paa u_i}$ are as follows,
\begin{eqnarray*}
(\paa_1,\paa_2,\paa_3)&=&(e_1,e_2,e_3)\frac{1}{3}
\begin{pmatrix}1&-1-z_3&-z_2\\1&2-z_3&-z_2\\1&-1+2z_3&2z_2
\end{pmatrix},\\
(e_1,e_2,e_3)&=&(\paa_1,\paa_2,\paa_3)
\begin{pmatrix}1&1&1\\-1&1&0\\ \frac{z_3-1}{z_2} & 
\frac{-z_3}{z_2} & \frac{1}{z_2}\end{pmatrix}.
\end{eqnarray*}
Unit field $e$, Euler field $E$ and multiplication $\circ$ 
are in canonical coordinates given by 
\begin{eqnarray*}
e=\sum_{j=1}^3 e_j,\quad E=\sum_{j=1}^3 u_je_j,\quad
e_i\circ e_j=\delta_{ij}e_i.
\end{eqnarray*}
In the coordinates $(z_1,z_2,z_3)$ they are given by 
\begin{eqnarray*}
e&=&3\paa_1,\quad E=z_1\paa_1+z_2\paa_2,\\
\paa_2\circ\paa_2&=& \frac{2-2z_3+2z_3^2}{3}\paa_1 + 
\frac{1-2z_3}{3}\paa_2+\frac{-z_3+z_3^2}{z_2}\paa_3,\\
\paa_2\circ\paa_3&=& \frac{-z_3+2z_2z_3}{3}\paa_1 
+ \frac{-z_2}{3}\paa_2 + \frac{-1+2z_3}{3}\paa_3,\\
\paa_3\circ\paa_3&=& \frac{2z_2^2}{3}\paa_1 
+\frac{z_2}{3}\paa_3.
\end{eqnarray*}

(b) The map
\begin{eqnarray*}
\sigma\mapsto \chi_\sigma:=
(C_3^{pure}\to C_3^{pure},\ \uuuu{u}\mapsto
(u_{\sigma^{-1}(1)},u_{\sigma^{-1}(2)},u_{\sigma^{-1}(3)}))
\end{eqnarray*}
is an isomorphism from $S_3$ to the group of deck transformations
of the covering $C_3^{pure}\to C_3^{conf}$.
The corresponding automorphisms
$\phi_\sigma:=f\circ \chi_\sigma \circ f^{-1}$ of 
$\C\times\C^*\times (\C-\{0,1\})$ are as follows
(see Theorem \ref{tb.1} (d) for $g_\sigma$):
\begin{eqnarray*}
\phi_{\id}(\uuuu{z})&=&\uuuu{z},\\
\phi_{(12)}(\uuuu{z})&=& (z_1,-z_2,g_{(12)}(z_3)),\\
\phi_{(13)}(\uuuu{z})&=& (z_1,z_2(1-z_3),g_{(13)}(z_3)),\\
\phi_{(23)}(\uuuu{z})&=& (z_1,z_2z_3,g_{(23)}(z_3)),\\
\phi_{(123)}(\uuuu{z})&=& (z_1,-z_2z_3,g_{(123)}(z_3)),\\
\phi_{(132)}(\uuuu{z})&=& (z_1,z_2(z_3-1),g_{(132)}(z_3)).
\end{eqnarray*}

(c) A universal covering of $\C\times\C^*\times (\C-\{0,1\})$
is given by the map 
\begin{eqnarray*}
\www{\pr}_3^{u,p}:\C\times \C\times \H&\to& \C\times\C^*\times
(\C-\{0,1\})\\
(z_1,\zeta,\tau)&\mapsto& (z_1,e^\zeta,T(\tau)),
\end{eqnarray*}
with $T$ as in Theorem \ref{tb.1} (d).
The group of deck transformations of the universal covering
\begin{eqnarray*}
\pr_3^{p,c}\circ f^{-1}\circ \www{\pr}_3^{u,p}:
\C\times\C\times\H
&\longrightarrow& \C\times \C^*\times (\C-\{0,1\})\\ 
&\stackrel{\cong}{\longrightarrow}& C_3^{pure}\\
&\stackrel{/S_3}{\longrightarrow}& C_3^{conf}
\end{eqnarray*}
is isomorphic to the group $\Br_3$.
An isomorphism $\psi$ 
\index{$\psi,\ \psi(\sigma_1),\ \psi(\sigma_2),\ \psi(\sigma^{mon})$} 
is determined by the images 
$\psi(\sigma_1)$ and $\psi(\sigma_2)$ of the elementary
braids $\sigma_1$ and $\sigma_2$. The following work
(recall the lift $\kappa:\H\to \C$ of the (multivalued) 
logarithm $\log:\C^*\dashrightarrow\C$ 
in Definition \ref{tb.2} and Lemma \ref{tb.3} (c)), 
\begin{eqnarray*}
\psi(\sigma_1):\C\times\C\times\H&\to& \C\times\C\times\H,\\ 
(z_1,\zeta,\tau)&\mapsto& (z_1,\zeta-\pi i,\tau-1)
\quad(\textup{recall }\tau-1=\mu_{(12)}(\tau)),\\
\psi(\sigma_2):\C\times\C\times\H&\to& \C\times\C\times\H,\\ 
(z_1,\zeta,\tau)&\mapsto& (z_1,\zeta+\kappa(\tau+1),\
\frac{\tau}{\tau+1})\\
&&\qquad (\textup{recall }\frac{\tau}{\tau+1}=\mu_{(13)}(\tau)).
\end{eqnarray*}
The isomorphism $\psi$ satisfies 
\begin{eqnarray*}
\psi(\sigma^{mon}):(z_1,\zeta,\tau)&\mapsto& 
(z_1,\zeta-2\pi i,\tau),\\
\psi((\sigma^{mon})^{-1}\sigma_1^2):(z_1,\zeta,\tau)&\mapsto&
(z_1,\zeta,\tau-2),\\
\psi(\sigma_2^2):(z_1,\zeta,\tau)&\mapsto&
(z_1,\zeta,\frac{\tau}{2\tau+1}),\\
\psi(\sigma_1\sigma_2^2):(z_1,\zeta,\tau)&\mapsto&
(z_1,\zeta-\pi i,\frac{-\tau-1}{2\tau+1}),\\
\psi(\sigma_2\sigma_1):(z_1,\zeta,\tau)&\mapsto&
(z_1,\zeta-\pi i+\kappa(\tau),\frac{\tau-1}{\tau}).
\end{eqnarray*}
Here recall
\begin{eqnarray*}
\mu_{(12)}^2=(\tau\to \tau-2),\quad 
\mu_{(13)}^2=(\tau\mapsto\frac{\tau}{2\tau+1}),\\
\mu_{(123)}=(\tau\mapsto \frac{\tau-1}{\tau})
=\textup{rotation of order 3 around }e^{2\pi i/6},\\
\mu(\begin{pmatrix}-1&-1\\2&1\end{pmatrix})
=\textup{rotation of order 2 around }\frac{-1+i}{2}.
\end{eqnarray*}
\end{theorem}

{\bf Proof:}
(a) It is rather obvious that $f$ is an isomorphism of complex
manifolds with inverse $f^{-1}$ as claimed.
The base changes between the two bases $(\paa_1,\paa_2,\paa_3)$
and $(e_1,e_2,e_3)$ of coordinate vector fields are obtained
by derivating $f(\uuuu{u})$ and $f^{-1}(\uuuu{z})$ suitably.
Unit field $e$, Euler field $E$ and multiplication $\circ$
in the new coordinates $\uuuu{z}=(z_1,z_2,z_3)$ can be 
calculated straightforwardly.

(b) Straightforward calculations.

(c) $\psi(\sigma_1)$ and $\psi(\sigma_2))$ are lifts to
$\C\times\C\times\H$ of the automorphisms 
$\phi_{(12)}$ and $\phi_{(13)}$ of 
$\C\times\C^*\times(\C-\{0,1\})$ because
\begin{eqnarray*}
\www{\pr}_3^{u,p}(\psi(\sigma_1)(z_1,\zeta,\tau)))
&=& \www{\pr}_3^{u,p}(z_1,\zeta-\pi i,\mu_{(12)}(\tau))\\
&=&(z_1,e^{\zeta-\pi i},T(\mu_{(12)}(\tau)))\\
&=& (z_1,-e^{\zeta},g_{(12)}(T(\tau)))\\
&=&(z_1,-z_2,g_{(12)}(z_3)),\\
\www{\pr}_3^{u,p}(\psi(\sigma_2)(z_1,\zeta,\tau)))
&=& \www{\pr}_3^{u,p}(z_1,\zeta+\kappa(\tau+1),
\mu_{(13)}(\tau))\\
&=&(z_1,e^{\zeta}T(\tau+1),T(\mu_{(13)}(\tau)))\\
&=& (z_1,e^{\zeta}(1-T(\tau)),g_{(13)}(T(\tau)))\\
&=&(z_1,z_2(1-z_3),g_{(13)}(z_3)).
\end{eqnarray*}
They satisfy the relation
\begin{eqnarray*}
\psi(\sigma_1)\psi(\sigma_2)\psi(\sigma_1)
=\psi(\sigma_2)\psi(\sigma_1)\psi(\sigma_2)
\end{eqnarray*}
because
\begin{eqnarray*}
\psi(\sigma_1)\psi(\sigma_2)\psi(\sigma_1):(z_1,\zeta,\tau)
&\stackrel{\sigma_1}{\longmapsto}& 
(z_1,\zeta-\pi i,\tau-1)\\
&\stackrel{\sigma_2}{\longmapsto}&
(z_1,\zeta-\pi i+\kappa(\tau),\frac{\tau-1}{\tau})\\
&\stackrel{\sigma_1}{\longmapsto}&
(z_1,\zeta-2\pi i+\kappa(\tau),\frac{-1}{\tau}),\\
\psi(\sigma_2)\psi(\sigma_1)\psi(\sigma_2):(z_1,\zeta,\tau)
&\stackrel{\sigma_2}{\longmapsto}& 
(z_1,\zeta+\kappa(\tau+1)),\frac{\tau}{\tau+1})\\
&\stackrel{\sigma_1}{\longmapsto}&
(z_1,\zeta-\pi i+\kappa(\tau+1),\frac{-1}{\tau+1})\\
&\stackrel{\sigma_2}{\longmapsto}&
(z_1,\zeta-\pi i+\kappa(\tau+1)+\kappa(\frac{\tau}{\tau+1}),
\frac{-1}{\tau})
\end{eqnarray*}
and $-2\pi i +\kappa(\tau) =-\pi i +\kappa(\tau+1) 
+\kappa(\frac{\tau}{\tau+1})$ by Theorem \ref{tb.3} (c) (iv).

Therefore $\psi$ is a group homomorphism from $\Br_3$ to
the group of deck transformations of the universal covering
$\C\times\C\times \H\to C_3^{conf}$. It is surjective because
$\psi(\sigma^{mon})$ is 
\begin{eqnarray*}
\psi(\sigma^{mon}) 
&=& (\psi(\sigma_1)\psi(\sigma_2)\psi(\sigma_1))^2:\\
(z_1,\zeta,\tau)&\mapsto& 
(z_1,\zeta-4\pi i+\kappa(\tau)+\kappa(\frac{-1}{\tau}),\tau)\\
&&= (z_1,\zeta-2\pi i,\tau)
\end{eqnarray*}
because
\begin{eqnarray*}
\kappa(\tau)+\kappa(\frac{-1}{\tau})
&=& \kappa(\tau)+\kappa(\frac{-1}{\tau}+2)+2\pi i\\
&=& \kappa(\tau)+\kappa(\frac{2\tau-1}{\tau})+2\pi i\\
&=& 2\pi i \quad\textup{ by Theorem \ref{tb.3} (c) (iv).}
\end{eqnarray*}
$\psi(\sigma_2^2)$ is as claimed because by Theorem \ref{tb.3}
(c) (iv) 
$$\zeta + \kappa(\tau+1)+\kappa(\frac{2\tau+1}{\tau+1})
=\zeta + \kappa(\tau+1)+\kappa(\frac{2(\tau+1)-1}{\tau+1})
=\zeta.$$
The formulas for $\psi((\sigma^{mon})^{-1}\sigma_1^2)$,
$\psi(\sigma_1\sigma_2^2)$ and $\psi(\sigma_2\sigma_1)$ 
are calculated straightforwardly.
\hfill$\Box$

\begin{remarks}\label{t9.2}
(i) The group of deck transformations of the universal covering
$\C\times\C\times\H\to C_3^{pure}$ splits into the product
\begin{eqnarray*}
\langle \psi(\sigma^{mon})\rangle\times 
\langle \psi((\sigma^{mon})^{-1}\sigma_1^2,\psi(\sigma_2^2)
\rangle \\
=\langle ((z_1,\zeta,\tau)\mapsto (z_1,\zeta-2\pi i,\tau))
\rangle \times \\
\langle((z_1,\zeta,\tau)\mapsto(z_1,\zeta,\mu_{(12)}^2(\tau)),
(z_1,\zeta,\tau)\mapsto(z_1,\zeta,\mu_{(13)}^2(\tau))\rangle.
\end{eqnarray*}
This fits to the splitting of $\Br_3^{pure}$,
\begin{eqnarray*}
\Br_3^{pure}=\langle\sigma^{mon}\rangle \times
\langle (\sigma^{mon})^{-1}\sigma_1^2,\sigma_2^2\rangle.
\end{eqnarray*}

(ii) The coordinate change 
$$\C^3\to\C^3,\quad  (u_1,u_2,u_3)\mapsto 
(z_1,z_2,z_4):=(u_1+u_2+u_3,u_2-u_1,u_3-u_1),$$
is linear. $f$ is the restriction to $C_3^{pure}$ of the 
composition of this coordinate change with the blowing up 
\index{blowing up} 
$$\C^3 \dashrightarrow \C\times \OO_{\P^{1}}(-1),\quad
(z_1,z_2,z_4)\dashrightarrow (z_1,z_2,\frac{z_4}{z_2})
=(z_1,z_2,z_3),$$
of $\C\times\{(0,0)\}$ in one of the two natural charts of 
$\C\times \OO_{\P^1}(-1)$. 

\begin{figure}[H]
\includegraphics[width=1.0\textwidth]{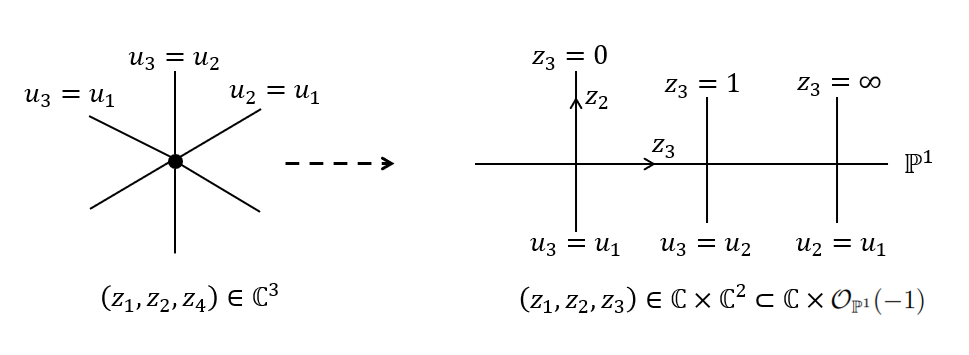}
\caption[Figure 9.1]{Blowing up of $\C\times\{(0,0)\}
\subset \C^3$}
\label{Fig:9.1}
\end{figure}

Also the multiplication is calculated only in this chart $\C^3$
with coordinates $(z_1,z_2,z_3)$. 
We refrain from calculating it in the 
other chart $\C^3$ with coordinates 
$(z_1,\www{z}_2,\www{z}_3)=(z_1,\frac{z_2}{z_4},z_4)$
as we do not really need it. 
Though in the chart considered, the hyperplane in the Maxwell
stratum where $u_2=u_1$ lies above $z_3=\infty$, so for that
hyperplane the other chart would be useful.
The hyperplanes where $u_3=u_1$ respectively $u_3=u_2$
lie above $z_3=0$ respectively $z_3=1$.

(iii) The multiplication on $\C\times\C\times\H\cong C_3^{univ}$
involves $\frac{\paa T(\tau)}{\paa\tau}$. 
Unit field and Euler field on it are $e=3\paa_1$ and 
$E=z_1\paa_1+\paa_{\zeta}$. 
\end{remarks}

\section{The manifolds as quotients and their covering maps}
\label{s9.2}

Let $(H_\Z,L,\uuuu{e})$ be a unimodular bilinear lattice
of rank 3 with triangular basis $\uuuu{e}$ and matrix
$S=S(\uuuu{x})=L(\uuuu{e}^t,\uuuu{e})^t\in T^{uni}_3(\Z)$
for some $\uuuu{x}\in\Z^3$.
We can and will restrict to a local minimum $\uuuu{x}\in C_i$
for some $i\in\{1,2,...,24\}$ as in Theorem \ref{t7.11}.
This theorem lists all 16 possible pairs
$((\Br_3)_{\uuuu{e}/\{\pm 1\}^3},
(\Br_3)_{S/\{\pm 1\}^3})$ of stabilizers, the case
$S(-l,2,-l)$ with parameter $l\in \Z_{\geq 3}$. 

Theorem \ref{t9.3} below gives in all 16 cases the manifolds
$C_3^{\uuuu{e}/\{\pm 1\}^3}$ and $C_3^{S/\{\pm 1\}^3}$.
Though this is coarse information.
The covering map of the normal covering
$C_3^{\uuuu{e}/\{\pm 1\}^3}\to C_3^{S/\{\pm 1\}^3}$
is also important.
But its description requires the whole discussion in the proof,
and therefore it is given only in the proof of Theorem
\ref{t9.3}. 

The group $(\Br_3)_{\uuuu{e}/\{\pm 1\}^3}$ acts on the third
factor $\H$ of $\C\times\C\times\H\cong C_3^{univ}$ via
its image in $PSL_2(\Z)$ under the homomorphism 
$\Br_3\to PSL_2(\Z)$ in Remark \ref{t4.15} (i). 
In all cases except $A_3$ and $A_2A_1$ this action is free
and the quotient of $\H$ by this action 
is a noncompact complex curve $C$. 
Furthermore, in all cases except $A_1^3$, $\HH_{1,2}$, $A_2A_1$, 
and $A_3$ the group $(\Br_3)_{\uuuu{e}/\{\pm 1\}^3}$ does not
contain a power of $\sigma^{mon}$. Then the quotient 
$\C\times\C\times\H/(\Br_3)_{\uuuu{e}/\{\pm 1\}^3}\cong
C_3^{\uuuu{e}/\{\pm 1\}^3}$ is 
$$\C\times (\textup{an affine }\C\textup{-bundle over }C).$$
Remark \ref{t9.4} (i) says what that means, and Remark
\ref{t9.4} (ii) states that it is isomorphic to
$\C\times\C\times C$. Remark \ref{t9.4} (iv) states
that any $\C^*$-bundle over a noncompact complex curve $C$ is 
isomorphic to $\C^*\times C$.

\begin{theorem}\label{t9.3}
Consider a local minimum $\uuuu{x}\in C_i\subset \Z^3$ for some 
$i\in\{1,2,...,24\}$ and the pseudo-graph $\GG_j$ with
$\GG_j=\GG(\uuuu{x})$. In the following table
the first and second column are copied from the table in
Theorem \ref{t7.11}.
The third and fourth column give the manifolds
$C_3^{\uuuu{e}/\{\pm 1\}^3}$ and 
$C_3^{S/\{\pm 1\}^3}$.
Here $\D:=\{z\in\C\,|\, |z|<1\}$ denotes the unit disk,
\index{$\D,\ \D^*,\ \D^{**}$}
$\D^*:=\D-\{0\}$, $\D^{**}:=\D-\{0,\frac{1}{2}\}$
and $\C^{***}:=\C-\{z\in\C\,|\, z^3=1\}$.
\index{$\C^{***}$}
\begin{eqnarray*}
\begin{array}{l|l|l|l}
 & \textup{sets} & 
C_3^{\uuuu{e}/\{\pm 1\}^3} 
& C_3^{S/\{\pm 1\}^3}\\ \hline 
\GG_1 & C_1\ (A_1^3)
& C_3^{pure} & 
C_3^{conf}\\
\GG_1 & C_2\ (\HH_{1,2})  
& \C\times\C^*\times\H & C_3^{conf} \\
\GG_2 & C_3\ (A_2A_1)
& \frac{\C\times\C^*\times (\C-\{0,1\})}
{\langle\phi_{(12)}\rangle} & 
C_3^{pure}/\langle\phi_{(12)}\rangle \\
\GG_2 & C_4\ (\P^1A_1),C_5
& \C^2\times (\C-\{0,1\}) & 
C_3^{pure}/\langle\phi_{(12)}\rangle \\
\GG_3 & C_6\ (A_3) 
& \C\times\frac{\C^*\times\C^{***}}
{\textup{group of order 3}} & 
\C\times\frac{\C^*\times\C^{***}}
{\textup{group of order 3}} \\
\GG_4 & C_7 \ (\whh{A}_2) 
& \C^2\times \C^{***} & 
\C\times\frac{\C^*\times \C^{***}}
{\textup{group of order 3}} \\
\GG_5 & C_8,\ C_9\ ((-l,2,-l))
&\C\times\C\times\D^*  & 
\C\times\C^*\times \D^* \\
\GG_6 & C_{10}\ (\P^2),\ C_{11},\ C_{12} 
& C_3^{univ} & \frac{\C\times\C^*\times\H }
{\textup{group of order 3}}\\
\GG_7 & C_{13}\ (\textup{e.g. }(4,4,8)) 
& C_3^{univ} & \frac{\C\times\C^*\times\H}
{\textup{group of order 2}} \\
\GG_8 & C_{14}\ (\textup{e.g. }(3,4,6)) 
& C_3^{univ} & \C\times\C^*\times \H \\
\GG_9 & C_{15},\ C_{16},\ C_{23},\ C_{24} 
& C_3^{univ} & \C\times\C^*\times\H \\
\GG_{10} & C_{17}\ (\textup{e.g. }(-2,-2,0)) 
& \C\times\C\times\D^* & \C\times\C^*\times\D^* \\
\GG_{11} & C_{18}\ (\textup{e.g. }(-3,-2,0)) 
& \C\times\C\times\D^* & \C\times\C^*\times\D^* \\
\GG_{12} & C_{19}\ (\textup{e.g. }(-2,-1,0)) 
& \C\times\C\times\D^{**} & 
\C\times\C^*\times\D^{**} \\
\GG_{13} & C_{20}\ (\textup{e.g. }(-2,-1,-1)) 
& \C\times\C\times\D^{**} & 
\C\times\C^*\times\D^{**} \\
\GG_{14} & C_{21},\ C_{22}
& \C\times\C\times \D^* & \C\times\C^*\times\D^* 
\end{array}
\end{eqnarray*}
In the case $\GG_2\,\&\, C_3\ (A_2A_1)$ the moduli spaces and the covering 
can also be presented as follows,
\begin{eqnarray*}
\begin{CD}
C_3^{\uuuu{e}/\{\pm 1\}^3}\cong 
@. \hspace*{0.3cm}\C\times \{(z_8,z_9)\in\C^*\times\C\,|\, z_8^3-4z_9^2\neq 0\}
@. \hspace*{0.5cm}(z_1,z_8,z_9)\\
@VV{3:1}V @VV{3:1}V @VV{3:1}V\\
C_3^{S/\{\pm 1\}^3}\cong 
@. \hspace*{0.3cm}\C\times\{(z_{10},z_9)\in\C^*\times\C\,|\, z_{10}-4z_9^2\neq 0\} 
@. \hspace*{0.5cm}(z_1,z_8^3,z_9)
\end{CD}
\end{eqnarray*}
In the case $\GG_3\,\&\, C_6\ (A_3)$ the moduli spaces and the covering
can also be presented as follows, 
\begin{eqnarray*}
\begin{CD}
C_3^{\uuuu{e}/\{\pm 1\}^3}\cong 
@. \hspace*{0.3cm}\C\times \{(z_{13},z_{14})\in\C^*\times\C\,|\, z_{14}^3-z_{13}^2\neq 0\}
@. \hspace*{0.5cm} (z_1,z_{13},z_{14})\\
@VV{4:1}V @VV{4:1}V @VV{4:1}V\\
C_3^{S/\{\pm 1\}^3}\cong 
@. \hspace*{0.3cm}\C\times\frac{\{(z_{15},z_{14})\in\C^*\times\C\,|\,z_{14}^3-z_{15}\neq 0\}}
{\langle (z_{15},z_{14})\mapsto (-z_{15},-z_{14})\rangle} 
@. \hspace*{0.5cm} [(z_1,z_{13}^2,z_{14})]
\end{CD}
\end{eqnarray*}
\end{theorem}

{\bf Proof:} 
{\bf The reducible case $\GG_1\,\&\, C_1\, (A_1)$:}\\
Here $(\Br_3)_{\uuuu{e}/\{\pm 1\}^3}= \Br_3^{pure}$
and $(\Br_3)_{\uuuu{x}/\{\pm 1\}^3}= \Br_3$.
Therefore $C_3^{\uuuu{e}/\{\pm 1\}^3}=C_3^{pure}$ and
$C_3^{\uuuu{x}/\{\pm 1\}^3}=C_3^{conf}$.
The normal covering $C_3^{pure}\to C_3^{conf}$ is known.

\medskip
{\bf The case $\GG_1\,\&\, C_2\, (\HH_{1,2})$:}\\
Here $(\Br_3)_{\uuuu{e}/\{\pm 1\}^3}= \langle(\sigma^{mon})^2
\rangle$ 
and $(\Br_3)_{\uuuu{x}/\{\pm 1\}^3}= \Br_3$.
Therefore $C_3^{\uuuu{x}/\{\pm 1\}^3}=C_3^{conf}$.
The equality 
$(\Br_3)_{\uuuu{e}/\{\pm 1\}^3}= \langle(\sigma^{mon})^2
\rangle$ and 
\begin{eqnarray*}
\psi((\sigma^{mon})^2): \C\times\C\times\H\to \C\times\C\times\H,
\ (z_1,\zeta,\tau)\mapsto (z_1,\zeta-4\pi i,\tau)
\end{eqnarray*}
show
\begin{eqnarray*}
C_3^{univ}\cong \C\times\C\times\H&\to& \C\times\C^*\times\H
\cong C_3^{\uuuu{e}/\{\pm 1\}^3},\\
(z_1,\zeta,\tau)&\mapsto& (z_1,e^{\zeta/2},\tau)=(z_1,
z_5,\tau).
\end{eqnarray*}
The quotient $(Br_3)_{S/\{\pm 1\}^3}/(\Br_3)_{\uuuu{e}/\{\pm 1\}^3} 
=\Br_3/\langle (\sigma^{mon})^2\rangle$ is by Remark
\ref{t4.15} (i) isomorphic to $SL_2(\Z)=\langle A_1,A_2
\rangle$ where $A_1=\begin{pmatrix}1&-1\\0&1\end{pmatrix}$,
$A_2=\begin{pmatrix}1&0\\1&1\end{pmatrix}$
and $\sigma_1,\sigma_2,\sigma^{mon}$ are mapped to 
$A_1,A_2,-E_2$. They induce the following deck transformations
of $\C\times\C^*\times\H \cong C_3^{\uuuu{e}/\{\pm 1\}^3}$,
\begin{eqnarray*}
\sigma_1:\ (z_1,z_5,\tau) &\mapsto & 
(z_1,(-i)z_5,\tau-1),\\
\sigma_2:\ (z_1,z_5,\tau) &\mapsto & 
(z_1,z_5\cdot e^{\kappa(\tau+1)/2},\frac{\tau}{\tau+1}),\\
\sigma^{mon}:\ (z_1,z_5,\tau) &\mapsto & 
(z_1,-z_5,\tau).
\end{eqnarray*}
Dividing out $\langle\sigma^{mon}\rangle$ first, one obtains
$\C\times\C^*\times\H$ with coordinates $(z_1,z_2,\tau)=(z_1,z_5^2,\tau)$
and an action of $PSL_2(\Z)$ on it, which is induced by the 
following action of the classes $[A_1]$ and $[A_2]$ in 
$PSL_2(\Z)$,
\begin{eqnarray*}
{}[A_1]\sim \sigma_1: (z_1,z_2,\tau)&\mapsto & 
(z_1,-z_2,\tau-1) =(z_1,-z_2,\mu_{(12)}(\tau)),\\
{}[A_2]\sim \sigma_2: (z_1,z_2,\tau)&\mapsto & 
(z_1,z_2T(\tau+1),\frac{\tau}{\tau+1})\\
&&=(z_1,z_2(1-T(\tau)),\mu_{(13)}(\tau)).
\end{eqnarray*}
The normal subgroup $\oooo{\Gamma(2)} =\langle [A_1^2],[A_2^2]
\rangle$ acts nontrivially only on the third factor $\H$
with quotient $\C-\{0,1\}$, so its quotient is $C_3^{pure}$.
The action of the quotient group $PSL_2(\Z)/\oooo{\Gamma(2)}
\cong S_3$ gives $C_3^{conf}$.

\medskip
We treat the cases $\GG_2\,\&\, C_4\, (\P^1A_1), C_5$
before the case $\GG_2\,\&\, C_3\, (A_2A_1)$.

\medskip
{\bf The reducible cases $\GG_2\,\&\, C_4\, (\P^1A_1), C_5$:}\\
By the proof of Theorem \ref{t7.11}, the group 
$(\Br_3)_{\uuuu{e}/\{\pm 1\}^3}= \langle\sigma_2^2,
(\sigma^{mon})^{-1}\sigma_1^2\rangle$ is the normal closure of
$\sigma_2^2$ in the group $(\Br_3)_{\uuuu{x}/\{\pm 1\}^3}= 
\langle \sigma_1,\sigma_2^2\rangle$, so especially a normal
subgroup.

The deck transformations $\psi(\sigma_2^2)$ and $
\psi((\sigma^{mon})^{-1}\sigma_1^2)$ of $C_3^{univ}$
in Theorem \ref {t9.1} (c) act nontrivially only on the
third factor $\H$, and there they act as the generators
$[A_2^2]$ and $[A_1^2]$ of $\oooo{\Gamma(2)}$. Therefore
\begin{eqnarray*}
C_3^{\uuuu{e}/\{\pm 1\}^3} \cong \C\times\C\times(\C-\{0,1\}) \quad
\textup{with coordinates }(z_1,\zeta,z_3).
\end{eqnarray*}
To obtain $C_3^{S/\{\pm 1\}^3}$ from $C_3^{\uuuu{e}/\{\pm 1\}^3}$, it is a priori sufficient to
divide out the action on $C_3^{\uuuu{e}/\{\pm 1\}^3}$ which the deck transformation
$\psi(\sigma_1)$ of $C_3^{univ}$ induces.
But it is easier to first divide out the action of 
$\sigma^{mon}\in\langle \sigma_1,\sigma_2^2\rangle
=(\Br_3)_{\uuuu{x}/\{\pm 1\}^3}$ on $C_3^{\uuuu{e}/\{\pm 1\}^3}$. 
The quotient is $C_3^{pure}\cong \C\times \C^*\times 
(\C-\{0,1\})$.
On this space the action of $\sigma_1$ is the action of 
$\phi_{(12)}$, which is a fixed point free involution.
Therefore $C_3^{S/\{\pm 1\}^3} = C_3^{pure}/\langle \phi_{(12)}\rangle$.

\medskip
{\bf The reducible case $\GG_2\,\&\, C_3\, (A_2A_1)$:}\\
Also here $(\Br_3)_{\uuuu{x}/\{\pm 1\}^3}= 
\langle \sigma_1,\sigma_2^2\rangle$, so 
$C_3^{S/\{\pm 1\}^3}=C_3^{pure}/\langle\phi_{(12)}\rangle$. 
Here $(Br_3)_{\uuuu{e}/\{\pm 1\}^3}= \langle \sigma_2^2,
(\sigma^{mon})^{-1}\sigma_1^2,\sigma^{mon}\sigma_1\rangle$. 
First $C_3^{univ}$ is divided by the action of the normal 
subgroup $\langle \sigma_2^2,(\sigma^{mon})^{-1}\sigma_1^2
\rangle$. The quotient is isomorphic to
$\C\times\C\times (\C-\{0,1\})$ with coordinates
$(z_1,\zeta,z_3)$ by the discussion in the cases 
$\GG_2\,\&\, C_4,C_5$.  

The deck transformation $\psi(\sigma^{mon}\sigma_1)$ on
$C_3^{univ}$ induces on this manifold the automorphism
$$(z_1,\zeta,z_3)\mapsto (z_1,\zeta-3\pi i,1-z_3).$$
If we divide out first the square of this action, the quotient
is isomorphic to $\C\times\C^*\times (\C-\{0,1\})$ with
coordinates $(z_1,z_6,z_3)$ with $z_6=e^{\zeta/3}$, 
and $\psi(\sigma^{mon}\sigma_1)$
induces on it the action of $\phi_{(12)}: 
(z_1,z_6,z_3)\mapsto (z_1,-z_6,1-z_3)$. We obtain the 
following diagram,
\begin{eqnarray*}
\begin{CD} 
(z_1,z_6,z_3) @>>> (z_1,z_6^3,z_3)=(z_1,z_2,z_3) \\
\C\times\C^*\times (\C-\{0,1\}) @>{3:1}>\text{covering}> 
\C\times\C^*\times (\C-\{0,1\}) \cong C_3^{pure} \\
@VV{/\langle\phi_{(12)}\rangle}V 
@VV{/\langle\phi_{(12)}\rangle}V \\
C_3^{\uuuu{e}/\{\pm 1\}^3} @>{3:1}>\text{covering}> 
C_3^{S/\{\pm 1\}^3} = C_3^{pure}/\langle \phi_{(12)}\rangle \\
\end{CD}
\end{eqnarray*}
Though $C_3^{\uuuu{e}/\{\pm 1\}^3}$ and 
$C_3^{S/\{\pm 1\}^3}$ can be presented in a better way.
Write $z_7:=z_3-\frac{1}{2}$. Then
$\phi_{(12)}:(z_1,z_6,z_7)\mapsto (z_1,-z_6,-z_7)$
and $\phi_{(12)}:(z_1,z_2,z_7)\mapsto 
(z_1,-z_2,-z_7)$. The following two diagrams belong together.
The upper diagram shows the sets,
the lower diagram shows the maps. 
The horizontal maps are isomorphisms. 
The sets on the right side in the upper diagram are 
better presentations of $C_3^{\uuuu{e}/\{\pm 1\}^3}$
and $C_3^{S/\{\pm 1\}^3}$ than the quotients on the left side.
See Corollary \ref{t9.8}. 
\begin{eqnarray*}
\begin{CD}
C_3^{\uuuu{e}/\{\pm 1\}^3}\cong 
\frac{\C\times\C^*\times 
(\C-\{\pm\frac{1}{2}\})}{\langle\phi_{(12)}\rangle}
@>{\cong}>> \C\times \{(z_8,z_9)\in\C^*\times\C\,|\, z_8^3-4z_9^2=0\}\\
@VV{3:1}V @VV{3:1}V \\
C_3^{S/\{\pm 1\}^3}\cong 
\frac{\C\times\C^*\times
(\C-\{\pm\frac{1}{2}\})}{\langle\phi_{(12)}\rangle}
@>{\cong}>> \C\times\{(z_{10},z_9)\in\C^*\times\C\,|\, z_{10}-4z_9^2=0\} 
\end{CD}
\end{eqnarray*}

\begin{eqnarray*}
\begin{CD}
 [(z_1,z_6,z_7)] @>>> (z_1,z_6^2,z_6^3z_7) @. 
 =(z_1,z_8,z_9)  \\
@VV{3:1}V @. @VV{3:1}V \\
[(z_1,z_6^3,z_7)] @. @. (z_1,z_8^3,z_9) \\
@|  @.  @|  \\
[(z_1,z_2,z_7)] @>>> (z_1,z_2^2,z_2z_7) @. 
=(z_1,z_{10},z_9) 
\end{CD}
\end{eqnarray*}
Geometrically, the horizontal maps are the restrictions of
the inversions on the level of 
$C_3^{\uuuu{e}/\{\pm 1\}^3}$ and $C_3^{S/\{\pm 1\}^3}$ 
of the blowing up described in Remark \ref{t9.2} (ii). 

\medskip
We treat the case $\GG_4\,\&\, C_7\, (\widehat{A}_2)$
before the case $\GG_3\,\&\, C_6\, (A_3)$.

\medskip
{\bf The case $\GG_4\,\&\, C_7\, (\widehat{A}_2)$:}\\
Here $(\Br_3)_{\uuuu{e}/\{\pm 1\}^3}= 
\langle \sigma_1^3,\sigma_2^3,\sigma_2\sigma_1^3\sigma_2^{-1}
\rangle$
and $(\Br_3)_{\uuuu{x}/\{\pm 1\}^3}= 
\langle \sigma_2\sigma_1,\sigma_1^3\rangle$.
The surjective homomorphism $\Br_3\to PSL_2(\Z)$ in
Remark \ref{t4.15} (i) catches the action of $\Br_3$ on the
third factor $\H$ of $\C\times\C\times\H\cong C_3^{univ}$.
It maps the subgroup  $\langle \sigma_1^3,\sigma_2^3,
\sigma_2\sigma_1^3\sigma_2^{-1}\rangle$ 
to the subgroup 
$\langle [A_1^3],[A_2^3],[A_2A_1^3A_2^{-1}]\rangle
=\oooo{\Gamma(3)}\subset PSL_2(\Z)$.
This group is isomorphic to $G^{free,3}$. 
The action of $\oooo{\Gamma(3)}$ on $\H$ is free.

The quotient $\H/\oooo{\Gamma(3)}$ is by Theorem \ref{tb.1} (c)
isomorphic to 
$$\P^1-\{\textup{the four vertices of a tetrahedron}\}
\cong \C-\{z_{11}\,|\, z_{11}^3=1\}=:\C^{***}.$$
Therefore the quotient $C_3^{\uuuu{e}/\{\pm 1\}^3}$ 
of $C_3^{univ}$ by 
$(\Br_3)_{\uuuu{e}/\{\pm 1\}^3}$ is isomorphic to
\begin{eqnarray*}
&&\C\times (\textup{an affine }\C\textup{-bundle over }
\C^{***})\\
&\cong& \C\times\C\times \C^{***}
\quad\textup{with coordinates }(z_1,\zeta,z_{11})
\end{eqnarray*}
by Remark \ref{t9.4} (ii).

In order to obtain $C_3^{S/\{\pm 1\}^3}$ it is because of 
$(\Br_3)_{\uuuu{x}/\{\pm 1\}^3} =\langle \sigma_2\sigma_1,\sigma_1^3\rangle$
a priori sufficient to divide out the action on 
$C_3^{\uuuu{e}/\{\pm 1\}^3}$ which the deck transformation $\psi(\sigma_2\sigma_1)$
of $C_3^{univ}$ induces. 
But it is easier to first divide out the action of 
$\sigma^{mon}\in\langle \sigma_2\sigma_1,\sigma_1^3\rangle$
on $C_3^{\uuuu{e}/\{\pm 1\}^3}$. The quotient is isomorphic to
$\C\times\C^*\times \C^{***}$
with coordinates $(z_1,z_2,z_{11})$. 

Now $\sigma_2\sigma_1$ acts because of $(\sigma_2\sigma_1)^3
=\sigma^{mon}$ as an automorphism of order three on the
$\C^*$-bundle $\C^*\times \C^{***}$.
To determine this automorphism, recall the action of 
$\sigma_2\sigma_1$ on $\C\times\C\times\H\cong C_3^{univ}$
in Theorem \ref{t9.1} (c),
\begin{eqnarray*}
\psi(\sigma_2\sigma_1):(z_1,\zeta,\tau)\mapsto
(z_1,\zeta-\pi i+\kappa(\tau),\frac{\tau-1}{\tau}).
\end{eqnarray*}
The M\"obius transformation 
$\mu_{(123)}=(\tau\mapsto \frac{\tau-1}{\tau})$ is elliptic
of order three with fixed point $\tau_0:=e^{2\pi i /6}$. It acts
on the tangent space $T_{\tau_0}\H$ by multiplication with
$e^{-2\pi i/3}$, because $\mu_{(123)}(\tau_0+\varepsilon)
\approx \tau_0+\varepsilon e^{-2\pi i/3}$ for small 
$\varepsilon\in\C$. The coordinate $z_{11}$ on $\C^{***}$
can be chosen so that $\sigma_2\sigma_1$ acts on $\C^{***}$
by $z_{11}\mapsto z_{11}e^{-2\pi i /3}$, especially
$\tau_0$ maps to $z_{11}=0$. On the $\C^*$-fiber
over $z_{11}=0$ $\sigma_2\sigma_1$ acts because of 
$T(\tau_0)=\tau_0$ by 
\begin{eqnarray*}
z_2=e^\zeta\mapsto e^{\zeta-\pi i+\kappa(\tau_0)}
=z_2(-1)T(\tau_0)=z_2 e^{-2\pi i /3}.
\end{eqnarray*}
By Remark \ref{t9.4} (iv) we can choose the trivialization
of the $\C^*$-bundle $\C^*\times \C^{***}$ over $\C^{***}$
such that the action of $\sigma_2\sigma_1$ on it is the action
$(z_2,z_{11})\mapsto (z_2e^{-2\pi i /3},z_{11}e^{-2\pi i /3})$.
We obtain the following diagram,
\begin{eqnarray*}
\begin{CD} 
(z_1,\zeta,z_{11}) @>>> (z_1,e^\zeta,z_{11})=(z_1,z_2,z_{11}) \\
C_3^{\uuuu{e}/\{\pm 1\}^3}\cong \C\times\C\times 
\C^{***}  @>{/\langle \sigma^{mon}\rangle}>> 
\C\times\C^*\times \C^{***}  \\
@. @VV{/\langle \sigma_2\sigma_1\rangle}V \\ 
\hspace*{2cm}C_3^{S/\{\pm 1\}^3}\cong @. 
\C\times\frac{\C^*\times \C^{***}}
{\langle (z_2,z_{11})\mapsto (z_2e^{2\pi i /3},z_{11}e^{2\pi i/3})
\rangle}  \\
\end{CD}
\end{eqnarray*}

\medskip
{\bf The case $\GG_3\,\&\, C_6\, (A_3)$:}\\
Here $(\Br_3)_{\uuuu{e}/\{\pm 1\}^3}= 
\langle (\sigma_1\sigma_2)^4,\sigma_1^3\rangle$
and $(\Br_3)_{\uuuu{x}/\{\pm 1\}^3}= 
\langle \sigma_1\sigma_2,\sigma_1^3\rangle$. 
Especially, $(\sigma_1\sigma_2)^4=\sigma_1\sigma_2\sigma^{mon}$
and
\begin{eqnarray*}
(\Br_3)_{\uuuu{e}/\{\pm 1\}^3}&\stackrel{3:1}{\supset}&
\langle \sigma_1^3,\sigma_2^{-1}\sigma_1^3\sigma_2
(=\sigma_1\sigma_2^3\sigma_1^{-1}),\sigma_2^3,(\sigma^{mon})^4
\rangle\\
&\supset& \langle\sigma_1^3,\sigma_2^{-1}\sigma_1^3\sigma_2,
\sigma_2^3 \rangle,
\end{eqnarray*}
and the group $\langle\sigma_1^3,\sigma_2^{-1}\sigma_1^3\sigma_2,
\sigma_2^3 \rangle$ maps under the surjective
homomorphism $\Br_3\to PSL_2(\Z)$ in Remark \ref{t4.15} (i) 
to the group 
$$\langle [A_1^3],[A_2^{-1}A_1^3A_2],[A_2^3]\rangle
=\oooo{\Gamma(3)}\subset PSL_2(\Z).$$
The action of $\oooo{\Gamma(3)}$ on $\H$ is free, and the
quotient is by Theorem \ref{tb.1} (c) isomorphic to $\C^{***}$.
The quotient of $C_3^{univ}$ by
$\langle \sigma_1^3,\sigma_2^{-1}\sigma_1^3\sigma_2,
\sigma_2^3\rangle$ is isomorphic to 
\begin{eqnarray*}
&&\C\times (\textup{an affine }\C\textup{-bundle over }\C^{***})
\\
&\cong& \C\times\C\times\C^{***}\quad\textup{with coordinates }
(z_1,\zeta,z_{11})
\end{eqnarray*}
by Remark \ref{t9.4} (ii). 
The quotient of $C_3^{univ}$ by 
$\langle \sigma_1^3,\sigma_2^{-1}\sigma_1^3\sigma_2,
\sigma_2^3,(\sigma^{mon})^4\rangle$ is isomorphic to 
\begin{eqnarray*}
\C\times\C^*\times\C^{***}\quad\textup{with coordinates }
(z_1,z_{12},z_{11})\textup{ with }z_{12}=e^{\zeta/4}.
\end{eqnarray*}
On this quotient $(\sigma_1\sigma_2)^4$ acts as an automorphism
of order three. We will determine the action now.

First, it acts on $\C\times\C^*\times\H$ with coordinates
$(z_1,z_{12},\tau)$ with $z_{12}=e^{\zeta/4}$ by
$$(\sigma_1\sigma_2)^4:(z_1,z_{12},\tau)\mapsto 
(z_1,z_{12}e^{-2\pi i 3/8}e^{\kappa(\tau+1)/4},\frac{-1}{\tau+1}).
$$
The M\"obius transformation $(\tau\mapsto \frac{-1}{\tau+1})$
is elliptic of order three with fixed point 
$\tau_1:=e^{2\pi i /3}$. It acts on the tangent space
$T_{\tau_1}\H$ by multiplication with $e^{-2\pi i/3}$ because
$\frac{-1}{\tau_1+\varepsilon}\equiv \tau_1+\varepsilon
e^{-2\pi i/3}$ for small $\varepsilon\in\C$. 
The coordinate $z_{11}$ on $\C^{***}$ can be chosen so that 
$(\sigma_1\sigma_2)^4$ acts on $\C^{***}$ by 
$z_{11}\mapsto z_{11}e^{-2\pi i/3}$, especially $\tau_1$ maps to 
$z_{11}=0$. 
On the $\C^*$-fiber over $z_{11}=0$ $(\sigma_1\sigma_2)^4$ acts
because of $\kappa(\tau_1+1)=\kappa(e^{2\pi i/6})=
\frac{2\pi i}{6}$ by
\begin{eqnarray*}
z_{12}=e^{\zeta/4}\mapsto e^{(\zeta-3\pi i+\kappa(\tau_1+1))/4}
=z_{12}e^{-2\pi i 3/8+2\pi i/24}=z_{12}e^{-2\pi i/3}.
\end{eqnarray*}
By Remark \ref{t9.4} (ii) we can choose the trivialization
of the $\C^*$-bundle $\C^*\times\C^{***}$ over $\C^{***}$ 
with coordinates $(z_{12},z_{11})$ such that the action of
$(\sigma_1\sigma_2)^4$ on it becomes the action 
$(z_{12},z_{11})\mapsto (z_{12}e^{-2\pi i/3},z_{11}e^{-2\pi i/3})$. 

$C_3^{S/\{\pm 1\}^3}$ is obtained by dividing out the action 
of $\sigma^{mon}$, because $(\Br_3)_{\uuuu{x}/\{\pm 1\}^3}
=\langle\sigma_1,\sigma_2^2\rangle$ is generated by
$(\Br_3)_{\uuuu{e}/\{\pm 1\}^3}$ and by $\sigma^{mon}$,
with $(\sigma^{mon})^4\in (\Br_3)_{\uuuu{x}/\{\pm 1\}^3}$.

We obtain the following diagram, 
\begin{eqnarray*}
\begin{CD} 
(z_1,z_{12},z_{11}) @>>> (z_1,z_{12}^4,z_{11})=(z_1,z_2,z_{11}) \\
\C\times\C^*\times \C^{***}  @>{4:1}>> 
\C\times\C^*\times \C^{***}  \\
@V{3:1}VV  @V{3:1}VV \\ 
\C\times\frac{\C^*\times \C^{***}}
{\langle (z_{12},z_{11})\mapsto (z_{12}e^{2\pi i /3},z_{11}e^{2\pi i/3})
\rangle} 
@>{4:1}>>
\C\times\frac{\C^*\times \C^{***}}
{\langle (z_2,z_{11})\mapsto (z_2e^{2\pi i /3},z_{11}e^{2\pi i/3})
\rangle}  \\
\cong C_3^{\uuuu{e}/\{\pm 1\}^3} @. 
\cong C_3^{S/\{\pm 1\}^3} \\
\end{CD}
\end{eqnarray*}
Though $C_3^{\uuuu{e}/\{\pm 1\}^3}$ and 
$C_3^{S/\{\pm 1\}^3}$ can be presented in a better way.
Write $\xi:=e^{2\pi i/3}$. The following two diagrams belong
together. The upper diagram shows the sets,
the lower diagram shows the maps. 
The horizontal maps are isomorphisms.  
In the line of $C_3^{\uuuu{e}/\{\pm 1\}^3}$ the set on the right
side is nicer than the set on the left side.
In the line of $C_3^{S/\{\pm 1\}^3}$ both sets are quotients. 
\begin{eqnarray*}
\begin{CD}
C_3^{\uuuu{e}/\{\pm 1\}^3}\cong 
\frac{\C\times\C^*\times \C^{***}}{\langle  
(z_1,z_{12},z_{11})\mapsto(z_1,z_{12}\xi,z_{11}\xi) \rangle}
@>{\cong}>> \C\times \{(z_{13},z_{14})\in\C^*\times\C\,|\, z_{14}^3-z_{13}^2=0\}\\
@VV{4:1}V @VV{4:1}V \\
C_3^{S/\{\pm 1\}^3}\cong 
\frac{\C\times\C^*\times \C^{***}}{\langle 
(z_1,z_2,z_{11})\mapsto (z_1,z_2\xi,z_{11}\xi) \rangle}
@>{\cong}>> \frac{\C\times\{(z_{15},z_{14})\in\C^*\times\C\,|\,z_{14}^3-z_{15}=0\}}
{\langle (z_1,z_{15},z_{14})\mapsto (z_1,-z_{15},-z_{14})\rangle} 
\end{CD}
\end{eqnarray*}

\begin{eqnarray*}
\begin{CD}
 [(z_1,z_{12},z_{11})] @>>> (z_1,z_{12}^3,z_{12}^2z_{11}) @. 
 =(z_1,z_{13},z_{14})  \\
@VV{4:1}V @. @VV{4:1}V \\
[(z_1,z_{12}^4,z_{11})] @. @. [(z_1,z_{13}^2,z_{14})] \\
@|  @.  @|  \\
[(z_1,z_2,z_{11})] @>>> [(z_1,z_2^{3/2},z_2^{1/2}z_{11})] @. 
=[(z_1,z_{15},z_{14})] 
\end{CD}
\end{eqnarray*}
We claim that geometrically, the horizontal maps are the 
restrictions of the inversions on the level of 
$C_3^{\uuuu{e}/\{\pm 1\}^3}$ and $C_3^{S/\{\pm 1\}^3}$ 
of the blowing up described in Remark \ref{t9.2} (ii).

\medskip
{\bf The cases $\GG_5\,\&\, C_8,C_9\, ((-l,2,-l))$:}\\
Here $(\Br_3)_{\uuuu{e}/\{\pm 1\}^3}= 
\langle (\sigma^{mon})^2\sigma_1^{-1}\sigma_2^{l^2-4}\sigma_1
\rangle$
and $(\Br_3)_{\uuuu{x}/\{\pm 1\}^3}= 
\langle \sigma^{mon},\sigma_1^{-1}\sigma_2\sigma_1\rangle$.
On the last factor $\H$ of $\C\times\C\times\H
\cong C_3^{univ}$, $\sigma_1^{-1}\sigma_2\sigma_1$ and
$\sigma_1^{-1}\sigma_2^{l^2-4}\sigma_1=
(\sigma_1^{-1}\sigma_2\sigma_1)^{l^2-4}$ act as parabolic
elements with fixed point $1\in\whh{\R}$, so freely, 
and $\sigma^{mon}$ acts trivially.
The quotients of $\H$ by the actions of 
$\sigma_1^{-1}\sigma_2^{l^2-4}\sigma_1$ and 
$\sigma_1^{-1}\sigma_2\sigma_1$ are both isomorphic to
$\D^*=\D-\{0\}$, where $\D=\{z\in\C\,|\,|z|<1\}$ is the
unit disk, and the covering is as follows, 
\begin{eqnarray*}
\H/(\textup{action of }\sigma_1^{-1}\sigma_2^{l^2-4}\sigma_1)
\cong \D^* &\to& \D^*\cong \H/(\textup{action of }
\sigma_1^{-1}\sigma_2\sigma_1),\\
z_{16}& \mapsto & z_{16}^{l^2-4}.
\end{eqnarray*}
Recall that $\psi(\sigma_2^2)$ acts trivially on the second
factor of $\C\times\C\times\H \cong C_3^{univ}$.
The total action of $\psi((\sigma^{mon})^2\sigma_1^{-1}
\sigma_2^{l^2-4}\sigma_1)$ on 
$\C\times\C\times \H\cong C_3^{univ}$ is as follows,
for even $l$:
\begin{eqnarray*}
(z_1,\zeta,\tau)&\stackrel{\sigma_1}{\mapsto}&
(z_1,\zeta-\pi i ,\tau-1)\\
&\stackrel{\sigma_2^{l^2-4}}{\mapsto}&
(z-1,\zeta-\pi i,\frac{\tau-1}{(l^2-4)(\tau-1)+1})\\
&\stackrel{\sigma_1^{-1}}{\mapsto}&
(z_1,\zeta,\frac{(\tau-1)+(l^2-4)(\tau-1)+1}
{(l^2-4)(\tau-1)+1})\\
&\stackrel{(\sigma^{mon})^2}{\mapsto}&
(z_1,\zeta-4\pi i, \frac{(l^2-3)(\tau-1)+1}
{(l^2-4)(\tau-1)+1}),
\end{eqnarray*}
and for odd $l$:
\begin{eqnarray*}
(z_1,\zeta,\tau)&\mapsto&
(z_1,\zeta-4\pi i+\kappa(\tau), \frac{(l^2-3)(\tau-1)+1}
{(l^2-4)(\tau-1)+1}),
\end{eqnarray*}
Because the action on the third factor $\H$ is free,
the action on the second factor $\C$ is not important for
the quotient up to isomorphism. The quotient is for even $l$
and for odd $l$
\begin{eqnarray*}
C_3^{\uuuu{e}/\{\pm 1\}^3} &\cong & \C\times (\textup{an affine }\C\textup{-bundle
over }\D^*)\\
&\cong & \C\times\C\times \D^*
\quad\textup{with coordinates }(z_1,\zeta,z_{16})
\end{eqnarray*}
by Remark \ref{t9.4} (ii).

Dividing out the action of $\sigma^{mon}$ gives 
\begin{eqnarray*}
C_3^{\uuuu{e}/\{\pm 1\}^3} \cong \C\times\C\times \D^* &\to & 
\C\times \C^*\times \D^* \\
(z_1,\zeta,z_{16})&\mapsto& (z_1,e^{\zeta},z_{16})=(z_1,z_2,z_{16}).
\end{eqnarray*}
On this quotient $\sigma_1^{-1}\sigma_2\sigma_1$ acts as a
cyclic automorphism of order $l^2-4$ because 
$(\sigma_1^{-1}\sigma_2\sigma_1)^{l^2-4} \in 
\langle \sigma^{mon},(\sigma^{mon})^2\sigma_1^{-1}
\sigma_2^{l^2-4}\sigma_1\rangle$. 
The action on the third factor $\D^*$ is free.
Therefore we can choose a posteriori the trivialization 
$\C\times\D^*$ with coordinates $(\zeta,z_{16})$ and
$\C^*\times\D^*$ with coordinates $(z_2,z_{16})$ so that we obtain
the following diagram,
\begin{eqnarray*}
\begin{CD} 
(z_1,\zeta,z_{16}) @>>> (z_1,e^\zeta,z_{16}) @.  \\
C_3^{\uuuu{e}/\{\pm 1\}^3}\cong \C\times\C\times \D^* 
@>{/\langle \sigma^{mon}\rangle}>> 
\C\times\C^*\times \D^* @. \hspace*{0.5cm}(z_1,z_2,z_{16})  \\
@.  @V{(l^2-4):1}VV 
@VV{/\langle \sigma_1^{-1}\sigma_2\sigma_1\rangle}V \\ 
\hspace*{2cm}C_3^{S/\{\pm 1\}^3}\cong 
@. \C\times\C^*\times \D^* @. 
\hspace*{0.5cm}(z_1,z_2,z_{16}^{l^2-4})  \\
\end{CD}
\end{eqnarray*}

\medskip
{\bf The cases $\GG_6\,\&\, C_{10}\, (\P^2),C_{11},C_{12}$:}\\
Here $(\Br_3)_{\uuuu{e}/\{\pm 1\}^3}= \{\id\}$
and $(\Br_3)_{\uuuu{x}/\{\pm 1\}^3}= \langle\sigma_2\sigma_1
\rangle$. Therefore
$$C_3^{\uuuu{e}/\{\pm 1\}^3} = C_3^{univ} \cong \C\times\C\times\H.$$
In order to obtain $C_3^{S/\{\pm 1\}^3}$ it is a priori sufficient to divide
out the action of $\sigma_2\sigma_1$. But it is easier to first
divide out the action of $\sigma^{mon}=(\sigma_2\sigma_1)^3$.
On the quotient $C_3^{\uuuu{e}/\{\pm 1\}^3} /\langle \sigma^{mon}\rangle
\cong \C\times\C^*\times \H$ with coordinates $(z_1,z_2,\tau)$
$\sigma_2\sigma_1$ acts as automorphism of order three.
Concretely,  
\begin{eqnarray*}
\begin{CD} 
@. (z_1,\zeta,\tau) @>>> (z_1,\zeta-\pi i+\kappa(\tau),
\frac{\tau-1}{\tau}) \\
(z_1,\zeta,\tau)\hspace*{0.5cm} @. 
C_3^{\uuuu{e}/\{\pm 1\}^3}\cong \C^2\times \H 
@>{\sigma_2\sigma_1}>> C_3^{\uuuu{e}/\{\pm 1\}^3}\cong
\C^2\times \H   \\
@VVV
@VV{/\langle \sigma^{mon}\rangle}V   
@VV{/\langle \sigma^{mon}\rangle}V  \\ 
(z_1,e^\zeta,\tau)\hspace*{0.5cm} 
@. \C\times\C^*\times \H
@>{\sigma_2\sigma_1}>> 
\C\times\C^*\times \H \\
@. 
(z_1,z_2,\tau) @>>> 
(z_1,z_2(-T(\tau)),\frac{\tau-1}{\tau})  \\
\end{CD}
\end{eqnarray*}
The action is fixed point free.
One reason is that $C_3^{S/\{\pm 1\}^3}$ is smooth.
Another reason is given by the following chain of arguments. 
The M\"obius transformation 
$\mu_{(123)}=(\tau\mapsto \frac{\tau-1}{\tau})$ is elliptic
of order three with fixed point $e^{2\pi i /6}$.
The quotient is again isomorphic to $\H$. On the fiber
of $\C^*\times\H$ over $e^{2\pi i /6}$, $\sigma_2\sigma_1$
acts by multiplication with 
$-T(e^{2\pi i /6})=-e^{2\pi i /6}=e^{-2\pi i/3}$, so fixed
point free. Also 
$$(-T(\tau))(-T(\frac{\tau-1}{\tau}))(-T(\frac{-1}{\tau-1}))
=-T(\tau)\frac{T(\tau)-1}{T(\tau)}\frac{1}{1-T(\tau)}=1$$
because of Theorem \ref{tb.1} (d). 
Therefore $\sigma_2\sigma_1$ acts on $\C\times\C^*\times\H$
fixed point free and has order three.

We obtain the following diagram, 
\begin{eqnarray*}
\begin{CD} 
(z_1,\zeta,\tau) @>>> (z_1,e^\zeta,\tau)=(z_1,z_2,\tau) \\
C_3^{\uuuu{e}/\{\pm 1\}^3}\cong \C\times\C\times \H 
@>{/\langle \sigma^{mon}\rangle}>> 
\C\times\C^*\times \H   \\
@.  @VV{/\langle \sigma_2\sigma_1\rangle}V  \\ 
C_3^{S/\{\pm 1\}^3}\cong
@. \C\times\frac{\C^*\times \H}
{\langle(z_2,\tau)\mapsto 
(z_2(-T(\tau)),\frac{\tau-1}{\tau})
\rangle}  \\
\end{CD}
\end{eqnarray*}
$C_3^{S/\{\pm 1\}^3}$ is $\C\times(\textup{a }\C^*
\textup{-orbibundle over }\H/\langle \mu_{(123)}\rangle)$ 
in the sense of Remark \ref{t9.5} with the only exceptional 
fiber over $[e^{2\pi i/6}]\in \H/\langle \mu_{(123)}\rangle$.

\medskip
{\bf The cases $\GG_7\,\&\, C_{13}\ (\textup{e.g. }(4,4,8))$:}\\
Here $(\Br_3)_{\uuuu{e}/\{\pm 1\}^3}= \{\id\}$
and $(\Br_3)_{\uuuu{x}/\{\pm 1\}^3}= \langle\sigma_2\sigma_1^2
\rangle$. Therefore
$$C_3^{\uuuu{e}/\{\pm 1\}^3} = C_3^{univ} \cong \C\times\C\times\H.$$
In order to obtain $C_3^{S/\{\pm 1\}^3}$ it is a priori sufficient to divide
out the action of $\sigma_2\sigma_1^2$. But it is easier to first
divide out the action of $\sigma^{mon}=(\sigma_2\sigma_1^2)^2$.
On the quotient $C_3^{\uuuu{e}/\{\pm 1\}^3} /\langle \sigma^{mon}\rangle
\cong \C\times\C^*\times \H$ with coordinates $(z_1,z_2,\tau)$
$\sigma_2\sigma_1^2$ acts as automorphism of order two.
Concretely, 
\begin{eqnarray*}
\begin{CD} 
@. (z_1,\zeta,\tau) @>>> (z_1,\zeta-2\pi i+\kappa(\tau-1),
\frac{\tau-2}{\tau-1}) \\
(z_1,\zeta,\tau)\hspace*{0.5cm} @. 
C_3^{\uuuu{e}/\{\pm 1\}^3}\cong \C^2\times \H 
@>{\sigma_2\sigma_1^2}>> C_3^{\uuuu{e}/\{\pm 1\}^3}\cong
\C^2\times \H   \\
@VVV
@VV{/\langle \sigma^{mon}\rangle}V   
@VV{/\langle \sigma^{mon}\rangle}V  \\ 
(z_1,e^\zeta,\tau)\hspace*{0.5cm} 
@. \C\times\C^*\times \H
@>{\sigma_2\sigma_1^2}>> 
\C\times\C^*\times \H \\
@. 
(z_1,z_2,\tau) @>>> 
(z_1,z_2(1-T(\tau)),\frac{\tau-2}{\tau-1})  \\
\end{CD}
\end{eqnarray*}
The action is fixed point free.
One reason is that $C_3^{S/\{\pm 1\}^3}$ is smooth.
Another reason is given by the following chain of arguments. 
The M\"obius transformation 
$(\tau\mapsto \frac{\tau-2}{\tau-1})$ is elliptic
of order two with fixed point $1+i$.
The quotient is again isomorphic to $\H$. On the fiber
of $\C^*\times\H$ over $1+i$, $\sigma_2\sigma_1^2$
acts by multiplication with 
$1-T(1+i)=1-2=-1$, so fixed
point free. Also 
\begin{eqnarray*}
(1-T(\tau))(1-T(\frac{\tau-2}{\tau-1}))
&=&(1-T(\tau))(1-T(\mu_{(123)}\mu_{(12)}))\\
&=&(1-T(\tau))(1-g_{(123)}g_{(12)}(T(\tau)))\\
&=&(1-T(\tau))(1-\frac{(1-T(\tau))-1}{1-T(\tau)})=1
\end{eqnarray*}
because of Theorem \ref{tb.1} (d). 
Therefore $\sigma_2\sigma_1$ acts on $\C\times\C^*\times\H$
fixed point free and has order two.

We obtain the following diagram, 
\begin{eqnarray*}
\begin{CD} 
(z_1,\zeta,\tau) @>>> (z_1,e^\zeta,\tau)=(z_1,z_2,\tau) \\
C_3^{\uuuu{e}/\{\pm 1\}^3}\cong \C\times\C\times \H 
@>{/\langle \sigma^{mon}\rangle}>> 
\C\times\C^*\times \H   \\
@.  @VV{/\langle \sigma_2\sigma_1^2\rangle}V  \\ 
C_3^{S/\{\pm 1\}^3}\cong
@. \C\times\frac{\C^*\times \H}
{\langle(z_2,\tau)\mapsto 
(z_2(1-T(\tau)),\frac{\tau-2}{\tau-1})
\rangle}  \\
\end{CD}
\end{eqnarray*}
$C_3^{S/\{\pm 1\}^3}$ is $\C\times(\textup{a }\C^*
\textup{-orbibundle over }
\H/\langle (\tau\mapsto \frac{\tau-2}{\tau-1})\rangle)$ 
in the sense of 
Remark \ref{t9.5} with the only exceptional fiber
over $[1+i]\in \H/\langle (\tau\mapsto \frac{\tau-2}{\tau-1})
\rangle$.

\medskip
{\bf The cases $\GG_8\,\&\, C_{14}$ and 
$\GG_9\,\&\, C_{15},C_{16},C_{23},C_{24}$:}\\
Here $(\Br_3)_{\uuuu{e}/\{\pm 1\}^3}= \{ \id\}$
and $(\Br_3)_{\uuuu{x}/\{\pm 1\}^3}= \langle\sigma^{mon}\rangle$.
Therefore
\begin{eqnarray*}
\begin{CD}
C_3^{\uuuu{e}/\{\pm 1\}^3} =C_3^{univ}\cong\C\times\C\times\H 
@>{/\langle \sigma^{mon}\rangle}>>
\C\times\C^*\times \H \cong C_3^{S/\{\pm 1\}^3}\\
(z_1,\zeta,\tau) @>>> (z_1,e^{\zeta},\tau)
=(z_1,z_2,\tau).
\end{CD}
\end{eqnarray*}

\medskip
{\bf The cases $\GG_{10}\,\&\, C_{17}
\ (\textup{e.g. }(-2,-2,0))$:}\\
Here $(\Br_3)_{\uuuu{e}/\{\pm 1\}^3}= \langle\sigma_2^2\rangle$
and $(\Br_3)_{\uuuu{x}/\{\pm 1\}^3}= \langle\sigma^{mon},
\sigma_2\rangle$.
The deck transformation $\psi(\sigma_2^2)$ of
$C_3^{univ}\cong \C\times\C\times\H$ acts nontrivially only on
the third factor $\H$, and there it acts as the parabolic
transformation $(\tau\mapsto \frac{\tau}{2\tau+1})$
with fixed point $0\in\whh{\R}$ and quotient 
$\H/\langle (\tau\mapsto \frac{\tau}{2\tau+1})\rangle
\cong\D^*$. Therefore
$$C_3^{\uuuu{e}/\{\pm 1\}^3} \cong \C\times\C\times \D^*
\quad\textup{with coordinates }(z_1,\zeta,z_{16}).$$
On this quotient $\sigma_2$ acts as fixed point free automorphism
of order two with the action
$\D^*\to \D^*$, $z_{16}\mapsto z_{16}^2=z_{17}$, on the last factor.
The quotient is isomorphic to $\C\times\C\times\D^*$
with coordinates $(z_1,\zeta,z_{17})$.
$\sigma^{mon}$ acts on the quotient with nontrivial action
$\zeta\mapsto \zeta-2\pi i$ only on the second factor $\C$.
We obtain the following diagram,
\begin{eqnarray*}
\begin{CD} 
(z_1,\zeta,z_{16}) @>>> (z_1,\zeta+\kappa(\tau+1),z_{16}^2)
@. \textup{with }z_{16}=[\tau] \\
C_3^{\uuuu{e}/\{\pm 1\}^3}\cong \C\times\C\times \D^* 
@>{/\langle \sigma_2\rangle}>> 
\C\times\C\times \D^* @. 
\hspace*{1cm}(z_1,\zeta,z_{17})  \\
@.  @VV{/\langle \sigma^{mon}\rangle}V  @VVV\\ 
\hspace*{2cm}C_3^{S/\{\pm 1\}^3}\cong 
@. \C\times\C^*\times \D^* @. 
\hspace*{1cm}(z_1,e^{\zeta},z_{17})  \\
\end{CD}
\end{eqnarray*}

\medskip
{\bf The cases $\GG_{11}\,\&\, C_{18}
\ (\textup{e.g. }(-3,-2,0))$:}\\
Here $(\Br_3)_{\uuuu{e}/\{\pm 1\}^3}=\langle\sigma_2^2\rangle$
and $(\Br_3)_{\uuuu{x}/\{\pm 1\}^3}=\langle \sigma^{mon},
\sigma_2^2\rangle $. The manifold $C_3^{\uuuu{e}/\{\pm 1\}^3}$ 
is the same quotient of $C_3^{univ}$ as in the cases 
$\GG_{10}\,\&\, C_{17}$, namely
$$C_3^{\uuuu{e}/\{\pm 1\}^3} \cong \C\times\C\times \D^*
\quad\textup{with coordinates }(z_1,\zeta,z_{16}).$$
$\sigma^{mon}$ acts nontrivially only on the second factor $\C$
of this manifold, by $\zeta\mapsto \zeta-2\pi i$.
We obtain the following diagram, 
\begin{eqnarray*}
\begin{CD}
C_3^{\uuuu{e}/\{\pm 1\}^3} \cong\C\times\C\times\D^* 
@>{/\langle\sigma^{mon}\rangle}>>
\C\times\C^*\times \D^* \cong C_3^{S/\{\pm 1\}^3}\\
(z_1,\zeta,z_{16}) @>>> (z_1,e^{\zeta},z_{16})
=(z_1,z_2,z_{16}).
\end{CD}
\end{eqnarray*}

\begin{figure}[H]
\includegraphics[width=0.9\textwidth]{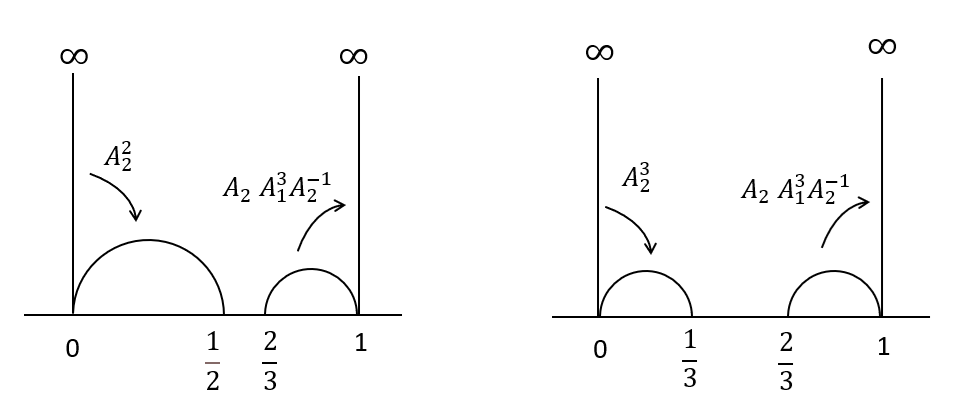}
\caption[Figure 9.2]{Fundamental domains in $\H$ for the
subgroups $\langle [A_2^2],[A_2A_1^3A_2^{-1}]\rangle$ (left) and
$\langle [A_2^3],[A_2A_1^3A_2^{-1}]\rangle$ (right) of 
$PSL_2(\Z)$}
\label{Fig:9.2}
\end{figure}

\medskip
{\bf The cases $\GG_{12}\,\&\, C_{19}
\ (\textup{e.g. }(-2,-1,0))$:}\\
Here $(\Br_3)_{\uuuu{e}/\{\pm 1\}^3}= \langle \sigma_2^2,
\sigma_2\sigma_1^3\sigma_2^{-1}\rangle$
and $(\Br_3)_{\uuuu{x}/\{\pm 1\}^3}= \langle\sigma^{mon},
\sigma_2^2,\sigma_2\sigma_1^3\sigma_2^{-1}\rangle$.
The group $\langle \sigma_2^2,\sigma_2\sigma_1^3\sigma_2^{-1}
\rangle$ acts on $\H$ as the group 
$\langle [A_2^2],[A_2A_1^3A_2^{-1}]\rangle\subset PSL_2(\Z)$.
Both generators are parabolic. Consider the hyperbolic polygon
$P$ whose relative boundary consists of the four arcs
$A(\infty,0)$, $A(0,\frac{1}{2})$, $A(\frac{2}{3},1)$,
$A(1,\infty)$, see the left picture in Figure \ref{Fig:9.2}.

Recall
\begin{eqnarray*}
A_1&=&\begin{pmatrix}1&-1\\0&1\end{pmatrix},\ 
A_2=\begin{pmatrix}1&0\\1&1\end{pmatrix},\
A_2^2=\begin{pmatrix}1&0\\2&1\end{pmatrix},\\
A_2A_1A_2^{-1}&=&\begin{pmatrix}2&-1\\0&1\end{pmatrix},\ 
A_2A_1^3A_2^{-1}=(A_2A_1A_2^{-1})^3
=\begin{pmatrix}4&-3\\3&-2\end{pmatrix}.
\end{eqnarray*}
$[A_2^2]$ has the fixed point $0\in\whh{\R}$ and maps 
$A(\infty,0)$ to $A(0,\frac{1}{2})$, and $[A_2A_1^3A_2^{-1}]$ 
has the fixed point $1\in\whh{\R}$ and maps $A(\frac{2}{3},1)$ to 
$A(1,\infty)$. By Theorem \ref{ta.2} (c) the hyperbolic polygon 
$P$ is a fundamental domain for the group $\langle [A_2^2],
[A_2A_1^3A_2^{-1}]\rangle$. Because of the two parabolic
fixed points and the euclidean boundary
$[\frac{1}{2},\frac{2}{3}]$ of $P$, the quotient
$\H/\langle [A_2^2],[A_2A_1^3A_2^{-1}]\rangle$ is isomorphic
to $\D-\{0,\frac{1}{2}\}=:\D^{**}$
where $z_{16}=0$ is the image of $0\in\whh{\R}$ and 
$z_{16}=\frac{1}{2}$ is the image of $1\in\whh{\R}$.  

As the action of $\langle\sigma_2^2,
\sigma_2\sigma_1^3\sigma_2^{-1}\rangle$ on the third factor 
$\H$ of $\C\times\C\times\H\cong C_3^{univ}$ is free,
the quotient of $\C\times\H$ with coordinates 
$(\zeta,\tau)$ is an affine $\C$-bundle over
$\H/\langle [A_2^2],[A_2A_1^3A_2^{-1}]\rangle$, so 
\begin{eqnarray*}
C_3^{\uuuu{e}/\{\pm 1\}^3} &\cong&\C\times (\textup{an affine }\C\textup{-bundle over }\D^{**})\\
&\cong& \C\times\C\times \D^{**}
\quad\textup{with coordinates }(z_1,\zeta,z_{16})
\end{eqnarray*}
by Remark \ref{t9.4} (ii). 
$\sigma^{mon}$ acts nontrivially only on the second factor
of this product, by the action $\zeta\to\zeta-2\pi i$.
Therefore 
\begin{eqnarray*}
\begin{CD}
C_3^{\uuuu{e}/\{\pm 1\}^3} \cong \C\times\C\times \D^{**}
@>{/\langle\sigma^{mon}\rangle}>>
\C\times\C^*\times \D^{**} \cong C_3^{S/\{\pm 1\}^3}\\
(z_1,\zeta,z_{16}) @>>> (z_1,e^{\zeta},z_{16})
=(z_1,z_2,z_{16}).\\
\end{CD}
\end{eqnarray*}

\medskip
{\bf The cases $\GG_{13}\,\&\, C_{20}
\ (\textup{e.g. }(-2,-1,-1))$:}\\
Here $(\Br_3)_{\uuuu{e}/\{\pm 1\}^3}= \langle \sigma_2^3,
\sigma_2\sigma_1^3\sigma_2^{-1}\rangle$
and $(\Br_3)_{\uuuu{x}/\{\pm 1\}^3}= \langle\sigma^{mon},
\sigma_2^3,\sigma_2\sigma_1^3\sigma_2^{-1}\rangle$.
This is analogous to the cases $\GG_{12}\,\&\, C_{19}$. 
The group $\langle \sigma_2^3,\sigma_2\sigma_1^3\sigma_2^{-1}
\rangle$ acts on $\H$ as the group 
$\langle [A_2^3],[A_2A_1^3A_2^{-1}]\rangle\subset PSL_2(\Z)$.
Both generators are parabolic. Consider the hyperbolic polygon
$P$ whose relative boundary consists of the four arcs
$A(\infty,0)$, $A(0,\frac{1}{3})$, $A(\frac{2}{3},1)$,
$A(1,\infty)$, see the right picture in Figure \ref{Fig:9.2}.

$[A_2^3]$ has the fixed point $0$ and maps $A(\infty,0)$
to $A(0,\frac{1}{3})$, and $[A_2A_1^3A_2^{-1}]$ has the 
fixed point $1$ and maps $A(\frac{2}{3},1)$ to 
$A(1,\infty)$. By Theorem \ref{ta.2} (c) the hyperbolic polygon 
$P$ is a fundamental domain for the group $\langle [A_2^3],
[A_2A_1^3A_2^{-1}]\rangle$.
As in the cases $\GG_{12}\,\&\, C_{19}$ we obtain
\begin{eqnarray*}
\begin{CD}
C_3^{\uuuu{e}/\{\pm 1\}^3} \cong \C\times\C\times \D^{**}
@>{/\langle\sigma^{mon}\rangle}>>
\C\times\C^*\times \D^{**} \cong C_3^{S/\{\pm 1\}^3}\\
(z_1,\zeta,z_{16}) @>>> (z_1,e^{\zeta},z_{16})
=(z_1,z_2,z_{16}).\\
\end{CD}
\end{eqnarray*}

{\bf The cases $\GG_{14}\,\&\, C_{21}, C_{22}$:}\\
Here $(\Br_3)_{\uuuu{e}/\{\pm 1\}^3}=\langle\sigma_2^3\rangle$
and $(\Br_3)_{\uuuu{x}/\{\pm 1\}^3}=\langle \sigma^{mon},
\sigma_2^3\rangle$.
The group $\langle \sigma_2^3\rangle$ acts on the third
factor $\H$ of $\C\times\C\times\H\cong C_3^{univ}$
as the group $\langle [A_2^3]\rangle\subset PSL_2(\Z)$.
The generator $[A_2^3]$ is parabolic. Therefore the quotient
$\H/\langle [A_2^3]\rangle$ is isomorphic to $\D^*$.
We obtain 
\begin{eqnarray*}
C_3^{\uuuu{e}/\{\pm 1\}^3} &\cong&\C\times (\textup{affine }\C\textup{-bundle over } \D^{*})\\
&\cong& \C\times\C\times \D^{*}
\quad\textup{with coordinates }(z_1,\zeta,z_{16}) 
\end{eqnarray*}
by Remark \ref{t9.4} (ii) and the diagram  
\begin{eqnarray*}
\begin{CD}
C_3^{\uuuu{e}/\{\pm 1\}^3} \cong \C\times\C\times \D^{*}
@>{/\langle\sigma^{mon}\rangle}>>
\C\times\C^*\times \D^{*} \cong C_3^{S/\{\pm 1\}^3}\\
(z_1,\zeta,z_{16}) @>>> (z_1,e^{\zeta},z_{16})
=(z_1,z_2,z_{16}).\\
\end{CD}
\end{eqnarray*}
\hfill$\Box$

\begin{remarks}\label{t9.4}
Let $X$ be a non-compact smooth complex curve.
\index{non-compact smooth complex curve}

(i) An {\it affine $\C$-bundle over $X$} 
\index{affine $\C$-bundle} is a holomorphic bundle
$p:Y\to X$ such that the following exist: An open covering
$(U_i)_{i\in I}$ of $X$, isomorphisms
$f_i:p^{-1}(U_i)\to\C\times U_i$, and for $i$ and $j\in I$
with $i\neq j$ and $U_i\cap U_j\neq\emptyset$ a holomorphic
map $f_{ij}\in\OO(U_i\cap U_j)$ with 
\begin{eqnarray*}
\C\times (U_i\cap U_j)\stackrel{f_i^{-1}}{\longrightarrow}
p^{-1}(U_i\cap U_j)\stackrel{f_j}{\longrightarrow}
\C\times (U_i\cap U_j),\\
(z,x)\mapsto (z+f_{ij}(x),x).
\end{eqnarray*}

(ii) $H^1(X,\OO)=0$ (e.g. \cite[26.1 Satz]{Fo77}).
This sheaf cohomology coincides with the \v{C}ech cohomology
$H^1((U_i)_{i\in I},\OO)$ with values in $\OO$
(e.g. \cite[12.8 Satz]{Fo77}).
Therefore any affine $\C$-bundle over $X$ is trivial, so 
isomorphic to the affine $\C$-bundle $\C\times X$.

(iii) Let $a:\C\times X\to\C\times X$ be an automorphism
of some finite order of $\C\times X$ as affine $\C$-bundle
over $X$, with a discrete set of fixed points in $X$,
such that $a$ acts as identity on the fiber of each fixed
point. Then $(\C\times X)/\langle a\rangle$ is again
an affine $\C$-bundle. The pull back of a trivialization
of this quotient bundle gives a new trivialization
of $\C\times X$ such that with respect to this new trivialization
$a$ acts nontrivially only on the base $X$,
$a:(z,x)\mapsto (z,a_2(x))$, where $a_2:X\to X$ is the
induced automorphism of the base $X$.

(iv) A {\it line bundle on $X$} 
\index{line bundle}
means a holomorphic rank 1 vector 
bundle. Any line bundle on $X$ is trivial, so isomorphic
to $\C\times X$ as vector bundle (e.g. \cite[30.3 Satz]{Fo77}).
Therefore any $\C^*$-bundle on $X$ is trivial, so isomorphic
to $\C^*\times X$.
\end{remarks}

\section{Partial compactifications of the F-manifolds}
\label{s9.3}

Theorem \ref{t9.3} gives all possible manifolds
$C_3^{\uuuu{e}/\{\pm 1\}^3}$ and $C_3^{S/\{\pm 1\}^3}$
as complex manifolds. 
The (long) proof describes also the coverings
$C_3^{\uuuu{e}/\{\pm 1\}^3}\to C_3^{S/\{\pm 1\}^3}$.
But the F-manifold structure and the Euler field
as well as possible partial compactifications are not treated
in Theorem \ref{t9.3}. They are subject of this section. 
Corollary \ref{t9.8} gives the Euler fields.
We refrain from writing down the multiplications in each
case. They can in principle be extracted from the construction
as quotients of $C_3^{univ}$ or as coverings of $C_3^{conf}$.

But we care about possible partial compactifications of 
$C_3^{\uuuu{e}/\{\pm 1\}^3}$ such that the F-manifold structure
extends. Lemma \ref{t9.7} says something conceptually about
such partial compactifications. Corollary \ref{t9.8} writes
them down in all cases and makes the Maxwell strata $\KK_2$ and
the caustics $\KK_3$ in them explicit. 

Especially interesting are the cases $A_3$ and $A_2A_1$.
In both cases $C_3^{\uuuu{e}/\{\pm 1\}^3}$ has a partial
compactification to $\C^3$ which is well known from 
singularity theory. To recover it starting from
$C_3^{\uuuu{e}/\{\pm 1\}^3}$ is especially in the case $A_3$
a bit involved. It requires in the case $A_3$ the notion
of a $\C^*$-orbibundle, which is given in the Remarks \ref{t9.5},
and the Example \ref{t9.6}.

\begin{remarks}\label{t9.5}
Let $\www{X}$ be a smooth complex curve, compact or not.
A {\it $\C^*$-orbibundle over $\www{X}$} 
\index{orbibundle} 
generalizes slightly
the notion of a $\C^*$-bundle over $\www{X}$. 
It is a complex two dimensional manifold $Y$ with a projection
$p:Y\to\www{X}$ and a $\C^*$-action which restricts to a 
transitive $\C^*$-action on each fiber and trivial stabilizer
on almost each fiber. There is a discrete subset 
$\Sigma\subset\www{X}$ such that the stabilizer of the 
$\C^*$-action is trivial on $p^{-1}(q)$ for 
$q\in\www{X}-\Sigma$ and nontrivial on $p^{-1}(q)$ for 
$q\in\Sigma$. 

Then $Y|_{\www{X}-\Sigma}$ is a $\C^*$-bundle over
$\www{X}-\Sigma$. For any $q\in \Sigma$ the restriction
$p^{-1}(\Delta)\subset Y$ for a small disk 
$\Delta\subset\www{X}$ with center $q$ and $\Delta\cap\Sigma
=\{q\}$ is isomorphic to one of the following models,
\index{$Y_{a,m}$} 
\begin{eqnarray*}
Y_{a,m}:=(\C^*\times\D)/\langle ((y_1,y_2)\mapsto 
(e^{2\pi i a/m}y_1,e^{2\pi i/m}y_2))\rangle,\\
\textup{with }\C^*\textup{-action }\C^*\times Y_{a,m}\to Y_{a,m},
\quad (t,[(y_1,y_2)])\mapsto [(ty_1,y_2)],\\
\textup{and projection }p_{a,m}:Y_{a,m}\to\D/\langle 
(y_2\mapsto e^{2\pi i /m}y_2)\rangle,\quad
[(y_1,y_2)]\mapsto [y_2]
\end{eqnarray*}
for some $m\in\Z_{\geq 2}$ and some $a\in\{1,2,...,m-1\}$
with $\gcd(a,m)=1$. 
The fiber over $[0]$ in this model is called 
{\it exceptional fiber of type $(a,m)$}. The stabilizer of 
the $\C^*$-action on it has order $m$. 

If $\www{X}$ is compact there is a notion of an Euler number
$\in\Q$ of a $\C^*$-orbibundle over $\www{X}$ which generalizes
the Euler number $\in\Z$ of a $\C^*$-bundle over $\www{X}$ 
\cite{OW71} (see also \cite{Do83}). 
If one takes out the center of an isolated quasihomogeneous
surface singularity one obtains a $\C^*$-orbibundle
over a smooth compact curve $\www{X}$ with negative Euler
number. Vice versa, any $\C^*$-orbibundle over a smooth
compact curve $\www{X}$ with negative Euler number arises in 
this way \cite{OW71}. 

We care about $\C^*$-orbibundles here mainly because of the case
$A_3$. The following example will be used in the proof of
the case $A_3$ in Corollary \ref{t9.8}.
\end{remarks}

\begin{example}\label{t9.6}
Consider $Y_1:=\C^2-\{(0,0)\}$ with coordinates 
$(z_{13},z_{14})$ and the $\C^*$-action
\begin{eqnarray*}
\C^*\times Y_1\to Y_1,\quad (t,z_{13},z_{14})\mapsto
(t^3 z_{13},t^2 z_{14}).
\end{eqnarray*}
It is a $\C^*$-orbibundle over $\www{X}=\P^1$ with projection
\begin{eqnarray*}
p_1:Y_1\to\P^1,\quad (z_{13},z_{14})\mapsto 
\frac{z_{14}^3}{z_{13}^2},
\end{eqnarray*}
and two exceptional fibers, the exceptional fiber 
$\{z_{14}=0\}$ over $0$ with stabilizer of order three 
and the exceptional fiber $\{z_{13}=0\}$ over $\infty$ of type 
$(1,2)$. 
We claim that the exceptional fiber $\{z_{14}=0\}$ over $0$ is 
of type $(1,3)$ and not of type $(2,3)$. 
We will show this below after studying the
$\C^*$-orbibundle $p_2:Y_2\to\P^1$ which is obtained from
$p_1:Y_1\to\P^1$ via
pull back with the $3:1$ branched covering
$b_2:\P^1\to\P^1, z\mapsto z^3$.

The $\C^*$-orbibundle $p_2:Y_2\to\P^1$ has only one 
exceptional fiber, the fiber of type $(1,2)$ over $\infty$. 
This $\C^*$-orbibundle is given by the following two charts
and their glueing.\\
First chart over $\C$:
\begin{eqnarray*}
p_2^{(1)}:Y_2^{(1)}=\C^*\times\C\to\C,&& (z_{12},z_{11})\mapsto 
z_{11},\\
\C^*\textup{-action}\quad  \C^*\times Y_2^{(1)}\to Y_2^{(1)},
&& (t,z_{12},z_{11})\mapsto (tz_{12},z_{11}),\\
Y_2^{(1)}\stackrel{3:1}{\longrightarrow} 
Y_1^{(1)}=p_1^{-1}(\C)=\C^*\times\C,&&
(z_{12},z_{11})\mapsto (z_{12}^3,z_{12}^2z_{11})=(z_{13},z_{14}).
\end{eqnarray*}
Second chart over $\P^1-\{0\}$:
\begin{eqnarray*}
p_2^{(2)}:Y_2^{(2)}=\C\times\C^*\to\P^1-\{0\},
&& (y_1,y_2)\mapsto \frac{y_2}{y_1^2},\\
\C^*\textup{-action}\quad  \C^*\times Y_2^{(2)}\to Y_2^{(2)},
&& (t,y_1,y_2)\mapsto (ty_1,t^2y_2),\\
Y_2^{(2)}\longrightarrow
Y_1^{(2)}=p_1^{-1}(\P^1-\{0\})=\C\times\C^*,&&
(y_1,y_2)\mapsto (y_1^3,y_2)=(z_{13},z_{14}),
\end{eqnarray*}
Glueing of the two charts:
\begin{eqnarray*}
Y_2^{(1)}-(p_2^{(1)})^{-1}(0)=\C^*\times\C^* &\to& 
Y_2^{(2)}-(p_2^{(2)})^{-1}(\infty)=\C^*\times\C^*,\\
(z_{12},z_{11})&\mapsto& (z_{12},z_{12}^2z_{11})=(y_1,y_2).
\end{eqnarray*}
The map $Y_2^{(1)}\to Y_1^{(1)}$, $(z_{12},z_{11})\mapsto
(z_{12}^3,z_{12}^2z_{11})$, shows that the chart $Y_1^{(1)}$
of the first $\C^*$-orbibundle $p_1:Y_1\to\P^1$ 
is a $\C^*$-orbibundle over $\C$ with exceptional fiber over 0
and that the $\C^*$-orbibundle is 
isomorphic to the model of type $(1,3)$ in Remark \ref{t9.4}
(v), namely
\begin{eqnarray*}
Y_1^{(1)}\cong (\C^*\times\C)/\langle ((z_{12},z_{11})\mapsto
(e^{2\pi i /3}z_{12},e^{2\pi i /3}z_{11}))\rangle.
\end{eqnarray*}

Consider a further $2:1$ branched covering 
$b_3:\P^1\to\P^1$, $\www{z}_{11}\to \www{z}_{11}^2=z_{11}$.
The pull back of $p_2:Y_2\to\P^1$ via $b_3$ gives a 
$\C^*$-bundle $p_3:Y_3\to\P^1$ over $\P^1$. 
We claim that its Euler number is $-1$, so it is isomorphic to 
$(\OO_{\P^1}(-1)-\textup{(zero section)})$. To prove this,
we consider two charts for $p_3:Y_3\to\P^1$ and their glueing.\\
First chart over $\C$:
\begin{eqnarray*}
p_3^{(1)}:Y_3^{(1)}=\C^*\times\C\to\C,
&& (z_{12},\www{z}_{11})\mapsto \www{z}_{11},\\
\C^*\textup{-action}\quad  \C^*\times Y_3^{(1)}\to Y_3^{(1)},
&& (t,z_{12},\www{z}_{11})\mapsto (tz_{12},\www{z}_{11}),\\
Y_3^{(1)}\stackrel{2:1}{\longrightarrow} 
Y_2^{(1)}=\C^*\times\C,&&
(z_{12},\www{z}_{11})\mapsto (z_{12},\www{z}_{11}^2)
=(z_{12},z_{11}).
\end{eqnarray*}
Second chart over $\P^1-\{0\}$:
\begin{eqnarray*}
p_3^{(2)}:Y_3^{(2)}=\C^*\times\C\to\P^1-\{0\},
&& (y_3,y_4)\mapsto \frac{1}{y_4},\\
\C^*\textup{-action}\quad  \C^*\times Y_3^{(2)}\to Y_3^{(2)},
&& (t,y_3,y_4)\mapsto (ty_3,y_4),\\
Y_3^{(2)}\longrightarrow
Y_2^{(2)}=\C\times\C^*,&&
(y_3,y_4)\mapsto (y_3y_4,y_3^2)=(y_1,y_2),
\end{eqnarray*}
Glueing of the two charts:
\begin{eqnarray*}
Y_3^{(1)}-(p_3^{(1)})^{-1}(0)=\C^*\times\C^* &\to& 
Y_3^{(2)}-(p_3^{(2)})^{-1}(\infty)=\C^*\times\C^*,\\
(z_{12},\www{z}_{11})&\mapsto& (z_{12}\www{z}_{11},
\frac{1}{\www{z}_{11}})=(y_3,y_4).
\end{eqnarray*}
The glueing map shows that $p_3:Y_3\to\P^1$ is the 
$\C^*$-bundle with Euler number $-1$. 
\end{example}

\begin{lemma}\label{t9.7}
Fix $m\in\Z_{\geq 2}$. 

(a) The following table defines four elements 
$\beta_1,\beta_2,\beta_3,\beta_4 \in \Br_3$, it gives their 
images $g_1,g_2,g_3,g_4$ in $PSL_2(\Z)$ under the 
homomorphism $\Br_3\to PSL_2(\Z)$ in Remark \ref{t4.15} (i),
and it defines four points $x_1,x_2,x_3,x_4\in\whh{\R}$,
\begin{eqnarray*}
\begin{array}{c|c|c|c|c}
i & 1 & 2 & 3 & 4 \\
\beta_i & \sigma_2^m & \sigma_1\sigma_2^m\sigma_1^{-1}
=\sigma_2^{-1}\sigma_1^m\sigma_2 & 
\sigma_1^{-1}\sigma_2^m\sigma_1=\sigma_2\sigma_1^m\sigma_2^{-1} 
& \sigma_1^m \\
g_i & [A_2^m] & [A_1A_2^mA_1^{-1}] = [A_2^{-1}A_1^mA_2] & 
[A_1^{-1}A_2^mA_1] = [A_2A_1^mA_2^{-1}] &[A_1^m] \\
x_i & 0 & -1 & 1 & \infty
\end{array} 
\end{eqnarray*}
The M\"obius transformation $g_i$ is parabolic with fixed point
$x_i$. The action of the braid $\beta_i$ on 
$\C\times\C\times\H\cong C_3^{univ}$ restricts to the action
of $g_i$ on $\H$. The quotient $\H/\langle g_i\rangle$
of that action is isomorphic to $\D^*=\D-\{0\}$. The quotient
$(\C\times\C\times\H)/\langle \beta_i\rangle$ is isomorphic
to $\C\times (\textup{an affine }\C\textup{-bundle over }\D^*)$.

(b) This quotient extends to 
$\C\times (\textup{an affine }\C\textup{-bundle over }\D)$
such that this is an F-manifold with Euler field and at each
point in the fiber over $0\in\D$ the germ is of type 
$I_2(m)A_1$ (compare Example \ref{t8.10} (ii)).

(c) (Concretization of a family of cases in part (b))
If $m$ is even then
\begin{eqnarray*}
(\C\times\C\times\H)/\langle \sigma_2^m\rangle 
& \stackrel{\cong}{\longrightarrow} & \C\times\C\times\D^* \\
{}[(z_1,\zeta,\tau)] &\mapsto & (z_1,\zeta,z_4) 
\quad\textup{with }z_4=z_4(\tau),
\end{eqnarray*}
and the partial compactification in part (b) is 
$\C\times\C\times\D$. 

(d) The semisimple F-manifold with Euler field
$\C\times\C^*\times (\C-\{0,1\})\cong C_3^{pure}$
extends to a semisimple F-manifold with Euler field
\begin{eqnarray*}
\C\times (\textup{the }\C^*\textup{-bundle over }\P^1
\textup{ which is }\OO_{\P^1}(-1)\textup{ minus the }0
\textup{-section})
\end{eqnarray*}
with Maxwell stratum the union of the fibers of 
$0,1,\infty\in\P^1$.
At each point in the Maxwell stratum the germ is of type
$I_2(2)A_1$ (compare Example \ref{t8.10} (ii)). 

(e) Consider $m\in\N$. The F-manifold with Euler field
\begin{eqnarray*}
\C\times\C^*\times\H \cong (\C\times\C\times\H)/\langle 
(\sigma^{mon})^m\rangle 
\cong C_3^{univ}/\langle (\sigma^{mon})^m\rangle
\end{eqnarray*}
extends to $\C\times\C\times\H$ if and only if $m\geq 2$.
If $m\geq 2$ then at each point of $\C\times\{0\}\times\H$
the germ of the F-manifold is irreducible.
\end{lemma}

{\bf Proof:}
(a) The equations 
$\sigma_2\sigma_1^m\sigma_2^{-1}=\sigma_1^{-1}\sigma_2^m\sigma_1$
and
$\sigma_2^{-1}\sigma_1^m\sigma_2=\sigma_1\sigma_2^m\sigma_1^{-1}$
are \eqref{4.14}. They show that $\beta_1,\beta_2,\beta_3$
and $\beta_4$ are all conjugate in $\Br_3$. 
Therefore also $g_1,g_2,g_3$ and $g_4$ are all conjugate in
$PSL_2(\Z)$. They are parabolic with the fixed points 
$x_1,x_2,x_3$ and $x_4$. The rest of part (a) is also obvious.

Part (d) is proved before the parts (b) and (c).

(d) $C_3^{pure}=\C^3-D_3^{pure}$ extends into the semisimple
F-manifold $\C^3$ of type $A_1^3$ with Maxwell stratum
$D_3^{pure}$, which consists of three hyperplanes which
intersect all three in the complex line 
$\{\uuuu{u}\in\C^3\,|\, u_1=u_2=u_3\}$. 
The isomorphism
$C_3^{pure}\cong \C\times\C^*\times(\C-\{0,1\})$ is the 
restriction to $C_3^{pure}$ of the blowing up of this complex 
line, 
\begin{eqnarray*}
\C\times \OO_{\P^1}(-1) &\longrightarrow& \C^3 \\
\cup & & \cup \\
\C\times \C\times\C & \longrightarrow & \C\times 
(\C\times\C -\{0\}\times\C^*) \\
\cup & & \cup \\
\C\times\C^*\times (\C-\{0,1\}) & \longrightarrow & 
\{(u_1+u_2+u_3,u_2-u_1,u_3-u_1)\in\C^3\,|\, \\
&& \hspace*{1cm}u_1\neq u_2\neq u_3\neq u_1\} \cong C_3^{pure}\\
(z_1,z_2,z_3) &\mapsto & (z_1,z_2,z_2z_3)\\
&& \hspace*{0.5cm}=(u_1+u_2+u_3,u_2-u_1,u_3-u_1)
\end{eqnarray*}
In Theorem \ref{t9.1} (a) and here only one of the two standard
charts of $\C\times \OO_{\P^1}(-1)$ 
is made explicit. This chart excludes the fiber over
$z_3=\infty\in\P^1$. The three hyperplanes in $D_3^{pure}$
give the fibers over $z_3=0,1,\infty\in\P^1$. 
Compare Remark \ref{t9.2} (ii). 

(b) and (c) Because $\beta_1,\beta_2,\beta_3$ and $\beta_4$ 
are all conjugate in $\Br_3$, it is for the proof of part (b)
sufficient to prove it for $\beta_1=\sigma_2^m$.
The action of $\sigma_2^m$ on 
$\C\times\C\times\H\cong C_3^{univ}$ is as follows,
\begin{eqnarray*}
\psi(\sigma_2^m):\ (z_1,\zeta,\tau)&\mapsto& (z_1,\zeta,
\frac{\tau}{m\tau+1})\quad\textup{for even }m,\\
\psi(\sigma_2^m):\ 
(z_1,\zeta,\tau)&\mapsto& (z_1,\zeta+\kappa(\tau+1),
\frac{\tau}{m\tau+1})\quad\textup{for odd }m.
\end{eqnarray*}
So, for even $m$ $\psi(\sigma_2^m)$ acts nontrivially only on
the third factor $\H$. Part (d) and this fact give the parts
(b) and (c) in the case $m=2$. 

For even $m$ the map 
$(\C\times\C\times\H)/\langle\sigma_2^m\rangle
\to (\C\times\C\times\H)/\langle\sigma_2^2\rangle$
is a branched covering of degree $\frac{m}{2}$, cyclically 
branched along the fiber over $[0]\in\H/\langle[A_2^2]\rangle$.
With the Examples \ref{t8.10} (ii) and \ref{t8.10} (v) 
(the case $I_2(m)A_1$) 
this shows for even $m$ part (b) for $\beta_1$ and part (c). 

For odd $m$ observe that $\kappa(\tau+1)$ has for $\tau\to 0$
along each hyperbolic line in $\H$ a limit, the limit
$\kappa(1)=0$. Therefore for odd $m$ we obtain a 2-valued map
$(\C\times\C\times\H)/\langle\sigma_2^m\rangle
\to (\C\times\C\times\H)/\langle\sigma_2^2\rangle$
which can be considered as a branched covering of the half-degree
$\frac{m}{2}$, cyclically branched along the fiber over
$[0]\in\H/\langle A_2^2\rangle$.
Again with the Examples \ref{t8.10} (ii) and \ref{t8.10} (v) 
this shows part (b) for $\beta_1=\sigma_2^m$. 

(e) For $m=1$ the multiplication in $\C\times\C^*\times\H
\cong C_3^{univ}/\langle \sigma^{mon}\rangle$ degenerates near
$\C\times\{0\}\times\H$ in the same way as the multiplication
in $\C\times\C^*\times(\C-\{0,1\})\cong C_3^{pure}$ 
degenerates near $\C\times\{0\}\times (\C-\{0,1\})$. 
Because of the denominator $z_2$ in the third summand 
on the right hand side of
\begin{eqnarray*}
\paa_2\circ\paa_2= \frac{2-2z_3+2z_3^2}{3}\paa_1 + 
\frac{1-2z_3}{3}\paa_2 + \frac{-z_3+z_3^2}{z_2}\paa_3,
\end{eqnarray*}
the multiplication does not extend holomorphically to
$\C\times\{0\}\times(\C-\{0,1\})$.

For $m\geq 2$ the (vertical) maps
\begin{eqnarray*}
\begin{CD}
C_3^{univ}/\langle (\sigma^{mon})^m\rangle  
@. \ \cong\ @. \C\times\C^*\times \H
@.  \hspace*{0.5cm}(z_1,z_{18},\tau) @. \\
@VVV @. @VVV @VVV @. \\
C_3^{univ}/\langle\sigma^{mon}\rangle 
@. \ \cong\ @. \C\times\C^*\times\H 
@.  \hspace*{0.5cm}(z_1,z_{18}^m,\tau)@. \ =(z_1,z_2,\tau)
\end{CD}
\end{eqnarray*}
are $m$-fold coverings.
Again, we can go from $\H$ to $\C-\{0,1\}$. 
The multiplication upstairs is obtained by pull back from the
multiplication downstairs via the $m$-fold covering,
\begin{eqnarray*}
\begin{CD}
@. \C\times\C^*\times (\C-\{0,1\})
@.  \hspace*{0.5cm}(z_1,z_{18},z_3)@. \\
@. @VVV @VVV @.\\
C_3^{pure} \cong\  @. \C\times\C^*\times(\C-\{0,1\}) 
@.  \hspace*{0.5cm}(z_1,z_{18}^m,z_3)@. \ =(z_1,z_2,z_3)
\end{CD}
\end{eqnarray*}
Write $\paa_{18}:=\frac{\paa}{\paa z_{18}}$. 
Because of $\paa_{18}\mapsto mz_{18}^{m-1}\paa_2$, the multiplication
is given by $e=3\paa_1$ and
\begin{eqnarray*}
\paa_{18}\circ\paa_{18}&=& m^2z_{18}^{2m-2}\frac{2-2z_3+2z_3^2}{3}\paa_1 
+ mz_{18}^{m-1}\frac{1-2z_3}{3}\paa_{18}\\
&+&  m^2z_{18}^{m-2}(-z_3+z_3^2)\paa_3,\\
\paa_{18}\circ\paa_3&=& mz_{18}^{m-1}\frac{-z_3+2z_{18}^mz_3}{3}\paa_1 
- z_{18}^m\frac{1}{3}\paa_{18} + mz_{18}^{m-1}\frac{-1+2z_3}{3}\paa_3,\\
\paa_3\circ\paa_3&=& \frac{2z_{18}^{2m}}{3}\paa_1 
+\frac{z_{18}^m}{3}\paa_3.
\end{eqnarray*}
The multiplication extends holomorphically to 
$\C\times\{0\}\times(\C-\{0,1\})$.
Moreover $\paa_{18}\circ\paa_3$ and 
$\paa_3\circ\paa_3$ vanish there, 
and also $\paa_{18}\circ\paa_{18}$ if $m\geq 3$.
If $m=2$ then $\paa_{18}\circ\paa_{18}=4(-z_3+z_3^2)\paa_3$
on $\C\times \{0\}\times (\C-\{0,1\})$. Therefore for any
$m\geq 2$ the algebra
$(T_p,\circ)$ for $p\in\C\times\{0\}\times (\C-\{0,1\})$ 
is irreducible, and thus the germ at $p$ of the F-manifold
is irreducible.
\hfill$\Box$

\bigskip
Recall that $C_3^{\uuuu{e}/\{\pm 1\}^3}$ is in all cases a 
semisimple F-manifold (with Euler field) with empty
Maxwell stratum. We consider the manifold in Theorem \ref{t9.3}
which is isomorphic to $C_3^{\uuuu{e}/\{\pm 1\}^3}$ with the
coordinates in the proof of Theorem \ref{t9.3}. 

Together Theorem \ref{t7.11}, Theorem \ref{t9.1} (a),
Theorem \ref{t9.3} and its proof, and Lemma \ref{t9.7}
allow in many cases to embed $C_3^{\uuuu{e}/\{\pm 1\}^3}$
into a bigger F-manifold $\oooo{C_3^{\uuuu{e}/\{\pm 1\}^3}}^p$ 
with nonempty Maxwell stratum
or nonempty caustic or both. This is subject of Corollary 
\ref{t9.8}.

\begin{corollary}\label{t9.8}
In many cases in Theorem \ref{t9.3}, 
$C_3^{\uuuu{e}/\{\pm 1\}^3}$ has a partial compactification, 
\index{partial compactification}
which we call $\oooo{C_3^{\uuuu{e}/\{\pm 1\}^3}}^p$, to an 
F-manifold with nonempty Maxwell stratum or nonempty 
caustic or both. 

(a) The exceptions are the cases $\GG_5,\GG_6,\GG_7,\GG_8$
and $\GG_9$. There $E=z_1\paa_1+\paa_\zeta$ and 
\begin{eqnarray*}
C_3^{\uuuu{e}/\{\pm 1\}^3}&=&\C\times\C\times\D^*
\quad\textup{in the case }\GG_5,\\
C_3^{\uuuu{e}/\{\pm 1\}^3}&=&C_3^{univ}
\quad\textup{in the cases }
\GG_6,\GG_7,\GG_8\textup{ and }\GG_9.
\end{eqnarray*}

(b) All cases except the cases in part (a) and the
cases $A_2A_1$ and $A_3$: The first table below
gives the manifold in Theorem \ref{t9.3} isomorphic to
$C_3^{\uuuu{e}/\{\pm 1\}^3}$ and the manifold isomorphic to
its extension $\oooo{C_3^{\uuuu{e}/\{\pm 1\}^3}}^p$. 
The second table gives the strata in the Maxwell stratum 
$\KK_2$ and in the caustic $\KK_3$ of 
$\oooo{C_3^{\uuuu{e}/\{\pm 1\}^3}}^p$. 
In the case $\HH_{1,2}$ at each point $p\in\KK_3$ 
the germ of the F-manifold is irreducible.
In all other cases with $\KK_3\neq\emptyset$, at each point
$p\in\KK_3$ the germ of the F-manifold is of type
$I_2(3)A_1$. 
\begin{eqnarray*}
\begin{array}{l|l|l|l}
 & \textup{sets} & 
C_3^{\uuuu{e}/\{\pm 1\}^3} 
& \oooo{C_3^{\uuuu{e}/\{\pm 1\}^3}}^p  \\ \hline 
\GG_1 & C_1\ (A_1^3)
& C_3^{pure} & \C^3 \\
\GG_1 & C_2\ (\HH_{1,2})  
& \C\times\C^*\times\H & \C\times\C\times\H \\
\GG_2 & C_4\ (\P^1A_1),C_5
& \C^2\times (\C-\{0,1\}) & \C^3 \\
\GG_4 & C_7 \ (\whh{A}_2) 
& \C^2\times \C^{***} & \C^3 \\
\GG_{10} & C_{17} 
& \C\times\C\times\D^* & \C\times\C\times\D \\
\GG_{11} & C_{18}
& \C\times\C\times\D^* & \C\times\C\times\D \\
\GG_{12} & C_{19} 
& \C\times\C\times\D^{**} & \C\times\C\times\D \\
\GG_{13} & C_{20} 
& \C\times\C\times\D^{**} & \C\times\C\times\D \\
\GG_{14} & C_{21},\ C_{22}
& \C\times\C\times \D^* & \C\times\C\times\D 
\end{array}
\end{eqnarray*}

\begin{eqnarray*}
\begin{array}{l|l|l|l}
 & \textup{sets} & 
\KK_2 & \KK_3  \\ \hline 
\GG_1 & C_1\ (A_1^3)
& D_3^{pure} & - \\
\GG_1 & C_2\ (\HH_{1,2})  
& - & \C\times\{0\}\times\H\\
\GG_2 & C_4\ (\P^1A_1),C_5
& \C^2\times \{0,1\} & - \\
\GG_4 & C_7 \ (\whh{A}_2) 
& - & \C^2\times\{z_5\,|\, z_5^3=1\} \\
\GG_{10} & C_{17} 
& \C\times\C\times\{0\} &  - \\
\GG_{11} & C_{18}
& \C\times\C\times\{0\} &  - \\
\GG_{12} & C_{19} 
& \C\times\C\times\{0\} & \C\times\C\times\{\frac{1}{2}\} \\
\GG_{13} & C_{20} 
& - & \C\times\C\times\{0,\frac{1}{2}\} \\
\GG_{14} & C_{21},\ C_{22}
& - & \C\times\C\times\{0\} 
\end{array}
\end{eqnarray*}
The Euler field is
\begin{eqnarray*}
\begin{array}{l|l|l|l}
 & A_1^3 & \HH_{1,2} & (\GG_2\,\&\, C_4,C_5), \GG_4, 
 \GG_{10},\GG_{11},\GG_{12},\GG_{13},\GG_{14} \\ \hline 
 E & \sum_{i=1}^3 u_ie_i & z_1\paa_1+\frac{1}{4}z_5\paa_5 & 
 z_1\paa_1+\paa_\zeta 
\end{array}
\end{eqnarray*}

(c) The case $A_2A_1$: We use the second presentation
in Theorem \ref{t9.1}.
\begin{eqnarray*}
C_3^{\uuuu{e}/\{\pm 1\}^3}&\cong& 
\C\times\{(z_8,z_9)\in\C^*\times\C\,|\, z_8^3-4z_9^2\neq 0\},\\
\oooo{C_3^{\uuuu{e}/\{\pm 1\}^3}}^p&\cong& \C^3,\\
\KK_2&\cong& 
\C\times\{(z_8,z_9)\in\C\times\C\,|\, z_8^3-4z_9^2= 0\},\\
\KK_3&\cong& \C\times\{0\}\times\C,\\
E&=& z_1\paa_1+\frac{2}{3}z_8\paa_8+z_9\paa_9.
\end{eqnarray*}
At each point $p\in\KK_3$ the germ of the F-manifold is of
type $I_2(3)A_1$. 

(d) The case $A_3$: We use the second presentation 
in Theorem \ref{t9.1}. 
\begin{eqnarray*}
C_3^{\uuuu{e}/\{\pm 1\}^3}&\cong& 
\C\times\{(z_{13},z_{14})\in\C^*\times\C\,|\, z_{14}^3-z_{13}^2
\neq 0\},\\
\oooo{C_3^{\uuuu{e}/\{\pm 1\}^3}}^p&\cong& \C^3,\\
\KK_2&\cong& \C\times\{0\}\times\C,\\
\KK_3&\cong& 
\C\times\{(z_{13},z_{14})\in\C\times\C\,|\,
z_{14}^3-z_{13}^2= 0\},\\
E&=& z_1\paa_1+\frac{3}{4}z_{13}\paa_{13}
+\frac{1}{2}z_{14}\paa_{14}.
\end{eqnarray*}
At each point $p\in \KK_3-\KK_2 =\KK_3-\C\times\{(0,0)\}$
the germ of the F-manifold is of type $I_2(3)A_1$. 
At each point $p\in \C\times\{(0,0)\}=\KK_3\cap\KK_2$ 
it is irreducible.

\begin{figure}[H]
\includegraphics[width=0.8\textwidth]{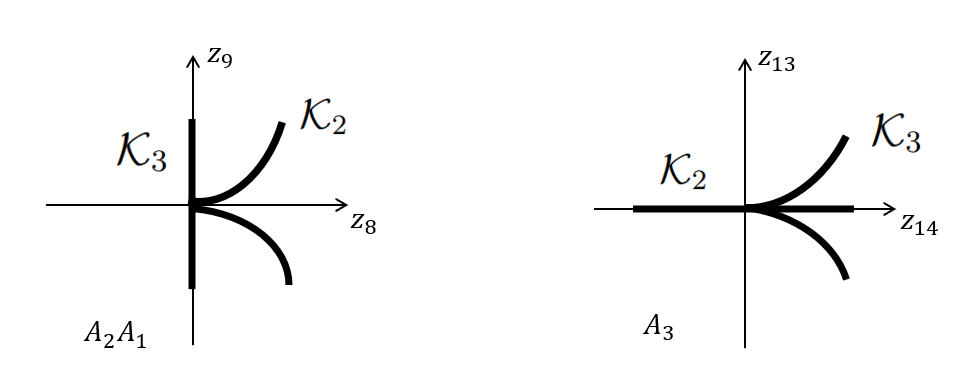}
\caption[Figure 9.3]{In the cases $A_2A_1$ (left) and $A_3$ 
(right) the restrictions of $\KK_2$ and $\KK_3$ to a 
hyperplane $\{z_1=\textup{constant}\}$ in the F-manifold
$\C^3$}
\label{Fig:9.3}
\end{figure}

\end{corollary}

{\bf Proof:}
(a) Here nothing has to be proved.
Though a remark on the cases $\GG_5\,\&\, C_8,C_9$: 
Here $(\Br_3)_{\uuuu{e}/\{\pm 1\}^3}=\langle (\sigma^{mon})^2
\sigma_1^{-1}\sigma_2^{l^2-4}\sigma_1\rangle$.
Without the factor $(\sigma^{mon})^2$ in the generator, 
Lemma \ref{t9.7} (b) would apply with $m=l^2-4$, and it 
would lead to an extension of the F-manifold to 
$\C\times\C\times\D$ with caustic $\C\times\C\times\{0\}$. 

(b) The case $\GG_1\,\&\, C_1\ (A_1^3)$: See
Example \ref{t8.10} (i).

The case $\GG_1\,\&\, C_2\ (\HH_{1,2})$: Lemma \ref{t9.7}
(e) applies with $m=2$. 

The cases $\GG_2\,\&\, C_4\ (\P^1A_1),C_5$: 
Recall $(\Br_3)_{\uuuu{e}/\{\pm 1\}^3}=\langle \sigma_2^2,
\sigma_1\sigma_2^2\sigma_1^{-1}\rangle$ and consider its
action on the third factor $\H$ of $\C\times\C\times\H
\cong C_3^{univ}$. The point $z_3=0$ in the closure of the
quotient $\C-\{0,1\}$ of $\H$ is the image of the
fixed point $0\in\whh{\R}$ of $\sigma_2^2$, 
the point $z_3=1$ is the image of the fixed point 
$-1\in\whh{\R}$ of $\sigma_1\sigma_2^2\sigma_1^{-1}$. 
By Lemma \ref{t9.7} (b)
the F-manifold extends to the fibers over $0$ and $1$ of 
$\C\times(\textup{an affine }\C\textup{-bundle over }\C)$
and has there its Maxwell stratum.
By Remark \ref{t9.4} (ii) one can choose the affine $\C$-bundle
over $\C$ to be $\C\times\C$. 

The case $\GG_4\,\&\, C_7\ (\whh{A}_2)$. 
Recall $(\Br_3)_{\uuuu{e}/\{\pm 1\}^3}=\langle \sigma_1^3,
\sigma_2^3, \sigma_2\sigma_1^3\sigma_2^{-1}\rangle$
and consider its action on the third factor $\H$ of 
$\C\times\C\times\H\cong C_3^{univ}$. The points $z_5$ 
with $z_5^3=1$ in the closure of the quotient $\C^{***}$ 
of $\H$ are the images of the fixed points $\infty,0$ and $1\in\whh{\R}$ of $\sigma_1^3,\sigma_2^3$ and 
$\sigma_2\sigma_1^3\sigma_2^{-1}$. By Lemma \ref{t9.7} (b)
the F-manifold extends to the fibers over $z_5$ with $z_5^3=1$ of
$\C\times(\textup{an affine }\C\textup{-bundle over }\C)$
and has there its caustic $\KK_3$, and at each point $p\in\KK_3$
the germ of the F-manifold is of type $I_2(3)A_1$. 
By Remark \ref{t9.4} (ii) one can choose the affine $\C$-bundle
over $\C$ to be $\C\times\C$. 

A remark on this case: Observe
$$(\sigma_1^3)(\sigma_2\sigma_1^3\sigma_2^{-1})(\sigma_2^3)
=(\sigma^{mon})^2(\sigma_2^{-1}\sigma_1^3\sigma_2)^{-1}.$$
Therefore $(\Br_3)_{\uuuu{e}/\{\pm 1\}^3}$ contains
$(\sigma^{mon})^{-2}\sigma_2^{-1}\sigma_1^3\sigma_2$. 
The point $z_5=\infty$ is the image of the fixed point of 
$A_2^{-1}A_1^3A_2$ which is the image in $PSL_2(\Z)$ of this
braid. 
But Lemma \ref{t9.7} (b) does not apply to $z_5=\infty$, 
because $(\Br_3)_{\uuuu{e}/\{\pm 1\}^3}$ does not contain 
$\sigma_2^{-1}\sigma_1^3\sigma_2$.

The cases $\GG_{10},\GG_{11},\GG_{12},\GG_{13},\GG_{14}$:
Recall from Theorem \ref{t7.11}
\begin{eqnarray*}
\begin{array}{l|l|l|l|l}
 & \GG_{10}\textup{ and }\GG_{11} & \GG_{12} & \GG_{13} & 
\GG_{14} \\
(\Br_3)_{\uuuu{e}/\{\pm 1\}^3} & \langle \sigma_2^2\rangle &
\langle \sigma_2^2,\sigma_2\sigma_1^3\sigma_2^{-1}\rangle & 
\langle \sigma_2^3,\sigma_2\sigma_1^3\sigma_2^{-1}\rangle &
\langle \sigma_2^3\rangle
\end{array}
\end{eqnarray*}
The point $0\in\D$ is the image of the fixed point $0\in\whh{\R}$
of $\sigma_2^2$ respectively $\sigma_2^3$, the point  
$\frac{1}{2}\in\D$ is in the cases $\GG_{12}$ and $\GG_{13}$
the image of the fixed point $-1\in\whh{\R}$ of the action of
$\sigma_2\sigma_1^3\sigma_2^{-1}$ on the third factor $\H$
of $\C\times\C\times\H\cong C_3^{univ}$. 
Lemma \ref{t9.7} (b) applies in all cases and yields an 
extension to an F-manifold 
$\C\times (\textup{an affine }\C\textup{-bundle over }\D)$
with the claimed Maxwell stratum and caustic. 
By Remark \ref{t9.4} (ii) one can choose the affine $\C$-bundle
over $\D$ to be $\C\times\D$. 

On the Euler field: For the case $A_1^3$ see Example \ref{t8.10}
(i).

In the case $\HH_{1,2}$ the moduli spaces is 
$C_3^{\uuuu{e}/\{\pm 1\}^3}\cong
\C\times\C^*\times\H$ with coordinates $(z_1,z_5,\tau)$
with $z_5=e^{\zeta/4}$. Therefore the Euler field 
$z_1\paa_1+\paa_\zeta$ on $C_3^{univ}$ pulls down to
$E=z_1\paa_1+\frac{1}{4}z_5\paa_5$.

In all other cases no power of $\sigma^{mon}$ is in 
$(\Br_3)_{\uuuu{e}/\{\pm 1\}^3}$, and 
\begin{eqnarray*}
C_3^{\uuuu{e}/\{\pm 1\}^3}&\cong& \C\times 
(\textup{an affine }\C\textup{-bundle over a noncompact curve }
C) \\ 
&\cong& \C\times\C\times C\quad\textup{with coordinates }
(z_1,\zeta,z_{3/11/16}).
\end{eqnarray*}
Therefore the Euler field $z_1\paa_1+\paa_\zeta$ on
$C_3^{univ}$ pulls down to $E=z_1\paa_1+\paa_\zeta$ on
$C_3^{\uuuu{e}/\{\pm 1\}^3}$.

(c) Recall the discussion of the case $A_2A_1$ in the proof of
Theorem \ref{t9.3}. Because $z_6=e^{\zeta/3}$, the Euler field
on $\C\times\C^*\times (\C-\{\pm\frac{1}{2}\})$ with coordinates
$(z_1,z_6,z_7)$ is $z_1\paa_1+\frac{1}{3}z_6\paa_6$.
Under the $2:1$ map
\begin{eqnarray*}
\C\times\C^*\times (\C-\{\pm\frac{1}{2}\})
&\to& 
\C\times \{(z_8,z_9)\in\C^*\times\C\,|\, z_8^3-4z_9^2=0\}\\
(z_1,z_6,z_7)&\mapsto& (z_1,z_6^2,z_6^3z_7)
=(z_1,z_8,z_9),
\end{eqnarray*}
it becomes
\begin{eqnarray*}
E=z_1\paa_1+\frac{1}{3}z_6(2z_6\paa_8+3z_6^2z_7\paa_9)
=z_1\paa_1+\frac{2}{3}z_8\paa_8+z_9\paa_9.
\end{eqnarray*}

Consider the isomorphism
\begin{eqnarray*}
C_3^{pure}\cong \C\times\C^*\times(\C-\{\pm\frac{1}{2}\})&\to&
\C\times\{(z_2,z_9)\in\C^*\times\C\,|\, z_2^2-4z_9^2=0\},\\
(u_1,u_2,u_3)\hspace*{2cm} 
(z_1,z_2,z_7) &\mapsto& (z_1,z_2,z_2z_7)
=(u_1+u_2+u_3,\\ && 
u_2-u_1,(u_3-u_1)-\frac{1}{2}(u_2-u_1)).
\end{eqnarray*}
The multiplication in the manifold on the right side extends
to $\C^3$ with Maxwell stratum the three hyperplanes
\begin{eqnarray*}
\C\times\{0\}\times\C\cup \C\times\{(z_2,z_9)\in\C^2\,|\, 
z_2^2-4z_9^2=0\}.
\end{eqnarray*}
The two-valued map
\begin{eqnarray*}
\C^3\dashrightarrow\C^3,\quad (z_1,z_8,z_9)\mapsto 
(z_1,z_8^{3/2},z_9) =(z_1,z_2,z_9)
\end{eqnarray*}
is outside of the branching locus $\{z_8=0\}$ locally an 
isomorphism of F-manifolds, and along 
$\{z_8=0\}$ it is branched of order $3/2$. 

Therefore the F-manifold $\C^3$ on the left side with 
coordinates $(z_1,z_8,z_9)$ has Maxwell stratum 
$\KK_2=\C\times\{(z_8,z_9)\in\C^2\,|\, z_8^3-4z_9^2=0\}$
and caustic $\KK_3=\C\times\{0\}\times\C$, and at each point
of the caustic the germ of the F-manifold is of type
$I_2(3)A_1$. 

(d) Recall the discussion of the case $A_3$ in the proof of
Theorem \ref{t9.3}. Because $z_{12}=e^{\zeta/4}$, 
the Euler field on $\C\times\C^*\times \C^{***}$ with coordinates
$(z_1,z_{12},z_{11})$ is $z_1\paa_1+\frac{1}{4}z_{12}\paa_{12}$.
Under the $3:1$ map
\begin{eqnarray*}
\C\times\C^*\times \C^{***}
&\to& 
\C\times \{(z_{13},z_{14})\in\C^*\times\C\,|\, z_{14}^3
-z_{13}^2\neq 0\}\\
(z_1,z_{12},z_{11})&\mapsto& (z_1,z_{12}^3,z_{12}^2z_{11})
=(z_1,z_{13},z_{14}),
\end{eqnarray*}
it becomes
\begin{eqnarray*}
E=z_1\paa_1+\frac{1}{4}z_{12}(3z_{12}^2\paa_{13}
+2z_{12}z_{11}\paa_{14})
=z_1\paa_1+\frac{3}{4}z_{13}\paa_{13}+\frac{1}{2}z_{14}\paa_{14}.
\end{eqnarray*}

In the proof of Theorem \ref{t9.3} the quotient of 
$C_3^{univ}$ by
\begin{eqnarray*}
\langle\sigma_1^3,\sigma_2^{-1}\sigma_1^3\sigma_2,\sigma_2^3,
(\sigma^{mon})^4\rangle
\stackrel{1:3}{\subset} (\Br_3)_{\uuuu{e}/\{\pm 1\}^3}
=\langle (\sigma_1\sigma_2)^4,\sigma_1^3\rangle
\end{eqnarray*}
was constructed as
\begin{eqnarray*}
\C\times\C^*\times\C^{***}\quad\textup{with coordinates }
(z_1,z_{12},z_{11})\\ 
\textup{ with } z_{12}=e^{\zeta/4},z_{11}=z_{11}(\tau).
\end{eqnarray*}
The three values $z_{11}\in\C$ with $z_{11}^3=1$ are the images
of the parabolic fixed points $\infty,-1,0\in\whh{\R}$
of the action of $\sigma_1^3,\sigma_2^{-1}\sigma_1^3\sigma_2,
\sigma_2^3$ on $\H$. 

By Lemma \ref{t9.7} (b) and Remark \ref{t9.4} (iv) 
the F-manifold $\C\times\C^*\times\C^{***}$ extends to 
$\C\times\C^*\times\C$ with caustic 
$\C\times\C^*\times\{z_{11}\,|\, z_{11}^3=1\}$, and at each
point $p$ in the caustic the germ of the F-manifold is
of type $I_2(3)A_1$. Under the extension
\begin{eqnarray*}
\C\times\C^*\times\C &\to& \C\times\C^*\times\C\\
(z_1,z_{12},z_{11}) &\mapsto& (z_1,z_{12}^3,z_{12}^2z_{11})
=(z_1,z_{13},z_{14})
\end{eqnarray*}
of the $3:1$ map above this caustic maps to the complex surface
\begin{eqnarray*}
\C\times\{(z_{13},z_{14}\in\C^*\times\C\,|\, 
z_{14}^3-z_{13}^2=0\},
\end{eqnarray*}
and the F-manifold extends to $\C\times\C^*\times\C$
with coordinates $(z_1,z_{13},z_{14})$ and with this surface
as caustic. 

We claim that the F-manifold
also extends to a fiber over $z_{14}=\infty$. 
Though proving this and understanding how this fiber glues in requires more care. The stabilizer 
$(\Br_3)_{\uuuu{e}/\{\pm 1\}^3}$ contains
$(\sigma^{mon})^4$ and
\begin{eqnarray*}
(\sigma_2^3)(\sigma_2^{-1}\sigma_1^3\sigma_2)(\sigma_1^3)
= (\sigma^{mon})^2 (\sigma_2\sigma_1^3\sigma_2^{-1})^{-1}
\end{eqnarray*}
so also $\sigma_2\sigma_1^6\sigma_2^{-1}$. 
The parabolic fixed point $1\in\whh{\R}$ of the action of
$\sigma_2\sigma_1^6\sigma_2^{-1}$ corresponds to $z_{11}=\infty$.
By Lemma \ref{t9.7} (b) the F-manifold $\C\times\C^*\times\C$
with coordinates $(z_1,z_{12},z_{11})$ extends to an F-manifold
with fiber over $z_{11}=0$, namely
\begin{eqnarray*}
\C\times (\textup{a certain }\C^*\textup{-orbibundle over }
\P^1)=:\C\times \www{Y}_2,
\end{eqnarray*}
and for a point in the fiber over $z_{11}=\infty$ the germ
of the F-manifold is of type $I_2(6)A_1$.

The fiber over $z_{11}=\infty$ is exceptional of type
$(1,2)$ (in the sense of Remark \ref{t9.5}),
because 
$$(\sigma^{mon})^4, 
(\sigma^{mon})^2\sigma_2\sigma_1^3\sigma_2^{-1},
\sigma_2\sigma_1^6\sigma_2^{-1}
\in (\Br_3)_{\uuuu{e}/\{\pm 1\}^3}.$$

The pull back of $\www{p}_2:\www{Y}_2\to\P^1$ via the 
$2:1$ branched covering $b_3:\P^1\to\P^1$, 
$\www{z}_{11}\to \www{z}_{11}^2=z_{11}$, is a $\C^*$-bundle
$\www{Y}_3$ with projection $\www{p}_3:\www{Y}_3\to\P^1$. 
We claim that its Euler number is $-1$. 

To see this, compare it with $C_3^{pure}$ as $\C^*$-bundle
with the standard action of $\C^*$. This is the
tautological bundle $\OO_{\P^1}(-1)$ minus its zero section,
so it has Euler number $-1$. Because of
the $2:1$ branched covering $b_3:\P^1\to\P^1$ and because
of $[PSL_2(\Z):\oooo{\Gamma(2)}]=6$ and
$[PSL_2(\Z):\oooo{\Gamma(3)}]=12$, the base space $\P^1$
of $\www{Y}_3$ is as a quotient of $\H\cup\{\textup{cusps}\}$ 
four times as large as the base space $\P^1$ of $C_3^{pure}$ as $\C^*$-bundle.
On the other hand, to obtain $C_3^{pure}$ from 
$C_3^{univ}$, we divided out $\sigma^{mon}$.
To obtain $\www{Y}_3$, we divided out only $(\sigma^{mon})^4$.
Therefore the Euler number of $\www{Y}_3$ is also $-1$. 

Now Example \ref{t9.6} applies. One sees
$\www{Y}_3=Y_3$ and $\www{Y}_2=Y_2$. 
The automorphism of order three in the proof of Theorem 
\ref{t9.3}, which $(\sigma_1\sigma_2)^4$
induces on the chart $Y_2^{(1)}$, 
$(z_{12},z_{11})\mapsto 
(e^{-2\pi i /3}z_{12},e^{-2\pi i/3}z_{11})$, 
extends to $Y_2$. It restricts because of $y_2=z_{12}^2z_{11}$ 
to the identity on the fiber over $z_{11}=\infty$ 
which is the fiber $\{0\}\times\C$ in the chart $Y_2^{(2)}$.
The quotient is $Y_1$. The quotient map is a $3:1$ 
branched covering with branching locus the fiber over
$z_{11}=\infty$ and the image $\{0\}\times\C^*\subset Y_1$
of this fiber. Therefore at each point in 
$\C\times (\textup{this image})=\C\times\{0\}\times\C^*$
the germ of the F-manifold is of type $I_2(2)A_1$.
Therefore $\C\times\{0\}\times\C^*$ 
is the Maxwell stratum of $\C\times Y_1$ as F-manifold.

Finally, the complement of $\C\times Y_1$ in $\C^3$ has
codimension 2. Therefore the multiplication extends to this
complement holomorphically, so the F-manifold extends to
$\C^3$ with coordinates $(z_1,z_{13},z_{14})$. 
Because this codimension 2 complement is in the closure
of the caustic of the F-manifold $\C\times Y_1$,
it is part of the caustic of $\C^3$. Part (d) is proved.
\hfill$\Box$

\chapter{Isolated hypersurface singularities}\label{s10}
\setcounter{equation}{0}
\setcounter{figure}{0}

An isolated hypersurface singularity 
gives rise to a unimodular bilinear lattice $(H_\Z,L)$ of some rank
$n\in\N$ and a $\Br_n\ltimes\{\pm 1\}^n$-orbit $\BB^{dist}$ of 
triangular bases. This chapter reports about these
structures and gives references. 
It does not contain new results.

The theory of isolated hypersurface singularities came
to life in the 1960ies. Milnor's book \cite{Mi69}
was a crucial early step. Soon several groups
around the world joined, especially Arnold and his 
students, and also the first author's teacher Brieskorn.

Nowadays there are several thorough and extensive 
presentations of the theory of isolated hypersurface
singularities, foremost the two volumes \cite{AGV85}
\cite{AGV88}, but also the more recent books
\cite{AGLV98} and \cite{Eb01}. A few years ago
Ebeling wrote a survey \cite{Eb20} on distinguished
bases for isolated hypersurface singularities.
But the sections below offer quite a lot of material
which is not covered in that survey. 

Section \ref{s10.1} introduces basic structures
induced by one isolated hypersurface singularity,
especially a unimodular bilinear lattice $(H_\Z,L)$,
but also the notion of a universal unfolding.
The base space $\MM$ of such an unfolding has to
be compared with a space $C_n^{\uuuu{e}/\{\pm 1\}^n}$.
The section finishes with a loose description of Arnold's
classification of singularities with modality
$\leq 2$.

Section \ref{s10.2} equips the pair $(H_\Z,L)$
with a natural $\Br_n\ltimes\{\pm 1\}^n$ orbit
$\BB^{dist}$ of distinguished bases.
Here the $\Z$-lattice bundles from section
\ref{s8.5} come from geometry and induce
this set $\BB^{dist}$. 
These triples $(H_\Z,L,\BB^{dist})$ are
rather special. The only triples in rank 2
and rank 3 which appear here are those of
$A_2$ and $A_3$. Section \ref{s10.2}
points at techniques and cites results
on the triples from singularities. 
It ends with the best studied cases,
the simple singularities and the 
simple elliptic singularities. 
It remains an interesting open problem to 
find properties which distinguish the 
triples $(H_\Z,L,\BB^{dist})$ for singularities
from all other triples.

Section \ref{s10.3} discusses results of the
first author and coauthors on the groups $G_\Z$ 
for singularities. A crucial notion is that
of an {\it Orlik block}. 
The most concrete results are those on the
singularities with modality $\leq 2$. 

Section \ref{s10.4} explains classical results
from the 1980ies on the monodromy groups 
and the vanishing cycles for singularities.
The situation for the even monodromy group $\Gamma^{(0)}$
is very satisfying, thanks to results of Ebeling.
In almost all cases it is determined by the pair
$(H_\Z,I^{(0)})$. The odd case is more complicated.
There the main work was done by  
Wajnryb, Chmutov and Janssen. 
It is also subject of section \ref{s10.4}. 

Section \ref{s10.5} comes to the few moduli spaces
$C_n^{\uuuu{e}/\{\pm 1\}^n}$ for singularities,
which had been studied, namely those for the
simple singularities and the simple elliptic singularities.
They turn out to be the complements of Maxwell stratum
and caustic in base spaces of certain global unfoldings.
In the case of the simple singularities, this is 
essentially due to Looijenga and Deligne 1974.
Hertling and Roucairol generalized it to
the simple elliptic singularities. 
The section offers key ideas of the proofs.

\section{Topology, Milnor lattice, unfolding, classification}\label{s10.1}

\begin{definition}\label{t10.1}
(E.g. \cite{Mi69}\cite{Lo84}\cite{AGV88}\cite{Eb01})
\index{isolated hypersurface singularity}\index{singularity}
An isolated hypersurface singularity is a holomorphic
function germ $$f:(\C^{m+1},0)\to(\C,0)$$ (with $m\in\Z_{\geq 0}$) 
with $f(0)=0$ and with an isolated singularity at 0, 
i.e. the Jacobi ideal
$J_f:=\bigl(\frac{\paa f}{\paa z_0},...,\frac{\paa f}{\paa z_m}\bigr)
\subset \C\{z_0,...,z_m\}$ has an isolated zero at 0.
Its Milnor number $n$ \index{Milnor number} is
\begin{eqnarray}\label{10.1}
n:=\dim\C\{z_0,...,z_m\}/J_f\in\N.
\end{eqnarray}
\end{definition}

Usually the Milnor number is called $\mu$, and the number of
variables is called $n$ \cite{AGV88} or $n+1$ \cite{Eb01}.
In order to be consistent with the other chapters of this book,
here the Milnor number is called $n$ and the number of variables
is called $m+1$. 

\begin{figure}
\includegraphics[width=0.5\textwidth]{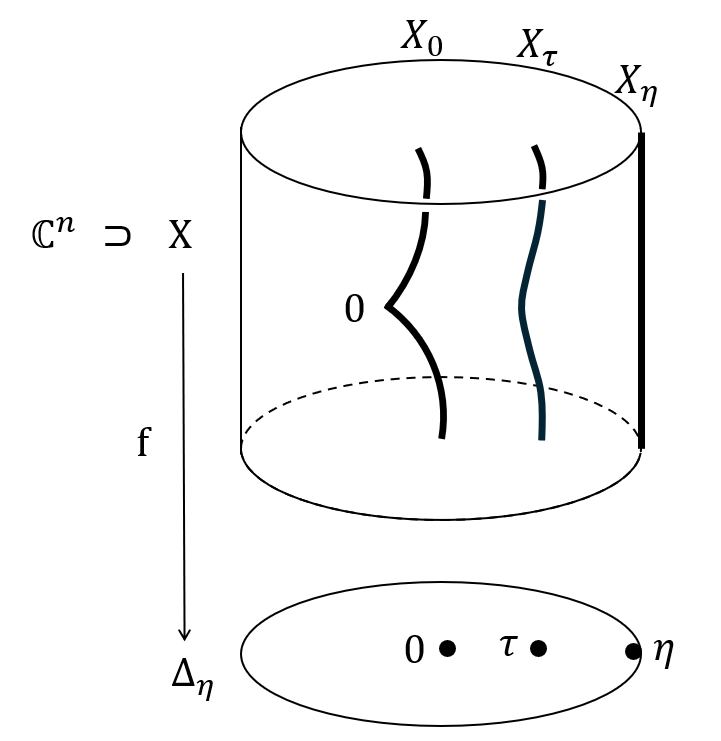}
\caption[Figure 10.1]{A schematic picture of a representative
$f:X\to\Delta_\eta$ of a singularity 
$f:(\C^{m+1},0)\to(\C,0)$}
\label{Fig:10.1}
\end{figure}

The first thing to do is to construct 
{\it a good representative} 
\index{good representative of $f$}
of $f$. One can choose and chooses 
a small positive number $\varepsilon>0$ such that the boundary 
of each ball $B_{\www{\varepsilon}}\subset \C^{m+1}$ 
with center at 0 and radius $\www{\varepsilon}\in(0,\varepsilon]$ 
intersects the fiber $f^{-1}(0)$ transversally.
Then one chooses a small positive number $\eta>0$ such that
each fiber $f^{-1}(\tau)$ for $\tau\in\Delta_\eta$
intersects the boundary of $B_\varepsilon$ transversally.
Here $\Delta_\eta:=\{\tau\in\C\, |\, |\tau|<\eta\}\subset\C$ 
is the disk of radius $\eta$ around 0. Then 
\begin{eqnarray}\label{10.2}
f:X\to\Delta_\eta\quad\textup{with}\quad
X:=B_\varepsilon\cap f^{-1}(\Delta_\eta)
\end{eqnarray}
is a good representative
of $f$ with fibers $X_\tau=f^{-1}(\tau)$ for $\tau\in\Delta_\eta$.
The restriction $f:X-X_0\to \Delta_\eta-\{0\}$ is a locally
trivial $C^\infty$ fiber bundle.
Only the 0-fiber $X_0$ is singular, and its only singularity
is at 0. The fibers $X_\tau$ for $\tau\in\Delta_\eta-\{0\}$ 
are called {\it Milnor fibers}. 
That all this works and is true can be found in 
\cite{Mi69}, \cite{Lo84}, \cite{AGV88} or \cite{Eb01}.
See Figure \ref{Fig:10.1} for a schematic picture.

The references \cite[(5.11)]{Lo84}, \cite[ch. 2]{AGV88} 
or \cite{Eb01} tell also the following.
Our signs and conventions are consistent with those in
\cite[ch. 2]{AGV88}. 

\begin{theorem}\label{t10.2}
Let $f:(\C^{m+1},0)\to(\C,0)$ be an isolated hypersurface singularity
with Milnor number $n\in\N$. Let $f:X\to\Delta_\eta$ be a good
representative as above. 

(a) Each regular fiber
$X_\tau$ with $\tau\in\Delta_\eta-\{0\}$ is homotopy equivalent to
a bouquet of $n$ $m$-spheres, so that its only 
(reduced if $m=0$) homology is the middle homology and that is
$$H_m^{(red)}(X_\tau,\Z)\cong\Z^n.$$
This $\Z$-module is called a {\sf Milnor lattice}.
\index{Milnor lattice} 

(b) For $\zeta\in S^1$ there is a canonical isomorphism 
\begin{eqnarray}\label{10.3}
H_m^{(red)}(X_{\zeta\eta},\Z)\cong H_{m+1}(X,X_{\zeta\eta},\Z)
=:Ml(f,\zeta)
\end{eqnarray}
of the Milnor lattice over a point $\zeta\eta\in\paa\Delta_\eta$
to the relative homology group with boundary in
$X_{\zeta\eta}$. The union $\bigcup_{\zeta\in S^1}Ml(f,\zeta)$
is a flat bundle over $S^1$ of $\Z$-lattices of rank $n$.

(c) An intersection form for Lefschetz thimbles
\index{intersection form for Lefschetz thimbles}
\index{Lefschetz thimble}
is well defined on relative homology groups with different 
boundary parts. It is for any $\zeta\in S^1$ a $(-1)^{m+1}$
symmetric unimodular bilinear form
\begin{eqnarray}\label{10.4}
I^{Lef}:Ml(f,\zeta)\times Ml(f,-\zeta)\to\Z
\end{eqnarray}
Let $\gamma_\pi$ be the
isomorphism $Ml(f,\zeta)\to Ml(f,-\zeta)$ by flat shift in
mathematically positive direction.
Then the classical Seifert form is given by
\begin{eqnarray}\label{10.5}
L^{sing}: Ml(f,\zeta)\times Ml(f,\zeta)\to\Z,\\
L^{sing}(a,b):=(-1)^{m+1}I^{Lef}(a,\gamma_{\pi} b).\nonumber
\end{eqnarray}
It has determinant $(-1)^{n(m+1)(m+2)/2}$. 
The classical monodromy $M^{sing}$ 
and the intersection form $I^{sing}$ on $Ml(f,\zeta)$ are given by
\begin{eqnarray}\label{10.6}
L^{sing}(M^{sing}a,b)&=&(-1)^{m+1}L^{sing}(b,a),\\
I^{sing}(a,b)&=& -L^{sing}(a,b)+(-1)^{m+1}L^{sing}(b,a)
\label{10.7}\\
&=&L^{sing}((M^{sing}-\id)a,b).\nonumber
\end{eqnarray}
The monodromy $M^{sing}$ is the monodromy operator of the flat
$\Z$-lattice bundle $\bigcup_{\zeta\in S^1}Ml(f,\zeta)
\cong \bigcup_{\zeta\in S^1}H_m(X_{\zeta\eta},\Z)$. 
The intersection form $I^{sing}$ is the natural intersection form
on $H_m(X_\tau,\Z)$. It is $(-1)^m$-symmetric. 

(d) The pair $(Ml(f,1),L^{sing})$ is up to a canonical
isomorphism independent of the choice of a good
representative $f:X\to\Delta_\eta$ of the function germ
$f$. So it is an invariant of the singularity $f$. 
\end{theorem}

\begin{definition}\label{t10.3}
Consider the situation in Theorem \ref{t10.2}. 
We define a unimodular bilinear lattice $(H_\Z,L)$ by 
\begin{eqnarray}\label{10.8}
H_\Z(f):=H_\Z&:=& Ml(f,1),\\
L(f):= L&:=&(-1)^{(m+1)(m+2)/2}\cdot L^{sing}.\label{10.9}
\end{eqnarray}
We call $L$ the {\it normalized} Seifert form. 
\index{normalized Seifert form}
Define $k\in\{0;1\}$ by
\begin{eqnarray}\label{10.10}
k:=m\textup{ mod}\, 2,\quad\textup{i.e. }\left\{\begin{array}{ll}
k:=0&\textup{if }m\textup{ is even,}\\
k:=1&\textup{if }m\textup{ is odd.}\end{array}\right.
\end{eqnarray}
Obviously, the monodromy $M$ and the bilinear form $I^{(k)}$
of the unimodular bilinear lattice $(H_\Z,L)$ 
are related to $M^{sing}$ and $I^{sing}$ by 
\begin{eqnarray}
M&=& (-1)^{m+1}M^{sing},\label{10.11}\\
I^{(k)}&=&(-1)^{m(m+1)/2}I^{sing}.\label{10.12}
\end{eqnarray}
We call $M$ the {\it normalized} monodromy. 
\index{normalized monodromy}
\end{definition}

\begin{remark}\label{t10.4}
Above the Seifert form $L^{sing}$ \index{Seifert form}
is defined with the help
of the isomorphism \eqref{10.3} and the intersection form
$I^{Lef}$ in \eqref{10.4}. An alternative is to derive
it from a {\it linking number} of cycles in different Milnor
fibers in the closely related fibration
$$\paa B_\varepsilon-f^{-1}(0) \to S^1, \quad z\mapsto 
\frac{f(z)}{|f(z)|}.$$
We prefer $I^{Lef}$ as the Lefschetz thimbles appear
implicitly anyway in the construction of distinguished
bases in section \ref{s10.2}. 
\end{remark}

The notion of an {\it unfolding} of an isolated hypersurface singularity 
$f$ goes back to Thom and Mather. 
An {\it unfolding} \index{unfolding} 
of $f:(\C^{m+1},0)\to(\C,0)$ is a holomorphic
function germ $F:(\C^{m+1}\times  {\MM},0)\to (\C,0)$ such that
$F|_{(\C^{m+1},0)\times\{0\}}=f$ and such that $( {\MM},0)$ 
is the germ of a complex manifold.
Its Jacobi ideal is 
$J_F:=(\frac{\paa F}{\paa z_i})\subset \OO_{\C^{m+1}\times  {\MM},0}$,
its critical space is the germ $(C,0)\subset (\C^{m+1}\times  {\MM},0)$
of the zero set of $J_F$ with the canonical complex structure.
The projection $(C,0)\to ( {\MM},0)$ is finite and flat of degree $n$.
A kind of Kodaira-Spencer map 
\index{Kodaira-Spencer map}
is the $\OO_{ {\MM},0}$-linear map
\begin{eqnarray}\label{10.13}
\aaa_C:\TT_{ {\MM},0}\to \OO_{C,0},\quad X\mapsto\www X(F)|_{(C,0)}
\end{eqnarray}
where $\www X$ is an arbitrary lift of $X\in\TT_{ {\MM},0}$ 
to $(\C^{m+1}\times  {\MM},0)$.

We will use the following notion of morphism between unfoldings. 
Let $F_i:(\C^{m+1}\times  {\MM}_i,0)\to(\C,0)$ for $i\in\{1,2\}$ be
two unfoldings of $f$ with projections 
$\pr_i:(\C^{m+1}\times  {\MM}_i,0)\to( {\MM}_i,0)$.
A {\it morphism} from $F_1$ to $F_2$ is a pair 
$(\Phi,\varphi)$ of map germs such that the following diagram
commutes,
\begin{eqnarray*}
\begin{xy}
\xymatrix{ (\C^{m+1}\times  {\MM}_1,0) \ar[r]^\Phi \ar[d]^{\pr_1} 
& (\C^{m+1}\times  {\MM}_2,0) \ar[d]^{\pr_2}\\
( {\MM}_1,0) \ar[r]^{\varphi} & ( {\MM}_2,0)  }
\end{xy}
\end{eqnarray*}
and
\begin{eqnarray*}
\Phi|_{(\C^{m+1},0)\times\{0\}}&=&\id,\\
F_1 &=& F_2\circ\Phi 
\end{eqnarray*}
hold. Then one says that $F_1$ {\it is induced} by $(\Phi,\varphi)$ from $F_2$.
An unfolding is {\it versal} if any unfolding is induced from it by a 
suitable morphism. A versal unfolding $F:(\C^{n+1}\times  {\MM},0)\to(\C,0)$ is
{\it universal} \index{universal unfolding}
if the dimension of the parameter space $( {\MM},0)$ is
minimal. Universal unfoldings exist by work of Thom and Mather.
More precisely, an unfolding is versal if and only if the map 
$\aaa_C$ is surjective, and it is universal if and only if the map $\aaa_C$
is an isomorphism (see e.g. \cite[ch. 8]{AGV85} for a proof). 
Observe that $\aaa_C$ is surjective/an isomorphism
if and only if its restriction to $0\in\MM$, the map
\begin{eqnarray}\label{10.14}
\aaa_{C,0}:T_0 {\MM}\to \OO_{\C^{m+1},0}/J_f
\end{eqnarray}
is surjective/an isomorphism. Therefore an unfolding
\begin{eqnarray}\label{10.15}
F(z_0,...,z_m,t_1,...,t_n)=F(x,t)=F_t(z)=f(z)+\sum_{j=1}^n m_jt_j,
\end{eqnarray}
with $( {\MM},0)=(\C^n,0)$ with coordinates $t=(t_1,...,t_n)$ 
where $m_1,...,m_n\in\OO_{\C^{m+1},0}$ represent a basis of the 
Jacobi algebra $\OO_{\C^{m+1},0}/J_f$, is universal.
In a versal unfolding the critical space $C$ is smooth.

Just as for $f$, one can also choose a good representative for an
unfolding $F$. In the following we suppose that such a
representative is chosen. Then $ {\MM}$ is an open subset of $\C^n$.

\begin{theorem}\label{t10.5}\cite[Theorem 5.3]{He02}
Let $f:(\C^{m+1},0)\to(\C,0)$ be an isolated hypersurface
singularity with Milnor number $n$. 
The base space $\MM$ of a good representative of a universal unfolding 
is a generically semisimple {\it F-manifold} with {\it Euler field}.
The multiplication comes from the multiplication on the right
side of the isomorphism ${\bf a}_C:\TT_{ {\MM},0}\to\OO_{C,0}$
in \eqref{10.13}, the unit field is $e={\bf a}_C^{-1}([1])$,
the Euler field is $E={\bf a}_C^{-1}([F])$. 
For each $t\in  {\MM}$, the eigenvalues of $E\circ$ are
the values of the critical points of $F_t$. 
\end{theorem}

The algebra $(T_0 {\MM},\circ)\cong (\OO_{\C^{m+1},0}/J_f)$ is irreducible.
Therefore the caustic $\KK_3\subset  {\MM}$ is not empty, so it is a
hypersurface. Also the Maxwell stratum $\KK_2\subset  {\MM}$ is not
empty, but a hypersurface. Within the caustic there is the 
{\it $\mu$-constant stratum} \index{$\mu$-constant stratum}
\begin{eqnarray}
\MM_\mu := \{t\in  {\MM}\,|\, F_t\textup{ has only one singularity }z^0,
F_t(z^0)=0\}.&&\label{10.16}
\end{eqnarray}
The only singularity which $F_t$ for $t\in \MM_\mu$ has, has
automatically Milnor number $n$ (because ${\bf a}_C$ is 
finite and flat of degree $n$). 
The $\mu$-constant stratum is a complex subvariety of ${\MM}$.
Its dimension is called the {\it modality} 
\index{modality} of the singularity $f$.

It is nontrivial, but true (see e.g. \cite[Theorem 2.2]{He11}) that
the lattices $(H_\Z,L)(F_t)$ for $t\in \MM_\mu$ glue to a local
system of $\Z$-lattices with Seifert forms (and monodromy $M$
and intersection forms $I^{(0)}$ and $I^{(1)}$) over $\MM_\mu$. 

Arnold classified all singularities with modality $\leq 2$
up to coordinate changes. We will describe them loosely in the
following theorem. We give normal forms for the simple singularities, 
the simple elliptic singularities and the hyperbolic singularities. 
Normal forms of the other singularities 
with modality $\leq 2$ can be found in \cite[II ch. 15.1]{AGV85}.
Wall classified also all singularities with modality $3$ \cite{Wa99}.

From a normal form $F_t(z)$ in $m+1$ variables $z_0,...,z_m$
and some parameters $t$ one obtains also a normal form
$F_t(z)+\sum_{j=m+1}^{m_2+1}z_j^2$ in $m_2+1$ variables.
In the following theorem, 
$(m\geq 2)$ after the symbol of a family means that the smallest
number of variables for which this family exists is $m=2$.
The meaning of $(m\geq 0)$ and $(m\geq 1)$ is analogous.
See also Remark \ref{t10.11}.

\begin{theorem}\label{t10.6} (Arnold, see \cite[II 15.1]{AGV85})

(a) \index{classification of singularities}\index{singularity}
The singularities with modality 0 are called {\sf simple singularities}.
\index{simple singularity}
There are two series and three exceptional cases (the index is the
Milnor number):
\begin{eqnarray*}
\begin{array}{c|c|c|c|c}
A_n\ (m\geq 0) & D_n\ (m\geq 1) & E_6\ (m\geq 1) 
& E_7\ (m\geq 1)& E_8\ (m\geq 1)\\ \hline 
z_0^{n+1} & z_0^{n-1}+z_0z_1^2 & z_0^4+z_1^3 
& z_0^3z_1+z_1^3 & z_0^5+z_1^3
\end{array}
\end{eqnarray*}

(b) The singularities with modality 1 are called 
{\sf unimodal singularities}. \index{unimodal singularity}
They fall into three groups of 1-parameter families: 

(i)  The three (1-parameter families of) 
{\sf simple elliptic singularities} \index{simple elliptic singularity}
$\www{E}_6$, $\www{E}_7$ and $\www{E}_8$, 
here the parameter is $\lambda\in \C-\{0,1\}$, 
\begin{eqnarray*}
\begin{array}{l|l|l}\textup{name} & n & normal form \\ \hline 
\www{E}_6 (m\geq 2) & 8 & f_\lambda=z_1(z_1-z_0)(z_1-\lambda z_0)-z_0z_2^2 \\
\www{E}_7\ (m\geq 1) & 9 & f_\lambda=z_0z_1(z_1-z_0)(z_1-\lambda z_0) \\
\www{E}_8\ (m\geq 1) & 10 & f_\lambda=z_1(z_1-z_0^2)(z_1-\lambda z_0^2)
\end{array}
\end{eqnarray*}
These normal forms are called {\sf Legendre normal forms};
they are from \cite[1.9]{Sa74}.

(ii) The (1-parameter families of)
{\sf hyperbolic singularities} \index{hyperbolic singularity}
$T_{pqr}$ with three discrete parameters
$p,q,r\in\Z_{\geq 2}$ with $p\geq q\geq r$ and 
$\frac{1}{p}+\frac{1}{q}+\frac{1}{r}<1$, 
here the continuous parameter is $t\in\C^*$,
\begin{eqnarray*}
z_0^p+z_1^q+z_2^r+tz_0z_1z_2.
\end{eqnarray*}

(iii) 14 (1-parameter families of) 
{\sf exceptional unimodal singularities}
\index{exceptional unimodal singularity}
with the following names (the index is the Milnor number),
\begin{eqnarray*}
m\geq 1:&& E_{12},\ E_{13},\ E_{14},\ Z_{11},\ Z_{12},\ Z_{13},\ 
W_{12},\ W_{13} \\
m\geq 2:&& Q_{10},\ Q_{11},\ Q_{12},\ S_{11},\ S_{12},\ U_{12}.
\end{eqnarray*}

(c) The singularities with modality 2 are called 
{\sf bimodal singularities}. \index{bimodal singularity}
They fall into three groups of 2-parameter families:

(i) The six (2-parameter families of) 
{\sf quadrangle singularities} \index{quadrangle singularity}
with the following names,
\begin{eqnarray*}
\begin{array}{l|l|l|l|l||l|l|l|l|l}
m\geq 1 & \textup{name} & E_{3,0} & Z_{1,0} & W_{1,0} &
m\geq 2 & \textup{name} & Q_{2,0} & S_{1,0} & U_{1,0} \\
                   & n  & 16      & 15      & 15 & 
                   & n  & 14      & 14      & 14 
\end{array}
\end{eqnarray*}

(ii) The eight (2-parameter families of) {\sf bimodal series}
\index{bimodal series}
with one discrete paramter $p\in\N$ and with the following names 
\begin{eqnarray*}
\begin{array}{l|l|l|l|l|l}
m\geq 1 & \textup{name} & E_{3,p} & Z_{1,p} & W_{1,p} & W_{1,p}^\sharp \\
                   & n  & 16+p    & 15+p    & 15+p    & 15+p \\ \hline
m\geq 2 & \textup{name} & Q_{2,p} & S_{1,p} & S_{1,p}^\sharp & U_{1,p} \\ 
                   & n  & 14+p    & 14+p    & 14+p           & 14+p 
\end{array}
\end{eqnarray*}

(iii) The 14 (2-parameter families of) {\sf exceptional bimodal 
singularities} \index{exceptional bimodal singularity}
with the following names (the index is the Milnor number), 
\begin{eqnarray*}
m\geq 1:&& E_{18},\ E_{19},\ E_{20},\ Z_{17},\ Z_{18},\ Z_{19},\ 
W_{17},\ W_{18} \\
m\geq 2:&& Q_{16},\ Q_{17},\ Q_{18},\ S_{16},\ S_{17},\ U_{16}.
\end{eqnarray*}
\end{theorem}

The simple singularities and the simple elliptic singularities can be
characterized in many different ways. Arnold found the following
characterizations.

\begin{theorem}\label{t10.7}
(a) (Arnold \cite{Ar73-1}) The simple singularities are the only 
isolated hypersurface singularities where $I^{(0)}$ on $H_\Z$ is
positive definite.

(Classical) Then $(H_\Z,I^{(0)})$ is a root lattice of the same type
as the name, $M$ is a Coxeter element, 
$\Gamma^{(0)}$ is the Weyl group, 
$\Delta^{(0)}$ is the set of roots, 
so especially $\Delta^{(0)}=R^{(0)}$. 

(b) (Arnold \cite{Ar73-2}) The simple elliptic singularities 
are the only isolated hypersurface singularities 
where $I^{(0)}$ on $H_\Z$ is positive semidefinite. 
\end{theorem}

\section{Distinguished bases}\label{s10.2}

Let $f:(\C^{m+1},0)\to(\C,0)$ be an isolated hypersurface
singularity. Theorem \ref{t10.2} and Definition \ref{t10.3}
gave a unimodular bilinear lattice $(H_\Z,L)$ which is
associated canonically to the singularity $f$. Here $H_\Z$
is a Milnor lattice of a fiber with value in $\R_{>0}$
and $L$ is the normalized Seifert form. 
In fact, the pair $(H_\Z,L)$ comes equipped naturally with a 
$\Br_n\ltimes\{\pm 1\}^n$ orbit $\BB^{dist}$ of 
triangular bases. They are called {\it distinguished bases}.
Also the set $\BB^{dist}$ is an invariant of the singularity 
$f$. This section gives the construction of the distinguished bases 
and lists several results on them. They have special properties.
There are many open questions around them. 

\begin{figure}
\includegraphics[width=0.7\textwidth]{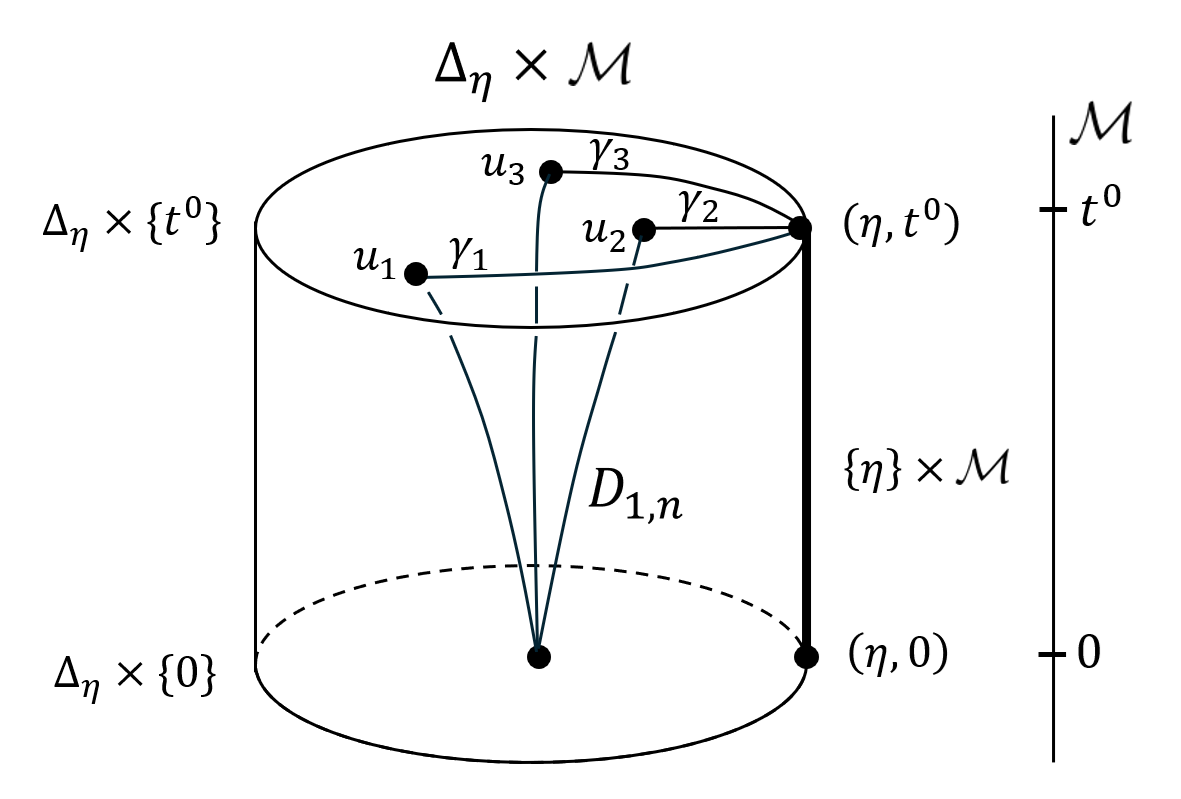}
\caption[Figure 10.2]{A schematic picture of 
$\Delta_\eta\times\MM\supset D_{1,n}$ for a representative
$F:\XX\to\Delta_\eta$ of a universal unfolding $F$ of a 
singularity $f$}
\label{Fig:10.2}
\end{figure}

Let $F:\XX\to\Delta_\eta$ be a good representative of
a universal unfolding $F$ of $f$, with 
$\XX=F^{-1}(\Delta_\eta)\cap (B_\varepsilon\times  {\MM})$
where $ {\MM}\subset \C^n$ is a sufficiently small open ball 
around 0 in $\C^n$. As above the critical space is
$C\subset \XX$. Its image under the map
$$(F,\pr_{{\MM}}):\XX\subset B_\varepsilon\times  
{\MM}\to \Delta_\eta\times  {\MM},
\quad (z,t)\mapsto (F(z,t),t)$$
is the {\it discriminant} $D_{1,n}$,
a hypersurface in $\Delta_\eta\times  {\MM}$. The projection
$D_{1,n}\to  {\MM}$ is finite and flat of degree $n$.
Of course, for that $ {\MM}$ must have been chosen small enough.
For any $t\in  {\MM}$ the intersection 
$D_{1,n}\times\{t\}\subset \Delta_\eta\times\{t\}$ is the
set of critical values of the critical points of $F_t$.
See Figure \ref{Fig:10.2}. 

Recall that the caustic $\KK_3\subset  {\MM}$ 
and the Maxwell stratum $\KK_2\subset  {\MM}$ are complex
hypersurfaces. Choose any point $t^0\in  {\MM}-(\KK_3\cup\KK_2)$,
so a generic point. Then $F_{t^0}:\XX_{t^0}\to\Delta_\eta$
with $\XX_{t^0}:=\XX\cap(B_\varepsilon\times\{t^0\})$ is a 
holomorphic function with $n$ $A_1$-singularities
and such that the critical values $u_1,...,u_n\in\Delta_\eta$
of these $A_1$-singularities are pairwise different. 

The union $\bigcup_{\tau\in \Delta_\eta-\{u_1,...,u_n\}}
H_m(F_t^{-1}(\tau),\Z)$ is a flat $\Z$-lattice bundle of
rank $n$. Near each critical value $u_j$ there is a
lattice vector, which is unique up to the sign
and which comes geometrically from the vanishing cycle
of the $A_1$ singularity in the fiber over $u_j$.
One chooses a distinguished system 
$(\uuuu{\gamma};\sigma)$ of paths $\gamma_1,...,\gamma_n$
from $\eta$ to $u_{\sigma(1)},...,u_{\sigma(n)}$ 
and pushes the $n$ vanishing
cycles along these paths to $H_m(F_t^{-1}(\eta),\Z)$. 
There is a canonical isomorphism
$$H_m(F_t^{-1}(\eta),\Z)\cong H_m(X_\eta,\Z)\cong Ml(f,1)=H_\Z.$$
One obtains a tuple $\uuuu{v}=(v_1,...,v_n)$ of cycles 
in $H_\Z$. Brieskorn \cite[Appendix]{Br70} 
first wrote a proof of part (a) of the following theorem
(see also \cite{AGV88}, \cite[Satz 5.5 in 5.4, and 5.6]{Eb01}, 
\cite[\S 1 1.5]{AGLV98}).

\begin{theorem}\label{t10.8}
(a) The resulting tuple $\uuuu{v}$ 
is a triangular basis of $(H_\Z,L)$ with 
$S:=L(\uuuu{v}^t,\uuuu{v})^t\in T^{uni}_n(\Z)$

(b) The $\Br_n\ltimes\{\pm 1\}^n$ orbit $\BB^{dist}$
of this basis $\uuuu{v}$ is an invariant of the singularity
$f$. Its elements are called {\sf distinguished bases}. 
\index{distinguished basis}
\end{theorem}

\begin{remarks}\label{t10.9}
(i) If one pushes the vanishing cycle above $u_{\sigma(j)}$
along the path $\gamma_j$ to the fiber $F_{t^0}^{-1}(\eta)$
over $\eta$, the union of cycles in the fibers above
$\gamma_j$ forms an $m+1$ dimensional cycle with boundary
in $F_{t^0}^{-1}(\eta)$, which gives a homology class
in $H_{m+1}(\XX_{t^0}, F_{t^0}^{-1}(\eta),\Z)$. 
It is a {\it Lefschetz thimble}. \index{Lefschetz thimble}
Its homology class is mapped by the isomorphism
$$H_{m+1}(\XX_{t^0}, F_{t^0}^{-1}(\eta),\Z) \cong 
H_m(F_{t^0}^{-1}(\eta),\Z)$$ 
to the homology class of its boundary, which is the
cycle in the fiber $F_{t^0}^{-1}(\eta)$. 
See Figure \ref{Fig:10.3}.

\begin{figure}
\includegraphics[width=0.7\textwidth]{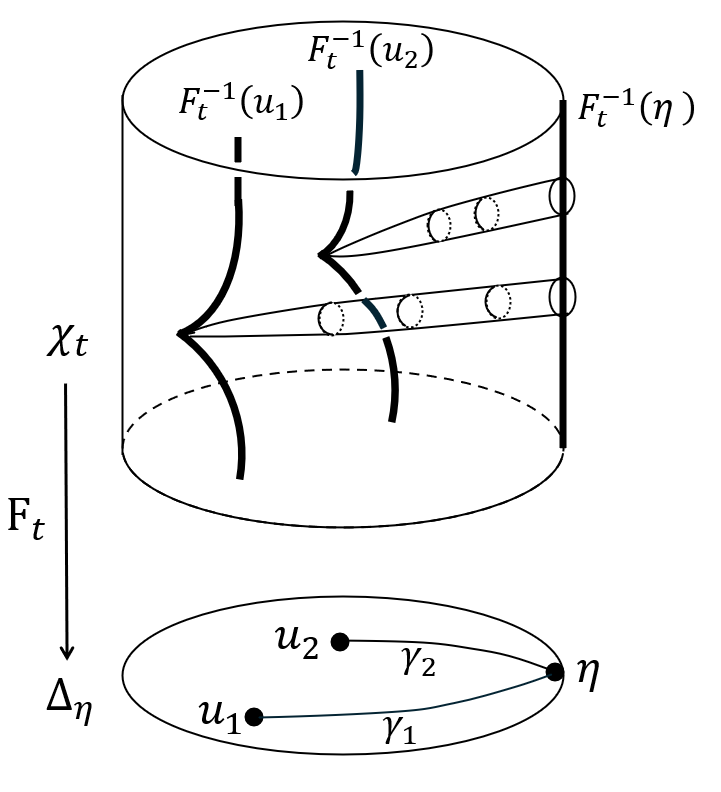}
\caption[Figure 10.3]{Two Lefschetz thimbles above paths
from $\eta$ to critical values $u_1$ and $u_2$}
\label{Fig:10.3}
\end{figure}

(ii) Recall $k=0$ if $m\equiv 0(2)$ and $k=1$ if
$m\equiv 1(2)$. The $\Z$-lattice bundle
$\bigcup_{\tau\in \Delta_\eta-\{u_1,...,u_n\}}
H_m(F_t^{-1}(\tau),\Z)$
is the $\Z$-lattice structure
$V_\Z^{(k)}((\uuuu{u},\eta),[(\uuuu{\gamma};\sigma)],\uuuu{v})$.
So, in the case of an isolated hypersurface singularity,
this bundle comes from geometry and is the source of the
triangular basis $\uuuu{v}$ in $(H_\Z,L)$. 

(iii) In \cite[5.6]{Eb01} the bilinear form $L^{sing}$ 
is called $v$. Korollar 5.3 (i) there claims that 
the matrix $V$ is upper triangular. But it is lower triangular. 

(iii) In \cite{Br83} Brieskorn considered the
larger set $\BB^*:=\bigcup_{g\in\Gamma^{(k)}}
g(\BB^{dist})$, possibly because he wanted an invariant
of the singularity $f$ and did not see that $\BB^{dist}$
is one. $\BB^{dist}$ is one because the set 
$\{\eta\}\times  {\MM}$ does not mix with the discriminant.
Therefore the $\Z$-lattice $H_m(F_t^{-1}(\eta),\Z)$
is canonically isomorphic to $H_m(X_\eta,\Z)$. 
Moving $t\in  {\MM}-(\KK_3\cup\KK_2)$ changes the tuple
$(u_1,...,u_n)$ of critical values and deforms the
distinguished system of paths. It does not do more. 
$\eta$ is fixed.
\end{remarks}

Now we start a report on specific properties of the
unimodular bilinear lattices $(H_\Z,L,\uuuu{v})$ 
with distinguished bases which
come from isolated hypersurface singularities.
We start with the Thom-Sebastiani sum.

\begin{theorem}\label{t10.10}
(a) (Definition) The {\sf Thom-Sebastiani sum}
\index{Thom-Sebastiani sum}
of a singularity $f=f(z_0,...,z_{m_1})$ in $m_1+1$ variables
and a singularity $g=g(z_{m_1+1},...,z_{m_1+m_2+1})$
in $m_2+1$ different variables is the singularity $f+g$. 

(b) Thom and Sebastiani \cite{ST71} observed that 
there is a canonical isomorphism
\begin{eqnarray}\label{10.17}
H_\Z(f)\otimes H_\Z(g)&\cong& H_\Z(f+g),
\end{eqnarray}
which respects the monodromies and the normalized
monodromies.
\begin{eqnarray}
M^{sing}(f)\otimes M^{sing}(g)&\cong& M^{sing}(f+g),\nonumber\\
M(f)\otimes M(g)&\cong& M(f+g),\label{10.18}
\end{eqnarray}
Deligne (1973?, cited in \cite{AGV88}) observed that it
respects the Seifert forms (up to a sign) and the 
normalized Seifert forms,
\begin{eqnarray}
L^{sing}(f)\otimes L^{sing}(g)
&\cong& (-1)^{(m_1+1)(m_2+1)}L^{sing}(f+g).\nonumber\\
L(f)\otimes L(g)&\cong& L(f+g).\label{10.19}
\end{eqnarray}
Any distinguished basis 
$(\delta^{f}_1,...,\delta^{f}_{n(f)})$ of $f$ 
and any distinguished basis 
$(\delta^{g}_1,...,\delta^{g}_{n(g)})$ of $g$
give rise to a basis
\begin{eqnarray}\label{10.20}
(\delta^{f}_1\otimes \delta^{g}_1,...,
\delta^{f}_1\otimes \delta^{g}_{n(g)},
\delta^{f}_2\otimes \delta^{g}_1,...,
\delta^{f}_2\otimes \delta^{g}_{n(g)},...,\\
\delta^{f}_{n(f)}\otimes \delta^{g}_1,...,
\delta^{f}_{n(g)}\otimes \delta^{g}_{n(g)})\nonumber
\end{eqnarray}
of $H_\Z(f)\otimes H_\Z(g)$, with the 
lexicographic order. Gabrielov \cite{Ga73} observed 
that its image under the isomorphism \eqref{10.17} 
is a distinguished basis of $H_\Z(f+g)$.
\end{theorem}

\begin{remark}\label{t10.11}
The Thom-Sebastiani sum $f+g$ of an isolated
hypersurface singularity $f=f(z_0,...,z_{m_1})$ and
an $A_1$-singularity $g=z^2_{m_1+1}+...+z^2_{m_1+m_2+1}$
is close to the singularity $f$. 
Here $f+g$ is called a {\it suspension} of $f$. 
The generator of the Milnor lattice 
$H_\Z(g)\cong \Z$ is unique up to a sign.
Therefore here the canonical isomorphism 
$H_\Z(f)\otimes H_\Z(g)\cong H_\Z(f+g)$ induces
an isomorphism $H_\Z(f)\cong H_\Z(f+g)$ which is
unique up to the sign. Therefore any bilinear form
and any endomorphism on $H_\Z(f)$ are also well-defined
on $H_\Z(f+g)$. 
The isomorphism \eqref{10.18} tells that the normalized
monodromies $M(f)$ and $M(f+g)$ coincide, because $M(g)=\id$. 
The isomorphism \eqref{10.19} tells that the normalized
Seifert forms $L(f)$ and $L(f+g)$ coincide, because 
$L(g)(\delta^g,\delta^g)=1$ where $H_\Z(g)=\Z\delta^g$. 
Gabrielov's result implies that $f$ and $f+g$ lead to the
same unimodular bilinear lattice $(H_\Z,L)$ with
$\Br_n\ltimes\{\pm 1\}^n$ orbit $\BB^{dist}$
of distinguished bases. 
\end{remark}

The technique which was most successfully applied for
the calculation of distinguished bases of isolated
hypersurface singularities, is due to A'Campo and
Gusein-Zade. It works for plane curve singularities, 
so functions in two variables. 
Together with the Thom-Sebastiani result above,
it can be applied to many singularities in more than
two variables. We will not describe the recipe here.
The following theorem will just describe roughly 
the result. For details see the given references
as well as \cite{AGV88}.
The recipe builds on resolution of plane
curve singularities. In concrete cases it leads to
many pictures of curves in the real plane and
an interesting way how to construct them and 
deal with them.

\begin{theorem}\label{t10.12} 
(A'Campo \cite{AC75-1}\cite{AC75-2} and Gusein-Zade 
\cite{Gu74-1}\cite{Gu74-2}, 
see also \cite[\S 2 2.1]{AGLV98}
\cite{BK96} \cite{LSh18})
Let $f$ be a plane curve singularity, 
\index{plane curve singularity}
so an isolated hypersurface singularity with $m=1$.
Suppose that $f=f_1...f_r$ is the decomposition
of $f$ into its irreducible branches and that
$f_j((\R^2,0))\subset (\R,0)$ for any $j$
(so $f$ is a complexification of a 
{\sf totally real singularity}).
Any {\sf totally real morsification} 
(see e.g. \cite{AGLV98}
for its definition) gives rise to three sets
$B^{\ominus},B^{0},B^{\oplus}\subset H_\Z$ 
with the following properties: 
If $\uuuu{v}=(v_1,...,v_n)$ is a list of the
elements of $B^{\ominus}\cup B^{0}\cup B^{\oplus}$
which lists first the elements of $B^{\ominus}$,
then the elements of $B^0$ and finally the elements
of $B^{\oplus}$, then $\uuuu{v}$ is a distinguished
basis, and its matrix $S$
takes the following shape,
\begin{eqnarray}\label{10.21}
S=L(\uuuu{v}^t,\uuuu{v})^t=E_n+ \left(\begin{array}{c|c|c}
0 & - & + \\ \hline
0 & 0 & -  \\ \hline
0 & 0 & 0\end{array}\right)\in T^{uni}_n(\Z),
\end{eqnarray}
more precisely, for $i<j$ 
\begin{eqnarray}\nonumber
-I^{sing}(v_i,v_j)=I^{(1)}(v_i,v_j)=L(v_j,v_i)=S_{ij}\\
=\left\{\begin{array}{ll}
0&\textup{if }v_i,v_j\in B^\ominus\textup{ or }
v_i,v_j\in B^0\textup{ or }v_i,v_j\in B^\oplus,\\
\leq 0&\textup{if }v_i\in B^\ominus, v_j\in B^0,\\
\geq 0&\textup{if }v_i\in B^\ominus, v_j\in B^\oplus,\\
\leq 0&\textup{if }v_i\in B^0, v_j\in B^\oplus.
\end{array}\right. \label{10.22}
\end{eqnarray}
\end{theorem}

Theorem \ref{t10.10} and Theorem \ref{t10.12} 
together allow to deal with singularities 
$f(z_0,...,z_m)$ with $m\geq 2$ if $f$ is an iterated 
Thom-Sebastiani sum of plane curve singularities.
The following Theorem \ref{t10.13} of Gabrielov is stronger.
It allows to construct a distinguished basis
of any singularity (in principle) from results
for singularities with less variables. 
Again we do not give the precise recipe, but only 
the rough result of it. For details see the given
references.

\begin{theorem}\label{t10.13}
(a) \cite{Ga79} (see also \cite[5.10]{Eb01} 
\cite[\S 2 2.4]{AGLV98})
Given an isolated hypersurface singularity
$f:(\C^{m+1},0)\to(\C,0)$, there is a recipe
of Gabrielov which allows to construct a 
certain distinguished basis of $f$ from the
following data:
a generic linear function $g:\C^{m+1}\to\C$ and 
a distinguished basis of $\www f:= f|_{g^{-1}(0)}$.
For the recipe one needs sufficient information
on the {\sf polar curve} \index{polar curve}
of $f$ with respect to $g$.
Let $S$ be the matrix of the distinguished 
basis of $f$, and 
let $\www S$ be the matrix of the
distinguished basis of $\www f$. 
Then especially 
\begin{eqnarray}\label{10.23}
\{S_{ij}\,|\, i,j\in\{1,..,n(f)\}\}\subset 
\{0,\pm 1\}\cup\{\pm\www S_{ij}\,|\, 
i,j\in \{1,...,n(\www f)\}\}.
\end{eqnarray}

(b) (Corollary of part (a)) Each isolated
hypersurface singularity has a distinguished basis
with matrix $S$ with entries only $0,\pm 1$, so 
$$\{S_{ij}\,|\, i,j\in\{1,..,n(f)\}\}\subset \{0,\pm 1\}.$$
\end{theorem}

The following theorem lists further properties of 
the unimodular bilinear lattice $(H_\Z,L,\uuuu{e})$
with a distinguished basis $\uuuu{e}$, which come
from a singularity $f$.

\begin{theorem}\label{t10.14}
Let $f:(\C^{m+1},0)\to(\C,0)$ be an isolated hypersurface
singularity with Milnor number $n\in\N$, and let 
$(H_\Z,L,\uuuu{e})$ be the induced unimodular bilinear
lattice with one chosen distinguished basis $\uuuu{e}$
and matrix $S=L(\uuuu{e}^t,\uuuu{e})^t\in T^{uni}_n(\Z)$.
Fix $k\in\{0;1\}$. 

(a) (Classical, e.g. \cite{Ga79}) 
$\Delta^{(k)}=\Gamma^{(k)}(e_j)$ for any element
$e_j$ of the distinguished basis, so $\Delta^{(k)}$
is a single orbit. 

(b) \cite{Ga74-1}\cite{La73}\cite{Le73}
The triple $(H_\Z,L,\uuuu{e})$ is irreducible.

(c) If $n\geq 2$ then there are $\delta_1,\delta_2\in 
\Delta^{(k)}$ with $I^{(k)}(\delta_1,\delta_2)=1$. 

(d) (Classical) The monodromy $M$ is quasiunipotent.
The sizes of the Jordan blocks are at most $m+1$. 
The sizes of the Jordan blocks with eigenvalue 1 are at most $m$.

(e) \cite{AC73} $\tr(M)=1$. 
\end{theorem}

\begin{remarks}\label{t10.15}
(i) The discriminant $D_{1,n}\subset\Delta_\eta\times M$ of
a good representative $F$ of a universal unfolding of $f$
is irreducible because it is the image of the smooth
critical set $C\subset\XX$. Its irreducibility leads
easily to part (a), $\Delta^{(k)}=\Gamma^{(k)}(e_j)$.
See e.g. \cite[Ch. 5 Satz 5.20]{Eb01}. 

(ii) Part (a) implies immediately part (b), the 
irreducibility of the triple $(H_\Z,L,\uuuu{e})$. 
This is the proof of part (b) in \cite{Ga74-1}.
The proof of part (b) in \cite{Le73} uses part (e). 

(iii) Part (c) follows from the fact that any
singularity with Milnor number $n\geq 2$ deforms to
the $A_2$ singularity. See e.g. \cite{Eb01}. 

(iv) The triple $(H_\Z,I^{(k)},\Delta^{(k)})$ with the
properties in the parts (a) and (c) and the property
that $\Delta^{(k)}$ generates $H_\Z$ and the further
property $\Delta^{(0)}\subset R^{(0)}$ if $k=0$
is a {\it vanishing lattice}. \index{vanishing lattice} 
This notion is defined 
in \cite{Ja83} and \cite{Eb84}. 
Here $\Gamma^{(k)}=\langle s_\delta^{(k)}\,|\, 
\delta\in\Delta^{(k)}\rangle$ is determined by the
triple. 

(v) In rank 2 and rank 3 the only unimodular bilinear
lattices whose monodromy has trace 1 are those of
$A_2$ and $A_3$. So, the condition $\tr(M)=1$ is
very restrictive. 

(vi) It is an interesting open question how to characterize
those matrices $S\in T^{uni}_n(\Z)$ which come
from distinguished bases of some isolated hypersurface
singularities. Theorem \ref{t10.13} (b) and Theorem
\ref{t10.14} give a number of necessary conditions.
But very probably these conditions are not sufficient. 
\end{remarks}

Finally we come to concrete results for the first singularities
in Arnold's classification. Recall that the matrices 
$S\in \SSS^{dist}\subset T^{uni}_n(\Z)$ of distinguished bases
are called {\it distinguished matrices}. 

Gabrielov gave distinguished matrices for all unimodal singularities
in \cite{Ga73} and \cite{Ga74-2}. Ebeling gave distinguished
matrices for all bimodal singularities in \cite{Eb83}.
Distinguished matrices for the simple singularities are classical. 

\begin{figure}
\includegraphics[width=0.8\textwidth]{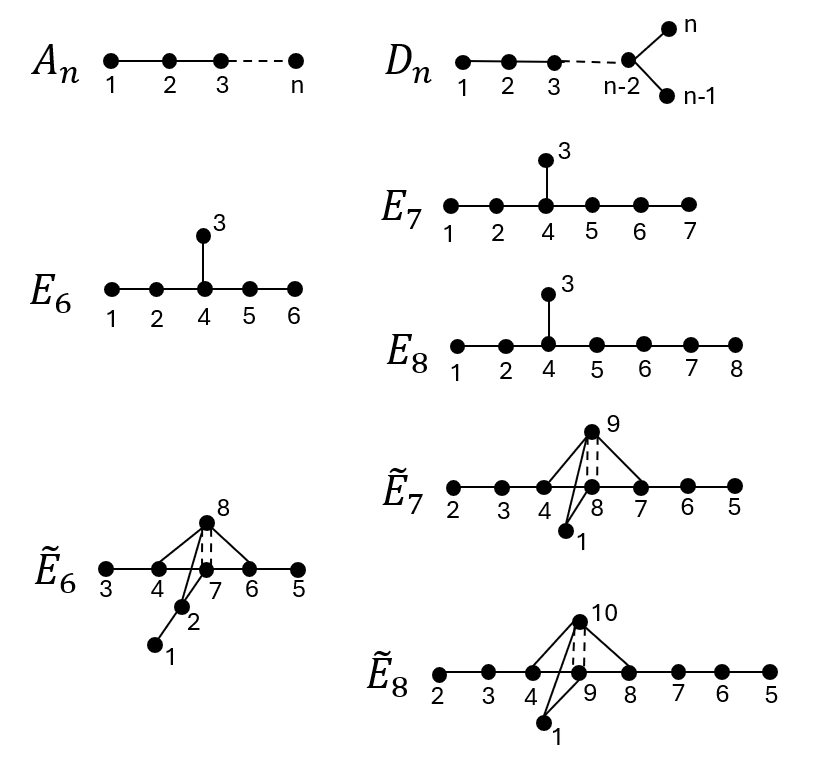}
\caption[Figure 10.4]{Coxeter-Dynkin diagrams of some
distinguished matrices for the simple and the simple elliptic
singularities}
\label{Fig:10.4}
\end{figure}

Figure \ref{Fig:10.4} encodes the information of matrices $S$
of distinguished bases for the simple singularities and
the simple elliptic singularities in {\it Coxeter-Dynkin 
diagrams}. This is more efficient than writing down 
series of matrices. See \cite[Ch. 5 Tabelle 5.3]{Eb01}
for the given Coxeter-Dynkin diagrams of the simple elliptic
singularities. In the given Coxeter-Dynkin diagrams for the 
simple singularities one can in fact change the numbering
arbitrarily.

\begin{definition}\label{t10.16}
Let $S=(S_{ij})\in T^{uni}_n(\Z)$. 
Its {\it Coxeter-Dynkin diagram} \index{Coxeter-Dynkin diagram}
$\textup{CDD}(S)$ is a graph with $n$ vertices which are 
numbered from $1$ to $n$ and with weighted edges which are
defined as follows: Between the vertices $i$ and $j$ with
$i<j$ one has no edge if $S_{ij}=0$ and an edge with
weight $S_{ij}$ if $S_{ij}\neq 0$. 
Alternatively, one draws $|S_{ij}|$ edges if $S_{ij}<0$ and
$S_{ij}$ dotted edges if $S_{ij}>0$. 
\end{definition}

The matrix $S$ and its Coxeter-Dynkin diagram 
$\textup{CDD}(S)$ determine one another.
A unimodular bilinear lattice $(H_\Z,L)$ with a 
triangular basis $\uuuu{e}$ with $S=L(\uuuu{e}^t,\uuuu{e})^t$
is reducible if and only if $CDD(S)$ consists of 
several components. 

Recall from Theorem \ref{t10.7} that the simple singularities
are the only singularities where $I^{(0)}$ on $H_\Z$ 
is positive definite and that the simple elliptic singularities
are the only singularities where $I^{(0)}$ on $H_\Z$ is 
positive semidefinite. The following Lemma applies.

\begin{lemma}\label{t10.17}
Let $S\in T^{uni}_n(\Z)$ and let $(H_\Z,L,\uuuu{e})$ 
be a unimodular bilinear lattice with triangular basis 
$\uuuu{e}$ with $S=L(\uuuu{e}^t,\uuuu{e})^t$.

(a) Suppose that $I^{(0)}$ is positive definite.
Then $S_{ij}\in\{0,\pm 1\}$. The sets $R^{(0)}$, 
$\BB^{dist}=\Br_n\ltimes\{\pm 1\}^n(\uuuu{e})$
and $\SSS^{dist}=\Br_n\ltimes\{\pm 1\}^n(S)$ are finite.

(b) Suppose that $I^{(0)}$ is positive definite.
Then $S_{ij}\in\{0,\pm 1,\pm 2\}$. The set 
$\SSS^{dist}=\Br_n\ltimes\{\pm 1\}^n(S)$ is finite.
\end{lemma}

{\bf Proof:} (a) $R^{(0)}$ is a discrete subset of the
compact set $\{\delta\in H_\R\,|\, I^{(0)}(\delta,\delta)=2\}$
and therefore finite. Therefore also the sets 
$\Delta^{(0)}$, $\BB^{dist}$ and $\SSS^{dist}$ are finite.
With $S+S^t$ also the submatrix
$\begin{pmatrix}2 & S_{ij} \\ S_{ij} & 2\end{pmatrix}$
(with $i<j$) is positive definite, so $S_{ij}\in\{0,\pm 1\}$.

(b) With $S+S^t$ also the submatrix
$\begin{pmatrix}2 & S_{ij} \\ S_{ij} & 2\end{pmatrix}$
(with $i<j$) is positive semidefinite, so $S_{ij}\in\{0,\pm 1,\pm 2\}$.
Therefore the set $\SSS^{dist}$ is finite.\hfill$\Box$ 

\bigskip
In section \ref{s3.3} we considered the unimodular bilinear
lattice $(H_\Z,L,\uuuu{e})$ with triangular basis $\uuuu{e}$
with arbitrary matrix $S=L(\uuuu{e}^t,\uuuu{e})^t\in T^{uni}_n(\Z)$
and posed the questions when the following two inclusions
are equalities,
\begin{eqnarray*}
\BB^{dist}&\subset& \{\uuuu{v}\in (\Delta^{(0)})^n\,|\, 
s_{v_1}^{(0)}...s_{v_n}^{(0)}=-M\},\hspace*{2cm}(3.3)\\
\BB^{dist}&\subset& \{\uuuu{v}\in (\Delta^{(1)})^n\,|\, 
s_{v_1}^{(1)}...s_{v_n}^{(1)}=M\}.\hspace*{2.4cm}(3.4)
\end{eqnarray*}
Chapter \ref{s7} was devoted to answer these questions in the
cases of rank 2 and 3. General positive answers were given in 
Theorem \ref{t3.7} and the Examples \ref{t3.22} (iv) and (v).
They embrace the cases of the simple singularities. 

References for the results in the next theorem are
given after the theorem.

\begin{theorem}\label{t10.18}
(a)

\begin{tabular}{l|l|l|l}
$f$ & $S_{ij}$ & $|\BB^{dist}|$ & 
$|\SSS^{dist}|$ \\ \hline 
ADE-singularity & in $\{0,\pm 1\}$ & finite & finite \\
simple elliptic sing. &  in $\{0,\pm 1,\pm 2\}$ & 
infinite & finite \\
any other singularity & unbounded & infinite & infinite 
\end{tabular}

\medskip
More details are given below in \eqref{10.33} and \eqref{10.47}.

(b) The cases of the simple singularities:
$\Delta^{(0)}=R^{(0)}$. 
The inclusion in \eqref{3.3} is an equality.

(c) The cases of the simple elliptic singularities:
$\Delta^{(0)}=R^{(0)}$. 
The inclusion in \eqref{3.3} is an equality.

(d) The cases of the hyperbolic singularities:
The inclusion in \eqref{3.3} is an equality.
\end{theorem}

\begin{remarks}\label{t10.19}
(a) The first two lines in the table in part (a) 
follow with the exception of
$|\BB^{dist}|=\infty$ 
for the simple elliptic singularities 
from Lemma \ref{t10.17}
and Theorem \ref{t10.7}. The fact
$|\BB^{dist}|=\infty$ for the simple elliptic singularities
is proved for example in \cite{Kl87} or in \cite{HR21}.
The third line in the table is due to \cite{Eb18}.

(b) $\Delta^{(0)}=R^{(0)}$ is classical.
The equality in \eqref{3.3} was first proved by Deligne
in a letter to Looijenga \cite{De74}.
It was rediscovered and generalized by Igusa and
Schiffler \cite{IS10} to all Coxeter groups, 
see the Theorems \ref{t3.6} and \ref{t3.7}.

(c) \cite[Korollar 4.4]{Kl83} says $\Delta^{(0)}=R^{(0)}$, 
\cite[VI 1.2 Theorem]{Kl87} gives equality in \eqref{3.3}.

(d) \cite[Theorem 1.5]{BWY19}.

(e) Theorem \ref{t7.1} (a) implies that in the case $A_2$
the inclusion in \eqref{3.4} is an equality.
Theorem \ref{t7.7} (g) implies that in the case $A_3$
the inclusion in \eqref{3.4} is not an equality, but that
it becomes an equality if one adds on the right hand side
the condition $\sum_{i=1}^n\Z v_i=H_\Z$. 
These two results for \eqref{3.4} and the results in 
Theorem \ref{t10.18} for \eqref{3.3} are the only results 
on the interesting question for which singularities 
the inclusions in \eqref{3.3} or \eqref{3.4} are equalities. 
\end{remarks}

\section[The groups $G_\Z$ for some singularities]
{Orlik blocks and the groups $G_\Z$ for the simple, 
unimodal and bimodal singularities}\label{s10.3}

The groups $G_\Z$ for the unimodular bilinear lattices
$(H_\Z,L)$ from isolated hypersurface singularities 
have been mainly studied by the first author and F. Gau{\ss} 
\cite{He92}\cite{He11}\cite{GH17}\cite{GH18}.
This section reports about results, 
especially for the singularities with modality $\leq 2$.
But it starts with some rather recent general results on Orlik blocks
\cite{He20} \cite{HM22-1} \cite{HM22-2} (which would have simplified
part of the work in the references above, if they had been
available at that time).

\begin{definition}\label{t10.20}
Consider a $\Z$-lattice $H_\Z$ of rank $n\in\N$ with an
automorphism $M$ which is of finite order and cyclic.
Here {\it cyclic} means that an element $e_1\in H_\Z$ with
$H_\Z=\bigoplus_{i=0}^{n-1}\Z M^{i-1}(e_1)$ exists.
The pair $(H_\Z,M)$ is called an {\it Orlik block}. \index{Orlik block}
\end{definition}

\begin{remarks}\label{t10.21} 
Consider an Orlik block $(H_\Z,M)$ with $H_\Z$ a
$\Z$-lattice of rank $n$. 

(i) The eigenvalues of $M$ are unit roots because $M$ is
of finite order. Cyclic implies regular. 
Finite order and regular together imply that each eigenvalue
has multiplicity one. Especially, $M$ is semisimple.
Denote by $\Ord(M)\subset\N$ the finite set of orders 
of the eigenvalues of $M$. 
Then $p_{ch,M}=\prod_{m\in\Ord(M)}\Phi_m$ is a
product of cyclotomic polynomials. The pair $(H_\Z,M)$ is
up to isomorphism determined by the set $\Ord(M)$
or by the polynomial $p_{ch,M}$. 

(ii) Lemma \ref{t5.2} (a) implies $\End(H_\Z,M)=\Z[M]$.
Consider the group
\begin{eqnarray*}
\Aut_{S^1}(H_\Z,M)&:=&\{g\in\Aut(H_\Z,M)\,|\, 
\textup{all eigenvalues of }g\textup{ are in }S^1\}\\
&=&\{p(M)\,|\, p(t)=\sum_{i=0}^{n-1}p_it^i\in\Z[t],\\
&& p(\lambda)p(\lambda^{-1})=1\textup{ for all eigenvalues }
\lambda\textup{ of }M\}.
\end{eqnarray*}
Here we used $\lambda^{-1}=\oooo{\lambda}$ for $\lambda$
eigenvalues of $M$.  

(iii) If $M$ is the monodromy of a unimodular bilinear lattice
$(H_\Z,L)$ then $\Aut_{S^1}(H_\Z,M)=G_\Z=G_\Z^{(0)}=G_\Z^{(1)}$
by Lemma \ref{t5.2} (b) (ii).

(iv) $\Aut_{S^1}(H_\Z,M)$ contains the finite group 
$\{\pm M^l\,|\, l\in\Z\}$. The main result Theorem 1.2 
in \cite{He20} gives necessary and sufficient conditions 
on the set $\Ord(M)$ or (equivalently) on 
the polynomial $p_{ch,M}$ for
\begin{eqnarray*}
\Aut_{S^1}(H_\Z,M)=\{\pm M^l\,|\, l\in\Z\}.
\end{eqnarray*}
\end{remarks}

Orlik \cite{Or72} conjectured that in the case of a 
{\it quasihomogeneous singularity} the pair $(H_\Z,M)$
decomposes in a specific way into a sum of Orlik blocks. 
Part of this conjecture was proved in \cite{HM22-1}\cite{HM22-2}.
Before stating the conjecture and the results precisely,
we define quasihomogeneous singularities and two special
families of quasihomogeneous singularities.

\begin{definition}\label{t10.22}
(a) An isolated hypersurface singularity 
$f:(\C^{m+1},0)\to(\C,0)$ is a {\it quasihomogoneous singularity},
\index{quasihomogeneous singularity}
if $f\in\C[z_0,...,z_m]$ is a {\it quasihomogeneous polynomial},
that means there exists a weight system $\uuuu{w}=(w_0,...,w_m)
\in (0,\frac{1}{2}]\cap \Q$ such that for any monomial 
$\prod_{i=0}^m z_i^{e_i}$ with nonvanishing coefficient in $f$ its
weighted degree $\deg_{\uuuu{w}}(\prod_{i=0}^m z_i^{e_i}):=
\sum_{i=0}^mw_ie_i$ is equal to $1$. 

(b) A {\it chain type singularity} \index{chain type singularity}
is a quasihomogeneous singularity of the
special shape
\begin{eqnarray}\label{10.24}
f=f(z_0,...,z_m)= z_0^{a_0+1}+\sum_{i=1}^mz_{i-1}z_i^{a_i}
\end{eqnarray}
for some $m\in\N$ and some $a_0,...,a_m\in\N$.
This quasihomogeneous polynomial has indeed an isolated singularity.

(c) A {\it cycle type singularity} \index{cycle type singularity}
is a quasihomogeneous singularity
of the special shape
\begin{eqnarray}\label{10.25}
f=f(z_0,...,z_m)=\sum_{i=0}^{m-1}z_i^{a_i}z_{i+1}+z_m^{a_m}z_0
\end{eqnarray}
for some $m\in\N$ and some $a_0,...,a_m\in\N$ which satisfy 
in the case of odd $m$ that neither $a_j=1$ for all even $j$ 
nor $a_j=1$ for all odd $j$. 
This quasihomogeneous polynomial has indeed an isolated singularity.
\end{definition}

The quasihomogeneous singularities form an important and especially
well studied subfamily of all isolated hypersurface singularities. 
Their monodromies are semisimple.

\begin{remarks}\label{t10.23}
Let $f\in\C[z_0,...,z_m]$ be a quasihomogeneous singularity.

(i) Milnor and Orlik \cite{MO70} gave a formula for the characteristic
polynomial $p_{ch,M^{sing}}$ of its monodromy. Recall
$M=(-1)^{m+1} M^{sing}$. 

(ii) The polynomial $p_{ch,M}$ has a unique decomposition
$p_{ch,M}=p_1p_2...p_l$ with $p_1\,|\, p_2\,|\, ...\,|\, p_l$
and $p_l$ the minimal polynomial of $M$ and $l\in\N$ suitable.

(iii) Orlik \cite[Conjecture 3.1]{Or72} conjectured that the pair $(H_\Z,M)$
has a decomposition into $l$ Orlik blocks with characteristic
polynomials $p_1,...,p_l$. 

(iv) His conjecture was proved in \cite[Theorem 1.3]{HM22-1} 
for all cycle type singularities and in \cite[Theorem 1.3]{HM22-2} 
for all iterated Thom-Sebastiani sums of chain type singularities 
and cycle type singularities.

(v) Furthermore, in \cite[Theorem 1.5]{HM22-3} it was proved that each
polynomial $p_j$ in  (ii) satisfies the sufficient
condition (I) in Theorem 1.2 in \cite{He20} for
\begin{eqnarray*}
&&\Aut_{S^1}(\textup{Orlik block for }p_j)\\
&=&\{\pm (\textup{monodromy of the Orlik block})^l\,|\, l\in\Z\}.
\end{eqnarray*}

(vi) Therefore if the singularity $f$ satisfies Orlik's conjecture,
one knows the automorphisms of the summands of a decomposition
of $(H_\Z,M)$ into Orlik blocks. This is a big step towards $G_\Z$.
Though in general, the decomposition is not left or right $L$-orthogonal.
Also, automorphisms might mix or exchange Orlik blocks.
So, there are still more steps to do towards $G_\Z$.

(vii) Many families of singularities with modality $\leq 2$ 
contain quasihomogeneous singularities, namely all families 
except the hyperbolic singularities $T_{pqr}$ 
and the eight bimodal series
$E_{3,p},Z_{1,p},W_{1,p},W_{1,p}^\sharp,Q_{2,p},S_{1,p},
S_{1,p}^\sharp,U_{1,p}$. 
Each family which contains at all a quasihomogeneous singularity, 
contains also a quasihomogeneous singularity 
which is an iterated Thom-Sebastiani sum of 
chain type singularities and cycle type singularities. 
Therefore Orlik's conjecture is true for all of these families
of singularities. In \cite{He92} this was proved in a case-by-case
work using Coxeter-Dynkin diagrams.

(viii) But Orlik blocks are not only useful for quasihomogeneous
singularities. Also in the families with modality $\leq 2$ which
do not contain quasihomogeneous singularities, $(H_\Z,M)$
usually contains an $M$-invariant sublattice of finite index in $H_\Z$,
which decomposes into Orlik blocks which each satisfy the sufficient
condition (I) in Theorem 1.2 in \cite{He20}. 
This fact was elaborated and used in many cases in \cite{GH17}
and \cite{GH18}. See for example Theorem \ref{t10.26} (ii) below. 

(ix) Orlik and Randell \cite{OR77} found for the chain type singularities an
$n$-th root of the monodromy, without using Theorem \ref{t3.26} (c).
The pair $(H_\Z,\textup{this }n\textup{-th root})$ is an Orlik 
block. Though they conjectured that a distinguished basis 
with matrix $S$ as in Theorem \ref{t3.26} (c) exists.
This conjecture was proved recently by Varolgunes \cite{Va23}.
\end{remarks}

Now we describe a part of the results in \cite{He11}\cite{GH17}
and \cite{GH18} on $G_\Z$ in the families of singularities with
modality $\leq 2$. We start with the table in the proof of
Theorem 8.3 in \cite{He11} of the characteristic polynomials
$p_{ch,-M}$ for all families of singularities with modality $\leq 2$
which contain quasihomogeneous singularities.
They can be extracted from the tables of spectral numbers in
\cite[13.3.4]{AGV88} or from \cite{He92}. 
The families consist of all singularities with modality $\leq 2$
except the hyperbolic singularities and the eight bimodal series.

\medskip
\begin{eqnarray*}
\begin{array}{ll|ll|ll}
 A_n   &  \frac{t^{n+1}-1}{t-1}  &  \www E_6 & \Phi_3^3\Phi_1^2  
& E_{3,0}  &  \Phi_{18}^2\Phi_6\Phi_2^2 \\
 D_n   &  (t^{n-1}+1)\Phi_2      &  \www E_7 & \Phi_4^2\Phi_2^3\Phi_1^2  
& Z_{1,0}  &  \Phi_{14}^2\Phi_2^3 \\
 E_6     &  \Phi_{12}\Phi_3          & \www E_8 & \Phi_6\Phi_3^2\Phi_2^2\Phi_1^2 
& Q_{2,0}  &  \Phi_{12}^2\Phi_4^2\Phi_3 \\
 E_7     &  \Phi_{18}\Phi_2          & & & W_{1,0}  &  
\Phi_{12}^2\Phi_6\Phi_4\Phi_3\Phi_2 \\
 E_8     &  \Phi_{30}                & & & S_{1,0}  &  \Phi_{10}^2\Phi_5\Phi_2^2 \\
 &   & & & U_{1,0}  &  \Phi_9^2\Phi_3 
\end{array}
\end{eqnarray*}

\begin{eqnarray*}
\begin{array}{ll|ll|ll}
E_{12}  &  \Phi_{42}           
& E_{13}  &  \Phi_{30}\Phi_{10}\Phi_2    
& E_{14}  &  \Phi_{24}\Phi_{12}\Phi_3 \\
Z_{11}  &  \Phi_{30}\Phi_6\Phi_2  
& Z_{12}  &  \Phi_{22}\Phi_2^2  
& Z_{13}  &  \Phi_{18}\Phi_9\Phi_2 \\
Q_{10}  &  \Phi_{24}\Phi_3  
& Q_{11}  &  \Phi_{18}\Phi_{6}\Phi_{3}\Phi_{2}  
& Q_{12}  &  \Phi_{15}\Phi_{3}^2 \\
W_{12}  &  \Phi_{20}\Phi_{5}  
& W_{13}  &  \Phi_{16}\Phi_{8}\Phi_{2} & \\ 
S_{11}  &  \Phi_{16}\Phi_{4}\Phi_{2} 
& S_{12}  &  \Phi_{13} 
& U_{12}  &  \Phi_{12}\Phi_{6}\Phi_{4}^2\Phi_{2}^2  
\end{array}
\end{eqnarray*}

\begin{eqnarray*}
\begin{array}{ll|ll|ll}
E_{18}   &  \Phi_{30}\Phi_{15}\Phi_3 
& E_{19}  &  \Phi_{42}\Phi_{14}\Phi_2
& E_{20}  &  \Phi_{66} \\
Z_{17}  &  \Phi_{24}\Phi_{12}\Phi_6\Phi_3\Phi_2  
& Z_{18}  &  \Phi_{34}\Phi_2^2 
& Z_{19}  &  \Phi_{54}\Phi_2 \\
Q_{16}  &  \Phi_{21}\Phi_3^2 
& Q_{17}  &    \Phi_{30}\Phi_{10}\Phi_{6}\Phi_{3}\Phi_2 
& Q_{18}  &  \Phi_{48}\Phi_{3} \\
W_{17}  &  \Phi_{20}\Phi_{10}\Phi_{5}\Phi_{2}
& W_{18}  &  \Phi_{28}\Phi_{7}
& \\
S_{16}  &  \Phi_{17}
& S_{17}  &  \Phi_{24}\Phi_{8}\Phi_{6}\Phi_{3}\Phi_2 
& U_{16}  &  \Phi_{15}\Phi_{5}^2 
\end{array}
\end{eqnarray*}

One sees that $p_{ch,-M}$ has no multiple roots
in the following cases: The simple singularities 
$A_n$, $D_{2n+1}$, $E_6,E_7,E_8$, so all simple singularities 
except the $D_{2n}$,
and 22 of the 28 exceptional unimodal
and bimodal singularities, namely all except 
$Z_{12},Q_{12},U_{12}$, $Z_{18},Q_{16}$ and $U_{16}$. 
Part (a) of the next theorem follows from the Remarks
\ref{t10.21} (iv) and \ref{t10.23} (iv), (v) and (vii). 

\begin{theorem}\label{t10.24}
(a) \cite[Theorems 8.3 and 8.4]{He11} In the cases
$A_n$, $D_{2n+1}$, $E_6,E_7,E_8$ and the 22 of the 28 exceptional unimodal
and bimodal singularities with the exceptions 
$Z_{12},Q_{12},U_{12},Z_{18},Q_{16}$ and $U_{16}$,
the pair $(H_\Z,M)$ is a single Orlik block and 
$$G_\Z=\{\pm M^l\,|\, l\in\Z\}.$$

(b) \cite[Theorem 4.1]{GH17} In the cases 
$D_{2n}$, $Z_{12},Q_{12},U_{12}$, $Z_{18},Q_{16}$ and $U_{16}$,
$(H_\Z,M)$ has a decomposition into two Orlik blocks and 
\begin{eqnarray*}
G_\Z=\{\pm M^l\, |\, l\in\Z\}\times U
\end{eqnarray*}
with
\begin{eqnarray*}
\begin{array}{l|l|l|l|l|l|l|l|l}
 & D_4 & D_{2n}\textup{ with }n\geq 3& 
 Z_{12} & Q_{12} & U_{12} & Z_{18} & Q_{16} & U_{16}\\ \hline
U\cong & S_3 & S_2 & \{\id\} & S_2 & S_3 & \{\id\} & S_2 & S_3
\end{array}
\end{eqnarray*}

(c) \cite[Theorem 3.1]{GH17} In the simple elliptic cases
$\Rad I^{(0)}\subset H_\Z$ has rank 2. The restriction
map $ G_\Z\to G_\Z|_{\Rad I^{(0)}}$ is called $\pr^{\Rad I^{(0)}}$. 
There is an exact non-splitting sequence
\begin{eqnarray*}
\{\id\}\to \ker\pr^{\Rad I^{(0)}} \to G_\Z
\stackrel{\pr^{\Rad I^{(0)}}}{\to} G_\Z|_{\Rad I^{(0)}}
\cong SL_2(\Z)\to \{\id\}
\end{eqnarray*}
with finite group $\ker \pr^{\Rad I^{(0)}}$. The sublattice
$(H_{\C,\neq 1}\cap H_\Z,M|_{..})$ has a left and right $L$-orthogonal 
decomposition into three Orlik blocks
$H^{(1)}\oplus H^{(2)}\oplus H^{(3)}$ with restricted monodromies
$M^{(1)}$, $M^{(2)}$ and $M^{(3)}$ and characteristic polynomials
in the following table,
\begin{eqnarray*}
p_{ch,-M^{(j)}}=\frac{t^{p_j}-1}{t-1}\quad\textup{with}\quad
\begin{array}{l|l|l|l}
 & \www{E}_6 & \www{E}_7 & \www{E}_8 \\ \hline 
(p_1,p_2,p_3) & (3,3,3) & (4,4,2) & (6,3,2) 
\end{array}
\end{eqnarray*}
Then $\ker \pr^{\Rad I^{(0)}}=U_1\rtimes U_2$, where\\
$U_2\cong S_3$ permutes the three Orlik blocks in the case of $\www{E}_6$,\\
$U_2\cong S_2$ permutes $H^{(1)}$ and $H^{(2)}$ in the case of $\www{E}_7$,\\
and $U_2=\{\id\}$ in the case of $\www{E}_8$.\\
The elements of $U_1$ are the extensions to $\Rad I^{(0)}$ by $\id$
of the following automorphisms of $H_{\C,\neq 1}\cap H_\Z$, 
\begin{eqnarray*}
\{(M^{(1)})^\alpha\times (M^{(2)})^\beta\times (M^{(3)})^{\gamma}\,|\, 
(\alpha,\beta,\gamma)\in\prod_{j=1}^3\{0,1,..,p_j-1\}\\
\textup{with}\quad 
\frac{\alpha}{p_1}+\frac{\beta}{p_2}+\frac{\gamma}{p_3}\equiv 0\mmod\Z\}
\end{eqnarray*}
\end{theorem}

\begin{remark}\label{t10.25}
Kluitmann \cite[III 2.4--2.6]{Kl87} calculated for the simple elliptic
singularities the group $\Aut(H_\Z,I^{(0)},M)$, which contains $G_\Z$. 
Comparison with part (c) above gives equality 
$\Aut(H_\Z,I^{(0)},M)= G_\Z$. Kluitmann described the group by an
exact sequence
\begin{eqnarray*}
\{\id\}\to \ker \pr^{\neq 1}&\to& \Aut(H_\Z,I^{(0)},M)\\
&\stackrel{\pr^{\neq 1}}{\to}& \Aut(H_{\C,\neq 1}\cap H_\Z,I^{(0)}|_{..},
M|_{..})\to\{\id\}.
\end{eqnarray*}
The group on the right side is finite of order 1296, 768 respectively 864 in 
the case $\www{E}_6,\www{E_7}$ respectively $\www{E}_8$.
The kernel $\ker \pr^{\neq 1}$ is isomorphic to 
$\Gamma(3),\Gamma(4)$ respectively 
$\Gamma(6)\subset SL_2(\Z)$ in the case $\www{E}_6,\www{E_7}$ 
respectively $\www{E}_8$.
\end{remark}

In the cases of the other families of singularities with modality $\leq 2$,
namely the hyperbolic singularities, 
the six families of quadrangle singularities 
and the eight bimodal series, we restrict here to the rough information
in the following tables and refer to \cite[Theorem 3.1]{GH17}
and \cite[Theorem 5.1 and Theorem 6.1]{GH18} for details.
The characteristic polynomial $p_{ch,-M}$ for the family $T_{pqr}$
of hyperbolic singularities 
with $p\geq q\geq r$ and $\frac{1}{p}+\frac{1}{q}+\frac{1}{r}<1$ is
\begin{eqnarray*}
p_{ch,-M}(t)=(t^p-1)(t^q-1)(t^r-1)(t-1)^{-1}.
\end{eqnarray*}
The next table from \cite[(5.1)]{GH18} gives 
the characteristic polynomial $p_{ch,-M}$
for the eight bimodal series as a product $b_1b_2$ respectively
$b_1b_2b_3$ in the case $Z_{1,p}$. The meaning of $m$ and $r_I$
is explained in Theorem \ref{t10.26}. 
\begin{eqnarray}\label{10.26}
\begin{array}{lllllll}
\textup{series} & n & b_1 & b_2 & b_3 & m & r_I\\   \hline 
W_{1,p}^\sharp & 15+p & \Phi_{12}& (t^{12+p}-1)/\Phi_1& - & 12 & 1\\
S_{1,p}^\sharp & 14+p & \Phi_{10}\Phi_2 & (t^{10+p}-1)/\Phi_1 & - & 10 & 1\\
U_{1,p} & 14+p & \Phi_9 & (t^{9+p}-1)/\Phi_1 & - & 9 & 1 \\
E_{3,p} & 16+p & \Phi_{18}\Phi_2 & t^{9+p}+1& - & 18 & 2\\
Z_{1,p} & 15+p & \Phi_{14}\Phi_2 & t^{7+p}+1 & \Phi_2 & 14 & 2\\
Q_{2,p} & 14+p & \Phi_{12}\Phi_4\Phi_3 & t^{6+p}+1 & - & 12 & 2\\
W_{1,p} & 15+p & \Phi_{12}\Phi_6\Phi_3\Phi_2 & t^{6+p}+1 & - & 12 & 2 \\
S_{1,p} & 14+p & \Phi_{10}\Phi_5\Phi_2 & t^{5+p}+1 & - & 10 & 2
\end{array}
\end{eqnarray}

\begin{theorem}\cite{He11}\cite{GH17}\cite{GH18}\label{t10.26}

(a) (i) Within the families of singularities with modality $\leq 2$ the
monodromy is not finite only in the hyperbolic singularities. 
It has one $2\times 2$ Jordan block with eigenvalue $-1$ for
the hyperbolic singularities. 

(ii) In the bimodal series, $H_\Z$ contains with index $r_I$ a
direct sum $H^{(1)}\oplus H^{(2)}$ respectively 
$H^{(1)}\oplus H^{(2)}\oplus H^{(3)}$ for $Z_{1,p}$
of Orlik blocks with characteristic polynomials
$b_1(-t)$ and $b_2(-t)$ and in the case $Z_{1,p}$ $b_3(-t)$. 

(iii) Within the bimodal series, the eigenvalue $\zeta:=e^{2\pi i /m}$
has multiplicity 2 exactly in the eight subseries with $m\,|\, p$, so
the subseries 
$E_{3,18p},Z_{1,14p},W_{1,12p},W_{1,12p}^\sharp,Q_{2,12p},S_{1,10p},
S_{1,10p}^\sharp,U_{1,9p}$.
In these cases $G_\Z$ contains automorphisms which act nontrivially
on the 2-dimensional eigenspace $H_{\C,\zeta}$ and which do not exist
for the other cases. 

(b) The quotient $G_\Z/\{\pm M^l\,|\, l\in\Z\}$ looks roughly
as follows in the families of singularities with modality $\leq 2$. 
\begin{eqnarray*}
\begin{array}{ll}
\textup{Singularity family} & 
G_\Z/\{\pm M^l\, |\, l\in\Z\} \\ \hline 
\textup{ADE-singularities} & 
\{\id\}\textup{ or }S_2\textup{ or }S_3  \\  
\textup{simple elliptic sing.} & 
\textup{a finite extension of }SL(2,\Z) \\
\textup{hyperbolic sing.} & \textup{a finite group} \\
\textup{exc. unimodal sing.} & \{\id\}\textup{ or }S_2\textup{ or }S_3  \\
\textup{exc. bimodal sing.} & \{\id\}\textup{ or }S_2\textup{ or }S_3  \\
\textup{quadrangle sing.} &  \textup{a triangle group}  \\
\textup{the 8 series, for }m\not|p & \textup{a cyclic finite group}  \\
\textup{the 8 subseries with }m|p & \textup{an infinite Fuchsian group}  
\end{array}
\end{eqnarray*}
(This table is taken from the introduction in \cite{GH18}.)
\end{theorem}

Finally, we come to the subgroup $G_\Z^{\BB}\subset G_\Z$. 

\begin{remarks}\label{t10.27}
In \cite{He11} a subgroup $G^{mar}\subset G_\Z$ is defined
for any isolated hypersurface singularity. It is called 
{\it $\mu$-constant monodromy group}. It comes from 
transversal monodromies along $\mu$-constant families
and from $\pm\id$. It has the property
\begin{eqnarray*}
G^{mar}\subset G_\Z^\BB \subset G_\Z
\end{eqnarray*}
(this is stated, though not really explained in 
\cite[Remark 3.4]{He11}). 
In \cite{He11}, \cite{GH17} and \cite{GH18} $G_\Z$ and 
$G^{mar}$ are determined for all singularities with modality $\leq 2$. 
It turns out that $G^{mar}=G_\Z^{\BB}=G_\Z$ for almost all of them,
namely all except the subseries with $m\,|\, p$ of the
eight bimodal series. 
For these subseries $G^{mar}\subsetneqq G_\Z$, and it is not
clear where in between $G_\Z^{\BB}$ is.
\end{remarks}

\section{Monodromy groups and vanishing cycles}\label{s10.4}

Consider the unimodular bilinear lattice $(H_\Z,L,\BB)$ 
with set of distinguished bases $\BB$ from an isolated
hypersurface singularity $f:(\C^{m+1},0)\to(\C,0)$.
The distinguished bases induce the even and odd monodromy
groups $\Gamma^{(0)}$ and $\Gamma^{(1)}$ and the sets 
$\Delta^{(0)}$ and $\Delta^{(1)}$ of even and odd vanishing
cycles. These two groups and these two sets have been studied
a lot and by many people. In the even case the pair $(H_\Z,I^{(0)})$ alone
determines $\Gamma^{(0)}$ and $\Delta^{(0)}$ in almost all cases.
This is satisfying. The precise results below are due to Ebeling.
In the odd case, the situation is more complicated. 
It was worked upon by Wajnryb, Chmutov and Janssen.
In the following we report first on the
even case and then on the odd case.

The lattice $(H_\Z,I^{(0)})$ with even intersection form is 
an a priori rather rough invariant. Lattices with even
intersection forms are well understood, due to work of Durfee,
Kneser, Nikulin, Wall and others. 

The following theorem of Ebeling was developed by him in
several papers, in greater and greater generality, until
the final state in \cite[(5.5) and (2.5)]{Eb84}. 
The version here is taken from 
\cite[Theorem 5.9]{Eb01}, as \cite{Eb84} does not cover 
the characterization of $\Delta^{(0)}$ for 
some special cases.

\begin{theorem}\label{t10.28} 
\cite[Theorem 5.9]{Eb01}
Suppose that the singularity $f$ is a hyperbolic singularity
of type $T_{pqr}$ with 
$(p,q,r)\in\{(2,3,7),(2,4,5),(3,3,4)\}$
(so these three triples are allowed, all other triples
$(p,q,r)$ with $p\geq q\geq r$ and 
$\frac{1}{p}+\frac{1}{q}+\frac{1}{r}<1$ are not allowed)
or any not hyperbolic singularity. 
Then 
\begin{eqnarray*}
\Gamma^{(0)}&=&O^{(0),*},\\
\Gamma^{(0)}_u&=&(\ker\tau^{(0)})_u
\stackrel{\ref{t5.2}\ (e)}{=}
T(\oooo{j}^{(0)}(\oooo{H_\Z}^{(0)})\otimes \Rad I^{(0)}),\\
\Gamma^{(0)}_s&=& \ker\oooo{\tau}^{(0)}\cap\ker\oooo{\sigma},\\
{}[O^{(0)}_s:\Gamma^{(0)}_s]&<&\infty,\\
\Delta^{(0)}&=&\{a\in R^{(0)}\,|\, I^{(0)}(a,H_\Z)=\Z\},
\end{eqnarray*}
(recall Lemma \ref{t6.2} (f) for $\oooo{\tau}^{(0)}$ and
Remark \ref{t6.3} (iii) for $\oooo{\sigma}$), 
the exact sequence
\begin{eqnarray*}
\{\id\}\to\Gamma^{(0)}_u\to\Gamma^{(0)}\to\Gamma^{(0)}_s\to\{\id\}
\end{eqnarray*}
splits non-canonically.
\end{theorem}

\begin{remarks}\label{t10.29}
(i) The characterizations of $\Gamma^{(0)}$ and $\Delta^{(0)}$
in the theorem show that they are determined 
by $(H_\Z,I^{(0)})$ alone. 

(ii) $\Gamma^{(0)}_u=(\ker\tau^{(0)})_u$ follows from
$\Gamma^{(0)}=O^{(0),*}\stackrel{Def.}{=}
\ker\tau^{(0)}\cap\ker\sigma$ and 
$O^{(0),Rad}_u\subset\ker\sigma$. 
The simple part $\Gamma^{(0)}_s$ has finite index in
$O^{(0)}_s$ because 
$\oooo{H_\Z}^{(0),\sharp}/\oooo{j}^{(0)}(\oooo{H_\Z}^{(0)})$
is a finite abelian group and therefore
$\ker\oooo{\tau}^{(0)}$ has finite index in 
$O^{(0)}_s$. 

(iii) We do not know a singularity where
$R^{(0)}\supsetneqq \Delta^{(0)}$. 
Equivalent: We do not know a singularity where a root
$a\in R^{(0)}$ with $I^{(0)}(a,H_\Z)\subsetneqq \Z$ exists.

(iv) The first idea of a proof of
$\Delta^{(0)}=\{a\in R^{(0)}\,|\, I^{(0)}(a,H_\Z)=\Z\}$
is due to Looijenga. It is given in 
\cite[Korollar 1 in 4.2]{Br83}.
This paper treats the exceptional unimodal singularities.
But it gives also many useful references and general facts.

(v) In the case of a hyperbolic singularity $T_{pqr}$
with $(p,q,r)\notin\{(2,3,7),(2,4,5),(3,3,4)\}$,
Gabrielov \cite{Ga74-2} determined $\Gamma^{(0)}$.
Then especially $[O^{(0)}_s:\Gamma^{(0)}_s|=\infty$
(see also \cite[(5.2)]{Eb84}).
For example the singularities $T_{2,7,7}$ and $T_{3,3,10}$
have isomorphic pairs $(H_\Z,I^{(0)})$, but the monodromy
groups for all hyperbolic singularities are pairwise different
\cite{Br83}. So in the case of the hyperbolic singularities,
the pair $(H_\Z,I^{(0)})$ does not determine the monodromy
group $\Gamma^{(0)}$.
See \cite[7.2]{BWY19} for a statement on $\Delta^{(0)}$
for the hyperbolic singularities. 
 
(vi) An example \cite[Example after Theorem 3.3]{Eb83}: 
The exceptional bimodal singularities $E_{18}$ and $Q_{18}$
have isomorphic pairs $(H_\Z,I^{(0)})$. Both are
isomorphic to the orthogonal sum 
$E_6\perp E_8\perp U\perp U$ where $E_6$ and $E_8$ mean the
root lattices and $U$ means the rank 2 
hyperbolic lattice, \index{hyperbolic lattice}
so with matrix $\begin{pmatrix}0&1\\1&0\end{pmatrix}$
for its even intersection form.
Therefore also $\Gamma^{(0)}$ and $\Delta^{(0)}$ are
isomorphic for these singularities. But Seifert form $L$
and monodromy $M$ are not isomorphic. 
The monodromy $M$ has order 30 
for $E_{18}$ and order 48 for $Q_{18}$. 
\end{remarks}

Now we come to the odd case.
The lattice $(H_\Z,I^{(1)})$ with odd intersection form 
has a basis $\uuuu{v}$ whose intersection matrix is
\begin{eqnarray*}
I^{(1)}(\uuuu{v}^t,\uuuu{v})= 
\begin{pmatrix}
0   &d_1&    &   &      &    &   & &      & \\
-d_1&0  &    &   &      &    &   & &      & \\
    &   &0   &d_2&      &    &   & &      & \\
    &   &-d_2&0  &      &    &   & &      & \\
    &   &    &   &\ddots&    &   & &      & \\
    &   &    &   &      &0   &d_l& &      & \\
    &   &    &   &      &-d_l&0  & &      & \\
    &   &    &   &      &    &   &0&      & \\
    &   &    &   &      &    &   & &\ddots& \\
    &   &    &   &      &    &   & &      &0
\end{pmatrix}
\end{eqnarray*} 
(at empty places put 0) where $d_1,...,d_l\in\N$ with 
$d_l\,|\, d_{l-1}\,|\, ...\,|\, d_1$. 
This is a linear algebra fact. The numbers $l$ and 
$d_1,...,d_l$ are unique. The last $\rk\Rad I^{(1)}$
vectors of $\uuuu{v}$ are a basis of the radikal
$\Rad I^{(1)}$ of $I^{(1)}$, so $2l+\rk\Rad I^{(1)}=n$. 
The basis $\uuuu{v}$ is called a {\it symplectic basis}. 
\index{symplectic basis}
Define the sublattice of rank $n$
\begin{eqnarray*}
H^{(d_1)}_\Z&:=& \{a\in H_\Z\,|\, I^{(1)}(a,H_\Z)\subset
d_1\Z\},
\end{eqnarray*}
and the finite abelian quotient group
\begin{eqnarray*}
H^{quot}&:=& H^{(d_1)}_\Z/2d_1H_\Z.
\end{eqnarray*}
Each automorphism in $O^{(1),Rad}$ respects
$H^{(d_1)}$ and $2d_1H_\Z$, so it acts on the quotient
$H^{quot}$. The induced finite group of automorphisms of $H^{quot}$ 
is called $O^{(1),quot}$, the homomorphism is called
$p^{quot}:O^{(1),Rad}\to O^{(1),quot}$. 
Each element of $O^{(1),quot}$ acts trivially on the
image of the radical $\Rad I^{(1)}$ in the quotient
$H^{quot}$. 

The following theorem cites four results of Chmutov in 
\cite{Ch82}:\\ 
(i) A relative characterization of $\Gamma^{(1)}$,\\
(ii) a relative characterization of $\Delta^{(1)}$,\\
(iii) that $\Gamma^{(1)}_s$ has finite index in $O^{(1),Rad}_s$,
\\
(iv) and a characterization of the subgroup
$\ker p^{quot}\subset \Gamma^{(1)}$.

\begin{theorem}\label{t10.30}
\cite[Theorem 1, Corollary of Proposition 1, Theorem 2,
Proposition 2]{Ch82}

(i) $\Gamma^{(1)}$ is the full preimage in $O^{(1),Rad}$
of $p^{quot}(\Gamma^{(1)})$ in $O^{(1),quot}$. 
Equivalent: 
\begin{eqnarray*}
\ker \bigl(p^{quot}:O^{(1),rad}\to O^{(1),quot}\Bigr)\subset\Gamma^{(1)}.
\end{eqnarray*}
So if one knowns
the image $p^{quot}(\Gamma^{(1)})$, one knows $\Gamma^{(1)}$. 

(ii) 
\begin{eqnarray*}
\Delta^{(1)}=\{a\in H_\Z&|& \textup{there exists }
b\in \Delta^{(1)} \textup{ with }\\
&& a-b\in 2H_\Z\textup{ and }I^{(1)}(a,H_\Z)=\Z\}.
\end{eqnarray*}
So if one knows the image of $\Delta^{(1)}$ in
$H_\Z/2 H_\Z$, one knows $\Delta^{(1)}$. 

(iii) 
\begin{eqnarray*}
\{g\in O^{(1),Rad}_s\,|\, g\textup{ acts trivially on the
quotient }\oooo{H_\Z}^{(1)}/2d_1\oooo{H_\Z}^{(1)}\}
\subset \Gamma^{(1)}_s.
\end{eqnarray*}
As the group on the left hand side has finite index in
$O^{(1),Rad}_s$, also $\Gamma^{(1)}_s$ has finite index in
$O^{(1),Rad}_s$. 

(iv)
\begin{eqnarray*}
\langle (s^{(1)}_a)^2\,|\, a\in H_\Z\rangle &=&
\ker\bigl(p^{quot}:O^{(1),rad}\to O^{(1),quot}\Bigr)
\stackrel{(i)}{\subset}\Gamma^{(1)}.
\end{eqnarray*}
\end{theorem}

\begin{remarks}\label{t10.31}
(i) In the case of a curve singularity (so $m=1$) the number
of components is $\rk\Rad I^{(1)}+1$. 
So the curve singularity is irreducible if and only if 
$I^{(1)}$ is nondegenerate. 

(ii) In the case of a curve singularity $d_1=...=d_l=1$.

(iii) Wajnryb \cite{Wa80} proved Theorem \ref{t10.30} for irreducible
curve singularities. Chmutov \cite{Ch81} generalized it to
all curve singularities. Here $d_1=1$ implies
$H^{(d_1)}=H_\Z$ and $H^{quot}=H_\Z/2H_\Z$. 

(iv) Wajnryb characterized in the case of irreducible
curve singularities also the image 
$p^{quot}(\Gamma^{(1)})\subset O^{(1),quot}$
and the image of $\Delta^{(1)}$ in $H_\Z/2H_\Z$. 
Also this was generalized by Chmutov first to all curve
singularities \cite{Ch81} and then to all singularities
\cite{Ch83}. 
\end{remarks}

First we discuss the characterization in \cite{Ch83} 
of the image
$$\Delta^{(1)}_{\F_2}:= \textup{image of }\Delta^{(1)}
\textup{ in }H_\Z/2H_\Z$$
of $\Delta^{(1)}$ in $H_\Z/2H_\Z$.
The following Lemma is elementary. It is essentially formulated
in \cite[3.1]{Ch83}, extending an observation of Wajnryb
\cite{Wa80}.

\begin{lemma}\label{t10.32}
Let $\www{H}_{\F_2}$ be an $\F_2$-vector space of dimension 
$\www{n}\in\Z_{\geq 2}$, let $\www{I}^{(1)}$ be an odd bilinear
form $\www{I}^{(1)}:\www{H}_{\F_2}\times \www{H}_{\F_2}\to\F_2$ 
on $\www{H}_{\F_2}$
({\sf odd} means here only $\www{I}^{(1)}(a,a)=0$ for $a\in\www{H}_{\F_2}$), 
and let $\uuuu{v}$ be a basis of $\www{H}_{\F_2}$. 

(a) (Definition) A {\sf quadratic form} $q:\www{H}_{\F_2}\to\F_2$
is a map with 
$$q(x+y)=q(x)+q(y)+\www{I}^{(1)}(x,y)\quad\textup{for all }x,y
\in \www{H}_{\F_2}.$$

(b) There is a unique quadratic form $q$ with $q(v_j)=1$ for
all $j\in\{1,...,n\}$. The value $q(x)$ for $x\in\www{H}_{\F_2}$ 
can be characterized as follows.
Write $x=\sum_{j=1}^{\www{n}} a_jv_j$ with $a_j\in\F_2$. Define a graph
with vertex set $V(x):=\{j\in\{1,...,n\}\,|\, a_j=1\}$ and edge set
$E(x):=\{(i,j)\in\{1,...,n\}^2\,|\, i<j,\www{I}^{(1)}(v_i,v_j)=1\}$.
Then $q(x)$ is the Euler characteristic of the graph $\mmod 2$ so 
$$q(x)\equiv |V(a)|-|E(a)|\ \mmod 2.$$
\end{lemma}

The next theorem collects results of Chmutov and Janssen. 
The main point is part (b). It characterizes $\Delta^{(1)}_{\F_2}$. 

\begin{theorem}\label{t10.33}
Suppose that the singularity $f$ is neither the singularity $A_n$
nor the singularity $D_n$. Choose a distinguished basis 
$\uuuu{v}\in\BB^{dist}$. Consider the quadratic form $q$ on
$H_\Z/2H_\Z$ with $q(v_j)=1$ for any $j\in\{1,...,n\}$
(it exists and is unique by Lemma \ref{t10.32}).

(a) \cite[7.1]{Ch83} The quadratic form $q$ on $H_\Z/2H_\Z$ is independent 
of the choice of the distinguished basis $\uuuu{v}\in\BB^{dist}$.

(b) \cite[Proposition 1]{Ch83} 
$$\Delta^{(1)}_{\F_2}=q^{-1}(1)-(\textup{image of }
\Rad I^{(1)}\textup{ in }H_\Z/2H_\Z).$$

(c) \cite[(5.1) and (5.2)]{Ja83} The subset 
$$(\Rad I^{(1)})^q:=\{a\in\Rad I^{(1)}\,|\, q(\textup{image of }a
\textup{ in }H_\Z/2H_\Z)=0\}$$ 
of the radical $\Rad I^{(1)}$ is either equal to $\Rad I^{(1)}$
or a subgroup of index 2 in $\Rad I^{(1)}$. In any case it is the
set $\{a\in\Rad I^{(1)}\,|\, a+b\in\Delta^{(1)}\}$, where
$b\in\Delta^{(1)}$ is an arbitrarily chosen odd vanishing cycle. 

(d) \cite[(5.3)]{Ja83} 
$$T(\oooo{j}^{(1)}(\oooo{H_\Z}^{(1)})\otimes (\Rad I^{(1)})^q)
\subset \Gamma^{(1)}_u.$$
\end{theorem}

\begin{remarks}\label{t10.34}
The singularities $A_n$ and $D_n$ play a special role 
in the odd case, in \cite{Ch83} and in \cite{Ja83}\cite{Ja85}.
They are the only singularities which have no deformation
to the singularity $E_6$. The groups $\Gamma^{(1)}$ and
the sets $\Delta^{(1)}$ for them had been determined by
Varchenko. Chmutov offers an account of this at the end of
\cite{Ch83}.
\end{remarks}

It remains to characterize the finite group
$p^{quot}(\Gamma^{(1)})$. The following holds for 
curve singularities (irreducible curve singularities without
$A_n$ and $D_n$ \cite[Theorem 4 (c)]{Wa80}, arbitrary curve
singularities \cite[Theorems 1 and 3]{Ch83}).

\begin{theorem}\label{t10.35}
Suppose that $f$ is a curve singularity. 
Recall that then $d_1=1$ and $H^{quot}=H_\Z/2H_\Z$.
The following three properties for an element 
$g\in O^{(1),quot}$ are equivalent:
\begin{list}{}{}
\item[(i)]
$g\in p^{quot}(\Gamma^{(1)})$.
\item[(ii)]
$g$ respects the quadratic form $q$ in Theorem \ref{t10.33}.
\item[(iii)]
$g$ maps $\Delta^{(1)}_{\F_2}$ to itself.
\end{list}
\end{theorem}

\begin{remarks}\label{t10.36}
(i) So, such an automorphism $g$ maps in the disjoint union
$$H_\Z/2H_\Z=(\textup{image of }\Rad I^{(1)})\,\dot\cup\,
\Delta^{(1)}_{\F_2}\,\dot\cup\,(\textup{the complement})$$
each of the three sets to itself, and $g$ is the identity on the
first set. 

(ii) Chmutov found a generalization of this to arbitrary singularities.
In Section 5 in \cite{Ch83} an equivalence relation for
the elements of $H^{quot}$ for an arbitrary singularity is
defined. Theorem 1 in \cite{Ch83} says that for an arbitrary
singularity an automorphism $g\in O^{(1),quot}$ is in
$p^{quot}(\Gamma^{(1)})$ if and only if it respects this 
equivalence relation. For details see \cite[Section 5]{Ch83}.
\end{remarks}

Janssen \cite{Ja83}\cite{Ja85} recovered a good part of
Chmutov's results. Whereas Chmutov worked mainly within 
the singularity cases, Janssen took a more abstract point of view.
He started with the notion of an {\it odd vanishing lattice}.
He classified odd vanishing lattices of $\F_2$ and over $\Z$.

\begin{definition}\label{t10.37}
Let $R=\Z$ or $R=\F_2$. 
An {\sf odd vanishing lattice} 
\index{odd vanishing lattice}\index{vanishing lattice}
over $R$ is a triple
$(V,I_V,\Delta_V)$ with the following properties.
$V$ is a free $R$-module of some rank $n_V\in\Z_{\geq 2}$.
$I_V:V\times V\to R$ is an odd bilinear form.
$\Delta_V$ is a subset of $V$ with three specific properties.
Define the group $\Gamma_{\Delta_V}:=
\langle s^{[1]}_a\,|\, a\in\Delta_V\rangle$ where
$s^{[1]}_a\in O(V,I_V)$ is as usual the transvection with
$s^{[1]}_a(b):=b-I_V(a,b)a$ for $b\in V$. 
The three specific properties:
\begin{list}{}{}
\item[(i)]
$\Delta_V$ is a single $\Gamma_{\Delta_V}$-orbit.
\item[(ii)]
$\Delta_V$ generates $V$ as $R$-module.
\item[(iii)]
There exist $a_1$ and $a_2\in\Delta_V$ with $I_V(a_1,a_2)=1$.
\end{list} 
\end{definition}

\begin{remarks}\label{t10.38}
(i) It is true that for any singularity with $n\geq 2$ the triple
$(H_\Z,I^{(1)},\Delta^{(1)})$ is an odd vanishing lattice
over $\Z$. 

(ii) Janssen classified all odd vanishing lattices over
$\F_2$ \cite[(4.8) Theorem]{Ja83} and over $\Z$ 
\cite[(7.8) Theorem]{Ja85}.
There are surprisingly few families, only seven families,
for $R=\F_2$ as well as for $R=\Z$

(iii) In the case of $R=\Z$, 
in each of the seven families, an odd vanishing lattice 
$(V,I_V,\Delta_V)$ is determined up to isomorphism by
the invariants $l,d_1,...,d_l$ and $\rk\Rad I_V$
of the pair $(V,I_V)$, and additionally in two of the
seven families by one more number.

(iv) Except for the $A_n$ and $D_n$ singularities,
for any singularity $f$ the triple $(H_\Z,I^{(1)},\Delta^{(1)})$
fits into one of the three families of odd vanishing lattices
with symbols $O^\sharp_1(d_1,...,d_l,\rho)$, 
$O^\sharp_0(d_1,...,d_l,\rho)$ and $O^\sharp(d_1,...,d_l,\rho,k_0)$
where $\rho=\rk \Rad I^{(1)}$ and $k_0$ is the additional number
\cite[(6.6) Theorem]{Ja83}\cite{Ja85}

(v) But even if one knows that for a singularity 
the triple $(H_\Z,I^{(1)},\Delta^{(1)})$ is for example of some
type $O^\sharp_1(d_1,...,d_l,\rho)$ where all invariants are 
determined by the pair $(H_\Z,I^{(1)})$, this does not mean
that the pair $(H_\Z,I^{(1)})$ determines $\Delta^{(1)}$ uniquely.
There might by an automorphism of the pair $(H_\Z,I^{(1)})$
which maps $\Delta^{(1)}$ to a different set. 
Then the pair $(H_\Z,I^{(1)})$ determines the triple 
$(H_\Z,I^{(1)},\Delta^{(1)})$ only up to isomorphism.

(vi) The singularities $A_6$ and $E_6$ exist both as
irreducible curve singularities, so the pairs
$(H_\Z,I^{(1)})$ have the same invariants $d_1=d_2=d_3=1$
and $\rk\Rad I^{(1)}=0$, so they are isomorphic.
But the sets $\Delta^{(1)}$, the groups $\Gamma^{(1)}$
and the quadratic forms $q$ on $H_\Z/2H_\Z$ are very different.
\end{remarks}

\begin{example}\label{t10.39}
In \cite[5.2 and 5.4]{DM94} several examples of pairs $(f_1,f_2)$ of
curve singularities with the following property are given.
The pairs $(H_\Z,L)$ of Milnor lattice and normalized
Seifert form are isomorphic for $f_1$ and $f_2$, but
$f_1$ and $f_2$ have distinct topological type. 
In all examples $f_1$ and $f_2$ are both reducible,
with two components each (so $\rk\Rad I^{(1)}=1$).

This leads naturally to the following questions
for each of these pairs of curve singularities.
We do no know their answers.
\begin{list}{}{}
\item[(A)] Are the triples $(H_\Z,I^{(1)},\Delta^{(1)})$
isomorphic?
\item[(B)] Are the triples $(H_\Z,L,\Delta^{(1)})$
isomorphic?
\item[(C)] Are the triples $(H_\Z,L,\BB^{dist})$
isomorphic?
\end{list}
\end{example}

\section[Moduli spaces]{Moduli spaces for the simple and the simple
elliptic singularities}\label{s10.5}

The only cases where $C_n^{\uuuu{e}/\{\pm 1\}^n}$ and
$C_n^{S/\{\pm 1\}^n}$ had been studied in detail up to now
are the cases from the simple singularities and the 
simple elliptic singularities, the cases from the simple
singularities by Looijenga \cite{Lo74} and (implicitly) 
Deligne \cite{De74},
all cases by Hertling and Roucairol \cite{HR21}.
Results on these cases are recalled and explained here.
First we discuss the simple singularities, then the
simple elliptic singularities.

The simple singularities in the normal forms in
Theorem \ref{t10.6} are quasihomogeneous polynomials.
They have universal unfoldings which are also quasihomogeneous
polynomials. Here the good representatives and the base spaces
$\MM$ can be chosen as algebraic and global objects.
We reproduce the universal unfoldings which are
given in \cite{Lo74}:
\index{unfolding of a simple singularity}
\begin{eqnarray}\label{10.27}
&&F^{alg}:\C^{m+1}\times \MM^{alg}\to\C\quad\textup{with}\quad 
\MM^{alg}=\C^n,\\
&&F^{alg}(z_0,...,z_m,t_1,...,t_n)=F^{alg}(z,t)=F^{alg}_t(z)
=f(z)+\sum_{j=1}^n t_jm_j \nonumber
\end{eqnarray}
with $f=F^{alg}_0$ and $m_1,...,m_n$ the monomials in the tables
\eqref{10.28} and \eqref{10.29},
\begin{eqnarray}\label{10.28}
\begin{array}{ccccccc}
\textup{name}  & m_1 & m_2 & m_3 & m_4 & ... & m_n \\ \hline 
A_n & 1 & z_0 & z_0^2 & z_0^3 & ... & z_0^{n-1} \\
D_n & 1 & z_1 & z_0 & z_0^2 & ... & z_0^{n-2} 
\end{array}
\end{eqnarray}
\begin{eqnarray}\label{10.29}
\begin{array}{ccccccccc}
\textup{name}  & m_1 & m_2 & m_3 & m_4 & m_5 
& m_6 & m_7 & m_8  \\ \hline
E_6 & 
1 & z_0 & z_1 & z_0^2 & z_0z_1 & z_0^2z_1 & &   \\
E_7 & 
1 & z_0 & z_1 & z_0^2 & z_0z_1 & z_0^3 & z_0^4 &  \\
E_8 & 
1 & z_0 & z_1 & z_0^2 & z_0z_1 & z_0^3 & z_0^2z_1 & z_0^3z_1 
\end{array}
\end{eqnarray}
One checks easily that the monomials form a basis 
of the Jacobi algebra $\OO_{\C^{m+1},0}/J_f$.
Therefore the unfolding $F^{alg}$ is indeed universal.
The following tables list the weights $\deg_{\bf w}t_j$ 
of the unfolding
parameters $t_1,...,t_n$ which make $F$ into a quasihomogeneous
polynomial of weighted degree 1. The weights are all positive.
The tables list also the Coxeter number $N_{Coxeter}$ of the
corresponding root system.
\begin{eqnarray*}
\begin{array}{l|l|l|l|l|l|l|l|l|l|l|l}
 & N_{Coxeter} & z_0 & z_1 & t_1 & t_2 & t_3 & t_4 & ... & 
 t_n \\
A_n & n+1 & \frac{1}{n+1} & & 1 & \frac{n}{n+1} & \frac{n-1}{n+1} 
& \frac{n-2}{n+1} & ... & \frac{2}{n+1} \\
D_n & 2(n-1) & \frac{1}{n-1} & \frac{n-2}{2(n-1)} & 1 & \frac{n}{2(n-1)}
& \frac{n-2}{n-1} & \frac{n-3}{n-1} & ... & \frac{1}{n-1} 
\end{array}
\end{eqnarray*}
\begin{eqnarray*}
\begin{array}{l|l|l|l|l|l|l|l|l|l|l|l}
 & N_{Coxeter} & z_0 & z_1 & t_1 & t_2 & t_3 & t_4 & t_5 & t_6 & t_7 & t_8 \\
E_6 & 12 & \frac{1}{4} & \frac{1}{3} & 1 & \frac{3}{4} 
& \frac{2}{3} & \frac{1}{2} & \frac{5}{12} & \frac{1}{6} & & 
\\[1mm] 
E_7 & 18 & \frac{2}{9} & \frac{1}{3} & 1 & \frac{7}{9} 
& \frac{2}{3} & \frac{5}{9} & \frac{4}{9} & \frac{1}{3} 
& \frac{1}{9} & \\[1mm]
E_8 & 30 & \frac{1}{5} & \frac{1}{3} & 1 & \frac{4}{5} 
& \frac{2}{3} & \frac{3}{5} & \frac{7}{15} & \frac{2}{5} 
& \frac{4}{15} & \frac{1}{15} 
\end{array}
\end{eqnarray*}
Because the weights are positive, $F_t$ for each $t$ 
satisfies that the sum of Milnor numbers of its singular
points is equal to $n$. It follows with little work
that $\MM$ is an F-manifold with 
$$\textup{unit field }e=\paa_1,\quad 
\textup{Euler field } 
E=\sum_{j=1}^n\deg_{\bf w}(t_j)t_j\paa_j,$$
and polynomial multiplication \cite[Theorem 4.3]{HR21}.

The positive weights and general properties of the 
Lyashko-Looijenga map lead easily to the following result
which was found independently by Looijenga and Lyashko
(but Lyashko published his version only much later).

\begin{theorem}\label{t10.40}
\cite{Lo74}\cite{Ly79}
The Lyashko-Looijenga map 
\begin{eqnarray}\label{10.30}
\LL^{alg}:\MM^{alg} \to \C[x]_n\cong\C^n
\end{eqnarray}
is a branched covering of degree
\begin{eqnarray}\label{10.31}
\deg \LL^{alg}=\frac{n!}{\prod_{j=1}^n\deg_{\bf w}t_j}.
\end{eqnarray}
It is branched along the caustic $\KK_3^{alg}\subset\MM^{alg}$
(at generic points of order 3)
and along the Maxwell stratum $\KK_2^{alg}\subset\MM^{alg}$
(at generic points of order 2). 
The restriction 
\begin{eqnarray}\label{10.32}
\LL^{alg}:\MM^{alg}-(\KK_3^{alg}\cup
\KK_2^{alg})\to \C[x]_n^{reg}\cong C_n^{conf}
\end{eqnarray}
is a covering.
\end{theorem}

\begin{remarks}\label{t10.41}
(i) It is known that 
$|\Gamma^{(0)}|=(N_{Coxeter})^n\prod_{j=1}^n\deg_{\bf w}t_j$
in the case of an ADE root lattice. Therefore
$$\deg \LL^{alg}= \frac{n!\cdot (N_{Coxeter})^n}
{|\Gamma^{(0)}|}.$$

(ii) The following table lists the numbers $|G_\Z|$ 
and $|\SSS^{dist}/\{\pm 1\}^n|$ for the simple singularities. 
Theorem \ref{t10.42} and \eqref{10.39} will tell
\begin{eqnarray*}
\deg\LL^{alg}=|\BB^{dist}/\{\pm 1\}^n|
=\frac{1}{2}\cdot |G_\Z|\cdot |\SSS^{dist}/\{\pm 1\}^n|.
\end{eqnarray*}
\begin{eqnarray}\label{10.33}
\begin{array}{lll}
  & |G_\Z| & |\SSS^{dist}/\{\pm 1\}^n| \\ 
\hline
A_n & 2(n+1) & (n+1)^{n-2}\\
D_4 & 36 & 9\\
D_n,n\geq 5 & 4(n-1) & (n-1)^{n-1}\\
E_6 & 24  & 2^7\cdot 3^3 = 3456\\
E_7 & 18  & 2\cdot 3^{10} = 118098\\
E_8 & 30  & 2\cdot 3^4\cdot 5^6 = 2531250
\end{array}
\end{eqnarray}
\end{remarks}

Recall from Definition \ref{t8.5} the open convex polyhedron
$F_n\subset C_n^{pure}\subset\C^n$ and its boundary
$W_n$, a union of walls. Also recall
$$C_n^{conf}\subset \C^n =\pr_n^{p,c}(F_n)\,\dot\cup\, 
\pr_n^{p,c}(W_n).$$
The union of walls $(\LL^{alg})^{-1}(\pr_n^{p,c}(W_n))$ 
upstairs in $\MM^{alg}$ contains $\KK_3^{alg}\cup \KK_2^{alg}$. 
The restriction
$$\LL^{alg}: (\LL^{alg})^{-1}(\pr_n^{p,c}(F_n))\to F_n$$
is an even covering. The components of the preimages
are called {\it Stokes regions}. Each of them is mapped
isomorphically to $F_n$. The number of Stokes regions is
$\deg \LL^{alg}$. 

We obtain a map \index{$LD$} 
\begin{eqnarray*}
LD:\{\textup{Stokes regions in }\MM^{alg}\}\to 
\BB^{dist}/\{\pm 1\}^n
\end{eqnarray*}
in the following way. For $t$ in one Stokes region choose
a standard system $(\uuuu{\gamma};\id)$ of paths and 
obtain as in Theorem \ref{t10.8} from $F_t$ a
distinguished basis up to signs. This works because
$\MM=\C^n$ is simply connected and therefore there is a 
canonical isomorphism
$$H_m(F_t^{-1}(\eta),\Z)\cong H_m(X_\eta,\Z)\cong H_\Z$$
for (depending on $t$) sufficiently large $\eta>0$. 
It is independent of the choice of $t$ within a fixed Stokes
region. The distinguished basis up to signs is
$LD(\textup{Stokes region})\in\BB^{dist}/\{\pm 1\}^n$. 

Choose one distinguished basis $\uuuu{e}\in\BB^{dist}$
for reference.
The map $LD$ is surjective because the restriction 
of $\LL^{alg}$ in \eqref{10.32} is a covering
and $\BB^{dist}/\{\pm 1\}^n$ is a single $\Br_n$ orbit.
Even more, the covering factorizes through the 
covering $\pr_n^{e,c}:C_n^{\uuuu{e}/\{\pm 1\}^n}\to C_n^{conf}$,
so it is a composition of two coverings, 
\begin{eqnarray}\label{10.34}
\MM^{alg}-(\KK_3^{alg}\cup\KK_2^{alg})
\stackrel{\pr_n^{alg,e}}{\longrightarrow} 
C_n^{\uuuu{e}/\{\pm 1\}^n}
\stackrel{\pr_n^{e,c}}{\longrightarrow} C_n^{conf}.
\end{eqnarray}
This follows from the construction of the spaces and the
compatibility of the braid group actions.
The main result in this section for the simple
singularities is the following.

\begin{theorem}\label{t10.42}\cite{Lo74}\cite{De74}
\cite[Theorem 7.1]{HR21}
Consider a simple singularity. The map
\begin{eqnarray}\label{10.35}
LD:\{\textup{Stokes regions in }\MM^{alg}\}\to 
\BB^{dist}/\{\pm 1\}^n
\end{eqnarray} 
is a bijection.
Equivalent: The covering 
\begin{eqnarray}\label{10.36}
\pr_n^{alg,e}: \MM^{alg}-(\KK_3^{alg}\cup\KK_2^{alg})
\to C_n^{\uuuu{e}/\{\pm 1\}^n}
\end{eqnarray} 
is an isomorphism.
\end{theorem}
 
Looijenga \cite{Lo74} asked whether the map $LD$ is a bijection 
and proved it in the case of the $A_n$ singularity. 
Deligne \cite{De74} answered the question for all simple
singularities positively in the following way. 
He calculated in a letter to Looijenga
for each simple singularity the number 
$|\BB^{dist}/\{\pm 1\}^n|$ and found equality with 
$\deg \LL^{alg}$. The map $LD$ is a surjection between
finite sets of same size, so it is a bijection.
This is one proof of Theorem \ref{t10.42}.

\cite{HR21} gives a different proof, which generalizes to
the simple elliptic singularities where $LD$ is a map between
infinite sets. Here we sketch part of this proof. 

We need a general result about symmetries of a singularity
\index{symmetries of a singularity} from \cite[13.2]{He02}.

Let $f:(\C^ {m+1},0)\to(\C,0)$ be an isolated hypersurface
singularity with Milnor number $n$.
Let $F:(\C^{m+1}\times \MM,0)\to(\C,0)$ 
be a universal unfolding with base space a germ 
$(\MM,0)=(\C^n,0)$. Denote by
\begin{eqnarray*}
\RR^f:=\{\varphi:(\C^{m+1},0)\to(\C^{m+1},0)&|&
\varphi\textup{ is a biholomorphic}\\
&&\textup{map germ with }f=f\circ\varphi\}
\end{eqnarray*}
the (very large) group of germs of coordinate changes
which leave $f$ invariant. 
Denote by
$$\Aut_\MM:=\Aut((\MM,0),\circ,e,E)$$
the group of automorphisms of $(\MM,0)$ as a germ of an
F-manifold with Euler field. It is a finite group
\cite[Theorem 4.14]{He02}. Consider a good representative
$F:\XX\to\Delta_\eta$ of the universal unfolding, whose 
base space $\MM\subset\C^n$ is a small ball in $\C^n$. 
The Lyashko-Looijenga map
$\LL:\MM\to\C[x]_n$ is holomorphic. The restriction
$\LL:\MM-(\KK_3\cup\KK_2)\to C_n^{conf}$ is locally 
biholomorphic. The components of the even smaller restriction
\begin{eqnarray}\label{10.37}
\LL:\LL^{-1}(\pr^{p,c}_n(F_n))\to F_n
\end{eqnarray}
are also now called {\it Stokes regions}, although
the restriction $\LL:(\textup{one Stokes region})\to F_n$
is injective, but in general not an isomorphism.
Nevertheless, as in the case of the simple singularities, 
we obtain a natural map
\begin{eqnarray}\label{10.38}
LD:\{\textup{Stokes regions in }\MM\}\to\BB^{dist}/\{\pm 1\}^n
\end{eqnarray}
by the construction in Theorem \ref{t10.8}
applied to $F_t$ for some $t$ in the Stokes region
together with a standard system of paths.

\begin{theorem}\label{t10.43}\cite[13.2]{He02}

(a) Each $\varphi\in\RR^f$
can be lifted to an automorphism of the unfolding.
It induces an automorphism $(\varphi)_{hom}\in G_\Z$
of $H_\Z$ and an automorphism $(\varphi)_{\MM}\in \Aut_\MM$. 
The group homomorphism $()_\MM:\RR^f\to\Aut_\MM$
is surjective.
The group homomorphism $()_{hom}:\RR^f\to G_\Z$ has finite
image. The image is a subgroup of $G_\Z^{\BB}$. 

(b) 
\begin{eqnarray*}
\begin{array}{cccc}
-\id\notin (\RR^f)_{hom}&\textup{ and }&
\ker()_\MM=\ker()_{hom}&\textup{ if }\mult(f)\geq 3.\\
-\id\in (\RR^f)_{hom}&\textup{ and }&
\ker()_\MM=\ker()_{hom}\times\{\pm \id\}&
\textup{ if }\mult(f)=2.
\end{array}
\end{eqnarray*}

(c) One can choose a finite subgroup $\RR^{finite}\subset\RR^f$
which lifts to a group of automorphisms of a (sufficiently) good
representative $F:\XX\to\Delta_\eta$ of the universal unfolding
and such that $()_{hom}:\RR^{finite}\to G_\Z$ is injective
with image $(\RR^f)_{hom}$. Then $(\RR^{finite})_\MM=\Aut_\MM$. 
This group acts on the set of Stokes regions of $\MM$.
The map $LD$ is equivariant, that means, for 
$\varphi\in\RR^{finite}$ the following diagram commutes,
\begin{eqnarray*}
\begin{CD}
\{\textup{Stokes regions in }\MM\} @>{LD}>> 
\BB^{dist}/\{\pm 1\}^n\\
@VV{(\varphi)_\MM}V  @VV{(\varphi)_{hom}}V \\
\{\textup{Stokes regions in }\MM\} @>{LD}>> 
\BB^{dist}/\{\pm 1\}^n\\
\end{CD}
\end{eqnarray*}
\end{theorem}

The last statement in part (a), $(\RR^f)_{hom}\subset G_\Z^{\BB}$,
and the last statement in part (c) that the map $LD$ is equivariant,
are not formulated in \cite[13.2]{He02}. But they follow easily
from the fact that any automorphism $\varphi\in \RR^{finite}$
lifts to an automorphism of the good representative of the 
universal unfolding.

We return to the simple singularities. Recall that the group $G_\Z$ is
\begin{eqnarray*}
G_\Z&=&\{\pm M^l\,|\, l\in\Z\}\times U\\
\textup{with}\quad U&\cong&
\left\{\begin{array}{ll}
\{\id\}&\textup{ in the cases }A_n,D_{2n+1},E_6,E_7,E_8,\\
S_2&\textup{ in the cases }D_{2n}\textup{ with }n\geq 3,\\
S_3&\textup{ in the case }D_4.
\end{array}\right.
\end{eqnarray*} 
In \cite[Ch. 8]{He11} as well as in \cite[5.1]{HR21} it is proved 
that the map $()_{hom}:\RR^{finite}\to G_\Z$ is an isomorphism
if $\mult(f)=2$ and that its image maps isomorphically to
$G_\Z/\{\pm \id\}$ if $\mult(f)=3$. 
As a conclusion from this discussion of symmetries of the simple
singularities, we obtain the following lemma
for the simple singularities. 

\begin{lemma}\label{t10.44}
The elements $(\varphi)_\MM$
for $\varphi\in\RR^{finite}$ form a group of deck transformations
for the covering 
$$\pr_n^{alg,S}:\MM^{alg}-(\KK_3^{alg}\cup \KK_2^{alg})
\to C_n^{S/\{\pm 1\}^n}$$
which is isomorphic to $G_\Z/\{\pm \id\}$. 
They are the restriction to $\MM^{alg}-(\KK_3^{alg}\cup\KK_2^{alg})$
of the elements of $\Aut_\MM$. Here $G_\Z=G_\Z^{\BB}$
and $\Aut_\MM\cong G_\Z/\{\pm\id\}$. 
\end{lemma}

On the other hand recall that the covering
$\pr_n^{e,S}:C_n^{\uuuu{e}/\{\pm 1\}^n}\to C_n^{S/\{\pm 1\}^n}$
is normal with group of deck transformations
\begin{eqnarray}\label{10.39}
\frac{(\Br_n)_{S/\{\pm 1\} ^n}}{(\Br_n)_{\uuuu{e}/\{\pm 1\}^n}}
\cong \frac{G_\Z^{\BB}}{Z((\{\pm 1\}^n)_S)}
\stackrel{\textup{here}}{=} \frac{G_\Z}{\{\pm \id\}}.
\end{eqnarray}
In the composition of coverings
\begin{eqnarray*}
\MM^{alg}-(\KK_3^{alg}\cup\KK_2^{alg})
\stackrel{\pr_n^{alg,e}}{\longrightarrow} C_n^{\uuuu{e}/\{\pm 1\}^n}
\stackrel{\pr_n^{e,S}}{\longrightarrow} C_n^{S/\{\pm 1\}^n}
\end{eqnarray*}
the second one is normal, and its group of deck transformations lifts to
a group of deck transformations of the composite covering
$\pr_n^{alg,S}=\pr_n^{e,S}\circ \pr_n^{alg,e}$. 
Because of the next lemma, the covering $\pr_n^{alg,e}$ is an isomorphism.

\begin{lemma}\label{t10.45}
The covering $\pr_n^{alg,S}$ is normal with group of deck transformations
the group $\Aut_\MM$.
\end{lemma}

The proof uses the following idea, which had been used first by
Jaworski \cite[Proposition 2]{Ja88} 
in the case of the Lyashko-Looijenga map for the simple
elliptic singularitiees. Consider two points $t^{(1)}$ and
$t^{(2)}\in\MM^{alg}$ with the same image 
$S/\{\pm 1\}^n=\pr_n^{alg,S}(t^{(1)})=\pr_n^{alg,S}(t^{(2)})
\in \SSS^{dist}/\{\pm 1\}^n$.
Consider in $C_n^{conf}$ a path from 
$\pr_n^{alg,c}(t^{(1)})=\pr_n^{alg,c}(t^{(2)})$ to a generic point of
$D_n^{conf}$. Then the lifts of this path to 
$\MM^{alg}-(\KK_3^{alg}\cup\KK_2^{alg})$ which start at $t^{(1)}$
or $t^{(2)}$ tend either both to $\KK_3^{alg}$
(if the relevant entry of the relevant matrix in the orbit of $S$ is $\pm 1$)
or both to $\KK_2^{alg}$
(if the relevant entry of the relevant matrix in the orbit of $S$ is $0$).
Therefore the covering $\pr_n^{alg,c}$ looks the same from
$t^{(1)}$ and $t^{(2)}$, so there is a deck transformation which maps
$t^{(1)}$ to $t^{(2)}$. It lifts to a deck transformation of the covering
$\pr_n^{alg,S}$. It extends to an automorphism of $\MM$ as F-manifold
with Euler field.

Now we turn to the simple elliptic singularities. The ideas above go through
and lead to Theorem \ref{t10.47}, 
which is analogoues to Theorem \ref{t10.42}.
A common point of the simple and simple elliptic singularities, which is
not shared by the other singularities, is that they have global algebraic
unfoldings. 

The 1-parameter families of simple elliptic singularities in the Legendre 
normal forms in Theorem \ref{t10.6} with parameter
$\lambda\in\C-\{0,1\}$ are quasihomogeneous polynomials 
except for $\lambda$ which has weight $\deg_{\bf w}(\lambda)=0$. 
They have global unfoldings which are also quasihomogeneous polynomials 
except for $\lambda$. \index{unfolding of a simple elliptic singularity}

\begin{eqnarray}
&&F^{alg}:\C^{m+1}\times \MM^{alg}\to\C \quad\textup{with}\quad
\MM^{alg}=\C^{n-1}\times(\C-\{0,1\}),\nonumber\\
&&F^{alg}(z_0,...,z_m,t_1,...,t_{n-1},\lambda)
=F^{alg}(z,t',\lambda)\\
&&=F^{alg}_{t',\lambda}(z) \nonumber 
=f_{\lambda}(z)+\sum_{j=1}^{n-1} t_jm_j \label{10.40}
\end{eqnarray}
with $f_{\lambda}=F^{alg}_{0,\lambda}$ and $m_1,...,m_{n-1}$ 
the monomials in table \eqref{10.41},
\begin{eqnarray}\label{10.41}
\begin{array}{cccccccccc}
\textup{name}  & m_1 & m_2 & m_3 & m_4 & m_5 
& m_6 & m_7 & m_8 & m_9 \\ \hline
\www E_6 &
1 & z_0 & z_1 & z_2 & z_0^2 & z_0z_1 & z_1z_2 & & \\
\www E_7 &
1 & z_0 & z_1 & z_0^2 & z_0z_1 & z_1^2 & z_0^2z_1 & z_0z_1^2 &  \\
\www E_8 & 
1 & z_0 & z_0^2 & z_1 & z_0^3 & z_0z_1 & z_0^2z_1 & z_1^2 & z_0z_1^2 
\end{array}
\end{eqnarray}
Let $\lambda:\H\to\C-\{0,1\},t_n\mapsto \lambda(t_n)$, be the standard
universal covering. For each of the three Legendre families of 
simple elliptic singularities, we will also 
consider the global family of functions
\begin{eqnarray}\label{10.42}
&&F^{mar}:\C^{m+1}\times \MM^{mar}\to\C, (z,t)
\mapsto F^{alg}(z,t',\lambda(t_m))\\
&&\textup{where}\quad \MM^{mar}=\C^{n-1}\times\H.\nonumber
\end{eqnarray}

The $\mu$-constant stratum within $\MM^{alg}$ is the 1-dimensional set
$\MM^{alg}_\mu=\{0\}\times(\C-\{0,1\})\subset\MM^{alg}$.
The global unfolding $F^{alg}$ is near any point 
$(0,\lambda)\in \MM^{alg}_\mu$ a universal unfolding of 
$f_\lambda=F_{0,\lambda}^{alg}$ \cite[Lemma 4.1]{HR21}. 
The weights $\deg_{\bf w}t_j$ for $j\in\{1,...,n-1\}$ are all positive.
The next table lists them and a rational number which will appear in
the degree of a Lyashko-Looijenga map. 
\begin{eqnarray*}
\begin{array}{l|l|l|l|l|l|l|l|l|l|l|l|l|l}
 & z_0 & z_1 & z_2 & t_1 & t_2 & t_3 & t_4 & t_5 & t_6 & t_7 & t_8 & t_9 
 & \frac{1}{2}\sum_{j=2}^{n-1}\frac{1}{\deg_{\bf w}t_j}\\
\www E_6 & \frac{1}{3} & \frac{1}{3} & \frac{1}{3} & 1 & \frac{2}{3} 
 & \frac{2}{3} & \frac{2}{3} & \frac{1}{3} & \frac{1}{3} & \frac{1}{3} & & 
 & \frac{27}{4} \\ 
\www E_7 & \frac{1}{4} & \frac{1}{4} & & 1 & \frac{3}{4} & \frac{3}{4} 
 & \frac{1}{2} & \frac{1}{2} & \frac{1}{2} & \frac{1}{4} & \frac{1}{4} &
 & \frac{25}{3} \\ 
\www E_8 & \frac{1}{6} & \frac{1}{3} & & 1 & \frac{5}{6} & \frac{2}{3} 
 & \frac{2}{3} & \frac{1}{2} & \frac{1}{2} & \frac{1}{3} & \frac{1}{3} 
 & \frac{1}{6} & \frac{101}{10}
\end{array}
\end{eqnarray*}
Because the weights are positive, $F^{alg}_{t',\lambda}$ for
each $(t',\lambda)$ satisfies that the sum of Milnor numbers of
its singular points is equal to $n$. It follows with little work that
$\MM^{alg}$ is an F-manifold with 
$$\textup{unit field }e=\paa_1,\quad
\textup{Euler field }E=\sum_{j=1}^{n-1}\deg_{\bf w}(t_j)t_j\paa_j$$
and polynomial multiplication \cite[Theorem 4.3]{HR21}
and also that $\MM^{mar}$ is an F-manifold with Euler field.

The Lyashko-Looijenga map for $\MM^{alg}$ was 
first studied by Jaworski \cite{Ja86}\cite{Ja88}. 
He found that the restriction of $\LL^{alg}:\MM^{alg}\to\C[x]_n$
to $\LL^{alg}:\MM^{alg}-(\KK_3^{alg}\cup\KK_2^{alg})\to C_n^{conf}$ 
is a finite covering. But contrary to Theorem \ref{t10.40} for the
simple singularities this is not at all trivial.
The map $\LL^{alg}$ is holomorphic, but not everywhere finite.
The 1-dimensional $\mu$-constant stratum $\MM^{alg}_\mu$
is mapped to $0\in\C[x]_n$. 

His methods did not allow Jaworski to calculate the degree
$\deg\LL^{alg}$. This was done by a huge effort in 
\cite[sections 5, 8, 9, 10]{HR21} by glueing into $\MM^{alg}$
suitable fibers over $\lambda\in\{0,1,\infty\}$ and extending
$\LL^{alg}$ to these fibers. 
A simpler way to calculate $\deg\LL^{alg}$ was recently
found by Takahashi and Zhang \cite{TZ23}. 

\begin{theorem}\label{t10.46}
(a) \cite{Ja86}\cite{Ja88} The restriction 
$$\LL^{alg}:\MM^{alg}-(\KK_3^{alg}\cup\KK_2^{alg})\to C_n^{conf}$$
of the Lyashko-Looijenga map $\LL^{alg}:\MM^{alg}\to \C[x]_n$ 
is a finite covering.

(b) \cite[Theorem 6.3]{HR21}\cite{TZ23} Its degree is
\begin{eqnarray}\label{10.43}
\deg\LL^{alg} =\frac{n!\cdot\frac{1}{2}\cdot 
\sum_{j=2}^{n-1}\frac{1}{\deg_{\bf w}t_j}}
{\prod_{j=2}^{n-1}\deg_{\bf w}t_j}.
\end{eqnarray}
\end{theorem}

The Stokes regions in $\MM^{alg}$ and in $\MM^{mar}$
are the components of the preimages in $\MM^{alg}$ 
and $\MM^{mar}$ of $\pr_n^{p,c}(F_n)\subset C_n^{conf}$.
The restriction of $\LL^{alg}$ respectively $\LL^{mar}$
to a Stokes region is an isomorphism to $\pr_n^{p,c}(F_n)$. 

We obtain a map \index{$LD$} 
$$LD:\{\textup{Stokes regions in }\MM^{mar}\}\to
\BB^{dist}/\{\pm 1\}^n$$
as for the simple singularities.
Though here we have to choose a reference point 
$(0,t_n^{ref})\in\MM^{mar}_\mu$. We choose it 
with $\lambda(t_n^{ref})=\frac{1}{2}$.
Then $H_\Z:=H_\Z(f_{1/2})$. 
Because the space $\MM^{mar}=\C^{n-1}\times\H$ is simply connected
there is a canonical isomorphism 
$$H_m(F_t^{-1}(\eta),\Z)\cong H_m(f_{1/2}^{-1}(\eta),\Z)
\cong H_\Z(f_{1/2})$$
for any $t\in\MM^{mar}$ and (depending on $t$) sufficiently large
$\eta>0$. We do not obtain a map $LD$ for the set of Stokes
regions in $\MM^{alg}$ because there we do not have such
canonical isomorphisms. The $\Z$-lattice bundle
$\bigcup_{\lambda\in\C-\{0,1\}}H_\Z(f_\lambda)$ has transversal
monodromy. 

Choose one distinguished basis $\uuuu{e}\in\BB^{dist}$ for
reference. The map $LD$ is surjective because the restriction
of $\LL^{mar}$ to a map
$$\LL^{mar}:\MM^{mar}-(\KK_3^{mar}\cup \KK_2^{mar})\to C_n^{conf}$$
is by Theorem \ref{t10.46} (a) a covering. Even more, it factorizes
in ways described by the following commutative diagram of coverings, 
\begin{eqnarray}\label{10.44}
\begin{xy}
\xymatrix{
\MM^{alg}-(\KK_3^{alg}\cup\KK_2^{alg}) \ar[dr]^{\textup{finite}}
\ar[rr]^{\textup{finite}} & & C_n^{conf} \\
 & C_n^{S/\{\pm 1\}^n} \ar[ur]^{\textup{finite}} &   \\
\MM^{mar}-(\KK_3^{mar}\cup\KK_2^{mar})  \ar[uu]^{\infty} 
\ar[ur]^{\infty} \ar[rr] & & 
C_n^{\uuuu{e}/\{\pm 1\}^n} \ar[ul]^{\infty} \ar[uu]^{\infty}
}
\end{xy}
\end{eqnarray}
This follows from the construction of the spaces and the compatibility
of the braid group actions. The main result in this section for the
simple elliptic singularities is the following.

\begin{theorem}\label{t10.47} \cite[Theorem 7.1]{HR21}
Consider a 1-parameter family of simple elliptic singularities. 
The map
\begin{eqnarray}\label{10.45}
LD:\{\textup{Stokes regions in }\MM^{mar}\}\to
\BB^{dist}/\{\pm 1\}^n
\end{eqnarray}
is a bijection. Equivalent: The covering
\begin{eqnarray}\label{10.46}
\pr_n^{mar,e}:\MM^{mar}-(\KK_3^{mar}\cup \KK_2^{mar})\to
C_n^{\uuuu{e}/\{\pm 1\}^n}
\end{eqnarray}
is an isomorphism.
\end{theorem}

The ideas of the proof are the same as in the proof of
Theorem \ref{t10.42} for the simple singularities in \cite{HR21}.
Recall that the group $G_\Z:=G_\Z(f_{1/2})$ is by Theorem 
\ref{t10.24} (c) the extension of the group
$G_\Z|_{\Rad I^{(0)}}\cong SL_2(\Z)$ by the finite group
$U_1\rtimes U_2$. The following results are proved in 
\cite{HR21} in a similar way as for the simple singularities.

\begin{lemma}\label{t10.48}
(a) $(\RR^{f_{1/2}})_{hom}=U_1\rtimes U_2$ if $\mult(f)=3$, and
$(\RR^{f_{1/2}})_{hom}=U_1\rtimes U_2\times\{\pm\id\}$ if $\mult(f)=2$.

(b) $G_\Z=G_\Z^{\BB}$.

(c) The covering $\pr_n^{mar,S}:\MM^{mar}-(\KK_3^{mar}\cup\KK_2^{mar})\to
C_n^{S/\{\pm 1\}^n}$ is a normal covering. The group of
deck transformations is isomorphic to 
$\Aut_{\MM^{mar}}:=\Aut(\MM^{mar},\circ,e,E)$ of $\MM^{mar}$
as F-manifold with Euler field and isomorphic to $G_\Z/\{\pm \id\}$.

(d) The covering $\pr_n^{e,S}:C_n^{\uuuu{e}/\{\pm 1\}^n}\to 
C_n^{S/\{\pm 1\}^n}$ is a normal covering with group of 
deck transformations isomorphic to $G_\Z/\{\pm \id\}$. 
\end{lemma}

The proof of part (c) uses Jaworski's idea which was described
after Lemma \ref{t10.45}. Now a matrix $S\in\SSS^{dist}$
can have the entries $0,\pm 1,\pm 2$. The lifts of a path
in $C_n^{conf}$ which tends to a generic point of $D_n^{conf}$ 
tend to $\KK_2^{mar}$ if the relevant entry is $0$,
to $\KK_3^{mar}$ if the relevant entry is $\pm 1$,
and to $\C^n\times\{0,1,\infty\}$ if the relevant entry
is $\pm 2$. 

Lemma \ref{t10.48} (c)+(d) and the diagram \eqref{10.44} of coverings 
show Theorem \ref{t10.47}.

\begin{remarks}\label{t10.49}
(i) The following table gives the degrees of the finite coverings in
the diagram \eqref{10.44} of coverings.
\begin{eqnarray}\label{10.47}
\begin{array}{llll}
 & \deg\pr_n^{alg,S}=6|U_1||U_2| & \deg\pr_n^{S,c}=|\SSS^{dist}/\{\pm 1\}^n|\\ 
\hline 
\www E_6 & 6\cdot 2\cdot 3\cdot 3^2 =324 & 3^7\cdot 5\cdot 7=76545\\
\www E_7 & 6\cdot 1\cdot 4\cdot 2^2 = 96 & 2^{13}\cdot 5^3\cdot 7 =7168000\\
\www E_8 & 6\cdot 1\cdot 6\cdot 1^2 = 36 & 2^7\cdot 3^8\cdot 7\cdot 101
=593744256
\end{array}
\end{eqnarray}
The surprising factor $101$ in the case of $\www{E}_8$ comes from
the term $\sum_{j=2}^{n-1}\frac{1}{\deg_{\bf w}t_j}$ in formula 
\eqref{10.43} for $\deg\LL^{alg}$. 

(ii) The number $\deg\LL^{alg}=\deg \pr_n^{alg,S}\cdot 
|\SSS^{dist}/\{\pm 1\}|$ was calculated in \cite{HR21} and
again in \cite{TZ23}. With the easily determined number
$\deg \pr_n^{alg,S}=6|U_1||U_2|$ it gives $|\SSS^{dist}/\{\pm 1\}^n|$.
Though this number $|\SSS^{dist}/\{\pm 1\}^n|$
was determined in the cases of $\www{E}_6$
and $\www{E}_7$ much earlier by Kluitmann, for $\www{E}_6$
in \cite{Kl83} and \cite[VII]{Kl87}, for $\www{E}_7$ in
\cite[VII]{Kl87}, by a long combinatorial and inductive procedure.
\end{remarks}

\begin{remarks}\label{t10.50}
In the cases of other singularities than the simple and simple 
elliptic singularities, there are a priori no natural global
unfoldings. A priori the base space of a universal unfolding
is a germ $(\MM,0)\cong(\C^n,0)$.

Consider a given $\mu$-homotopy class of singularities.
In \cite{He11} a global moduli space $\MM^{mar}_\mu$ 
for {\it right equivalence classes} of {\it marked singularities}
is constructed. Locally it is isomorphic to the $\mu$-constant
stratum in the base space of a universal unfolding of one singularity.
So one can consider it as a global $\mu$-constant stratum for all
singularities in one $\mu$-homotopy class. 

In order to obtain global moduli spaces of dimension $n$,
it would be good to thicken $\MM^{mar}_\mu$ to an $n$-dimensional
manifold which is near each point of $\MM^{mar}_\mu$ the base space
of a universal unfolding and then glue this manifold to the corresponding
manifold $C_n^{\uuuu{e}/\{\pm 1\}^n}$ where $\uuuu{e}\in\BB^{dist}$
is one chosen distinguished basis. 

By Theorem \ref{t10.42} and \ref{t10.47} this works for the
simple singularities and the simple elliptic singularities.
It gives the manifold $\MM^{alg}$ for the simple singularities and
the manifold $\MM^{mar}$ for the simple elliptic singularities. 

For other singularities this is an open project.  
\end{remarks}

\begin{appendix}

\chapter{Tools from hyperbolic geometry}\label{sa}
\setcounter{equation}{0}
\setcounter{figure}{0}

\renewcommand{\theequation}{\mbox{A.\arabic{equation}}}

\renewcommand{\thefigure}{\mbox{A.\arabic{figure}}}

The upper half plane 
\index{upper half plane}\index{hyperbolic plane}
$\H=\{z\in\C\,|\, \Im(z)>0\}$ together
with its natural metric (whose explicit form we will not need)
is one model of the hyperbolic plane. In the study of 
the monodromy groups $\Gamma^{(0)}$ and $\Gamma^{(1)}$ 
for the cases with $n=2$ or $n=3$, we will often encounter
subgroups of $\Isom(\H)$. The theorem of Poincar\'e-Maskit
\cite{Po82}\cite{Ma71} allows under some conditions to show
for such a group that it is discrete, to find a fundamental domain
and to find a presentation. Three special cases of this
theorem, which will be sufficient for us, are formulated
in Theorem \ref{ta.2}.

Before, we collect basic facts and set up some notations
in the following Remarks and Notations \ref{ta.1}.

Subgroups of $\Isom(\H)$ arise here in two ways.
Either they come from groups of real $2\times 2$ matrices.
This is covered by the Remarks and Notations \ref{ta.1} (v).
Or they come from the action of certain groups of
real $3\times 3$ matrices on $\R^3$ with an indefinite
metric. This will treated in Theorem \ref{ta.4}.

\begin{remarksandnotations}\label{ta.1}
(Some references for the following material are 
\cite{Fo51}\cite{Le66}\cite{Be83})

(i) Let $n\in\N$. Recall the notions of the 
free group $G^{free,n}$ with $n$ generators
and of the free Coxeter group $G^{fCox,n}$ 
with $n$ generators from Definition \ref{t3.1}.
%

(ii) $\widehat{\C}:=\C\cup\{\infty\}$,
$\widehat{\R}:=\R\cup\{\infty\}$. 
The hyperbolic lines \index{hyperbolic line}
in $\H$ are the parts in $\H$ of 
circles and euclidean lines which meet $\R$ orthogonally.
For any $z_1,z_2\in\H\cup\widehat{\R}$ with $z_1\neq z_2$,
denote by $A(z_1,z_2)$ 
the part between $z_1$ and $z_2$ of the unique hyperbolic
line whose closure in $\H\cup\widehat{\R}$ contains $z_1$
and $z_2$. Here $z_i\in A(z_1,z_2)$ if $z_i\in\H$, but
$z_i\notin A(z_1,z_2)$ if $z_i\in\widehat{\R}$.
Such sets are called {\it arcs}. \index{arc}

(iii) We simplify the definition of a polygon in \cite{Ma71}.
A {\it hyperbolic polygon} \index{hyperbolic polygon}
$P$ is a contractible open subset
$P\subset\H$ whose relative boundary in $\H$ consists
of finitely many arcs $A_1=A(z_{1,1},z_{1,2}),...,
A_m=A(z_{m,1},z_{m,2})$ (Maskit allows countably many arcs).
The arcs and the points are numbered
such that one runs through them in the order
$A_1,...,A_m$ and $z_{1,1},z_{1,2},z_{2,1},z_{2,2},...,
z_{m,1},z_{m,2}$ if one runs mathematically positive on the
euclidean boundary of $P$ in $\H\cup\widehat{\R}$.
For $A_i$ and $A_{i+1}$ (with $A_{m+1}:=A_1$) there are three
possibilities:
\begin{list}{}{}
\item[(a)] $z_{i,2}=z_{i+1,1}\in\H;$ then this point is called
a {\it vertex} of $P$;
\item[(b)] $z_{i,2}=z_{i+1,1}\in\widehat{\R}$;
\item[(c)] $z_{i,2}\in\widehat{\R},z_{i+1,1}\in\widehat{\R},
z_{i,2}\neq z_{i+1,1}$; then the part of $\widehat{\R}$
between $z_{i,2}$ and $z_{i+1,1}$ (moving from smaller to 
larger values) is in the euclidean boundary of $P$
between $A_i$ and $A_{i+1}$.
\end{list}
In the second and third case $A_i\cap A_{i+1}=\emptyset$.
A polygon has no vertices if and only if all arcs $A_1,...,A_m$
are hyperbolic lines, and if and only if all points
$z_{1,1},...,z_{m,2}\in\widehat{\R}$. 

(iv) Denote \index{$GL_2^{(\pm 1)}(\R)$}
\begin{eqnarray*}
Gl^{(-1)}_2(\R)&:=&\{A\in GL_2(\R)\,|\, \det A=-1\},\\
Gl^{(\pm 1)}_2(\R)&:=& \{A\in GL_2(\R)\,|\, \det A=\pm 1\}
=SL_2(\R)\cup GL^{(-1)}_2(\R),
\end{eqnarray*}
and analogously $Gl^{(-1)}_2(\Z)$, $Gl^{(\pm 1)}_2(\Z)$.
Recall 
$$A^t\begin{pmatrix}0&1\\-1&0\end{pmatrix}A=
\det A\cdot \begin{pmatrix}0&1\\-1&0\end{pmatrix}
\quad\textup{for }A\in Gl^{(\pm 1)}_2(\R).$$

(v) The following map $\mu:Gl^{(\pm 1)}_2(\R)\to\Isom(\H)$
is a surjective group homomorphism with kernel
$\ker\mu=\{\pm E_2\}$, 
\begin{eqnarray*}
\mu(A)=\Bigl(z\mapsto \frac{za+b}{cz+d}\Bigr)
&&\textup{ if }A=\begin{pmatrix}a&b\\c&d\end{pmatrix}
\in SL_2(\R),\\
\mu(A)=\Bigl(z\mapsto \frac{\oooo{z}a+b}{c\oooo{z}+d}\Bigr)
&&\textup{ if }A=\begin{pmatrix}a&b\\c&d\end{pmatrix}
\in GL^{(-1)}_2(\R).
\end{eqnarray*}
$\mu(A)$ for $A\in SL_2(\R)$ is orientation preserving and is
called a {\it M\"obius transformation}. 
\index{M\"obius transformation}
\index{$\mu:Gl^{(\pm 1)}_2(\R)\to\Isom(\H)$} 
$\mu(A)$ for $A\in GL^{(-1)}_2(\R)$ is orientation reversing.

If $A\in SL_2(\R)-\{\pm E_2\}$, there are three possibilities:
\begin{list}{}{}
\item[(a)] $|\tr(A)|<2$; then $A$ has a fixed point in $\H$
(and the complex conjugate number is a fixed point in $-\H$)
and is called {\it elliptic}.
\index{elliptic M\"obius transformation}
\item[(b)] $|\tr(A)|=2$; then $A$ has one fixed point in 
$\widehat{\R}$ and is called {\it parabolic}.
\index{parabolic M\"obius transformation}
\item[(c)] $|\tr(A)|>2$; then $A$ has two fixed points in
$\widehat{\R}$ and is called {\it hyperbolic}.
\index{hyperbolic M\"obius transformation}
\end{list}
If $A=\begin{pmatrix}a&b\\c&d\end{pmatrix}
\in GL^{(-1)}_2(\R)$ with $\tr(A)=0$ then $\mu(A)$ is
a reflection along the hyperbolic line 
\begin{eqnarray*}
\{z\in\H\,|\, z=\mu(A)(z)\}
=\{z\in\H\,|\, 0=cz\oooo{z}-2a\Ree(z)-b\}\\
\Bigl(= \{z\in\H\,|\, 0=(z-\frac{a}{c})(\oooo{z}-\frac{a}{c})
-\frac{1}{c^2}\}\quad\textup{if }c\neq 0\Bigr).
\end{eqnarray*}
\end{remarksandnotations}

The theorem of Poincar\'e-Maskit starts with a hyperbolic 
polygon $P$ whose relative boundary in $\H$ consists of
arcs $A_1,...,A_m$, with an involution 
$\sigma\in S_m$ and with elements $g_1,...,g_m\in \Isom(\H)$
with $g_i(A_i)=A_{\sigma(i)}$ and $g_{\sigma(i)}=g_i^{-1}$. 
Under some additional conditions, it states that 
the group $G:=\langle g_1,...,g_m\rangle\subset\Isom(\H)$
is discrete, that $P$ is a fundamental domain
(i.e. each orbit of $G$ in $\H$ meets the relative closure
of $P$ in $\H$, no orbit of $G$ in $\H$ meets $P$ in more
than one point), and it gives a complete set of relations
with respect to $g_1,...,g_m$ of $G$. 
Poincar\'e \cite{Po82} had the case when $\H/G$ is compact,
Maskit \cite{Ma71} generalized it greatly. In \cite{Ma71}
the relative boundary of $P$ in $\H$ may consist of
countably many arcs.  The following theorem singles out three
special cases, which are sufficient for us.
Remark \ref{ta.3} (ii) illustrates them with pictures.

\begin{theorem}\label{ta.2}
\index{Poincar\'e-Maskit theorem}
\cite{Ma71} Let $P\subset \H$ be a hyperbolic polygon whose
relative boundary in $\H$ consists of arcs $A_1,...,A_m$
with $A_i=A(z_{i,1}z_{i,2})$ where one runs through these
arcs and these points in the order $A_1,...,A_m$ and
$z_{1,1},z_{1,2},z_{2,1},z_{2,2},...,z_{m,1},z_{m,2}$
if one runs mathematically positive on the euclidean 
boundary of $P$ in $\H\cup\widehat{\R}$.

(a) Let $I\subset\{1,...,m\}$ be the set of indices such that
$z_{i,2}$ is a vertex, so $z_{i,2}=z_{i+1,1}\in\H$,
with $A_{m+1}:=A_1$ and $z_{m+1,i}:=z_{1,i}$ ($I$ may be empty).
Suppose that at a vertex $z_{i,2}$ the arcs $A_i$ and $A_{i+1}$
meet at an angle $\frac{\pi}{n_i}$ for some number 
$n_i\in\Z_{\geq 2}$. For $i\in \{1,...,m\}$ let 
$g_i\in \Isom(\H)$ be the reflection along the hyperbolic line
which contains $A_i$.

The group $G:=\langle g_1,...,g_m\rangle\subset\Isom(\H)$ is
discrete, $P$ is a fundamental domain, and the set of 
relations 
$$g_1^2=...=g_m^2=\id,\quad (g_ig_{i+1})^{n_i}=\id
\quad\textup{for }i\in I,$$
form a complete set of relations. Especially, if $I=\emptyset$ then
$G$ is a free Coxeter group with generators $g_1,...,g_m$. 

(b) Let $P\subset\H$ have no vertices.
Choose on each hyperbolic line $A_i$ a point $p_i$,
and let $g_i$ be the elliptic element with fixed point $p_i$
and rotation angle $\pi$.

The group $G:=\langle g_1,...,g_m\rangle\subset\Isom(\H)$ is
discrete, $P$ is a fundamental domain, and the set of 
relations $g_1^2=...=g_m^2=\id$
a complete set of relation, so $G$ is a free Coxeter group 
with generators $g_1,...,g_m$. 

(c) Let $P\subset\H$ have no vertices.
Suppose that $m$ is even. Suppose
$z_{2i-1,2}=z_{2i,1}$ for $i\in\{1,2,...,\frac{m}{2}\}$.
Let $g_i$ for $i\in\{1,2,...,\frac{m}{2}\}$ be the 
parabolic element with fixed point $z_{2i-1,2}$ which
maps $A_{2i-1}$ to $A_{2i}$.

The group $G:=\langle g_1,...,g_m\rangle\subset\Isom(\H)$ is
discrete, $P$ is a fundamental domain, and the group 
$G$ is a free group with generators $g_1,...,g_m$.
\end{theorem}

\begin{remarks}\label{ta.3}
(i) The Cayley transformation \index{Cayley transformation}
$\widehat{\C}\to\widehat{\C}$,
$\bigl(z\mapsto \frac{z-i}{z+i}\bigr)$ maps the upper half
plane $\H$ to the unit disk $\D=\{z\in\C\,|\, |z|<1\}$.
It leads to the unit disk model of the hyperbolic plane.
The hyperbolic lines in this model are the parts in $\D$
of circles and euclidean lines which intersect $\partial\D$ 
orthogonally. 

(ii) The following three pictures illustrate Theorem \ref{ta.2}
in the unit disk model instead of the upper half plane model.

\begin{figure}[H]
\includegraphics[width=1.0\textwidth]{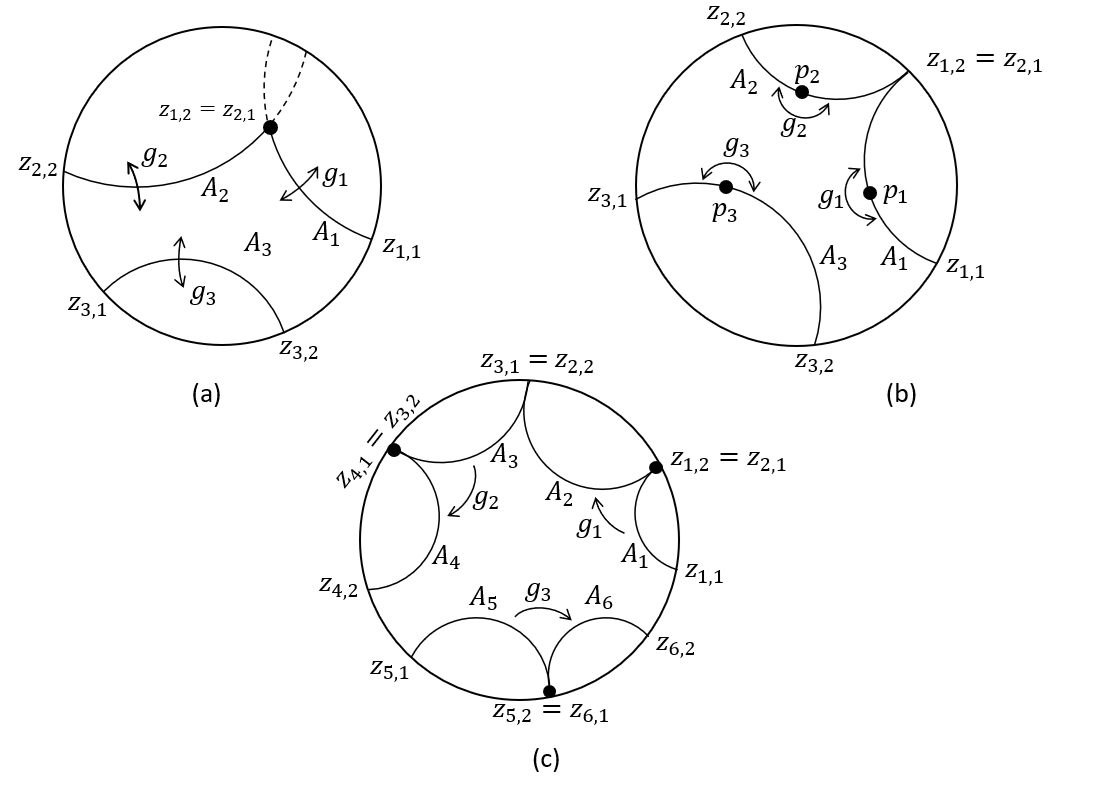}
\caption[Figure A.1]{Three pictures for Theorem \ref{ta.2}}
\label{Fig:A.1}
\end{figure}

\end{remarks}

The surjective group homomorphism 
$\mu:Gl^{(\pm 1)}_2(\R)\to\Isom(\H)$ in Remark \ref{ta.1} (v)
shows how to go from groups of real $2\times 2$ matrices
to subgroups of $\Isom(\H)$. The next theorem shows how to go
from groups of certain real $3\times 3$ matrices to subgroups
of $\Isom(\H)$. It is classical. But as we need some details,
we prefer to explain these details and not refer to
some literature.

\begin{theorem}\label{ta.4}
Let $(H_\R,I^{[0]})$ be a 3-dimensional real vector space with
a symmetric bilinear form $I^{[0]}$ with signature $(+--)$.

(a) (Elementary linear algebra) 
A vector $v\in H_\R-\{0\}$ is called {\sf positive} if 
$I^{[0]}(v,v)>0$, {\sf isotropic} if $I^{[0]}(v,v)=0$,
{\sf negative} if $I^{[0]}(v,v)<0$. 
The positive vectors \index{positive vector}
form a (double) cone $\KK\subset H_\R$, 
\index{cone}\index{$\KK\subset H_\R$}
the isotropic vectors \index{isotropic vector}
and the vector 0 form its boundary,
the negative vectors \index{negative vector} 
form its complement.
The orthogonal hyperplanes $(\R\cdot v)^\perp$ satisfy the 
following:
\begin{list}{}{}
\item[(i)]
$(\R\cdot v)^\perp\cap \KK\neq\emptyset$ if $v$ is negative.
\item[(ii)]
$(\R\cdot v)^\perp\cap\oooo{\KK}=\R\cdot v$ if $v$ is isotropic.
\item[(iii)]
$(\R\cdot v)^\perp\cap\oooo{\KK}=\{0\}$ if $v$ is positive.
\end{list}
$\KK/\R^*$ denotes the lines in $\KK$, i.e. the 1-dimensional
subspaces.

(b) (Basic properties of $\Aut(H_\R,I^{[0]})$)
Let $\sigma:\Aut(H_\R,I^{[0]})\to\{\pm 1\}$ be the spinor
norm map (see Remark \ref{t6.3} (iii)). The group 
$\Aut(H_\R,I^{[0]})$ is a real 3-dimensional Lie group
with four components. The components are the fibers 
of the group homomorphism
\begin{eqnarray*}
(\det,\sigma):\Aut(H_\R,I^{[0]})\to \{\pm 1\}\times \{\pm 1\}.
\end{eqnarray*}
$-\id\in \Aut(H_\R,I^{[0]})$ has value $(\det,\sigma)(-\id)=(-1,1)$.
If $v$ is positive then $(\det,\sigma)(s^{(0)}_v)=(-1,1)$.
If $v$ is negative then $(\det,\sigma)(s^{(0)}_v)=(-1,-1)$. 
An isometry $g\in \Aut(H_\R,I^{[0]})$ maps each of the two
components of the cone $\KK$ to itself if and only if $g$ is
in the two components which together form the kernel of $\det\cdot\sigma$,
so if $\det(g)\sigma(g)=1$.

(c) Choose a basis $\uuuu{f}=(f_1,f_2,f_3)$ of $H_\R$ with
$$I^{[0]}(\uuuu{f}^t,\uuuu{f})=
\begin{pmatrix}0&0&1\\0&-2&0\\1&0&0\end{pmatrix}.$$

(i) The map \index{$\Theta: Gl^{(\pm 1)}_2(\R)\to \Aut(H_\R)$} 
\begin{eqnarray*}
\Theta: Gl^{(\pm 1)}_2(\R)\to \Aut(H_\R),\quad
\begin{pmatrix}a&b\\c&d\end{pmatrix}\mapsto
\Bigl(\uuuu{f}\mapsto \uuuu{f}
\begin{pmatrix}a^2&2ab&b^2\\ac&ad+bc&bd\\c^2&2cd&d^2
\end{pmatrix}\Bigr),
\end{eqnarray*}
is a group homomorphism with kernel
$\ker\Theta=\{\pm E_2\}$ and image
$\ker(\det\cdot\sigma:\Aut(H_\R,I^{[0]})\to\{\pm 1\})$.

(ii) The map \index{$\vartheta:\H\to (H_\R-\{0\})/\R^*$} 
\begin{eqnarray*}
\vartheta:\H\to (H_\R-\{0\})/\R^*,\quad 
z\mapsto\R^*(z\oooo{z}f_1+\Ree(z)f_2+f_3),
\end{eqnarray*}
is a bijection $\vartheta:\H\to \KK/\R^*$. 

(iii) For $A\in Gl^{(\pm 1)}_2(\R)$, the automorphism
$\vartheta\circ\mu(A)\circ\vartheta^{-1}:\KK/\R^*\to \KK/\R^*$
coincides with the action of $\Theta(A)$ on $\KK/\R^*$.

(iv) The natural maps
\begin{eqnarray*}
\begin{array}{ccccc}
\Aut(H_\R,I^{[0]})/\{\pm \id\} & \longleftarrow & 
\ker(\det\cdot\sigma) & \longrightarrow & \Isom(\H)\\
\{\pm B\} & \longleftarrow & B\ , \ \Theta(A) & \longmapsto &
\mu(A) \end{array}
\end{eqnarray*}
are group isomorphisms.

(v) For each hyperbolic line $l$ there is a negative vector
$v\in H_\R$ with $\vartheta(l)=((\R\cdot v)^\perp\cap\KK)/\R^*$.

(vi) Let $\sigma_l\in \Isom(\H)$ be the reflection along 
a hyperbolic line $l$, and let $v$ be as in (v). The action
of $s^{(0)}_v$ on $\KK/\R^*$ coincides with 
$\vartheta\circ\sigma_l\circ\vartheta^{-1}$.

(vii) Let $\delta_p\in\Isom(\H)$ be the elliptic element
with fixed point $p\in\H$ and order 2 (so rotation angle $\pi$).
Let $v\in H_\R$ be a positive vector with $\R\cdot v=\vartheta(p)$.
The action of $s^{(0)}_v$ on $\KK/\R^*$ coincides with
$\vartheta\circ \delta_p\circ\vartheta^{-1}$.
\end{theorem}

{\bf Proof:}
(a) and (b) are elementary and classical, their proofs are skipped.

(c) (i) Start with a real 2-dimensional vector space $V_\R$
with basis $\uuuu{e}=(e_1,e_2)$ and a skew-symmetric bilinear
form $I^{[1]}$ on $V_\R$ with matrix 
$I^{[1]}(\uuuu{e}^t,\uuuu{e})=\begin{pmatrix}0&1\\-1&0\end{pmatrix}$.

The tensor product $V_\R\otimes V_\R$ comes equipped with an induced
symmetric bilinear form $\www{I}^{(0)}$ via
(here $I$ and $J$ are finite index sets and $a_i,b_i,c_j,d_j\in V_\R$)
\begin{eqnarray*}
\www{I}^{(0)}(\sum_{i\in I}a_i\otimes b_i,\sum_{j\in J}c_j\otimes d_j)
=\sum_{i\in I}\sum_{j\in J}I^{[1]}(a_i,c_j)I^{[1]}(b_i,d_j).
\end{eqnarray*}
An element $g\in\Aut(H_\R)$ with 
$g\uuuu{e}=\uuuu{e}A$ and $A\in Gl^{(\pm 1)}_2(\R)$ respects
$I^{[1]}$ in the following weak sense:
\begin{eqnarray*}
I^{[1]}(g(v_1),g(v_2))&=&\det A\cdot I^{[1]}(v_1,v_2),\\
\textup{because }A^t\begin{pmatrix}0&1\\-1&0\end{pmatrix}A
&=&\det A\cdot \begin{pmatrix}0&1\\-1&0\end{pmatrix}.
\end{eqnarray*}
It induces an element 
$\www{\Theta}(g)\in\Aut(V_\R\otimes V_\R,\www{I}^{(0)})$ via
\begin{eqnarray*}
\www{\Theta}(g)(\sum_{i\in I}a_i\otimes b_i)=
\sum_{i\in I}g(a_i)\otimes g(b_i).
\end{eqnarray*}

The symmetric part $\www{H}_\R\subset V_\R\otimes V_\R$ of the
tensor product has the basis
$\www{\uuuu{f}}=(\www{f}_1,\www{f}_2,\www{f}_3)
=(e_1\otimes e_1,e_1\otimes e_2+e_2\otimes e_1,e_2\otimes e_2)$.
One sees
\begin{eqnarray*}
\www{I}^{(0)}(\www{\uuuu{f}}^t,\www{\uuuu{f}})
=\begin{pmatrix}0&0&1\\0&-2&0\\1&0&0\end{pmatrix}.
\end{eqnarray*}
From now on we identify 
$(\www{H}_\R, \www{I}^{(0)}|_{\www{H}_\R},\www{\uuuu{f}})$
with $(H_\R,I^{[0]},\uuuu{f})$. 

For an element $g\in \Aut(V_\R)$ with 
$g\uuuu{e}=\uuuu{e}A$ with 
$A=\begin{pmatrix}a&b\\c&d\end{pmatrix}\in Gl^{(\pm 1)}_2(\R)$,
the automorphism $\www{\Theta}(g)$ on $V_\R\otimes V_\R$
restricts to an automorphism of the symmetric part $H_\R$ with
matrix 
\begin{eqnarray*}
\www{\Theta}(g)\uuuu{f}=\uuuu{f}
\begin{pmatrix}a^2&2ab&b^2\\ac&ad+bc&bd\\c^2&2cd&d^2
\end{pmatrix}.
\end{eqnarray*}
This fits to $\Theta$. 
It shows especially 
$\Theta(A)\in\Aut(H_\R,I^{[0]})$. 

The kernel of $\Theta$ is $\{\pm E_2\}$.
The Lie groups $Gl^{(\pm 1)}_2(\R)$ and 
$\Aut(H_\R,I^{[0]})$ are real 3-dimensional.
The Lie group $Gl^{(\pm 1)}_2(\R)$ has two components.
$(\det,\sigma)(\Theta(\begin{pmatrix}1&0\\0&-1\end{pmatrix}))
=(\det,\sigma)(s^{(0)}_{f_2})=(-1,-1)$. 
Therefore the image of $\Theta$ consists of the two components
of $\Aut(H_\R,I^{[0]})$ which together form the kernel of
$\det\cdot \sigma$. This finishes the proof of part (i).

(ii) Define $\www{\vartheta}(z):=z\oooo{z}f_1+\Ree(z)f_2+f_3$.
It is a positive vector because
\begin{eqnarray*}
I^{[0]}(\www{\vartheta}(z),\www{\vartheta}(z))
=z\oooo{z}\cdot 1-2(\Ree(z))^2+1\cdot z\oooo{z}=2(\Im(z))^2>0.
\end{eqnarray*}
It is easy to see that $\vartheta$ is a bijection from $\H$
to $\KK/\R^*$. 

(iii) In fact, $\www{\vartheta}(z)$ is the symmetric part of
$$(ze_1+e_2)\otimes (\oooo{z}e_1+e_2)
=(\uuuu{e}\begin{pmatrix}z\\1\end{pmatrix})\otimes
(\uuuu{e}\begin{pmatrix}\oooo{z}\\1\end{pmatrix})
\in V_\C\otimes V_\C.$$
For $A=\begin{pmatrix}a&b\\c&d\end{pmatrix}
\in Gl^{(\pm 1)}_2(\R)$ 
\begin{eqnarray*}
&&\www{\vartheta}(\mu(A)(z))\\
&=& \Bigl(\textup{symmetric part of }
(\uuuu{e}\begin{pmatrix}\mu(A)(z)\\1\end{pmatrix})\otimes
(\uuuu{e}\begin{pmatrix}\mu(A)(\oooo{z})\\1\end{pmatrix}\Bigr)\\
&=& |cz+d|^{-2}\Bigl(\textup{symmetric part of }
(\uuuu{e}A\begin{pmatrix}z\\1\end{pmatrix})\otimes
(\uuuu{e}A\begin{pmatrix}\oooo{z}\\1\end{pmatrix})\Bigr)\\
&=& |cz+d|^{-2}\uuuu{f}
\begin{pmatrix}a^2&2ab&b^2\\ac&ad+bc&bd\\c^2&2cd&d^2
\end{pmatrix}
\begin{pmatrix}z\oooo{z} \\ \Ree(z) \\ 1\end{pmatrix}\\
&=& |cz+d|^{-2}\Theta(A)(\uuuu{f})
\begin{pmatrix}z\oooo{z}\\ \Ree(z)\\1\end{pmatrix}\\
&=& |cz+d|^{-2}\Theta(A)(\www{\vartheta}(z)).
\end{eqnarray*}
This shows part (iii). Part (iv) follows from part (iii).

(v) A hyperbolic line $l$ is the fixed point set of a
reflection $\mu(A)$ for a matrix 
$A=\begin{pmatrix}a&b\\c&-a\end{pmatrix}\in Gl^{(-1)}_2(\R)$, so
\begin{eqnarray*}
l&=& \{z\in\H\,|\, z=\mu(A)(z)\}
=\{z\in \HH\,|\, 0=cz\oooo{z}-2a\Ree(z)-b\}.
\end{eqnarray*}
Observe
\begin{eqnarray*}
cz\oooo{z}-2a\Ree(z)-b=I^{[0]}(\uuuu{f}
\begin{pmatrix}-b\\a\\c\end{pmatrix},
\uuuu{f}\begin{pmatrix}z\oooo{z}\\ \Ree(z)\\1\end{pmatrix}).
\end{eqnarray*}
Therefore
\begin{eqnarray*}
\vartheta(l)&=& 
\Bigl(\R(-bf_1+af_2+cf_3)^\perp\cap\KK\Bigr)/\R^*.
\end{eqnarray*}

(vi) and (vii) are clear now.\hfill$\Box$

\begin{remarks}\label{ta.5}
(i) In the model $\KK/\R^*$ of the hyperbolic plane,
the isometries of $\H$ are transformed to linear
isometries of $(H_\R,I^{[0]})$.
The hyperbolic lines in $\H$ are transformed to 
linear hyperplanes in $H_\R$ (modulo $\R^*$) which
intersect $\KK$. 

(ii) If one chooses an affine hyperplane in $H_\R$ 
which intersects one component of $\KK$ in a disk,
this disk gives a new disk model of the hyperbolic plane,
which is not conformal to $\H$ and $\D$
(angles are not preserved), but where the hyperbolic
lines in $\H$ are transformed to the segments in the
new disk of euclidean lines in the affine hyperplane
which intersect the disk. The following picture
sketches three hyperbolic lines in the unit disk $\D$ 
and in the new disk which is part of the cone $\KK$.
\begin{figure}[H]
\includegraphics[width=0.7\textwidth]{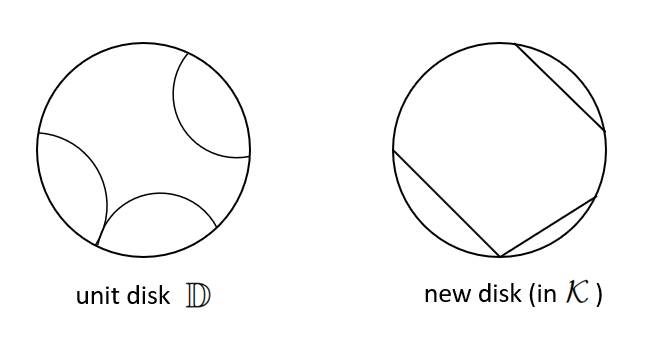}
\caption[Figure A.2]{Two disk models of the hyperbolic plane}
\label{Fig:A.2}
\end{figure}

\end{remarks}

\chapter{The first congruence subgroups}\label{sb}
\setcounter{equation}{0}
\setcounter{figure}{0}

\renewcommand{\theequation}{\mbox{B.\arabic{equation}}}

\renewcommand{\thefigure}{\mbox{B.\arabic{figure}}}

For the 3-dimensional covering spaces in section
\ref{s9.2}, we need some classical facts on the quotient map
$\H\to\H/\Gamma(n)$ for $n\in\{2,3\}$. In the case $n=2$
we need additionally a lift of the multivalued logarithm
on $\C-\{0,1\}\cong\H/\Gamma(2)$ to a univalued holomorphic
function on $\H$. 

For $n\in\Z_{\geq 2}$ the congruence subgroup 
\index{congruence subgroup}\index{$\Gamma(n)$}
$\Gamma(n)$ is the group
\begin{eqnarray*}
\Gamma(n):=\{A\in SL_2(\Z)\,|\, A\equiv E_2\mmod n\}.
\end{eqnarray*}
For $n=2$ we have $-E_2\in\Gamma(2)$, and $\Gamma(2)/\{\pm E_2\}$
embeds into $PSL_2(\Z)$. For $n\geq 3$ we have 
$-E_2\notin \Gamma(n)$,
and $\Gamma(n)$ embeds into $PSL_2(\Z)$. In any case
we write $\oooo{\Gamma(n)}\subset PSL_2(\Z)$ for the image
in $PSL_2(\Z)$. The following theorem collects classical facts
(see e.g. \cite[5.4 and 5.7.3]{La09}).

\begin{theorem}\label{tb.1}
(a) Fix $n\in\Z_{\geq 2}$. The group $\oooo{\Gamma(n)}$ is a 
normal subgroup of finite index in $PSL_2(\Z)$. It does not
contain elliptic elements. It acts properly discontinuously
and freely on $\H$. The quotient $\H/\Gamma(n)$ is the 
complement of finitely many points in a smooth compact complex
curve. The finitely many missing points correspond to the
$\oooo{\Gamma(n)}$ orbits in $\Q\cup\{\infty\}$,
which is the set of fixed points of the parabolic elements
in $\oooo{\Gamma(n)}$. The number of the missing points is
$[PSL_2(\Z):\oooo{\Gamma(n)}]\cdot n^{-1}$. 

(b) Fix $n\in\{2,3,4,5\}$. Then $\H/\Gamma\cong \P^1\C-
\{\textup{finitely many points}\}$. 
\begin{eqnarray*}
\begin{array}{cccc}
n & PSL_2(\Z)/\oooo{\Gamma(n)} & [PSL_2(\Z):\oooo{\Gamma(n)}] 
& |\{\textup{missing points}\}| \\ \hline 
2 & \cong S_3 & 6 & 3 \\
3 & \cong A_4 & 12 & 4 \\
4 & \cong S_4 & 24 & 6 \\
5 & \cong A_5 & 60 & 12 
\end{array}
\end{eqnarray*}
There are a natural homeomorphism $h:\P^1\C\to S^2\subset\R^3$
and embeddings of the vertices of a 
\index{platonic solid}
tetrahedron, an octahedron
and an icosahedron into $S^2$ such that the image
$h(\{\textup{missing points}\})$ is the set of vertices
of the tetrahedron in the case $n=3$, the set of vertices
of the octahedron in the case $n=4$ and the set of vertices
of the icosahedron in the case $n=5$. The group
$PSL_2(\Z)/\oooo{\Gamma(n)}$ for $n\in\{3,4,5\}$
acts on $S^2$ as group of rotations which fix the corresponding
platonic solid. 
For $n=2$ $\H/\oooo{\Gamma(2)}\cong\C-\{0,1\}$
(for the group $PSL_2(\Z)/\oooo{\Gamma(2)}$ see part (d)).

For $n\in\{2,3,4,5\}$ the composition 
\begin{eqnarray*}
T:\H\to \H/\oooo{\Gamma(n)}\stackrel{\cong}{\longrightarrow}
\P^1\C-\{\textup{missing points}\}\hookrightarrow \C 
\end{eqnarray*}
is a Schwarzian triangle function.
\index{$T:\H\to \P^1\C-\{\textup{points}\}$}
\index{Schwarzian triangle function}
It maps the degenerate hyperbolic triangle 
\index{hyperbolic triangle}
with vertices
$0,1,\infty$ and angles $0,0,0,$ to a spherical triangle
with angles $\frac{2\pi}{n},\frac{2\pi}{n},\frac{2\pi}{n}$.
It maps the degenerate hyperbolic triangle with vertices 
$i,e^{2\pi i/6},\infty$ with angles $\frac{\pi}{2},
\frac{\pi}{3},0$ to a spherical triangle 
\index{spherical triangle}
with angles
$\frac{\pi}{2},\frac{\pi}{3},\frac{\pi}{n}$. 

(c) In the case $n=3$, there is an isomorphism 
$\P^1\C\cong\C\cup\{\infty\}$ such that the 24 spherical
triangles with angles 
$\frac{\pi}{2},\frac{\pi}{3},\frac{\pi}{3}$ in $\P^1\C$
have the following images in 
$\C-\{0,1\}\subset \C\cup\{\infty\}$.
Figure \ref{Fig:b.1} shows a subdivision of each of the 
two visible  faces of a tetrahedron into 6 triangles. 
Projecting the tetrahedron to $S^2$ gives 24 spherical triangles.
Their boundaries consist of altogether six great circles in
$S^2$.
Figure \ref{Fig:b.2} shows the images of the six great circles in
$\C$. Three are circles, three are lines (where the point 
$\infty$ is missing).
Figure \ref{Fig:b.3} restricts to showing the preimages of the
edges of the tetrahedron.
The preimages in $\C\cup\{\infty\}$ of the vertices 
of the tetrahedron are the points 
$\{\infty,1,e^{2\pi i/3},e^{-2\pi i/3}\}$.
Therefore 
\begin{eqnarray*}
\H/\oooo{\Gamma(3)}\cong \C-\{z\,|\, z^3=1\}.
\end{eqnarray*}

\begin{figure}
\parbox{4cm}{\includegraphics[width=0.3\textwidth]{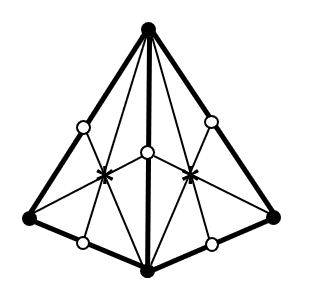}}
\parbox{6cm}{
$\bullet$ \hspace*{0.3cm} 4 vertices, \\
$\circ$ \hspace*{0.3cm} 5 (of 6) center points of edges,\\
$*$ \hspace*{0.3cm} 2 (of 4) center points of faces.}
\caption[Figure B.1]{A tetrahedron and a subdivision of each 
face into six triangles, only two faces shown}
\label{Fig:b.1}
\end{figure}

\begin{figure}
\includegraphics[width=0.8\textwidth]{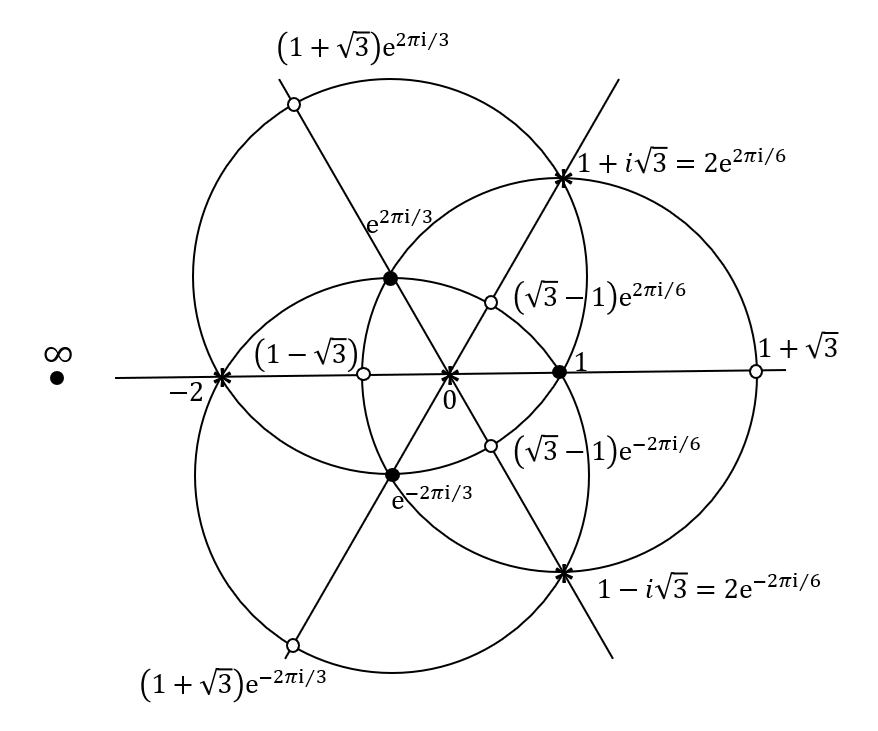}
\caption[Figure B.2]{Images in $\C\cup\{\infty\}$ of 24 spherical triangles}
\label{Fig:b.2}
\end{figure}

\begin{figure}
\includegraphics[width=0.4\textwidth]{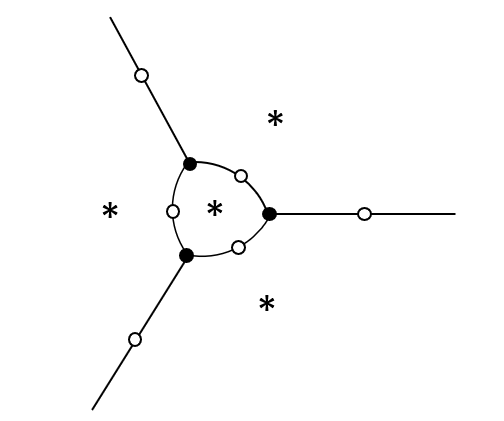}
\caption[Figure B.3]{The edges of the tetrahedron in Figure B.2}
\label{Fig:b.3}
\end{figure}

(d) The case $n=2$. The group of M\"obius transformations of 
$\C\cup\{\infty\}$ which fix the set $\{0,1,\infty\}$
is called {\sf anharmonic group} 
\index{anharmonic group} 
$G^{anh}$. The first four
entries of the following table fix an isomorphism
$g_\bullet:S_3\to G^{anh}$ with 
$g_{\sigma_1}g_{\sigma_2}=g_{\sigma_1\sigma_2}$. 

\begin{eqnarray*}
\begin{array}{llllllll}
\sigma & g_\sigma & g_\sigma(z) & 
(g_\sigma(1),g_\sigma(0),g_\sigma(\infty)) & \mu_\sigma & 
\mu_\sigma(\tau) & \mu_\sigma & \textup{fp} \\
\hline
\id & g_{\id} & z & (1,0,\infty) & 
\mu_{\id} & \tau & \textup{identity} & \H\\ 
(1\, 2) & g_{(12)} & 1-z & (0,1,\infty) & 
\mu_{(12)} & \tau-1 & \textup{parabolic} & \infty \\
(1\, 3) & g_{(13)} & \frac{z}{z-1} & (\infty,0,1) & 
\mu_{(13)} & \frac{\tau}{\tau+1} & 
\textup{parabolic} & 0 \\
(2\, 3) & g_{(23)} & \frac{1}{z} & (1,\infty,0) & 
\mu_{(23)} & \frac{-1}{\tau} & 
\textup{rotation} & i \\
(1\, 2\, 3) & g_{(123)} & \frac{z-1}{z} & (0,\infty,1) & 
\mu_{(123)} & \frac{\tau-1}{\tau} & 
\textup{rotation} & e^{2\pi i/6} \\
(1\, 3\, 2) & g_{(132)} & \frac{1}{1-z} & (\infty,1,0) & 
\mu_{(132)} & \frac{-1}{\tau-1} &
\textup{rotation} & e^{2\pi i/6}
\end{array}
\end{eqnarray*}

The Schwarzian triangle function
$T:\H\to \H/\oooo{\Gamma(2)}\cong 
\P^1\C-\{\textup{missing points}\}$ 
\index{Schwarzian triangle function}
\index{$T:\H\to \P^1\C-\{\textup{points}\}$}
can here be chosen as 
\begin{eqnarray*}
T=g_{(12)}\circ 
(\textup{the classical }\lambda\textup{-function}):
\H\to\C-\{0,1\}
\end{eqnarray*}
(also the $\lambda$-function itself works, but we prefer
$T$ as it maps $0,1,\infty$ to $0,1,\infty$). 
The following table gives the values under $T$ of some points
or sets.

\begin{eqnarray*}
\begin{array}{l|l||l|l}
\textup{point} & T(\textup{point})
& \textup{set} & T(\textup{set}) \\ \hline
0 & 0 &  i\R_{>0}=A(0,\infty) & \R_{<0}\\
1 & 1 & 1+i\R_{>0}=A(\infty,1) & \R_{>1}\\
\infty & \infty & A(1,0) & (0,1)\\
i & -1 & \frac{1}{2}+i\R_{>0}=A(\frac{1}{2},\infty) & 
\frac{1}{2}+i\R\\
\frac{1}{2}(1+i) & \frac{1}{2} & A(1,-1) & S^1-\{1\}\\
e^{2\pi i /6} & e^{2\pi i /6} & A(2,0) & 1+(S^1-\{-1\}) 
\end{array}
\end{eqnarray*}

The hyperbolic polygon with boundary 
$A(1,0)\cup A(0,-1)\cup A(-1,\infty)\cup A(\infty,1)$
is a fundamental domain for the action of $\oooo{\Gamma(2)}$
on $\H$. It consists of 12 degenerate hyperbolic triangles
with angles $\frac{\pi}{2},\frac{\pi}{3},0$. 
They are shown in the left picture in Figure B.4.
The right picture in Figure \ref{Fig:b.4} shows their images in
$\C-\{0,1\}$. The images are the connected components of the
complement of the two circles and the lines $\R$ and $i\R$. 

\begin{figure}
\includegraphics[width=1.0\textwidth]{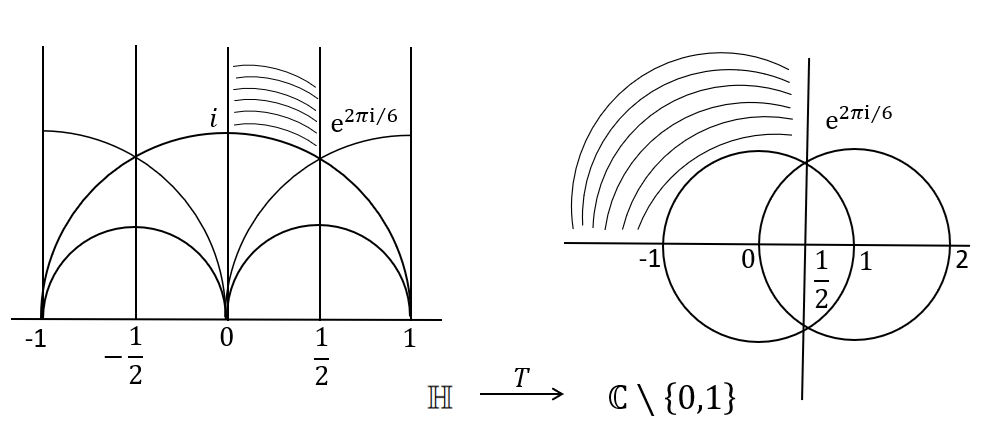}
\caption[Figure B.4]{$T$ maps 12 hyperbolic triangles 
to 12 spherical triangles}
\label{Fig:b.4}
\end{figure}

The six M\"obius transformations in the fifth and sixth 
column of the table above are representatives of the classes
in the quotient group $\mu(SL_2(\Z))/\mu(\Gamma(2))\cong S_3$,
with 
\begin{eqnarray*}
T(\mu(A)\mu_\sigma(\tau))=T(\mu_\sigma\mu(A)(\tau))
=g_\sigma(T(\tau))\quad\textup{for }\sigma\in S_3,
A\in\Gamma(2).
\end{eqnarray*}
The seventh and eighth column give properties of them
({\sf fixed point} is abbreviated {\sf fp}). 
\end{theorem}

{\bf Proof:}
See e.g. \cite[5.4 and 5.7.3]{La09}. \hfill$\Box$

\begin{definition}\label{tb.2}
The logarithm $\ln:\C-\{0\}\dashrightarrow\C$
is multivalued. As $T:\H\to\C-\{0,1\}$ is a universal
covering of $\C-\{0,1\}$, the logarithm has univalued
lifts $\H\to\C$. We denote by $\kappa:\H\to\C$ the lift
\index{lift of the logarithm}\index{$\kappa:\H\to\C$}
with values in $\R_{>0}$ on $1+i\R_{>0}$
(this set maps under $T$ to $\R_{>1}$).
\end{definition}

The geometry of this lift is described in part (c) of Theorem 
\ref{tb.3}. The parts (a) and (b) help to make part (c) more
transparent. They contain basic well known facts.

\begin{theorem}\label{tb.3}
(a) For each cusp $c\in \Q\cup\{\infty\}$ of $PSL_2(\Z)$,
the stabilizer $(\mu(SL_2(\Z)))_c$ of $c$ in $\mu(SL_2(\Z))$
is infinite cyclic with a parabolic generator.
For example, for $c\in\{\infty,0,1\}$ generators are as follows,
\begin{eqnarray*}
\begin{array}{l|l}
c & \textup{generator of }(\mu(SL_2(\Z))_c \\ \hline
\infty &\mu_{(12)}=\mu(\begin{pmatrix}1&-1\\0&1\end{pmatrix})\\
0      &\mu_{(13)}=\mu(\begin{pmatrix}1&0&\\1&1\end{pmatrix})\\
1      &\mu(\begin{pmatrix}2&-1\\1&0\end{pmatrix})
\end{array}
\end{eqnarray*}
The cusps form a single $\mu(SL_2(\Z))$ orbit.
If $c=\mu(B)(\infty)$ then 
$(\mu(SL_2(\Z)))_c = 
\langle\mu(B)\mu_{(12)}\mu(B)^{-1}\rangle$.

(b) For each cusp $c\in \Q\cup\{\infty\}$ of $\oooo{\Gamma(2)}$,
the stabilizer in $\mu(\Gamma(2))$ is infinite cyclic with
a parabolic generator. The generator is the square of a  
generator of $(\mu(SL_2(\Z)))_c$. For example, for 
$c\in\{\infty,0,1\}$ generators are as follows,
\begin{eqnarray*}
\begin{array}{l|l}
c & \textup{generator of }(\mu(\Gamma(2))_c \\ \hline
\infty &\mu_{(12)}^2=\mu(\begin{pmatrix}1&-2\\0&1\end{pmatrix})\\
0      &\mu_{(13)}^2=\mu(\begin{pmatrix}1&0&\\2&1\end{pmatrix})\\
1      &\mu(\begin{pmatrix}3&-2\\2&-1\end{pmatrix})
=\mu_{(12)}^{-2}\mu_{(13)}^{-2}
\end{array}
\end{eqnarray*}
The cusps form three $\Gamma(2)$ orbits, the orbits of
$0$, $1$ and $\infty$. If $c\in\Q\cup\{\infty\}$ is in the
orbit of $c_0\in\{0,1,\infty\}$ with 
$c=\mu(B)(c_0)$ for some $B\in\Gamma(2)$ then  
$(\mu(\Gamma(2)))_c = 
\mu(B)(\mu(\Gamma(2)))_{c_0} \mu(B)^{-1}$.

(c) The lift $\kappa:\H\to\C$ of the logarithm 
$\ln:\C-\{0\}\dashrightarrow\C$ satisfies the following.
\index{$\kappa:\H\to\C$}

(i) The very fact that it is a lift with $T$ says
\begin{eqnarray*}
T(\tau)=\exp(\kappa(\tau))\quad\textup{for }\tau\in\H.
\end{eqnarray*}
\begin{eqnarray*}
\begin{xy}
\xymatrix{\H \ar[drr]_\kappa \ar[rr]^T && \C-\{0,1\} 
\ar[d]_{\ln} \ar[dr]_\id & \\
 & & \C \ar[r]_\exp & \C-\{0\}}
\end{xy}
\end{eqnarray*}
The following table shows the values under $\kappa$ of some
sets.
\begin{eqnarray*}
\begin{array}{l|l|l|l}
\textup{set} & 1+i\R_{>0}=A(\infty,1) & A(1,0) & 
i\R_{>0} = A(0,\infty) \\ \hline 
\kappa(\textup{set}) & \R_{>0} & \R_{<0} & \pi i +\R
\end{array}
\end{eqnarray*}

(ii) Let $\nu$ be the (because of 
$\mu(\Gamma(2))=\langle\mu_{(12)}^2,\mu_{(13)}^2\rangle
\cong G^{free,2}$ well defined) group homomorphism
\begin{eqnarray*}
\nu:\mu(\Gamma(2))\to (2\pi i \Z,+)\quad\textup{with }
\nu(\mu_{(12)}^2)=2\pi i,\ 
\nu(\mu_{(13)}^2)=-2\pi i.
\end{eqnarray*}
Then
\begin{eqnarray*}
\kappa(\mu(B)(\tau)) = \kappa(\tau) + \nu(B)\quad
\textup{for }B\in\Gamma(2),\tau\in \H.
\end{eqnarray*}
Especially
\begin{eqnarray*}
\kappa(\tau-2)&=& \kappa(\tau)+2\pi i,\\
\kappa(\frac{\tau}{2\tau+1})&=& \kappa(\tau)-2\pi i,\\
\kappa(\frac{3\tau-2}{2\tau-1})&=& \kappa(\tau).
\end{eqnarray*}
For $B\in\Gamma(2)$ and $m\in\Z$ 
\begin{eqnarray*}
\kappa(\mu(B)(\tau))&=& \kappa(\tau)+2\pi i m
\quad\textup{if }\mu(B)\textup{ is conjugate in }\Gamma(2) 
\textup{ to }\mu_{(12)}^{2m},\\
\kappa(\mu(B)(\tau))&=& \kappa(\tau)-2\pi i m
\quad\textup{if }\mu(B)\textup{ is conjugate in }\Gamma(2) 
\textup{ to }\mu_{(13)}^{2m},\\
\kappa(\mu(B)(\tau))&=& \kappa(\tau)
\quad\textup{if }\mu(B)\textup{ is conjugate in }\Gamma(2) 
\textup{ to }\mu(\begin{pmatrix}3&-2\\2&-1\end{pmatrix})^{m}.
\end{eqnarray*}

(iii) Consider a cusp $c=\mu(B)(1)$ in the $\mu(\Gamma(2))$
orbit of $1$, so $B\in\Gamma(2)$. If $\tau$ moves along any
geodesic in $\H$ to $c$, then $\kappa(\tau)\to \nu(B)$.

(iv) Finally for $\tau\in\H$ 
\begin{eqnarray*}
\kappa(\frac{2\tau-1}{\tau})&=& -\kappa(\tau),\\
0&=&\pi i-\kappa(\tau-1)+\kappa(\tau)
+\kappa(\frac{\tau-1}{\tau}).
\end{eqnarray*}
\end{theorem}

{\bf Proof:}
(a) and (b) are well known.

(c) (i) $T(\tau)=\exp(\kappa(\tau))$ 
holds by definition of a lift
of the logarithm with $T$. The table on the values of $\kappa$
follows from the table in Theorem \ref{tb.1} on the values of
$T$ and from properties of the logarithm.

(ii) Because of $\mu(\Gamma(2))=
\langle \mu_{(12)}^2,\mu_{(13)}^2 \rangle\cong G^{free,2}$,
it is sufficient to show
\begin{eqnarray*}
\kappa(\tau-2)=\kappa(\tau)+2\pi i\quad\textup{and}\quad
\kappa(\frac{\tau}{2\tau+1})=\kappa(\tau)-2\pi i.
\end{eqnarray*}

First consider a path from $z_0\in \R_{>1}\subset
\C-\{0,1\}$ to $z_0$ which goes once counterclockwise around 0
along a circle with center 0.
The logarithm takes up the summand $2\pi i$. The lift to
a path in $\H$ via $T:\H\to\C-\{0,1\}$, which
starts at the preimage $\tau_0\in 1+i\R_{>0}$ of $z_0$,
ends at $\tau_0-2$. 
Therefore $\kappa(\tau-2)=\kappa(\tau)+2\pi i$.

Now consider a path from $z_1\in (0,1)\subset
\C-\{0,1\}$ to $z_1$ which goes once clockwise around 0
along a circle with center 0
The logarithm takes up the summand $-2\pi i$. The lift to
a path in $\H$ via $T:\H\to\C-\{0,1\}$, which
starts at the preimage $\tau_1\in A(-1,0)$ of $z_1$, 
ends at $\frac{\tau_1}{2\tau_1+1}\in A(1,0)$. 
Therefore $\kappa(\frac{\tau}{2\tau+1})=\kappa(\tau)-2\pi i$.

(iii) The logarithm is univalued near 1 with value 
$\ln(1)=0$. Therefore if $\tau\to 1$ along any geodesic
in $\H$ then $\kappa(\tau)\to 0$. Together with
$\kappa(\mu(B))(\tau)=\kappa(\tau)+\nu(B)$ for 
$B\in\Gamma(2)$, this shows 
$$\kappa(t)\to\nu(B)\quad\textup{if}\quad \tau\to\mu(B)(1).$$

(iv) Observe $\mu(\begin{pmatrix}2&-1\\1&0\end{pmatrix})
=\mu_{(12)}^{-1}\mu_{(123)}$ and 
\begin{eqnarray*}
&&\exp\left(\kappa(\frac{2\tau-1}{\tau})+\kappa(\tau)\right)
= T\left(\mu_{(12)}^{-1}\mu_{(123)}(\tau)\right)
T(\tau)\\
&=& g_{(12)}^{-1}(g_{(123)}(T(\tau))) T(\tau)
=g_{(23)}(T(\tau)) T(\tau) = 1,
\end{eqnarray*}
so $\kappa(\frac{2\tau-1}{\tau})+\kappa(\tau)
=2\pi im$ for some $m\in\Z$. If $\tau\to 1$ along a geodesic
in $\H$ then $\kappa(\frac{2\tau-1}{\tau})+\kappa(\tau)
\to 0+0=0$, so $m=0$.

Similarly
\begin{eqnarray*}
&&\exp\left(\pi i-\kappa(\tau-1)+\kappa(\tau)+ 
\kappa(\frac{\tau-1}{\tau})\right)\\
&=& (-1)T(\mu_{(12)}(\tau))^{-1} T(\tau) 
T(\mu_{(13)}\mu_{(12)}(\tau))\\
&=& (-1)g_{(12)}(T(\tau))^{-1}T(\tau)g_{(13)}g_{(12)}(T(\tau))\\
&=& (-1)\frac{1}{1-T(\tau)} T(\tau) 
\frac{T(\tau)-1}{T(\tau)} = 1,
\end{eqnarray*}
so $\pi i-\kappa(\tau-1)+\kappa(\tau)+ 
\kappa(\frac{\tau-1}{\tau})=2\pi i m$ for some $m\in \Z$.
If $\tau\in 1+i\R_{>0}=A(\infty,1)$, then
\begin{eqnarray*}
\tau-1\in i\R_{>0},\quad\textup{so }\kappa(\tau-1)
\in\pi i+\R,\\
\kappa(\tau)\in \R_{>0},\\
\frac{\tau-1}{\tau}\in A(1,0),\quad\textup{so }
\kappa(\frac{\tau-1}{\tau})\in \R_{<0},\\
\textup{so }\pi i-\kappa(\tau-1)+\kappa(\tau)+ 
\kappa(\frac{\tau-1}{\tau})\in\R,
\end{eqnarray*}
so $m=0$.\hfill$\Box$

\chapter{Quadratic units via continued fractions}\label{sc}
\setcounter{equation}{0}
\setcounter{figure}{0}

\renewcommand{\theequation}{\mbox{C.\arabic{equation}}}

\renewcommand{\thefigure}{\mbox{C.\arabic{figure}}}

The purpose of this appendix is to prove 
Lemma \ref{tc.1} with two statements on the units in 
certain rings of algebraic integers.
The proof of Lemma \ref{tc.1} will be given after the proof
of Theorem \ref{tc.6}.

A convenient and very classical tool to prove this lemma is the
theory of continued fractions as best approximations of 
irrational numbers, applied to the case of quadratic irrationals
which are algebraic integers. 

Theorem \ref{tc.4} below cites standard results on the 
continued fractions of real irrationals. 
It is prepared by the Definitions \ref{tc.2} and \ref{tc.3}

Lemma \ref{tc.5} provides the less well known formulas
\eqref{c.2} and \eqref{c.3} for the case of a quadratic
irrational. Theorem \ref{tc.6} describes the unit group
$\Z[\alpha]^*$ where $\alpha$ is a quadratic irrational
and an algebraic integer, in terms of the continued fractions
of $\alpha$.

\begin{lemma}\label{tc.1}
(a) Let $x\in \Z_{\geq 3}$, and let
$\kappa_{1/2}:=\frac{x}{2}\pm \frac{1}{2}\sqrt{x^2-4}$
be the zeros of the polynomial $t^2-xt+1$, so
$\kappa_1+\kappa_2=x,\kappa_1\kappa_2=1,\kappa_1^2=x\kappa_1-1$.
Then
\begin{eqnarray*}
\Z[\kappa_1]^*=\left\{\begin{array}{ll}
\{\pm\kappa_1^l\,|\, l\in\Z\}&\textup{ if }x\in\Z_{\geq 4},\\
\{\pm (\kappa_1-1)^l\,|\, l\in\Z\}&\textup{ if }x=3.
\end{array}\right. 
\end{eqnarray*}
$\kappa_1$ has norm 1. If $x=3$ then $(\kappa_1-1)^2=\kappa_1$,
and $\kappa_1-1$ has norm $-1$.

(b) Let $x\in\Z_{\geq 2}$, and let 
$\lambda_{1/2} =x^2-1\pm x\sqrt{x^2-2}$ be the zeros of the polynomial
$t^2-(2x^2-2)t+1$, so 
$\lambda_1+\lambda_2=2x^2-2, \lambda_1\lambda_2=1,
\lambda_1^2=(2x^2-2)\lambda_1-1$. Then 
\begin{eqnarray*}
\Z[\sqrt{x^2-2}]^*=\left\{\begin{array}{ll}
\{\pm \lambda_1^l\,|\, l\in\Z\}&\textup{ if }x\geq 3,\\
\{\pm (1+\sqrt{2})^l\,|\, l\in\Z\}&\textup{ if }x=2.
\end{array}\right.
\end{eqnarray*}
$\lambda_1$ has norm 1. If $x=2$ then $(1+\sqrt{2})^2=\lambda_1$,
and $1+\sqrt{2}$ has norm $-1$.
\end{lemma}

Theorem \ref{tc.4} is mainly taken from several theorems in
\cite[1.2 and 1.3]{Ai13}, but with part (b) 
from \cite[I 2.]{Ca65}. It is preceded by two definitions.
According to \cite[5.9 Lagrange's Theorem]{Bu00},
this part (b) is orginally due to Lagrange 1770.
In fact, we will not use this part (b), 
but we find it enlightening.

\begin{definition}\label{tc.2}
Let $\theta\in\R-\Q$ be an irrational number. 
\index{irrational number}

(a) Define recursively sequences $(a_n)_{n\geq 0}$, 
$(\theta_n)_{n\geq 0}$, 
$(p_n)_{n\geq -1},(q_n)_{n\geq -1}$, 
$(r_n)_{n\geq 0}$ as follows:
\begin{eqnarray*}
\theta_0&:=&\theta, \\
a_0&:=& \lfloor \theta_0\rfloor\in\Z,\\
\theta_n&:=& \frac{1}{\theta_{n-1}-a_{n-1}}\in\R_{>1}-\Q
\quad\textup{for }n\in\N,\\
a_n&:=&\lfloor \theta_n\rfloor\in\N\quad\textup{for }n\in\N,\\
(p_{-1},p_0,q_{-1},q_0)&:=& (1,a_0,0,1),\\
p_n&:=& a_np_{n-1}+p_{n-2}\in\Z\quad\textup{for }n\in\N,\\
q_n&:=& a_nq_{n-1}+q_{n-2}\in\N\quad\textup{for }n\in\N,\\
r_n&:=& \frac{p_n}{q_n}\in\Q\quad\textup{for }n\in\Z_{\geq 0}.
\end{eqnarray*}
$\theta_n$ and $a_n$ are defined for all $n\in\N$, because each
$\theta_{n-1}$ is in $\R-\Q$, so $\theta_{n-1}-a_{n-1}\in(0,1)$. 

(b) Following \cite[Notation 2.]{Ca65} define
\begin{eqnarray*}
\|\theta\|:=\min (\theta-\lfloor\theta\rfloor,
\lceil\theta\rceil-\theta)\in(0,\frac{1}{2})
\end{eqnarray*}
\end{definition}

We are interested especially in the case when 
$\theta\in\R-\Q$ is a quadratic irrational.
We recall some notations for this case.

\begin{definition}\label{tc.3}
Let $\theta\in\R-\Q$ be a quadratic irrational, 
\index{quadratic irrational} 
i.e. $\dim_\Q \Q[\theta]=2$. The other root of the minimal
polynomial of $\theta$ is called $\theta^{conj}$, so 
$\theta+\theta^{conj}=:\www{a_0}\in\Q$ and 
$-\theta\theta^{conj}=:d_0\in\Q$.  
For any $\alpha=a+b\theta\in\Q[\theta]$ with $a,b\in\Q$
write $\alpha^{conj}:=a+b\theta^{conj}$. It is the algebraic
conjugate of $\alpha$. The algebra 
homomorphism 
$$\NN:\Q[\theta]\to\Q,\quad \alpha\mapsto \alpha\alpha^{conj},$$
is the norm map. The number $\alpha$ is called
\index{reduced number}
{\it reduced} if $\alpha>1$ and $\alpha^{conj}\in (-1,0)$.
Recall that $\alpha$ is an algebraic integer if and only if
$\alpha+\alpha^{conj}\in\Z$ and $\NN(\alpha)\in\Z$ and that
in this case $\alpha$ is a unit in $\Z[\alpha]$ if and only
if $\NN(\alpha)\in\{\pm 1\}$.  
\end{definition}

\begin{theorem}\label{tc.4} (Classical)
In the situation of Definition \ref{tc.2} the following holds.

(a) \cite[1.2]{Ai13}
$a_0\in\Z$, $a_n\in\N$ for $n\in\N$. For $n\in\Z_{\geq 0}$
the rational number $r_n$ is 
\begin{eqnarray*}
r_n=a_0+\frac{1}{a_1+\frac{1}{\ddots \frac{1}{a_{n-1}+\frac{1}{a_n}}}}
=:[a_0,a_1,...,a_n].
\end{eqnarray*}
It is called {\sf partial quotient} or {\sf continued fraction}
\index{partial quotient}\index{continued fraction}\index{convergent}
or {\sf $n$-th convergent} of $\theta$. 
These numbers approximate $\theta$,
\begin{eqnarray*}
r_0<r_2<r_4<...<\theta<...<r_5<r_3<r_1,\\
|\theta-r_n|<\frac{1}{q_n^2}.
\end{eqnarray*}
This allows to write $\theta=[a_0,a_1,...]$ as an infinite continued fraction.
The numerator $p_n$ and the denominator $q_n$ of $r_n$ are coprime,
\begin{eqnarray*}
\gcd(p_n,q_n)&=&1,\quad\textup{and more precisely }\\
p_nq_{n-1}-p_{n-1}q_n&=&(-1)^{n-1}\textup{ for }n\in\Z_{\geq 0}.
\end{eqnarray*}
The denominators grow strictly from $n=1$ on,
\begin{eqnarray*}
1=q_0\leq q_1<q_2<q_3<... .
\end{eqnarray*}
(b) \cite[I 2.]{Ca65}
The partial quotients $r_n$ are in the following precise sense
the only best approximations of $\theta$: 
\begin{eqnarray*}
|p_n-q_n\theta|&=&\|q_n\theta\|\quad\textup{for }n\in\N,\\
\|q_{n+1}\theta\|&<&\|q_n\theta\|\quad\textup{for }n\in\N,\\
\|q\theta\|&\geq& \|q_n\theta\|\quad\textup{for }n\in\Z_{\geq 0}
\textup{ and }q\in\N\textup{ with }q<q_{n+1},\\
|p_0-q_0\theta|&=&\|q_0\theta\|>\|q_1\theta\|\quad 
\textup{if }q_1>1 (\iff a_1>1),\\
|p_0-q_0\theta|&\in& (\frac{1}{2},1)\textup{ and }
|p_0-q_0\theta|>\|q_1\theta\|\quad\textup{if }q_1=1
(\iff a_1=1).
\end{eqnarray*}
In any case 
\begin{eqnarray*}
|p_{n+1}-q_{n+1}\theta|<|p_n-q_n\theta|\quad\textup{for }n\in\Z_{\geq 0}.
\end{eqnarray*}

(c) \cite[Theorem 1.19]{Ai13}
The partial quotients $r_n$ are also in the following precise sense the
only best approximations of $\theta$:
A rational number $\frac{p}{q}$ with $p\in\Z,q\in\N$ and 
$\gcd(p,q)=1$ satisfies
\begin{eqnarray*}
|\theta-\frac{p}{q}|<\frac{1}{2q^2}&\Longrightarrow&
(p,q)=(p_n,q_n) \quad\textup{for a suitable}\quad n\in\Z_{\geq 0}.
\end{eqnarray*}

(d) \cite[Theorem 1.17 and Proposition 1.18]{Ai13}
The continued fraction is {\sf periodic}, 
\index{periodic continued fraction} 
i.e. there exist
$k_0\in\Z_{\geq 0}$ and $k_1\in\N$ with
\begin{eqnarray*}
a_{n+k_1}=a_n\textup{ for }n\geq k_0,
\end{eqnarray*}
if and only if $\theta$ is a quadratic irrational, i.e.
$\dim_\Q\Q[\theta]=2$. Then one writes
$[a_0a_1...]=[a_0a_1...a_{k_0-1}\oooo{a_{k_0}
a_{k_0+1}...a_{k_0+k_1-1}}]$. 
Furthermore, then $k_0$ can be chosen as 0 
if and only if $\theta$ is reduced. In this case the continued 
fraction $[a_0a_1...]$ is called {\sf purely periodic}.
\index{purely periodic continued fraction}
\end{theorem}

Lemma \ref{tc.5} fixes useful additional 
observations for the case of a quadratic irrational $\theta$. 
These observations are used in the proof of Theorem \ref{tc.6}.
It considers an algebraic integer 
$\alpha\in\R-\Q$ which is a quadratic irrational.
We are interested in the group $\Z[\alpha]^*$ of units in
$\Z[\alpha]$. 
Theorem \ref{tc.6} shows how to see a generator of this group
(and a quarter of its elements) in the continued fractions
of a certain reduced element $\theta$ in $\Z[\alpha]$.
Theorem \ref{tc.6} is not new. 
For example, \cite[Theorem 8.13]{Bu00} gives its main part. 
But the proof here is more elegant than what we found 
in the literature. 

\begin{lemma}\label{tc.5}
Let $\theta\in\R-\Q$ be a quadratic irrational which is
reduced. Let $[\oooo{a_0a_1...a_{k-1}}]$ be its purely periodic
continued fraction of some minimal length $k\in\N$. 
We consider the objects in Definition \ref{tc.2} for this
$\theta$. Then 
\begin{eqnarray}\label{c.1}
\theta_{n+k}=\theta_n\quad\textup{for }n\in\Z_{\geq 0}.
\end{eqnarray}
$\theta_m$ is reduced for $m\in\{0,1,...,k-1\}$, and its 
purely periodic continued fraction is 
$[\oooo{a_m...a_{k-1}a_0...a_{m-1}}]$. 
Write 
\begin{eqnarray*}
\www{a_0}&:=&\theta+\theta^{conj}\in\Q_{>0},\qquad 
d_0:=-\NN(\theta)=-\theta\theta^{conj}\in\Q_{>0},\\
\beta&:=&p_{k-1}-q_{k-1}\theta\in \Q[\theta]-\Q.
\end{eqnarray*}
Then for $n\in\Z_{\geq -1}$
\begin{eqnarray}\label{c.2}
&&p_{n+k}-q_{n+k}\theta= \beta\cdot (p_n-q_n\theta)
\end{eqnarray}
and for $n\in\Z_{\geq 0}$
\begin{eqnarray}
\theta_n= \frac{-p_{n-2}p_{n-1}+p_{n-2}q_{n-1}\www{a_0}
+q_{n-2}q_{n-1}d_0+(-1)^n\theta}{\NN(p_{n-1}-q_{n-1}\theta)}.
\label{c.3}
\end{eqnarray}
\end{lemma}

{\bf Proof:} The natural generalization of the notation
$[a_0,a_1,...,a_m]$ to numbers $a_0\in\R,a_1,...,a_m\in\R_{>0}$
gives for $n\in\Z_{\geq 0}$ 
\begin{eqnarray*}
\theta&=&[a_0,a_1,...,a_{n-1},\theta_n]
=\frac{\theta_np_{n-1}+p_{n-2}}{\theta_nq_{n-1}+q_{n-2}},
\end{eqnarray*}
see Proposition 1.9 in \cite{Ai13}. 
One concludes that the continued fraction of $\theta_n$
is purely periodic, that $\theta_n=\theta_m$ if 
$n=kl+m$ with $l\in\Z_{\geq 0}$ and $m\in\{0,1,...,k-1\}$,
and that its continued fraction is 
$[\oooo{a_m...a_{k-1}a_0...a_{m-1}}]$. Therefore 
$\theta_n$ is reduced. Recall 
$p_{n-1}q_{n-2}-p_{n-2}q_{n-1}=(-1)^n$. Inverting the 
equation above gives 
\begin{eqnarray*}
\theta_n&=&
\frac{\theta q_{n-2}-p_{n-2}}{\theta(-q_{n-1})+p_{n-1}}\\
&=&\frac{(-p_{n-2}+q_{n-2}\theta)(p_{n-1}-q_{n-1}\theta^{conj})}
{\NN(p_{n-1}-q_{n-1}\theta)}\\ 
&=&\frac{-p_{n-2}p_{n-1}+p_{n-2}q_{n-1}\www{a_0}
+q_{n-2}q_{n-1}d_0+(-1)^n\theta}
{\NN(p_{n-1}-q_{n-1}\theta)}.
\end{eqnarray*}
The formula $\theta=\theta_{k}=\frac{\theta q_{k-2}-p_{k-2}}
{\theta(-q_{k-1})+p_{k-1}}$ shows
\begin{eqnarray*}
(1,-\theta)
\begin{pmatrix}p_{k-1}&p_{k-2}\\q_{k-1}&q_{k-2}\end{pmatrix}
= (p_{k-1}-q_{k-1}\theta)(1,-\theta)=\beta(1,-\theta).
\end{eqnarray*}
The inductive definition of $p_n$ and $q_n$ shows
\begin{eqnarray*}
\begin{pmatrix}p_n& p_{n-1}\\q_n& q_{n-1}\end{pmatrix}
=\begin{pmatrix}a_0&1\\1&0\end{pmatrix}
\begin{pmatrix}a_1&1\\1&0\end{pmatrix}...
\begin{pmatrix}a_n&1\\1&0\end{pmatrix}.
\end{eqnarray*}
With the periodicity $a_n+k=a_n$ we obtain
\begin{eqnarray*}
(1,-\theta)
\begin{pmatrix}p_{n+k}&p_{n-1+k}\\q_{n+k}&q_{n-1+k}\end{pmatrix}
&=&(1,-\theta)
\begin{pmatrix}p_{k-1}&p_{k-2}\\q_{k-1}&q_{k-2}\end{pmatrix}
\begin{pmatrix}p_{n}&p_{n-1}\\q_{n}&q_{n-1}\end{pmatrix}\\
&=& \beta(1,-\theta)
\begin{pmatrix}p_n&p_{n-1}\\q_n&q_{n-1}\end{pmatrix}.
\end{eqnarray*}
This gives formula \eqref{c.2}. \hfill$\Box$

\begin{theorem}\label{tc.6}
Let $\alpha\in\R-\Q$ be a quadratic irrational and an
algebraic integer.\index{algebraic integer}
\index{quadratic algebraic integer}

(a) There are a unique sign $\varepsilon_\alpha\in\{\pm 1\}$
and a unique number $n_\alpha\in\Z$ such that 
$\theta:=\varepsilon_\alpha \alpha+n_\alpha$ is reduced.
Then $\Z[\alpha]=\Z[\theta]$,  and any reduced element 
$\www\theta\in\Z[\alpha]$ with $\Z[\alpha]=\Z[\www\theta]$ 
satisfies $\www\theta=\theta$. 
We consider the objects in Definition
\ref{tc.2} for this $\theta$. We define
\begin{eqnarray*}
\www{a_0}&:=&\theta+\theta^{conj}\in\N,\qquad 
d_0:=-\NN(\theta)=-\theta\theta^{conj}\in\N,\\
\beta&:=&p_{k-1}-q_{k-1}\theta\in \Z[\theta]-\Z.
\end{eqnarray*}
as in Lemma \ref{tc.5}. Then $a_0=\www{a_0}$ and 
$d_0\in\{1,2,...,a_0\}$. 

(b) Then $\beta$ is a unit and generates together with $-1$
the unit group $\Z[\alpha]^*$, the $l$-th power of $\beta$ is 
$\beta^l=p_{lk-1}-q_{lk-1}\theta$ for $l\in\Z_{\geq 0}$,
and 
\begin{eqnarray*}
\{\pm \beta^l\,|\, l\in\Z\}&=&\Z[\alpha]^*,\\
\{\beta^l\,|\, l\in\Z_{\geq 0}\}
&=&\Z[\alpha]^*\cap\{p_{n-1}-q_{n-1}\theta\,|\, n\in\Z_{\geq 0}\}.
\end{eqnarray*}
The element $\beta$ is uniquely characterized by the following
properties: 

(i) $-1$ and $\beta$ generate the unit group $\Z[\alpha]^*$, 

(ii) $|\beta|<1$, 

(iii) $\beta=p-q\theta$ with 
$p\in\Z,q\in\N$ (namely $p=p_{k-1},q=q_{k-1}$). 
\end{theorem}

{\bf Proof:}
(a) Choose $\varepsilon_\alpha\in\{\pm 1\}$ such that
$\varepsilon_\alpha(\alpha-\alpha^{conj})>0$.
Then choose $n_\alpha\in\Z$ such that 
$\varepsilon_\alpha\alpha^{conj}+n_\alpha\in(-1,0)$.
Define $\theta:=\varepsilon_\alpha\alpha+n_\alpha$.
Then $\theta^{conj}=\varepsilon_\alpha\alpha^{conj}+n_\alpha\in(-1,0)$.
Also $\theta>\theta^{conj}$ and 
$\theta\theta^{conj}\in\Z-\{0\}$. This shows 
$\theta>1$, so $\theta$ is reduced. 
Also $\Z[\alpha]=\Z[\theta]$ is clear. 

Any reduced element $\www{\theta}\in\Z[\alpha]$ with 
$\Z[\alpha]=\Z[\theta]$ has the shape 
$\www{\theta}=\www{\varepsilon_\alpha}\alpha+\www{n_\alpha}$ with 
$\www{\varepsilon_\alpha}\in\{\pm 1\}$ and $\www{n_\alpha}\in\Z$.
The sign $\www{\varepsilon_\alpha}$ is because of
$\www{\theta}>1>0>\www{\theta}^{conj}$ the unique sign with 
$\www{\varepsilon_\alpha}(\alpha-\alpha^{conj})>0$, so
$\www{\varepsilon_\alpha}=\varepsilon_\alpha$.  
Now $\www{n_\alpha}$ is the unique integer with 
$\varepsilon_\alpha\alpha^{conj}+n_\alpha\in(-1,0)$,
so $\www{n_\alpha}=n_\alpha$. Therefore $\www{\theta}=\theta$.

We have $a_0=\lfloor \theta\rfloor=\theta+\theta^{conj}=\www{a_0}$
and $d_0=-\theta\theta^{conj}\leq a_0$, 
both because $\theta^{conj}\in(-1,0)$. 

(b) We apply Lemma \ref{tc.5}. It tells us which of the 
elements $p_n-q_n\theta$ for $n\in\Z_{\geq 0}$ are units,
in the following way. 

Consider $n\in\Z_{\geq 0}$ and write $n=lk+m$ with 
$l\in\Z_{\geq 0}$ and $m\in\{0,1,...,k-1\}$. 
Recall that $\theta_n$ is reduced, that
$\theta_n=\theta_m$ because of formula \eqref{c.1} 
and that $\theta_m=\theta$ only for
$m=0$, because for $m\in\{1,...,k-1\}$ the purely periodic
continued fractions of $\theta$ and $\theta_m$ differ.
Recall also from formula \eqref{c.2} that 
$$p_{n-1}-q_{n-1}\theta=\beta^l(p_{m-1}-q_{m-1}\theta).$$

If for some $n\in\Z_{\geq 0}$
$\NN(p_{n-1}-q_{n-1}\theta)\in\{\pm 1\}$, then by formula
\eqref{c.3} $\theta_n$ satisfies $\Z[\alpha]=\Z[\theta]
=\Z[\theta_n]$. The uniqueness of $\theta$ in part (a) implies
that then $\theta_n=\theta$, so $m=0$. 
Therefore for $n\in \Z_{\geq 0}-k\Z_{\geq 0}$, 
$\NN(p_{n-1}-q_{n-1}\theta)\notin\{\pm 1\}$,
so then $p_{n-1}-q_{n-1}\theta$ is not a unit. 

On the other hand, if $n=kl$, so $m=0$, then 
$\theta_n=\theta$, and formula \eqref{c.3} tells
$\NN(p_{n-1}-q_{n-1}\theta)=(-1)^n$, so 
$p_{n-1}-q_{n-1}\theta$ is a unit. In fact, formula
\eqref{c.2} tells $p_{n-1}-q_{n-1}\theta=\beta^l$.
We see
\begin{eqnarray}\label{c.4}
\{p_{n-1}-q_{n-1}\theta\,|\, n\in\Z_{\geq 0}\}
\cap \Z[\theta]^* =\{\beta^l\,|\, l\in\Z_{\geq 0}\}.
\end{eqnarray}

It remains to see that $-1$ and $\beta$ generate $\Z[\theta]^*$.

By Dirichlet's unit theorem \cite[Ch. 2 4.3 Theorem 5]{BSh73},
the set $\Z[\alpha]^*$ is as a group isomorphic to
$\{\pm 1\}\times \Z$. It has two generators $\pm\www\beta$
with $|\pm\www\beta|<1$. They are the unique elements in
$\Z[\alpha]^*$ with maximal absolute value $<1$. 
One of them has the shape $p-q\theta$ with $q\in\N$.
This is called $\www\beta$. Then also $p\in\N$,
because $|\www\beta|=|p-q\theta|<1$ and $q\theta>1$. 

{\bf 1st case,} $\theta\in(1,2)$: Then $a_0=d_0=1$, and
$\theta=\frac{1+\sqrt{5}}{2}$ is the golden section with
$\theta^2=\theta+1$. This case is well known. Here
the continued fraction of $\theta$ is purely periodic with 
period $\oooo{1}$ of length one, because $a_0=1$ and 
\begin{eqnarray*}
\theta_1&=&(\theta_0-a_0)^{-1}=\theta=\theta_0.
\end{eqnarray*} 
Here $\beta=1-\theta=-\theta^{-1}=\theta^{conj}$.
It is well known that 
$$\Z[\alpha]^*=\Z[\theta]^*=\{\pm\theta^l\,|\, l\in\Z\}
=\{\pm \beta^l\,|\, l\in\Z\}$$ 

{\bf 2nd case,} $\theta>2$: 
$\www{\beta}=p-q\theta$ is a unit, so $\pm 1=\NN(p-q\theta)$.
Also $|\www{\beta}|<1$ and $\theta>2$ imply $p\geq 2q$. Therefore
\begin{eqnarray*}
|\frac{p}{q}-\theta|&=& \frac{1}{q(p-q\theta^{conj})}
=\frac{1}{q(p+q|\theta^{conj}|)}\\
&<& \frac{1}{q(2q+0)}=\frac{1}{2q^2}.
\end{eqnarray*}
By Theorem \ref{tc.4} (c) $n\in\Z_{\geq 0}$ with 
$(p,q)=(p_n,q_n)$ exists. By \eqref{c.4} 
$\www{\beta}$ is a power of $\beta$, so $\www\beta=\beta$. 
\hfill$\Box$ 

\bigskip
The parts (a) and (b) in the following proof of Lemma \ref{tc.1}
serve also as examples for Theorem \ref{tc.6}.

\bigskip
{\bf Proof of Lemma \ref{tc.1}:}
(a) Here 
\begin{eqnarray*}
x\geq 3\Rightarrow 2x>5\Rightarrow x^2-4&>&x^2-2x+1\\
\Rightarrow \sqrt{x^2-4}&\in& (x-1,x),\\
\Rightarrow \theta=\kappa_1-1
=\frac{x-2}{2}+\frac{1}{2}\sqrt{x^2-4}
&\in& (x-\frac{3}{2},x-1)\\
\textup{and}\quad \theta^{conj}=\kappa_2-1
=\frac{x-2}{2}-\frac{1}{2}\sqrt{x^2-4}&\in& (-1,-\frac{1}{2}).
\end{eqnarray*}
Observe
\begin{eqnarray*}
\theta+\theta^{conj}=x-2,\quad \theta\theta^{conj}=-x+2.
\end{eqnarray*}
Therefore
\begin{eqnarray*}
\theta_0&=&\theta,\quad 
a_0=\lfloor \theta_0\rfloor =x-2,\\
\theta_1&=&(\theta_0-a_0)^{-1}
= \frac{\theta^{conj}-(x-2)}{(\theta-(x-2))
(\theta^{conj}-(x-2))} \\
&=& \frac{-\theta}{-x+2}=\frac{\theta}{x-2},\\
a_1&=&\lfloor \theta_1\rfloor = 1, \\
\theta_2&=& (\theta_1-a_1)^{-1} 
=\frac{x-2}{\theta-(x-2)}
=\frac{(x-2)(\theta^{conj}-(x-2))}{-x+2}=\theta,\\
\theta&=&[\oooo{x-2,1}].
\end{eqnarray*}
The continued fraction of $\theta$ is purely periodic 
with period $\oooo{x-2,1}$ of length 2 if $x\geq 4$ and
purely periodic with period $\oooo{1}$ if $x=3$. 
The norm of $p-q\theta$ for $p,q\in\Z$ is 
\begin{eqnarray*}
\NN(p-q\theta):=(p-q\theta)(p-q\theta^{conj})
=p^2-q(p+´q)(x-2)\in\Z.
\end{eqnarray*}
It is $\pm 1$ if and only if $p-q\theta$ is a unit. 
\begin{eqnarray*}
\begin{array}{lllll}
n & 0 & 1 & 2 & 3 \\
a_n & x-2 & 1 & x-2 & 1 \\
(p_n,q_n) & (x-2,1) & (x-1,1) & (x^2-2x,x-1) & (x^2-x-1,x) \\
\NN(p_n-q_n\theta) & -x+2 & 1 & -x+2 & 1
\end{array}
\end{eqnarray*}

If $x=3$ then $\beta=p_0-q_0\theta=1-\theta$ in the notation of 
Lemma \ref{tc.5}, so $\Z[\theta]^*$ is generated by $(-1)$ and 
$\beta$ or $-\beta^{-1}=\theta=\kappa_1-1$, so 
$$\Z[\theta]^*=\{\pm (1-\theta)^l\,|\, l\in\Z\}
=\{\pm(\kappa_1-1)^l\,|\, l\in\Z\}.$$ 
This is also consistent with the 1st case in the proof of part (b) 
of Theorem \ref{tc.6}. 

If $x\geq 4$ then $\beta=p_1-q_1\theta=x-1-\theta$ in the notation of
Corollary \ref{tc.4}, so $\Z[\theta]^*$ is generated by $(-1)$ and 
$\beta$ or $x-1-\theta^{conj}=\kappa_1$, so 
$$\Z[\theta]^*=\{\pm(x-1-\theta)^l\,|\, l\in\Z\}
=\{\pm\kappa_1^l\,|\, l\in\Z\}.$$

This proves part (a) of Lemma \ref{tc.1}.

(b) The case $x=2$ is treated separately and first.
This case is well known. 
Then $\theta=1+\sqrt{2}\in(2,3)$, 
$\theta^{conj}=1-\sqrt{2}\in (-1,0)$, so $a_0=2$.
The continued fraction of $\theta$ is purely periodic with period
$\oooo{2}$ of length one, because $a_0=2$ and 
\begin{eqnarray*}
\theta_1=(\theta_0-a_0)^{-1}=(\sqrt{2}-1)^{-1}=\theta_0.
\end{eqnarray*}
The element 
$$p_0-q_0\theta=2-\theta=1-\sqrt{2}=\theta^{conj}$$ 
is a unit. This and Theorem \ref{tc.6} (b) show 
$$\Z[\alpha]^*=\Z[\theta]^*=\{\pm\theta^l\,|\, l\in\Z\}
=\{\pm (1+\sqrt{2})^l\,|\, l\in\Z\}.$$

Now we treat the cases $x\geq 3$. Here 
\begin{eqnarray*}
x\geq 3\Rightarrow x^2-2>x^2-x+\frac{1}{4}
\Rightarrow \sqrt{x^2-2}>x-\frac{1}{2},\\
\Rightarrow \theta=(x-1)+\sqrt{x^2-2}
\in(2x-\frac{3}{2},2x-1)\\
\textup{and}\quad \theta^{conj}=(x-1)-\sqrt{x^2-2}
\in(-1,-\frac{1}{2}).
\end{eqnarray*}
Observe
\begin{eqnarray*}
\theta+\theta^{conj}=2x-2,\quad 
\theta\theta^{conj}=-2x+3.
\end{eqnarray*}
Therefore
\begin{eqnarray*}
\theta_0&=&\theta,\quad a_0=\lfloor \theta_0\rfloor=2x-2,\\
\theta_1&=&(\theta_0-a_0)^{-1}=...=
\frac{\theta}{2x-3}\in(1,2),\quad a_1=1,\\
\theta_2&=&(\theta_1-a_1)^{-1}=...=
\frac{\theta-1}{2}\in (x-2,x-1),\quad a_2=x-2,\\
\theta_3&=&(\theta_2-a_2)^{-1}=...=
\frac{\theta-1}{2x-3}\in(1,2),\quad a_3=1,\\
\theta_4&=&(\theta_3-a_3)^{-1}=...=
\theta=\theta_0,\\
\theta&=&[\oooo{2x-2,1,x-2,1}].
\end{eqnarray*}
The continued fraction of $\theta$ is 
purely periodic with period $\oooo{2x-2,1,x-2,1}$ of 
length four. 
The norm of $p-q\theta$ is
$$\NN(p-q\theta)=(p-q\theta)(p-q\theta^{conj})
=p^2+q^2-q(p+q)(2x-2).$$
It is $\pm 1$ if and only if $p-q\theta$ is a unit.
\begin{eqnarray*}
\begin{array}{lllll}
n & 0 & 1 & 2 & 3 \\
a_n & 2x-2 & 1 & x-2 & 1 \\
(p_n,q_n) & (2x-2,1) & (2x-1,1) & (2x^2-3x,x-1) & 
(2x^2-x-1,x) \\
\NN(p_n-q_n\theta) & -2x+3 & 2 & -2x+3 & 1
\end{array}
\end{eqnarray*}
We conclude with Theorem \ref{tc.6} (and with 
the notation of Lemma \ref{tc.5}) that 
$$\beta=p_3-q_3\theta=(2x^2-x-1)-x\theta
=(x^2-1)-x\sqrt{x^2-2}=\lambda_2$$
is together with $(-1)$ a generator of
$\Z[\sqrt{x^2-2}]^*=\Z[\theta]^*$. Therefore also 
$\lambda_1$ together with $(-1)$ is a generator of
$\Z[\sqrt{x^2-2}]^*$. This proves part (b) of Lemma \ref{tc.1}.
\hfill$\Box$ 

\begin{remark}\label{tc.7}
In the situation of Theorem \ref{tc.6}, 
Satz 9.5.2 in \cite{Ko97} tells that the unit group 
$\Z[\theta]^*$ is generated by $-1$ and $q_{k-2}+q_{k-1}\theta$.
This is consistent with Theorem \ref{tc.6} because of the 
following. Here
\begin{eqnarray*}
\theta&=&\frac{\theta_kp_{k-1}+p_{k-2}}{\theta_kq_{k-1}+q_{k-2}}
=\frac{\theta p_{k-1}+p_{k-2}}{\theta q _{k-1}+q_{k-2}},\\
\textup{so}\quad 
0&=& q_{k-1}\theta^2-(p_{k-1}-q_{k-2})\theta-p_{k-2}, \\
\textup{but also}\quad
0&=&\theta^2-\www{a_0}\theta-d_0,\\
\textup{so}\quad 
a_0&=&\www{a_0}=\frac{p_{k-1}-q_{k-2}}{q_{k-1}},\quad
d_0=\frac{p_{k-2}}{q_{k-1}},\\
p_{k-1}-q_{k-1}\theta^{conj} 
&=&p_{k-1}-q_{k-1}(a_0-\theta)
=q_{k-2}+q_{k-1}\theta.
\end{eqnarray*}
\end{remark}

\chapter{Powers of quadratic units}\label{sd}
\setcounter{equation}{0}
\setcounter{figure}{0}

\renewcommand{\theequation}{\mbox{D.\arabic{equation}}}

\renewcommand{\thefigure}{\mbox{D.\arabic{figure}}}

The following definition and lemma treat powers of units of norm 1 
\index{power of a quadratic unit}\index{quadratic unit}
in the rings of integers of quadratic number fields.
Though these powers appear explicitly only in Lemma \ref{td.2}
(c). Lemma \ref{td.2} will be used in the proof of Theorem 
\ref{t5.18}.

\begin{definition}\label{td.1}
(a) Define the polynomials $b_l(a)\in\Z[a]$ for $l\in\Z_{\geq 0}$
by the following recursion.
\begin{eqnarray}
b_0:=0,\quad b_1:=1, \quad b_l:=ab_{l-1}-b_{l-2}
\quad\textup{for }l\in\Z_{\geq 2}.\label{d.1}
\end{eqnarray}

(b) Define for $l\in\Z_{\geq 0}$ the polynomial $r_l\in\Z[a]$ 
and for $l\in\N$ the 
rational functions $q_{0,l},q_{1,l},q_{2,l}\in\Q(t)$,
\index{$q_{0,l},\ q_{1,l},\ q_{2,l}\in\Q(t)$} 
\begin{eqnarray*}
r_0&:=&0,\\
r_l&:=& -ab_l+2b_{l-1}+2\quad\textup{for }l\in\N,\\
q_{0,l}&:=& \frac{b_l-b_{l-1}}{b_l},\\
q_{1,l}&:=& \frac{b_l-b_{l-1}-1}{r_lb_l},\\
q_{2,l}&:=& q_{0,l}-2q_{1,l}.
\end{eqnarray*}

(c) A notation: For two polynomials $f_1,f_2\in\Z[a]$,
$(f_1,f_2)_{\Z[a]}:=\Z[a]f_1+\Z[a]f_2\subset\Z[a]$ denotes
the ideal generated by $f_1$ and $f_2$. 

Remark: If $(f_1,f_2)_{\Z[a]}=\Z[a]$ then for any integer
$c\in\Z$ $\gcd(f_1(c),f_2(c))=1$. 

(d) For $a\in\Z_{\leq -3}\cup\Z_{\geq 3}$ define
$\kappa_a:=\frac{a}{2}+\frac{1}{2}\sqrt{a^2-4}$ and 
$\kappa_a^{conj}:=\frac{a}{2}-\frac{1}{2}\sqrt{a^2-4}$ 
as the zeros of the polynomial $t^2-at+1$, so that 
$\kappa_a+\kappa_a^{conj}=a$,
$\kappa_a\kappa_a^{conj}=1$, $\kappa_a^2=a\kappa_a-1$. 
They are algebraic integers and units with norm 1.
\end{definition}

The following table gives the first twelve of the polynomials
$b_l(a)$. 
The software Maxima \cite{Maxima22} claims that the factors in the
products are irreducible polynomials as polynomials in
$\Q[a]$. We will not use this claim.

\begin{eqnarray*}
b_0 &=& 0\\
b_1 &=& 1\\
b_2 &=& a\\
b_3 &=& (a-1)(a+1)\\
b_4 &=& a(a^2-2)\\
b_5 &=& (a^2-a-1)(a^2+a-1)\\
b_6 &=& (a-1)a(a+1)(a^2-3)\\\
b_7 &=& (a^3-a^2-2a+1)(a^3+a^2-2a-1)\\
b_8 &=& a(a^2-2)(a^4-4a^2+2)\\
b_9 &=& (a-1)(a+1)(a^3-3a-1)(a^3-3a+1)\\
b_{10} &=& a(a^2-a-1)(a^2+a-1)(a^4-5a^2+5)\\
b_{11} &=& (a^5-a^4-4a^3+3a^2+3a-1)(a^5+a^4-4a^3-3a^2+3a+1)
\end{eqnarray*}

\begin{lemma}\label{td.2}
(a) For any $l\in\N$
\begin{eqnarray}
&&1= b_l^2-ab_lb_{l-1}+b_{l-1}^2
=b_l^2-b_{l+1}b_{l-1},\label{d.2}\\
&&(b_{l-1},b_l)_{\Z[a]}=
(b_l-b_{l-1},b_l)_{\Z[a]}=\Z[a],\label{d.3}\\
&&r_l=\left\{\begin{array}{ll}
(2-a)(b_{(l+1)/2}+b_{(l-1)/2})^2&\textup{ for }
l\textup{ odd},\\
(2-a)(a+2)b_{l/2}^2&\textup{ for }l\textup{ even,}\\
 &(\textup{also }l=0)
\end{array}\right. \label{d.4}\\
&&(r_{l-1}/(2-a),r_l/(2-a))_{\Z[a]}=\Z[a],\label{d.5}
\end{eqnarray}
and 
\begin{eqnarray}
q_{2,l}=1-\frac{r_{l-1}/(2-a)}{r_l/(2-a)}.\label{d.6}
\end{eqnarray}

(b) For $a\in\Z_{\leq -3}$
\begin{eqnarray*}
b_l(a)\in(-1)^{l-1}\N\quad\textup{for}\quad l\geq 1,\\
b_1(a)=1,\quad b_2(a)=a,\quad |b_2(a)+b_1(a)|=|a|-1,\\
|b_l(a)|>2|b_{l-1}(a)|\geq |b_{l-1}(a)|+1\quad\textup{for}\quad
l\geq 2,\\
|b_{l+1}(a)+b_{l}(a)|>|b_{l}(a)+b_{l-1}(a)|
\quad\textup{for}\quad l\geq 1.
\end{eqnarray*}
For $a\in\Z_{\geq 3}$
\begin{eqnarray*}
b_l(a)>0\quad\textup{for}\quad l\geq 1,\\
b_1(a)=1,\quad b_2(a)=a,\quad b_2(a)+b_1(a)=a+1,\\
b_l(a)>2b_{l-1}(a)\quad\textup{for}\quad l\geq 1,\\
b_{l+1}(a)+b_l(a)>2(b_l(a)+b_{l-1}(a))
\quad\textup{for}\quad l\geq 1.
\end{eqnarray*}

(c) Consider $a\in\Z_{\leq -3}\cup\Z_{\geq 3}$ and $l\in\N$. Then
\begin{eqnarray}
\kappa_a^l&=& b_l(a)\kappa_a-b_{l-1}(a),\label{d.7}\\
\kappa_a&=& (1-q_{0,l}(a)) + (q_{0,l}(a)-r_l(a)q_{1,l}(a))
\kappa_a^l, \label{d.8}\\
\kappa_a^l&=&\frac{2-r_l(a)}{2}+\frac{1}{2}\sqrt{r_l(a)(r_l(a)-4)},
\label{d.9}
\end{eqnarray}
so $\kappa_a^l$ is a zero of the polynomial
$t^2-(2-r_l(a))t+1$. 
\end{lemma}

{\bf Proof:}
(a) The recursive definition \eqref{d.1} of $b_l$ shows
immediately the equality of the middle and right term
in \eqref{d.2}, it shows
$$b_l^2-ab_lb_{l-1}+b_{l-1}^2
=b_{l-1}^2-ab_{l-1}b_{l-2}+b_{l-2}^2,$$
and it shows $b_1^2-ab_1b_0+b_0^2=1$.
This proves \eqref{d.2}. It implies \eqref{d.3}. 

The sequence $(r_l)_{l\in\N}$ satisfies the recursion
\begin{eqnarray*}
r_0=0,\quad r_1=2-a,\quad 
r_l=ar_{l-1}-r_{l-2}+2(2-a)\quad\textup{for }l\geq 2.
\end{eqnarray*}
For $l=2$ one verifies this immediately. For $l\geq 3$
it follows inductively with \eqref{d.1},
\begin{eqnarray*}
r_l&=&-a(ab_{l-1}-b_{l-2})+2(ab_{l-2}-b_{l-3})+2\\
&=& a(-ab_{l-1}+2b_{l-2}+2)-(-ab_{l-2}+2b_{l-3}+2)+2(2-a)\\
&=&ar_{l-1}-r_{l-2}+2(2-a).
\end{eqnarray*}
For $l=1$ and $l=0$ \eqref{d.4} is obvious. 
For odd $l=2k+1\geq 3$ as well as even $l=2k\geq 2$, 
one verifies \eqref{d.4} inductively with this recursion 
and with \eqref{d.2}, for odd $l=2k+1\geq 3$:
\begin{eqnarray*} 
&&(2-a)(b_{k+1}+b_k)^2-ar_{2k}+r_{2k-1}-2(2-a)\\
&=&(2-a)[ ((ab_k-b_{k-1})+b_k)^2]\\
&&-(2-a)[a(a+2)b_{k}^2-(b_k+b_{k-1})^2+2]\\
&=&(2-a)[2b_k^2-2ab_kb_{k-1}+2b_{k-1}^2-2]\\
&=& 0 \quad (\textup{with }\eqref{d.2}),
\end{eqnarray*}
for even $l=2k\geq 2$:
\begin{eqnarray*}
&& (2-a)(a+2)b_k^2-ar_{2k-1}+r_{2k-2}-2(2-a)\\
&=& (2-a)[(a+2)b_k^2]\\
&&-(2-a)[a(b_k+b_{k-1})^2-(a+2)b_{k-1}^2+2]
\\
&=&(2-a)[2b_k^2-2ab_kb_{k-1}+2b_{k-1}^2-2]\\
&=&0 \quad (\textup{with }\eqref{d.2}).
\end{eqnarray*}
\eqref{d.5} claims for $k\geq 0$ 
\begin{eqnarray*}
((b_{k+1}+b_k)^2,(a+2)b_k^2)_{\Z[a]} =\Z[a]\\
\textup{and}\quad ((a+2)b_{k+1}^2,(b_{k+1}+b_k)^2)_{\Z[a]}=\Z[a].
\end{eqnarray*}
The following {\bf claim} is basic: For $f_1,f_2,f_3\in\Z[a]$
\begin{eqnarray*}
(f_1,f_3)_{\Z[a]}=(f_2,f_3)_{\Z[a]}=\Z[a]\quad\Rightarrow\quad
(f_1f_2,f_3)_{\Z[a]}=\Z[a].
\end{eqnarray*}
To see this claim consider $1=\alpha_1f_1+\alpha_2f_3$,
$1=\beta_1f_2+\beta_2f_3$. Then
\begin{eqnarray*}
1&=& (\alpha_1f_1+\alpha_2f_3)(\beta_1f_2+\beta_2f_3)\\
&=& \alpha_1\beta_1f_1f_2+\alpha_2\beta_2f_3^3
+ \alpha_1\beta_2f_1f_3+\alpha_2\beta_1f_2f_3.
\end{eqnarray*}
The claim and \eqref{d.3} show that for \eqref{d.5} it is
sufficient to prove $$(a+2,b_{k+1}+b_k)_{\Z[a]}=\Z[a].$$
This follows inductively in $k$ with
$$b_{k+1}+b_k=(a+2)b_k-(b_k+b_{k-1})\quad\textup{and}\quad
b_1+b_0=1.$$
Finally, we calculate $q_{2,l}$:
\begin{eqnarray*}
q_{2,l}&=& 
(r_lb_l)^{-1}(r_l(b_l-b_{l-1})-2(b_l-b_{l-1}-1))\\
&=& 1+(r_lb_l)^{-1}(-(-ab_l+2b_{l-1}+2)b_{l-1})-2b_l+2b_{l-1}+2)\\
&=& 1+(r_lb_l)^{-1}(b_l(ab_{l-1}-2)-2(b_{l-1}^2-1))\\
&=& \left\{\begin{array}{ll}
1+(r_l)^{-1}((ab_{l-1}-2)-2b_{l-2})
\quad(\textup{with \eqref{d.2}})&\textup{ for }l\geq 2,\\
1&\textup{ for }l=1\end{array}\right. \\
&=& 1-r_l^{-1}r_{l-1}.
\end{eqnarray*}

(b) All inequalities and signs follow inductively with \eqref{d.1}.

(c) \eqref{d.7} is true for $l=1$. It follows inductively in
$l$ with the following calculation, which uses 
$\kappa_a^2=a\kappa_a-1$. 
\begin{eqnarray*}
\kappa_a^{l+1}&=& \kappa_a(b_l(a)\kappa_a-b_{l-1}(a))\\
&=& b_l(a) (a\kappa_a-1)-b_{l-1}(a)\kappa_a\\
&=& (ab_l(a)-b_{l-1}(a))\kappa_a-b_l(a)\\
&=&b_{l+1}(a)\kappa_a-b_l(a).
\end{eqnarray*}
The right hand side of \eqref{d.8} is
$$\frac{b_{l-1}(a)}{b_l(a)} + \frac{1}{b_l(a)}\kappa_a^l$$
which is $\kappa_a$ by inverting \eqref{d.7}.
Writing $\kappa_a=\frac{a}{2}+\frac{1}{2}\sqrt{a^2-4}$ gives
for $\kappa_a^l$
$$\kappa_a^l=\frac{ab_l(a)-2b_{l-1}(a)}{2}
+\frac{b_l(a)}{2}\sqrt{a^2-4}.$$
One verifies that this equals the right hand side of \eqref{d.9}.
\hfill$\Box$

\end{appendix}


\printindex

\end{document}